\documentclass{gtpart}
 \usepackage{pdfsync,hyperref}
 \usepackage{amssymb}
 \usepackage{latexsym}
 \usepackage{amsmath,amsfonts}
 \usepackage{amscd}
 \usepackage[mathscr]{eucal}

\usepackage{times,tikz}
\usetikzlibrary{patterns}
\def\basepic{
\draw [thick] (0,0) -- (0.2,0) .. controls (0.9,0.0) and (1.1,0.4) .. (1.8,0.4) --
(2.2,0.4) .. controls (2.9,0.4) and (3.1,0) .. (3.8,0) -- (4,0);
\draw [thick] (0,3) -- (0.2,3) .. controls (0.9,3.0) and (1.1,2.6) .. (1.8,2.6) --
(2.2,2.6) .. controls (2.9,2.6) and (3.1,3) .. (3.8,3) -- (4,3);
\draw [thick] (0,1) arc[start angle = -90, end angle=90, x radius=1, y radius=0.5];
\draw [thick] (4,2) arc[start angle = 90, end angle=270, x radius=1, y radius=0.5];
}
\def\blueR{
\draw [line width=0.8pt,densely dotted](4,2.1) arc[start angle = 90, end angle=270, x radius=1.1, y radius=0.6];
}
\def\blueRL{
\draw [line width=0.8pt,densely dotted](3.2929,1.9596) arc[start angle = 130, end angle=230, x radius=1.1, y radius=0.6];
}
\def\blueRR{
\draw [line width=0.8pt,densely dotted](4,2.1) arc[start angle = 90, end angle=130, x radius=1.1, y radius=0.6];
\draw [line width=0.8pt,densely dotted](4,0.9) arc[start angle = 270, end angle=230, x radius=1.1, y radius=0.6];
}
\def\blueL{
\draw [line width=0.8pt,densely dotted](0,0.9) arc[start angle = -90, end angle=90, x radius=1.1, y radius=0.6];
}
\def\blueLR{
\draw [line width=0.8pt,densely dotted](0.7071,1.0404) arc[start angle = -50, end angle=50, x radius=1.1, y radius=0.6];
}
\def\blueLL{
\draw [line width=0.8pt,densely dotted](0,2.1) arc[start angle = 90, end angle=50, x radius=1.1, y radius=0.6];
\draw [line width=0.8pt,densely dotted](0,0.9) arc[start angle = 270, end angle=310, x radius=1.1, y radius=0.6];
}
\def\redM{\draw [dashed] (2,0.4) -- (2,2.6);}
\def\redL{\draw [dashed] (0.7071,0.1) -- (0.7071,1.0404);
  \draw [dashed] (0.7071,1.9696) -- (0.7071,2.9);}
\def\redR{\draw [dashed] (3.2929,0.1) -- (3.2929,1.0404);
  \draw [dashed] (3.2929,1.9696) -- (3.2929,2.9);}
\def\greyM{
  \draw [pattern=north east lines] (1.5,0.4) rectangle (2.5,2.6);
  \draw (2,1.5) node {S};
}
\def\yjpic#1{\mathord{\;\begin{tikzpicture}[scale=0.5,baseline=19]\basepic#1\end{tikzpicture}\;}}
\def\baseOpic{
\draw [thick] (0,0.4) .. controls (1,0.4) and (1,0) .. (2,0)
.. controls (3,0) and (3,0.4) .. (4,0.4);
\draw [thick] (0,2.6) .. controls (1,2.6) and (1,3) .. (2,3)
.. controls (3,3) and (3,2.6) .. (4,2.6);
\draw [thick] (2,1.5) circle [x radius=1, y radius=0.5];
}
\def\Oblue{\draw [line width=0.8pt,densely dotted] (2,1.5) circle [x radius=1.1, y radius=0.6];}
\def\OredM{\draw [dashed] (2,0) -- (2,1); \draw [dashed] (2,2) -- (2,3);}
\def\OredL{\draw [dashed] (0.2,0.4) -- (0.2,2.6);}
\def\OredR{\draw [dashed] (3.8,0.4) -- (3.8,2.6);}
\def\Ogrey{
  \draw [pattern=north east lines] (1.5,0) rectangle (2.5,1);
  \draw [pattern=north east lines] (1.5,2) rectangle (2.5,3);
  \draw (2,0.5) node {S};
}
\def\Opic#1{\mathord{\;\begin{tikzpicture}[scale=0.5,baseline=19]\baseOpic#1\end{tikzpicture}\;}}
\def\OblueL{\draw [line width=0.8pt,densely dotted] (2,2.1) arc
 [start angle=90, end angle=270,x radius=1.1, y radius=0.6];}
\def\OblueR{\draw [line width=0.8pt,densely dotted] (2,0.9) arc
 [start angle=-90, end angle=90,x radius=1.1, y radius=0.6];}

  \DeclareFontFamily{OT1}{pzc}{}
  \DeclareFontShape{OT1}{pzc}{m}{it}{<-> s * [1.050] pzcmi7t}{}
  \DeclareMathAlphabet{\mathpzc}{OT1}{pzc}{m}{it}

 \usepackage{xy}
 \xyoption{all}
\usepackage{graphicx}
\usepackage{calc}

\newlength{\wcwidth}
\newlength{\wcheight}
\newcommand{\widecheck}[1]{\ensuremath{
\settowidth{\wcwidth}{#1}
\settoheight{\wcheight}{#1}
\addtolength{\wcheight}{1pt}
\makebox[0cm][l]{%
\raisebox{\depth+\wcheight}[0cm][0cm]{%
\scalebox{-1}{$\widehat{\hphantom{#1}}$}}}#1
\rule{0pt}{\wcheight+2.5pt}}}

 \renewcommand{\theequation}{{\thesection.\arabic{equation}}}
 \newcommand{\p}{\partial}

\theoremstyle{plain}
\newtheorem*{thm*}{Theorem}
\newtheorem{thm}{Theorem}[section]\newtheorem{prop}[thm]{Proposition}
\newtheorem{lemma}[thm]{Lemma}\newtheorem{cor}[thm]{Corollary}

\theoremstyle{definition}
\newtheorem{defn}[thm]{Definition}

\theoremstyle{remark}

\newtheorem{remarks}[thm]{Remark}

 \newcommand{\bth}{\begin{*thm}}
 \newcommand{\ethm}{\end{*thm}}
 \newcommand{\bco}{\begin{*cor}}
 \newcommand{\eco}{\end{*cor}}
 \newcommand{\bcj}{\begin{*conj}}
 \newcommand{\ecj}{\end{*conj}}
 \newcommand{\bpr}{\begin{*prop}}
 \newcommand{\epr}{\end{*prop}}
 \newcommand{\bprs}{\begin{**prop}}
 \newcommand{\eprs}{\end{**prop}}
 \newcommand{\ble}{\begin{*lemma}}
 \newcommand{\ele}{\end{*lemma}}
 \newcommand{\bles}{\begin{**lemma}}
 \newcommand{\eles}{\end{**lemma}}
 \newcommand{\bsl}{\begin{*sublemma}}
 \newcommand{\esl}{\end{*sublemma}}
 \newcommand{\bre}{\begin{*rem}}
 \newcommand{\ere}{\end{*rem}}
 \newcommand{\bres}{\begin{**rem}}
 \newcommand{\eres}{\end{**rem}}
 \newcommand{\bnt}{\begin{*nota}}
 \newcommand{\ent}{\end{*nota}}
 \newcommand{\bnts}{\begin{**nota}}
 \newcommand{\ents}{\end{**nota}}
 \newcommand{\bde}{\begin{*defn}}
 \newcommand{\ede}{\end{*defn}}
 \newcommand{\bdes}{\begin{**defn}}
 \newcommand{\edes}{\end{**defn}}

\DeclareMathAlphabet\eusm{U}{eus}{m}{n}
\def\makebb#1{\expandafter\def\csname bb#1\endcsname{{\mathbb{#1}}}\ignorespaces}
\def\makerm#1{\expandafter\def\csname rm#1\endcsname{{\rm #1}}\ignorespaces}
\def\makebf#1{\expandafter\def\csname bf#1\endcsname{{\bf #1}}\ignorespaces}
\def\makegr#1{\expandafter\def\csname gr#1\endcsname{{\mathfrak{#1}}}\ignorespaces}
\def\makescr#1{\expandafter\def\csname scr#1\endcsname{{\mathscr{#1}}}\ignorespaces}
\def\makecal#1{\expandafter\def\csname cal#1\endcsname{{\cal #1}}\ignorespaces}
\def\makeudl#1{\expandafter\def\csname udl#1\endcsname{{\underline{#1}}}\ignorespaces}
\def\doLetters#1{%
  #1A #1B #1C #1D #1E #1F #1G #1H #1I #1J #1K #1L #1M
  #1N #1O #1P #1Q #1R #1S #1T #1U #1V #1W #1X #1Y #1Z}
\def\doletters#1{%
  #1a #1b #1c #1d #1e #1f #1g #1h #1i #1j #1k #1l #1m
  #1n #1o #1p #1q #1r #1s #1t #1u #1v #1w #1x #1y #1z}
\doLetters\makebb\doLetters\makecal\doLetters\makerm\doLetters\makebf\doLetters\makescr
\doletters\makerm \doletters\makebf \doLetters\makegr
\doletters\makegr \doLetters\makeudl \doletters\makeudl

 \newcommand{\op}{\operatorname}
\def\co{\colon\thinspace}

\def\HM{\mathop{\textit{HM} } \nolimits}
\def\Hf{\mathop{\textit{HF} }\nolimits}

\def\CM{\mathop{\textit{CM} } \nolimits}

 \def\pf{\noindent{\it Proof.}\enspace}

 \def\Hom{\mathop{\mathrm{Hom\, }}\nolimits}

 \def\dim{\mathop{\mathrm{dim}}\nolimits}
 
\def\smx{\textsc{x}}
\def\smk{\textsc{k}}
\def\smz{\textsc{z}}

\def\scf{{\mbox{\normalsize\textsc{f}}}}
\def\scp{\textsc{p}}

\def\w{\mathrm{w}}
\def\m{\mathrm{m}}
 \def\ff{\mathpzc{f}}
 \def\ss{\mathpzc{s}}
\def\uu{\mathpzc{u}}
\def\mm{\mathpzc{m}}
\def\hata{\op{\hat{a}}}

 \def\gr{\mathop{\mathrm{gr}}\nolimits}

 \def\End{\mathop{\mathrm{End}}\nolimits}
 
 \def\ker{\mathop{\mathrm{Ker}}\nolimits}

\def\Met{\mathop{Met}\nolimits}
  
 \def\pr{\mathop{\mathrm{pr}}\nolimits}

 \def\red{\textrm{red}}
\def\ev{\textrm{even}}
\def\odd{\textrm{odd}}

\def\c{\mathpzc{c}}
\def\e{\mathpzc{e}}
\def\ir{\textit{r}}
\def\r{\mathpzc{r}}
\def\b{\mathpzc{b}}
\def\k{\mathpzc{k}}
\def\q{\mathpzc{q}}
\def\n{\mathpzc{n}}
\def\o{\mathpzc{o}}
\def\p{\mathpzc{p}}
 \def\G{\scriptsize \textsc{G}}

\def\M{\tiny \mathcal{M}}
\def\B{\mathcal{B}}
\def\C{\mathcal{C}}
\def\N{\mathcal{N}}
\def\V{\mathcal{V}}
 \def\dist{\mathop{\mathrm{dist}}\nolimits}
 \def\cl{\mathop{\mathrm{cl}}}

 \def\ps{\frac{\partial}{\partial s}}
 \def\ul{\underline}
\def\td{\tilde}

\def\ulh{\ul{\op{h}}}
 \def\Spin{\mathop{\mathrm{Spin}}\nolimits}

 \newcommand{\epf}{\hfill$\Box$\medbreak}
 
\def\BTitem#1\ETitem{\begin{equation}\hss\left\{\mkern-20mu\parbox{0.9\hsize}{%
\begin{itemize}#1\end{itemize}}\right.\end{equation}}
\def\BTTitem#1\ETTitem{\begin{equation*}\tag{\Ttag}\hss\left\{\mkern-20mu\parbox{0.9\hsize}{%
\begin{itemize}#1\end{itemize}}\right.\end{equation*}}

\mathcode`r="7072

\begin{document}

 \title{ \(\Hf =\HM\) \hskip 4pt  \(V\): \\
Seiberg-Witten-Floer homology and handle additions}

\author{Cagatay Kutluhan} 
 \givenname{Cagatay }
 \surname{Kutluhan}
 \address{Department of Mathematics\\University of Buffalo
}
 \email{kutluhan@buffalo.edu}
 
 \author{Yi-Jen Lee} 
 \givenname{Yi-Jen}
 \surname{Lee}
 \address{Institute of Mathematical Sciences\\the Chinese University
   of Hong Kong
}
 \email{yjlee@ims.cuhk.edu.hk}
 
 \author{Clifford Henry Taubes}
 \givenname{Clifford Henry}
 \surname{Taubes}
 \address{Department of Mathematics\\Harvard University
}
 \email{chtaubes@math.harvard.edu}


\arxivreference{}
\arxivpassword{}

 \begin{abstract}
This is the last of five papers \cite{KLT1, KLT2, KLT3, KLT4} 
that construct an isomorphism between
the Seiberg-Witten Floer homology and the Heegaard Floer homology of a
given compact, oriented 3-manifold. See Theorem \ref{cor:main} below
for a
precise statement. As outlined in \cite{KLT1}, this isomorphism is given as a
composition of three isomorphisms.
In this article, we establish the third isomorphism that relates the
Seiberg-Witten Floer homology on the auxiliary manifold with the
appropriate version of Seiberg-Witten Floer homology on the original
manifold. This constitutes Theorem 4.1 in \cite{KLT1}, re-stated in a
more refined form as Theorem \ref{thm:main} below.
The tool used in the proof is a filtered variant of connected sum
formula for Seiberg-Witten Floer homology, in special cases where one
of the summand manifolds is \(S^1\times S^2\) (referred to as
``handle-addition'' in all five articles in this series). Nevertheless, the
arguments leading to the afore-mentioned connected sum formula
are general enough to establish a connected sum formula in the
wider context of Seiberg-Witten Floer homology 
with non-balanced perturbations. This is stated as Proposition
\ref{prop:conn-1} below. Although what is asserted in this proposition
has been known to experts for some time, a detailed proof has not
appeared in the literature, and therefore of some independent interest.
 \end{abstract}

\maketitle

\tableofcontents

\section{Introduction}
To summarize what was done in the predecessors to this article,
\cite{KLT1}-\cite{KLT4} cited above: the first article in this series
outlined a program for a proof of Theorem \ref{cor:main}, based on a
concatenation of three isomorphisms. 
The first isomorphism (Theorem 2.3 in \cite{KLT1}) relates a
version of embedded contact homology on an auxiliary manifold to
the Heegaard Floer homology on the original, and was accomplished in
\cite{KLT2}-\cite{KLT3}. 
The second isomorphism (Theorem 3.4 in \cite{KLT1})
relates the relevant version of the embedded contact homology on the
auxiliary manifold and a version of the Seiberg-Witten Floer
homology on this same manifold. This was established in \cite{KLT4}. 
This last installment of the \(\HM =\Hf\) series contains the proof of
the third isomorphism, stated as Theorem 4.1 in \cite{KLT1}. Part of the content of this paper are
drawn from unpublished details of the proof of
the second author's Corollary 8.4 in \cite{L}, which describes the
behavior of certain Seiberg-Witten Floer homology under handle addition. 

\subsection{The main theorem and an outline of proof}\label{sec:1.1}

Let $M$ be a closed, connected, and oriented $3$-manifold. 
Given a $\Spin^c$ structure $\grs$ on $M$,
P. B. Kronheimer and T. S. Mrowka defined in \cite{KM} three
flavors of Seiberg-Witten Floer homology, 
\(\widehat{\HM}_*\), \(\overline{\HM}_*\), and \(\widecheck{HM}_*\), 
modelling on three different
versions of \(S^1\)-equivariant homologies.
These homology groups have the structure of modules over the graded ring 
\[
\mathbf{ A}_\dag(M):= \mathbb{Z}[U]\otimes\textstyle\bigwedge^*(H_1(M;\mathbb{Z})/\op{Tors}),
\]
where \(U\) has degree \(-2\) and elements in
\(H_1(M;\mathbb{Z})/\op{Tors}\) has degree \(-1\). These modules are
graded by an affine space over \(\bbZ/c_\grs\bbZ\), where
\(c_\grs\in 2 \bbZ^{\geq 0}\) is the divisibility of \(c_1(\grs)\),
the first Chern class of the $\Spin^c$ structure $\grs$.
Moreover, as \(\mathbf{ A}_\dag(M)\)-modules, these three flavor of
Seiberg-Witten Floer homologies fit into a long exact
sequence modelling on the fundamental exact sequence of
\(S^1\)-equivariant Floer homologies. (Cf. Equation (3.4) in
\cite{KM}). 
\begin{equation}\label{HM-sq1}
\cdots\widehat{\HM}\to \overline{\HM}\to \widecheck{HM}\to \cdots
\end{equation}
This is called {\em the first fundamental exact sequence
  of \(\HM\)} in this article. In \cite{L},
the second author defined a fourth flavor of Seiberg-Witten Floer
homology \(\widetilde{\HM}_*\) with the same module structure and relative
grading. (It was originally denoted by \(\op{HM}^{tot}\) in
\cite{L}; given here as Definition \ref{def:HM-tilde}). The definition
models on the ordinary homology of an
\(S^1\)-space.  As such, it fits into a second long exact sequence
together with \(U\, \widehat{\HM}_*\) and \(\widehat{\HM}_*\). This is
referred to as the {\em second fundamental exact sequence of \(\HM\)};
see Lemma \ref{lem:1.11} below. 

In this article, we regard these four flavors of \(\HM\) as a system, in the
order of \(\widehat{\HM}_*\), \(\overline{\HM}_*\), and
\(\widecheck{HM}_*\), \(\widetilde{\HM}_*\). They are denoted
collectively by \(\mathring{\HM}_*\).

 As will be detailed in the upcoming Section \ref{sec:2}, the
 Seiberg-Witten Floer homologies (also referred to as the monopole
 Floer homology in this article)
 \(\mathring{\HM}_*\) depends on the cohomology class of the
 perturbation form \(\varpi\) in addition to the \(\Spin^c\) structure \(\grs\).
 One may also define a monopole Floer homology with 
 local coefficients \(\Gamma\) compatible with \(\grs\) and \([\varpi]\).
Of particular interest to us is the case when the perturbation is
``balanced''; in this case \(\Gamma\) may be taken to be \(\bbZ\).
 These are denoted by $\mathring{\HM}_* (M,\grs,c_b)$; and this is the
 variant of monopole Floer homology to be equated with the Heegaard
 Floer homology \(\Hf^\circ_*\), in Theorem \ref{cor:main}
 below. This is, in a sense, the strongest possible statement of
 equivalence between \(\HM\) and \(\Hf\), as the monopole Floer
 homology \(\overline{\HM}\neq 0\) and \(\widehat{\HM}\neq
 \widecheck{HM}\) only in the balanced case. The equivalence between
 other versions of \(\HM\) and \(\Hf\) may be deduced from this case
 through the use of local coefficients. It is also worth mentioning
 that a coarser version of Seiberg-Witten-Floer homology,
 \(\HM_\bullet\), defined by taking a completion of the Floer complex with respect to
 grading, frequently appears in \cite{KM} and other literature. In
 this article we work exclusively with the original version, \(\mathring{\HM}_*\).

The upcoming Theorem \ref{thm:main} relates $\mathring{\HM}_*
(M,\grs,c_b)$ with two filtered variants of monopole Floer homology.
The first was introduced in \cite{L}, originally denoted by
\(\op{HMT}^\circ\) therein. Here, the label \(\circ\) stands, in
specific order, for \(-,\infty, +, \wedge\). The fact that they appear
in the superscript (instead of the top) of the notation, and the order
in which they appear, reflects the nature of their definition. The
latter is done following the algebraic framework of Ozsvath-Szabo in
\cite{OS1}. The second of these two variants was introduced in
\cite{KLT4} (Cf. also Section 4 of \cite{KLT1} for a brief
summary). They are denoted by \(\textrm{H}_*^\circ(Y)\) in
\cite{KLT1}, and by \(\textrm{H}^\circ_{SW}\) in \cite{KLT4}. The
construction of both these filtered monopole Floer homologies is based
on the same general framework, which we describe in Section
\ref{sec:3} below. This framework always produces four flavors of
Floer homologies, labeled by \(\circ=-,\infty, +, \wedge\); and they are
related by two fundamental long exact sequences parallel to those
appearing in the Heegaard Floer theory,  cf. Equation
(\ref{eq:hf-fund-sq}) below. (To be more precise,
only the first three flavors appeared in \cite{KLT1} and \cite{KLT4},
but it shall become clear in Section \ref{sec:3} that the
afore-mentioned general construction actually gives rise to a fourth
flavor). The basic ingredient of this construction consists
of a triple of data: A certain \(\Spin^c\) 3-manifold \(Y_Z\), a
closed 2-form \(w\) on \(Y_Z\) used to define a monotone perturbation
to the Seiberg-Witten equations, and a special 1-cycle \(\gamma\)
embedded in \(Y_Z\) useful for defining a filtration on the
associated monopole Floer complex. Further constraints on the choice of this
triple are given in Section \ref{sec:Aa)}. 

The triple
that enters the definition of \(\op{HMT}^\circ\) is what was denoted
by \((\ul{M}, *d\ul{\ff}, \ul{\gamma})\) in \cite{L}. Here, \(\ul{M}\) is
constructed from \(M\) by adding a 1-handle
\footnote{Cf. Item (7) in Section \ref{sec:notation}}
along the extrema of
\(\ff\), the latter being a Morse function giving rise to the Heegaard
diagram used to define \(\Hf^\circ\). Denote this 1-handle by
\(\mathcal{H}_0\). What was denoted by \(\ul{\ff}\)
is an \(S^1\)-valued harmonic Morse function obtained by a natural extension of
\(\ff\). The 1-circle \(\ul{\gamma}\) is related to the path
\(\gamma_z\subset M\) used by Ozsvath-Szabo to define a filtration on
the Heegaard Floer complex. The triple used for the
definition of \(\textrm{H}^\circ_{SW}\) in \cite{KLT4} was denoted by
\((Y, w, \gamma^{(z_0)})\) in \cite{KLT4} and \cite{KLT2}. The
3-manifold \(Y\) is obtained from \(\ul{M}\) by attaching additional
1-handles along pairs of index 1 and index 2 critical points of \(\ff\).
The two-form \(w\) on \(Y\) is constructed from a natural extension of
\(*d\ff\). The 1-cycle \(\ul{\gamma}\) in \(\ul{M}\) becomes the
1-cycle \(\gamma^{(z_0)}\) in \(Y\) after the handle-attachment. The
precise definitions of \(\op{HMT}^\circ\) and
\(\textrm{H}^\circ_{SW}\) may be found in Section \ref{sec:3.8}. By
construction, \(\op{HMT}^\circ\) and
\(\textrm{H}^\circ_{SW}\) are respectively \(\mathbf{A}_\dag
(\ul{M})\)- and \(\mathbf{A}_\dag(Y)\)-modules, and each is equipped
with a pair of fundamental exact sequences parallel to (\ref{eq:hf-fund-sq}).

Let \(\textsc{g}\) denote the number of 1-handles added to \(\ul{M}\) in order
to obtain \(Y\) and denote these handles by \(\mathcal{H}_\grp\),
\(\grp\in \Lambda\), where the label set 
\(\Lambda\) is an ordered set consisting of \(\G\) elements. Recall
that \(\gamma_z\subset M\) is defined so that \(\partial\gamma _z\)
is the attaching 0-cycle of \(\mathcal{H}_0\). As described in
\cite{L}, the path \(\gamma_z\) determines a decomposition of \(\ul{M}\)
as a connected sum \(\ul{M}\simeq M\# (S^1\times S^2)\) (cf. \cite{L}
Equation (15)), and hence a splitting 
\begin{equation}\label{eq:H_Y-split0}
H_1(\ul{M};\bbZ)\simeq H_1(M;\bbZ)\oplus H_1(S^1\times S^2;\bbZ),
\end{equation}
with the second summand generated by \([\ul{\gamma}]\in
H_1(\ul{M};\bbZ)\). 
Correspondingly, this determines a factorization of the algebra \begin{equation}\label{A-mod-1}
\mathbf{A}_\dag (\ul{M})\simeq \mathbf{A}_\dag
(M)\otimes _{\bbZ[U]}\mathbf{A}_\dag(S^1\times S^2)=\mathbf{A}_\dag
(M)\otimes \textstyle\bigwedge^*H_1(S^1\times S^2;\bbZ).
\end{equation}
The last factor above, \(\textstyle\bigwedge^*H_1(S^1)=\textstyle\bigwedge^*H_1(S^1\times
S^2;\bbZ)\) has a natural action on its dual algebra
\(\textstyle\bigwedge^*(H^1(S^1))\). The latter is regarded as a
graded \(\bbZ\)-algebra generated by two
elements, one of degree 0 and the other of degree \(1\). This was denoted by \(\hat{V}\) in
\cite{KLT1} and by \(H_*(S^1)\) in the rest of this article. (Cf. Item (6)
of Section \ref{sec:notation} below). For this reason we shall use the shorthand
\(H_{-*}(S^1)\) for the factor \(\textstyle\bigwedge^*H_1(S^1\times
S^2;\bbZ)\) in (\ref{A-mod-1}), and the aforementioned dual action is
implied whenever we refer to ``the
\(H_{-*}(S^1)\) action on \(H_{*}(S^1)\)'' below. 

The auxiliary manifold \(Y\) may be decomposed as a connected sum of
\(\ul{M}\) and \(\G\) copies of \(S^1\times S^2\), one for each of the
1-handles \(\mathcal{H}_\grp\), in a similar manner: 
For each \(\grp\in \Lambda\), we fix an arc \(\lambda_\grp\) in \(M\) connecting the
attaching 0-cycle of  \(\mathcal{H}_\grp\). Let \(S_\grp\) denote the
boundary sphere of a small tubular neighborhood of \(\lambda_\grp\), and
use the same notation for the corresponding sphere in \(Y\). The precise description of \(\lambda_\grp\)
and \(S_\grp\) is given in Part
1 of Section \ref{sec:De)} below. Now split
\(Y\) along these spheres \(S_\grp\) to get the aforementioned
connected sum, and use this to define a splitting 
\begin{equation}\label{eq:H_Y-split1}
H_1(Y;\bbZ)\simeq H_1(\ul{M};\bbZ)\oplus\bigoplus _{\grp\in \Lambda}H_1((S^1\times S^2)_\grp;\bbZ),
\end{equation}
where \((S^1\times S^2)_\grp\) denotes the copy of \(S^1\times S^2\)
coming from \(\mathcal{H}_\grp\). 
This in turn determines a factorization \begin{equation}\label{A-mod-2}\mathbf{A}_\dag (Y)\simeq \mathbf{A}_\dag
(\ul{M})\otimes _{\bbZ[U]}\bigotimes_{\grp\in \Lambda}\mathbf{A}_\dag((S^1\times S^2)_\grp)=\mathbf{A}_\dag
(\ul{M})\otimes H_{-*}(S^1)^{\G}
\end{equation}
like (\ref{A-mod-1}).


The main theorem of this article relates the three versions of
monopole Floer homologies: \(\mathring{\HM}(M, \grs, c_b)\),
\(\op{HMT}^\circ\), and \(\mathrm{H}^\circ_{SW}=\op{H}^\circ(Y)\).

\begin{thm}\label{thm:main}

{\bf (1)} Use \(\op{HMT}^\circ\boxtimes
H_{*}(S^1)^{\boxtimes\G}\) to denote the exterior tensor product of
the \(\mathbf{A}_\dag(\ul{M})\)-module \(\op{HMT}^\circ\) and 
\(\G\) copies of the \(H_{-*}(S^1)\)-module \(H_*(S^1)\). With respect to the factorization (\ref{A-mod-2}), there exists a system of 
isomorphisms of 
\(\mathbf{A}_\dag(Y)=\mathbf{A}_\dag(\ul{M})\otimes
H_{-*}(S^1)^{\otimes \G}\)-modules 
\[
\op{H}^\circ(Y)\stackrel{\simeq}{\longrightarrow}\op{HMT}^\circ\boxtimes
H_{-*}(S^1)^{\boxtimes\G}, \quad \circ=-, \infty, +, \wedge \]
which 
  preserves the relative
gradings and is natural with respect to the fundamental long exact
sequences on both sides.

{\bf (2)} The \(H_{-*}(S^1)\)-factor of the factorization 
\(\mathbf{A}_\dag(\ul{M})=\mathbf{A}_\dag(M)\otimes H_{-*}(S^1)\)
in (\ref{A-mod-1}) acts trivially on \(\op{HMT}\).  
Regarding \(\op{HMT}\) as an \(\mathbf{A}_\dag(M)\)-module in this
manner, 
there exists a system of 
 isomorphisms of
\(\mathbf{A}_\dag(M)\)-modules from
\[\text{\(\op{HMT}^\circ\), \(\circ=-, \infty, +, \wedge\) respectively to
\(\mathring{\HM} \,( M, \grs, c_b)\), \(\circ=\wedge, -, \vee, \sim\),}\]
that preserves the relative
gradings and is natural with respect to the fundamental long exact
sequences on both sides. 
 \end{thm}

 The proof of this theorem is given in Section
\ref{sec:filtered-conn}. The remainder of this section gives a brief outline of this proof.

Given how \(Y\) is constructed from \(\ul{M}\), and \(\ul{M}\) in turn
from \(M\), it is little surprise that the preceding theorem is a
consequence of certain filtered variant of connected sum
formula for Seiberg-Witten-Floer homologies. See Propositions \ref{prop:conn-f}
in Section \ref{sec:filtered-conn}. The first steps of the proof of this
formula, via understanding the chain maps on Seiberg-Witten Floer
complexes induced by cobordisms
associated to the connected sum, lead to a connected sum formula
for Seiberg-Witten Floer homologies {\em sans filtration}. This is
stated as Proposition \ref{prop:conn-1} below. 

The more essential part of the proof, which also constitutes the major
technical component of this article, consists of an extension of the
framework defining \(\op{HMT}^\circ\) and \(\mathrm{H}^\circ(Y)\) to
the context of cobordisms and their assciated chain maps. The
analytical foundation of such an extension is provided in Sections
\ref{sec:B}-\ref{sec:D} of this article.

The proof of part (2) of Theorem \ref{thm:main} also involves some
homological algebra computation that turns out to be a 
manifestation of so-called ``Koszul duality''.  An elementary account
of the relevant part of this story is given in Section \ref{sec:4}.
This algebraic machinery expresses all four flavors of the balanced
monopole Floer homology, \(\mathring{\HM}(M, \grs, c_b)\) in terms of
a balanced monopole Floer complex of the first flavor, \(\widehat{\CM}_*
(M,\grs, c_b)\). Meanwhile, the filtered connected sum formula
previously mentioned expresses all four flavors of \(\op{HMT}^\circ\)
in terms of a monopole Floer complex with ``negative monotone''
perturbation, \(\CM_*(M,\grs, c_-)\). See Proposition \ref{prop:KM-ES} below.
These two monopole Floer
complexes are linked via a chain-level variant of the 
following result of Kronheimer-Mrowka's:
\begin{thm}[\cite{KM} Theorem 31.5.1]
Suppose \(c_1(\grs)\) is not torsion. Then \[\widehat{\HM}_*(M, \grs, c_b)\simeq\HM_* (M, \grs, c_-).\]
\end{thm}
The right hand side of the preceding isomorphism refers to the
monopole Floer homology for negative monotone perturbations. A brief
account of this variant of monopole Floer homology can be found in
Section \ref{sec:2.3}. The construction of both \(\op{H}^\circ_{SW}\)
and \(\op{HMT}^\circ\) are based on negative monotone monopole Floer
complexes.

More on the motivation for various constructions in the article may be
found in \cite{L}.

\begin{remarks}\label{rmk:suture}
With the hind-sight gained from Juhasz's \cite{Ju} and
Kronheimer-Mrowka's \cite{KMs} definitions of sutured Floer homologies, we feel that \(\op{HMT}^\circ\) are best interpreted as variants of sutured Floer
homology. In particular, \(\HM\,( M(1), \grs(1))=\widehat{\op{HMT}}
(\underline{M}, \ul{\grs})\) in terms of the notation in \cite{Ju}, \cite{KMs}
and \cite{L}. 
From this point of view, Theorem \ref{thm:main} (2) may be viewed as a
re-interpretation of monopole Floer homology of closed 3-manifolds as
(generalized) sutured Floer homology. In particular, the
\(\circ =\wedge\) variant of this statement is a Seiberg-Witten
analog of Prop 2.2 in \cite{Ju}, where the hat-version of the Heegaard
Floer homology is re-interpreted as a sutured Floer homology. See also
Theorem 1.6 announced in \cite{CGHH} for an \(ECH\) analog (of the
\(\circ=\wedge\) variant). We hope to discuss this in more detail
elsewhere. (See also the end of Part 4 in Section \ref{sec:Da)}).
\end{remarks}

\subsection{Relating Heegaard and Seiberg-Witten Floer homologies}

With all said and done, the main result of this articles combines with
those in \cite{KLT1}-\cite{KLT4} to reach our ultimate goal:
\begin{thm} \label{cor:main}
Let \(M\) be a closed, oriented 3-manifold, and \(\grs\) be a
\(\Spin^c\) structure on \(M\). Then there exists a system of
isomorphisms from \(\Hf^\circ_*(M, \grs)\), \(\circ=-, \infty, +,
\wedge\), respectively to \(\mathring{\HM}_*(M, \grs, c_b)\),
\(\circ=\wedge, -, \vee, \sim\),  as \(\bbZ/c_\grs\bbZ\)-graded 
\(\mathbf{A}_\dag(M)\)-modules, which is natural with respect to the
fundamental exact sequences of the Heegaard and monopole Floer
homologies.   
\end{thm}

The result summaries the relation between the Heegaard and
monopole Floer homlogies, which has been conjectured since the inception
of Heegaard Floer theory. See for example Conjecture 1.1 in \cite{OS2}
, I.3.12 in \cite{KM}, Conjecture 1 in \cite{KMan}, and Conjecture 1.1 in \cite{L}.

As the Heegaard Floer homology \(\Hf^\circ\) makes no other
appearances for the rest of this article, the reader is referred to
\cite{OS1} and \cite{OS2} for its definition and properties. In
particular, the fundamental exact sequences relating its four flavors
take the following form:
\begin{equation}\label{eq:hf-fund-sq}\begin{split}
& \cdots \to \Hf^-\to \Hf^\infty\to \Hf^+\to \cdots\\
&\cdots \to \Hf^-\stackrel{U}{\to }\Hf^-\to \widehat{\Hf}\to \cdots
\end{split}
\end{equation}
\paragraph{\it Proof of Theorem \ref{cor:main}.} An outline of the
proof is already given in \cite{KLT1}. To summarize, by combining the
two parts of Theorem \ref{thm:main}, one has (Cf. Theorem 4.1 in \cite{KLT1}):
\begin{equation}\label{compute:H(Y)}
\op{H}^\circ(Y)\simeq\mathring{\HM} (M, \grs, c_b)\boxtimes
H_{-*}(S^1)^{\boxtimes\G}, 
\end{equation}
as modules over the algebra \(\mathbf{A}_\dag(M)\otimes
H_{-*}(S^1)^{\otimes\G}\). Here, the \(\mathbf{A}_\dag(M)\otimes
H_{-*}(S^1)^{\otimes\G}\) structure on \(\op{H}^\circ(Y)\) comes from
the latter's \(\mathbf{A}_\dag(Y)\)-module structure via the
inclusion \begin{equation}\label{eq:A-inclusion}\mathbf{A}_\dag(M)\otimes
1 \otimes H_{-*}(S^1)^{\otimes\G}\hookrightarrow\mathbf{A}_\dag(M)\otimes
H_{-*}(S^1)\otimes H_{-*}(S^1)^{\otimes\G}\stackrel{i_{sum}}{\simeq}\mathbf{A}_\dag(Y)
\end{equation}
with respect to the factorization combining 
(\ref{A-mod-1}) and (\ref{A-mod-2}). 

It is asserted in Theorem 3.4 of \cite{KLT1} and proven in Theorem 1.5
of 
\cite{KLT4} that the left hand side of (\ref{compute:H(Y)}),
\(\op{H}^\circ(Y)\), is isomorphic to what was called
``\(ech^\circ\)'' as \(\mathbf{A}_\dag(Y)\)-modules. The
\(ech^\circ\) chain complex, as well as a (particular choice of)
\(\mathbf{A}_\dag(Y)\) action on it, is explicitly described in
\cite{KLT2}-\cite{KLT3}. A computation based on this explicit
description yields: 
\begin{prop}\label{prop:ech-compute} (Cf. also Theorem 2.4 of \cite{KLT1})
There is a system of isomorphisms 
\[ech^\circ\simeq\Hf^\circ (M, \grs)\boxtimes
H_{*}(S^1)^{\boxtimes\G}\]
as modules over  \(\mathbf{A}_\dag(M)\otimes
H_{-*}(S^1)^{\otimes\G}\), which preserves relative gradings and is
natural with respect to the fundamental exact sequences on both sides.
Here, the \(\mathbf{A}_\dag(M)\otimes
H_{-*}(S^1)^{\otimes\G}\) structure on \(ech^\circ \) also refers to
the one induced from
the latter's \(\mathbf{A}_\dag(Y)\)-module structure via the same 
inclusion (\ref{eq:A-inclusion}).
\end{prop}

The proof of this proposition involves some details of \cite{KLT3}'s
description 
of the \(\mathbf{A}_\dag(Y)\) actions on \(ech^\circ\), as well as
some 
particular choice of the arcs \(\lambda_{\grp}\) used to define the
factorization (\ref{A-mod-2}), and will be postponed to Section
\ref{pf:prop1.5}.

The assertion of the theorem is a direct consequence of the
composition of the three isomorphisms from the preceding proposition,
(\ref{compute:H(Y)}), and Theorem 1.5 of \cite{KLT4} (Theorem 3.4 of \cite{KLT1}).
\epf

\subsection{Some notations and conventions}\label{sec:notation}

 Throughout the remainder of this
paper, section numbers, equation numbers, and other references from
\cite{KLT1}-\cite{KLT4} are distinguished from those in this paper by the use of the
appropriate Roman numeral as a prefix.  For example,
`Section II.1' refers to Section 1 in \cite{KLT2}.  
In addition, the following conventions are used:
\begin{enumerate}
\item  As in \cite{KLT1}-\cite{KLT4}, we use
\(c_{0}\) to denote a constant in $(1, \infty)$ whose value
is independent of all relevant parameters.  The value of
\(c_{0}\) can increase between subsequent appearances.  
\item  As in \cite{KLT1}-\cite{KLT4}, we denote by $\chi$ a fixed,
non-increasing function on $\mathbb{R}$ that equals \(1\) on a neighborhood
of $(-\infty, 0]$ and equals \(0\) on a neighborhood of $[1, \infty)$.
\item When left unspecified, the modules, chain complexes and homologies in this
  article are over the coefficient ring \(\bbK\), which can be taken
  to be \(\bbZ\) as was done in \cite{KLT1}-\cite{KLT4}. Using a
  separate notion serves to distinguish different roles
  the abelian group \(\bbZ\) plays in this article, e.g. as the group of
  deck transformations versus the coefficient ring of the chain
  complexes. 
\item The term ``module'' in this article refers to either a left module or
a right module.  Thus, both the monopole Floer homology and monopole
Floer cohomology are said to have a module structure over the ring
\(H^*(BS^1)\). Note in contrast that in \cite{KM}, a ``module'' refers
specifically to a left module.
Moreover, what appears as \(U_\dag\) in \cite{KM} is denoted by \(U\)
in this article for simplicity, since we focus on Floer homology as
opposed to cohomology.
\item The definition of Floer complexes in this article often depends
  on several parameters, yet there are chain homotopies relating the
  Floer complexes with the values of some of the parameters
  changed. In the interest of simplicity, these parameters are usually
  left unspecified in our notation for the Floer complexes unless necessary.
\item Due to geometric motivations (cf. \cite{GKM}), we
view \(H_*(S^1)\) and \(H^*(BS^1)\) both as free commutative
differential graded algebras with zero differential and a single
generator, where the odd 
generator \(y\) for \(H_1(S^1)\) has degree 1, while the even generator
\(u\) for \(H^*(BS^1)\) has degree \(-2\). In this section
commutativity and the commutator \([\cdot, \cdot]\) are meant in the
{\em graded} sense. In particular, what is called an ``anti-chain map''
in \cite{KM} is in our terminology an odd chain map. If necessary, we use 
notations \([\cdot, \cdot]_{\text{odd}}\) or \([\cdot,
\cdot]_{\text{even}}\) to emphasize the parity of the commutator.  When \(H_*(S^1)\)
is written as a polynomial algebra in \(y\), \(\bbZ[y]\),
\(H_{-*}(S^1)\) is often written as \(\bbZ[\partial_y]\), to reflect
the action of \(H_{-*}(S^1)\) on \(H_*(S^1)\). 
\item In this article as well as its prequels, a ``1-handle''
  frequently refer to \([0,1]\times S^2\), and ``attaching a 1-handle to a
  3-manifold'' refers to a 0-dimensional surgery on the
  3-manifold.
\item[(8)] In the context of fiber bundles over a fixed base manifold,
  \(\ul{\bbF}\) typical stands for a trivial bundle with fibers
  \(\bbF\). 
\end{enumerate}

\paragraph{Acknowledgments.} The authors are supported in part by
grants from the National Science Foundation. The first author is
supported by a National Science Foundation Postdoctoral Research
Fellowship under Award No.~DMS-1103795. The second author
thanks Harvard University for hosting her during multiple
visits through the long course of working on this project. She also
wish to thank T. Mrowka and P. Ozsvath for suggesting the general form
of  the connected sum formula, Proposition \ref{prop:conn-1}
below. A similar statement under different assumptions, via a
different and more involved route of proof, is to appear in \cite{BMO}.

\section{Elements of Seiberg-Witten Floer theory}\label{sec:2}
\setcounter{equation}{0}

This subsection reviews some backgrounds on Seiberg-Witten-Floer
theory, with the book \cite{KM} as the definitive
reference. By way of this, we introduce some notation and
terminlogy used in the rest of this article, some of which differ from
those in \cite{KM}. We focus mostly on the special cases involved in
the proof of Theorem \ref{thm:main}, leaving the general details
for the reader to consult \cite{KM}. Many notions here have analogs in, e.g. \cite{KLT4},
\cite{LT}, which work with similar settings.

\subsection{Seiberg-Witten equations on 3-manifolds}\label{sec:2.1}  
Let \(M\) be a
closed, oriented, Riemannian 3-manifold. Fix a \(\Spin^c\)-structure \(\grs\)
on \(M\) and let \(\bbS\) denote its associated spinor bundle. 
We call a pair, \((\bbA, \Psi )\), consisting of a Hermitian connection on \(\det\,
(\bbS)\) and a section of \(\bbS\) a (Seiberg-Witten) {\em configuration}. 
The gauge group \(C^{\infty }(M; U (1))\) acts on the space of configurations
in the following fashion:  Let \(\hat{u}: M \to  U (1)\).  Then \(\hat{u}\) sends
a configuration, \((\bbA, \Psi )\), to \((\bbA - 2\hat{u}^{-1}d\hat{u}, \hat{u}\Psi )\).   Two solutions
obtained one from the other in this manner are said to be
\textit{gauge equivalent}.  Note that this \(C^{\infty }(M; U (1))\) action
is free except at pairs of the form \((\bbA, \Psi=0)\); these are
called {\em reducible} configurations. Configurations which are not
reducible is {\em irreducible}.

In the most general form, the 3-dimensional Seiberg-Witten equations ask that  a configuration
\((\bbA, \Psi )\) obey
\begin{equation}\label{eq:SW-3d}
\begin{cases}
B_{\bbA } - \Psi^{\dag }\tau\Psi + i\varpi-\grT= 0  & \textit{  and    }\\
D^{\bbA}\Psi  -\grS= 0 , & \end{cases}
\end{equation}
where \(B_\bbA\) denotes the Hodge dual of the curvature form of
\(\bbA\), \(D^\bbA\) denotes the Dirac operator, and the quadratic
term \(\Psi^{\dag }\tau\Psi \) is as in Section 1.2 of
\cite{LT}. \(\varpi\) is a closed 2-form, and the pair \((\grT,
\grS)\) is
a small perturbation arising as the formal gradient of a gauge-invariant function of \((\bbA, \Psi
)\). This is called a tame perturbation in \cite{KM}, and is in general needed
to guarantee the transversality properties necessary for the
definition of Seiberg-Witten-Floer homology. See Chapters 10 and 11
in \cite{KM}. In the simplest case, \((\grT, \grS)\) may be taken to
be of the form 
\begin{equation}\label{eq:grTS}
(\grT, \grS)=(2i*d\mu, 0) 
\end{equation}
for a smooth 1-form \(\mu\) taken from
a Banach space called \(\Omega\) in \cite{KLT4}. This may be assumed
to be a subspace of the Banach space of tame perturbations in
Chapter 11.6 in \cite{KM}, and hence inherits the 
so-called ``\(\mathcal{P}\)-norm'' from \cite{KM}. This norm
bounds the norms of the derivatives of \(\mu\) to any given order.

Irreducible solutions to (\ref{eq:SW-3d})
may exist only when the cohomology class
\([\varpi]=2\pi c_1(\det\, \bbS)\). In this case the Seiberg-Witten equations
(\ref{eq:SW-3d}) is said to have {\em balanced perturbation}, while it
is said to have {\em exact perturbation} when \([\varpi]=0\). The
cases when \([\varpi]=2rc_1(\det\, \bbS)\) is said to be
{\em monotone}: when \(r>\pi\) it is said to be {\em negative
  monotone}, and when \(r<\pi\) it is said to be {\em positive
  monotone}. Note that when \(c_1(\det\,\bbS)\) is torsion, the notation
of balanced, exact, and positive or negative monotone perturbations
are equivalent. We work in the negative monotone case with nontorsion
\(c_1(\det\,\bbS)\) for most part of
this paper where all Seiberg-Witten solutions are irreducible. 
Note in contrast that in the closely related series of
articles \cite{T1}-\cite{T5}, \(\varpi\) is taken to be \(da\) for
a contact 1-form \(a\), which is an exact perturbation.

This said, unless otherwise specified, from now on we set 
\begin{equation}\label{eq:2rw}
\varpi=2rw
\end{equation}
for a closed 2-form \(w\) in the cohomology class of \(c_1(\det\, \bbS)\) and a
real number \(r>\pi\). When \(c_1(\det\, \bbS)\) is torsion, we always
set \(w\equiv 0\). Otherwise, the particulars of \(w\) for the proof
of our main theorem \ref{thm:main} are described in Section \ref{sec:Aa)}.

To make contact with the notation in \cite{KLT4}, write 
\begin{equation}\label{def:E-3d}
\det\,(\bbS)=E^{2}\otimes K^{-1}
\end{equation}
with \(K\to M\) being a fixed
complex line bundle. Fix a smooth connection, \(A_K\), on \(K^{-1}\).
Where \(w\) is nowhere vanishing (such as over the stable Hamiltonian
manifold \(Y\) in \cite{KLT4}), \(K^{-1}\) is typically given by \(\ker\, 
(*w)\subset TM\) 
and \(E\) the \(i|w|\)-eigen-bundle of the Clifford action by
\(w\). More constraints on the choice of \(K\) and \(A_K\) will be
specified along the way through the rest of this article. 

With \(A_K\) chosen, let \(A\) denote the connection on
the \(E\)-summand corresponding to \(\bbA\), and write
\(\Psi=\sqrt{2r}\, \psi\). In this case, perturbations of the form
(\ref{eq:grTS}) suffice for our purpose. Since 
the Riemannian metric and a connection on \(E\) determine a \(\Spin^c\)
connection on \(\bbS\), 
we often consider the equivalent equations for \((A,
\psi)\) of the form
\begin{equation}\label{eq:(A.4)}
\begin{cases}
B_{A} - r\, (\psi^{\dag}\tau\psi -i*w) +   \frac{1}{2} B_{A_K} - i *d\mu = 0, & \\
D_{A}\psi = 0, &
\end{cases}
\end{equation}
where 
\(D_A=D^\bbA\), \(B_{A}\) is the Hodge star of the
curvature 2-form of \(A\) and \(B_{A_K}\)
 denotes the Hodge star of the curvature 2-form for the connection,
\(A_{K}\).

Given a Hermitian line bundle \(V\to M\), we use \(\op{Conn}\, (V)\) to
denote the space of Hermitian connections on \(V\). 
The equations in (\ref{eq:SW-3d}) 
are the variational equations of the functional \(\gra\)
of \((A, \psi)\in \op{Conn}\, (E) \times C^{\infty}(M;\bbS)\), 
given by
\begin{equation}\label{eq:(A.5)}
\gra =  \frac{1}{2} \grc\grs - r \textsc{w} + \gre_{\mu} + r\int _{M}\psi^\dag
D_A\psi,
\end{equation}
where the notation is as follows:  The functions \(\grc\grs\) and
\(\textsc{w}\) are defined using a chosen reference 
connection on \(E\).  Let \(A_{E}\) denote the latter. 
With \(A\) written as \[
A=A_{E}+\hata_{A}, 
\]
then \(\textsc{w}\) and \(\grc\grs\) are given by
\begin{equation}\label{eq:(A.6)}
\textsc{w} = i\int_{M}\hata_A\wedge w\quad 
\text{and}\quad
\grc\grs = -\int_{M}\hata_A\wedge d\hata_A - 2\int_{Y_Z}\hata_A\wedge \, \big(F_{A_E}+\frac{1}{2}F_{A_K}\big).
\end{equation}
What is denoted by \(\gre_{\mu}\) is the integral over
\(M\) of  \(i\mu \wedge F_{A}\).  The functionals \(\gra\),
\(\textsc{w}\) and \(\grc\grs\) in general is not invariant under the
\(C^{\infty}(M;U(1))\) action on \(\op{Conn}\, (\det\,\bbS)\times
C^{\infty}(M;\bbS)\), however their differentials descend to
the orbit space. These differentials are henceforth denoted by 
\(d\gra\), \(d(\grc\grs)\), \(\cdots\), etc. 

To define the Seiberg-Witten Floer homology in general, \cite{KM} takes a
real blow-up of the space \(\op{Conn}\, (\det\,\bbS)\times
C^{\infty}(M;\bbS)=:\mathcal{C}(M)\) along the set of reducibles (Cf. Chapter 6 of \cite{KM}). This blown-up space
is denoted as \(\mathcal{C}^\sigma(M, \grs)\) therein and has a free
\(C^\infty(M, U(1))\)-action. (Cf. \cite{KM} p.115). The vector field
dual to \(d\gra\) extends to
\(\mathcal{C}^\sigma\),  which is then used
to define the Seiberg-Witten equations. (Cf. \cite{KM} Section 6.2).

A solution \(\grc\) to the Seiberg-Witten equations or its
corresponding gauge equivalence class \([\grc]\in \mathcal{C}^\sigma(M,
\grs)/C^\infty(M, U(1))\) is said to be {\em non-degenerate} when certain differential operator 
\(\grL_\grc\) has trivial kernel. The explicit form of this operator is given for
irreducible solutions of (\ref{eq:SW-3d}) in
(\ref{eq:(B.30)}) below. In general, this notion of nodegeneracy
arises from the interpretation of \([\grc]\) as a zero of the 1-form \(d\gra\) on \(\mathcal{C}^\sigma(M,
\grs)/C^\infty(M, U(1))=:\mathcal{B}^\sigma(M)\). With the metric and \(\varpi\) fixed, a choice of \((\grT, \grS)\), (or in the case of
(\ref{eq:(A.4)})), of \(\mu\)) such that all solutions to (\ref{eq:SW-3d})
or (\ref{eq:(A.4)})  are non-degenerate is said in what follows to be
 {\em suitable.} 
In the negative monotone case with nontorsion \(c_1(\det\,\bbS)\), a suitable choice
for \(\mu\) can be found with $\mathcal{P}$-norm bounded by any
given positive number.  (Cf. e.g. (1.18) in \cite{KLT4} and references
therein). Otherwise, especially when reducible solutions exist, a
suitable pair \((\grT, \grS)\) is typically of more general form than
that of (\ref{eq:grTS}). Nondegenerate gauge-equivalence classes of reducible
Seiberg-Witten solutions are further classified into the ``stable''
and ``unstable'' types in \cite{KM}.

\subsection{Seiberg-Witten equations on 4-dimensional cobordisms}\label{sec:2.2}
Let \(Y_-\), \(Y_+\) be closed oriented 3-manifolds. In this paper
\(X\) will denote a simple cobordism from \(Y_-\) to \(Y_+\) of
the following sort: \(X\) is 
an oriented complete 4-manifold equipped with the extra structure
listed below.
\BTitem\label{(A.9a,11)}
\item 
 There is a proper function \(s\co X \to \mathbb{R}\)
 with non-degenerate critical points with at most one single critical
 value, \(0\).
\item  
 There exists an orientation preserving diffeomorphism between
the \(s < 0\)  part of \(X\)  and \((-\infty, 0)\times Y_{-}\)  that identifies \(s\)  with
the Eulidean coordinate on the \((-\infty, 0)\)  factor.
\item  
 There exists an orientation preserving diffeomorphism between
the \(s >0\)  part of \(X\)  and \((0, \infty)\times Y_{+}\)  that identifies \(s\)  with
the Eulidean coordinate on the \((0, \infty)\) factor. 
\item There is an even class in \(H^{2}(X;\bbZ)\)  that restricts to
  the \(s < 0\) and \(s > 0\)  parts of \(X\)  as the respective
  \(Y_{-}\)  and \(Y_{+}\)  versions of
\(c_{1}(\det\, (\bbS))\).
\ETitem
The diffeomorphism in the second bullet of (\ref{(A.9a,11)}) is used, often
implicitly, to identify the \(s < 0\) part of \(X\) with
\((-\infty, 0) \times Y_{-}\); and the
diffeomorphism in the third bullet of (\ref{(A.9a,11)}) is likewise used to
identify the \(s > 0\) part with \((0, \infty) \times
Y_{+}\). Fix a class satisfying the last bullet of (\ref{(A.9a,11)})
and denote it also by \(c_{1}(\det\, (\bbS))\).

Assume that the Riemannian metric on \(X\) satisfies the following: 
\BTitem\label{eq:(A.12,15a)}
\item There exists \(L \geq 100\)  such that the metric on the
\(s \leq -L \) and \(s \geq L \) parts 
of \(X\) are identified by the embeddings in the second and third
bullets of (\ref{(A.9a,11)}) with the respective product metrics on \((-\infty,
-L] \times Y_{-}\)  and \([L, \infty)\times Y_{+}\).
\item  The metric pulls back from the \(|s | 
\in[L  - 8, L]\)  part of \(X\)  via the embeddings from
the second and third bullets of (\ref{(A.9a,11)}) as the quadratic form
\(ds^{2} + \grg\)  with \(\grg\)  being an \(s\)-dependent metric on either \(Y_{-}\) 
or \(Y_{+}\)  as the case may be.
\ETitem
The chosen metric on \(X\) is used to write \(\bigwedge^{2}
T^*X\) as \(\Lambda^{+} \oplus\Lambda^{-}\) with \(\Lambda^{+}\)
denoting the bundle of self-dual 2-forms and with
\(\Lambda^{-}\) denoting the corresponding bundle of
anti-self dual 2-forms.  A given 2-form \(\grw\) is written with respect to
this splitting as \(\grw = \grw^{+} + \grw^{-}\). 
 
Use the metric to define the notion of a
$\Spin^c$-structure on \(X\).
 It follows from the last bullet in (\ref{(A.9a,11)}) that there is a
$\Spin^c$ structure that restricts to the \(s
\leq -2\) and \(s \geq 2\) parts of \(X\) as the given
$\Spin^c$ structures from \(Y_{-}\)
and \(Y_{+}\), and has its first Chern class equal to \(c_{1}(\det\, (\bbS))\).  Fix such a
$\Spin^c$ structure and use \(\bbS^{+}\) and \(\bbS^{-}\)
to denote the respective bundles of self-dual and anti-self dual
spinors.  

The Seiberg-Witten equations on \(X\) are equations for a
pair \((\bbA, \Psi)\) with \(\bbA\) being a Hermitian
connection on the line bundle \(\det\, (\bbS^{+})\) and
with \(\Psi\) being a section of \(\bbS^{+}\).  It takes the following
general form:
\begin{equation}\label{eq:(A.14)}
F_{\bbA}^{+} -(\Psi^{\dag}\tau\Psi- i \varpi_{X}) -\grT^+=0
\quad \text{and}\quad \mathcal{D}_{\bbA}^+\Psi-\grS^+=0,
\end{equation}
where the notation uses $F_{\mathbb{A}}$ to denote the
curvature 2-form of \(\bbA\), and it uses \(\Psi^{\dag}\tau\Psi\) to denote the bilinear
map from $\mathbb{S}^{+}$ to $i\Lambda^{+}$ that is defined using the Clifford
multiplication.  Meanwhile, $\mathcal{D}_{\mathbb{A}}^+\co \Gamma
(\bbS^+)\to \Gamma(\bbS^-)$ and $\mathcal{D}_{\mathbb{A}}^-\co \Gamma
(\bbS^-)\to \Gamma(\bbS^+)$ are the 4-dimensional Dirac operators on \(X\) defined
by the metric and the chosen connection \(\bbA\).  What is denoted by 
\(\varpi_X\) is a self-dual 2-form satisfying the following list for
some \(L'\geq L\):
 \BTitem\label{eq:(A.13a)}
\item  
 The pull-back of \(\varpi_{X}\)  from the \(s
< -L' \) part of \(X\)  via the embedding from the
second bullet of (\ref{(A.9a,11)}) is twice the self dual part
of a closed 2-form \(\varpi_-\) on \(Y_-\).  
\item The pull-back of \(\varpi_{X}\) from the \(s >
L' \) part of \(X\)  via the embedding from the third bullet of (\ref{(A.9a,11)}) is twice the self dual part of a closed
2-form \(\varpi_+\) on \(Y_{+}\). 
\ETitem
 The pair \((\grT^+, \grS^+)\) is the 4-dimensional analog of \((\grT,
 \grS)\) in (\ref{eq:SW-3d}); see (24.2) in \cite{KM}. 

We denote \(X_c:=s^{-1}([-L'-1, L'+1])\subset  X\) and call it 
the ``compact piece'' of \(X\). Each connected component of \(X-X_c\)
is called an {\em end} of \(X\). The diffeomorphisms in
(\ref{(A.9a,11)}) identifies each end with a product \((-\infty, -L')\times
M\) or \(M\times (L', \infty)\) for a connected oriented manifold
\(M\); in the first case it is said to be an {\em negative end}, and
in the second case an {\em positive end}. In either case we call this
end the {\em \(M\)-end} of \(X\). ``The negative end of \(X\)'' refers
to \(s^{-1}(-\infty, -L'-1)\simeq (-\infty, -L'-1)\times Y_-\), and ``the
positive end of \(X\)'' refers to \(s^{-1}(L'+1, \infty)\simeq (L'+1,
\infty)\times Y_+\).
\paragraph{\it Caveat.} 
What is denoted as \(X_c\) in this article was denoted by \(X\) in
\cite{KM}. Correspondingly, the noncompact manifold \(X\) in this
article was denoted as \(X^*\) in \cite{KM}.

An important
 special case is when (\ref{eq:(A.14)}) is defined on a {\em product
   cobordism}. By this we mean that \(X=\bbR\times M\) for a closed
 oriented \(\Spin^c\) 3-manifold \(M\), with the function \(s\) as the
 Euclidean coordinate of the \(\bbR\) factor; the Riemannian metric on \(X\) is
 the product of the affine metric on \(\bbR\) and the Riemannian
 metric on \(M\), and both \(\varpi_X\) and \((\grT^+,
 \grS^+)\) are invariant under the natural \(\bbR\)-action on
 \(\bbR\times M\). 
Thus, the conditions in the first bullet of (\ref{eq:(A.12,15a)}) and in (\ref{eq:(A.13a)}) may be
paraphrased as saying that the \(s^{-1}[L', \infty)\) and
\(s^{-1}(-\infty, -L']\) part of the Seiberg-Witten equations on \(X\) are
those of product cobordisms. As explained in \cite{KM}, Clifford action by \(ds\) over product
cobordisms may be used to identify
\(\bbS^+\simeq \bbS^-\). Meanwhile, both are the pull-back of a spinor
bundle \(\bbS\) over \(M\).  In this way, (\ref{eq:(A.14)}) may be
re-written as a gradient flow equation of the action functional
\(\gra\), cf. (IV.1.20). The gradient vector field here is \((-1)\)
times the left-hand side of (\ref{eq:SW-3d}), with (\ref{eq:(A.14)})'s 
\(\varpi_X=2\varpi^+\), and \(\grT^+\), \(\grS^+\) induced
respectively from the \(\grT\) and \(\grS\) in (\ref{eq:SW-3d}).

A solution \(\grd = (\bbA, \Psi)\) to (\ref{eq:(A.14)}) is said to be an
{\em  instanton} if the constant \(s \leq -L\) pull-backs
converge as \(s \to -\infty \) to a pair that can be written
as \((\bbA_{-},\Psi_{-})\), with \((\bbA_{-},\Psi_{-})\) being a
solution to (\ref{eq:SW-3d}) on
\(Y_{-}\); and if the constant \(s \geq L\) pull-backs converge
as \(s \to \infty \) to a pair \((\bbA_{+}, \Psi_{+})\), with \((\bbA_{+}, \Psi_{+})\) being a solution to
(\ref{eq:SW-3d}) on \(Y_{+}\).  If \(\grd\) is an instanton then the convention
in what follows will be to say that the respective \(s \to-\infty \) and \(s \to \infty \) limits of \(\grd\) are
\((\bbA_{-},  \Psi_{-})\) and \((\bbA_{+},  \Psi_{+})\). As in 
the 3-dimensional case, \cite{KM} define a real ``blow-up'' of the space
\(\mathcal{C}_{loc}(X):=\op{Conn}\, \,( \det\bbS^+)\times C^\infty(X, \bbS^+)\), this denoted by
\(\mathcal{C}^\sigma_{loc}(X)\) below. To describe
\(\mathcal{C}^\sigma_{loc}(X)\) in more detail, consider the tautological
bundle \(C^\infty(X, \bbS^+)-\{0\}\) over the sphere \(\bbU (C^\infty(X, \bbS^+)):=(C^\infty(X, \bbS^+)-\{0\})/\bbR^+\), and
let \(\Gamma^\sigma(X;\bbS^+)\) denote the \(\bbR^{\geq 0}\)-bundle
associated to this principal \(\bbR^+\)-bundle. Then
\(\mathcal{C}_{loc}^\sigma (X):=\op{Conn}\, \,( \det\bbS^+)\times
\Gamma ^\sigma (X, \bbS^+)\). Alternatively,
\[\mathcal{C}_{loc}^\sigma
(X)=\bigcap_{l\in\bbZ^+}\mathcal{C}_{l,loc}^\sigma (X),\] 
where \(\mathcal{C}_{l,loc}^\sigma (X)\) is the \(L^2_{l,loc}\)-variant of
\(\mathcal{C}_{loc}^\sigma (X)\) defined in p.464 of 
\cite{KM}. We denote an element in \(\Gamma ^\sigma (X, \bbS^+)\) in
the form of 
 \(\Psi^\sigma =(\pmb{\Psi}, \Psi)\), where \(\pmb{\Psi}\in
 \bbU (C ^\infty(X, \bbS^+))\), and \(\Psi\) is 
 in the fiber of the bundle \(\Gamma^\sigma(X;\bbS^+)\) over \(\pmb{\Psi}\).

The 4-dimensional Seiberg-Witten equations
(\ref{eq:(A.14)}) may be generalized to elements in
\(\mathcal{C}^\sigma_{loc}(X)\), and hence also the notion of
an instanton. (Cf. \cite{KM} Equation (6.5).) An (generalized) instanton has its \(s\to
-\infty\) and \(s\to\infty\) limits in \(\mathcal{C}^\sigma(Y_-)\) and
\(\mathcal{C}^\sigma(Y_+)\) respectively, in the sense explained in
p.219 in \cite{KM}. The 4-dimensional Seiberg-Witten equation is
invariant under the actions of the gauge group \(C^\infty(X; U(1))\). An instanton \((\bbA,
\Psi^\sigma )\), or a gauge-equivalence class of instantons,  is said to be {\em reducible} when \(\Psi\equiv 0\); otherwise
it is {\em irreducible}. 
 
The perturbation \((\grT^+, \grS^+)\) is introduced in (\ref{eq:(A.14)}) so that a certain operator
that is associated to any given instanton solutions to (\ref{eq:(A.14)}) is
Fredholm with trivial cokernel.  Cf. Chapter 24.3 of \cite{KM} in
general and Equation (1.21) in \cite{KLT4} for a special case closely related to this
article. Instanton solutions with this
property are said to be non-degenerate.  
We call perturbation term {\em suitable} when all
instanton solutions to the corresponding version of (\ref{eq:(A.14)}) are
non-degenerate.  A suitable perturbation can be found for
(\ref{eq:(A.14)}) with norm bounded by any given positive number.  The relevant
norm is also called the $\mathcal{P}$-norm.  As in the case with
elements in \(\Omega\), the $\mathcal{P}$-norm of a perturbation term
bounds the norms of its derivatives to all orders.  

Just as in the 3-dimensional case, 
the 4-dimensional cobordisms relevant to Theorem \ref{thm:main} are
equipped with  \(\varpi_X\) and pairs \((\grT^+, \grS^+)\) of the following form:
\[\label{eq:(A.14b)}
\varpi_X=2rw_X\quad \text{and}\quad  (\grT^+, \grS^+)=(i \grw_{\mu}^{+}, 0).
\] 
Here, \(\grw_{\mu}\) is a 2-form of the form \(d\big(\chi(L + s) \mu_{-} +
\chi(L  -s) \mu_{+}\big)\) for some 1-forms \(\mu_-\), \(\mu_+\) on
\(Y_-\), \(Y_+\) respectively. However, in the case of a product
cobordism \(X = \bbR \times M\), we take \(\grw_{\mu} = d\mu_{-}
=d\mu_{+}\). Meanwhile, \(w_X\) is a self-dual 2-form constrained by the
properties listed in (\ref{eq:(A.13b)}) below, among others. 
These constraints involve another constant denoted by \(L_{tor}\) below. The latter is
no smaller than \(L+4\). The constraints use \(X_{tor}\) to denote the
union of the components of the \(|s|>0\) part of \(X\) where \(c_1(\det\,\bbS)\)
is torsion. 
\BTitem\label{eq:(A.13b)}
\item  
 The pull-back of \(w_{X}\)  to each
constant \(s\)  slice of \(X\)  is a closed 2-form whose de
Rham cohomology class is that of
\(c_{1}(\det\, (\bbS))\) .
\item  
 The embedding from the first bullet of
(\ref{eq:(A.13a)}) pulls back \(w_{X}\)  from the \(s
< -L \) part of \(X-X_{tor}\)  as twice the self dual part of the
\(Y_{-}\)  version of the 2-form \(w\).  The embedding from the second
bullet of (\ref{eq:(A.13a)}) pulls back \(w_{X}\) from the \(s >
L \) part of \(X-X_{tor}\) as twice the self dual part of the
\(Y_{-}\)  version of the 2-form \(w\). The 2-form \(w_X\)
is identically zero on any component of the \(|s|>L_{tor}\) part of \(X_{tor}\).
\ETitem

Similarly to the 3-dimensional case, the 4-dimensional Seiberg-Witten
equations may be rewritten in terms of the pair \((A, \psi)\in
\op{Conn}\, (E)\times C^\infty(\bbS^+)\) that is obtained from the pair \((\bbA,
\Psi)\in \op{Conn}\, (\det\bbS^+)\times C^\infty(\bbS^+)\) via the same formulae
as those in the previous subsection. This requires an extension of
\(K\) and \(A_K\) from the ends \(s^{-1}[L', \infty)\cup
s^{-1}(-\infty, -L']\). Constraints on such choices will be introduced
in subsequent sections as needs arise; typically where \(\varpi_X\)
is nowhere vanishing, \(E\) is chosen to be the \(i|\varpi_X|\)-eigenbundle
under the Clifford action of \(\varpi_X\) on \(\bbS^+\).

\subsection{The monopole Floer chain complex}\label{sec:2.3}

Fix a closed, oriented, connected Riemannian 3-manifold \(M\) and a \(\Spin^c\)
structure \(\grs\) on it. We first give in Part 1 below a precise definition of the
monopole Floer complexes involved in the proof of the Theorem
\ref{cor:main}, the main objective of this series of articles. 
Sketches of how they  generalize to other cases are provided in Parts 2 and 3 of
this subsection. 

\paragraph{Part 1: Non-torsion \(c_1(\grs)\), positive/negative
  monotone \(\varpi\).}  Suppose for now that \(\grs\) has non-torsion first
Chern class, and \(\varpi\), \((\grT, \grS)\) are as in (\ref{eq:2rw})
and (\ref{eq:grTS}) respectively, with \(r>\pi\). Fix also a complex Hermitian line
bundle $K \to M$ as specified in Section \ref{sec:2.1} above. 
 The {\em spectral flow} function on \(\op{Conn}\,  (E) \times
C^{\infty}(M; \bbS)\) is defined initially on the complement of a certain codimension 1
subvariety just as in Section 1e in \cite{KLT4} using a chosen Hermitian connection
on \(E\) and a suitably generic section of \(\bbS\).  As such it is
locally constant and integer-valued.  The definition can be extended
to the whole of \(\op{Conn}\,  (E) \times C^{\infty}(M;
\bbS)\) as explained in Sections \ref{sec:Bf)} and \ref{sec:Bh)}
below.  This spectral flow
function is denoted by \(\grf_{s}\).   It suffices for now to know only
that this extended function \(\grf_{s}\) has integer values and that the
functions \[\grc\grs^\grf:=\grc\grs  - 4\pi ^{2}\grf_{s}\quad
\text{and} \quad \gra^\grf:=\gra + 2\pi
(r  - \pi )\grf_{s}\] 
are invariant under the action of
\(C^{\infty}(M; U(1))\) on \(\op{Conn}\,  (E) \times C^{\infty}(M; \bbS)\) that
has a \(\hat{u} \in C^{\infty}(M;U(1))\) sending \((A, \psi)\) to \((A
-\hat{u}^{-1}d\hat{u}, \hat{u}\psi)\).  By
way of comparison, \(\gra\), \(\grf_{s}\), and \(\grc\grs\) are not
invariant under this action.   (The notions \(\gra^\grf\) and
\(\grc\grs^\grf\) can be generalized to be defined over the blown-up
configuration space; cf. e.g. \cite{KM} Equation (16.4). The arguments
in the proof of Lemma 16.4.4 therein show that this generalization is
also invariant under gauge actions). 

Denote by \(\mathcal{Z}_{w,r}\) the set of gauge-equivalence classes of
solutions to the corresponding (\ref{eq:(A.4)}). (This was denoted by
a slightly different notation, \(\mathcal{Z}_{SW,r}\), in
\cite{KLT4}). It is well-known that in this case, for a generic choice
of \(r\), \(\mu\), this set \(\mathcal{Z}_{w,r}\) consists of
finitely many, nondegenerate irreducible
elements. (Cf. e.g. (IV.1.18) and references therein, ignoring the
``holonomy nondegenerate'' condition there for the moment). Assume
this to be the case. Consider next the 4-dimensional Seiberg-Witten
equations on the product cobordism \(\bbR\times M\), with
\(w_X=2w^+\), and \(\mu_-=\mu_+=\mu\). Here, \(w\) is used to denote
the pull-back of the 2-form \(w\) on \(M\) under the projection of
\(\bbR\times M\) to its second factor. Given an instanton \(\grd\) on
this product cobordism with \(s\to-\infty\) and \(s\to \infty\) limits
given respectively by representatives of \(\grc_-\), \(\grc_+\) in
\(\mathcal{Z}_{w,r}\). The differential operator in
(IV.1.21) has a Fredholm extension, whose index we denote by
\(\imath_\grd\). By \cite{APS}, in this case
\begin{equation}\label{eq:ind-SF}
\imath_\grd=\grf_s(\grc_+)-\grf_s(\grc_-).
\end{equation}

Let \(\mathcal{M}_k(\grc_-, \grc_+)\) denote the space 
of gauge equivalence classes of such instantons with \(\imath_\grd=k\). These spaces are
\(k\)-dimensional manifolds with a free \(\bbR\)-action when the
perturbation term in the Seiberg-Witten equations is suitable and
\(k>0\). In particular, the monotonicity assumption guarantees that
\(\mathcal{M}_1(\grc_-, \grc_+)/\bbR\) consists of finitely many
elements. With a coherent orientation chosen (this amounts to choices
of preferred elements of \(\Lambda(\grc)\) for all \(\grc\in
\mathcal{Z}_{w,r}\) in the language of \cite{KM}), each element in
\(\mathcal{M}_1(\grc_-, \grc_+)/\bbR\) is assigned a sign.

Fix a ring \(\bbK\), which can be taken to be \(\bbZ\) for the rest of this
article. The chain module for the monopole (or alternatively,
Seiberg-Witten) Floer chain group is the free \(\bbK\)-module
generated by \(\mathcal{Z}_{w,r}\), denoted by
\(\bbK(\mathcal{Z}_{w,r})\) below. The spectral flow function
\(\grf_s\) descends to define a relative \(\bbZ/c_\grs\bbZ\)-grading
on this module, where \(c_\grs\in 2\bbZ\) is the divisibility 
of the first Chern class of the \(\Spin^c\) structure \(\grs\). The differential \(\partial_{w, r}\) of ths monopole Floer
complex in this situation is the endomorphism of
\(\bbK(\mathcal{Z}_{w,r})\) given by the rule
\begin{equation}\label{eq:differential}
\grc_1\mapsto\sum_{\grc_2\in \mathcal{Z}_{w,r}}\op{w}(\grc_1, \grc_2)\, \grc_2,
\end{equation}
where \[\op{w}(\grc_1, \grc_2)=\sum_{\grd\in
  \mathcal{M}_1(\grc_1, \grc_2)/\bbR}\op{sign} (\grd)=\chi(
\mathcal{M}_1(\grc_1, \grc_2)/\bbR).\] 
The afore-mentioned properties of \(\mathcal{Z}_{w,r}\) and \(
\mathcal{M}_1(\grc_1, \grc_2)/\bbR\) for suitable monotone
perturbations guarantee that this homomorphism is well-defined, and it
is of degree 1 according to (\ref{eq:ind-SF}). A
typical gluing argument shows that \(\partial_{w,
  r}^2=0\). (Cf. e.g. \cite{KM} Sections 19 and 22). The
homology of the above monopole Floer complex is the monopole Floer homology,
or alternatively, the Seiberg-Witten-Floer homology of the negative
monotone perturbation (\ref{eq:2rw}). This is denoted
as \(\HM_*(M, \grs, c_-)\) below. The monopole Floer homology for
positive monotone perturbation forms, still assuming that
\(c_1(\grs)\) is nontorsion, is defined in the same way. 

\paragraph{Part 2: Local coefficients. } One may also associate monopole Floer homologies for
more general Seiberg-Witten equations (\ref{eq:SW-3d}). The
construction of monopole Floer complexes in Part 1 may fail to work
due mainly to two reasons: \begin{itemize}
\item[(1)] With balanced perturbations, the generating set of the chain group,
\(\mathcal{Z}\), namely the set of gauge equivalence classes of solutions to
(\ref{eq:SW-3d}), may contain reducible elements. (Recall that
\(\mathcal{Z}=\mathcal{Z}_{w,r}\) in the previous part, which consists of
finitely many irreducible elements. 
\item[(2)] The space \(\mathcal{M}_1(\grc_1, \grc_2)/\bbR\) might
  contain infinitely many elements, making the coefficients appearing in (\ref{eq:differential})'s
  formula for the Seiberg-Witten differential, \(\op{w} (\grc_1,
  \grc_2)\), undefined. 
\end{itemize}
The second issue above can be dealt with by working with monopole
Floer complexes with more general coefficients (as opposed to Part 1's
\(\bbZ\)-coefficient monopole Floer complex). Cf. \cite{KM} Section
22.6. 

Assume for simplicity that the perturbation \(\varpi\) in (\ref{eq:SW-3d}) is non-balanced, so
that the issue (1) above can be ignored: namely, with a generic
perturbation 
\(\mathcal{Z}\) will still consist of finitely many,
nondegenerate, irreducible elements. 
Fix a local system \(\Gamma\) in the sense described in
\cite{KM}. This assigns to every \(\grc\in\mathcal{Z}\) a group
\(\Gamma(\grc)\) and for each relative homotopy class \(z\) of paths between
\(\grc_1, \grc_2\in \mathcal{Z}\subset \mathcal{B}(M)\), a homomorphism \(\Gamma(z)\co
\Gamma(\grc_1)\to \Gamma(\grc_2)\). The monopole Floer chain complex with
local coefficient \(\Gamma\), \((C, \partial)\), has \(C:=\bigoplus_{\grc\in
  \mathcal{Z}}\Gamma(\grc)\) as its chain module. 
As for its differential \(\partial\), regard each \(\grd\in
\mathcal{M}(\grc_1, \grc_2)\) as path in
\(\mathcal{B}(M)\) and let \(\mathcal{M}_z(\grc_1,
\grc_2)\subset\mathcal{M}(\grc_1, \grc_2)\) be the subspace consisting
of elements of relative homotopy class \(z\). Refining
(\ref{eq:differential}), the following formula
defines the associated differential \(\partial\in \op{End} (C)\)
\begin{equation}\label{eq:differential-loc}
\partial=\sum_{\grc_1, \grc_2\in \mathcal{Z}}\sum_{z\in\pi_1
  \mathcal{B}(M;\grc_1, \grc_2)}\op{w}(\grc_1,
\grc_2; z)\,\Gamma (z),
\end{equation}
where \begin{equation}\label{eq:differential-lo}\op{w}(\grc_1, \grc_2;z)=\sum_{\grd\in
  \mathcal{M}_{1,z}(\grc_1, \grc_2)/\bbR}\op{sign} (\grd)=\chi(
\mathcal{M}_{1,z}(\grc_1, \grc_2)/\bbR),
\end{equation}
and \(\pi_1
  \mathcal{B}(M;\grc_1, \grc_2)\) denotes the space of relative
  homotopy classes of instantons with \(\grc_1\) and
  \(\grc_2\) as its \(s\to -\infty\) and \(s\to \infty\) limits respectively. 
Typical compactness results can be used to ensure that each
coefficient \(\op{w}(\grc_1,
\grc_2; z)\) is finite. (See e.g. \cite{KM} Theorem 8.1.1 and
Proposition 16.1.4.) Though (\ref{eq:differential-loc}) may have
infinitely many non-vanishing terms, the sum may be well-defined when \(\Gamma\)
is chosen to satisfy certain completeness 
conditions depending on the choice of \(\grs\) and \([\varpi]\). See
Definition 30.2.2 in \cite{KM}. We call a local system \(\Gamma\)
satisfying this completeness condition {\em \((\grs,
  [\varpi])\)-complete} (as opposed to ``\(c\)-complete'' in
\cite{KM}).  There is also a more stringent notion of completeness
which depends only on the cohomology class \([d\gra]\in H^1(\mathcal{B}(M);\bbZ)\) due originally to
Novikov. This sort of local system is said to be ``strongly
\(c\)-complete'' in \cite{KM}; see Definition 30.2.4 therein. We call
such \(\Gamma \) {\em strongly \((\grs, [\varpi])\)-complete}
instead. We shall not encounter local systems other than \(\bbZ\)  except in
Proposition \ref{prop:conn-1} (b) below, which is not directly relevant to the proof of
Theorem \ref{thm:main}. The interested reader is therefore referred to
\cite{KM} for more details on the definition of monopole Floer homology with local coefficients. A brief
summary in alternative language may also be found in the last section
of \cite{LT}. In the monotone case discussed in Part 1 or the balanced
case in the upcoming Part 3, the (strong) \((\grs,
  [\varpi])\)-completeness condition is met for all coefficients, and
  sum (\ref{eq:differential-lo}) has finitely many non-vanishing
  terms. 



\paragraph{Part 3: Balanced perturbations.}
 We now briefly describe how Issue (1)  in Part 2 is dealt with in the balanced case. 
 For details, 
see Chapters VI and VIII in \cite{KM}. As
already mentioned in Section \ref{sec:2.1}, \cite{KM}
considered the extension of (\ref{eq:SW-3d}) to
\(\mathcal{C}^\sigma\). The set of gauge equivalence classes of solutions to this
extended Seiberg-Witten equation is denoted by \(\grC\). Suppose that
the perturbation to the Seiberg-Witten equation is suitable. The subsets
of irreducible, unstable reducible, and stable reducible elements are
respectively denoted by \(\grC^o\), \(\grC^u\), \(\grC^s\). (In the
nonbalanced situation previously considered,
\(\grC=\grC^o=\mathcal{Z}_{r,w}\)). The first three flavors of
monopole Floer homology as defined in \cite{KM} use different
combinations of \(\grC^o\), \(\grC^u\), \(\grC^s\) to generate the
chain groups: Set
\[\begin{split}
& C^o=\bbK(\grC^o), \quad C^u=\bbK(\grC^u), \quad C^s=\bbK(\grC^s),
\quad \text{and let}\\
& \hat{C}=C^o\oplus C^u, \quad \bar{C}=C^s\oplus C^u,
\quad \check{C}= C^o\oplus C^s.
\end{split}
\]

Meanwhile, the operator in (IV.I.21) has a Fredholm generalization for
paths \(\grd(s)\) in \(\op{Conn}\,  (E)\times C^\infty(M, \bbS)\) with
\(s\to\infty\) or \(s\to -\infty\) limits that are nondegenerate elements in 
\(\mathcal{C}^\sigma\). (See Sections 14.4 and 22.3 in \cite{KM}). The index of this
operator is also denoted by \(\imath_\grd\) below; and it may be used
to generalize the spectral flow function \(\grf_s\) to the set of
nondegenerate elements in \(\mathcal{C}^\sigma\). This in turn defines
a relative \(\bbZ/c_\grs\)-grading, \(\op{gr}\), on the modules \(C^o\),
\(C^u\), \(C^s\). The chain modules \(\hat{C}\), \(\bar{C}\),
\(\check{C}\) are also \(\bbZ/c_\grs\)-graded according to the
following rule:
\[\begin{split}
\hat{C}=\bigoplus_j\hat{C}_j, \quad \hat{C}=\bigoplus_j\hat{C}_j,
\quad \hat{C}=\bigoplus_j\hat{C}_j, \quad \text{where}\\ 
\hat{C}_j=C^o_j\oplus C^u_j, \quad \bar{C}_j=C^s_j\oplus C^u_{j+1},
\quad \check{C}_j= C^o_j\oplus C^s_j.
\end{split}
\]
Note that the \(\bar{C}\) chain module above is graded by a modified
grading \(\overline{\op{gr}}\), related to \(\op{gr}\) via Equation
(22.15) in \cite{KM}. 
To define the differentials, define homomorphisms
\(\partial^\sharp_\natural\co C^\sharp\to C^\natural\) via rules
similar to (\ref{eq:differential}) or (\ref{eq:differential-loc}) by counting irreducible instantons
with \(\imath_\grd=1\) whose  \(s\to -\infty\) and \(s\to\infty\)
limits are in \(\grC^\sharp\) and \(\grC^\natural\) respectively; see \cite{KM} Equation (22.8) for the
precise formulae. Here, \(\sharp\) and \(\natural \) may stand for one
of the labels \(u\), \(o\), \(s\); however, due to the way 
\(\grC^u\), \(\grC^s\), \(\grC^o\) are defined, only the homomorphisms
\(\partial^o_o\), \(\partial^o_s\), \(\partial^u_o\), \(\partial^u_s\)
are nontrivial. Meanwhile, there are homomorphisms \(\bar{\partial}^\sharp_\natural\co
C^\sharp\to C^\natural\), and with \(\sharp\) and \(\natural \) denoting
either the label \(u\) or \(s\), by counting reducible instantons whose  \(s\to -\infty\) and \(s\to\infty\)
limits are in \(\grC^\sharp\) and \(\grC^\sharp\) respectively, with
\(\bar{\op{gr}}\) differing by \(-1\). 
If the \(\Spin^c\)-structure and \([\varpi]\)
satisfy monotonicity condition, then
the differentials for the
complexes, \(\hat{\partial}\co \hat{C}\to \hat{C}\),
\(\bar{\partial}\co \bar{C}\to \bar{C}\), \(\check{\partial}\co
\check{C}\to \check{C}\) are defined in terms of these
homomorphisms via Equation (22.7) and Definition 22.1.3 in
\cite{KM}. To give some examples, \(\check{\partial}\co C^o\oplus C^s\to
C^o\oplus C^s\), \(\hat{\partial}\co C^o\oplus C^u\to
C^o\oplus C^u\) are respectively written in block form as:

\begin{equation}\label{eq:hat-d}\left[\begin{array}{cc}
\partial^o_o &-\partial^u_o\bar{\partial}^s_u\\
\partial^o_s &\bar{\partial}^s_s-\partial^u_s\bar{\partial}^s_u
\end{array}\right]; \quad \left[\begin{array}{cc}
\partial^o_o &-\partial^u_o\\
-\bar{\partial}^s_u\partial^o_s &-\bar{\partial}^u_u-\bar{\partial}^s_u\partial^u_s
\end{array}\right].
\end{equation}

The gluing theorems in \cite{KM} show that
\(\hat{\partial}^2\), \(\bar{\partial}^2\), \(\check{\partial}^2\) are
indeed all 0. When the perturbation is balanced, such as in the
statement of Theorem \ref{thm:main}, the homology of these chain
complexes \((\mathring{C}_*, \mathring{\partial}_*)\), namely the
corresponding monopole Floer homology, is denoted \(\mathring{\HM}_* \,(
M,
\grs, c_b)\) for \(\circ=\wedge, -, \vee\).

The aforementioned homomorphisms
\(\partial^\sharp_\natural\), \(\bar{\partial}^\sharp_\natural\) are
also used to define chain maps (denoted \(i\co \bar{C}\to \check{C}\),
\(j\co \check{C}\to \hat{C}\), \(p\co\hat{C}\to \bar{C}\) in
\cite{KM}) that do not
define a short exact sequence, but their induced maps
on homologies do, this being the first of the fundamental exact sequences
referred to in Theorem \ref{thm:main}. See Proposition 22.2.1 in
\cite{KM}.

\paragraph{Part 4: Notation and other remarks. } 
When specificality is desired, the notation 
\[\mathring{C}_*(M, \grs,
[\varpi];\Gamma)=\mathring{C}_*(M, \grs,
\varpi;\Gamma), \quad \circ=\wedge, -, \vee,\] is used 
to denote the monopole Floer complex
corresponding to the cylindrical version of (\ref{eq:(A.14)}) with an \((\grs,
[\varpi])\)-complete local coefficients \(\Gamma\), and \(\mathring{\HM}_*\,
(M,\grs, [\varpi];\Gamma)\) is used to denote  the corresponding monopole Floer
homology. (The Floer chain complex \(\mathring{C}_*(M, \grs,
\varpi;\Gamma)\) does depend on the choice of \(\varpi\), not just its
cohomology class, though its associated Floer homology only depends on
the cohomology class \([\varpi]\). The notation \(\mathring{C}_*(M, \grs,
[\varpi];\Gamma)\) is adopted when the specific 
representative \(\varpi\) of \([\varpi]\) is irrelevant). 
In particular, when \([\varpi]=2\pi c_1(\det \bbS)\), \(\mathring{C}_*(M, \grs,
[\varpi];\Gamma)\) and \(\mathring{\HM}_*\,(M,\grs, [\varpi];\Gamma)\)
are also respectively denoted by \(\mathring{C}_*\,
(M,\grs, c_b;\Gamma)\) and \(\mathring{\HM}_*\,
(M,\grs, c_b;\Gamma)\). The coefficient \(\Gamma\) is dropped from the
notation when it is \(\bbZ\), or not important. 
The following (admittedly sloppy) convention will be adopted for the rest of this article:
Since the Floer complexes \((\check{C},
\check{\partial})=(\hat{C},\hat{\partial})=(C^o,\partial^o_o)\) when
the perturbation is non-balanced, we
use \(\CM\) or \((\CM,\partial )\) to denote the one complex in this case. When we wish to emphasize
the \(\Spin^c\) manifold and/or cohomology class of perturbation etc used to define the
monopole Floer complex,  these data are added to the above expression
in parentheses such as \(\CM_* (M, \grs, [\varpi] )\) or \(\big(\CM_* (M, \grs),
\partial  _*(M, \grs)\big)\).

As final remarks to this subsection, note that in \cite{KM} there is an
equivalent, geometric version of grading for the monopole Floer
complexes in terms of homotopy classes
of oriented 2-plane fields. This is briefly described in Part 1 of Section
\ref{sec:6.1} below, and denoted by 
  \(\bbJ(M)\) therein. A very brief description of this in the special
  cases relevant to this article will appear in Part 1 of Section
  \ref{sec:6.1}. Meanwhile, the signs \(\op{sign}(\grd)\)
  assigned according to the rules in \cite{KM} depend on a choice of {\em
    homology orientation} of \(M\). See Definition 22.5.2 in
  \cite{KM}.

\subsection{Cobordism-induced maps between monopole Floer
  complexes}\label{sec:2.4}

Instantons on cobordisms \(X\) described in Section \ref{sec:2.2} are used to define maps 
between the monopole Floer
complexes. Details of the construction of these maps are given in Chapter VII of
\cite{KM} for cobordisms \(X\) between {\em connected} 3-manifolds
\(Y_-\) and \(Y_+\), even though properties of moduli spaces of
Seiberg-Witten instantons on more general \(X\), where \(Y_-\) may be
disconnected, are also established therein. In particular, taking
\(X\) to be a product cobordism \(\bbR\times M\), this construction is used to define
chain maps from \(\mathring{C}\) back to itself, \(\circ=\wedge, -,
\vee\), that induce the \(\mathbf{A}_\dag\)-module structure on the
corresponding Floer homology. 
Another application of these cobordism-induced chain maps is to define
chain homotopies between monopole Floer complexes \(\mathring{C}\)
associated to different metrics and \((\grT, \grS)\). See
e.g. the proof for Corollary 23.1.6 in \cite{KM} and its variants. According to the
conventions set forth in Section \ref{sec:notation}, this justifies
our notation for the monopole Floer complex, \(\mathring{C} (M, \grs,
[\varpi];\Gamma)\). 
In fact,
this type of arguments show that \(\mathring{C}\) with positively
proportional ``period class'' (\cite{KM} p.591) are chain homotopic to
each other. (See Theorem 31.4.1 in \cite{KM}). This in turn justifies using
the notation \(\CM\, (M, \grs, c_-)\) for any negatively monotone,
nonbalanced 
perturbation, according to our convention. 

The rest of this subsection is divided into four parts. In the first three
parts we review some basic elements in the construction to the afore-mentioned
cobordism-induced maps. The last part contains a 
generalization of \cite{KM}'s construction to certain simple
cobordisms between possibly disconnected manifolds, in order to accommodate our needs in
Section \ref{sec:6}. 

\paragraph{Part 1: Moduli spaces and their compactifications.}
Fix a \(\Spin^c\)
structure \(\grs_X\) on \(X\) which restricts to the \(s\leq -2\) and
the 
\(s\geq 2\) part of \(X\) respectively as \(\Spin^c\) structures
\(\grs_{-}\) on \(Y_-\) and \(\grs_+\) on \(Y_+\). Fix also a
self-dual two form \(\varpi_X\) on \(X\) satisfying (\ref{eq:(A.13a)})
and a suitable pair \((\grT^+, \grS^+)\). Let \(c_{\grs_X}\)
 denote the divisibility of \(c_1(\grs_X)\). This number divides both
 \(c_{\grs_-}\) and \(c_{\grs_+}\). Assume that \(Y_\pm\) are both
 connected in this part. 

Consider instantons \(\grd\)
 defined from (\ref{eq:(A.14)}) with 
 representatives of \(\grc_-\) and \(\grc_+\) respectively as its
 \(s\to -\infty \) and \(s\to \infty\)
 limits.  
 The index of the Fredholm operator that entered the definition of
 nondegeneracy for instantons is denoted by \(\imath_\grd\). This
 generalizes the notion of index in the case of product cobordisms described in
 the previous subsection; and it depends only on the relative homotopy
 class of \(\grd\). See again Chapter 24 of \cite{KM}.
Let \(\mathcal{M}_k(X; \grc_-, \grc_+)\) denote the space of gauge
equivalence classes of such 
instantons with \(\imath_\grd=k\). When \(\grc_-\in
\grC^\sharp(Y_-)\), \(\grc_+\in \grC^\flat (Y_+)\) are both reducible,
let \(\mathcal{M}_k^{\text{red}}(X;
\grc_-, \grc_+)\subset \mathcal{M}_k(X; \grc_-, \grc_+)\) be the
subspace consisting of reducible instantons. Note that \(\mathcal{M}_k^{\text{red}}(X;
\grc_-, \grc_+)=\mathcal{M}_k(X; \grc_-, \grc_+)\) in the cases when the pair
\((\sharp, \flat)\) is \((u, u)\), \((s,s)\), or \((s,u)\). When
\((\grT^+, \grG^+)\) is suitable, \(\mathcal{M}_k^{\text{red}}(X;
\grc_-, \grc_+)\) is a smooth manifold with dimension respectively \(k\), \(k\), \(k+1\), \(k-1\) in
the cases when the pair \((\sharp, \flat)\) is \((u, u)\), \((s,s)\), \((s,u)\),
or \((u,s)\). The moduli space \(\mathcal{M}_k(X; \grc_-, \grc_+)\) is a
\(k\)-dimensional manifold consisting purely of irreducible
instantons in the case when at least one of \(\grc_-\) or \(\grc_+\) is
irreducible, while it is a \(k\)-manifold with boundary  \(\partial\mathcal{M}_k(X; \grc_-, \grc_+)=\mathcal{M}_k^{\text{red}}(X;
\grc_-, \grc_+)\) in the case when \((\sharp, \flat)=(u,s)\). 

All the spaces \(\mathcal{M}_k(X; \grc_-, \grc_+)\) and \(\mathcal{M}_k^{\text{red}}(X; \grc_-, \grc_+)\) are given
orientations according to the rules specified in \cite{KM}. This
depends on a choice of what is called a ``homological orientation'' of
\(X\) {\em as a cobordism} in \cite{KM}. (See Definition 3.4.1 in \cite{KM}). Let \(\mathcal{M}_{k,z}(X;
\grc_-, \grc_+)\) and \(\mathcal{M}_{k,z}^{\text{red}}(X; \grc_-,
\grc_+)\) respectively be subspaces of \(\mathcal{M}_k(X; \grc_-, \grc_+)\) and \(\mathcal{M}_k^{\text{red}}(X; \grc_-, \grc_+)\)
consisting of instantons with relative homotopy class \(z\). (Given
\(\grc_-\), \(\grc_+\), and \(z\), the spaces \(\mathcal{M}_{k,z}(X;
\grc_-, \grc_+)\) (resp. \(\mathcal{M}_{k,z}^{\text{red}}(X; \grc_-,
\grc_+)\)) are empty for all \(k\in \bbZ\) except one. This is denoted \(\mathcal{M}_{z}(X;
\grc_-, \grc_+)\) (resp. \(\mathcal{M}_{z}^{\text{red}}(X;
\grc_-, \grc_+)\)) below.)  All the
moduli spaces introduced above lie in the orbit
space of \(\mathcal{C}^\sigma_{loc}(X)\) under the gauge action by
\(C^\infty(X, U(1))=:\mathcal{G}_{loc}(X)\). This orbit
space is denoted as 
\(\mathcal{B}^\sigma_{loc}(X)\). Let
\(\mathcal{M}_k(X)\subset\mathcal{B}^\sigma _{loc}(X)\) denote the union of
all spaces \(\mathcal{M}_{z}(X;
\grc_-, \grc_+)\) and  \(\mathcal{M}_{z}^{\text{red}}(X; \grc_-,
\grc_+)\) with dimension less or equal to \(k\), for all \(\grc_-\in
\grC(Y_-)\), \(\grc_+\in \grC (Y_+)\), and \(z\in \pi_0(\mathcal{B}^\sigma _{loc}(X))\).

It follows from \cite{KM}'s Section 13.6 that the embeddings \(\M (X)=\bigcup_k \M_k
(X)\hookrightarrow\mathcal{B}^\sigma _{loc}(X)\) and \(\M
(X)\hookrightarrow\mathcal{B}^\sigma _{l,loc}(X)\) factor respectively
through 
subspaces \(\mathcal{B}^\sigma (X)\subset
\mathcal{B}^\sigma_{loc}(X)\) and \(\mathcal{B}^\sigma _l(X)\subset
\mathcal{B}^\sigma_{l, loc}(X)\), described below. These subspaces are homotopy equivalent to
\(\mathcal{B}^\sigma _{loc}(X)\) but are sometimes more convenient to
work with. In particular, \(\mathcal{B}^\sigma _l(X)\) has the virtue
of carrying a Banach manifold structure. Let 
\[\begin{split}
\mathcal{B}^\sigma (X)& :=\bigcup_{\grc_-\in \mathcal{B}^\sigma (Y_-)}\bigcup_{\grc_+\in \mathcal{B}^\sigma (Y_+)}\mathcal{B}^\sigma (X; \grc_-,
\grc_+),\\
\mathcal{B}^\sigma _l(X)& :=\bigcup_{\grc_-\in \mathcal{B}^\sigma _l(Y_-)}\bigcup_{\grc_+\in \mathcal{B}^\sigma _l(Y_+)}\mathcal{B}^\sigma _l(X; \grc_-,
\grc_+),\\
\end{split}
\]
where \(\mathcal{B}^\sigma (X; \grc_-,
\grc_+)=\bigcap_l \mathcal{B}^\sigma _l(X; \grc_-,
\grc_+)\subset \mathcal{B}^\sigma_{loc}(X)\), and \(\mathcal{B}^\sigma _l(X; \grc_-,
\grc_+)\subset \mathcal{B}^\sigma _{l,loc}(X)\) is defined as
follows. Let \(c_\pm=(\bbA_\pm, (\pmb{\Psi}_\pm, \Psi_\pm))\in
\mathcal{C}^\sigma _l(Y_\pm)\) be respectively representatives of
\(\grc_\pm \in \mathcal{B}^\sigma _l(Y_\pm)\), and use the same notation \((\bbA_\pm,
(\pmb{\Psi}_\pm, \Psi_\pm))\) to denote the corresponding
\(\bbR\)-invariant element in \(\mathcal{C}^\sigma
_{l,loc}(\bbR\times Y_\pm)\). Using the diffeomorphisms in
(\ref{(A.9a,11)}) to identify connected components of \(X-X_c\) with
subdomains of \(\bbR\times Y_+\)  or \(\bbR\times Y_+\), let \(\mathcal{C}^\sigma _l(X; c_-,
c_+)\subset \mathcal{C}^\sigma _{l,loc}(X)\) be the subspace
consisting of \((\bbA, (\pmb{\Psi}, \Psi))\in \mathcal{C}^\sigma
_{l,loc}(X)\) such that \(\bbA-\bbA_+, \Psi-\Psi_+\) are both
\(L^2_l\) on the positive end of \(X\), and \(\bbA-\bbA_-, \Psi-\Psi_-\) are both
\(L^2_l\) on the negative end of \(X\). Let \(\mathcal{B}^\sigma _l(X; \grc_-,
\grc_+)\subset \mathcal{B}^\sigma _{l,loc}(X)\) be the subspace
consisting of elements represented by elements in \(\mathcal{C}^\sigma _l(X; c_-,
c_+)\subset \mathcal{C}^\sigma _{l, loc}(X)\). By construction,
\(\mathcal{B}^\sigma (X)\), \(\mathcal{B}^\sigma _l(X)\) come equipped
with maps
\[\begin{split}
\Pi^\partial& =\Pi^{-\infty}\times \Pi^{\infty}\co \mathcal{B}^\sigma
(X)\to \mathcal{B}^\sigma(Y)\times \mathcal{B}^\sigma(Y),\\
\Pi^\partial& =\Pi^{-\infty}\times \Pi^{\infty}\co \mathcal{B}^\sigma_l
(X)\to \mathcal{B}^\sigma_l(Y)\times \mathcal{B}^\sigma_l(Y )
\end{split}
\]
sending \((\bbA, (\pmb{\Psi}, \Psi))\) to \((\bbA_-, (\pmb{\Psi}_-,
\Psi_-))\times (\bbA_+, (\pmb{\Psi}_+, \Psi_+))\).

Let \(\mathcal{M}_k^+(X; \grc_-,
\grc_+)\), \(\mathcal{M}_{k,z}^+(X; \grc_-,
\grc_+)\) be respectively the compactification of \(\mathcal{M}_k(X; \grc_-,
\grc_+)\) and \(\mathcal{M}_{k,z}(X; \grc_-,
\grc_+)\) by adding (parametrized) ``broken trajectories'' as described in
\cite{KM}'s Definition 24.6.1 and Theorem 24.6.2. In
Definition 24.6.9 of \cite{KM}, a surjective map \(\grr\) from
\(\mathcal{M}_{k,z}^+(X; \grc_-,\grc_+)\) to a smaller
compactification, \(\bar{\mathcal{M}}_{k,z}(X, \grc_-,
\grc_+)\subset \mathcal{B}^\sigma _{loc}(X)\), was introduced. Both
compactifications  \(\mathcal{M}_{k,z}^+(X; \grc_-,
\grc_+)\) and \(\bar{\mathcal{M}}_{k,z}(X; \grc_-,
\grc_+)\) are ``spaces stratified by manifolds'' in the sense of
\cite{KM}'s Definition 16.5.1. (Cf. \cite{KM} Propositions 24.6.8 and
24.6.10). For brevity, we refer to such spaces
simply as ``stratified manifolds'' in this article. 
By definition, \(\mathcal{M}_{k,z}(X; \grc_-,
\grc_+)\) is the top dimensional stratum of both  \(\mathcal{M}_{k,z}^+(X; \grc_-,
\grc_+)\), \(\bar{\mathcal{M}}_{k,z}(X; \grc_-,
\grc_+)\), and each \(\bar{\mathcal{M}}_k(X; \grc_-,
\grc_+)\) embeds in \(\mathcal{B}^\sigma (X)\subset \mathcal{B}^\sigma _{loc}(X)\) through the 
stratified manifold \[\begin{split}\mathcal{M}(X)
  &=\bigcup_k\mathcal{M}_k(X)\subset \mathcal{B}^\sigma (X),\\ 
\emptyset \subset & 
\cdots\mathcal{M}_{k-1}(X)\subset\mathcal{M}_{k}(X)\cdots\subset 
\mathcal{M}(X).
\end{split}
\] 
Meanwhile, the map \(\grr\) sends strata of \(\mathcal{M}^+_k(X; \grc_-,
\grc_+)\) to strata of \(\bar{\mathcal{M}}_k(X; \grc_-,
\grc_+)\) (not necessarily of the same dimension), and restricts to an
isomorphism on the top stratum. The moduli
spaces of reducible instantons \(\mathcal{M}_{k}^{\text{red}}(X; \grc_-,
\grc_+)\) are compactified similarly. 

\paragraph{Part 2. Integrating cochains on stratified manifolds.}
Generalizing the formula for the differential of monopole Floer
complex, (\ref{eq:differential-lo}), the purported maps between
monopole Floer complexes have coefficients given in terms of ``integrals'' of
the form \(\langle\mathpzc{u}, \mathcal{M}\rangle \), where
\(\mathpzc{u}\in \op{C}(\mathcal{B}^\sigma (X);\bbK)\),
\((\op{C}(\mathcal{B}^\sigma (X);\bbK), \delta)\) being a suitable version of
cochain complex for \(\mathcal{B}^\sigma (X)\),
\(H(\op{C}(\mathcal{B}^\sigma (X);\bbK))=H^*(\mathcal{B}^\sigma
(X);\bbK)\), and
\(\mathcal{M}\subset \mathcal{B}^\sigma (X)\) is a 
compactified moduli space of the types described in the Part
1. Explicit formulae for these maps are given below; see 
(\ref{eq:m-coeff}) and thereabouts.
Before proceeding to explain the possible choices of \((\op{C}(\mathcal{B}^\sigma
(X);\bbK), \delta)\) and the definition of the integrals \(\langle\mathpzc{u},
\mathcal{M}\rangle \) associated to them, we make a few motivational remarks. 

Ideally, the stratification structure of 
the relevant \(\mathcal{M}=\mathcal{M}_k\)
is sufficiently simple, e.g. it is a manifold with corners such that
\begin{equation}\label{m-w-corner}\partial
  \mathcal{M}_k=\mathcal{M}_{k-1},\quad \partial\mathcal{M}_{k-1}=0.
\end{equation}
Defined from broken trajectories, the lower-dimensional strata of
\(\mathcal{M}\) typically have an explicit description in the manner of  
\cite{KM}'s Propositions 24.6.8 and 24.6.10. Thus, when (\ref{m-w-corner})
holds, Stokes' theorem of the form 
\[
\langle\delta\mathpzc{v}, \mathcal{M}\rangle =\langle \mathpzc{v}, \partial\mathcal{M}\rangle=\langle \mathpzc{v}, \mathcal{M}_{k-1}\rangle
\]
can be invoked to derive various essential identities for the associated
cobordism maps between monopole Floer complexes. For example, this type of arguments are used to show that when \(\mathpzc{u}\) is closed, the
associated map \(\mathring{m}[\mathpzc{u}]\) is a chain map, and thus
it induces maps between the corresponding monopole Floer homology groups. Moreover, the
induced maps on Floer homologies depend only on the cohomology class of
\(\mathpzc{u}\), \([u]\in H^k(\mathcal{B}^\sigma (X);\bbK)\), rendering
the specific choice of \((\op{C}(\mathcal{B}^\sigma
(X);\bbK), \delta)\) irrelevant on the homological level.

With suitable \((\grT^+, \grS^+)\), (\ref{m-w-corner}) indeed holds in the non-balanced case, when all the
relevant 
Seiberg-Witten solutions are irreducible. Though the moduli spaces one
encounters may in general have more complicated stratification, it
was shown in \cite{KM} (e.g. Theorem 24.7.2 therein) that in most settings of interest, the
stratification is still simple enough so that (\ref{m-w-corner}) holds in a formal
sense (see \cite{KM}'s Lemma 21.3.1 for a precise statement). Thus,
via a suitable variant of Stokes' theorem (see \cite{KM}'s Equation
(21.4)), the arguments sketched above still apply, leading to the
desired identities.

Returning to the issue of choosing 
\(\op{C}(\mathcal{B}^\sigma (X);\bbK)\), a simplest option is the
de-Rham complex: taking \(\mathpzc{u}\) to be a differential \(k\)-form on
\(\mathcal{B}^\sigma (X)\), its restriction to \(\mathcal{M}\subset
\mathcal{B}^\sigma (X)\) or any stratum of \(\mathcal{M}\) is
well-defined, and the ``integral'' \[\langle\mathpzc{u},
\mathcal{M}\rangle=\int_{\cal
  M}\mathpzc{u}=\int_{\mathcal{M}_k\backslash\mathcal{M}_{k-1}}\mathpzc{u}\]
is literally the integral of \(\mathpzc{u}\) over
\(\mathcal{M}\). This however only works for \(\bbK=\bbR\). To be able
to work with more general \(\bbK\), in particular \(\bbK=\bbZ\),
\cite{KM} chooses to work with particular types of C\v{e}ch cochain
complexes \((C^*(\mathcal{U};\bbK), \delta)\), where  \(\mathcal{U}\)
is an open cover of \(\mathcal{B}^\sigma (X)\) satisfying certain
transversality conditions relative to
\(\mathcal{M}\subset \mathcal{B}^\sigma (X)\). It was shown that such
covering \(\mathcal{U}\) exists and any two of them have a common
refinement. See Chapter 21 of
\cite{KM}. The exposition in \cite{KM} focuses on maps between
monopole Floer {\em homology groups} instead of their underlying chain
maps between monopole Floer {\em complexes}. As mentioned previously, the former depends only
on the cohomology class \([\mathpzc{u}]\in H
(C^*(\mathcal{U};\bbK))=H^*(\mathcal{B}^\sigma (X);\bbK)\); thus, in
\cite{KM} the specific choice of the
covering \(\mathcal{U}\) and the cochain \(\mathpzc{u}\) representing
\([\mathpzc{u}]\in H^*(\mathcal{B}^\sigma (X);\bbK)\) is typically
left unspecified. In this article however, 
specific maps between monopole Floer {\em complexes} do play a role,
and the 
cochains \(\mathpzc{u}\) used to define these maps need to be
specified. This shall be done without reference to the covering \(\mathcal{U}\), as
there is no natural choice for the latter. Instead, in the upcoming Remark we introduce a
notion of equivalence (depending on \(\mathcal{M}(X)\subset \mathcal{B}^\sigma
(X)\)) among cochains possibly from different choices of underlying chain complexes
\(\op{C}^*(\mathcal{B}^\sigma (X);\bbK)\) for
\(H^*(\mathcal{B}^\sigma (X);\bbK)\). The maps
\(\mathring{m}[\mathpzc{u}]\) between monopole Floer complexes depend
only on the equivalence class of \(\mathpzc{u}\). Abusing terminology,
what is called an ``\(k\)-cochain \(\mathpzc{u}\) on \(\mathcal{B}^\sigma (X)\)'' in
this article typically refers to {\em any representative} \(\mathpzc{u}\in
\op{C}^*(\mathcal{B}^\sigma (X);\bbK)\) of a given
equivalence class (relative to a fixed \(\mathcal{M}(X)\)). 

\begin{remarks}\label{rem:cochain}
Let \(\mathcal{M}\) be a finite dimensional compact oriented stratified manifold embedded in a metric
space \(\mathcal{B}\). Suppose \(\mathcal{U}\) is an open covering of
\(\mathcal{B}\) transverse to \(\mathcal{M}\) in the sense defined in
\cite{KM}'s Chapter 21. As explained in \cite{KM}, 
the transversality condition on \(\mathcal{U}\) makes it possible 
to associate to each C\v{e}ch cochain \(\mathpzc{u}\in
C^k(\mathcal{U}; \bbK)\) a well-defined cohomology class on the
\(k\)-dimensional stratum of \(\mathcal{M}\), 
\[[\mathpzc{u}]\in \check{H}^k(\mathcal{M}_k,
\mathcal{M}_{k-1};\bbK)\simeq H^k_c(\mathcal{M}_k\backslash
\mathcal{M}_{k-1};\bbK),\] and the
value of \(\langle \, \mathpzc{u}, \grM \, \rangle\) for each stratum
\(\grM\) of \(\mathcal{M}\) is given in terms of
this cohomology class. See p. 408 of \cite{KM}. To rephrase the
constructions in \cite{KM},  we introduce a  cochain complex
\((C_{\cal M}^*, \delta_{\cal M})\)
defined as follows: Let \(C_{\cal M}^k=C_{\cal M}^{k;\bbK}:=H^k(\mathcal{M}_k,
\mathcal{M}_{k-1};\bbK)\), and let \(\delta_{\cal M}\co H^k(\mathcal{M}_k,
\mathcal{M}_{k-1};\bbK)\to H^{k+1}(\mathcal{M}_{k+1},
\mathcal{M}_{k};\bbK)\) be the connecting map in the long exact
sequence for the triple \((\mathcal{M}_{k+1},
\mathcal{M}_{k}, \mathcal{M}_{k-1})\). (The fact that \(\delta^2_{\cal
  M}=0\) is inessential in this article and we leave its verification to the reader). Use \([\mathpzc{u}]_{\cal M}\in C_{\cal M}^k\) to denote the
cohomology class of \(\mathpzc{u}\) in \(H^k(\mathcal{M}_k,
\mathcal{M}_{k-1};\bbK)\)  in the preceding expression. Then
by construction, 
\[[\delta\mathpzc{u}]_{\cal M}=\delta_{\cal M} [\mathpzc{u}]_{\cal
  M}.\] Let  \((C^{\cal
  M}_*, \partial _{\cal M})\) denote the dual chain complex of \((C_{\cal M}^*,
\delta_{\cal M})\). There is a canonical basis \(\{\mu^k_\alpha \}_\alpha \)
for \(C_{\cal M}^k\), with \(\alpha \) indexing all the connected
\(k\)-dimensional strata \(\grM_\alpha \) of \({\cal M}\), and
\(\mu^k_\alpha \) generating \(H^k (\grM_\alpha ,
\mathcal{M}_{k-1};\bbK)=\bbK\subset H^k(\mathcal{M}_k,
\mathcal{M}_{k-1};\bbK)\). Then duals of \(\mu^k_\alpha \), denoted
\([\grM_\alpha ]\) below, then form a corresponding basis for \(C^{\cal
  M}_k\).  This is used to define a notion of ``fundamental class''
for stratified manifolds: Given  a \(k\)-dimensional stratum \(\grM\) of
\(\mathcal{M}\), let
\[[\grM]:=\sum_\beta[\grM_\beta]\in C^{\cal
  M}_k,\] where \(\grM_\beta\) are the
connected components of \(\grM=\bigcup_\beta\grM_\beta\). We say that
\(\mathcal{M}'\subset\mathcal{M}\) is a \(k\)-dimensional {\em stratified
submanifold} of \({\cal M}\) if \({\cal M}'\) is a \(k\)-dimensional
stratified manifold whose strata are strata of \({\cal M}\). Given
such \(\mathcal{M}'\), let \[[\mathcal{M}']:=[\mathcal{M}'\backslash 
\mathcal{M}_{k-1}]\in C^{\cal M}_k.\] Then \(\langle \mathpzc{u}, \grM\rangle\) (resp. \(\langle
\mathpzc{u}, \mathcal{M}'\rangle\)) simply denotes the pairing of
\([\mathpzc{u}]_{\cal M}\in C_{\cal M}^*\) and \([\grM]\in C^{\cal
  M}_*\) (resp. \([\mathcal{M}']\in C^{\cal
  M}_*\)), and \cite{KM}'s version of   ``Stokes' theorem'' states: 
\[
\langle \delta\mathpzc{v}, \mathcal{M}'\rangle=\langle \delta_{\cal
  M}[\mathpzc{v}]_{\cal M}, [\mathcal{M}']\rangle=\langle
[\mathpzc{v}]_{\cal M}, \partial_{\cal M}[\mathcal{M}']\rangle.
\]
(Cf. \cite{KM}'s Equations (21.3) and (21.4)). 
In the case when \(\mathcal{M}'\) is a manifold with corners,
\(\partial_{\cal M}[\mathcal{M}']=[\partial\mathcal{M}']\) and the
right hand side of the preceding formula equals \(\langle
\mathpzc{v}, \partial\mathcal{M}'\rangle\), reducing the formula to the
usual Stoke's theorem. As noted in \cite{KM}, the compactness of
\(\mathcal{M}\) ensures the finiteness of the integrals \(\langle
\mathpzc{u}, \mathcal{M}'\rangle\), even though \(C^{\cal M}_*\) may
have infinite rank. 

Now suppose \(\mathpzc{u}\) is a differential \(k\)-form on
\(\mathcal{B}\). Since \(\mathpzc{u}\) restricts to a
closed form on any \(k\)-dimensional submanifold, its also determines
an element \([\mathpzc{u}]_{\cal M}\in C^{k;\bbR}_{\cal M}=H^k(\mathcal{M}_k,
\mathcal{M}_{k-1};\bbR)\). With \([\mathpzc{u}]_{\cal M}\) for
differential forms so defined,
one has \[\begin{split}
& \delta_{\cal M}[\mathpzc{u}]_{\cal
  M}=[d\mathpzc{u}]_{\cal M} \in C^{k+1;\bbR}_{\cal M} \quad \text{
and}\\ 
&  \langle [\mathpzc{u}]_{\M}, [\grM]\rangle=\int_\grM
\uu
\end{split}
\] for any \(k\)-dimensional stratum \(\grM\) of \(\mathcal{M}\).

Fix \(\mathcal{M}\subset \mathcal{B} \) and \(\bbK\). Let \(\mathpzc{u}\) be a \(k\)-cochains in one of the models for
\(\op{C}^*(\mathcal{B};\bbK)\) described above, namely, it is a C\v{e}ch cochain \(\mathpzc{u}\in
C^k(\mathcal{U};\bbK)\) for an arbitray open cover  \(\mathcal{U}\)
transverse to \(\mathcal{M}\), or when \(\bbK=\bbR\), it can be a differential
\(k\)-form on \(\mathcal{B}\). Let \(\mathpzc{u}'\) be another
\(k\)-cochains in a possibly different model of
\(\op{C}^*(\mathcal{B};\bbK)\). We say that the two ``\(k\)-cochains
on \(\mathcal{B}\)'', \(\uu\) and \(\uu'\), are  {\em equivalent on 
\(\mathcal{M}\),} (or simply ``equivalent'' if  the 
\(\mathcal{M}\) being referred to is clear) if \([\uu]_{\M}=[\uu']_{\M}\in C^{k;\bbK}_{\M}\).
(In other words, \(\uu\) and \(\uu'\) evaluate identically on all
\(k\)-dimensional strata of \(\M\).) 
To keep notations simple, we usually omit the subscript \({\cal
  M}\) from \(\delta _{\cal M}\) or \(\partial_{\cal M}\) below. 
\end{remarks}

Now let \(\mathcal{B}^\sigma (X)\) be as in Part 1, namely the orbit
space of \(\mathcal{C}^\sigma (X)\) under gauge group actions. Let  
\(\mathpzc{u}\) be a \(k\)-cochain on \(\mathcal{B}^\sigma (X)\) in
the sense just explained. For each fixed
 \(\Spin^c\)-structure, introduce homomorphisms \[\begin{split}
m^\sharp_\natural[\mathpzc{u}](X, \grs_X)\co
 C^\sharp(Y_-,\grs_-)\to C^\natural (Y_+, \grs_+)  &\quad  \text{for \(\sharp=o,
 u\), \(\natural=o, s\), } \\ \bar{m}^\sharp_\natural[\mathpzc{u}]( X,
 \grs_X)\co C^\sharp(Y_-,\grs_-)\to C^\natural (Y_+, \grs_+) & \quad \text{for \(\sharp=u,
 s\), \(\natural=u, s\)}
\end{split}
\]
respectively by the rules
\begin{equation}\label{eq:m-coeff}
\begin{split}
& \grC^\sharp\ni \grc_-\mapsto\sum_{\grc_+\in\grC^\natural} \langle \, \mathpzc{u},
\mathcal{M}_k(X; \grc_-, \grc_+)\, \rangle\, \grc_+;\\ & \grC^\sharp\ni \grc_-\mapsto\sum_{\grc_+\in\grC^\natural} \langle \, \mathpzc{u},
\mathcal{M}^{\text{red, k}}(X; \grc_-, \grc_+)\, \rangle \, \grc_+,
\end{split}
\end{equation}
where \(\mathcal{M}^{\text{red}, k}(X; \grc_-, \grc_+):=\mathcal{M}_{\ul{k}}^{\text{red}}(X;
\grc_-, \grc_+)\) with \(\ul{k}:=k, k, k-1, k+1\) respectively in
the cases when the pair \((\sharp, \flat)\) is \((u, u)\), \((s,s)\), \((s,u)\),
or \((u,s)\). (In other words, \(\mathcal{M}^{\text{red}, k}(X;
\grc_-, \grc_+)\) stands for the moduli space of reducible instantons
of dimension \(k\)). Note that  (the interior of) all \(\mathcal{M}_k(X; \grc_-,
\grc_+)\), \(\mathcal{M}^{\text{red}, k}(X; \grc_-, \grc_+)\), \(\mathcal{M}_z(X; \grc_-,
\grc_+)\), \(\mathcal{M}^{\text{red}}_z(X; \grc_-, \grc_+)\) are
strata or stratified submanifolds of \(\mathcal{M}(X)\subset
\mathcal{B}^\sigma (X)\). By the preceding Remark, the maps \(m^\sharp_\natural[\mathpzc{u}](X,
 \grs_X)\), \(\bar{m}^\sharp_\natural[\mathpzc{u}](X,
 \grs_X)\) depend only on the class \([\uu]_{\M}\in C^*_{\M}\). 

Once 
in place, the homomorphisms \(m^\sharp_\natural[\mathpzc{u}](X, \grs_X)\), \(\bar{m}^\sharp_\natural[\mathpzc{u}](X, \grs_X)\),
\(\partial^\sharp_\natural\,( Y_\pm,\grs_\pm)\),
\(\bar{\partial}^\sharp_\natural\,( Y_\pm,\grs_\pm)\) can be assembled
according to the formulae in Equation (25.5) and Definition 25.3.3 of \cite{KM}
into homomorphisms \[\mathring{m}[\mathpzc{u}] (X, \grs_X) \co \mathring{C}_*(Y_-,
\grs_-)\to \mathring{C}_*(Y_+, \grs_+)\] 
for \(\circ=\vee, -,
\wedge\). For example, for \(\mathpzc{u}\in C^k(\mathcal{U};\bbK)\),  \(\hat{m}[\mathpzc{u}]\co C^o(Y_-)\oplus
C^u(Y_-)\to C^o(Y_+)\oplus C^u(Y_+)\) is given in block form as:

\begin{equation}\label{def:m-map}\left[\begin{array}{cc}
m^o_o[\mathpzc{u}] &m^u_o[\mathpzc{u}]\\
(-1)^k\bar{m}^s_u[\mathpzc{u}]\partial^o_s-\bar{\partial}^s_um^o_s[\mathpzc{u}]
&(-1)^k\bar{m}^u_u[\mathpzc{u}]+(-1)^k\bar{m}^s_u[\uu] \partial
^u_s-\bar{\partial}^s_um^u_s [\mathpzc{u}]
\end{array}\right].
\end{equation}

The gluing theorems in Section 24.7 of \cite{KM} show
that when \(\mathpzc{u}\) is closed, these are chain maps,  with both \(\mathring{C}_*(Y_-,\grs_-)\) and
\(\mathring{C}_*(Y_+, \grs_+)\) regarded as chain complexes with
relative \(\bbZ/c_{\grs_X}\)-grading. As remarked in Section
\ref{sec:2.3}, gradings on \(\mathring{C}_*(Y_-,\grs_-)\) and
\(\mathring{C}_*(Y_+, \grs_+)\) are alternatively described in
\cite{KM} by
\(\bbJ(Y_-)\) and \(\bbJ(Y_+)\), the geometrically defined grading
sets \(\bbJ(Y_-)\) and \(\bbJ(Y_+)\). A cobordism \(X\) determines a relation \(\sim_X\) between
the grading sets \(\bbJ(Y_-)\) and \(\bbJ(Y_+)\) mentioned in Section \ref{sec:2.3}.
\begin{remarks}\label{rem:generalized-m}
In subsequent discussions, we make use of cobordism maps
\(\mathring{m}[\uu]\) associated to more general cochains than those
described above. (See in particular Part 3 of Section
\ref{sec:A-module} below.) Note that  the formula
(\ref{def:m-map}) defining \(\mathring{m}[\uu]\) assembles \(m^\#_\flat[\uu]\),
\(\bar{m}^\#_\flat[\uu]\), \(\partial ^\#_\flat\),
\(\bar{\partial}^\#_\flat\) in the particular manner specified in
\cite{KM}, 
so that desirable properties for \(\mathring{m}[\uu]\) may be obtained
by applying the Stokes' theorem for integrands of the form
\(\grr^*\uu\) on stratified submanifolds of \(\M^+(X)\), with \(\uu\in \op{C}(\B^\sigma _{loc}(X);\bbK)\). In other words, the
integrals defining \(\mathring{m}[\uu](X)\) factors through integrals
over the small compactified moduli space \(\bar{\M}(X) \). The more
general maps 
\(\mathring{m}[\uu](X)\) that we shall encounter are constructed by mapping \(\M^+(X)\) to a
larger space (typically a bundle over \(\B^\sigma _{loc}(X)\)), and
considering integrals of pull-backs of cochains on the latter larger
space over \(\M^+(X)\). To correctly assemble these integrals so as to
make the Stokes' theorem useful, the formula defining such
\(\mathring{m}[\uu](X)\) generalizes that given in \cite{KM} (as exemplified in
(\ref{def:m-map})) by replacing terms of the form
\(\bar{\partial}^s_um^\#_s [\mathpzc{u}]\) or \(m^u_\flat
[\mathpzc{u}]\bar{\partial}^s_u\) in \cite{KM}'s formulas with a sum,  in which it appears as
the first terms. In the notation of Part 2 of
Section \ref{sec:A-module}, the terms in the this sum take the general
form of \(\bar{n}^s_u[\op{u}_+]m^\#_s [\mathpzc{u}']\) or
\(m^u_\flat[\mathpzc{u}']\bar{n}^s_u[\op{u}_+]\), where \(\op{u}_\pm\)
are cochains on \(\B^\sigma (Y_\pm)\), and \(\deg (\uu')+\deg
(\op{u}_\pm)=\deg (\uu)-1\). In particular, 
\(\bar{\partial}^s_u=\bar{n}^s_u[1]\) in this notation. 
\end{remarks}

In order for the maps between Floer homologies induced by these chain maps to
behave well when composing cobordisms (exemplified by
Proposition 23.2.2 in \cite{KM}), 
one works with the assembled maps 

\[
\mathring{m}[\mathpzc{u}] \,( X)=\sum
_{\grs_X}\mathring{m}[\mathpzc{u}] \,( X, \grs_X)\co \bigoplus_{\grs_-}\mathring{C}_*(Y_-,
\grs_-)\to \bigoplus_{\grs_+}\mathring{C}_*(Y_+, \grs_+), \quad
\circ=\vee, -, \wedge,\]
where the direct sum \(\bigoplus_{\grs_\pm}\) is over the set of all
\(\Spin^c\) structures on \(Y_\pm\), and \(\grs_X\) runs through all
\(\Spin^c\) structures on \(X\).
As explained in Remark 24.6.6 in \cite{KM}, 
there can be infinitely many \(\grs_X\) to sum over for a fixed pair of \(\grs_-\),
\(\grs_+\). This necessitates the replacement of the chain complexes \(\mathring{C}_*(Y_-,
\grs_-)\), \(\mathring{C}_*(Y_+, \grs_+)\) in the preceding expression
by their ``grading-completed'' variants, \(\mathring{C}_\bullet(Y_-,
\grs_-)\), \(\mathring{C}_\bullet (Y_+, \grs_+)\) (Cf. Definition
3.1.3 and paragraphs around (30.1) in \cite{KM}). The cobordisms
relevant to our proof of Theorem \ref{thm:main} however have \(H^2(X,
Y_-)=0\), and this is why we  may 
use the pre-completion Floer complexes \(\mathring{C}_*\)
as the domain and target of \(\mathring{m}[\mathpzc{u}]\).

\paragraph{Part 3: Local coefficients.} 
The values  \(\langle \, \mathpzc{u},\mathcal{M}_k(X; \grc_-,
\grc_+)\, \rangle\), \(\langle \, \mathpzc{u},\mathcal{M}^{\text{red},k}(X; \grc_-,
\grc_+)\, \rangle\)  in (\ref{eq:m-coeff}) are finite only if the
moduli spaces \(\mathcal{M}_k(X; \grc_-,
\grc_+)\) , \(\mathcal{M}^{\text{red},k}(X; \grc_-,
\grc_+)\) have certain compactness properties. The
standard compactness arguments can be adapted to work with
nonvanishing \(\varpi_X\), when the perturbation form \(\varpi_X\) can be written as
\begin{equation}\label{pert-cobord}
\varpi_X=2\omega^+
\end{equation}
for some closed 2-form \(\omega\) on \(X\). We assume that \(\varpi_X\) satisfies (\ref{pert-cobord}) throughout
this article.  As with the monopole Floer
complex in Section \ref{sec:2.3},
the coefficients in
(\ref{eq:m-coeff}) are finite only when the cohomology classes \(c_1[\grs_X]\) and 
\([\omega]\) are related by certain constraints. 
A generalization of \cite{KM}'s  
Lemma 25.3.1 (making use the modified energy bounds from Section 29.1
therein)  guarantees that these 
constraints are met when 
\begin{equation}\label{eq:w_X}\begin{split}
& \text{\(\varpi_X=2rw_X\) for
\(r\neq 0\) and a \(w_X\) satisfying (\ref{eq:(A.13b)}), }\\
&\text{and when
\(X_{tor}\neq s^{-1}(\bbR-\{0\})\). }
\end{split}
\end{equation}
For more general pairs of \(c_1[\grs_X]\) and 
\([\omega]\), cobordism maps \(\mathring{m}[\mathpzc{u}]\)  
may still be well-defined for suitable local
coefficients. Let \(\Gamma_X\) be an ``\(X\)-morphism''
between local systems \(\Gamma_-\) on \(\mathcal{B}^\sigma(Y_-)\) and \(\Gamma_+\)
on \(\mathcal{B}^\sigma(Y_+)\) in the sense of \cite{KM}'s Definition
23.3.1. To each relative homotopy class \(z\in \pi_0
(\mathcal{B}^\sigma (X; \grc_-, \grc_+))\), \(\Gamma_X\) assigns an 
isomorphism \(\Gamma_{X}(z)\co \Gamma_-(\grc_-)\to
\Gamma_+(\grc_+)\). One then generalizes the homomorphisms (of
\(\bbZ\)-modules)
\(m^\sharp_\natural[\mathpzc{u}]( X, \grs_X)\), \(\bar{m}^\sharp_\natural[\mathpzc{u}]( X, \grs_X)\)   given by (\ref{eq:m-coeff}) 
to   \[\begin{split}
m^\sharp_\natural[\mathpzc{u}]( X, \grs_X; \Gamma_X) \co &
 C^\sharp(Y_-,\grs_-;\Gamma _-)\to C^\natural (Y_+, \grs_+;\Gamma _+),
 \quad \sharp=o,
 u, \quad \natural=o, s;\\
\bar{m}^\sharp_\natural[\mathpzc{u}]( X, \grs_X; \Gamma_X) \co &
 C^\sharp(Y_-,\grs_-;\Gamma _-)\to C^\natural (Y_+, \grs_+;\Gamma _+),
 \quad \sharp=
 u,s \quad \natural=u, s;
\end{split}
\] 
these are defined respectively by the formulae 
\begin{equation}\label{eq:m-coeff-loc}
\begin{split}
m^\sharp_\natural[\mathpzc{u}]( X, \grs_X; \Gamma_X) & =\sum_{\grc_-\in\grC^\sharp}\sum_{\grc_+\in\grC^\natural} \sum_{z\in \pi_0
(\mathcal{B}^\sigma (X; \grc_-, \grc_+))}\langle \, \mathpzc{u},
\mathcal{M}_{k,z}(X; \grc_-, \grc_+)\, \rangle\, \Gamma_X(z),\\
\bar{m}^\sharp_\natural[\mathpzc{u}](X, \grs_X; \Gamma_X) & =\sum_{\grc_-\in\grC^\sharp}\sum_{\grc_+\in\grC^\natural} \sum_{z\in \pi_0
(\mathcal{B}^\sigma (X; \grc_-, \grc_+))}\langle \, \mathpzc{u},
\mathcal{M}_{z}^{\text{red}, k} (X; \grc_-, \grc_+)\, \rangle\, \Gamma_X(z)
\end{split}
\end{equation}
where \(\mathcal{M}_{k,z}(X; \grc_-, \grc_+)\subset \mathcal{M}_{k}(X;
\grc_-, \grc_+)\), \(\mathcal{M}_{z}^{\text{red}, k} (X; \grc_-,
\grc_+)\subset \mathcal{M}^{\text{red}, k} (X; \grc_-, \grc_+)\) are
the subspaces consisting of elements with relative homotopy class \(z\). 
These \(m^\sharp_\natural\)  \(\bar{m}^\sharp_\natural\) are assembled in the same manner (e.g. (\ref{def:m-map}) for
\(\hat{m}\)) into the cobordism maps \[\mathring{m}[\mathpzc{u}](X,
\grs_X;\Gamma_X)\co \mathring{C}(Y_-,\grs_-;\Gamma _-)\to
\mathring{C}(Y_+, \grs_+;\Gamma _+)\quad \circ =\vee, -, \wedge.\]
Again, for the sums in (\ref{eq:m-coeff-loc}) to be well-defined, \(\Gamma_X\) and
\(\Gamma_{\pm}\) need to satisfy certain completeness conditions
depending on \(\grs_X\) and \(\varpi_X\). Here we limit ourselves to
some general remarks; 
more details will be provided on a case-by-case basis as
occasions arise. See also Section 25.3 in \cite{KM}, which contains some discussion on the case with
\(\varpi_X =0\). 
\begin{remarks}\label{rem:2.1}
In the more formal language of \cite{LT}'s Section 6.1,
  where a ``local system'' in Floer theory is described as a functor,
  an ``\(X\)-morphism'' from \(\Gamma_-\) to \(\Gamma_+\) is a natural
  transformation that intertwines the fundamental-groupoid structure
  on both sides. That is to say, it satisfies the composition law in
  \cite{KM}'s Equation (23.7).  (In \cite{KM}, \(\pi_0
(\mathcal{B}^\sigma (X; \grc_-, \grc_+))\) is denoted as \(\pi(\grc_-,
X, \grc_+)\) and an element in \(\mathcal{B}^\sigma (X; \grc_-,
\grc_+)\) is called an ``\(X\)-path''). For each pair \(\grc_-,
\grc_+\), the fundamental groups
\(\pi_1\mathcal{B}^\sigma (Y_-)\simeq H^1(Y_-;\bbZ)\) and
\(\pi_1\mathcal{B}^\sigma (Y_+)\simeq H^1(Y_+;\bbZ)\) act
respectively from the right and from the left on \(\pi_0
(\mathcal{B}^\sigma (X; \grc_-, \grc_+))\) through ``concatenation of paths''. Meanwhile, 
\begin{equation}\label{eq:pi_0-H_2}
\pi_0
(\mathcal{B}^\sigma (X; \grc_-, \grc_+))\simeq(j^{*})^{-1}
(c_1(\grs_X))\subset H^2(X, \partial X;\bbZ)
\end{equation}
in the following relative long
exact sequence: 
\[
\cdots\to H^1(\partial X;\bbZ)\stackrel{\delta}{\to }H^2(X, \partial X;\bbZ)\stackrel{j^*}{\to }H^2(X;
\bbZ)\stackrel{i^*}{\to }H^2(\partial X;\bbZ)\to \cdots.
\]
Note that \((j^{*})^{-1} (c_1(\grs_X))\) is 
an affine space under the abelian group \(\op{Im} (\delta )=\ker
(j^*)\). Under the identification (\ref{eq:pi_0-H_2}), the \(\pi_1\mathcal{B}^\sigma (Y_\pm)\simeq
H^1(Y_\pm;\bbZ)\) actions on \(\pi_0
(\mathcal{B}^\sigma (X; \grc_-, \grc_+))\)  respectively factor
through the aforementioned \(\delta H^1(\partial
X;\bbZ)\)-action on \((j^{*})^{-1} (c_1(\grs_X))\) under \(\delta_\pm:=\delta\circ
i_\pm\), where \(i_\pm\co H^1(Y_\pm;\bbZ)\hookrightarrow H^1(\partial
X;\bbZ)\) denotes the inclusion. 
The following simple consequences of the above observations will be
useful in this article: 
\begin{itemize}
\item  When \(\delta_\pm\co
H^1(Y_\pm;\bbZ)\to  \op{Im} \delta\) are both isomorphisms, any local
system \(\Gamma_-\) on \(\mathcal{B}^\sigma (Y_-)\) determines a local
system \(\Gamma_+\) on \(\mathcal{B}^\sigma (Y_+)\) and a unique (modulo automorphisms of \(\Gamma_-\), \(\Gamma_+\))
\(X\)-morphism \(\Gamma_X\) from \(\Gamma_-\) to
\(\Gamma_+\). Conversely, any local
system \(\Gamma_+\) on \(\mathcal{B}^\sigma (Y_+)\) also determines a local
system \(\Gamma_-\) on \(\mathcal{B}^\sigma (Y_-)\) and a unique 
\(X\)-morphism \(\Gamma_X\) from \(\Gamma_-\) to
\(\Gamma_+\).  In this case \(\pi_0(\mathcal{B}^\sigma (X, \grc_-,
\grc_+))\) is an affine space under both the actions of
\(\pi_1(\mathcal{B}^\sigma (Y_-))\) and \(\pi_1(\mathcal{B}^\sigma
(Y_+))\) and a choice of an element \(z_0\in \pi_0(\mathcal{B}^\sigma (X, \grc_-,
\grc_+))\) induces isomorphisms \(\iota^\pm_{z_0}\co \pi_1(\mathcal{B}^\sigma
(Y_\pm))\to \pi_0(\mathcal{B}^\sigma (X, \grc_-,
\grc_+))\) as \(\pi_1(\mathcal{B}^\sigma
(Y_\pm))\)-spaces. 
\item 
It was explained in
\cite{LT} that the ``\((\grs, [\varpi])\)-completeness'' condition for a
local system \(\Gamma\) in \(\mathring{C}(M, \grs, [\varpi];\Gamma )\)
is determined by the class \([\varpi]|_{\ker
  c_1(\grs)}\); in particular, when \([\varpi]|_{\ker
  c_1(\grs)}=0\) any \(\Gamma\) (including \(\bbZ\)) is \((\grs,
[\varpi])\)-complete. In the more general setting of cobordisms, the
cobordism map \(\mathring{m}[\mathpzc{u}](X;\Gamma_X)\) is
well-defined via (\ref{def:m-map}) when \(\Gamma_\pm\) are
respectively \((\grs_\pm, [\varpi_\pm])\)-complete, {\em and} an additional
completeness condition depending on the class
\([\omega]|_{\ker c_1(\grs_X)}\) is satisfied. (Here, \(c_1(\grs_X)\), \([\omega]\) are
both viewed as homomorphisms from \(H^2(X, \partial X)\) to \(\bbZ\)
via the Poincar\'e-Lefschetz duality). In particular, this additional
completeness condition is vacuous when \([\omega]|_{\ker c_1(\grs_X)}=0\).
Thus, the cobordism map \(\mathring{m}[u](X)\) is well-defined with
coefficient \(\bbZ\) via (\ref{eq:m-coeff}) when \([\omega]=2r
c_1(\grs_X)\) for \(r\in \bbR\), the setting relevant to the proof of
Theorem \ref{thm:main}. 
\end{itemize}
\end{remarks}

\paragraph{Part 4: Disconnected \(Y_-\) or \(Y_+\).} 
Suppose \(X_i\), \(i=1, \ldots ,k\) are respectively cobordisms from
\(Y_i^+\) to \(Y_i^+\), where all \(Y_i^\pm\) are connected. Then
\(X:=\coprod_iX_i\) may be viewed as a cobordism from \(Y_-:=\coprod_i
Y^-_i\) to \(Y_+:=\coprod_i Y_i^+\). The cobordism map \(\hat{m}[\uu]\)
introduced in Parts 2 and 3 above has a straightforward
generalization in this setting: Let \(\hat{C}(Y_\pm):=\bigotimes_{i=1}^k
\hat{C}(Y_i^\pm)\). Observe that in this case \(\B^\sigma
(X)=\prod_{i=1}^k\B^\sigma (X_i)\), and so given cochains \(\uu_i\in
\op{C}^*(\B^\sigma (X_i))\) (in the sense explained in Part 2) and
\(X_i\)-morphisms \(\Gamma_{X_i}\) from \(\Gamma_i^-\) to \(\Gamma
_i^+\) for
each \(i\), one has a cochain 
\(\uu:=\prod_i\uu_i\in \prod_i\op{C}^*(\B^\sigma
(X_i))=\op{C}^*(\B^\sigma (X))\) and an \(X\)-morphism \(\Gamma_X\)
from \(\Gamma_-:=\Pi_i\Gamma_i^-\) to \(\Gamma_+:=\Pi_i\Gamma_i^+\). Meanwhile, a set of local systems
\(\Gamma _i^-\) for each \(Y^+_i\) Define \(\hat{m}[\uu](X; \Gamma _X)\co
\hat{C}(Y_-; \Gamma _+)\to \hat{C}(Y_+; \Gamma _-)\) as 
\begin{equation}\label{def:m-disconn}
\hat{m}[\uu](X; \Gamma _-):=\bigotimes_{i=1}^k\hat{m}[\uu_i](X_i;\Gamma _{X_i})\co \bigotimes_{i=1}^k
\hat{C}(Y_i^-;\Gamma_i^- )\to \bigotimes_{i=1}^k
\hat{C}(Y_i^+; \Gamma_i^+). 
\end{equation}
The proof of  Theorem
\ref{thm:main} also requires maps associated to more general
cobordisms. For this purpose, it suffices to consider the \(\hat{m}\)
variant of the chain map for cobordisms \(X\)
satisfying the following constraint:

\begin{equation}\label{eq:assumption-X}
\begin{split}
&\text{At most one of \(Y_-\) or \(Y_+\) is disconnected, in which
  case it consists of }\\
&\text{two components. Moreover, at most one end of \(X\) is associated }\\
& \text{with balanced perturbation. }
\end{split}
\end{equation}

Assume that one of \(Y_-\), \(Y_+\)
is of the form \(Y_\sqcup=Y_1\sqcup Y_2\) for connected \(Y_1\) and \(Y_2\),
while the other is connected. Take \(Y_-=Y_\sqcup\) for example, since
the case where \(Y_+=Y_\sqcup\) is entirely parallel. Given the
self-dual 2-form \(\varpi_X\) described in (\ref{eq:(A.13a)}),  We shall always take \(Y_2\) to be the
only end of \(X\) possibly associated with a  balanced perturbation. 
Thus, \(\grC (Y_\sqcup)=\grC (Y_1)\times\grC(Y_2)=\grC^{oo}\sqcup
\grC^{ou}\sqcup \grC^{os}\), with \(\grC^{oo}\), \(\grC^{ou}\),
\(\grC^{os}\) denoting \(\grC ^o(Y_1)\times\grC^o(Y_2)\), \(\grC^o
(Y_1)\times\grC^u(Y_2)\), \(\grC ^o(Y_1)\times\grC^s(Y_2)\) respectively. 
Let \(C^{oo}(Y_\sqcup)=\bbK(\grC^{oo})=\CM (Y_1)\otimes C^o(Y_2)\),
\(C^{ou}(Y_\sqcup)=\bbK(\grC^{ou})=\CM (Y_1)\otimes C^u(Y_2)\),
\(C^{os}(Y_\sqcup)=\bbK(\grC^{os})=\CM (Y_1)\otimes C^s(Y_2)\).

In these cases we have the analogs of \(m^\sharp_\natural \) in
\cite{KM}, this being the homomorphisms \(m^{o \sharp}_\natural\co
\CM (Y_1)\otimes C^{\sharp}(Y_2)\to C^\natural (Y_+)\) (or in the case
where \(Y_\sqcup=Y_+\), 
\(m^\sharp_{o \natural }\co C^\sharp(Y_-)\to \CM (Y_1)\otimes
C^{\natural}(Y_2)\)), with 
\(\sharp\) standing for \(o\) or \(u\); and with the label \(\natural \) standing for \(o\) or \(s\).
Meanwhile, the analogs of \(\bar{m}^\sharp_\natural \) are all
trivial, since by (\ref{eq:assumption-X}) there are no reducible
instantons on \(X\).

As the condition (\ref{eq:assumption-X}) implies that
\(\hat{C}(Y_\#)=\CM (Y_\#)\), \(\hat{C} (Y_\sqcup)=C^{oo}\oplus
C^{ou}\), the maps \(\bar{m}=0\),  \(\hat{m}[\mathpzc{u}]\co\CM (Y_\#)\to C^{oo}\oplus
C^{ou}\) and \(\hat{m}[\mathpzc{u}]\co C^{oo}\oplus C^{ou}\to \CM
(Y_\#)\) respectively take the following simple form:
\begin{equation}\label{eq:m-V}
\left[\begin{array}{cc}
m^{oo}_o &m^{ou}_o
\end{array}\right], \quad 
\left[\begin{array}{c}
m_{oo}^o \\-(1\otimes\bar{\partial}^s_u(Y_2))\circ m_{os}^o
\end{array}\right].
\end{equation}

Further properties of the Floer complex \(\hat{C}(Y_\sqcup)\) and the
 maps \(\hat{m}\) associated to cobordisms \(X\) satisfying
 (\ref{eq:assumption-X}) will be discussed in 
Section \ref{sec:6.1}.

\paragraph{\it Caveat.} This Part assumes implicitly
that the \(X\)-morphisms and
local coefficients involved satisfy appropriate completeness
conditions in the sense of Remark \ref{rem:2.1}. While we forgo
general discussions of 
this issue; it will be addressed for the special cases in Sections
\ref{sec:6}.

\subsection{\(\mathbf{A}_\dag\)-module actions and geometric cochains}\label{sec:A-module}

In this subsection we introduce some useful cochains
\(\mathpzc{u}\) on \(\mathcal{B}^\sigma (X)\). They are described in
terms of differential forms on \(\mathcal{B}^\sigma (X)\) or
\(\mathcal{B}^\sigma _{loc}(X)\). (Note that a differential form on the latter
induces a 
corresponding differential form on the former via pulling back the embedding \(\mathcal{B}^\sigma (X)\hookrightarrow
\mathcal{B}^\sigma _{loc}(X)\); since \(\M(X)\subset \mathcal{B}^\sigma (X)\hookrightarrow
\mathcal{B}^\sigma _{loc}(X)\), they are equivalent on  \(\M(X)\) in the
sense of Remark \ref{rem:cochain}.) To a connected \(d\)-dimensional
submanifold of \(X\), we associate an element of
\(\Omega^{2-d}(\mathcal{B}^\sigma (X))\). There are many possible
choices of this differential form, but its equivalence class in
\(C^{*;\bbR}_{\M}\) will be fixed. To work with more general \(\bbK\),
this class is often replaced by a cohomologous element from
\(C^{*;\bbZ}_{\M}\subset C^{*;\bbR}_{\M}\). 
We then describe 
\(\mathbf{A}_\dag\)-module actions on monopole Floer complexes and
related chain homotopy maps as maps
\(\mathring{m}[\mathpzc{u}]\) associated to product cobordisms
\(X=\bbR\times M\) and cochains \(\mathpzc{u}\) of this type. 

The significance of such geometrically constructed cochains is that
the Seiberg-Witten cobordism maps \(\hat{m}[\uu]\) have natural counterparts
in invariants (some yet to be rigorously defined) constructed from
counting pseudo-holomorphic curves; in the latter case, the cobordism
maps are constructed from submanifolds in \(X\).  

Let \(X\) be a \(\Spin^c\) 4-manifold described by (\ref{(A.9a,11)})
and (\ref{eq:(A.12,15a)}), and let \(\mathpzc{E}=\{M_i\}_i\) be the
set of connected oriented \(\Spin^c\) manifolds indexing the ends of \(X\).

Fix a self-dual two form \(\varpi_X\) on \(X\) satisfying (\ref{eq:(A.13a)})
and a suitable pair of \((\grT^+, \grS^+)\). Let \(\mathcal{M}(X)\) be the
stratified manifold of instanton solutions to
(\ref{eq:(A.14)}) introduced in Part 1 of the last subsection, with stratification
\(\emptyset \subset
\mathcal{M}_0(X)\subset\cdots\mathcal{M}_k(X)\subset
\mathcal{M}_{k+1}(X)\subset \mathcal{M}(X)\) as before. 

Fix a hermitian line bundle \(K\) on \(X\) and a smooth connection
\(A_K\) on \(K^{-1}\), and write 
\begin{equation}\label{def:E-4d}
\det (\bbS^+)=E^2\otimes
K^{-1};
\end{equation}
namely, a 4-dimension version of (\ref{def:E-3d}). Let
\(A\in \op{Conn}(E)\) denote the unitary connection induced from
\(\bbA\in \op{Conn}(\det \bbS^+)\). As mentioned previously in the end of
Section \ref{sec:2.2}, both \((\bbA, \Psi)\in \op{Conn}(\det \bbS^+)\times C^\infty(\bbS^+)\) and its corresponding \((A, \psi)\in \op{Conn}
(E)\times C^\infty(\bbS^+) \) are used to denote  an
element in \(\mathcal{C}(X)\). At this point \(K\) is not assumed to
be related to \(\varpi_X\). 
In the case when the factorization (\ref{def:E-3d}) or
(\ref{def:E-4d}) arises from a
splitting \(\bbS\) or \(\bbS^+=E\oplus E\otimes K^{-1}\), we write \(\psi=(\alpha
, \beta )\), where \(\alpha , \beta \) are respectively the
\(E\)- and the \(E\otimes K^{-1}\)-component of \(\psi\) under the
decomposition.

\paragraph{Part 1. Cocycles on \(\mathcal{B}^\sigma (X)\)
                   from closed \(d\)-submanifolds in \(X\).}
The cocycles in this Part are constructed from differential forms on
\(\mathcal{B}^\sigma _{loc}(X)\). As mentioned previously, they induce
differential forms on \(\mathcal{B}^\sigma (X)\) and we shall use the
same notations for forms on \(\mathcal{B}^\sigma _{loc}(X)\) and their
corresponding forms on \(\mathcal{B}^\sigma (X)\). Alternatively, one
may define the forms on  \(\mathcal{B}^\sigma (X)\) by carrying out
parallel arguments using \(\mathcal{C}^\sigma (X)\) in place of
\(\mathcal{C}^\sigma _{loc}(X)\). 

\paragraph{(a)} {\em  When \(d=0\).} To a point \(x\in X\) we
associate an integral 2-cocycle   \([\mathpzc{e}]_{\M (X)}\in
C^{2;\bbZ}_{\M (X)}\) as follows. 
Consider the subgroup
\[\mathcal{G}_{x,loc}\subset C^\infty(X, U(1)):=\mathcal{G}_{loc}(X)\] consisting of maps \(u\co X\to
U(1)\) with \(u(x)=1\). Then
\[\tilde{\mathcal{B}}^\sigma_{x,loc}(X):=\mathcal{C}^\sigma_{loc}(X)/\mathcal{G}_{x, loc}\]
admits a free \(U(1)=\mathcal{G}_{loc}(X)/\mathcal{G}_{x, loc}\)-action, and
\(\mathcal{B}^\sigma_{x, loc}(X)\) is the orbit space of this action. Let
\[\pi_x\co \tilde{\mathcal{B}}^\sigma_{x, loc}(X)\to \mathcal{B}^\sigma_{loc}(X)\]
denote the quotient map of this action. We use \(\vartheta\in \Omega
^1( \tilde{\mathcal{B}}^\sigma_{x, loc}(X))\) to denote a Thom form of the 
\(U(1)\)-fibration \(\pi_x\co \tilde{\mathcal{B}}^\sigma_{x, loc}(X)\to
\mathcal{B}^\sigma_{loc}(X)\), so that
\[
d\vartheta=\pi_x^*\mathpzc{e},
\] 
\(\mathpzc{e}\in\Omega^2(\mathcal{B}^\sigma_{loc}(X))\) being an
Euler form. Choose \(\vartheta\) so that 
it defines a principal \(U(1)\)-connection on \(
\tilde{\mathcal{B}}^\sigma_{x, loc}(X)\), now regarded as a principal
\(U(1)\)-bundle. 
In this setting  \((\pi_x)_*:=(\pi_x^*)^{-1}\) is
well-defined at \(d\vartheta \), and we formally write 
\(\mathpzc{e}=(\pi_x)_*(d\vartheta)\). Let 
 \(\mathcal{E}_x\) be the hermitian  line bundle associated to the principal \(U(1)\)-bundle
\(\tilde{\mathcal{B}}^\sigma_{x, loc}(X)\). The latter is identified with the
(\(U(1)\)-) fiber
product 
\[
\mathcal{E}_x(X):=\tilde{\mathcal{B}}^\sigma_{x, loc}(X)\times_{U(1)}E_x=(\tilde{\mathcal{B}}^\sigma_{x, loc}(X)\times
E_x)/\text{diagonal \(U(1)\)-action},
\]
where \(E_x\simeq \bbC\) is the fiber of the bundle \(E\) over \(x\in
X\), equipped with the
\(\mathcal{G}_{loc}(X)/\mathcal{G}_{x, loc}(X)=U(1)\)-action. 
 Then  \(\mathpzc{e}\) has an alternative interpretation as 
\(\frac{i}{2\pi } \) times the curvature form of the unitary connection
associated to \(\vartheta\) on \(\mathcal{E}_x\). 

The following alternative interpretation of
\({\cal E}_x(X)\) will come in handy later: Let
\(\pmb{\pi}\co \pmb{\cal
  E}(X)\to X\times {\cal B}^\sigma
_{loc}(X)\) be the ``universal family'' (described below) for the bundle \(\pi_E\co E\to
X\); then 
\[
\mathcal{E}_x(X)=\pmb{\cal
  E}(X)|_{\{x\}\times \mathcal{B}^\sigma _{loc}(X)}.
\]
The bundle \(\pmb{\cal
  E}(X)\) is constructed in the following manner. Consider the hermitian line bundle \(\pi_E\times \op{Id}\co E\times
\mathcal{C}^\sigma _{loc}(X)\to X\times
\mathcal{C}^\sigma _{loc}(X)\). This bundle is equipped with a
tautological unitary connection \(\tilde{A}\) characterized by the
following property: \(\tilde{A}|_{X\times\{\bbA,
  (\pmb{\Psi}, \Psi))\}}=A\) for all \((\bbA, (\pmb{\Psi},
\Psi))\in \mathcal{C}^\sigma _{loc}(X)\), and
\(\tilde{A}|_{\{x\}\times \mathcal{C}^\sigma _{loc}(X)}\) is trivial
for each \(x\in X\). In the case when \(E\subset \bbS^+\) is a summand
of a splitting of \(\bbS^+\), the bundle \(\pi_E\times \op{Id}\co E\times
\mathcal{C}^\sigma _{loc}(X)\to X\times
\mathcal{C}^\sigma _{loc}(X)\) also carries a 
tautological section \(\tilde{\alpha }\), characterized by the
property that  \(\tilde{\alpha }|_{X\times\{(\bbA,
  (\pmb{\Psi}, \Psi=\sqrt{2r}(\alpha , \beta ))\}}=\alpha \). Let \(\pmb{\cal
  E}(X)\) be the quotient of  \(E\times
\mathcal{C}^\sigma _{loc}(X)\) by the diagonal
\(\mathcal{G}_{loc}(X)\) action. The map  \(\pi_E\times \op{Id}\co E\times
\mathcal{C}^\sigma _{loc}(X)\to X\times
\mathcal{C}^\sigma _{loc}(X)\) then
descends to define a hermitian line bundle 
\[
\pmb{\pi}\co \pmb{\cal
  E}(X)\to X\times {\cal B}^\sigma_{loc}(X),
\] 
and \(\tilde{\alpha}\) (when defined), \(\tilde{A}\) descend respectively to define a tautological section and a
tautological unitary connection on \(\pmb{\cal
  E}(X)\), also denoted \(\tilde{\alpha}\) and \(\tilde{A}\)
below. Let \(\bar{X}\supset X\) denote the compactification of \(X\)
over which the diffeomorphisms in (\ref{(A.9a,11)}) extend to define a 
diffeomorphism between \(([-\infty, -L')\times Y_-) \sqcup ((L',
\infty]\times Y_-)\) and \(\bar{X}-X_c\). When restricted to \(X\times \mathcal{B}^\sigma (X)\subset
X\times \mathcal{B}^\sigma _{loc}(X)\), the bundle \(\pmb{\cal
  E}(X)|_{X\times \mathcal{B}^\sigma (X)}\) extends to define a bundle
over \(\bar{X}\times \mathcal{B}^\sigma (X)\), denoted 
\[
\pmb{\pi}\co \pmb{\cal
  E}(\bar{X})\to \bar{X}\times {\cal B}^\sigma(X)
\] 
below. The tautological
section and connection, \(\tilde{\alpha}\) (when defined) and \(\tilde{A}\), extend
over \(\pmb{\cal
  E}(\bar{X})\) and will be denoted by the same notation. 

Restricting the tautological connection \(\tilde{A}\) to \(\pmb{\cal
  E}(X)|_{\{x\}\times \mathcal{B}^\sigma _{loc}(X)}=\mathcal{E}_x(X)\), one
has a unitary connection on \(\mathcal{E}_x(X)\). Let \(\tilde{\vartheta}\)
denote the corresponding principal \(U(1)\)-connection on
\(\tilde{\cal B}_{x,loc}(X)\), i.e. the principal \(U(1)\)-bundle
associated to \(\mathcal{E}_x(X)\), and let 
\begin{equation}\label{def:theta}
\theta:=(\pi_x)_*
(\vartheta-\tilde{\vartheta})\in \Omega^1(\mathcal{B}^\sigma _{loc}(X)).
\end{equation}

The form \(\vartheta\) (and consequently its associated \(\mathpzc{e}\)) is far from
unique. However, as mentioned in Remark \ref{rem:cochain},
we are only interested in \(\mathpzc{e}\)'s equivalence class rel
\(\mathcal{M}(X)\) or \(\vartheta\)'s equivalence class rel
\(\pi_x^{-1}\mathcal{M}(X)\), where
\(\pi_x^{-1}\mathcal{M}(X)\subset\tilde{\cal B}^\sigma _{x, loc}(X)\) is
viewed a stratified manifold with stratification \(\emptyset \subset
\cdots\pi_x^{-1}\mathcal{M}_{k+1}(X)\subset
\pi_x^{-1}\mathcal{M}_k(X)\cdots\subset
\pi_x^{-1}\mathcal{M}(X)\). For this purpose it suffices to describe \(\vartheta|_{\pi_x^{-1}\mathcal{M}_1(X)}\). 

We say that the connection \(\vartheta\) is {\em integral} over
\(\mathcal{M}_1(X)\) if it is induced from 
a trivialization
\(\rho_\vartheta \co \ul{\bbC} \stackrel{\sim }{\to }
\mathcal{E}_x|_{\mathcal{M}_1(X)}\), where \(\ul{\bbC} \) denotes the
trivial  \(\bbC\)-bundle \(\mathcal{M}_1(X)\times \bbC\). Note that
conversely 
\(\rho _\vartheta\) is uniquely determined by \(\vartheta \) modulo
constant \(U(1)\) actions. We also use \(\rho _\vartheta \) to
denote the associated trivialization \(\ul{U(1)} \stackrel{\sim }{\to
}\tilde{\cal B}^\sigma _{x, loc}(X)|_{\mathcal{M}_1(X)}\). We require
that 
\begin{equation}\label{int-con}
\text{\(\vartheta\) be {\em integral} over
\(\mathcal{M}_1(X)\). }
\end{equation}
Such \(\vartheta \) exists since there is no
obstruction to trivializing \(U(1)\)-bundles over \(1\)-complexes. 
 Since the boundary of each 2-dimensional stratum
\(\grM\) of \(\mathcal{M}(X)\) lies in \(\mathcal{M}_1(X)\), a choice
of such  \(\vartheta\) determines a
well-defined relative euler class  for the \(U(1)\)-bundle
\(\tilde{\mathcal{B}}^\sigma_{x, loc}(X)|_{\mathcal{M}_2(X)\backslash \mathcal{M}_1(X)}\) (or equivalently, for 
\(\mathcal{E}_x|_{\mathcal{M}_2(X)\backslash \mathcal{M}_1(X)}\)).
This class in \(H^2(\mathcal{M}_2(X), \mathcal{M}_1(X); \bbZ)\) is by
definition the equivalence class \([\mathpzc{e}]_{\M (X)}\in
C^{2;\bbZ}_{\mathcal{M}(X)}\). 

Two connections \(\vartheta_1\) and \(\vartheta_2\)  that are both integral
over \(\mathcal{M}_1(X)\) differ by
\(\vartheta_2-\vartheta_1=\pi_x^*df\) on \(\mathcal{M}_1(X)\), where \(f\) is a map \(f\co
\mathcal{M}_1(X)\to U(1)=\bbR/\bbZ\) (such \(f\) is unique modulo
constant maps). 
Thus, \([(\pi_x)_*(\vartheta_2-\vartheta _1)]_{\M(X)}\in
C^{0,\bbR}_{\M(X)}\) is a closed element. 
We say that \(\vartheta _1\),
\(\vartheta _2\) are {\em \(\delta\)-cohomologous} if \([df]=0\in
H^1(\M_1(X);\bbZ)\). In this case \(f\) factors through a
map \(\tilde{f}\co
\mathcal{M}_1(X)\to \bbR\), and the restriction \(\tilde{f}|_{\mathcal{M}_0(X)}\) defines a
class 
\([\tilde{f}]_{\cal M}\in C^{0;\bbR}_{\mathcal{M}(X)}\). We have
\([(\pi_x)_*(\vartheta_2-\vartheta_1)]_{\M(X)}=\delta [\tilde{f}]_{\mathcal{M}(X)}\) for \(\delta\)-cohomologous
  \(\vartheta _1\) and \(\vartheta _2\), and hence \([\e]_{\M(X)}\)
    depends only on the \(\delta\)-cohomology class of  \(\vartheta\). 

\paragraph{\it Convention.} When we wish to emphasize the choice of \(x\), we add a subscript
\(x\) to the forms \(\vartheta \), \(\tilde{\vartheta}\), \(\theta \), \(\e\) introduced
above. E.g. \(\vartheta _x\) denotes the \(\vartheta \) associated to
\(x\). 

\paragraph{(b)} {\em When \(d=1\).} Let \(\gamma\subset X\) be an
embedded oriented circle in the interior of \(X\). To such a \(\gamma
\), we associate a real 1-cocycle \([\theta_\gamma ]_{\M(X)}\in
C^{1;\bbR}_{\M (X)}\). Modifying \([\theta_\gamma ]_{\M(X)}\),
\(\gamma\) is also associated an integral 1-cocycle \([\uu_\gamma
]_{\M(X)}\) cohomologous to \([\theta_\gamma ]_{\M(X)}\). Let 
\[\op{hol}_\gamma \co \mathcal {B}^\sigma_{loc}(X)\to U(1)=\bbR/\bbZ\] be the map
sending an element \(\grd\in\mathcal{C}^\sigma_{loc}(X)\) to the holonomy
of \(A\in \op{Conn}(E)\) associated to
\(\grd\). Let 
\[
\theta_\gamma:=d\op{hol}_\gamma \in \Omega^1(\B^\sigma_{loc}(X)).
\] 
This is an integral closed 1-form on \(\mathcal {B}^\sigma_{loc}(X)\) and defines
a class \([\theta_\gamma ]_{\M(X)}\in C^{1,\bbR}_{\M (X)}\). For the
purpose of defining cobordism maps, it is
often desirable to replace \([\theta_\gamma ]_{\M(X)}\in
C^{1,\bbR}_{\M (X)}\) with a cohomologous element
\[\uu_{\gamma}=[\theta_\gamma ]_{\M(X)}-\delta\varepsilon_\gamma, \]
\(\varepsilon_\gamma\in C^{0;\bbR}_{\M (X)}\), so that
\(\uu_{\gamma}\in C^{1,\bbZ}_{\M (X)}\subset C^{1,\bbR}_{\M
  (X)}\). We call \(\uu_\gamma\) an 
``integral correction'' of \([\theta_\gamma ]_{\M(X)}\). 
A choice of integral correction \(\uu_\gamma \) is equivalent  to a
choice of lifting, 
\[
h_\gamma \co \M_0(X)\to\bbR
\] 
for \(\op{hol}_\gamma |_{\M_0 (X)}\co \M_0(X)\to
\bbR/\bbZ\), as \(\delta\varepsilon_\gamma=\delta[h_\gamma ]_{\M(X)}\). Different
choices of \(h_\gamma\) defer by elements in
\(C^{0;\bbZ}_{\M (X)}\).

\paragraph{Part 2: Product cobordisms and \({\bf A}_\dag(M)\)-actions.}
In this part we apply the construction in Part 1 to the case of product cobordisms. Let
\(X=\bbR\times M\), \(M\) being a closed connected \(\Spin^c\)
3-manifold. The cocycles \(\e\), \(\theta_\gamma\) on
\(\mathcal{B}^\sigma (X)\) described below, loosely speaking,  will
take the form of pull-backs from corresponding cocycles on
\(\mathcal{B}^\sigma (M)\). The latter cocycles are chosen to represent generators of the
cohomological algebra 
\begin{equation}\label{coho-B}
\begin{split}
H^*(\mathcal{B}^\sigma
(M);\bbZ)& =H^*(\mathcal{B}^\sigma_{loc}
(M);\bbZ)\\
& \simeq H^*(\bbC P^\infty;\bbZ)\otimes
H^*(H^1(M;\bbR)/H^1(M;\bbZ);\bbZ)=\mathbf{A}_\dag(M).
\end{split}
\end{equation}
Let \(U\in H^2(\bbC
P^\infty;\bbZ)\) be the generator of the polynomial algebra \(H^*(\bbC
P^\infty;\bbZ)\), and let \(\{\grt_i\}_i\) be a basis of \(H_1(M;\bbZ)/\op{Tors}\simeq
H^1(H^1(M;\bbR)/H^1(M;\bbZ);\bbZ)\). We use the same notations \(U\),
\(\grt_i\) to denote the corresponding generating elements of the
algebra (\ref{coho-B}). We shall introduce 2-cocycles
\(\mu_U\) representing \(U\) and 1-cocycles  \(\mu_{\grt_i}\)
representing \(\grt_i\), and the cobordism maps associated to pull-backs
of these
cocycles are referred to respectively as \(U\)-actions or
\(\grt_i\)-actions on the monopole Floer complex \(\mathring{C}(M)\). Together they generate the \(\mathbf{A}_\dag(M)\)-actions on
\(\mathring{C}(M)\). The choice of \(\mu_U\) depends on a choice of a
point \(p\in M\), while the choice of \(\mu_{\grt_i}\) depends on the
choice of an embedded circle \(\gamma_i\subset M\) representing \(\grt_i\).

Before proceeding, we make some preparatory remarks on \(\mathcal{B}^\sigma
_{loc}(X)\) and its variants in the case 
\(X=\bbR\times M\). 
As explained in \cite{KM}, by a unique continuation theorem
\(\M(\bbR\times M)\)
falls in a smaller blown-up configuration space
 \(\mathcal{B}^\tau _{loc}(\bbR\times M)=\op{Conn} (\det \bbS^+)\times
 \Gamma^\tau(\bbR\times M, \bbS^+)\subset \mathcal{B}^\sigma_{loc}
 (\bbR\times M)\),  which is often more
 convenient to work with. Here, \(\Gamma^\tau(\bbR\times M,
 \bbS^+)\subset \Gamma^\sigma(\bbR\times M, \bbS^+)\) consists of
 elements \((\bbA, (\pmb{\Psi}, \Psi))\) such that
 \(\Psi|_{\{s\}\times M}\neq 0\) for all \(s\in \bbR\). By
 construction, there exists for each \(s\in \bbR\) a map 
\[
\Pi^s\co \mathcal{B}^\tau_{loc}(\bbR\times M)\to \mathcal{B}^\sigma
(M)
\]
which is defined by restricting \(\bbA\) and \(\Psi\) to \(\{s\}\times M\subset X\). When restricted to \(\mathcal{B}^\tau(\bbR\times M):=\mathcal{B}^\tau
 _{loc}(\bbR\times M)\cap \mathcal{B}^\sigma (X)\), \(\Pi^s\) has
 well-defined limits as \(s\to\pm\infty\), 
\[\Pi^{\pm\infty}\co \mathcal{B}^\tau(\bbR\times M)\to
\mathcal{B}^\sigma (M).
\]
An element in \(\grd\in \mathcal{B}^\tau_{loc}(\bbR\times M)\)
defines a path \(\grd(s) \) in \(\mathcal{B}^\sigma (M)\): \(s\in
\bbR\mapsto \Pi^s\grd\in \mathcal{B}^\sigma (M)\). Conversely, a
 path \(\grd(\cdot)\co \bbR\to \mathcal{B}^\sigma (M)\) {\em together} with a
 \(\nabla_{s}\), the latter being the
 \(\frac{\partial}{\partial s}\)-component of an \(A\in
 \op{Conn}(E)\), 
 determines a \(\grd\in \mathcal{B}^\tau_{loc}(\bbR\times M)\). Denote
 the \(\nabla_s\) associated to \(\grd\) by \(\nabla_s^\grd\). This
 corresponds to the second term in \cite{KM}'s (4.10), and is 
 a lift of the  vector field  \(\frac{\partial}{\partial s}\) on the
 base \(\bbR\times M\) to the total space of the bundle \(E\). 
 
As \(\M (\bbR\times M)\subset \mathcal{B}^\tau(\bbR\times M)\), the
cocycles introduced in Part 1 may be defined using
\(\mathcal{B}^\tau(\bbR\times M)\) in place of
\(\mathcal{B}^\sigma_{loc}(\bbR\times M)\). 
 
Translations on \(\bbR\times M\) induce an 
 \(\bbR\)-action on \(\mathcal{B}^\sigma _{loc}(\bbR\times M)\) or 
 \(\mathcal{B}^\tau(\bbR\times M)\) in the following manner: For each
 \(a\in \bbR\), let 
\[
\tau_a\co \bbR\times M\to \bbR\times M
\] 
denote the map sending \((s, p)\in  \bbR\times M\) to
 \((s+a, p)\). For each \(\grd\in \mathcal{B}^\tau_{loc}(\bbR\times
M)\), its associated \(\nabla_s^\grd\) defines a lift of \(\tau_a\) to
a bundle automorphism of \(E\) (or equivalently of \(\bbS^+\)), denoted \(\tau_a^\grd\)
below. Let  
\[
\tau_a^{\cal B}\co \mathcal{B}^\sigma _{loc}(\bbR\times
M)\to \mathcal{B}^\sigma _{loc}(\bbR\times M)
\]
 be the map that sends
\(\grd\) to the pull-back of \(\grd\) (as a gauge-equivalence class of
\(\op{Conn}(\det\bbS^+)\times \Gamma^\sigma (\bbS^+)\)) via 
 \(\tau_{-a}^\grd\). Use the same notation, \(\tau_a^{\cal B}\),  to
 denote the similarly defined map from 
 \(\mathcal{B}^\tau(\bbR\times M)\) to itself. In particular,
 \(\tau_a^{\cal B}\) sends \(\grd(s)\) to \(\grd (s+a)\).  
Let 
\[\bar{i}\co
 \mathcal{B}^\sigma (M)\to \mathcal{B}^\tau(\bbR\times M)
\] 
be the embedding that sends a
 \(\grc\in  \mathcal{B}^\sigma (M)\) to \(\bbR\)-invariant element
 \(\grd_\grc\) with \(\grd_\grc(s)=\grc\) \(\forall s\in \bbR\). The
 fixed point set of the \(\bbR\)-action on \(\mathcal{B}^\sigma _{loc}(\bbR\times
M)\) is the image of \(\bar{i}\), and the action is free on the rest
of \(\mathcal{B}^\sigma _{loc}(\bbR\times
M)\). 

The \(\bbR\)-actions \(\tau_a^{\cal B}\) 
 preserve the subspace \(\mathcal{M}(\bbR\times M)\subset
 \mathcal{B}^\tau(\bbR\times M)\), together with  all of  its strata. 
The fixed point set of the aforementioned \(\bbR\)-action on
\(\mathcal{M}(\bbR\times M)\) is 
\[\mathcal{M}_0(\bbR\times M)\simeq \grC (M)\subset\mathcal{B}^\sigma (M),\] 
and the action is free on all higher dimensional strata
of \(\mathcal{M}(\bbR\times M)\). Thus, the orbit space
\(\mathcal{N}_k(M):=(\mathcal{M}_{k+1}(\bbR\times M)\backslash
\mathcal{M}_k(\bbR\times M))/\bbR\) is a \(k\)-dimensional manifold. 
 As explained in Section 16.1 of \cite{KM}, the spaces
\(\mathcal{N}_k(M)\) are compactified into a stratified manifold \(\mathcal{N}_k^+(M)\) by
adding ``(unparametrized) broken trajectories'', and the quotient
map \[q_\bbR\co \mathcal{M}_{k+1}(\bbR\times M)\backslash
\mathcal{M}_k(\bbR\times M)\stackrel{/\bbR}{\longrightarrow} \mathcal{N}_k(M)\]
extends to a
map, also denoted \(q_\bbR\),  between the stratified manifolds
\((\mathcal{M}_{k+1}(\bbR\times M)/\mathcal{M}_k(\bbR\times M))^+\) and
\(\mathcal{N}_k^+(M)\). (Recalling that each
\(\mathcal{M}_{k+1}(\bbR\times M)/\mathcal{M}_k(\bbR\times M)\) is a disjoint union of
moduli spaces of the form \(\mathcal{M}_z(\grc_-, \grc_+)\) or
\(\mathcal{M}^{\red}_z(\grc_-, \grc_+)\), the space
\((\mathcal{M}_{k+1}(\bbR\times M)/\mathcal{M}_k(\bbR\times M))^+\) above denotes the
disjoint union of their respective compactifications \(\mathcal{M}_z^+(\grc_-, \grc_+)\) or
\(\mathcal{M}^{\red\, +}_z(\grc_-, \grc_+)\). Correspondingly,
\(\mathcal{N}_k^+(M)\) is a disjoint union of compactified spaces of
the form \(\mathcal{N}_z^+(\grc_-, \grc_+):=(\mathcal{M}_z(\grc_-, \grc_+)/\bbR)^+\) or
\(\mathcal{N}^{\red\, +}_z(\grc_-, \grc_+):=
(\mathcal{M}^{\red}_z(\grc_-, \grc_+)/\bbR)^+\)).

\paragraph{(a)} {\em The \(U\)-map.} Fix \(p\in M\) and let
\(x=(0,p)\in \bbR\times M=X\).   Let 
\[
\pi_x\co\tilde{\mathcal{B}}_x^\tau(\bbR\times M):=\pi_x^{-1}\mathcal{B}^\tau(\bbR\times
M)\to \mathcal{B}^\tau(\bbR\times M)
\]
be the the principal
\(U(1)\)-bundle obtained by pulling back  \(\pi_x\co \tilde{\mathcal{B}}_{x,loc}^\sigma (\bbR\times M)\to
\mathcal{B}_{loc}^\sigma (\bbR\times M)\) via the embedding
\(\mathcal{B}^\tau(\bbR\times M)\hookrightarrow
\mathcal{B}_{loc}^\sigma (\bbR\times M)\).  Define
\(\pi_x\co\tilde{\mathcal{B}}^\tau_{x}(\bbR\times M)\to
\mathcal{B}^\tau(\bbR\times M)\) similarly.  Let \(\tilde{\cal B}^\sigma
 _p(M)\) be the 3-dimensional analog of \(\tilde{\mathcal{B}}_x^\sigma
 (X)\); namely,  \(\tilde{\cal B}^\sigma _p(M):={\cal C}^\sigma (M)/{\cal
  G}_p(M),\) with \({\cal G}_p(M)\subset C^\infty(M, U(1))\) being the
subgroup that consists of maps with value 1 at \(p\in M\). Then by
construction, the map \(\Pi^0\) lifts to a map \(\tilde{\Pi}^0\), that
fits into the following commutative diagram: 
\[\begin{CD}
\tilde{\cal B}^\tau_x(\bbR\times M) @>\tilde{\Pi}^0>>  \tilde{\cal B}^\sigma_p(M) \\
@V\pi_xVV @V\pi_pVV\\
{\cal B}^\tau(\bbR\times M) @>\Pi^0>>  {\cal B}^\sigma(M).
\end{CD}
\]
Regard \(\pi_p\co \tilde{\cal B}^\sigma _p(M)\to {\cal B}^\sigma (M)\) as a
principal \(U(1)\)-bundle and let \(\vartheta '_p\in
\Omega^1(\tilde{\mathcal{B}}_p^\sigma (M))\) denote a principal
\(U(1)\)-connection on this bundle. We will choose the principal
\(U(1)\) connection \(\vartheta\) from Part 1(a) to be of
the form \[\vartheta=(\tilde{\Pi}^0)^*\vartheta'_p.\] 
 By the unique continuation theorem (cf. \cite{KM}),  
\(\Pi^0|_{\mathcal{M}(\bbR\times M)}\) is an isomorphism, and we choose \(\vartheta'_p\)
to be integral over  \(\Pi^0\mathcal{M}_1(\bbR\times M)\), so that \(\vartheta \)
meets the integrability requirement (\ref{int-con}). Given  \(p\in M\), the 2-cocycle
\(\mu_U\) on \(\B^\sigma (M)\) used to define the \(U\)-action
is the Euler form \(\e'_p=(\pi_p)_* d\vartheta '_p\) of the bundle
\(\tilde{\B}^\sigma _p(M)\). It is straightforward to verify that indeed
\([\e'_p]=U\in H^*(\B^\sigma (M);\bbZ)\).
Let \[\mathring{U}_p:=\mathring{m}[\mathpzc{e}](\bbR\times M)\co
\mathring{C}(M)\to \mathring{C}(M),\] where
\(\mathpzc{e}=(\pi_x)_*(d\vartheta )=\Pi_0^*\e'_p\) as before.  We
call this degree \(-2\) cobordism chain map {\em the \(U\)-map
  associated to \(p\)} on the monopole Floer complex \(\mathring{C}(M)\).

It is desirable to express \(U_p\) in terms of integrals over the
unparametrized moduli spaces, 
\(\mathcal{N}_k(M)\), in a way similar to the formulae 
(\ref{eq:differential-loc}), (\ref{eq:differential-lo}) for the
differential \(\partial\) of the 
monopole Floer complex. For this purpose we digress to make some preparatory
observations. 

Let \((a, b)\in
 \bbR\times \bbR\mapsto\tau_a\times \tau_b^{\cal B}\) be the product
 \(\bbR\times \bbR\)-action on \((\bbR\times M)\times \mathcal{B}^\sigma
 _{loc}(\bbR\times M)\), i.e. the base space of the bundle \(\pmb{\cal
   E}(\bbR\times M)\), and use the same notation to denote the lift
 via \(\tilde{A}\) of
 this \(\bbR\times \bbR\)-action to the total space, \(\pmb{\cal
   E}(\bbR\times M)\). By
 construction, the tautological \(\tilde{A}\) and \(\tilde{\alpha}\)
 (when defined)  on \(\pmb{\cal E}(\bbR\times M)\) are invariant under pull-back of the anti-diagonal
 \(\bbR\)-action; namely, 
\[
(\op{Id}\times
 \tau_a^{\cal B})^*\tilde{\alpha}=(\tau_a\times
 \op{Id})^*\tilde{\alpha},\,\,  \text{ and
 similarly for \(\tilde{A}\)}.
\] 
Let 
 \(\bbR\grd\subset\mathcal{B}^\sigma _{loc}(\bbR\times M) \) denote
 the \(\bbR\)-orbit through a \(\grd\in \mathcal{B}^\sigma
 _{loc}(\bbR\times M)\) and let \(\hat{p}:=\bbR\times \{p\}\subset
 \bbR\times M\) denote the \(\bbR\)-orbit through \(x=(0,p)\). Then
 the aforementioned anti-diagonal \(\bbR\)-action on  \(\pmb{\cal
   E}(\bbR\times M)\) defines a bundle isomorphism \(\iota_\Delta \) between \(\pmb{\cal
   E}(\bbR\times M)|_{\hat{p}\times \{\grd\}}\simeq E|_{\hat{p}}\) 
 and \(\pmb{\cal
   E}(\bbR\times M)|_{\{x\}\times
   \bbR\grd}\simeq\mathcal{E}_x(X)|_{\bbR\grd}\), and parallel transports
 via \(\tilde{A}\) along the two paths \(\hat{p}\times \{\grd\}\)
 and \(\{x\}\times \bbR\grd\) in \((\bbR\times M)\times
 \mathcal{B}^\sigma _{loc}(\bbR\times M)\) are identified under
 \(\iota_\Delta \). Note that the connection
 \(\tilde{A}|_{\hat{p}\times \{\grd\}}\) on
 \(\pmb{\cal E}(\bbR\times M)|_{\hat{p}\times \{\grd\}}\simeq E|_{\hat{p}}\) is
 precisely the restriction of \(\nabla_s^\grd\) to
 \(E|_{\hat{p}}\). This is identified via \(\iota_\Delta\) with  the connection
 \(\tilde{A}|_{\{x\}\times \bbR\grd}\) on  \(\pmb{\cal E}(\bbR\times
 M)|_{\{x\}\times \bbR\grd}\simeq{\cal E}_x(X)|_{\bbR\grd}\), which 
 corresponds to \(\tilde{\vartheta}|_{\bbR\grd}\) on
 \(\tilde{\mathcal{B}}^\sigma _{x,loc}(\bbR\times
 M)|_{\bbR\grd}\). (Recall the definition of \(\tilde{\vartheta}\)
 from Part 1(a).) Observe, by the way, that  \({\cal E}_x(X)\simeq \pmb{\cal E}(\bbR\times
 M)|_{\{x\}\times \mathcal{B}^\sigma _{loc}(\bbR\times M)}\) admits an
 \(\bbR\)-action \(a\in \bbR\mapsto \op{Id}\times \tau_a^{\cal B}\),
 and the associated \(\bbR\)-action on \(\tilde{\mathcal{B}}^\sigma
 _{x,loc}(\bbR\times M)\) is precisely the lift of the \(\bbR\)-action
 on \(\mathcal{B}^\sigma _{loc}(\bbR\times M)\) via
 \(\tilde{\vartheta}\). Namely, denoting the lift of \(\tau_a^{\cal
   B}\) by the same notation, we have the commutative diagram: 
\[
\begin{CD}
\tilde{\cal B}^\sigma _{x,loc}(\bbR\times M) @>\tau_a^{\cal B}>> \tilde{\cal B}^\sigma _{x,loc}(\bbR\times M)  \\
@V\pi_xVV @V\pi_xVV\\
{\cal B}^\sigma _{loc}(\bbR\times M) @>\tau_a^{\cal B}>>{\cal B}^\sigma _{loc}(\bbR\times M).
\end{CD} 
\]

Let \(\bar{\bbR}:=[-\infty,
 \infty]\supset\bbR\), and so in the present setting
 \(\bar{X}=\bar{\bbR}\times M\). Suppose \(\grd\in
 \mathcal{B}^\tau(\bbR\times M)\). Then by definition, \(\tau_a^{\cal B}(\grd)\)
 converges as \(a\to \pm \infty\) (in the subspace topology of \(\mathcal{B}^\tau(\bbR\times
 M)\subset \mathcal{B}^\sigma _{loc}(\bbR\times M)\)) respectively
 to 
\[
\tau_{\pm \infty}^{\cal B}(\grd)=\bar{i} (\Pi^{\pm\infty}\grd).
\] 
Let \(\bar{\bbR}\grd\subset \mathcal{B}^\tau(\bbR\times M)\) denote
 \(\{\tau_s^{\cal B}(\grd)\}_{s\in \bar{\bbR}}\). Thus, 
the paths \(\hat{p}\times \{\grd\}\)
 and \(\{x\}\times \bbR\grd\) in \((\bar{R}\times M)\times
 \mathcal{B}^\tau(\bbR\times M)\) extend respectively to arcs
 \((\bar{\bbR}\times \{p\})\times \{\grd\}\) and \(\{x\}\times
 \bar{\bbR}\grd\). The previously introduced bundle isomorphism
 \(\iota_\Delta\) extends to define a bundle isomorphism 
 \[
\iota_\Delta \co \pmb{\cal E}(\bar{\bbR}\times M)|_{(\bar{\bbR}\times\{p\})\times
   \{\grd\}}\simeq E|_{\bar{\bbR}\times\{p\}}\stackrel{\sim}{\to }\pmb{\cal
   E}(\bar{\bbR}\times M)|_{\{x\}\times
   \bar{\bbR}\grd}\simeq\mathcal{E}_x(\bar{\bbR}\times
 M)|_{\bar{\bbR}\grd}.
\] 
The assumption that \(\grd\in \mathcal{B}^\tau(\bbR\times
 M)\) also ensures that parallel transport via \(\nabla_s^\grd\) along
 \(\bar{\bbR}\times\{p\}\) gives  a well-defined unitary holonomy map
 \[
\op{hol}_{\hat{p}}^E(\grd)\in \op{Hom} (E|_{(-\infty, p)},
E|_{(\infty, p)})\simeq \op{Hom} (\pmb{\cal E}|_{(-\infty, p)\times \{\grd\}},\pmb{\cal
E}|_{(\infty, p)\times \{\grd\}}).
\]
 As \(\iota_\Delta^*\) preserves \(\tilde{A}\),
the holonomy of \(\tilde{A}\) along \(\{x\}\times
   \bar{\bbR}\grd\) also gives a well-defined unitary element agreeing
   with 
   \(\iota_\Delta \circ \op{hol}_{\hat{p}}^E(\grd)\circ \iota_\Delta ^{-1}\) in \(\op{Hom} (\mathcal{E}_x(\bbR\times
 M)|_{\{\bar{i}(\grc_-)\}}, \mathcal{E}_x(\bbR\times
 M)|_{\{\bar{i}(\grc_+)\}})\), where
 \(\grc_\pm:=\Pi^{\pm\infty}(\grd)\). The space of unitary elements in \(\op{Hom} (\mathcal{E}_x(\bbR\times
 M)|_{\{\bar{i}(\grc_-)\}}, \mathcal{E}_x(\bbR\times
 M)|_{\{\bar{i}(\grc_+)\}})\) is precisely 
 \(\tilde{\cal B}_x^\tau(\bbR\times M)|_{\{\bar{i}(\grc_-)\}}\times _{U(1)}\tilde{\cal
     B}_x^\tau(\bbR\times M)|_{\{\bar{i}(\grc_+)\}}\simeq \tilde{\cal
       B}^\sigma _p(M)|_{\{\grc_-\}}\times _{U(1)}\tilde{\cal
       B}^\sigma _p(M)|_{\{\grc_+\}}\). This is the fiber over
     \((\grc_-, \grc_+)\) of the
     \(U(1)\)-bundle: 
\[
\begin{split}
\tilde{\cal
  B}^\sigma _p(M)\times _{U(1)}\tilde{\cal
  B}^\sigma _p(M)& =(\tilde{\cal B}^\sigma _p(M)\times \tilde{\cal
  B}^\sigma _p(M))/\text{diagonal \(U(1)\) actions}\\
& \qquad \quad \stackrel{\pi_{p-p}}{\longrightarrow}{\cal B}^\sigma (M)\times {\cal B}^\sigma (M),
\end{split}
\]
where \(\pi_{p-p}\) is the quotient map by the residual
\(U(1)\)-action.  Let 
\[
\op{hol}_{\hat{p}}\co    {\cal B}^\tau(\bbR\times M)\to \tilde{\cal
       B}^\sigma _p(M)\times _{U(1)}\tilde{\cal
       B}^\sigma _p(M)
\]
be the map that sends \(\grd\) to the element in \( \tilde{\cal
  B}^\sigma _p(M)\times _{U(1)}\tilde{\cal
  B}^\sigma _p(M)\) corresponding to to \(\iota_\Delta \circ\op{hol}^E_{\hat{p}}(\grd)\circ
\iota_\Delta ^{-1}\). This map is a lift of
     the map \(\Pi^\partial=\Pi^{-\infty}\times \Pi^{\infty}\co {\cal
     B}^\tau(\bbR\times M)\to B^\sigma (M)\times  B^\sigma (M)\)
in the sense that \(\pi_{p-p}\circ\op{hol}_{\hat{p}}=\Pi^\partial\).
Meanwhile, letting \(\tilde{\Pi}^a=\tilde{\Pi}^0\circ\tau_a\), 
the map \(\tilde{\Pi}^{\partial}:=\tilde{\Pi}^{-\infty}\times \tilde{\Pi}^{\infty}\co \tilde{\cal
  B}^\tau _x(\bbR\times M)\to \tilde{\cal
  B}^\sigma _p(M)\times \tilde{\cal
  B}^\sigma _p(M)\) is in turn a lift of \(\op{hol}_{\hat{p}}\) under the
quotient map 
\[
\pi_\Delta \co \tilde{\cal
  B}^\sigma _p(M)\times \tilde{\cal
  B}^\sigma _p(M) \to \tilde{\cal
  B}^\sigma _p(M)\times _{U(1)}\tilde{\cal  B}^\sigma _p(M).
\]

A choice of \(\vartheta '_p\) determines a
principal \(U(1)\)-connection on the bundle \(\pi_{p-p}\co \tilde{\cal
  B}^\sigma _p(M)\times _{U(1)}\tilde{\cal
  B}^\sigma _p(M)\to {\cal
  B}^\sigma (M)\times {\cal
  B}^\sigma (M)\), which we denote by \(\vartheta '_{p-p}\). Since
\(\vartheta '_p\) is integral over \(\Pi_0\M_1(\bbR\times M)\), it also
determines a trivialization \(\rho_{p-p}\co
\ul{U(1)}\stackrel{\sim}{\to} \tilde{\cal
  B}^\sigma _p(M)\times _{U(1)}\tilde{\cal
  B}^\sigma _p(M)|_{\Pi_0\M_1(\bbR\times M)\times \Pi_0\M_1(\bbR\times
  M)}\) of the bundle \(\pi_{p-p}\) over \(\Pi_0\M_1(\bbR\times M)\times \Pi_0\M_1(\bbR\times
  M)\). As \(\pi_{p-p} (\M (\bbR\times M))\subset \grC(M)\times
  \grC(M)\subset \Pi_0\M_1(\bbR\times M)\times \Pi_0\M_1(\bbR\times
  M)\), combining the
trivialization \(\rho_{p-p}\) with \(\tilde{\Pi}^\partial\), we get a map 
\[
\op{h}_{\hat{p}}\co \M(\bbR\times M)\to U(1)=\bbR/\bbZ.
\] 
Observe that the maps \(\tilde{\Pi}^\partial\), \(\Pi^\partial\),
\(\Pi^\infty\), \(\tilde{\Pi}^\infty\), \(\Pi^{-\infty}\),
\(\tilde{\Pi}^{-\infty}\),  \(\op{h}_{\hat{p}}\) are all invariant
under the respective \(\bbR\)-actions on their domains, and therefore descend to define
maps from the orbit spaces under the \(\bbR\)-actions. Our convention
is to denote the
corresponding maps from \(\mathcal{B}^\tau (\bbR\times M)/\bbR\),
\(\tilde{\mathcal{B}}^\tau (\bbR\times M)/\bbR\), or
\(\M (\bbR\times M)/\bbR\) by adding underlines to the notations. For example,
\(\op{h}_{\hat{p}}=\ul{\op{h}}_{\hat{p}}\circ q_\bbR\). By construction, we have 
\[\begin{split}
(\op{hol}_{\hat{p}})^*\vartheta'_{p-p}& =-d\op{h}_{\hat{p}}\quad
\text{and} \\
(\ul{\op{hol}}_{\hat{p}})^*\vartheta'_{p-p} & =-d\ul{\op{h}}_{\hat{p}}\quad
\text{over \({\cal N}_k(M)\).}
\end{split}
\]
Let \(\rho_{\vartheta '_p}\co
\ul{\bbC}\to \tilde{\cal B}_p(M)|_{\Pi_0\M_1(\bbR\times M)}\) be a trivialization
inducing \(\vartheta'_p\), and use the notation to denote the
associated trivialization of \({\cal E}_p(M)|_{\Pi_0\M_1(\bbR\times M)}\),
\({\cal E}_p(M)\) being the hermitian
line bundle associated to \(\tilde{\cal B}_p(M)\). Using \(\iota_\Delta \) to identify \(E|_{(\pm\infty,
  p)}\) respectively with \({\cal E}_p(M)|_{\{\grc_\pm\}}\), we have:
\[
e^{2\pi i\op{h}_{\hat{p}}(\grd)}=(\rho_{\vartheta'_p})^{-1}\circ\op{hol}_{\hat{p}}^E(\grd)\circ
\rho_{\vartheta'_p}\in \bbC^*. 
\]
Meanwhile, given \(\grd\in {\cal B}^\tau (\bbR\times M)\) and an
arbitrary \(\tilde{\grd}\in \pi_x^{-1}(\grd)\), 
\[
\int_{\bbR\tilde{\grd}}\vartheta
=\int_{\bbR\grd}\theta=:-h_{\hat{p}}(\grd)\in \bbR
\]
where \(\theta\) is as in Part 1(a)'s  (\ref{def:theta}), with 
\(x\)  set to be \(\dot{p}:=(0, p)\in \bbR\times M\). In particular, when \(\grd\in
\M_1 (\bbR\times M)\), 
\[
\op{h}_{\hat{p}}(\grd)=h_{\hat{p}}(\grd)\quad\mod
\bbZ.\]
Like \(\op{h}_{\hat{p}}\), the function \(h_{\hat{p}}(\grd)\co \M
(\bbR\times M)\to \bbR\) is invariant under the \(\bbR\)-action on
\(\M\), and hence induces a function \(\ul{h}_{\hat{p}}\co
\N(M):=\bigcup_k \N_k(M)\to \bbR\), with
\[
\text{\(\ul{\op{h}}_{\hat{p}}=\ul{h}_{\hat{p}}\) \(\mod \bbZ\)
over \(\mathcal{N}_0(M)\).}
\]
This function can be used to write
\begin{equation}\label{theta-h}
[\theta_{\dot{p}}]_{\M (\bbR\times M)}=-\sum_{\ul{\grd}\in \N_0(M)}\ul{h}_{\hat{p}}(\ul{\grd})\mu^1_{\ul{\grd}},
\end{equation}
where \(\{\mu^1_{\ul{\grd}}\}_{\ul{\grd}\in \N_0(M)}\) is the
canonical basis for \(C_{\M(\bbR\times M)}^1\). 
One may extend by continuity both 
\(\ul{\op{h}}_{\hat{p}}\) and \(\ul{h}_{\hat{p}}\) to the spaces of broken trajectories
 \(\N^+_k(M)\) as respectively  \(\bbR/\bbZ\)-
 and \(\bbR\)-valued functions. With this extension we have
\(\ul{h}_{\hat{p}}(\pmb{\grd})=\sum_{\ul{\grd}_i}\ul{h}_{\hat{p}}(\ul{\grd}_i)\)
for \(\pmb{\grd}=\{\ul{\grd}_i\}_i\in \N^+_k(M)\), and 
 \(\ul{\op{h}}_{\hat{p}}=\ul{h}_{\hat{p}}\) \(\mod \bbZ\)
over the 0-dimensional strata \((\mathcal{N}_k^+(M))_0\). Consequently, 
\[
\op{u}_p:=[d\ul{\op{h}}_{\hat{p}}]_{\N^+_1(M)}-\delta[\ul{h}_{\hat{p}}]_{\N^+_1(M)}\in
C^{1;\bbZ}_{\N_1^+(M)}\subset C^{1;\bbR}_{\N_1^+(M)};
\]
namely, \([\op{u}_{p}]_{\N^+_1(M)}\) is an integral correction
of
\([-(\ul{\op{hol}}_{\hat{p}})^*\vartheta'_{p-p}]_{\N^+_1(M)}=[d\ul{\op{h}}_{\hat{p}}]_{\N^+_1(M)}\).

We next express \(\mathring{U}_p\) in terms of integrals of \(\op{u}_p\)
over
\(\mathcal{N}_1^+(M)\). According to (\ref{eq:m-coeff}) and (\ref{eq:m-coeff-loc}), the coefficients of
\(\mathring{U}_p \) take the form of 
\[\begin{split}\langle \, \mathpzc{e}, 
\mathcal{M}_{2,z}(X; \grc_-, \grc_+)\, \rangle & =\langle \, \mathpzc{e}, 
\bar{\mathcal{M}}_{2,z}(X; \grc_-, \grc_+)\, \rangle\quad \text{ or}\\ 
\langle \, \mathpzc{e},
\mathcal{M}^{\text{red, 2}}_z(X; \grc_-, \grc_+)\, \rangle & =\langle \, \mathpzc{e},
\bar{\mathcal{M}}^{\text{red, 2}}_z(X; \grc_-, \grc_+)\, \rangle ,
\end{split}
\] 
where \(X=\bbR\times M\). Let
\(\bar{\mathcal{M}}\) be one of the compactified moduli spaces 
\(\bar{\mathcal{M}}_{2,z}(X; \grc_-, \grc_+)\) or
\(\bar{\mathcal{M}}^{\text{red, 2}}_z(X; \grc_-, \grc_+)\) named above. This
is a 2-dimensional
stratified submanifold of \(\mathcal{M}(X)\). 
Let
\(\grM=\bar{\cal M}\backslash \mathcal{M}_1(X)\) denote the
top-dimensional stratum of \(\bar{\cal M}\), and let \(\mathcal{M}^+\)
be the larger compactification of \(\grM\) by adding (parametrized)
broken trajectories. The latter carries a stratification of the form 
\[
\emptyset\subset (\mathcal{M}^+)_0\subset (\mathcal{M}^+)_1\subset (\mathcal{M}^+)_2=\mathcal{M}^+.
\]
Meanwhile, \(\grM\) consists of \(\bbR\)-orbits; let
\(\grN:=\grM/\bbR\subset \mathcal{N}_1(M)\), and use \(\mathcal{N}^+\subset \mathcal{N}^+_1(M)\) to
denote the compactification of \(\grN\) by adding unparametrized
broken trajectories. It is stratified as \(\emptyset\subset
(\mathcal{N}^+)_0\subset (\mathcal{N}^+)_1=\mathcal{N}^+\). The strata
of \(\mathcal{M}_z^+(\grc_-, \grc_+)\), 
\(\mathcal{M}^{\red\, +}_z(\grc_-, \grc_+)\), \(\bar{\mathcal{M}}_z(\grc_-, \grc_+)\), 
\(\bar{\mathcal{M}}^{\red}_z(\grc_-, \grc_+)\), \(\mathcal{N}_z^+(\grc_-, \grc_+)\) or
\(\mathcal{N}^{\red\, +}_z(\grc_-, \grc_+)\) are described in
\cite{KM}'s (24.27), (24.28) and Proposition 24.6.10. Applied to the
case 
under discussion, this entails: 
\BTitem\label{M-strata}
\item \((\bar{\mathcal{M}})_0=\bar{i}(\{\grc_-, \grc_+\}) \subset \bar{i}(\grC(M))\). 
 \item \((\mathcal{N}^+)_0\) consists of finitely
many once- or twice-broken trajectories. We denote such a broken
trajectory in the form 
\(\pmb{\grd}=(\grd_i)_i\), where each \(\grd_i\in \mathcal{N}_0(M)\)
and \(i\in \{1,2\}\) or \(\{1,2,3\}\).
\item The 1-dimensional strata of \(\bar{\M}\) consists of
  \(\bbR\)-orbits in \(\mathcal{B}^\sigma (X)\). More precisely,
  \((\bar{\M})_1\backslash (\bar{\M})_0=\bigcup_{\pmb{\grd}=(\grd_i)_i\in
    (\mathcal{N}^+)_0}\bigcup_i\bbR\grd_i\). 
\item \((\M^+)_1\) is a union of two parts: 
\[
(\M^+)_1=\grr^{-1}(\bar{\M})_0\cup \grr^{-1} ((\bar{\M})_1\backslash
(\bar{\M})_0). 
\]
In the first part, \(\grr^{-1}\bar{\grc}\simeq (\N^+)_1\) for each
\(\bar{\grc}\in(\bar{\M})_0\subset \grC(M)\), while on the second part,
\(\grr\) restricts to an  isomorphism from \(\grr\co \grr^{-1}((\bar{\M})_1\backslash
(\bar{\M})_0)\) to \((\bar{\M})_1\backslash (\bar{\M})_0\).
\ETitem
Our strategy  to compute \(\langle \, \mathpzc{e},
\bar{\mathcal{M}}\rangle\) is to 
introduce a map
\(\tilde{\varsigma }\co \mathcal{M}^+ \to \tilde{\mathcal{B}}^\sigma _{x,loc}(X)\)
so that the following diagram commutes: 
\begin{equation}\label{CD-M+}\begin{CD}
\mathcal{M}^+@>\tilde{\varsigma} >>\tilde{\mathcal{B}}^\sigma_{x,loc}(X)\\
@ V\grr VV @ V\pi_x VV\\
\bar{\cal M}  @>\varsigma>> {\cal B}^\sigma  _{loc}(X),
\end{CD}
\end{equation} 
where \(\varsigma \co \mathcal{M}_2(X)\to
\mathcal{B}^\sigma  _{loc}(X)\) denotes  the embedding. 
This means that \(\tilde{\varsigma}|_{\grM}\) is then a lift the embedding
\(\varsigma|_\grM\) under \(\pi_x\). We choose this lifting so that
\(\tilde{\varsigma}(\grM)\subset \tilde{\cal B}^\tau (\bbR\times
M)\subset \tilde{\mathcal{B}}^\sigma_{x,loc}(X)\)
is tangent to the \(\bbR\)-action. Such a choice is specified in turn
by a lift \(\tilde{\varsigma}_{\cal N}\) of \(\grN\subset {\cal
  N}_1(M)\) to \(\tilde{\cal N}_1(M)\). As an extension of
\(\tilde{\varsigma}|_\grM\), \(\tilde{\varsigma }\)'s image is also
tangent to the \(\bbR\)-action on \(\tilde{\cal B}^\tau (\bbR\times
M)\). 
With \(\tilde{\varsigma }\) chosen, we then write 
\begin{equation}\label{value}
\begin{split}& \langle\mathpzc{e}, \bar{\cal M}\rangle=\langle\mathpzc{e},
\grM\rangle \\
& \quad =\langle\tilde{\varsigma }^*\mathpzc{e},
\grM\rangle =\langle\tilde{\varsigma }^*\mathpzc{e},
\mathcal{M}^+\rangle \\
& \quad =\langle\tilde{\varsigma }^*\vartheta, \partial [\mathcal{M}^+]\rangle=\langle\tilde{\varsigma }^*\vartheta, [(\mathcal{M}^+)_1]\rangle.
\end{split}
\end{equation}
using \cite{KM}'s Theorem 24.7.2 and Lemma
21.3.1. 
By (\ref{M-strata}), the last term above is written as a sum
\[\begin{split}
\langle\tilde{\varsigma }^*\vartheta, [(\mathcal{M}^+)_1]\rangle &=\langle\tilde{\varsigma }^*\vartheta, \grr^{-1}(\bar{\M})_0\rangle+\langle\tilde{\varsigma }^*\vartheta, \grr^{-1} (\bar{\M}_1\backslash
\bar{\M}_0)\rangle\\
& =\langle \tilde{\varsigma }^*_{\N}(\ul{\tilde{\Pi}}^-)^*\vartheta'
_p, \N^+\rangle-\langle \tilde{\varsigma }^*_{\N}(\ul{\tilde{\Pi}}^+)^*\vartheta'
_p, \N^+\rangle+\langle\theta_{\dot{p}}, \bar{\M}_1\backslash
\bar{\M}_0\rangle\\
& =-\langle (\ul{\op{hol}}_{\hat{p}})^*\vartheta'_{p-p}, \N^+\rangle-\sum_{\pmb{\grd}\in
    (\mathcal{N}^+)_0}\op{sign}(\pmb{\grd})\, 
  \ul{h}_{\hat{p}}(\pmb{\grd}).\\
\end{split}
\]
Note that the first term in the last line above is independent of the
choice of \(\tilde{\varsigma }_{\N}\); it is also independent of the
choice of \(\vartheta '_p\), since by (\ref{CD-M+}), \(\Pi_{\N}\grN\)
lies in a fiber of \(\tilde{\B}^\sigma (M)\times
_{U(1)}\tilde{\B}^\sigma (M)\), over which \(\vartheta'_p\) is the
standard \(U(1)\)-invariant volume form generating
\(H^1(U(1);\bbZ)\). To summarize, we have 
\begin{equation}\label{e-hol}
\begin{split}
\langle\mathpzc{e}, \bar{\cal M}\rangle
& = \langle \op{u}_{p}, \N^+\rangle;\\
\langle d\ul{\op{h}}_{\hat{p}}, \N^+\rangle& =\langle\mathpzc{e}-\delta[\theta_{\dot{p}}]_{\M(\bbR\times M)}, \bar{\cal M}\rangle.
\end{split}
\end{equation}

\paragraph{(b)} {\em The \(\grt_i\)-maps.} Given an embedded oriented
circle \(\ul{\gamma }\subset M\), one may define a
3-dimensional counterpart of the 1-form \(\theta_\gamma\) in Part 1(b)
above: Let 
\[
\op{h}'_{\ul{\gamma }}\co \B^\sigma (M)\to U(1)=\bbR/\bbZ 
\]
be the map
sending a \(\grd\in \B^\sigma (M)\) to the holonomy along \(\ul{\gamma}\)
of \(A\in
\op{Conn}(E)\) associated to \(\grd\), and set 
\[
\theta'_{\ul{\gamma}}=d\op{h}'_{\ul{\gamma }}\in \Omega^1(\B^\sigma (M)).
\]
By construction \(\theta'_{\ul{\gamma}}\) is closed and its cohomology
class \([\theta'_{\ul{\gamma}}]\in H^1(\mathcal{B}^\sigma (M);\bbZ)\) 
equals \([\ul{\gamma}]\in H_1(M;\bbZ)/\op{Tors}\) under the
isomorphism \[H_1(M;\bbZ)/\op{Tors}\simeq
H^1(H^1(M;\bbR)/H^1(M;\bbZ);\bbZ)\simeq H^1(\mathcal{B}^\sigma
(M);\bbZ).
\]
Let \(\dot{\gamma}:=\{0\}\times \ul{\gamma}\subset
\bbR\times M\). Then \(\theta_{\dot{\gamma
  }}=\Pi_0^*\theta'_{\ul{\gamma}}\) and \([\theta_{\dot{\gamma
  }}]_{\M(X)}\in C^{1;\bbR}_{\M(X)}\). Let 
\[
\uu_{\dot{\gamma
  }}=[\theta_{\dot{\gamma  }}]_{\M(X)}-\delta [h_{\dot{\gamma
  }}]_{\M(X)}\in C^{1;\bbZ}_{\M(X)}
\]
be an integral correction of \([\theta_{\dot{\gamma
  }}]_{\M(X)}\in C^{1;\bbR}_{\M(X)}\), as described in Part 1(b). 
In the present case, \(\M_0(X)\stackrel{\Pi_0}{\simeq}\grC(M)\subset\B^\sigma (M)\), and
the function \(h_{\dot{\gamma }}\co \M_0(X)\to \bbR\) in Part 1(b)
takes the form of \(\Pi_0^* h'_{\ul{\gamma }}\), where
\[
h'_{\ul{\gamma }}\co  \grC(M)\to \bbR
\]
 is a lift of
\(\op{h}'_{\ul{\gamma }}|_{\grC(M)}\co \grC(M)\to
\bbR/\bbZ\). Noting that the strata of \(\M(X)\) are \(\bbR\)-spaces
in this product cobordism case, we have 
\[\begin{split}
h_{\dot{\gamma }}& =\Pi_0^* h'_{\ul{\gamma }}=\Pi_s^*
h'_{\ul{\gamma }}=h_{\{s\}\times\ul{\gamma }};\\
[\theta_{\dot{\gamma  }}]_{\M(X)}&=[\Pi_0^*\theta'_{\ul{\gamma  }}]_{\M(X)}=[\Pi_s^*\theta'_{\ul{\gamma  }}]_{\M(X)}=[\theta_{\{s\}\times\ul{\gamma }}]_{\M(X)}
\end{split}
\]
for all \(s\in [-\infty, \infty]\). Thus, for our purpose
\(\dot{\gamma}\) may be taken to be \(\{s\}\times \ul{\gamma
}\subset \bbR\times M\) for arbitrary \(s\).

Let \(\grt:=[\ul{\gamma }]\in H_1(M;\bbZ)/\op{Tors}\) and let 
\[
\textsc{p}_{\grt}\co \hat{\cal
  B}^\sigma_{\grt}(M)\to \mathcal{B}^\sigma (M)
\] 
be the \(\bbZ\)-covering of
\(\mathcal{B}^\sigma(M)\) with \(\pi_1(\hat{\mathcal{B}}^\sigma _\grt)\subset
\pi_1(\mathcal{B}^\sigma (M))\) being the kernel of the map \(\grt\co
\pi_1(\mathcal{B}^\sigma (M))\simeq H^1(M;\bbZ)\to \bbZ\). The function
\(\op{h}'_{\ul{\gamma}}\co \mathcal{B}^\sigma (M)\to
U(1)=\bbR/\bbZ\) lifts to an \(\bbR\)-valued function 
\[
\hat{h}_{\ul{\gamma}}\co \hat{\mathcal{B}}^\sigma _{\grt}(M)\to\bbR.
\] 
This lift \(\hat{h}_{\ul{\gamma}}\)
is unique modulo addition by constant, \(\bbZ\)-valued
functions, and it can be fixed by choosing a base point
\(\hat{\grc}_0\in \hat{\mathcal{B}}^\sigma _{\grt}(M)\) with
\(\op{h}'_{\ul{\gamma }}(\grc_0)=0\mod \bbZ\), where 
\(\grc_0:=\textsc{p}_{\grt}(\hat{\grc}_0)\). Namely, let 
\(\hat{h}_{\ul{\gamma}}\) be such that
\[
\hat{h}_{\ul{\gamma}}(\hat{\grc}_0)=0\in \bbR.
\] 
Given an arc \(\grb\co [0,1]\to \B^\sigma  (M)\), the difference
\(\hat{h}_{\ul{\gamma}}(\hat{\grb}(1))-\hat{h}_{\ul{\gamma}}(\hat{\grb}(0))\)
    takes the same value for any lift \(\hat{\grb}\co
[0,1]\to \hat{\B}_\grt^\sigma  (M)\) of \(\grb\); we denote this
value by \(\Delta_\grb\hat{h}_{\ul{\gamma}}\in \bbR\). It depends only
on \(\op{h}'_{\ul{\gamma}}\) and not on the choice of the lift
\(\hat{h}_{\ul{\gamma}}\). In particular,
any \(\ul{\grd}\in \N (M)\) defines an arc \(\grb_{\ul{\grd}}\) in
\(\B^\sigma  (M)\), and we adopt the short-hand 
\(\Delta_{\ul{\grd}}\hat{h}_{\ul{\gamma}}:=\Delta_{\grb_{\ul{\grd}}}\hat{h}_{\ul{\gamma}}\). This value only depends on the relative homotopy class of
\(\ul{\grd}\). Observing
  that \(\int_{\bbR\ul{\grd}}\theta _{\dot{\gamma
    }}=\Delta_{\ul{\grd}}\hat{h}_{\ul{\gamma}}\), we have 
\[
[\theta _{\dot{\gamma }}]_{\M(\bbR\times
  M)}=\sum_{\ul{\grd}\in\N_0(M)}
(\Delta_{\ul{\grd}}\hat{h}_{\ul{\gamma}}) \mu^1_{\ul{\grd}}\in
C^{1;\bbR}_{\M(\bbR\times M)}.
\]
An integral correction \(\uu_{\dot{\gamma }}\) of \([\theta _{\dot{\gamma }}]_{\M(\bbR\times
  M)}\) can be written in a similar fashion by replacing the function
\(\hat{h}_{\ul{\gamma}}\co \hat{\mathcal{B}}^\sigma
_{\grt}(M)\to\bbR\) in the preceding discussion by a modified function 
\[
\op{x}_{\ul{\gamma}}\co \hat{\mathcal{B}}^\sigma
_{\grt}(M)\to\bbR,
\]
where \(\op{x}_{\ul{\gamma}}=\hat{h}_{\ul{\gamma}}-\textsc{p}_\grt^*\varepsilon'_{\ul{\gamma}}\)
for a function \(\varepsilon'_{\ul{\gamma}}\co
\mathcal{B}^\sigma(M)\to\bbR\) satisfying 
\(\varepsilon'_{\ul{\gamma}}|_{\grC(M)}=h'_{\ul{\gamma }}\).

Returning to the subject of \(\mathbf{A}_\dag\) actions, take  \(\ul{\gamma }=\gamma_i\subset M\) to be one that
represents \(\grt_i\in H_1(M;\bbZ)/\op{Tors}\). The
{\em 
\(\grt_i\)-map associated to \(\gamma_i\)} is defined to be 
\[
\mathring{\grm}_{\grt_i}=\mathring{\grm}_{\gamma_i}:=\mathring{m}[\uu_{\dot{\gamma
  }_i}](\bbR\times M)\co \mathring{C} (M)\to \mathring{C}(M).
\]
This corresponds to the 1-cocycle
\(\mu_{\grt_i}=(\textsc{p}_{\grt_i})_*(d\op{x}_{\ul{\gamma}})\) on
\(\B^\sigma (M)\). 

It will also be handy to introduce an 
analog of Part (a)'s \(\op{u}_p\) (cf. (\ref{e-hol})): Given
\(\ul{\grd}\in \mathcal{N}(M)\), let \(\Delta
_{\ul{\grd}}\op{x}_{\gamma _i}\) be defined in the same way as \(\Delta
_{\ul{\grd}}\hat{h}_{\gamma _i}\) above. Let \(\op{u}_{\gamma_i}\)
denote the function on \(\mathcal{N}(M)\) that sends each
\(\ul{\grd}\in \mathcal{N}(M)\) to
\(\Delta_{\ul{\grd}}\op{x}_{\gamma _i}\). Note that
\(\op{u}_{\gamma_i}\) is \(\bbZ\)-valued, and hence defines a class in
\(C^{0;\bbZ}_{\mathcal{N}(M)}\), denoted by the same notation. The
coefficients appearing in the formula for \(\grm_{\grt_i}\) then may
be re-expressed as integrals of \(\op{u}_{\gamma_i}\) over
\(\mathcal{N}(M)\):
\[
\langle \uu_{\dot{\gamma }_i}, \M\rangle=\langle \op{u}_{\gamma_i}, \N\rangle,
\]
where \(\M\) is a 1-dimensional stratum in \(\M (\bbR\times M)\) and
\(\N=\M/\bbR\) is the corresponding stratum in \(\N (M)\). In 
general, we use the notation \(\mathring{n}[\op{u}]: =\mathring{m}[\uu]\)
when \(X=\bbR\times M\) is a product cobordism and the coefficients in
the formula for 
\(\mathring{m}[\uu]\) may be expressed as integrals of \(\op{u}\in
C^k_{\mathcal{N}(M)}\) over \(\mathcal{N}(M)\) in the way described
above. E.g. we write 
\begin{equation}\label{m-n}
\mathring{\grm}_{\gamma _i}=\mathring{n}[\op{u}_{\gamma _i}]; \quad
\mathring{U}_p=\mathring{n}[\op{u}_p].
\end{equation}


\paragraph{Part 3: Cochains on \(\B^\sigma (X)\) from noncompact
  \(d\)-submanifolds of \(X\).} In this part we consider
\(d\)-submanifolds in \(X\) that are ``asymptotically cylindrical'' in
the sense described below, and use them to define cochains on
\(\B^\sigma (X)\) (or more generally, on various bundles over
\(\B^\sigma (X)\)) in a manner similar to Part 1. These cochains are
often useful for defining chain homotopy equivalences between Floer
complexes, as will be demonstrated by examples. 

\paragraph{(a)} {\em When \(d=1\).} Let \(M_1, M_2\in \mathpzc{E}\)
label two ends of \(X\), allowing \(M_1=M_2\). We say that an oriented connected 1-submanifold  \(\lambda\subset
X\)  is a {\em path from \(p_1\in M_1\) to \(p_2\in M_2\)} if  \(\lambda\cap
(X-X_c)\) consists of two connected components of the following form:
the first component is \((-\infty,
L)\times\{p_1\}\subset(-\infty, L)\times M_1\) or \((L, \infty)\times
\{-p_1\}\subset (L, \infty)\times M_1\) under the diffeomorphisms in
(\ref{(A.9a,11)}), depending on whether \(M_1\) is a negative end or a
positive end, and the second component is \((-\infty,
L)\times\{-p_2\}\subset(-\infty, L)\times M_2\) or \((L, \infty)\times
\{p_2\}\subset (L, \infty)\times M_2\) under the diffeomorphisms in
(\ref{(A.9a,11)}). We
shall 
define a 1-cochain \([\theta_\lambda]_{\M(X)}\in C^{1;\bbR}_{\M(X)}\)
and its integral correction \([\kappa _\lambda]_{\M(X)}\in
C^{1;\bbZ}_{\M(X)}\), beginning by introducing generalizations of
notions such as \(\op{hol}_{\hat{p}}\), \(\vartheta _{p-p}\), \(\pi_{p-p}\) etc. previously encountered
in Part 2(a). 

Fix choices of \(\vartheta '_{p_1}\in\Omega^1(
\tilde{\B}^\sigma _{p_1}(M_1))\), \(\vartheta '_{p_2}\in\Omega^1(
\tilde{\B}^\sigma _{p_1}(M_2))\) as described in Part 2(a) and  note that \(\tilde{\B}^\sigma _{p_1}(M_1)\times_{U(1)}\tilde{\B}^\sigma
_{p_2}(M_2)\) is a principal \(U(1)\)-bundle over \(\B^\sigma
(M_1)\times \B^\sigma (M_2)\), and \(\vartheta'_{p_1}\),
\(\vartheta'_{p_2}\) together define a principal \(U(1)\)-connection
on this bundle, which we denote by \(\vartheta'_{p_2-p_1}\): Consider
the commutative diagram 
\begin{equation}\label{pr-B}
\xymatrixcolsep{7pc}
\xymatrix{
& \tilde{\mathcal{B}}^\sigma _{p_1}(M_1)\times
\tilde{\mathcal{B}}^\sigma _{p_2}(M_2)\ar@{->}[d]^{\pi_{\Delta
  }}\ar@{->}[ld]_{\tilde{\pr}_i} \ar@/^5pc/[dd]^{\pi_{p_1}\times \pi _{p_2}}\\
\tilde{\mathcal{B}}_{p_i}^\sigma (M_i)
\ar@{->}[d]^{\pi_{p_i}} & \tilde{\mathcal{B}}^\sigma _{p_1}(M_1)\times_{U(1)}
\tilde{\mathcal{B}}^\sigma _{p_2}(M_2) \ar@{->}[l]^{\pr'_i} \ar@{->}[ld]^{\ul{\pr}_i}
\ar@{->}[d]^{\pi_{p_2-p_1}}\\
\B^\sigma (M_i) &\B^\sigma (M_1)\times \B^\sigma (M_2) \ar@{->}[l]^{\pr_i} 
}
\end{equation}
for \(i=1\) or \(2\), where \(\pr_i\) denotes projecting to the \(i\)-th
factor and \(\pi_\Delta\) denotes quotienting by the diagonal
\(U(1)\) action. Then 
\begin{equation}\label{def:theta12}
\vartheta'_{p_2-p_1}=(\pi_\Delta)_!(\tilde{\pr}_1^*\vartheta '_{p_1}\wedge\tilde{\pr}_2^*\vartheta'
_{p_2})=(\pr_2')^*\vartheta '_{p_2}=-(\pr_1')^*\vartheta '_{p_1},
\end{equation}
where \((\pi_\Delta
)_!\co \Omega^2(\tilde{\mathcal{B}}^\sigma _{p_1}(M_1)\times
\tilde{\mathcal{B}}^\sigma _{p_2}(M_2)\to \Omega^1(\tilde{\mathcal{B}}^\sigma _{p_1}(M_1)\times_{U(1)}
\tilde{\mathcal{B}}^\sigma _{p_2}(M_2))\) denotes integrating over the
fibers of \(\pi_\Delta\). Let \(\tilde{\B}^\sigma _\lambda (X)\),
\(\tilde{\tilde{\B}}^\sigma _\lambda (X)\) be pull-back bundles defined by
the following commutative diagram:
\begin{equation}\label{diag:B-lambda}
\xymatrixcolsep{5pc}
\xymatrix{
\tilde{\tilde{\B}}^\sigma _\lambda (X)
\ar@{->}[r]^{\tilde{\tilde{\Pi}}^\partial_\lambda}\ar@{->}[d]^{\tilde{\pi}} \ar@/^-3pc/[dd]_{\tilde{\pi}_\lambda}
& \tilde{\mathcal{B}}^\sigma _{p_1}(M_1)\times
\tilde{\mathcal{B}}^\sigma _{p_2}(M_2)\ar@{->}[d]^{\pi_\Delta}\ar@/^5pc/[dd]^{\pi_{p_1}\times \pi _{p_2}}\\
\tilde{\mathcal{B}}_{\lambda }^\sigma (X) \ar@{->}[r]^{\tilde{\Pi}^\partial_\lambda}
\ar@{->}[d]^{\pi_{\lambda }} & \tilde{\mathcal{B}}^\sigma _{p_1}(M_1)\times_{U(1)}
\tilde{\mathcal{B}}^\sigma _{p_2}(M_2) 
\ar@{->}[d]^{\pi_{p_2-p_1}}\\
\B^\sigma (X) \ar@{->}[r]^{\Pi^\partial_\lambda}& \B^\sigma
(M_1)\times \B^\sigma (M_2)\\
\mathcal{M}(X)\ar@{^(->}[u]\ar@{->}[r]^{\Pi^\partial_\lambda} & \grC
(M_1)\times \grC (M_2)\ar@{^(->}[u],
}
\end{equation}
where \(\Pi^\partial_\lambda :=\Pi^{M_1}\times
\Pi^{M_2}\), and \(\Pi^{M_i}:=\Pi^{\pm\infty}|_{M_i\subset Y_\pm}\).
For \(i=1,2\), let \(\rho _{\vartheta '_{p_i}}\co \ul{U(1)}\to \tilde{\cal B}^\sigma
(M_i)|_{\Pi_0\M_1(\bbR\times M_i)}\) be a trivialization of the
\(U(1)\)-bundle \(\tilde{\cal B}^\sigma
(M_i)\) over \(\Pi_0 \M_1(\bbR\times M_i)\subset {\cal B}^\sigma
(M_i)\). Over \(\grC
(M_1)\times \grC (M_2)\subset  \Pi_0\M_1(\bbR\times M_1)\times
\Pi_0\M_1\subset (\bbR\times M_2)\subset \B^\sigma
(M_1)\times \B^\sigma (M_2)\), the \(U(1)\times U(1)\)-bundle
\(\pi_{p_1}\times \pi_{p_2}\co \tilde{\mathcal{B}}^\sigma _{p_1}(M_1)\times
\tilde{\mathcal{B}}^\sigma _{p_2}(M_2)\to \B^\sigma
(M_1)\times \B^\sigma (M_2)\) is equipped with a trivialization
\(\rho_{\vartheta '_{p_1}}\times \rho_{\vartheta '_{p_1}}\). This
trivialization factors through a trivialization, \(\rho_{p_2-p_1}\),
of the \(U(1)\)-bundle, \(\pi_{p_2-p_1}\co \tilde{\mathcal{B}}^\sigma _{p_1}(M_1)\times_{U(1)}
\tilde{\mathcal{B}}^\sigma _{p_2}(M_2) \to \B^\sigma
(M_1)\times \B^\sigma (M_2)\) 
over  \(\grC(M_1)\times \grC (M_2)\subset \B^\sigma
(M_1)\times \B^\sigma (M_2)\) and a trivializtion, \(\rho_\Delta \), of the
\(U(1)\)-bundle \(\pi_\Delta \co \tilde{\mathcal{B}}^\sigma _{p_1}(M_1)\times
\tilde{\mathcal{B}}^\sigma _{p_2}(M_2)\to \tilde{\mathcal{B}}^\sigma _{p_1}(M_1)\times_{U(1)}
\tilde{\mathcal{B}}^\sigma _{p_2}(M_2)\) over \[
\tilde{\grC}:=\pi_{p_2-p_1}^{-1}(\grC 
(M_1)\times \grC (M_2))\subset \tilde{\mathcal{B}}^\sigma _{p_1}(M_1)\times_{U(1)}
\tilde{\mathcal{B}}^\sigma _{p_2}(M_2).
\] 
The trivializations \(\rho_{p_2-p_1}\) and \(\rho _\Delta \) above are
compatible respectively with \(\vartheta '_{p_2-p_1}\) and
\(\tilde{\op{pr}}_1^*\vartheta '_{p_1}+\tilde{\op{pr}}_2^*\vartheta'
_{p_2}\), which are in turn integral respectively over
\(\grC(M_1)\times \grC (M_2)\) and \(\tilde{\grC}\) as \(\vartheta'
_{p_i}\) satisfy (\ref{int-con}). 
All the trivializations above are
determined by \(\vartheta'_{p_1}\) and \(\vartheta'_{p_2}\) modulo
constant \(U(1)\)-maps. 

Identify
the hermitian line
bundle associated to the principal \(U(1)\)-bundle \(\pi_\Delta \co \tilde{\mathcal{B}}^\sigma _{p_1}(M_1)\times
\tilde{\mathcal{B}}^\sigma _{p_2}(M_2) \to \tilde{\mathcal{B}}^\sigma _{p_1}(M_1)\times_{U(1)}
\tilde{\mathcal{B}}^\sigma _{p_2}(M_2)\) with the bundle \(\op{Hom}\, \big(\ul{\pr}^*_1\mathcal{E}_{p_1}(M_1),
\ul{\pr}^*_2\mathcal{E}_{p_2}(M_2)\big)\), and use \(E_{p_i}\) to denote the
fiber of the bundle \(E\to M_i\) at \(p_i\in M_i\). Given
\(\ul{\grd}\in \mathcal{C}^\sigma (X)\), let 
\(\op{hol}_\lambda^E(\ul{\grd})\in \Hom (E_{p_1}, E_{p_2}) \) denote
the holonomy along \(\lambda \) of the \(A\in \op{Conn}(E)\)
associated to \(\ul{\grd}\). Observing that given a
\(\tilde{\grd}\in\tilde{\cal B}^\sigma _\lambda (X)\), the value 
\(\op{hol}_\lambda^E(\ul{\grd})\) is identical for 
all representatives
\(\ul{\grd}\in \mathcal{C}^\sigma (X)\) of \(\tilde{\grd}\), this then
defines a map, also denoted \(\op{hol}_\lambda\), from \(\tilde{\cal
  B}^\sigma _\lambda(X)\) to \(\tilde{\mathcal{B}}^\sigma _{p_1}(M_1)\times
\tilde{\mathcal{B}}^\sigma _{p_2}(M_2)\) that fits into the
following commutative diagram: 
\begin{equation}
\label{CD-ttB}\xymatrixcolsep{5pc}
\xymatrix{
\tilde{\tilde{\B}}^\sigma _\lambda (X) \ar@{->}[r]^{\op{hol}_\lambda}
\ar@<.5ex>[d]^{\tilde{\pi}}& \tilde{\mathcal{B}}^\sigma _{p_1}(M_1)\times
\tilde{\mathcal{B}}^\sigma _{p_2}(M_2)\ar@{->}[d]^{\pi_\Delta}\\
\tilde{\mathcal{B}}_{\lambda }^\sigma (X)
\ar@{->}[r]^{\tilde{\Pi}^\partial_\lambda}\ar@{->}[ur]^{\op{hol}_\lambda}
& \tilde{\mathcal{B}}^\sigma _{p_1}(M_1)\times_{U(1)}
\tilde{\mathcal{B}}^\sigma _{p_2}(M_2).
}
\end{equation}
Since \(\Pi^\partial_\lambda (\M (X))\subset \grC(M_1)\times \grC
(M_2)\subset \B^\sigma(M_1)\times \B^\sigma (M_2)\) and \(\tilde{\Pi}^\partial_\lambda (\tilde{\M}(X))\subset
\tilde{\grC}\subset \tilde{\mathcal{B}}^\sigma _{p_1}(M_1)\times_{U(1)}
\tilde{\mathcal{B}}^\sigma _{p_2}(M_2)\), we may compose 
\(\op{hol}_\lambda|_{\tilde{\M}(X)}\) with the trivialization
\(\rho_\Delta \) to get a map \(\op{h}_\lambda\co
\tilde{\M}(X)\subset\tilde{\mathcal{B}}_{\lambda }^\sigma (X)\to
U(1)=\bbR/\bbZ\). Let 
\[
\vartheta_\lambda:=d\op{h}_\lambda,
\] 
a closed 1-form on \(\tilde{\cal M}(X)\). Note that
\(\vartheta_\lambda \) depends only on \(\vartheta '_{p_1}\), \(\vartheta'
_{p_2}\), not the choices of \(\rho_{\vartheta '_{p_1}}\),
\(\rho_{\vartheta '_{p_2}}\).  Let \(\tilde{\Pi}^{M_i}_\lambda :=\op{pr}'_i\circ
\tilde{\Pi}^\partial_\lambda \) and observe that both  \(\vartheta _\lambda
\) and 
\begin{equation}\label{Pi_lambda}
(\tilde{\Pi}^\partial_\lambda)^*\vartheta'_{p_2-p_1}=(\tilde{\Pi}^{M_2}_\lambda )^*\vartheta'
_{p_2}=-(\tilde{\Pi}^{M_1}_\lambda )^*\vartheta '_{p_1}
\end{equation}
define principal \(U(1)\)-connections on the bundle \(\pi_\lambda \co
\tilde{\M}(X)\to \M (X)\). Thus, 
\begin{equation}\label{vartheta-theta}
\vartheta _\lambda -(\tilde{\Pi}^\partial_\lambda )^*\vartheta'_{p_2-p_1}=\pi_\lambda^*\theta_\lambda
\end{equation}
for a 1-form \(\theta_\lambda\) on \(\M (X)\), and correspondingly, a \([\theta _\lambda ]_{\M(X)}\in C^{1;\bbR}_{\M(X)}\). Note that
\(\theta_\lambda \) does not depend on the choice of either
\(\vartheta '_{p_1}\) or \(\vartheta '_{p_2}\), since varying the choice
of either changes \(\vartheta _\lambda
\) and \((\tilde{\Pi}^\partial_\lambda )^*\vartheta'_{p_2-p_1}\) by the
same amount. 

As observed in Remark \ref{rem:generalized-m}, with \(\theta _\lambda
\) constructed from forms on the bundle \(\tilde{\B}^\sigma _\lambda
(X)\), the cobordism map \(\mathring{m}[\theta _\lambda ]\) is
defined by a generalization of the formula in \cite{KM}. Let
\(m^\#_\flat[\theta _\lambda ]\), \(\bar{m}^\#_\flat[\theta _\lambda
]\) be defined as a sum of integrals in the usual way, i.e. by
(\ref{eq:m-coeff-loc}), the explicit formula for \(\hat{m}[\theta
_\lambda ](X)\), generalizing (\ref{def:m-map}), is given below: 
\begin{equation}\label{def:m-lambda}
\left[\begin{array}{cc}
m^o_o[\theta_\lambda] &m^u_o[\theta_\lambda]\\
\hat{m}^o_u[\theta _\lambda ]
&\hat{m}^u_u[\theta _\lambda ]
\end{array}\right],
\end{equation}
where 
\begin{equation}\label{def:m-lambda0}
\begin{split}
\hat{m}^o_u[\theta _\lambda
]&
:=-\bar{m}^s_u[\theta_\lambda]\partial^o_s-\bar{\partial}^s_um^o_s[\theta_\lambda]+\bar{n}^s_u[d\ul{\op{h}}_{\hat{p}_2}]\, 
m^o_s[1],\\ 
\hat{m}^u_u[\theta _\lambda ] & :=-\bar{m}^u_u[\theta_\lambda]-\bar{m}^s_u[\theta _\lambda ] \partial
^u_s-\bar{\partial}^s_um^u_s
[\theta_\lambda]+\bar{n}^s_u[d\ul{\op{h}}_{\hat{p}_2}]\, m^u_s[1]
\end{split}
\end{equation}
when \(p_1\in Y_-\) and
\(p_2\in Y_+\); when \(p_1, p_2\) are both in \(Y_+\), then the 
\(\bar{n}^s_u[d\ul{\op{h}}_{\hat{p}_2}]\)'s in the formulas 
 above are replaced by
\(\bar{n}^s_u[d\ul{\op{h}}_{\hat{p}_2}]-\bar{n}^s_u[d\ul{\op{h}}_{\hat{p}_1}]\). When \(p_1, p_2\) both
belongs to \(Y_-\), \(\hat{m}[\theta _\lambda]\) is given by
(\ref{def:m-map}). (In this case it is \(\check{m}[\theta _\lambda ]\)
that gains additional terms.)

\paragraph{\it Example.} Let \(M\) be connected, and take
\(X=\bbR\times M\) to be a product cobordism. Let 
\(\lambda\subset \bbR\times M\) be the graph of a path
\(\bar{\lambda}(\cdot)\co \bbR\to M\) that sends the \((-\infty, -L')\subset
\bbR\) to \(p_1\in M\) and \((L',\infty)\subset \bbR\) to \(p_2\in
M\). Let \(\bar{\M}\) be as in Part 2(a). 
We choose a lifting \(\tilde{\varsigma }_\lambda\) of the embedding \(\varsigma \co\bar{\M}\to \B^\sigma (X)\) in a way parallel to (\ref{CD-M+}), namely,
such that the following diagram commutes: 
\begin{equation}
\label{CD-M+1}
\begin{CD}
\mathcal{M}^+@>\tilde{\varsigma}_\lambda
>>\tilde{\mathcal{B}}^\sigma_\lambda (X)@>\tilde{\Pi}^\partial_\lambda
>>\tilde{\B}^\sigma _{p_1}(M)\times _{U(1)}\tilde{\B}^\sigma _{p_2}(M)\\
@ V\grr VV @ V\pi_\lambda  VV @ V\pi_{p_2-p_1}VV\\
\bar{\cal M}  @>\varsigma>> {\cal B}^\sigma  (X)  @>\Pi^\partial_\lambda
>>\B^\sigma(M)\times\B^\sigma (M). 
\end{CD}
\end{equation}
As observed previously, over \((\varsigma\circ
\Pi^\partial_\lambda)\bar{\M}\subset  \grC(M_1)\times
\grC(M_2)\subset \B^\sigma(M)\times\B^\sigma (M)\), the bundle
$\pi_{p_2-p_1}\co \tilde{\B}^\sigma _{p_1}(M)\times_{U(1)}\tilde{\B}^\sigma _{p_2}(M)$ is trivialized by
$\rho_{p_2-p_1}$. This induces a trivialization of its  pull-back bundle
$\pi_\lambda\co \tilde{\B}^\sigma_\lambda (X)\to \B^\sigma(X)$ 
(via \((\tilde{\Pi}^\partial_\lambda )^*\)) over
\(\bar{\M}\stackrel{\varsigma }{\hookrightarrow}\B^\sigma(X)\). Choose
\(\tilde{\varsigma }_\lambda \) to be constant with respect to this
trivialization. Then $(\tilde{\Pi}^\partial_\lambda)^*\vartheta'
_{p_2-p_1}$ vanishes over
\(\tilde{\varsigma}_\lambda((\M^+)_1\backslash(\M^+)_0)=\pi_\lambda
^{-1}(\bar{\M}_1\backslash\bar{\M}_0)\). As 
$\vartheta _\lambda \in \Omega^1(\B^\sigma (X))$ is closed by
construction, arguing as in (\ref{value}) and the subsequent
discussions, again using \cite{KM}'s Theorem 24.7.2, Lemma
21.3.1 and (\ref{M-strata}), we have: 
\begin{equation}
\label{value1}
\begin{split}
0 = & \langle\tilde{\varsigma }^*_\lambda (d\vartheta _\lambda ), \mathcal{M}^+\rangle \\
& =\langle\tilde{\varsigma }^*_\lambda
\vartheta_\lambda , \partial [\mathcal{M}^+]\rangle=\langle\tilde{\varsigma }^*_\lambda \vartheta_\lambda , [(\mathcal{M}^+)_1]\rangle\\
&=\langle\tilde{\varsigma }^*_\lambda \vartheta_\lambda , \grr^{-1}(\bar{\M})_0\rangle+\langle\tilde{\varsigma }^*_\lambda \vartheta_\lambda , \grr^{-1} (\bar{\M}_1\backslash
\bar{\M}_0)\rangle\\
& =\langle \tilde{\varsigma }^*_\lambda (\tilde{\Pi}^\partial_\lambda )^*\vartheta'
_{p_2-p_1}, \{\bar{i}\grc_-, \bar{i}\grc_+\}\times\N^+\rangle+\langle\theta_\lambda , \bar{\M}_1\backslash
\bar{\M}_0\rangle\\
& =\langle \tilde{\varsigma }^*_\lambda (\tilde{\Pi}^{M_1}_\lambda
)^*\vartheta'_{p_1},\{\bar{i}\grc_-\}\times \N^+\rangle-\langle \tilde{\varsigma }^*_\lambda (\td{\Pi}_\lambda ^{M_2})^*\vartheta'_{p_2}, \{\bar{i}\grc_+\}\times \N^+\rangle\\
& \qquad\quad  +\langle[\theta_\lambda ]_{\M
  (X)}, \partial[\bar{\M}]\rangle\\
&=\langle \mathpzc{e}_{\dot{p}_1}, \bar{\M}\rangle-\langle \mathpzc{e}_{\dot{p}_2}, \bar{\M}\rangle+\Big\langle[\theta_\lambda ]_{\M
  (\bbR\times M)} -[\theta_{\dot{p}_1}]_{\M
  (\bbR\times M)}+[\theta_{\dot{p}_2}]_{\M
  (\bbR\times M)}, \partial[\bar{\M}]\Big\rangle, 
\end{split}
\end{equation}
(To see the last two lines in the preceding expression,
recall (\ref{e-hol}) and (\ref{Pi_lambda}).) Summarizing, we have 
\begin{equation}\label{e-diff}\begin{split}
\mathpzc{e}_{\dot{p}_2}-\mathpzc{e}_{\dot{p}_1} & =\delta\uu_{\lambda }\in C^{2;\bbZ}_{\M
  (\bbR\times M)}, \quad \text{where}\\
\uu_{\lambda }& :=[\theta_\lambda ]_{\M
  (\bbR\times M)} -[\theta_{\dot{p}_1}]_{\M
  (\bbR\times M)}+[\theta_{\dot{p}_2}]_{\M(\bbR\times M)}\in C^{1;\bbR}_{\M
  (\bbR\times M)}. 
\end{split}
\end{equation}
Note that \(\uu _{\lambda }\) in fact has integral coefficients,
i.e. \(\uu _\lambda \in C^{1;\bbZ}_{\M(\bbR\times M)}\subset C^{1;\bbR}_{\M(\bbR\times M)}\). To see this, recall (\ref{theta-h}) and write 
\[
\uu _{\lambda }=\sum_{\ul{\grd}\in
  \N_0(M)}\Big(\int_{\bbR\ul{\grd}}\theta_\lambda
+\ul{h}_{\hat{p}_1}(\ul{\grd})-\ul{h}_{\hat{p}_2}(\ul{\grd})\Big) \mu^1_{\ul{\grd}}. 
\]
Meanwhile, for a \(\ul{\grd}\) in \(\M (\grc_-, \grc_+)/\bbR\), 
\[\begin{split}
& \int_{\bbR\ul{\grd}}\theta_\lambda \quad\mod \bbZ\\
&\quad 
=\int_{\bbR\ul{\grd}}\tilde{\varsigma}_\lambda^*(d\op{h}_\lambda)\quad\mod
\bbZ\\
&\quad =\op{h}_\lambda(\td{\varsigma}
  _\lambda(\bar{i}(\grc_+)))-\op{h}_\lambda(\td{\varsigma }_\lambda(\bar{i}(\grc_-)))
  \in \bbR/\bbZ
\end{split}
\]
Let \(\grd\in \C^\tau (\bbR\times M)\) represent  an element in
\(\td{\pi}_\lambda ^{-1}q_\bbR^{-1}(\ul{\grd})\) and let \(\gamma
\subset \bar{\bbR}\times M\) be the loop formed by the union of four
arcs \((\bar{\bbR}\times \{p_1\})\cup(\bar{\bbR}\times \{-p_2\})\cup (\{\infty\}\times \ul{\lambda
})\cup(\{-\infty\}\times (-\ul{\lambda}))\) in \(\bar{\bbR}\times M\),
where  \(\ul{\lambda }\subset M\) is the closure of the image of the path  \(\lambda(\cdot)\co \bbR\to M\). Then 
\[
\op{h}_\lambda(\td{\varsigma}
  _\lambda(\bar{i}(\grc_+)))-\op{h}_\lambda(\td{\varsigma
  }_\lambda(\bar{i}(\grc_-)))+\ul{\op{h}}_{\hat{p}_1}(\ul{\grd})-\ul{\op{h}}_{\hat{p}_2}(\ul{\grd})
=-\frac{i}{2\pi}\ln \big(\op{hol}_\gamma ^E(\grd))\big)=0\in \bbR/\bbZ;
\] 
and hence the coefficients in \(\uu_\lambda \), \(\int_{\bbR\ul{\grd}}\theta_\lambda
-\ul{h}_{\hat{p}_1}(\ul{\grd})+\ul{h}_{\hat{p}_2}(\ul{\grd})\in \bbZ\).

 
Let \(\mathring{\op{K}}_\lambda=\mathring{m}[\uu_\lambda]\co
\mathring{C}(M)\to \mathring{C}(M)\), a degree \(-1\) map defined in
the same manner as \(\mathring{m}[\theta_\lambda]\), namely as was in
(\ref{def:m-lambda}) and (\ref{def:m-lambda0}) with
\(\bar{n}^s_u[d\ul{\op{h}}_{\hat{p}_i}]\) there replaced by
\(\bar{n}^s_u[\op{u}_{p_i}]=(\bar{U}_{p_i})^s_u\). It follows
from (\ref{e-diff}) and \cite{KM}'s Proposition 25.3.4 that 
\begin{equation}\label{def:opK}
\mathring{U}_{p_2}-\mathring{U}_{p_1}=[
\mathring{\op{K}}_\lambda, \mathring{\partial}]. 
\end{equation}
Namely, \(\mathring{\op{K}}_\lambda\) defines a chain homotopy equivalence
between the two \(U\)-maps \(\mathring{U}_{p_2}\) and
\(\mathring{U}_{p_1}\). 


The arguments in the preceding example generalizes readily to
cobordisms \(X\) of the types considered in Section
\ref{sec:2.4}. 
Note that the diagram (\ref{CD-M+1}) and the first three
lines of (\ref{value1}) hold in general. When \(X\) is not a product
cobordism, the fourth line of (\ref{value1}) has a simple modification
by replacing its first term by the more general \(\langle \tilde{\varsigma }^*_\lambda (\tilde{\Pi}^\partial_\lambda )^*\vartheta'
_{p_2-p_1}, \grr^{-1}(\bar{\mathcal{M}})_0\rangle\), where
\(\grr^{-1}(\bar{\mathcal{M}})_0\) fibers over
\((\bar{\mathcal{M}})_0\), with fibers consisting of 1-dimensional strata of
\(\mathcal{N}^+(M_1)\) or \(\mathcal{N}^+(M_2)\). 

The map \(\hat{\op{K}}_\lambda \) in the Example has an analog in this
setting, which we denote by the same notation:
\begin{equation}\label{def:K_lambda}
\begin{split}
\mathring{\op{K}}_\lambda  (X)& :=\mathring{m}[\theta _\lambda](X)+
 \mathring{\Theta}_{p_2} * \mathring{m}[1](X)-\mathring{\Theta }_{p_1}
 * \mathring{m}[1](X). 
\end{split}
\end{equation}
In the above, 
\[
\mathring{\Theta}_{p_i}:=\mathring{m}[\theta
_{\dot{p}_i}](\bbR\times M_i)
=-\mathring{n}[\ul{h}_{\hat{p}}](M_i). 
\]
 is used to denote both an endomorphism
on 
\(\mathring{C}(M_i)\)  and its associated endomorphism, 
\(\mathring{\Theta}_{p_i}\otimes 1\), on \(\mathring{C}(Y_\pm)\). Meanwhile, \(\mathring{\Theta}_{p_i} * \mathring{m}[1](X)\)
denotes either the composition \(\mathring{\Theta}_{p_i} \, \mathring{m}[1](X)\) or
\(\mathring{m}[1](X) \, \mathring{\Theta}_{p_i} \), depending on whether
\(p_i\in Y_-\) or \(p_i\in Y_+\). Note that 
while \(\mathring{m}[\theta _\lambda
](X)\) is defined for coefficient ring
\(\bbK=\bbR\), arguments similar to those in the preceding example
show that \(\mathring{\op{K}}_\lambda  \) is in fact defined for
coefficient ring \(\bbK=\bbZ\).  
The arguments also give rise to an
analog of the identity (\ref{def:opK}):
\begin{equation}\label{eq:K-lambda}
\mathring{U}_{p_2}*
\mathring{m}[1](X)-\mathring{U}_{p_1}* \mathring{m}[1](X)=[
\mathring{\op{K}}_\lambda, \mathring{\partial}].
\end{equation}
 
\begin{remarks}
Instead of the formula given in 
(\ref{def:K_lambda}), it is possible to express \(\mathring{\op{K}}_\lambda \) as 
\[
\mathring{\op{K}}_\lambda =\mathring{m}[\uu_\lambda ](X),
\] 
with \(\uu_\lambda \in C^{1;\bbZ}_{\M(X)}\), in a way
parallel to (\ref{e-diff}). This often yields cleaner formulae in
later discussions but is less practical, being not as concrete as
(\ref{def:K_lambda}). In what follows we alternate between these two
equivalent description of \(\mathring{\op{K}}_\lambda \), depending on
which is more convenient in the context. 
\end{remarks}

\paragraph{(b)} {\em When \(d=2\).} For each \(M_i\in
\mathpzc{E}\), let \(\gamma _i\subset M_i\) be an embedded
(oriented) 
circle or the empty set. 
Let \(\Sigma\subset X\) be an
embedded oriented surface {\em asymptotic to \(\{\gamma_i\}_{i\in
    \mathpzc{E}}\)} in the following sense: 
 \(\Sigma \cap
(X-X_c)\) is the union of connected components of the following form:
under the diffeomorphisms in
(\ref{(A.9a,11)}), for each \(M_i\) there is a component 
\((-\infty, L')\times(-\gamma_i)\subset(-\infty, L')\times M_i\) if
\(M_i\) is a negative end, and it is 
\((L', \infty)\times\gamma_i\subset (L', \infty)\times M_i\) if
\(M_i\) is a positive end. Let \(F_\Sigma\co
\mathcal{B}^\sigma (X)\to \bbR\) be the function sending a \(\grd\in\mathcal{B}^\sigma
(X)\) to 
\[F_\Sigma(\grd): =\int_{\Sigma}\frac{iF_A}{2\pi}, \] where \(A\in
\op{Conn}(E)\) is the connection associated to an arbitrary
representative of \(\grd\). The function \(F_\Sigma \) depends only on the
relative homology class of \(\Sigma\): for another embedded surface
\(\Sigma'\) asymptotic to the end \(\{\gamma_i\}_{i\in \mathpzc{E}}\), 
\[F_{\Sigma'}-F_\Sigma =\langle c_1(\grs)-c_1(K^{-1}),
[\Sigma'-\Sigma]\rangle/2.\] 
Let \(\theta '_{\gamma_i}\in \Omega ^1(\B^\sigma (M_i))\),
\(\op{h}'_{\gamma_i}\co \B^\sigma (M_i)\to \bbR/\bbZ\) be as defined
in Part 2(b) if \(\gamma _i\neq \emptyset\), and let
\(\theta'_{\gamma_i}:=0\) and \(\op{h}'_{\gamma_i}=0\) if \(\gamma
_i=\emptyset\). Then 
\begin{equation}\label{dF_Sig}
d F_\Sigma=\sum_{i\in \mathpzc{E}}(\Pi^{M_i})^*\theta'_{\gamma_i}. 
\end{equation}
Thus, an integral correction of \([d F_\Sigma]_{\M
  (X)}=\delta[F_\Sigma]_{\M(X)}\in C^{1;\bbR}_{\M (X)}\) takes the
form of \(\delta \textsc{f}_\Sigma\), where
\(\textsc{f}_\Sigma\in C^{0;\bbZ}_{\M(X)}\subset
C^{0;\bbR}_{\M(X)}\) is given by 
\begin{equation}\label{scF_S}
\textsc{f}_\Sigma:=[F_\Sigma ]_{\M(X)}-\sum_{i\in
  \mathpzc{E}}\big[(\Pi^{M_i})^*h'_{\gamma_i}\big]_{\M(X)}, 
\end{equation}
where \(h'_{\gamma_i}\co \grC(M_i)\to \bbR\) is as in Part
2(b). 

 \paragraph{\it Example.} 
Take \(X=\bbR\times M\) again to be the
product cobordism,  
and let \(\Sigma\) be
such that \(s\mapsto \Sigma\cap (\{s\}\times M)\) forms a homotopy between the circles
\(\gamma_-, \gamma_+\subset M\), both representing the element
\(\grt_i\in H_1(M;\bbZ)/\op{tors}\). 
Applying Equations (\ref{scF_S}) and (\ref{dF_Sig}) to this setting,
and recalling from Part 2(b) the definition and properties of \(\uu_{\dot{\gamma }}\),
we have 
\begin{equation}\label{dF}
\delta \textsc{f}_\Sigma=\uu_{\dot{\gamma }_+}-\uu_{\dot{\gamma}_-}. 
\end{equation}
By Proposition 25.3.4 of \cite{KM}, this implies that 
\begin{equation}\label{eq:dF}
\mathring{\grm}_{\gamma_-}-\mathring{\grm}_{\gamma_+}=\big[\mathring{\partial},
\mathring{m}[\textsc{f}_\Sigma](\bbR\times M)\big].
\end{equation}
Namely, \(\mathring{m}[{\textsc{f}}_\Sigma]\) defines a chain homotopy
equivalence between the two \(\grt_i\) maps  \(\mathring{\grm}_{\gamma_-}:=\mathring{m}[\uu_{\dot{\gamma }_-}]\) and \(\mathring{\grm}_{\gamma_+}:=\mathring{m}[\uu_{\dot{\gamma }_+}]\).

\bigbreak 

The preceding example also generalizes readily. When \(X\) is not a
product cobordism, the identities (\ref{dF}) and (\ref{eq:dF}) have
respectively the following analog: 
\begin{equation}\label{eq:dF-X}
\begin{split}
\delta \textsc{f}_\Sigma  = & \sum _{i\in \mathcal{E}}[(\Pi^{M_i})^*\mu_{\gamma _i}];\\
-\sum_i\mathring{m}[1](X) *\mathring{\grm}_{\gamma_i}& =\Big[\mathring{\partial},
\mathring{m}[\textsc{f}_\Sigma](X)\Big], 
\end{split}
\end{equation}
where \(\mathring{m}[1](X) *\mathring{\grm}_{\gamma_i}\) denotes the
composition map \(\mathring{m}[1](X) \mathring{\grm}_{\gamma_i}\) when
\(\gamma_i\subset Y_+\), and it denotes
\(-\mathring{\grm}_{-\gamma_i}\mathring{m}[1](X) \) when \(\gamma
_i\subset Y_-\). 

\begin{remarks}\label{rmk:S_1-prod}
In view of (\ref{def:m-disconn}), the actions \(\mathring{U}_p\) and
\(\mathring{\grm}_{\grt_i}\) defined in Part 2 above extend to the case when \(M\) is not necessarily connected, and together they define a
\(\mathbf{A}_\dag(M):=\bbK[U]\otimes
\bigwedge^*H_1(M;\bbZ)/\text{Tors}\) action associated to each choice of \(p\) and
\(\{\grt_i\}_i\) for possibly disconnected \(M\). 
These more general \(\mathring{U}_p\) and
\(\mathring{\grm}_{\grt_i}\) are chain maps as in the connected case;
in fact it follows as a straightforward consequence of the case for connected 3-manifolds, already 
verified in \cite{KM} in the process of defining the
\(\mathbf{A}_\dag \) actions on the monopole Floer homology \(\widehat{\HM}\). 
The arguments in Part 3(a)
show that in this more general setting, \(\mathring{U}_{p_1}\) and \(\mathring{U}_{p_2}\)
are chain homotopy equivalent when \(p_1\) and \(p_2\) belong to the
same connected component of \(M\), but {\em not} if \(p_1\), \(p_2\)
lie on different components of \(M\). In fact, the cohomology classes
\([\mu_{p_1}], [\mu_{p_2}]\in H^*(\B^\sigma (M);\bbZ)=H^*(\B^\sigma
(M_1); \bbZ) \otimes H^*(\B^\sigma
(M_2); \bbZ)\) are independent. Generalizing the definition of \(\mathring{U}_p\), given
a 0-cycle \({\bf p}\) consisting of finitely many signed points
\(p_i\) in \(M\), let 
\[
\mathring{U}_{\bf p}:=\sum_i\op{sign}(p_i) \, \mathring{U}_{p_i}. 
\] 
Suppose \(M=M_\sqcup:=M_1\sqcup M_2\) consists of two connected
components \(M_1\) and \(M_2\), and so \(\B^\sigma
(M_\sqcup)=\B^\sigma (M_1)\times \B^\sigma (M_2)\). Suppose \(p_i\in
M_i\) for \(i=1,2\). Then the spaces \( \tilde{\mathcal{B}}^\sigma _{p_1}(M_1)\times
\tilde{\mathcal{B}}^\sigma _{p_2}(M_2)\), \( \tilde{\mathcal{B}}^\sigma _{p_1}(M_1)\times_{U(1)}
\tilde{\mathcal{B}}^\sigma _{p_2}(M_2)\) in the second column of the
diagram (\ref{pr-B}) are respectively \(U(1)\times U(1)\)- and
\(U(1)\)-bundles over \(\B^\sigma
(M_\sqcup)\), and we abbreviate them respectively as 
\(\td{\td{\B}}^\sigma _{p_1, p_2}(M_\sqcup)\) and 
\(\tilde{\mathcal{B}}^\sigma _{p_2-p_1}(M_\sqcup) \). 
It is worth noting that \(\tilde{\mathcal{B}}^\sigma
_{p_2-p_1}(M_\sqcup) \) is of the same homotopy
type as
\(\mathcal{B}^\sigma (M_1\# M_2)\), \(M_1\# M_2\) being the connected
sum of \(M_1\) of \(M_2\) along \(p_1\) and \(p_2\). While the Floer complex \(\hat{C}(M_\sqcup)\) in
Part 2 above  (heuristically) reflects the topology of 
\(\mathcal{B}^\sigma(M_\sqcup)\), the connected sum theorem in Section
\ref{sec:6} relates the Floer complex \(\hat{C}(M_\#)\) (associated
with \(\mathcal{B}^\sigma(M_\#)\)) not directly to
\(\hat{C}(M_\sqcup)\), but to a ``Floer complex associated with \(\tilde{\mathcal{B}}^\sigma
_{p_2-p_1}(M_\sqcup)\)''. Using the description of  \(\tilde{\mathcal{B}}^\sigma
_{p_2-p_1}(M_\sqcup) \) as an \(S^1\)-bundle over
\(\mathcal{B}^\sigma(M_\sqcup)\), the latter complex is constructed 
using what was called the ``algebraic \(S^1\)-bundle'' operation in
\cite{L}, described in more detail in Section \ref{sec:4} below. 
The ingredients of this construction consist of a
chain-complex for the orbit space of the \(S^1\)-action,
endowed with a ``\(U\)-map'' associated to its Euler class. 
The Euler class
of the bundle \(\tilde{\mathcal{B}}^\sigma
_{p_2-p_1}(M_\sqcup) \) is \(\pr^*_2\e_{\dot{p}_2}-\pr^*_1\e_{\dot{p}_1}\); so in
the setting under discussion, these are 
\(\hat{C}(M_\sqcup)\), endowed with \(U\)-map 
\begin{equation}\label{eq:U-cup}
\hat{U}_{\sqcup}:=1\otimes\hat{U}_{p_2}-\hat{U}_{p_1}\otimes 1=\hat{U}_{p_2-p_1}. 
\end{equation}
The precise definition of (the hat-flavor of ) ``the Floer
complex for \(\tilde{\mathcal{B}}^\sigma
_{p_2-p_1}(M_\sqcup) \)'' is then what is called
\(S_{\hat{U}_\sqcup}\hat{C}_*(M_\sqcup)\) in Part 3 of
Section \ref{sec:6.1}. There, for any given \(p\in M_\sqcup=M\) 
we also introduce an associated \(U\)-map on this Floer
complex. Two such \(U\)-maps
associated to different points \(p, p'\in M\) are chain
homotopy-equivalent  even if \(p\), \(p'\) belong to different
 connected components of \(M\).   (Cf. Lemma \ref{lem:U-diff} below.) 
\end{remarks}

\paragraph{Part 4. \(\mathbf{A}_\dag\)-actions under large \(r\)
  perturbations.}  Let \(M\) be connected and let \(Q\) denote one of the generating elements
of \(U\) or \(\grt_i\) of \(\mathbf{A}_\dag(M)\), \(U\), \(\grt_i\)
being as defined in the beginning of Part 2. 
In the non-balanced setting discussed in  \cite{KLT4}, a particular
choice of \(p\) and \(\gamma_{i}\)'s was made for the case when \(M\)
is the auxiliary manifold \(Y\) in Theorem \ref{thm:main} (cf. Part 7 of Section
IV.1.3), and the associated \(U\)-maps and
\(\grt_i\)-maps were defined concretely. In this part we relate the
description therein with the more general and abstract construction
given in Part 2 above. The same arguments can be used to re-interprete
the type of cobordism maps  in
Parts 1 and 3 under large \(r\) perturbations in a manner similar to
\cite{KLT4}. Details will be provided for some particular 3-manifolds and
cobordisms (including \(Y\) and the product cobordism \(\bbR\times Y\)
from \cite{KLT4} as special cases) in Section \ref{sec:3} below. 

 In the context of 
\cite{KLT4} as well those to be discussed in Section 3, the spinor bundle \(\bbS\) on \(M\) splits as \(E\oplus E\otimes
K^{-1}\), and hence also a splitting of \(\bbS^+\) on \(\bbR\times
M\), which we denote by the same notation.  As pointed out in
Part 2(a), in this case the tautological section \(\td{\alpha }\) on \(\pmb{\cal
  E}(\bbR\times M)\) or \(\pmb{\cal E}(\bar{\bbR}\times M)\) is
well-defined. 

The non-balanced assumption implies that there are no reducible
Seiberg-Witten solutions, leading to significant simplications. To
name a few:  this allows one to replace all the blowup space
\({\cal B}^\sigma\) occurring in last part by the space \(\mathcal{B}\)
before blowing-up. 
It also implies that \(\bar{U}\)- and
\(\bar{\grt}_i\)-maps are trivial, and \(\hat{U}=\check{U}=:U\),
\(\hat{\grm}_{\grt_i}=\check{\grm}_{\grt_i}:=\grm_{\grt_i}\). Morever,
the relevant moduli spaces are manifolds with
corners in this setting; namely (\ref{m-w-corner}) holds and \([\partial\M]=\partial[\M]\). 

The generating set of the relevant Floer complex, \(\grC(M)\), in \cite{KLT4} is
denoted by \(\mathcal{Z}=\mathcal{Z}_{SW,r}\). For large \(r\), this
is a finite set, and its elements are all represented by elements of
the form \((A, (\alpha , \beta ))\in
\op{Conn}(E)\times \Gamma(E\oplus  E\otimes K^{-1})\) with
\(\alpha^{-1}(0)\) consisting of finitely many points in \( M\). This
makes it possible to choose the 
point \(p\in M\) used to define the \(U\)-map and the embedded circles
\(\gamma _i\) used to define the \(\grt_i\)-maps  to be
mutually disjoint and to all lie in the complement of \(\alpha
^{-1}(0)\subset M\). Write the map, \(\grm_Q\), associated to each \(Q\) in a form
similar to (\ref{eq:differential}) and (\ref{eq:differential-loc}):  (\(\grm_U:=U\))
\[
\grm_Q=\sum_{\grc_1, \grc_2\in \mathcal{Z}}\sum_{z\in\pi_1
  \mathcal{B}(M;\grc_1, \grc_2)}\op{w}_Q(\grc_1,
\grc_2; z)\,\Gamma (z);
\]
and in the monotone case, let \(\op{w}_Q(\grc_1,
\grc_2)=\sum_z\op{w}_Q(\grc_1,\grc_2; z)\). The discussion in the rest
of this Part works for both \(\op{w}_Q(\grc_1,\grc_2; z)\) and
\(\op{w}_Q(\grc_1,\grc_2)\), but for simplicity only the latter will be mentioned. 

\paragraph{(a)} {\em The \(U\)-map associated to \(p\in M\).} In the formulation of Part 2(a),
the coefficients of the \(U_p\)-map are  given by 
\[
\op{w}_U(\grc_1,\grc_2)=\langle \e , \bar{\M}_2(\grc_1,
\grc_2)\rangle. 
\]
This is the Euler number of the bundle \(\mathcal{E}(\bbR\times
M)|_{\bar{\M}_2(\grc_1, \grc_2)}\) relative to the trivialization
\(\rho_{\vartheta }|_{\partial\bar{\M}_2(\grc_1, \grc_2)\subset
  \M_1(\bbR\times M)}\). In comparison, 
Section IV.1.3's \(\op{w}_U(\grc_1,
\grc_2)\) 
is taken to be the signed count of elements in
\(\mathcal{M}_{2,p}(\grc_1, \grc_2)\), where \(\mathcal{M}_{k,p}(\grc_1, \grc_2)\subset  \mathcal{M}_{k}(\grc_1,
\grc_2)\) consists of elements \(\grd\in\mathcal{M}_{k}(\grc_1,
\grc_2)\) represented by some \((A, (\alpha , \beta ))\in \op{Conn}(E)\times \Gamma(\bbS^+)\) with \(\alpha \) 
vanishing at \(x=(0,p)\in X=\bbR\times M\). Suitable genericity
assumptions on \((\grT, \grS)\) and \(p\) were imposed so that for all
\(\grc_1, \grc_2\in \mathcal{Z}\), 
\(\mathcal{M}_{k,p}(\grc_1, \grc_2)=\emptyset\) for \(k<2\) and
\(\mathcal{M}_{2,p}(\grc_1, \grc_2)\) consists of finitely many
regular points. Let \(\td{\alpha }_x\in \Gamma (\mathcal{E}_x (\bar{\bbR}\times M))\) be
the section obtained by restricting the tautological section
\(\td{\alpha }\) to \(\pmb{\cal E}|_{\{x\}\times
  \mathcal{B}(\bar{\bbR}\times M)\subset X\times \mathcal{B}(\bbR\times
  M)}=\mathcal{E}_x (\bbR\times M)\). Then 
the space \(\mathcal{M}_{k,p}(\grc_1,
\grc_2)\) is precisely the zero locus of the section \(\td{\alpha
}_x\) on \(\mathcal{M}_{k}(\grc_1,
\grc_2)\subset\mathcal{B}_{loc}(\bbR\times M)\). The fact that \(\mathcal{M}_{k,p}(\grc_1, \grc_2)=\emptyset\)
\(\forall \grc_1, \grc_2\) for \(k<2\) implies that \(\tilde{\alpha}_x\) is nowhere-vanishing on \(\M_1(\bbR\times M)\), and hence
\(\td{\alpha }_x/|\td{\alpha }_x|\) defines a trivialization of
\(\mathcal{E}_x(\bar{\bbR}\times M)|_{\partial\bar{\M}_2 (\grc_1,
  \grc_2)\subset\M_1(\bbR\times M)}\), and the Euler number of the
complex line bundle \(\mathcal{E}_x(\bar{\bbR}\times M)|_{\bar{\M}_2 (\grc_1,
  \grc_2)}\) relative to this trivialization is precisely the Euler
characteristic of \(\mathcal{M}_{2,p}(\grc_1, \grc_2)=\bar{\M}_2 (\grc_1,
  \grc_2)\cap \td{\alpha }_x^{-1}(0)\), namey, the value of
  \(\op{w}_U(\grc_1, \grc_2)\) defined in \cite{KLT4}. This agrees
  with the expression from Part 2(a) if \(\vartheta \) therein is {\em
    chosen so
  that}  \(\td{\alpha }_x/|\td{\alpha }_x|\) is constant with respect to
  the trivialization \(\rho_\vartheta \) on \(\M_1(\bbR\times M)\). As
  observed in Part 1(a), the cocycle \(\e\in C^{2;\bbZ}_{\M(\bbR\times
    M)}\) depends on the \(\delta \)-cohomology class of \(\vartheta
  \), which in turn depends on the class \([\theta _p]_{\M (\bbR\times
    M)}\in C^{2;\bbZ}_{\M(\bbR\times
    M)}\). The aforementioned choice in the large-perturbation setting
  is natural in the sense that under
  proper setup, one expects
 \begin{equation}\label{large-r-U}
[\theta _p]_{\M (\bbR\times
    M)}\to 0 \quad \text{and} \quad [\ul{h}_{\hat{p}}]_{\N^+_1(M)}\to
  0 \quad \text{as \(r\to \infty\), }
\end{equation}
which in turn is based on the expection that, roughly speaking, 
\begin{equation}\label{large-r-flat}
\text{\(|\nabla_A\alpha |\to 0\) pointwise away from \(\alpha
  ^{-1}(0)\) as \(r\to \infty\)};
\end{equation}
or, put in another way, a variant of \cite{Ta1}'s Proposition 4.1 holds. A weak version of
the latter in the setting of Section \ref{sec:3} is provided in Lemma
\ref{lem:B.6}. 

To see how (\ref{large-r-U}) would follow from (\ref{large-r-flat}), recall (\ref{theta-h}) and note that as
\(\bar{\M}_1(\bbR\times M)\cap \tilde{\alpha }^{-1}(0)=\emptyset\), \(|\alpha |\big|_{\hat{p}}\) is
nowhere-vanishing for all \(\grd\in \in \M_1(\bbR\times M)\). Let \((A, (\alpha ,
\beta ))\in \mathcal{C}(\bbR\times M)\) is a representative of the
aforementioned 
\(\grd\), and use \(\hat{A}_\alpha \) to denote the connection defined
on \((\bbR\times M)\backslash\alpha ^ {-1}(0)\) satisfying \(\nabla_{\hat{A}_\alpha
}(\alpha /|\alpha |)=0\). 
Thus, for \(\grd\in \M_1(\bbR\times M)\)  
\[\begin{split}
h_{\hat{p}}(\grd)& =-\int_{\bbR\grd}\theta
_p= -\int_{\hat{p}}(\iota_\Delta
)^*\theta_{\dot{p}}\\
& =\frac{i}{2\pi}\int_{\hat{p}}(A-\hat{A}_\alpha  )\to 0 \quad \text{as
  \(r\to 0\)}
\end{split}
\]
if (\ref{large-r-flat}) holds, and (\ref{large-r-U}) follows as a consequence.

\paragraph{(b)} {\em The \(\grt_i\)-map associated to \(\gamma
  _i\subset M\).} According to Part 2(b), 
\begin{equation}\label{def:w_t}
\begin{split}
\op{w}_{\grt_i}(\grc_1, \grc_2)& =\langle\uu_{\dot{\gamma }_i},
\bar{\M}_1(\grc_1, \grc_2)\rangle\\
& =\sum_{\ul{\grd}\in\N_0(\grc_1,
  \grc_2)}\op{sign} (\ul{\grd}) \, (\Delta_{\ul{\grd}}\op{x}_{\gamma_i}) \\
& =
\langle\theta_{\dot{\gamma }_i}, \bar{\M}_1(\grc_1, \grc_2)\rangle
-\langle \Pi_0^*h'_{\gamma _i},  \partial\bar{\M}_1(\grc_1, \grc_2)\rangle.
\end{split}
\end{equation}
If \(\gamma _i\) lies on the
complement of \(\alpha ^{-1}(0)\) for all \((A, (\alpha  , \beta ))\)
representing elements \(\grc\) in \(\grC(M)\),  there is a natural choice of
\(h'_{\gamma _i}\co \grC(M)\to \bbR\) among the \(\bbZ^{\grC(M)}\)-many possible
lifts of \(\op{h}'_{\gamma _i}|_{\grC(M)}\), leading to a natural
choice of \(\uu_{\gamma_i}\). Namely, one sets 
\begin{equation}\label{large-r-h'}
h'_{\gamma _i}(\grc)=\frac{i}{2\pi}\int_{\gamma _i}(A-\hat{A}_\alpha)
\end{equation}
in this case. With this choice of \(h'_{\gamma _i}\), the
corresponding \(\op{x}_{\gamma _i}\) satisfies
\[
\text{ \(\op{x}_{\gamma _i}(\grc)=\op{hol}_{\gamma _i}(\hat{A}_\alpha )\)
 mod \(\bbZ\) \(\forall \grc\in \grC(M)\), }
\]
where \(\op{hol}_{\gamma
   _i}(\hat{A}_\alpha )\in U(1)=\bbR/\bbZ\) denotes the holonomy of
 \(\hat{A}_\alpha \) along \(\gamma _i\). 
Such \(\gamma _i\) can be found in large \(r\)-perturbation settings when a suitable variant of
(\ref{large-r-flat}) holds; in fact with such \(\gamma _i\), 
\begin{equation}\label{large-r-t}
h'_{\gamma _i}\to
0 \quad \text{ as \(r\to \infty\). }
\end{equation}
This may indeed be arranged in the setting of this series
of articles. In various
parts of \cite{KLT4} as well as in latter parts of this article
(e.g. (\ref{eq:(A.8)})), choices of \(\op{x}_{\gamma_i}\) were  made
via explicit formulae, and  (\ref{large-r-t}) in this context is given
a precise reformulation in terms of \(\op{x}_{\gamma_i}\) in Lemma
\ref{lem:A.1}. 

In part 7 of Section IV.1.3, the integer
\(\Delta _{\ul{\grd}}\op{x}_{\gamma_i}\) in (\ref{def:w_t})
is given an alternative description as the algebraic intersection number between
\(\alpha^{-1}(0)\) and the cylinder \(\bbR\times \gamma_i\subset
\bbR\times M\). To relate this with the definition in Part 2(b), note
that by the choice of \(\gamma _i\), the section \(\alpha
|_{\bar{\bbR}\times\gamma _i\subset \bar{\bbR}\times M}\) of the bundle \(E|_{\bar{\bbR}\times \gamma _i}\) is
nowhere-vanishing over the boundary of the cylinder \(\partial(\bar{\bbR}\times
\gamma _i)\subset\{-\infty,
\infty\}\times M\), and the aforementioned intersection number agrees
with the relative Chern number of \(E|_{\bar{\bbR}\times \gamma _i}\)
relative to the trivialization over \(\partial (\bar{\bbR}\times\gamma _i)\) defined by \(\alpha /|\alpha |\). This relative Chern
number in turn can be expressed as 
\[
\int_{\bar{\bbR}\times \gamma _i}\frac{i}{2\pi} F_{\hat{A}_\alpha
}=\Delta _{\ul{\grd}}\op{x}_{\gamma_i}. 
\]

\section{Filtered monopole Floer homologies}\label{sec:3}
\setcounter{equation}{0}

The algebraic recipe for Ozsvath-Szabo's definition of the four flavors
of Heegaard Floer homologies, labeled by the superscripts \(-,
\infty, +, \wedge\), was summarized abstractly in section 4 of
\cite{L}. In this section, we explain how the same recipe may be
applied in the Seiberg-Witten context to define analogs of
Ozsvath-Szabo's Floer homologies. These intermediate Floer homologies
play a pivotal role in the proof of Theorem \ref{thm:main}. 

\subsection{Motivation and sketches of construction}\label{sec:3.1}

The afore-mentioned recipe hinges on the existence of certain filtration
on a Floer chain complex with local coefficients in the group ring
\(\bbK[\bbZ]=\bbK[U, U^{-1}]\), with \(U\) corresponding to
the generator \(1\in\bbZ\). This Floer complex with local cofficients
constitutes the \(\infty\)-flavor of the Ozsvath-Szabo construction,
while the ``filtration'' refers to the filtration of the coefficient
ring \(\bbK[U,U^{-1}]\) by submodules
\[
\cdots U\, \bbK[U]\subset \bbK[U]\subset U^{-1}\bbK[U]\subset\cdots\subset\bbK[U, U^{-1}].
\]
If the differential of the \(\infty\)-flavor of the Floer complex
preserves this filtration, then it induces a filtration on the
\(\infty\)-flavor Floer complex by \(\bbK[U]\)-subcomplexes, which are all
isomorphic via multiplication by powers of \(U\).
This defines the \(-\)-flavor Floer complex. With these two basic
flavors in place, the \(+\)- and the \(\wedge\)-flavors are defined so
that they fit into short exact sequences (see (\ref{eq:fund-short})
below) inducing what are called the {\em fundamental
exact sequences} of corresponding Floer
homologies. In \cite{L}, the existence of such filtration is
attributed to the existence of what was termed a 
``semi-positive 1-cocycle''. The 1-cocycle used here refers to the
cocycle 
that defines the local system on the \(\infty\)-flavor of Floer
complex. The ``semi-positivity'' condition serves to guarantee 
that the differential is filtration-preserving. Note that the
\(\infty\)-flavor of Floer homology depends only on the cohomology
class of this cocycle. 
The
other three flavors of Ozsvath-Szabo's construction depend on
the choice the cocycle that defines the semi-positivity condition. 

Section 4.2 of \cite{L} provides some examples where this recipe may be
applied. Section 6 of the same
article sketched how such semi-positive 1-cocycles might arise
in certain versions of Seiberg-Witten Floer theory associated to 
equations of the form of (\ref{eq:(A.4)}). In particular, choosing the metric
and 2-form \(w\) in (\ref{eq:(A.4)}) to reflect the data that
go into the definition of Heegaard Floer homology provides a
bridge to relate the Heegaard and Seiberg-Witten Floer homologies. 

To elaborate, the local system underlying the Seiberg-Witten
analog of Ozsvath-Szabo construction is closely related to what was
denoted \(\Gamma_\eta\) in \cite{KM} (Cf. Example in the end of their
Section 22.6), where \(\eta\) is a singular
1-cycle in a certain 3-manifold \(\ul{M}\). Use \([(\bbA, \Psi)]\in \mathcal{B}^\sigma\) to
denote the gauge equivalence of \((\bbA, \Psi)\).
In \cite{KM}, this local system associates to each
point on \(\mathcal{B}^\sigma\) the ``fiber'' \(\bbR\), and to each
path \(\{[(\bbA(\tau), \Psi(\tau))]\}_\tau\) from 
\([(\bbA_-, \Psi_-)]\) to \([(\bbA_+, \Psi_+)]\), an isomorphism
\(\bbR^\times\subset \op{End}(\bbR)\) between the fibers over the end
points. The latter isomorphism is given by multiplication by the real number
\begin{equation}\label{hol-cocycle}
e^{\frac{i}{2\pi}\int_\tau\int_\eta\frac{d}{d\tau}\bbA(\tau)}.
\end{equation}
Note that the exponent is the difference of the holonomy of \(\bbA_-\)
along the cycle \(\eta\subset M\) from that of \(\bbA_+\), and it defines a
real 1-cocycle in \(\mathcal{B}^\sigma\). 
Meanwhile, as only points in \(\grC\subset
\mathcal{B}^\sigma\) and paths constituting the sets \(\mathcal{M}_1(\grc_-,
\grc_+)\), \(\grc_-\), \(\grc_+\in \grC\) enter the definition of a
monople Floer complex, it suffices to
consider the holonomy difference of paths corresponding to 
 elements in \(\mathcal{M}_1(\grc_-,\grc_+)\).  
The observation leading to \cite{L}'s construction of filtered monopole
Floer homologies (in the sense of Ozsvath-Szabo) is the following: 
\begin{equation}\label{L-observ}
\begin{split}
&\text{for monopole
Floer complexes associated to certain \(\varpi\) in the form of}\\
& \text{ 
(\ref{eq:2rw}) with large \(r\) and certain choice of \(\eta\), 
the value of the afore-mentioned
}\\
&\text{ holonomy difference is very close to a non-negative integer. }
\end{split}
\end{equation}
(Cf. also (\ref{large-r-t}).)

Associating to each
element in \(\mathcal{M}_1(\grc_-,\grc_+)\) its corresponding
integer, one has a (partially defined) integer 1-cocycle on
\(\mathcal{B}^\sigma\) with which one may define a Floer complex with
more refined local coefficients  than \(\Gamma_\eta\). We denote the
latter local system by \(\Lambda_\eta\). It replaces the fibers \(\bbR\) over
\(\grC\) of \(\Gamma_\eta\) by the group ring \(\bbK[\bbZ]=\bbK[U, U^{-1}]\); and
it replaces the isomorphism induced by an element \(\grd\) in
\(\mathcal{M}_1(\grc_-,\grc_+)\) between these fibers, namely
(\ref{hol-cocycle}), by \(U^n\) where \(n\) denotes the
afore-mentioned non-negative integer associated to \(\grd\). The fact that
\(n\geq 0\) in all cases has the following consequence: Use
the corresponding monopole
Floer complex with local coefficients, \(\Lambda_\eta\), as the {\em \(\infty\)-flavor Floer complex}. There
is filtration on this chain complex, \(\CM_* (\ul{M}; \Lambda_\eta)\), by subcomplexes of
\(\bbK[U]\)-modules. This can be used to define the other three flavors of Floer
complexes.

The program described in \cite{L} assumes various plausible
conjectures and assertions that come from an extension of  
the geometric picture in the last author's work relating the
Seiberg-Witten and Gromov invariant for closed 4-manifolds
(Cf. \cite{T}, \cite{Ta1}). A proof of these conjectures 
constitute a major part of the technical hurdle for
implementing the program in \cite{L}. 
The difficulties arise because the 2-form \(\varpi\) in \cite{L} must
have zeros. 

In this series of articles \cite{KLT1}-\cite{KLT4}, the road block to
the approach in \cite{L} is circumvented by a
modification of \cite{L}'s outline. Very roughly, the manifold
\(\ul{M}\) in \cite{L} is replaced by the manifold denoted by \(Y\) in
\cite{KLT2}. This is obtained from \(\ul{M}\) by adding further
1-handles along the zeros of \(w\) on \(\ul{M}\). The 2-form \(w\) 
extends into \(Y\) as a no-where vanishing closed 2-form, which we also denote by \(w\). Over the middle of the added 1-handle, this
\(w\) approximates \(da\) for a certain contact form \(a\), and as the
special 1-cycle \(\eta\) (denoted \(\ul{\gamma}\) therein) lies away
from the zeros of \(w\) on \(\ul{M}\), this 1-cycle also embeds in
\(Y\).  This was denoted by \(\gamma^{(z_0)}\) in
\cite{KLT1}-\cite{KLT4}. The technical challenge in this new approach
involves, among other things, the analog of 
(\ref{L-observ}) for the monopole Floer complex associated to \(Y\),
\(w\), and \(\eta=\gamma^{(z_0)}\). Some of these technical issues are
dealt with in \cite{KLT4}. Those that remain are dealt with in
Section \ref{sec:B}-\ref{sec:D} of this article. 

In Section \ref{sec:Aa)} below, we specify the class of 3-manifolds,
denoted \(Y_Z\) therein, 
together with the 2-form \(w\) on it and the 1-cycle \(\eta\) for
which positively results of the kind (\ref{L-observ}) hold. Section
\ref{sec:Ac)} describes the sort of cobordisms \(X\) for which the companion
statements hold. Cf. Propositions \ref{prop:A.2}, \ref{prop:A.3},
\ref{prop: A.7}, \ref{prop:A.9}, \ref{prop: A.11}.          
The remaining subsections give precise statements of
the desired positivity results. The formulation here involves an
``cut-off'' version of the connection \(A\) (called \(\hat{A}\)), so that in
place of (\ref{L-observ}), its associated holonomy difference is
integer-valued. (Cf. Lemma \ref{lem:A.1}).
The conditions on \(Y_Z\) and \(X\)
are introduced more for technical convenience rather than essential
reasons, and the statements in Sections \ref{secAe)}-\ref{sec:Ah)} may conceivably hold for
more general 3-manifolds and 4-dimensional cobordisms.

\subsection{The 3-manifold \(Y_Z\)}\label{sec:Aa)}

Let \(Z\) denote a given connected, oriented closed 3-manifold; and let
\(Y_{Z}\) denote the manifold that is obtained from \(Z\) by
attaching a 1-handle at a chosen pair of points, denoted \((p_0,
p_3)\) below.  In the proof of the main theorem \ref{thm:main}, \(Z\) 
is taken to be either \(S^{3}\), the manifold \(M\) in the statement
of Theorem \ref{thm:main}, or a manifold that
is obtained from \(M\) by attaching some number of 1-handles.
 Although \(Y_{Z}\) is diffeomorphic to the connected sum of
\(Z\) and $S^{1} \times S^{2}$, it is
viewed for the most part as $Z_{\delta}\cup \mathcal{H}_{0}$ with $\mathcal{H}_{0}$ the
attached 1-handle and with $Z_{\delta}$ being the
complement of a pair of coordinate balls about the chosen points
\(p_{0}\) and \(p_{3}\) in \(Z\).   The manifold
\(Y_{Z}\) has a distinguished embedded loop that crosses the
handle $\mathcal{H}_{0}$ once.  This loop is denoted by
$\gamma$.  The three parts of this subsection say more about the
geometry of \(Y_{Z}\) near $\mathcal{H}_{0}$,
near $\gamma$, and in general.

\paragraph{Part 1:}  The geometry of \(Y_{Z}\) near
$\mathcal{H}_{0}$ is just like that given in Section
II.1a.  By way of a reminder, the description of the geometry requires
the a priori specification of constants $\delta_{*}\in (0, 1) $
and \(R > -100\ln \delta_{*}\).  Also needed are coordinate charts
centered on \(p_{0}\) and \(p_{3}\).  The latter are
used to identify respective neighborhoods of these points with balls of
radius $10\delta_{*}$ in
$\mathbb{R}^{3}$.  The pull-back of the standard
spherical coordinates on $\mathbb{R}^{3}$ gives
spherical coordinate functions on the neighborhood of
\(p_{0}\), these denoted by \((r_{+}, (\theta_{+}, \phi_{+}))\).  There are
corresponding coordinate functions for the neighborhood of
\(p_{3};\) these are denoted in what follows by \((r_{-}, (\theta_{-}, \phi_{-}))\).

The handle $\mathcal{H}_{0}$ is diffeomorphic to the
product of an interval with \(S^{2}\).  The interval factor
is written as \([-R - 7\ln \delta_{*}, R +
7\ln \delta_{*}]\) and \(u\) is used to denote the
Euclidean coordinate for this interval.  The spherical coordinates for
the \(S^{2}\) factor are written as \((\theta, \phi)\).
 The handle $\mathcal{H}_{0}$ is attached to the
coordinate balls centered on \(p_{0}\) and \(p_{3}\)
as follows:  Delete the \(r_{+} < e^{ -2R} (7\delta_{*})^{ -1}\) part
of the
coordinate ball centered on \(p_{0}\) and the corresponding
part of the coordinate ball centered on \(p_{3}\).  Having
done so, identify $\mathcal{H}_{0}$ with the respective
\(r_{+} \in [e^{ -2R}(7\delta_{*})^{ -1},7\delta_{*}]\) and \(r_{ -} \in
[e^{ -2R}(7\delta_{*})^{ -1}, 7\delta_{*}]\) parts of these coordinate balls with
$\mathcal{H}_{0}$ by writing

\begin{equation}\label{eq:(A.1)}
\begin{split}
&(r_{+} = e^{-R+u},
(\theta_{+} = \theta, \phi_+ =\phi))  \quad  \textit{and }\\
& (r_{-} =e^{-R-u}, (\theta_{-} = \pi  -\theta, \phi_{ -} = \phi) ).
\end{split}
\end{equation}

The handle $\mathcal{H}_{0}$ has a distinguished closed
2-form, this being  \(\frac{1}{2} \sin \theta d\theta d\phi\).  This 2-form is nowhere zero on
the constant \(u\) cross-sectional spheres and thus orients these spheres.
 Granted this orientation, then   
$\frac{1}{2}\sin \theta d\theta  d\phi$ has integral 2 over constant \(u\)
sphere. 

  \paragraph{Part 2:}  The loop $\gamma$ intersects
$\mathcal{H}_{0}$ as the \(\theta = 0\) arc.  Thus it has
geometric intersection number \(1\) with each \(u =\) constant sphere.  This
loop is oriented so that the corresponding algebraic intersection
number is \(+1\).  A tubular neighborhood of $\gamma$ is specified with
a diffeomorphism to the product of \(S^{1}\) and a disk
about the origin in $\mathbb{C}$.  The latter is denoted by
\(D_{\gamma}\) and its complex coordinate is denoted by \(z\).
 The diffeomorphism identifies the \(z = 0\) circle in
 $S^{1}{\times} D $ with $\gamma$.  The circle
\(S^{1}\) is written in what follows as $\mathbb{R}/(
\ell_\gamma\mathbb{Z} )$ with \(\ell_\gamma> 0\) being a chosen constant.  The affine coordinate for
$\mathbb{R}/(\ell  _\gamma\bbZ )$ is denoted by \(t\).  The product structure on such a
neighborhood is constrained where it intersects
$\mathcal{H}_{0}$ by the requirement that the
\(\mathcal{H}_{0}\) coordinate \(u\) on the intersection depend
only on \(t\).  A neighborhood with these coordinates is fixed once and
for all; it is denoted by \(U_{\gamma}\).   

\paragraph{Part 3:} 
Use the Mayer-Vietoris principle to write the second homology of \(Y_{Z}\) as
\begin{equation}\label{eq:(A.2)}
 H_{2}(Y_{Z}; \bbZ) =H_{2}(Z; \bbZ) \oplus H_{2}(\mathcal{H}_{0}; \bbZ).
\end{equation}
The convention in what follows is to take the generator of
$H_{2}(\mathcal{H}_{0}; \bbZ)$ to be
the class of any cross-sectional sphere with the orientation given by
the 2-form \(\sin\theta d\theta d\phi\).  Fix a class in \(H^{2}(Y_{Z};
\bbZ)\) which has even pairing with the classes in
\(H_{2}(Y_{Z}; \bbZ)\) and pairing \(2\) with
the generator of the $H_{2}(\mathcal{H}_{0};\bbZ)$ summand in (\ref{eq:(A.2)}).  This class is denoted in what follows
by $c_{1}(\det\, (\mathbb{S}))$, and it is necessarily non-torsion by
the above assumption.  

There is a corresponding,
closed 2-form on \(Y_{Z}\) whose de Rham cohomology class is
that of $c_{1}(\det\, (\mathbb{S}))$.  In particular, there
are forms \(w\) of this sort satisfying the following additional
constraints:
\BTitem\label{eq:(A.3a)}
\item The form restricts to $\mathcal{H}_{0}$
as \(\frac{1}{2\pi}\sin\theta \, d\theta \, d\phi \); 
\item The form restricts to \(U_{\gamma}\)
as  \(\frac{i}{2\pi} g( |z |)\, dz \wedge d\bar{z}\)
  with \(g\) denoting a strictly positive function.  
\item There is a closed 1-form on \(Y_{Z}\), typically denoted by
  \(\upsilon\) below, with the following properties:
\begin{itemize}
\item[a)]   It has non-negative wedge product with \(w\).
\item[b)]   It restricts to \(U_{\gamma}\)  as \(dt\), and restricts to $\mathcal{H}_{0}$  as
\(\textsc{h}(u) \, du\)  with \(\textsc{h}(u) > 0 \)
for all \(u\).
\end{itemize}
\ETitem
Fix such a 2-form as the perturbation form \(w\) in
(\ref{eq:(A.4)}). 

The metric on \(Y_Z\) is chosen to satisfy the following constraints: 
\BTitem\label{eq:(A.3b)}
\item The metric appears on
$\mathcal{H}_{0}$  as the product metric
 of an \(S^{2}\)-independent metric on the
interval \([-R  - \ln (7\delta_{*}), R +\ln (7\delta_{*})] \) and the round metric
\(d\theta^{2} + \sin^{2}\theta d\phi^{2} \) on the \(S^{2}\) factor.  Meanwhile, the curvature 2-form of
\(A_{K}\)  on $\mathcal{H}_{0}$ is \(\frac{i}{2\pi} \sin\theta\, d\theta\, d\phi\).
\item The metric appears on \(U_{\gamma}\) as \(dt^{2} +\, g(|z|)  \, dz \otimes d\bar{z}\)
with \(g\) being the function in the second bullet of (\ref{eq:(A.3a)}).
Meanwhile, \(A_{K}\)  has  holonomy \(1\)  on \(\gamma\)  and its curvature 2-form on
\(U_{\gamma}\)  is \(iw\).
\ETitem
Many of the lemmas and propositions in the rest of this section depend
implicitly on the radius of \(D_{r}\) and on the
injectivity radius of the Riemannian metric.  They also depend
implicitly on the norms of \(w\), the
curvature of \(A_{K}\), the Riemannian curvature, and the norms of
their derivatives up to some order less than \(10\).  

There are suitable choices for
\(\mu\) with positive but small as desired $\mathcal{P}$-norm that
vanish on $\mathcal{H}_{0}\cup U_{\gamma}$ .  This last property is not a direct
consequence of an explicit assertion in \cite{KM} but it follows
nonetheless from their constructions.
 
The reference connection \(A_{E}\) is chosen constrained only to the extent
that it is flat on $\mathcal{H}_{0}$ and is flat with
holonomy 1 on \(U_{\gamma}\).  

The function on  \(\op{Conn}\, (E) \times C^{\infty}(Y_{Z}; \bbS)\) of
central concern in what follows is the analog here of the function that
is defined in (IV.1.16).  This function is denoted by \(\textsc{x}\).
 The definition requires the a priori choice of a smooth function
$\wp: [0, \infty) \to[0, \infty)$  which is
non-decreasing, obeys \(\wp(x) = 0\) for \(x <  \frac{7}{16}\)
and \(\wp(x) = 1\) for \(x \geq   \frac{9}{16}\).  As in \cite{KLT4}, it proves convenient to choose \(\wp\) so that its
derivative, \(\wp'\), is bounded by \(2^{10} (1  -
\wp)^{3/4}\).  The definition of \(\textsc{x}\) uses the
fact that \(w\) is nowhere zero on \(U_{\gamma}\).  In
particular, Clifford multiplication by \(*w\) on
\(U_{\gamma}\) splits \(\bbS \) over
\(U_{\gamma}\) as the direct sum of eigenbundles.  This
splitting is \begin{equation}\label{split-S}
\bbS = E \oplus E\otimes K^{ -1}
\end{equation}
with
the convention being that \(*w\) acts as \(+i|w|\) on \(E\).  A given section
\(\psi\) of \(\bbS\) is written with respect to this splitting over
\(U_{\gamma}\) as a pair denoted by \(|w|^{1/2}(\alpha,\beta)\).

Granted this notation, use \(\wp\) with a given pair \(\grc = (A, \psi)
\in \op{Conn}\,  (E) \times C^{\infty}(Y_{Z}; \bbS)\) to define the connection
\begin{equation}\label{eq:(A.7)}
\hat{A} = A -\frac{1}{2}  \wp(|\alpha|^{2})|\alpha|^{ -2} (\bar{\alpha}\nabla  _A\alpha-\alpha\nabla_{A} \bar{\alpha} ),
\end{equation}
on \(E|_{U_\gamma}\).  The salient point is that the connection
\(\hat{A}\) is flat on the part of \(U_\gamma\) where \(|\alpha|^2>\frac{9}{16}\)
 (this is where \(\wp=1\)) and the \(A\)-derivative of \(\alpha
 /|\alpha |\) is zero on this same part of \(U_\gamma \). 
  This can be seen from the following formulas: \BTitem\label{(3.9)v2}
\item \(F_{\hat{A}}=(1-\wp(|\alpha|^2))\, F_A+\wp' (|\alpha|^2)\nabla_A\alpha \wedge\nabla_A\bar{\alpha }\); 
 \item \(\nabla_{\hat{A}}\alpha =(1-\wp(|\alpha|^2))\nabla_A\alpha +\wp (|\alpha|^2)\, d(\ln |\alpha
   |) \, \alpha \).
\ETitem
Meanwhile, the connections \(\hat{A} \)and \(A\) are equal where
\(|\alpha |^2\leq \frac{7}{16}\) (this is where \(\wp=0\)).

With \(\hat{A}\) understood, then 
the value of the function \(\textsc{x}=\textsc{x}_\gamma\) on the
given configuration \(\grc= (A, \psi)
\in \op{Conn}\,  (E) \times C^{\infty}(Y_{Z}; \bbS)\) is defined by
rule whereby
\begin{equation}\label{eq:(A.8)}
 \textsc{x}(\grc) =   \frac{i}{2\pi}\int_\gamma (\hat{A}-A_E).
\end{equation}
\begin{remarks}
To relate with the general discussion in Part 2(b) of Section
\ref{sec:A-module}, note that \(\hat{A}\) from (\ref{eq:(A.7)})
agrees with the connection \(\hat{A}_\alpha \) over \(\gamma\); and so
setting \(\varepsilon _\gamma '(\grc)=\int_\gamma (\hat{A}-A)\) for
\(\grc\in \B^\sigma (M)\) would meet the requirement that  \(\varepsilon _\gamma
'|_{\grC(M)}=h'_\gamma \) when \(h'_\gamma \) is given by (\ref{large-r-h'}). 
Meanwhile, the
reference connection \(A_E\) plays the role of 
the base point \(\hat{\grc}_0\in \hat{\B}^\sigma (M)\) in Section
\ref{sec:A-module} in the following sense: Let \((A_0, (\alpha _0,
\beta _0))\) be an arbitrary representative of \(\hat{\grc_0}\) and
fore any \(\hat{\grc}\in \hat{\B}^\sigma _{\grt}(M)\), let \((A,
(\alpha , \beta ))\) be an arbitray representative of
\(\hat{\grc}\). Then 
\(\op{x}_\gamma (\hat{\grc})\), as defined in Section
\ref{sec:A-module}'s Part 2(a), equals 
\[\begin{split}
\op{x}_\gamma (\hat{\grc})& =\hat{h}_\gamma
(\hat{\grc})-\Pi_0^*\varepsilon '_\gamma (\hat{\grc})\\
& =\frac{i}{2\pi}\Big( \int_\gamma (A-A_0)-\int_\gamma  (A-\hat{A})\Big)
=\frac{i}{2\pi}\int_\gamma (\hat{A}-A_0).
\end{split}
\]
The last term above equals \((\ref{eq:(A.7)})\) when \(A_0=A_E\). Note
that \(\hat{\grc}_0\) and \(A_E\) are required to satisfy consist
constraints; namely, both \(\op{hol}_\gamma (A_0)=0\) \(\mod \bbZ\) and
\(\op{hol}_\gamma (A_E)=0\) \(\mod \bbZ\). 
\end{remarks}

The following lemma supplies a fundamental observation about
\(\textsc{x}.\)
\begin{lemma}\label{lem:A.1}  
If the conditions in (\ref{eq:(A.3a)}), (\ref{eq:(A.3b)}) hold, then
there exists \(\kappa > \pi \)  with the
following significance:  Fix \(r > \kappa  \) and a 1-form \(\mu \in \Omega  \) with
$\mathcal{P}$-norm less than 1.  The function
\(\textsc{x}\)  has only integer values on the solutions to the
corresponding \((r,  \mu)\)-version of (\ref{eq:(A.4)}). 
\end{lemma}
This lemma is proved in \S\ref{sec:Bc)}.

\subsection{4-dimensional cobordisms}\label{sec:Ac)}

This subsection describes in general terms the
sorts of cobordisms that are considered. 

To start, let \(Z_{-}\) and \(Z_{+}\) denote two
versions of the manifold \(Z\) and let \(Y_{-}\) and
\(Y_{+}\) denote the respective \(Z = Z_{-}\) and \(Z
= Z_{+}\) versions of \(Y_{Z}\).  There is no need
to assume that either \(Y_{-}\) or \(Y_{+}\) is
connected, but if not, then the handle $\mathcal{H}_0$
is attached to the same connected component.  Use
$\gamma_{-}$ to denote the \(Y_{-}\)
version of the curve $\gamma$ and use $\gamma_{+}$ to
denote the \(Y_{+}\) version.  The corresponding versions of
\(U_{\gamma}\) are denoted in what follows by \(U_{-}\) and \(U_{+}\).  
 
Of interest here is a smooth, oriented, 4-dimensional manifold \(X\)
with the properties listed below, in addition to those in (\ref{(A.9a,11)}): 
\BTitem\label{eq:(A.9b)}
\item  
 There exists an embedding of \(\bbR \times [-R  -\ln (7\delta_{*}),  R +\ln (7\delta_{*}) ] \times S^{2}\)
 into \(X\)  that pulls back \(s\)  as the Euclidean coordinate on the \(\bbR\)-factor.  Moreover, the
composition of this embedding with the diffeomorphism in the second
bullet identifies the \(s < 0\)  part with $(-\infty,0) \times \mathcal{H}_0$  in \((-\infty,0) \times Y_{-}\); and the composition with
the diffeomorphism from the third bullet identifies the \(s> 0\)  part with 
$(0, \infty) \times\mathcal{H}_0$  in \((0, \infty) \times Y_{+}\).
\item  
 There exists an embedding of \(\bbR \times
S^{1}\)  into \(X\)  that pulls back \(s\) 
as the Euclidean coordinate on the \(\bbR\)-factor.
 Moreover, the composition of this embedding with the  diffeomorphism
in the second bullet identifies the \(s < 0\)  part of
\(\bbR \times S^{1}\)  with \((-\infty, 0) \times  \gamma_{-}\);
 and the composition with the diffeomorphism from the third
bullet identifies the \(s > 0\)  part of \(\bbR \times S^{1}\)  with \((0, \infty)\times \gamma_{+}\). 
\ETitem

The image in \(X\) of the embedding of \(\bbR
\times [-R  - \ln (7\delta_{*}),  R +
\ln (7\delta_{*}) ] \times S^{2}\)
from the first bullet above is denoted by \(U_{0}\).

 The notation used in the next constraint has \(C\) denoting the image in \(X\)
of $\mathbb{R}\times S^{1}$  as described by the
second bullet of (\ref{eq:(A.9b)}).  This constraint requires that the
\(\gamma_{-}\) and \(\gamma_{+}\) versions of  \(\ell_\gamma\) are equal. 

\BTitem\label{eq:(A.10)}
\item  
 There exists \(\ell_\gamma > 0\)  and  a diffeomorphism of a neighborhood of
\(C\)  to the product of \(\bbR \times \bbR/( \ell_\gamma\bbZ)\)  with a disk about the origin in
$\mathbb{C}$.  This disk is denoted by \(D\).
\item  
 The diffeomorphism identifies the Euclidean coordinate on
\(\bbR \times \bbR/(\ell_\gamma\bbZ) \times D\)  with \(s\). 
\item  
 The respective \(s < 0\)  and \(s >0\)  parts of the neighborhood are in  \((-\infty, 0)
\times U_{-}\)  and in \((0 , \infty)\times U_{+}\).   Moreover, the diffeomorphism
on these parts of the neighborhood respects the respective splittings
of \(U_{-}\)  and \(U_{+}\)  as
\((-\infty, 0)\times \bbR/(\ell _\gamma\bbZ) \times D\)  and \((0 ,
\infty)\times \bbR/(\ell _\gamma   \bbZ) \times D\).
\ETitem
 
By way of an explanation, a diffeomorphism of this sort exists if the
co-normal bundle to \(C\) in \(X\) has a nowhere zero section that restricts to
the \(s < 0\) part of \(X\) as the real part of the
$\mathbb{C}$-valued 1-form \(dz\) along \(\gamma_{-}\) and
restricts to the \(s > 0\) part of \(X\) as the real part of the
$\mathbb{C}$-valued 1-form \(dz\) along \(\gamma_{+}\).
 The tubular neighborhood in (\ref{eq:(A.10)}) is denoted in what follows as
\(U_{C}\).  The diffeomorphism in (\ref{eq:(A.10)}) is used, often
implicitly, to identify \(U_{C}\) with \(\bbR \times\bbR/(\ell _\gamma\bbZ) \times D\).
 
In addition to those listed in (\ref{eq:(A.13b)}),
the 2-form \(w_X\) to use in the Seiberg-Witten equations is required
to satisfy the following additional constraint: 
\begin{equation}\label{eq:(A.13c)}
\begin{split}
 & \text{The pull-back of \(w_{X}\)  to
\(U_{0}\)  via the embedding from the fourth bullet}\\
&\text{of
(\ref{eq:(A.9b)}) is twice the self dual part of 
\(\frac{1}{2} \sin \theta d\theta d\phi   \) }\\
&\text{and its pull-back to
\(U_{C}\)  via the embedding in (\ref{eq:(A.10)})}\\
&\text{is twice the
self dual part of \(\frac{i}{2\pi} g(|z|)\,  dz \wedge d \bar{z}\).}
\end{split}
\end{equation}

Meanwhile, the metric on \(X\) is required to satisfy the following
constraints in addition to those in (\ref{eq:(A.12,15a)}):

\BTitem\label{eq:(A.15b)}
\item  
 The metric pulls back from \(U_{0}\)  via the
embedding of the first bullet of (\ref{eq:(A.9b)}) as the product metric defined
by the Euclidean metric on the \(\bbR\)-factor and an
\(\bbR\)-independent product metric on the \([-R  -\ln (7\delta_{*}),  R +\ln (7\delta_{*}) ] \times S^{2}\)  factor.
\item  
 The metric pulls back from \(U_{C}\)  via the
embedding in (\ref{eq:(A.10)}) as the product metric given by the quadratic form
\(ds^{2} + dt^{2} +g(|z|) \, dz \otimes d\bar{z}\).  
\ETitem

Extensions to \(U_{C}\) and \(U_{0}\) of the
\(Y_{-}\) and \(Y_{+}\) versions of the line
bundles \(K\) and \(E\) and their connections \(A_{K}\) and
\(A_{E}\) are needed for what follows.  There is no
obstruction to making these extensions.  Even so, it is necessary to
constrain \(A_{K}\) and \(A_{E}\) on
\(Y_{-}\) and \(Y_{+}\) so that extended versions
of \(A_{K}\) and \(A_{E}\) on \(U_{C}\cup U_{0}\)
 exist with the curvature of the extended
version of \(A_{K}\) pulling back via the embeddings from the
first bullet of (\ref{eq:(A.9b)}) and (\ref{eq:(A.10)}) as 
 \(\sin\theta\, d\theta\, d\phi \) and \(g(|z|) \, dz \wedge d\bar{z}\).  Meanwhile, the pull-backs of the curvature of \(A_{E}\)
via these embeddings is zero.  Extensions with this property are
assumed implicitly.     

The definitions in
\cite{KM} are sufficiently flexible so as to allow for the following:  For
any given $r > \pi $ , there are suitable perturbation
terms for (\ref{eq:(A.14)}) with positive but as-small-as desired
$\mathcal{P}$-norm that vanish on \(U_{C}\) and on the image
of $\bbR \times \mathcal{H}_0$ via the embedding map from the first bullet of (\ref{eq:(A.9b)}).   

 With regards to notation and conventions, the propositions and
lemmas that follow refer only to (\ref{eq:(A.14)}).  Even so, all assertions
still hold for the versions with an extra perturbation term if the
perturbation term has $\mathcal{P}$-norm bounded by  \(e^{-r^2}\)
or has small, \(r\)-independent $\mathcal{P}$-norm and vanishes on
\(U_{C}\) and on the image of $\bbR \times\mathcal{H}_0$ via the embedding from the first bullet
of (\ref{eq:(A.9b)}).  Proofs of the propositions and lemmas will likewise refer
only to (\ref{eq:(A.14)}).  The modifications that are needed to deal with the
extra perturbation terms are straightforward and so left to the
reader.

The second set of constraints require the choice of constants
\(\c \geq 1\) and \(\r \geq 1\).  By way of notation, one of the upcoming
constraints uses the embeddings from the second and third bullets of
(\ref{(A.9a,11)}) to write \(w_{X}\) on the \(|s | 
\in[L  - 4, L]\) part of \(X\) as \(w_{X} =  ds \wedge*w_{*} + w_{*}\) with
\(w_{*}\) denoting a closed, \(s\)-dependent 2-form on
\(Y_{-}\) or \(Y_{+}\), and with \(*\) here denoting
the Hodge star for the metric \(\grg\) in the second bullet of
(\ref{eq:(A.12,15a)}). 

\BTitem\label{eq:(A.16)}
\item[1)] The constant \(L\)  in (\ref{eq:(A.12,15a)})
is less than \(\c\). The constant \(L_{tor}\) in (\ref{eq:(A.13b)}) is
equal to \(\c\ln \r\). 
\item[2)] 
 The norm of the Riemannian curvature tensor and those of its
covariant derivatives up to order 10 are less than \(\r^{1/\c}\) on the \(s \in [-L, L]\)  part of \(X\).
\begin{itemize}
\item[a)]   The injectivity radius is larger than \(\r^{-1/\c}\) on the \(s \in [-L, L]\)  part of \(X\).
\item[b)]   The metric volume of the \(s\)-inverse image in
\(X\)  of any unit interval is bounded by \(\c\).
\end{itemize}
\item[3)]  The metric \(\grg\)  from (\ref{eq:(A.12,15a)})'s second bullet
obeys \(|\frac{\partial }{\partial s}\grg| \leq   \r^{1/\c}\).
 \item[4)] The norm of \(w_{X}\)  is bounded by \(\c\).  The norms of its
  covariant derivatives to order 10 are bounded by \(\r^{1/\c}\) on the \(s \in [-L, L]\)  part of \(X\).
\begin{itemize}
 \item[a)] The 2-form \(w_{X}\)  is closed on the
 \(|s | \leq L  - 4\)  part of \(X\). 
 \item[b)] Use the
 embeddings from the second and third bullets of (\ref{(A.9a,11)}) to write \(w_{X}\)  on the \(|s| \in [L  - 4, L]\)  parts of \(X-X_{tor}\)  as \(w_{X}= ds 
 \wedge *w_{*} + w_{*}\). Then \(\frac{\partial}{\partial s}w_{*} = d\b\)  where \(\b\)  is a smooth,
 \(s\)-dependent 1-form on the relevant components of \(Y_{-}\)  or
 \(Y_{+}\)  with \(\int_{(X-X_{tor})\cap|s|^{-1}([L  -4, L])}|\b|^2
 <\r^{-1/\c}\).
 \item[c)] The 2-form \(w_{X}\)  is closed on the components of the 
 \(L-4\leq|s | \leq L_{tor} - 4\)  part of \(X_{tor}\). 
\item[d)] Use the
embeddings from the second and third bullets of (\ref{(A.9a,11)}) to write \(w_{X}\)  on the \(|s| \in [L_{tor}  - 4, L_{tor}]\)  parts of \(X_{tor}\)  as \(w_{X}= ds 
\wedge *w_{*} + w_{*}\). Then \(\frac{\partial}{\partial s}w_{*} = d\b\)  where \(\b\)  is a smooth,
\(s\)-dependent 1-form on the relevant components of \(Y_{-}\)  or
\(Y_{+}\)  with \(\int_{X_{tor}\cap|s|^{-1}([L_{tor}  -4,
  L_{tor}])}|\b|^2 <\r^{-1/\c}\).
\end{itemize}
\item[5)]  There is a smooth, closed 1-form on \(X\), denoted by
\(\upsilon_{X}\) below, with norm bounded by \(\c\) and such that:
\begin{itemize}
\item[a)] The pull-back of \(\upsilon_{X}\) 
to \((-\infty, -L] \times Y_{-}\)  and to \([L, \infty)  \times Y_{+}\)  via the
embeddings from the second and third bullets of (\ref{(A.9a,11)}) is an
\(s\)-independent 1-form on \(Y_{-}\)  and \(Y_{+}\).
\item[b)] The pull back of \(\upsilon_{X}\) 
 to \(U_{C}\) via the embedding from (\ref{eq:(A.10)}) is
\(dt\)  and its pull-back to  \(U_{0}\) via the embedding from the first bullet of
(\ref{eq:(A.9b)}) is \(\textsc{h}(u) du\)  with \(\textsc{h}(\cdot) \geq\c^{-1}\).
\item[c)]   \(*(ds \wedge \upsilon_{X}  \wedge w_{X}) \geq -
  \r^{-1/\c}\) on the \(|s | \in [L  - 4, \infty)\) part of \(X\).
\end{itemize}
\ETitem 

Note that item 4) of the preceding constraints ensures that the condition
(\ref{pert-cobord}) holds. 
\begin{defn}
The metric and \(w_{X}\) on \(X\) are said to be {\em \((\c, \r)\)-compatible} when
one of the following conditions are met:  
\BTitem\label{eq:(A.17)}
\item  
 The space \(X = \bbR \times Y_{Z}\);  the metric has the form
\(ds^{2} + \grg\)  with \(\grg\)  being an
\(s\)-independent metric on \(Y_{Z}\); and the
2-form \(w_{X}\)  is the \(s\)-independent
form \(ds \wedge *w + w\).   Moreover, there exists a closed 1-form on \(Y_{Z}\),
denoted by \(\upsilon\) below, that restricts
to \(U_{\gamma}\)  as \(dt\), and restricts to $\mathcal{H}_0$  as
\(\textsc{h}(u) \, du\)
 with \(\textsc{h}(\cdot) >\c^{-1}\) , and is such that
\(\upsilon \wedge w \geq - \r^{-1/\c}\).
\item  The metric and \(w_{X}\)  obey the constraints
in (\ref{eq:(A.12,15a)}), (\ref{eq:(A.13b)}), (\ref{eq:(A.13c)}), (\ref{eq:(A.15b)}) and (\ref{eq:(A.16)}).
\ETitem
\end{defn}
By way of a look ahead, the notion of \((\c, \r)\)-compatibility is invoked
below with \(\r\) given by the constant \(\textrm{r}\) in (\ref{eq:(A.14)}).

\subsection{Positivity on cobordisms}\label{secAe)}

 An analog of the connection that is defined in (\ref{eq:(A.7)}) plays a
role in what follows. This connection is denoted in what follows by \(\hat{A}\).
 To define it, keep in mind that \(w_{X} \neq 0\) on \(U_{C}\) and so
 Clifford multiplication by \(w_{X}^{+}\) on \(\bbS^{+}\) over \(U_{C}\) or
$(-\infty, -2] \times \mathcal{H}_0$ or $[2, \infty) \times \mathcal{H}_0$ splits
\(\bbS^{+}\) as a direct sum of eigenbundles, this
written as \begin{equation}
\bbS^{+} = E \oplus (E\otimes K^{-1})
\end{equation}
with it understood that
\(w_{X}\) acts as multiplication by \(i|w_{X}|\) on the left
most summand (namely, \(E\)).  (This
splitting is the analog of the splitting in (\ref{split-S})). A section, \(\psi\), of \(\bbS\) is written with
respect to this splitting over \(U_{C}\) as \[\psi  =|w_{X}|^{1/2}
(\alpha, \beta).\]
  Meanwhile, \(\bbA\) is written as
\(A_{K} + 2A\) with \(A\) being a connection on \(E\).  Granted this
notation, write \(\hat{A}\) using the formula in (\ref{eq:(A.7)}) with it understood that
the covariant derivatives of \(\alpha\) that appear have non-zero
pairing with the vector field \(\frac{\partial}{\partial s}\).
This connection is flat where \(|\alpha |^2>\frac{9}{16}\) and
\(\alpha /|\alpha |\) is \(\hat{A}\)-covariantly constant. Meanwhile, \(\hat{A}\)
is equal to \(A\) where \(|\alpha |^2\leq \frac{7}{16}\). The formulas
for the curvature of \(\hat{A}\) and the \(\hat{A}\)-covariant
derivative of \(\alpha \) is given in (\ref{(3.9)v2}) with it understood that \(F_A\) and \(\nabla_A\alpha \) now have components that have non-zero pairing with \(\frac{\partial}{\partial s}\).
 
With a look ahead at the upcoming propositions, note that the
integral of \(iF_{\hat{A}}\) over \(C\) is proved to be well defined when \((A,
\psi)\) is an instanton solution to (\ref{eq:(A.14)}).  This is proved using
integration by parts to express the integral of \(iF_{\hat{A}}\) as
the difference between integrals of the  $i \mathbb{R}$-valued 1-form \(\hat{A}
- A_{E}\) over respective \(s \gg 1\) and \(s \ll -1\) slices of \(C\). 
 
The first proposition below concerns the integral of 
\(iF_{\hat{A}}\) on \(C\) when \(X\), its metric, and the 2-forms
\(w_{X}\) and \(\grw_{\mu}\) define the product cobordism.
 \begin{prop}\label{prop:A.2}   
Assume that \(X\), the metric, and \(w_{X}\) can be used to define a
product cobordism once \(\mu\) is chosen.   Assume
in addition that \(Y_{Z}\)  has a closed 1-form, \(\upsilon_{\diamond}\),  such that
\(\upsilon_{\diamond} \wedge w \geq 0\), whose restriction to \(U_{\gamma}\)  is \(dt\)  and whose
restriction to $\mathcal{H}_{0}$  is \(\textsc{h} \, du\) 
with \(\textsc{h}\)  being a strictly positive function of \(u\).  
 Given \(\c \geq 1\), there exists $\kappa> \pi $  with the following significance:  Fix
\(r \geq \kappa  \) and \(\mu \in \Omega  \) with either
$\mathcal{P}$-norm bounded by \(e^{-r^2}\) or with $\mathcal{P}$-norm bounded by  \(1\)
 but vanishing on $\bbR \times(\mathcal{H}_0\cup  U_{\gamma})$.   Let
\(\grc_{-}\)  and \(\grc_{+}\)  denote
solutions to the \((r, \mu)\)-version of (\ref{eq:(A.4)}) on
\(Y_{Z}\)  with \(\gra(\grc_{-})  -
\gra(\grc_{+})  \leq   r^{2-1/\c}\).   Suppose that \(\grd = (\bbA, \psi)\)  is an
instanton solution to the corresponding version of (\ref{eq:(A.14)}) on \(X\) 
with \(s \to -\infty  \) limit
\(\grc_{-}\)  and \(s \to \infty \) limit \(\grc_{+}\).  Then  \(i \int _CF_{\hat{A}}\geq 0\).
\end{prop}

Proposition \ref{prop:A.2} is a special case of the next proposition which concerns
the integral of \(iF_{\hat{A}}\) on \(C\) when the relevant data does
not necessarily define the product cobordism.  
 
\begin{prop}\label{prop:A.3}   
Assume that \(X\)  and
\(w_{X}\)  obey the conditions in Sections \ref{sec:Ac)}, 
and that the metric on \(X\)  obeys (\ref{eq:(A.12,15a)}) and (\ref{eq:(A.15b)}).
Then there exists $\kappa  > \pi $ such that given any \(\c \geq\kappa \), there exists \(\kappa_{\c}\) 
with the following property:  Fix \(r \geq\kappa_{c}\)  and assume that the metric and
\(w_{X}\)  are \((\c, \r = \textup{r})\)-compatible data. Fix \(\mu_{-}\)  and
\(\mu_{+}\)   from the respective \(Y_{-}\)-  and \(Y_{+}\)-versions
of \(\Omega\)  with either $\mathcal{P}$-norm less than \(e^{-r^2}\)
or with $\mathcal{P}$-norm less than 1 but vanishing
on the respective \(Y_{-}\)-  and \(Y_{+}\)-versions of $\mathcal{H}_0 \cup 
U_{\gamma}$.  Let \(\grc_{-}\) and \(\grc_{+}\)  denote solutions to
the \((r, \mu_{-})\)-version of (\ref{eq:(A.4)}) on \(Y_{-}\)  and \((r, \mu_{+})\)-version of (\ref{eq:(A.4)}) on
\(Y_{+}\)  with \(\gra(\grc_{-})-\gra(\grc_{+}) \leq   r^{2-1/\c}\).
If \(\grd = (\bbA, \psi)\) is an instanton
solution to (\ref{eq:(A.14)}) with \(s \to -\infty  \) limit
\(\grc_{-}\)  and \(s \to \infty\)
 limit \(\grc_{+}\), then \(i\int_C F_{\hat{A}}\geq 0\).
\end{prop}

Proposition \ref{prop:A.3} is proved in Section \ref{sec:Cb)}.

\subsection{The bound for \(\gra(\grc_-)-\gra (\grc_+)\) in
Proposition \ref{prop:A.3}}\label{sec:Af)}

Proposition \ref{prop:A.3} concerns only those instanton solutions to
(\ref{eq:(A.14)}) that obey the added constraint \(\gra(\grc_{-})-\gra(\grc_{+}) \leq   r^{2-1/\c}\).
The two propositions that are stated momentarily are used to
guarantee that this constraint is met in the cases that are relevant to
the body of this paper.  What follows sets the stage for the first
proposition.

\begin{defn}\label{def:c-tight}
Fix \(\c > 1\).  The metric on \(Y_{Z}\) and the
2-form \(w\) are said to define {\em \(\c\)-tight} data when there exists a
positive, \(\c\)-dependent constant  with the following
significance:  Use the metric, the 2-form \(w\), a choice of \(r\) greater
than this constant and a chosen 1-form from \(\Omega\) with
$\mathcal{P}$-norm less than 1 to define (\ref{eq:(A.4)}).  If \(\grc\) is a solution,
then \(|\gra^{\grf}(\grc)| < 
r^{2-1/\c}\). 
\end{defn}
\begin{prop}\label{prop:A.4}
Let \(Y_{Z}\) denote a compact, oriented Riemannian 3-manifold with a chosen
Riemannian metric and a $\Spin^c$-structure with
non-torsion first Chern class.  Let \(w\) denote a harmonic
2-form on \(Y_{Z}\)  whose de Rham class is this first
Chern class.  Assume that \(w\)  has non-degenerate zeros on any
component of \(Y_{Z}\)  where it is not identically
zero.  Then the metric and \(w\)  define a \(\c\)-tight data
set  if \(\c\) is sufficiently large. 
\end{prop}
This proposition is proved in Section \ref{sec:Bh)}.

This notion of being \(\c\)-tight is used in the second of the promised
propositions.  To set the stage for this one, suppose that \(X\) is a
cobordism of the sort that is described in Section \ref{sec:Ac)}.  Fix a metric
on \(X\) and the auxiliary data as described in (\ref{eq:(A.12,15a)}),
(\ref{eq:(A.13b)}), (\ref{eq:(A.13c)}), 
and let \(\grd = (A, \psi)\) denote an instanton solution to a given $r
> \pi $ version of (\ref{eq:(A.14)}).  Use \(\grc_{-}\) and
\(\grc_{+}\) to denote the respective \(s \to-\infty \) and \(s \to \infty \) limits of \(\grd\).  Associated
to \(\grd\) is a certain first order, elliptic differential operator, this
being the operator that is depicted in (IV.1.21) when \(X\) is the product
cobordism.  The operator in the general case is written using slightly
different notation in (2.61) of \cite{T3}.  This operator has a natural
Fredholm incarnation when the respective \(Y_{-}\) version
of \(\grf_{s}\) is constant on a neighborhood of \(\grc_{-}\) and the \(Y_{+}\) version is constant on
a neighborhood of \(\grc_{+}\).  Use \(\imath_{\grd}\)
to denote the corresponding Fredholm index.  By way of a relevant
example, \(\imath_{\grd}\) is equal to \(\grf_{s}(\grc_{+}) -\grf_{s}(\grc_{-})\) when \(X\) and the associated data
define the product cobordism.  Section \ref{sec:Cg)} associates an integer,
\(\imath_{\grd+}\), to \(\grd\) which is defined
without preconditions on \(\grc_{-}\) and \(\grc_{+}\).
 The latter is equal to the maximum of \(\imath_{\grd}\) and
0 in the case when \(\imath_{\grd}\) can be defined.

\begin{prop}\label{prop:A.5}   
Assume  that \(X\)  obeys the conditions in Sections \ref{sec:2.2} and
\ref{sec:Ac)}, 
that the metric on \(X\) obeys (\ref{eq:(A.12,15a)}) and (\ref{eq:(A.15b)})  for a given
\(L> 100\), and that \(w_X\) obeys the conditions in
(\ref{eq:(A.13a)}) and (\ref{eq:(A.13b)})
for a given \(L_*\geq L+4\).  Then there exists $\kappa  >\pi $ such
that for any given \(\c \geq \kappa \), there exists
\(\kappa_{\c}\)  with the following significance:
 Suppose that the respective pairs of metric and version of \(w\) 
on \(Y_{-}\)  and  \(Y_{+}\) define \(\c\)-tight data.  Fix \(r > \kappa_{\c}\)  and fix
\(\mu_{-}\)  and \(\mu_{+}\)  from the respective \(Y_{-}\)  and
\(Y_{+}\)  versions of \(\Omega\) with $\mathcal{P}$-norm less than 1 so as to define
(\ref{eq:(A.14)}) on \(X\) .  Let \(\grd\)  denote an instanton solution to
these equations with \(\imath_{\grd+} \leq \c
\).  Use \(\grc_{-}\)  and \(\grc_{+}\)  to denote the respective \(s \to-\infty  \) and \(s \to \infty  \) limits of
\(\grd\).  Then \(\gra(\grc_{-})  - \gra(\grc_{+})< r^{2-1/\c} \).
\end{prop}
 
This proposition is proved in Section \ref{sec:Cg)}.

\subsection{The cases when \(Y_Z\) is from \(\{M\sqcup(S^1\times S^2), Y\}\), 
  \(\{Y_k\}_{k=0, \ldots, \G}\), or \(\{Y_k\sqcup(S^1\times
  S^2)\}_{k=0, \ldots,\G-1}\)}\label{sec:Ag)}

In what follows, the notation \(Y_Z\) stands, in addition to the
manifold itself, also implicitly for its associated metric and two-form \(w\)
from Part 3 of Section \ref{sec:Aa)}.

 The body of this article is concerned with \(2\G + 3\) specific versions of \(Y_Z\), these
being as follows: The first manifold of interest is \(M\) and \(S^1
\times S^2\) and the second is the
manifold \(Y\) form Section II.1. The next \(\G + 1\) manifolds are labeled as \(\{Y_k\}_{k=0,...\G}\) with a
given \(k \in  \{0, ..., \G\}\) version being the manifold that is obtained from \(M\) by
attaching the handle \(\mathcal{H}_0\) as directed in Part 2 of Section II.1a and attaching \(k\) of the
handles from the set \(\{\mathcal{H}_\grp\}_{\grp\in \Lambda}\) as directed in Part 1 of Section II.1a. Note in this regard
that \(Y\) and \(Y_{\G}\) are the same manifold, endowed with
different metric and 2-form \(w\). Their main difference is the
behavior of \(w\) over the attached handles \(\mathcal{H}_\grp\)-- for
\(Y\) it
approximates certain standard contact form
(cf. (\ref{eq:(D.47)}) below), while for \(Y_{\G}\) it is harmonic
(cf. Proposition \ref{prop:A.6}).
 The last \(\G\) manifolds of interest are the disjoint unions of the various \(k \in  \{0, ..., \G -1\}\)
versions of \(Y_k\) and \(S^1\times S^2\).

  \paragraph{Part 1:}  Let \(Y_{Z}\) denote the disjoint union
of \(M\) and \(S^{1} \times S^{2}\).   To
see about the constraints in Sections \ref{sec:Aa)}, 
take \(Z\) to be the disjoint union of \(M\) and \(S^{3}\).  The handle
$\mathcal{H}_{0}$ is attached to \(S^{3}\) so as to obtain \(S^{1} \times S^{2}\).
 Write \(S^{1}\) as $\bbR/(2\pi \bbZ)$ and
let \(t\) denote the corresponding affine coordinate.  Use the spherical
coordinates \((\theta, \phi)\) for \(S^{2}\).  The loop
\(\gamma\) is the \(\theta = 0\) circle in \(S^{1}\times S^{2}\).  

 To see about \(w\) and the metric, consider first their appearance on
\(S^{1} \times S^{2}\).  Take the 2-form
\(w\) on \(S^{1} \times S^{2}\) to be  
 \(\sin\theta\, d\theta\, d\phi \) and the metric to be \(\textsc{h}\,
 dt^{2} + d\theta^{2} +\sin^{2}\theta d\phi^{2}\) with
\(\textsc{h}\) denoting a positive constant.   If the first Chern class
of \(\det\,(\bbS|_{M})\) is torsion, take \(w = 0\)
on \(M\) and take any smooth metric.  If the first Chern class of
\(\det\,(\bbS|_{M})\) is not torsion, 
take a metric on \(M\) such that the
associated harmonic 2-form with de Rham cohomology class that of \(c_{1}(\det\,(\bbS|_{M}))\) has non-
degenerate zeros. Take \(w\) in this case to be this same harmonic 2-form. By way of a
parenthetical remark, a sufficiently generic metric on \(M\) will have
this property. See, for example, \cite{Ho} for a proof that such is
the case.

 The data just described obeys the conditions in Section \ref{sec:Aa)}.
 Use Proposition \ref{prop:A.4} to see that this data is also \(\c\)-tight for a
suitably large version of \(\c\).

 \paragraph{Part 2:}  Let \(Y_{Z}\) denote the manifold \(Y\) that is
described in Section II.1.  Suffices it to say for now that \(Y\) is
obtained from \(Y_{0}\) via a surgery that attaches some
positive number of 1-handles to the \(M_{\delta}\) part of
\(Y_{0}\).  This number is denoted by \(\textsc{g}\).
 
The 2-form \(w\) is described in Section II.1e.  See also Part 3 of Section
IV.1a.  Let \(b_{1}\) denote the first Betti number of \(M\).
 Part 2 of  Section II.1d  describes a set of \(b_{1} + 1\)
closed integral curves of the kernel of \(w\) that have geometric
intersection number 1 with each \(u\) = constant 2-sphere in
$\mathcal{H}_{0}$.  One of these curves intersects
$\mathcal{H}_{0}$ as the \(\theta = 0\) arc.  This is the
curve \(\gamma^{(z_0)}\) in the notation from Part 2 of Section III.1a.  Use the latter for
$\gamma$.  It follows from what is said in (II.1.5) and Part 2 of
Section II.1d that the $\gamma$ has a tubular neighborhood with
coordinates as described in Section \ref{sec:Aa)} such that the 2-form \(w\) has the
desired appearance.  Section II.1e and (IV.1.5) describe a closed
1-form on \(Y\) that can be used to satisfy the requirements in the third
bullet in (\ref{eq:(A.3a)}).  This 1-form is denoted by \(\upsilon_{\diamond}\).  

 A set of Riemannian metrics on \(Y\) are described in Part 5 of Section
IV.1a that have the desired form on $\mathcal{H}_{0}$.
 Although not stated explicitly, a metric of the sort that is
described in Part 5 of Section IV.1a can be chosen so that it has the
desired behavior on some small radius tubular neighborhood of
$\gamma$.   Note that the set of metrics under consideration are
obtained from the choice of an almost complex structure on the kernel
of a 1-form \(\hat{a}\) given in (IV.1.6).  These almost complex structures are
taken from the set $\mathcal{J}_{ech}$ that is described
in Theorem II.A.1 and Section III.1c.  None of the conclusions in \cite{KLT2}-\cite{KLT4} are
compromised if the almost complex structure from
$\mathcal{J}_{ech}$ is chosen near $\gamma$  so that
the metric obeys the constraints in (\ref{eq:(A.3b)}).  To be sure, the chosen
almost complex structure must have certain genericity properties to
invoke the propositions and theorems in these papers.  These
genericity results are used to preclude the existence of certain
pseudoholomorphic subvarieties in $\mathbb{R}\times Y$ .  An
almost complex structure giving a metric near $\gamma$ that obeys
(\ref{eq:(A.3b)}) is not generic.  Even so, the subvarieties that must be excluded
can be excluded using a suitably almost generic almost complex
structure from the subset described in $\mathcal{J}_{ech}$
that give a metric that is described by (\ref{eq:(A.3b)}) near $\gamma$.  What
follows is the key observation that is used to prove this: The curves
to be excluded have image via the projection from $\mathbb{R}\times Y$
 that intersects the complement of small radius
neighborhoods of $\gamma$.  A detailed argument for the existence of
the desired almost complex structures from
$\mathcal{J}_{ech}$  amounts to a relatively
straightforward application of the Sard-Smale theorem along the lines
used in the proof of Theorem 4.1 in \cite{HT2}. 

It follows from Lemma IV.2.5 and Proposition IV.2.7 that the metric just
described together with \(w\) define a \(\c\)-tight data set on \(Y\) for a suitably large
choice for \(\c\).  
 
\paragraph{Part 3:}  
This part of the subsection considers the case when \(Y_Z\) is some \(k \in  \{0,
  ..., \G\}\) version of \(Y_k\). As noted previously, the manifold \(Y_k\) is
obtained from \(M\) by attaching the handle \(\mathcal{H}_0\) in the manner that is
described in Part 2 of Section II.1a and attaching \(k\) of the handles
from the set \(\{\mathcal{H}_\grp\}_{\grp\in \Lambda}\) as described in Part 1 of Section II.1a. Part 3
in Section II.1a defines a subset \(M_\delta  \subset M\) and the constructions of both
\(Y\) and \(Y_k\) identify \(M_\delta \) as a subset of both. The
curve \(\gamma^{(z_0)}\) that was introduced above in Part 2 sits in the latter part of \(Y\) and
so it can be viewed using this identification as a curve in
\(Y_k\).  Use this \(Y_k\) incarnation of \(\gamma^{(z_0)}\) for the curve $\gamma$.  

 The proposition that follows says what is needed with regards to the
2-form \(w\) and the metric to use on \(Y_k\).
 
\begin{prop}\label{prop:A.6}
Fix \(k\in\{0, \cdots, \G\}\).
There exists a nonempty set of
Riemannian metrics on \(Y_k\) with the following two properties:  Let \(w\)  denote
the metric's associated harmonic 2-form with de Rham cohomology
class that of \(c_{1}(\det\,(\bbS))\).   Then
\(w\)  has non-degenerate zeros.  Moreover, the metric and \(w\) obey the
conditions in Section \ref{sec:Aa)}. 
\end{prop}

This proposition is proved in Section \ref{sec:Db)}.  

The set of metrics in 
Proposition \ref{prop:A.6} is denoted by \(\Met\) in what follows.  Take the metric on
\(Y_k\) from this set and take \(w\) to be the associated
harmonic 2-form with de Rham cohomology class that of
\(c_{1}(\det\, (\bbS))\).  Proposition \ref{prop:A.4} asserts that
the resulting data set is \(\c\)-tight for a suitably large choice of \(\c\).

\paragraph{Part 4:} This part of the subsection discusses the case when
\(Y_Z\)  is the disjoint union of some \(k\in 0, \cdots, \G-1\) version of \(Y_k\)   and 
\(S^{1} \times S^{2} \). The metric on \(Y_{0}\)   comes from the
Proposition \ref{prop:A.6}'s set \(\Met\) , and the 2-form \(w\)
on \(Y_k\)  is the
corresponding harmonic 2-form with de Rham cohomology class that of
 \(c_{1} (\det\,( \mathbb{S} ))\).
Any smooth metric can be chosen for a given
\(S^{1} \times S^{2}\) component. The class \(c_{1} (\det\,( \mathbb{S} ))\)
is taken equal to zero on each \(S^{1} \times S^{2}\)  component and
this understood, the 2-form  \(w\)  is
identically zero on each such component. 

 What is said in Proposition \ref{prop:A.4} implies that the
resulting data set is  \(\c\)-tight for a
suitably large choice of  \(\c\).

\subsection{Cobordisms with \(Y_+\) and \(Y_-\) either \(Y\),
  \(M\sqcup (S^1\times S^2)\), from \(\{Y_k\}_k\), 
or \(\{Y_k \sqcup (S^1\times S^2)\}_k\) }
\label{sec:Ah)}
 
The first proposition concerns the product cobordisms when \(Y_Z\) is one of the manifolds
from the set \(Y\), \(M\sqcup (S^1 \times S^2)\), \(\{Y_k \}_{k=0, \ldots, \G}\) or \(\{Y_k \sqcup (S^1 \times  S^2)\}_{k=0,
  \ldots, \G-1}\). 
The subsequent propositions concern certain cobordisms of the sort described in
Section \ref{sec:Ac)} with \(Y_+\) and \(Y_-\) as follows:
\begin{itemize}
\item One is \(Y\) and the other is \(Y_{\G}\).
\item One is \(Y_k\) and the other is \(Y_{k-1} \sqcup (S^1 \times
  S^2)\) for some \(k\in \{1,...,\G \}\).
\item One is \(Y_0\) and the other is \(M\sqcup (S^1 \times S^2)\).
\end{itemize}
These propositions assume implicitly that
the metric and version of \(w\) on these manifolds are those supplied by
the relevant part of Section \ref{sec:Ag)}.  In particular, the metric
and \(w\) on \(M\sqcup (S^{1} \times S^{2})\) is
described by Part 1 of Section \ref{sec:Ag)}, and this data on \(Y\) is described in
Part 2 of Section \ref{sec:Ag)}.  Meanwhile, the metric on the
relevant \(k\in \{0, \ldots, \G\}\) version of \(Y_k\) is
from the set \(\Met\) and \(w\) is the associated harmonic 2-form with de Rham
cohomology class that of \(c_{1}(\det\, (\bbS))\). 
 
\begin{prop}\label{prop: A.7}   
Let \(Y_{Z}\) denote either \(M\sqcup  (S^{1} \times S^{2})\) or \(Y\)
or some \(k\in \{0, \ldots, \G\}\) version of \(Y_k\) with the 2-form \(w\)  and metric as
described in the preceding paragraph.  Given \(\imath \geq 0\),
there exists $\kappa > \pi $  with the following significance:  Fix
any \(r > \kappa  \) and
a 1-form \(\mu \in \Omega  \) with either $\mathcal{P}$-norm less than \(e^{-r^2}\)
or $\mathcal{P}$-norm less than 1 but vanishing on
\(\mathcal{H}_0\cup U_\gamma \).    Use this data with the metric and
\(w\) to define the product cobordism \(X = \bbR \times
Y_{Z}\) as prescribed in Section 2.2.  Suppose that \(\grc_{-}\) and \(\grc_{+}\)  are solutions to the \((r,
\mu)\)-version of (\ref{eq:(A.4)}) on \(Y_{Z}\)  with \(|\grf_{s}(\grc_{+}) -\grf_{s}(\grc_{-})| \leq \imath\);
and suppose that \(\grd\)  is an instanton solution to (\ref{eq:(A.14)}) on
\(X\)  with \(s \to -\infty  \) limit equal to \(\grc_{-}\)  and \(s
\to \infty \) limit equal to \(\grc_{+}\).  Then \(\textsc{x}(\grc_{+}) \geq\textsc{x}(\grc_{-})\).    
\end{prop}
 
\pf This follows directly from
Propositions \ref{prop:A.2} and \ref{prop:A.4} given what is said in Section \ref{sec:Ag)} about \(w\) and
the metric. \epf
 

The next proposition describes cobordisms between \(Y_{0}\)
and \(M\sqcup (S^{1} \times S^{2})\) of the sort that obey the conditions in Section \ref{sec:Ac)}.  
 
\begin{prop}\label{prop:A.8}  
Take the metric on \(M\sqcup (S^{1}\times S^{2})\) and harmonic 2-form \(w\) to be as described in Part 1 of Section \ref{sec:Ag)}. The metric on \(M\sqcup (S^{1}\times S^{2})\) determines a
corresponding set of metrics in the \(Y_0\) version of \(\Met\). Choose a metric
from this set and take \(w\) on \(Y_0\) to be the associated harmonic 2-form
with de Rham class \(c_{1}(\det\, (\bbS))\).
Denote one of \(Y_{0}\)  or \(M\sqcup (S^{1}\times S^{2})\) by
\(Y_{-}\)  and the other by \(Y_{+}\).  There exists a cobordism that obeys the
conditions in Section \ref{sec:Ac)} and the conditions in the list below.  This
list uses \(X\)  to denote the cobordism manifold.
\begin{itemize}
\item  The function \(s\) on \(X\) has exactly one critical
point.  This critical point has index 3 when \(Y_{-} =
Y_{0}\)  and index 1 when \(Y_{-} = M \sqcup (S^{1} \times S^{2})\).   
\item  There is a metric on \(X\)  with an associated self-dual
2-form that are \((\c, \r)\)-compatible if \(L\) and \(\c\) are sufficiently large and if  $\r >\pi
$.
\end{itemize}
\end{prop}
 
This proposition is proved in Section \ref{sec:Dd)}.

 The next proposition uses \(C\) to denote the cylinder in Proposition
\ref{prop:A.8}'s cobordism that is described by the first bullet
of (\ref{eq:(A.9b)}).  The proposition also reintroduces the notation in (\ref{(3.9)v2}).

\begin{prop}\label{prop:A.9}  
Take \(w\)  and the metric on \(Y_{0}\) and on \(M\sqcup (S^{1}\times S^{2})\)  to be as described in
Proposition \ref{prop:A.8}.   Denote one of \(Y_{0}\) 
or \(M\sqcup (S^{1} \times S^{2})\)  by \(Y_{-}\)  and the
other by \(Y_{+}\).  Take the cobordism space
\(X\), the metric on \(X\),  and the associated self-dual 
2-form
\(w_{X}\)  to be as described by Proposition \ref{prop:A.8}. Given \(k \geq 0\), there exists 
$\kappa >\pi $  with the following significance:  Fix \(r >\kappa  \) and 1-forms \(\mu_{-}\) 
and \(\mu_{+}\)  from the \(Y_{-}\)-  and \(Y_{+}\)-versions
of \(\Omega\) with either $\mathcal{P}$-norm less than \(e^{-r^2}\)
or with $\mathcal{P}$-norm less than 1 but vanishing
on the \(Y_{-}\)-  and \(Y_{+}\)-versions of $\mathcal{H}_0\cup
U_{\gamma}$.  Let \(\grd  = (A, \psi)\) denote an instanton solution to the resulting version of (\ref{eq:(A.14)})
with \(\imath_{\grd+} \leq k\). Then \(i \int_CF_{\hat{A}}   \geq 0\).
\end{prop}
\pf The proposition follows
directly from Propositions \ref{prop:A.3}, \ref{prop:A.5} and \ref{prop:A.8} given what is said in
Section \ref{sec:Ag)} about the respective \(Y_{0}\) and \(M  \sqcup
(S^{1} \times  S^{2})\) metrics and
versions of \(w\). \epf
 
The next set of propositions are analogs of Propositions \ref{prop:A.8} and \ref{prop:A.9} in the case when
one of \(Y_-\) and \(Y_+\) is some \(k\in \{1,\cdots, \G\}\) version
of \(Y_k\) and the other is \(Y_{k-1}  \sqcup(S^{1} \times  S^{2})\), or 
when one is \(Y\) and the other is \(Y_{\G}\).
The propositions that follow assume that \(c_1(\det\, \bbS)\) on each
\(k\in \{0, \ldots, \G\}\) version of \(Y_k\)
vanishes on the cross-sectional spheres in any \(\grp\in\Lambda\)
version of \(\mathcal{H}_\grp\) and that it has pairing
2 with the cross-sectional spheres in \(\mathcal{H}_0\). This class is
also assumed to be zero on the \(S^1\times S^2\) component of any \(k\in \{0, \ldots, \G\}\) version of \(Y_{k-1}  \sqcup(S^{1} \times  S^{2})\). Meanwhile, its 
restriction to the \(H_2(M;\bbZ)\) summand from the associated Mayer-Vietoris sequence for the various \(k\in \{0, \ldots, \G\}\) versions of \(H_2(Y_k;\bbZ)\) is assumed to be independent of \(k\).

\begin{prop}\label{prop:A.10a}  
There exists, for each \(k\in \{0, \ldots, \G\}\),  a subset to be denoted by
\(\op{Met}\,( Y_k)\) in the \(Y_k\) version of \(\Met\) with the
following significance: Let \(\op{Met}(Y_0)\) denote the subset from
Proposition \ref{prop:A.8}. For each \(k\in \{1, \ldots, \G\}\), 
take a metric from an open subset of \(\op{Met}(Y_{k-1})\)  and a
metric on \(S^1\times S^2\) to define a
metric on  \(Y_{k-1}\sqcup(S^1\times S^2)\). Take \(w\) on
\(Y_{k-1}\sqcup(S^1\times S^2)\) to be the associated harmonic 2-form
with de Rham class \(c_1(\det\, \bbS)\). The chosen metric determines
a corresponding subset of metrics \(\op{Met}\,( Y_k)\subset \Met\). Take a metric from the latter subset and take \(w\) to be the
associated harmonic 2-form with de Rham class
\(c_1(\det(\bbS))\). Take \(Y_-\) to be one of \(Y_k\) and
\(Y_{k-1}\sqcup(S^1\times S^2)\), and  take \(Y_+\) to be the other. There exists a cobordism that obeys the 
conditions in Section \ref{sec:Ac)} 
and the conditions listed below. This list uses \(X\) to denote the cobordism manifold.\begin{itemize}
\item The function \(s\) on \(X\) has precisely 1 critical point. This critical point has index 3
when \(Y_+\) has the \(S^1\times S^2\) component and it has index 1
when \(Y_-\) has the \(S^1\times S^2\) component.
\item There is a metric on \(X\) with an associated self-dual 2-form that are \((\c, \r)\)-compatible if
\(L\), \(\c\) and \(\r > \pi \).
\end{itemize}
\end{prop}

This proposition is proved in Section \ref{sec:De)}.  

The next proposition considers the case when one of \(Y_-\) and \(Y_+\) is \(Y\) and the other is \(Y_{\G}\).
 \begin{prop}\label{prop:A.10b}  
Take \(w\) and the metric on \(Y\) to be as described in the opening paragraphs of this subsection. Take the metric on \(Y_{\G}\) from a certain non-empty subset of
\(\op{Met}(Y_{\G})\) and take \(w\) on \(Y_{\G}\) to be the associated harmonic
2-form with de Rham class that of \(c_1(\det (\bbS))\). Take \(Y_-\) to be one of \(Y\) and \(Y_{\G}\) and take \(Y_+\) to be the other. There exists a
cobordism that obeys the conditions in Section \ref{sec:Aa)} and the
conditions listed below. This list uses \(X\) to denote the cobordism
manifold.
\begin{itemize}
\item The function \(s\) on \(X\) has no critical points.
\item There is a metric on \(X\) with an associated self-dual 2-form that are \((\c, \r)\)-compatible if \(L\), \(\c\) and \(\r > \pi \).
\end{itemize}
\end{prop}

The proof of Proposition \ref{prop:A.10b} is in Section \ref{sec:Df)}.

 The upcoming proposition uses \(C\) to denote the cylinder in Proposition
\ref{prop:A.10a} and \ref{prop:A.10b}'s cobordism that is described by the first bullet
of (\ref{eq:(A.9b)}).  Notation from (\ref{(3.9)v2}) 
is also used.
 
\begin{prop}\label{prop: A.11}  
Let \(X\) denote one of the cobordism manifolds that are described in
Propositions \ref{prop:A.10a} and \ref{prop:A.10b} with
\(c_1(\det(\bbS))\) and the 2-form and metrics on \(Y_-\), \(Y_+\) and \(X\) as
described therein.
 Given \(\imath\geq 0\), there exists $\kappa >\pi $  with the following property:  Fix \(r >
\kappa  \) and 1-forms \(\mu_{-}\) and \(\mu_{+}\)  from the \(Y_{-}\)-  and \(Y_{+}\)-versions
of \(\Omega\)  with either $\mathcal{P}$-norm less than \(e^{-r^2}\)
or with $\mathcal{P}$-norm less than 1 but vanishing on the \(Y_{-}\)  and \(Y_{+}\) 
versions of $\mathcal{H}_{0} \cup U_{\gamma}$.  Let \(\grd  = (A, \psi)\)
 denote an instanton solution to the resulting version of (\ref{eq:(A.14)})
with \(\imath_{\grd+} \leq \imath\). Then \(i   \int_CF_{\hat{A}}\geq 0\).
\end{prop} 

\pf The proposition follows
directly from Propositions \ref{prop:A.3}, \ref{prop:A.5},
\ref{prop:A.10a}  and \ref{prop:A.10b} given what is said in
Section \ref{sec:Ag)}. 
\epf

\subsection{Filtered Floer homologies and
  filtration-preserving chain maps}\label{sec:3.8}

This subsection is divided into two parts. In the first part, we associate to each triple \((Y_Z, w, \gamma)\)
described in Section \ref{sec:Aa)} a system of filtered monopole Floer
homologies \(\HM^\circ (Y_Z, rw; \Lambda_\gamma)\), for \(\circ=-, \infty,
+, \wedge\) and \(r>\pi\), in the manner described in Section
\ref{sec:3.1}. Recall the constraint on the cohomology class \([w]\)
from Part 3 of Section \ref{sec:Aa)}. Together with the first bullet of
(\ref{eq:(A.3a)}), this implies that \(\CM_* (Y_Z, rw;
\Lambda_\gamma)\) is associated with a negative-monotone, non-balanced
perturbation. For reasons that will become clear momentarily, we use
\(\CM^\circ (Y_Z, \langle w\rangle; \Lambda_\gamma)\) to denote \(\CM_* (Y_Z, rw;
\Lambda_\gamma)\) for \(r\gg \pi\) and similarly for its
homology. (The notation \(\langle w\rangle \) stands for the ray
\(\bbR^+[w]\subset H^2(Y_Z; \bbR)\). This includes, as special cases, the triple
\((\underline{M}, \underline{w}, \underline{\gamma})\) in \cite{L}
(\(\underline{M}\) is denoted by \(Y_0\) in this
article), and the triple \((Y, w, \gamma^{(z_0)})\) in Section II.1a.  

In the second part, a filtration-preserving chain
map from \(\CM^- (Y_-, \langle w_-\rangle;\Lambda_{\gamma_-})\) to \(\CM^- (Y_+, \langle w_+\rangle;
\Lambda_{\gamma_+})\) is associated to each triple \((X, \varpi_X, C)\)
described in Section \ref{sec:Ac)}. To explain the notation, \(X\) is a cobordism from
the 3-manifold \(Y_-\) to \(Y_+\), \(\varpi_X\) is a self-dual 2-form
on \(X\) related to \(w_-\) and \(w_+\) as prescribed by
(\ref{eq:(A.13a)}). What is denoted by \(C\) signifies an embedded surface in
\(X\), with ends \(\gamma_-\subset Y_-\) and \(\gamma_+\subset Y_+\);
see the second bullet of (\ref{eq:(A.9b)}). 

\paragraph{Part 1:}
To accomplish this task, begin by introducing the (partially
defined) integral 1-cocycle on \(\mathcal{B}^\sigma(Y_Z)\) defining
\(\Lambda_\gamma\). This local system associates each \(\grc\in
\mathcal{Z}_{w, r}\) the group algebra \(\bbK[\bbZ]=\bbK[U,
U^{-1}]\). To each \(\grd\in \mathcal{M}_1(\grc_-, \grc_+)\) it
associates \(U^{n(\grd)}\in \op{End}(\bbK[U,U^{-1}])\), where \(n(\grd)=\smx
(\grc_+)-\smx (\grc_-)\). Here, \(\smx\) is the ``modified holonomy
function'' given in (\ref{eq:(A.8)}). Lemma \ref{lem:A.1} asserts that
\(n(\grd)\in \bbZ\) for \(\grc_-, \grc_+\in \mathcal{Z}_{w,
  r}\). Following the recipe in Section \ref{sec:3.1}, we then set
\((\CM^\infty, \partial^\infty)\) to be the monopole Floer complex with twisted
coefficients:
\[\begin{split}
\CM^\infty& =\bbK[U,
U^{-1}]\,( \mathcal{Z}_{w,r}),\\
 \partial^\infty\grc_- &=\sum_{\grc_+\in \mathcal{Z}_{w,r}}\sum_{\grd\in
  \mathcal{M}_1(\grc_-, \grc_+)/\bbR}\op{sign}(\grd)\,
U^{n(\grd)}\grc_+\, \, \quad \text{for
    \(\grc_-\in \mathcal{Z}_{w,r}\).}
\end{split}
\]
The monotonicity condition guarantees that the sum here is finite.
The \(\op{sign}(\grd)\in \{\pm 1\}\) in the preceding expression
is assigned according to the orientation convention laid out in \cite{KM}.

One may regard \(\CM^\infty\) as a chain complex over
\(\bbK\), generated by
\(\hat{\mathcal{Z}}_{w,r}=\mathcal{Z}_{w,r}\times \bbZ\). The
generating set \(\hat{\mathcal{Z}}_{w,r}\) lies in
\(\tilde{\mathcal{B}}^\sigma=\mathcal{C}^\sigma/\mathcal{G}_\gamma\),
a \(\bbZ\)-covering of \(\mathcal{B}^\sigma\). Here,
\(\mathcal{G}_\gamma\subset\mathcal{G}\) consists of smooth maps
\(u\co Y_Z\to S^1\), with \(\deg (u|_\gamma)=0\). Multiplication by
\(U^n\) then corresponds to a deck-transformation on this
\(\bbZ\)-covering, and the condition on \(c_1(\det \bbS)\) set forth
in Part 3 of Section \ref{sec:Aa)} then implies that \(\deg
U=-2\). 
The grading set of \(\hat{\mathcal{Z}}_{w,r}\) is an affine
space over \(\bbZ/c_Z\bbZ\), where \(c_Z\in 2\bbZ\) is the gcd of the
values of \(c_1(\det\bbS)\) on \(H_2(Z;\bbZ)\) according the splitting
(\ref{eq:(A.2)}). 
\begin{remarks}\label{rmk:3.15}
(a) Here, we use the same notation \(U\) for the map on monopole Floer
complexes described in Part 2 of Section \ref{sec:2.4} and deck
transformation here. This is because for the kind of \(Y_Z\)
considered in this article, they turn out identical by the arguments
for the last bullet of Proposition IV.7.6. 

(b) The way grading on a monopole Floer complex with local coefficients
follow the way they are defined in some literature, e.g. what is
called a Floer-Novikov complex
\cite{L1}. The book \cite{KM} does not seem to contain an explicit
discussion on the grading of Floer complex with local coefficients.
\end{remarks}  
Suppose that \((Y_Z, w)\) define \(\c\)-tight data for \(\c>1\)
(cf. Definition \ref{def:c-tight}). 
Take \(X\) to be the product cobordism \(\bbR\times Y_Z\), \(w_X=w+ds\wedge*w\), and
\(C=\bbR\times\gamma\subset X\). Let \(\grd\in \mathcal{M}(\grc_-,
\grc_+)\) be as in Proposition \ref{prop:A.3}. In this case,
\(i\int_CF_{\hat{A}}=2\pi(\smx(\grc_+)-\smx (\grc_-))\) and 
Proposition \ref{prop:A.3} asserts that one has \(n (\grd)\geq 0\). Thus,
\[\CM^-=\bbK(\mathcal{Z}\times \bbZ^{\geq 0})\subset\CM^\infty\] 
is a
subcomplex of \(\bbK[\bbZ^{\geq 0}]=\bbK[U]\)-modules. One may then introduce \[\CM^+=\CM^\infty/\CM^-,\quad
\widehat{CM}=\CM^-/U\CM^-.\] 
The resulting short exact sequences 
\begin{equation}\label{eq:fund-short}
\begin{split}0& \to
\CM^-\to \CM^\infty\to \CM^+\to 0\quad \text{and}\\
\quad 0& \to U\CM^-\to \CM^-\to
\widehat{\CM}\to 0
\end{split}
\end{equation}
induce the {\em fundamental exact sequences} on
the homologies. As the \(\bigwedge^*H_1(Y_Z)/\op{Tors}\)-action on the
monopole Floer complexes commute with \(U\), the exact sequences above
preserve the \(\mathbf{A}_\dag\)-module structure. 

In Section \ref{sec:Ag)}, the assumption that \((Y_Z, w)\) is
\(\c\)-tight is verified for the particular manifolds listed 
therein.  In particular, 
\BTitem\label{eq:nota-H}
\item when \(Y_Z=Y\) and its assciated \(\grs\), \(w\), \(\gamma\) and metric are
  as in Part 2 of Section~\ref{sec:Ag)}, \[\HM^\circ(Y, \langle w\rangle;
  \Lambda_\gamma)=\op{H}^\circ (Y)=\op{H}_{SW}^\circ\]
in the notation of \cite{KLT1, KLT4}.
\item when \(Y_Z=Y_k\), \(k\in \{0, \cdots, \G\}\) and its
  assciated \(\grs\), \(w\), \(\gamma\) and metric are in Part 3 of Section
\ref{sec:Ag)},  the corresponding \(\HM^\circ(Y_Z, \langle w\rangle;
  \Lambda_\gamma)\) are instrumental in the proof for Theorem
  \ref{thm:main}. Recalling that \(Y_0\) and its
  assciated \(\grs\), \(w\), \(\gamma\) and metric are respectively
  what was denoted by 
  \(\underline{M}\) \(\ul{\grs}\), \(\ul{w}\), \(\ul{\gamma}\) in
  \cite{L}, we observe that
\[
\HM^\circ(Y_0, \langle w\rangle; \Lambda_\gamma)=\op{HMT}^\circ
\]
introduced in \cite{L}.
\ETitem

Note that \(\CM^\circ(Y_Z, \langle w\rangle; \Lambda_\gamma)\), \(\HM^\circ(Y_Z, \langle w\rangle;
\Lambda_\gamma)\) introduced above implicitly depend on \(r\) and
\((\grT, \grS)\). According to the convention set forth in Section
\ref{sec:notation}, this is permissible if there are chain homotopies
between the monopole Floer complexes associated with different parameters
preserving the \(\mathbf{A}_\dag\)-module structure. This is
justified by combining the arguments proving Proposition IV.1.4 with
what is said in the upcoming Part 2.

\paragraph{Part 2:}
We now consider chain maps induced by (non-product) cobordisms \(X\) described in
Section \ref{sec:Ac)}. To begin, we introduce an \(X\)-morphism from
\(\Lambda_{\gamma_-}\) to \(\Lambda_{\gamma_+}\). (See Definition 23.3.1
in \cite{KM} for ``\(X\)-morphism''). This is done in a way similar to the definition of
\(\Gamma_C\) in Equation (23.8) in \cite{KM}. In \cite{KM}, a
``cobordism'' from \(Y_-\) to \(Y_+\) refers to a compact 4-manifold
with boundary \(Y_+\sqcup (-Y_-)\). This corresponds to the compact
part of our \(X\), denoted by \(X_c=s^{-1}([-L_{tor}, L_{tor}])\). The
surface \(C\cap X_c\) plays the role of the singular 2-chain \(\nu\)
in (23.8) of \cite{KM}. It has boundary \(\gamma_+-\gamma_-\), with
\(\gamma_+\simeq\gamma\simeq \gamma_-\). Given
\(\grc_-\in\mathcal{Z}_{w_-, r}(Y_-)\) and
\(\grc_+\in\mathcal{Z}_{w_+, r}(Y_+)\), let \(\grd\) denote an element
in \(\mathcal{B}^\sigma(X)\) with \(s\to -\infty\) limit \(\grc_-\)
and \(s\to \infty\) limit \(\grc_+\). Then \(\Gamma_C\) is an
isomorphism from \(\Gamma_{\gamma_-}(\grc_-)\simeq \bbR\) to \(\Gamma_{\gamma_+}(\grc_+)\simeq \bbR\)
given by multiplication by \(e^{\frac{i}{2\pi}\int_CF_{\bbA}}\). 
The analog of \(\Gamma_C\) in our setting, denoted \(\Lambda_C\) below, is
given by an homomorphism from
\(\Lambda_{\gamma_-}(\grc_-)\simeq\bbK[U, U^{-1}]\) to
\(\Lambda_{\gamma_+}(\grc_+)\simeq\bbK[U, U^{-1}]\) for each
pair \(\grc_-\), \(\grc_+\). This is
given by multiplication with \(U^{n(\grd)}\), where \begin{equation}\label{eq:n_cob}
n(\grd)=\frac{i}{2\pi}\int_CF_{\hat{A}}=\smx_{\gamma_+}(\grc_+)-\smx_{\gamma_-}(\grc_-),
\end{equation}
the right most equality being a consequence of the Stokes' theorem. This is again an integer according to Lemma
\ref{lem:A.1}. With \(\Lambda_C\) in place, given a \(k\)-cochain \(\mathpzc{u}\in
\op{C}^k(\mathcal{B}^\sigma (X);\bbK)\) in the notation of Section \ref{sec:2.4}, we
define the map \[\begin{split}m^\infty[u](X, \langle w_X\rangle;
  \Lambda_C)\co & \CM^\infty (Y_-)=\bbK[U,
U^{-1}]\, (\mathcal{Z}_{w_-,r})\\
&\quad \qquad  \to  \CM^\infty (Y_+)=\bbK[U,U^{-1}]\, (\mathcal{Z}_{w_+,r})
\end{split}
\] by the following rule:
\[\begin{split}
\mathcal{Z}_{w_-,r}(Y_-)\ni \grc_- &\mapsto \sum_{\grc_+\in \mathcal{Z}_{w_+,r}(Y_+)}\sum_{i}
 \big\langle\mathpzc{u}, \mathcal{M}_k(X, \grc_-, \grc_+)\big\rangle\,
 U^{n(\grd_i)}\, \grc_+,
\end{split}
\]
where \(i\) runs through each connected component of
\(\mathcal{M}_k(X, \grc_-, \grc_+)\); and for every \(i\), \(\grd_i\) is an element in
the corresponding connected component.
In order for the sum on the right hand side to be well defined, we
assume that \(H^2(X, Y_-)=0\) and \(w_X\) satisfies (\ref{eq:w_X}).

To proceed, suppose \((Y_-, w_-)\), \((Y_+, w_+)\) are \(\c\)-tight
and consider \(C(X, \langle w_X\rangle; \Lambda_C)|_{ \CM^-(Y_-)}\). Suppose
furthermore that \((X, w_X)\) satisfies the conditions in Propositions
\ref{prop:A.3} and \ref{prop:A.5}. By these propositions, the integers
\(n(\grd_i)\) in (\ref{eq:n_cob}) are non-negative, implying that the
image of \(C(X, \langle w_X\rangle; \Lambda_C)|_{ \CM^-(Y_-)}\) under
\(m^\infty\) lies in
\(\CM^-(Y_+)\). Use \[m^-[u](X, \langle w_X\rangle; \Lambda_C)\co \CM^-(Y_-)\to
\CM^-(Y_+)\] 
to denote this map. It is straightforward to see that both
\(m^\infty\) and \(m^-\) are chain maps, given that \(\CM^\infty\) is
a variant of monopole Floer complexes, and the non-negativity of the
integers \(n(\grd)\) appearing in the formulae for \(\partial^\infty\)
and \(m^\infty\). These then induce homomorphisms between the
respective homologies: 
\[
\HM_*(X, \langle w_X\rangle; \Lambda_C)\co  \HM^\circ(Y_-, \langle w_-\rangle; \Lambda_{\gamma_-})\to  \HM^\circ(Y_+, \langle w_+\rangle; \Lambda_{\gamma_+})
\]
for \(\circ=-, \infty\). Like those in Part 1, these maps preserve the \(\mathbf{A}_\dag
\)-module structure.

\section{Some homological algebra}\label{sec:4}
\setcounter{equation}{0}

As mentioned in Section 1,  the purpose of this section is to review
the algebraic background for the upcoming Proposition
\ref{prop:KM-ES}. The latter is used to relate the fomula for monopole
Floer homology of a connected sum, given in Proposition \ref{prop:conn-f} below, in
terms of the monopole Floer homology with balanced perturbation that
appears in Theorem \ref{thm:main} and Theorem \ref{cor:main}.
This computation turns out to be 
a simplest manifestation of the so-called ``Koszul duality'', well-known in certain
circles. For a sampling of literature on this subject, see
e.g. \cite{C}, \cite{Ko}, \cite{GKM}. The variant most relevant to
this article is discussed in \cite{GKM}, which relates the ordinary
chain complex of an \(S^1\)-space, equipped with an
\(H_*(S^1)\)-module structure capturing the \(S^1\)-action, with the
\(S^1\)-equivariant chain complexes of the same space, which are
naturally endowed with \(H^*(BS^1)\)-module structures. We need
however only a small portion of the full machinery in \cite{GKM}.  
Thus, in this section we give a self-contained though elementary 
exposition of the relevant part of this story, tailored
to our needs.

\subsection{Terminology and conventions}

By a {\em modules over \(H^*(BS^1)\)} we mean a chain complex with a module
  structure over \(\bbK[u]\), where \(u\) acts as a chain map of degree
\(-2\). The prime examples of such modules in this article are the
monopole Floer complexes. In parallel, a {\em module over
  \(H_*(S^1)\)} stands for a chain complex with a module structure over
\(\bbK[y]\), where \(y\) acts as a degree 1 chain map. An example that
appears later is the chain complex to compute the monopole Floer
homology of a connected sum, see (\ref{eq:conn-sum}) in Proposition \ref{prop:conn-1}. Meanwhile, 
a graded homology module \(H_*\) will be viewed as 
 a chain complex with zero differentials. 
 We use capital letters \(U\), \(Y\) to
denote the chain maps corresponding to the action of \(u\), \(y\).

\begin{defn}
A {\em morphism} from one module over \(H^*(BS^1)\) to another is a
\(\bbK\)-chain map which commutes with \(U\)-actions. 
Morphisms between
\(H_*(S^1)\)-modules are defined similarly, with \(Y\) replacing \(U\). 
We shall also often encounter a weaker notion: a {\em
  p-morphism} 
 between two \(H^*(BS^1)\)-modules is a \(\bbK\)-chain map
which commutes with \(U\)-actions up to \(\bbK\)-chain homotopy. 
\end{defn}

\subsection{From \(H^*(BS^1)\)-modules to \(H_*(S^1)\)-modules}\label{sec:S_U}
Given a module \((C, \partial _C)\) over
\(H^*(BS^1)\), we define the module \(S_U(C)\) over
\(H_*(S^1)=\bbK[y]\) as follows:
\begin{equation}\label{def:S_U}
\big(S_U(C), S_U(\partial_C)\big)=(C\otimes \bbK[y],
\, \partial_C\otimes\jmath +U
\otimes y),
\end{equation}
where the homomorphism \(\jmath\co \bbK[y]\to \bbK[y]\) is defined by
\[\jmath \,( a+by)=a-by\quad \text{for \(a, b\in \bbK\)},\]
and the \(y\)-action is simply the multiplication \(1\otimes
y\). (\(\jmath\) was denoted by \(\sigma\) in \cite{L}. Cp. Equation (5.1) therein). 

To see that \(S_U(C)\) is indeed a chain complex, note that the condition
\(S_U(\partial_C)^2=0\) is equivalent to the pair of identities
\(\partial_C^2=0\), and \([\partial_C, U]=0\).

\begin{lemma}\label{lem:u-y} A
p-morphism \(\Phi\) between two \(H^*(BS^1)\)-modules \((C_{(1)}, \partial_{(1)})\),
\((C_{(2)}, \partial_{(2)})\) induces 
an \(H_*(S^1)\)-module morphism \(S_U(\Phi)\) between
\(S_{U_{(1)}}(C_{(1)})\) and \(S_{U_{(2)}}(C_{(2)})\), where \(U_{(1)}\),
\(U_{(2)}\) denote the \(u\)-action on \(C_{(1)}\), \(C_{(2)}\)
respectively. Furthermore, 
\begin{itemize}
\item \(S_U(\Phi)\) is injective if \(\Phi\) is
injective, and it is surjective if \(\Phi\) is surjective;
\item Let  \(\Phi'\) be another p-morphism of \(H^*(BS^1)\)-modules
  from \((C_{(1)}, \partial_{(1)})\) to
  \((C_{(2)}, \partial_{(2)})\). Then \(\Phi+\Phi'\) is a
  p-morphism as well, and 
\begin{equation*}
S_U(\Phi+\Phi')=S_U(\Phi)+S_U(\Phi').
\end{equation*}
\item Let  \(\Psi\) be a p-morphism of \(H^*(BS^1)\)-modules
  from \((C_{(2)}, \partial_{(2)})\) to \((C_{(3)}, \partial_{(3)})\). Then \(\Psi\circ\Phi\) is a
  p-morphism as well, and 
\begin{equation}\label{eq:S-compose}
S_U(\Psi\circ\Phi)=S_U(\Psi)\circ S_U(\Phi).
\end{equation}
\end{itemize}
\end{lemma}
\pf As a p-morphism, \(\Phi\)
satisfies both \begin{equation}\label{eq:S_UPhi}\begin{split}
\Phi\, \partial_{(1)}-(-1)^{\deg (\Phi)}\partial_{(2)}\, \Phi& =0 \qquad \text{and}\\
\Phi \, U_{(1)}-U_{(2)}\, \Phi & =K _\Phi\,\partial_{(1)}+(-1)^{\deg (\Phi)}\partial_{(2)}\, K_\Phi,
\end{split}
\end{equation}
for a \(\bbK\)-linear  homomorphism \(K_\Phi\).  
This is equivalent to the identity \begin{equation}\label{eq:S-chainmap}
S_U(\Phi)\,
S_U(\partial_{(1)})-(-1)^{\deg (\Phi)}S_U(\partial_{(2)})\,
S_U(\Phi)=0,
\end{equation}
where \(S_U(\Phi)\co C_{(1)}\otimes \bbK[y]\to  C_{(2)}\otimes \bbK[y]\) is
defined as
\begin{equation}\label{def:p-morphism}
S_U(\Phi)= \Phi\otimes \, \jmath^{\deg\, (\Phi)}+K_\Phi\otimes y. \quad 
\end{equation}
This verifies that \(S_U(\Phi)\) is a chain map. Moreover, since the
\(y\)-action on \(S_U(C_{(1)})\), \(S_U(C_{(2)})\) is 
multiplication by \(1\otimes y\), it is immediate that \(S_U(\Phi)\)
commutes with the \(y\)-actions on both sides. The claim that \(S_U\)
preserves injectivity and surjectivity can be checked directly from
the definition of \(S_U(\Phi)\).

Since the construction of \(S_U(\Phi)\) is linear, the second item in
the statement of the Lemma is obvious. 

To verify the third bullet about the composition of p-morphisms, let 
\[
S_U(\Psi)= \Psi\otimes \, \jmath^{\deg\, (\Psi)}+K_\Psi\otimes y.
\]
Then (\ref{eq:S-compose}) is straightforward to verify, given that
\[
K_{\Psi\circ\Phi}=K_\Psi\circ\Phi+(-1)^{\deg (\Psi)}\Psi\circ K_\Phi.
\]
The fact that \(\Psi\circ\Phi\) is a p-morphism follows directly
from (\ref{eq:S-chainmap}) and its analog for \(\Psi\).

The first bullet may be directly verified after writing out the
definition of \(S_U(\Phi)\) explicitly. More is said in the proof of
Lemma \ref{lem:y-u} below. 
\epf

\begin{remarks}\label{rem:S_U-htpy}
Given a p-morphism \(\Phi\), 
the \(H_*(S^1)\)-morphism \(S_U(\Phi)\) given in
(\ref{def:p-morphism}) apparently depends on the choice of the degree 
\((\deg \Phi  +1)\) map \(K_\Phi
\). By (\ref{eq:S_UPhi}), two different choices of \(K_\Phi
\), say \(K_\Phi\) and \(K'_\Phi\), differ by a chain map:
\[
[\partial_C, K_\Phi -K'_\Phi]:=(K _\Phi- K'_\Phi)
\,\partial_{(1)}+(-1)^{\deg (\Phi)}\partial_{(2)}\, (K_\Phi -K'_\Phi) =0. 
\]
We say that \(K_\Phi\) and \(K'_\Phi\) are {\em homotopic}
if there exists a degree \(\deg (\Phi)\) linear map \(Z_\Phi \co C_{(1)}\to C_{(2)}\) such that 
\[
K_\Phi -K'_\Phi=[\partial_C, Z_\Phi]:=\partial_{(2)} Z_\Phi -Z_\Phi (-1)^{\deg (\Phi)}\partial_{(1)}. 
\]
Let \(S_U(\Phi)\), \(S_U(\Phi)'\) respectively denote the versions of
\(S_U(\Phi)\) defined using \(K_\Phi\) and
\(K'_\Phi\). They are chain-homotopic when \(K_\Phi\) and \(K'_\Phi\) are
  homotopic:
\begin{equation}\label{diff:S_UPhi}
S_U(\Phi)-S_U(\Phi)'=([\partial_C, Z_\Phi])\otimes y=[S_U(\partial_C),
-Z_\Phi \otimes y]. 
\end{equation}
(Keep in mind that in our notation, \([\cdot, \cdot]\) stand for a
commutator in a graded sense.) Thus, for a given \(\Phi \), the
homology \(H_*(S_U(\Phi);\bbK)\) depends only on the (relative)
homotopy class of \(K_\Phi \). 
\end{remarks}

\begin{defn}
Two \(H^*(BS^1)\)-modules \((C_{(1)}, \partial_{(1)})\), \((C_{(2)}, \partial_{(2)})\) are said to be {\em
  p-homotopic} if there exist p-morphisms \(\Phi\co C_{(1)}\to
C_{(2)}\), \(\Psi\co C_{(2)}\to C_{(1)}\), and \(H_1\co C_{(1)}\to C_{(1)}\), \(H_2\co  C_{(2)}\to
C_{(2)}\), such that 
\[
\Psi\circ\Phi -\op{Id}_{(1)}=[\partial_{(1)}, H_1], \quad 
\Phi\circ\Psi-\op{Id}_{(2)}=[\partial_{(2)}, H_2].
\]
They are said to be {\em homotopic} if \(\Phi\), \(\Psi\), \(H_1\),
\(H_2\) are  morphisms. The notion of two \(H_*(S^1)\)-modules being {\em
  homotopic} is defined similarly.
\end{defn}
\begin{lemma}\label{lem:S-htpy}
Suppose two \(H^*(BS^1)\)-modules \((C_{(1)}, \partial_{(1)})\), \((C_{(2)}, \partial_{(2)})\) are  {\em
  p-homotopic} via p-morphisms \(\Phi\co C_{(1)}\to
C_{(2)}\), \(\Psi\co C_{(2)}\to C_{(1)}\) as above. Then the
\(H_*(S^1)\)-modules \(S_U(C_{(1)}\)), \(S_U(C_{(2)})\) are homotopic
via the maps \(S_U(\Phi)\) and \(S_U(\Psi)\).
\end{lemma}
\pf By assumption, there exist \(H_1\), \(H_2\) such that \(\Phi\), \(\Psi\) satisfy:
\[
\Psi\circ\Phi -\op{Id}_{(1)}=[\partial_{(1)}, H_1], \quad 
\Phi\circ\Psi-\op{Id}_{(2)}=[\partial_{(2)}, H_2].
\]
We need to verify the identities:
\begin{equation*}\begin{split}
S_U(\Psi)\circ S_U(\Phi)-\op{Id}_{(1)}& =[S_U(\partial_{(1)}),
S_U(H_1)]; \\ 
S_U(\Phi)\circ S_U(\Psi)-\op{Id}_{(2)}& =[S_U(\partial_{(2)}),
S_U(H_2)];\\
[S_U(\Phi), Y]&=0;\\ 
 [S_U(\Psi), Y]& =0.
\end{split}
\end{equation*}
It suffices to verify the first and the third identities, since the
second and the fourth are entirely parallel. 

To verify the first identity, use (\ref{eq:S-compose}) and the fact
that \(\Psi\circ\Phi -\op{Id}_{(1)}=[\partial_{(1)}, H_1]\) to reduce it to 
\[
S_U(\op{Id}_{(1)})=\op{Id}. 
\]
This holds by taking \(\Psi=\op{Id}_{(1)}\) and \(K_\Psi=0\) in
(\ref{def:p-morphism}).

To verify the third identity, simply plug in the definition of
\(S_U(\Phi)\) and \(Y=I\otimes y\). \epf 

\subsection{From \(H_*(S^1)\)-modules to \(H^*(BS^1)\)-modules}
First, introduce the \(\bbK[u]\)-modules:
\begin{equation}\label{defn:V}
\begin{split}
V^- &:=u\bbK[u],\, \\
V^\infty  &:= \bbK[u, u^{-1}], \, \\
V^+  &:= \bbK[u, u^{-1}]/u \bbK[u],\\
V^{\wedge} &:= \bbK[u]/u \bbK[u].\\
\end{split}
\end{equation} 
These modules by definition fit into the short exact sequences:
\begin{gather}
0\to V^-\stackrel{i_V}{\hookrightarrow }V^\infty\to V^+\to 0, \label{eq:V-seq1}\\
0\to V^-\stackrel{u}{\to } V^-\to V^{\wedge }\to 0.\label{eq:V-seq2}
\end{gather}
 We shall frequently view these four modules as a system, and write
them collectively as \(V^\circ\). The same convention applies to the
various system of modules we construct out of these four below. 

\begin{defn}\label{def:EC} (\cite{J, L})
Given a module \((C, \partial _C)\) over
\(H_*(S^1)\), we define the following system of modules over \(H^*(BS^1)=\bbK[u]\):
\begin{equation}\label{def:E}
\big(E^\circ_Y(C), E_Y(\partial_C)\big):=(C\otimes V^\circ, \, \partial_C\otimes 1+Y\otimes u)
\quad \text{for \(\circ=-, \infty, +, \wedge\)}, 
\end{equation}
where the \(u\)-action is the multiplication \(1\otimes u\).
\end{defn}
The fact that \(E_Y(\partial_C)^2=0\) again follows directly from the
definition of \(H_*(S^1)\)-modules: \(\partial_C^2=0\),
\([Y, \partial_C]=0\), \(Y^2=0\). 
By taking tensor product of
(\ref{eq:V-seq1}), (\ref{eq:V-seq2}) with \(C\),
one has the following corresponding short exact sequences of \(\bbK[u]\)-modules:
\begin{gather}
0\to E_Y^-(C)\stackrel{\op{Id}\otimes i_V}{\longrightarrow}E^\infty_Y(C)\to E^+_Y(C)\to 0, \label{eq:CV-seq1}\\
0\to E_Y^-(C)\stackrel{\op{Id}\otimes u}{\longrightarrow} E_Y^-(C)\to E_Y^{\wedge }(C)\to 0.\label{eq:CV-seq2}
\end{gather}
It is also straightforward to verify
that the maps in the above exact sequences commute with
\(E_Y(\partial_C)\), and therefore
they induce the following long exact sequences of \(H^*(BS^1)\)-modules associated to
\((C, \partial _C)\): 

\begin{gather}\label{eq:E-sq1}
\cdots\to H_*(E_Y^-(C))\stackrel{i_{V_*}}{\longrightarrow}H_*(E_Y^\infty(C))\to H_*(E_Y^+(C))\stackrel{\delta_{V*}}{\longrightarrow}
H_{*-1}(E_Y^-(C))\to \cdots\\
\cdots\to H_*(E_Y^-(C))\stackrel{u}{\to }H_*(E_Y^-(C))\to H_*(E_Y^\wedge(C))\to
H_{*-1}(E_Y^-(C))\to \cdots
\label{eq:E-sq2}
\end{gather}

We call (\ref{eq:E-sq1}), (\ref{eq:E-sq2}) the (resp. first and
second) {\em fundamental exact
  sequences} 
for the \(H_*(S^1)\)-module \(C\).
For convenience of later reference, we denote the short exact
sequences of \(H^*(BS^1)\)-modules (\ref{eq:CV-seq1}),
(\ref{eq:CV-seq2}) by \(\bbE_Y(C)\), \(\textbf{E}_Y(C)\) respectively.
Correspondingly, the long exact sequences (\ref{eq:E-sq1}),
(\ref{eq:E-sq2}) are denoted by \(H(\bbE_Y(C))\), \(H(\textbf{E}_Y(C))\).
It is straightforward to verify the assertion in the following Lemma
and so we leave it to the reader to check that
\begin{lemma}\label{lem:y-u} 
A morphism \(\phi\) between two
\(H_*(S^1)\)-modules \((C_{(1)}, \partial_{(1)})\),
\((C_{(2)}, \partial_{(2)})\) induces a system of
\(H^*(BS^1)\)-module morphisms
\[\begin{split}
E^\circ(\phi) \co E^\circ_{Y_{(1)}}(C_{(1)}) & \to
E^\circ_{Y_{(2)}}(C_{(2)})\\
\phi & \mapsto\phi\circ 1
\end{split}
\] 
for \(\circ=-, \infty, +, \wedge\), where
\(Y_{(1)}\), \(Y_{(2)}\) denote the \(y\)-actions on \((C_{(1)}, \partial_{(1)})\),
\((C_{(2)}, \partial_{(2)})\) respectively. Moreover:
\begin{itemize}
\item \(E^\circ(\phi)\) is injective if \(\phi\) is injective; it is
  surjective if \(\phi\) is surjective.
\item Let  \(\phi'\) be another morphism of \(H_*(S^1)\)-modules
  from \((C_{(1)}, \partial_{(1)})\) to
  \((C_{(2)}, \partial_{(2)})\). Then \(\phi+\phi'\) is an \(H_*(S^1)\)-morphism as well, and 
\begin{equation*}
E^\circ (\phi+\phi')=E^\circ(\phi)+E^\circ (\phi').
\end{equation*}
\item Let  \(\psi\) be another morphism of \(H_*(S^1)\)-modules
  from \((C_{(2)}, \partial_{(2)})\) to
  \((C_{(3)}, \partial_{(3)})\). Then \(\psi\circ\phi\) is an \(H_*(S^1)\)-morphism as well, and 
\begin{equation}\label{eq:E-compose}
E^\circ(\psi\circ\phi)=E^\circ(\psi)\circ E^\circ(\phi).
\end{equation}
\item The system of morphisms \(E^\circ(\phi)\) combine to define morphisms of short exact sequences of
\(H^*(BS^1)\)-modules \[
\bbE(\phi)\co \bbE_Y(C_{(1)})\to
\bbE_Y(C_{(2)})\quad \text{and} \quad \textbf{E}(\phi)\co \textbf{E}_Y(C_{(1)})\to
\textbf{E}_Y(C_{(2)}).\] 
Correspondingly, their induced
maps on homologies \(H_*(E_Y(\phi))\) combine to define morphisms of
long exact sequences of 
\(H^*(BS^1)\)-modules \[\begin{split}
& H(\bbE_Y(\phi))\co H(\bbE_Y(C_{(1)}))\to
H(\bbE_Y(C_{(2)}))\quad \text{and} \\
\quad & H(\textbf{E}_Y(\phi))\co H(\textbf{E}_Y(C_{(1)}))\to
H(\textbf{E}_Y(C_{(2)})).
\end{split}
\] 
\end{itemize}
\end{lemma}
\pf The proofs are straightforward; thus we shall say no more than making 
 the following remarks: Both \(E_Y^\circ\) and \(S_U\) preserve injectivity and surjectivity
due to the same reason, namely they can be written in polynomial form (in \(u\) and \(y\)
respectively, which defines a filtration), where their 0-th order term
takes the form of a tensor product of the
original morphism and an automorphism.
This in turn implies that both of them takes short exact sequences
to short exact sequences.
\epf

\begin{lemma}\label{lem:E-htpy}
Let \(C_{(1)}\), \(C_{(2)}\) denote homotopic
\(H_*(S^1)\)-modules. Then \(E^\circ_Y(C_{(1)})\),
\(E^\circ_Y(C_{(2)})\) are homotopic \(H^*(BS^1)\)-modules.
\end{lemma}
\pf By assumption, there exist morphisms \(\Phi\co C_{(1)}\to
C_{(2)}\), \(\Psi\co C_{(2)}\to
C_{(1)}\) and \(H_1\co C_{(1)}\to
C_{(1)}\), \(H_2\co C_{(2)}\to C_{(2)}\) such that \begin{equation}\label{S-htpy}
\Psi\circ\Phi -\op{Id}_{(1)}=[\partial_{(1)}, H_1], \quad 
\Phi\circ\Psi-\op{Id}_{(2)}=[\partial_{(2)}, H_2].
\end{equation}
Lemma \ref{lem:y-u} claims that \(E^\circ(\Phi)\co E^\circ_Y(C_{(1)})\to 
E^\circ_Y(C_{(2)})\), \(E^\circ(\Psi)\co E^\circ_Y(C_{(2)})\to 
E^\circ_Y(C_{(1)})\) are systems of morphisms. Meanwhile, the desired
identities are:
\begin{equation}\label{eq:E-htpy}
E^\circ(\Psi)\circ E^\circ(\Phi )-\op{Id}_{(1)}=[E_Y(\partial_{(1)}), E^\circ(H_1)], \quad 
E^\circ(\Phi)\circ E^\circ(\Psi)-\op{Id}_{(2)}=[E_Y(\partial_{(2)}), E^\circ(H_2)].
\end{equation}
We shall only verify the first identity, since the second is
similar. For this purpose, apply \(E^\circ\) to the first identity in
(\ref{S-htpy}), then apply Lemma \ref{lem:y-u} and subtract the first
line of (\ref{eq:E-htpy}) to the resulting identiy. This leads to
\[
E^\circ(\op{Id}_{(1)})-\op{Id}=[Y, H_1]\otimes u.
\]
This is true because of the definition of \(E^\circ\) and the fact
that \(H_1\) is a morphism. 
\epf

\subsection{Koszul duality} 

The functors \(S_U\) and \(E^-\) may be
viewed as inverses of each other in the following sense.
\begin{prop}\label{prop:koszul}
{\rm\bf (a)} \, \, Let \((C, \partial_C)\) be an \(H^*(BS^1)\)-module.
Then there is a system of isomorphisms of \(H^*(BS^1)\)-modules:
\begin{equation}\label{eq:ES}
H_*(E^\circ _YS_U(C))\simeq H_*(C\otimes _{\bbK[u]}V^\circ).
\end{equation}
Moreover, these isomorphisms have the
following naturality properties: 
\begin{itemize}
\item[(i)] they are 
natural with respect to p-morphisms of \(H^*(BS^1)\)-modules, and
\item[(ii)] they combine to define isomorphisms of long exact
  sequences of \(H^*(BS^1)\)-modules 
\[
H(\bbE_YS_U(C))\simeq H(C\otimes _{\bbK[u]}\bbV),\quad 
  \text{and}\quad H(\textbf{E}_YS_U(C))\simeq
  H(C\otimes_{\bbK[u]}\textbf{V}).\] 
Here, \(H(C\otimes _{\bbK[u]}\bbV)\) and \(H(C\otimes_{\bbK[u]}\textbf{V})\) respectively
  denote the long exact sequence induced by the short exact sequences
  of \(H^*(BS^1)\)-modules:
\[\begin{split}
0\to C\otimes _{\bbK[u]}V^-\to C\otimes _{\bbK[u]}V^\infty\to
C\otimes _{\bbK[u]}V^+\to 0\quad \text{and}\\ 0\to C\otimes
_{\bbK[u]}V^-\stackrel{1\otimes u}{\longrightarrow }C\otimes
_{\bbK[u]}V^-\to C\otimes _{\bbK[u]}V^\wedge\to 0.
\end{split}
\]
\end{itemize}

{\rm\bf (b)}\,\, Let \((C, \partial_C)\) be an
\(H_*(S^1)\)-module. Then there is an isomorphism of \(H_*(S^1)\)-modules:
\[
H_*(S_UE_Y^-(C))\simeq H_*(C).
\] 
\end{prop}
\pf 
{\em Part (a):} Written out explicitly,
\[\begin{split}
E^\circ_YS_U(C) &=C\otimes \bbK[y]\otimes  V^\circ, \\
E_YS_U(\partial _C) &=\partial_C\otimes  \jmath\otimes 1+U\otimes y\otimes
1+1\otimes y\otimes u.
\end{split}
\]
View this as a filtered complex by the total degree in the \(C\otimes
V^\circ\) factor.
Then the \(E_1\)-term of the associated spectral sequence is simply 
\begin{equation}\label{eq:E1}
C\otimes  \bbK\{y\}\otimes  V^\circ/\big((U\otimes y\otimes
1+1\otimes y\otimes u )(C\otimes  \bbK\{1\}\otimes  V^\circ)\big)\, \simeq  C\otimes_{\bbK[u]}V^\circ,
\end{equation}
with differential \(d_1\) given by \(-\partial_C\). This 
spectral sequence degenerates at \(E_2\), and we have 
\(H_*(E^\circ_YS_U(C))\simeq H_*(C\otimes _{\bbK[u]}V^\circ)\) as
claimed. As the \(u\)-action on \(E^\circ_YS_U(C)\) is \(1\otimes
1\otimes u\) and the \(u\)-action on \(C\) is \(U\), the quotient in
(\ref{eq:E1}) shows that the isomorphism preserves the \(\bbK[u]\)-module
structure. 
Property (ii) also follows immediately from this
computation.  

On the other hand, given a
p-morphism \(\Phi\) between two \(H^*(BS^1)\)-modules \((C_{(1)}, \partial_{(1)})\),
\((C_{(2)}, \partial_{(2)})\), by Lemmas \ref{lem:u-y} and
\ref{lem:y-u}, there is a corresponding system of morphisms of
\(H^*(BS^1)\)-modules \(E^\circ_YS_U(\Phi)\). The naturality property
(i) follows from the fact that these morphisms preserve the filtration.

\noindent{\em Part (b):} Written out explicitly, 
\[\begin{split}
S_UE_Y^-(C) &=C\otimes u\bbK[u]\otimes \bbK[y], \\
S_UE_Y^-(\partial_C) & =\partial_C\otimes 1\otimes \sigma+Y\otimes
u\otimes\sigma+1\otimes u\otimes y.
\end{split}
\]
Filtrate by the total degree of the factor \(C\otimes u\bbK[u]\) as in
the previous part.
Then the \(E_1\) term is \[C\otimes u\bbK[u]\otimes R\{y\}/\big((1\otimes
u\otimes y)(C\otimes u\bbK[u]\otimes R\{1\})\big)\simeq C,\] 
on which \(d_1\) acts as \(-\partial_C\). The spectral sequence again
degenerates at \(E_2\), yielding the claimed isomorphism
\(H_*(S_UE_Y^-(C))\simeq H_*(C).\) To see that the module structure
agree, note that a cycle in the \(E_1\) term given by an element
\(-z_1\in C\), \(\partial_Cz_1=0\) corresponds to a cycle in
\(S_UE_Y^-(C)\) of the form \(Z_0\otimes 1+z_1\otimes u\otimes y\),
where \(Z_0\in   C\otimes u\bbK[u]\) satisfies 
\[-((Y\otimes u)Z_0)\otimes 1+((1\otimes u)Z_0)\otimes
y-(Yz_1)\otimes u^2\otimes y=0.
\]
In other words, \(Z_0=-(Yz_1)\otimes u\), and the cycle in
\(S_UE_Y^-(C)\) has the form \(-(Yz_1)\otimes u\otimes 1+z_1\otimes
u\otimes y\). The \(y\)-action \(1\otimes 1\otimes y\) takes this element to \(-(Yz_1)\otimes
u\otimes y\), while the element corresponding to \(-Yz_1\in C\) in the
\(E_1\) term is \(-(Yz_1)\otimes
u\otimes y\) as well, since \(Y^2=0\).
\epf
\begin{remarks}\label{rmk:ES}
(a) Spelled out explicitly, (\ref{eq:ES}) says that 
\(H_*(E_Y^-S_U(C))\simeq H_*(C)\), and \(H_*(E_Y^\infty S_U(C))\) is the localization of
\(H_*(C)\) as a \(\bbK[u]\)-module. On the other hand, note that since
\(V^\wedge=\bbK[u]/u\bbK[u]\), \(E^\wedge S_U(\partial_C)\) reduces to
  \(S_U(\partial_C)\otimes 1\), and therefore \(H_*(E_Y^\wedge
  S_U(C))\simeq H_*(S_U(C))\).

(b) The constructions \(E_Y^-,
E_Y^\infty, E_Y^+\) above are directly copied from J. Jones's
formulation of the ``co-Borel'', ``Tate'', and Borel (the usual)
versions of equivariant homologies \cite{J}.  
It is proved in \cite{GKM} that \(S_U\), and \(E_Y\) induce isomorphisms
of derived categories.
\end{remarks}
\begin{remarks}\label{rem:period_cpx}
As stated, the general \(H^*(BS^1)\)-module \((C, \partial_C)\) in
the present section is assumed to be \(\bbZ\)-graded. For our application
however, results in this section are
typically applied to monopole chain complexes \(\mathring{C}\). These
are only  relatively graded, and the grading group is 
\(\bbZ\) only when \(c_1(\grs)\) is torsion, in other cases it is 
\(\bbZ/c_\grs\), where \(c_\grs\in 2\bbZ\). These chain complexes
\(\mathring{C}\) are also equipped with a canonical absolute
\(\bbZ/2\)-grading. (Cf. \cite{KM}). Nevertheless, we observe that
such a monopole chain complex \(\mathring{C}\) can be alternatively 
interpreted as a \(\bbZ\)-graded
chain complex \((C, \partial)\) with periodicity \(c_\grs\in 2 \bbZ\), namely, 
\((C_k, \partial_k)=(C_{k+c_\grs}, \partial_{k+c_\grs})\) \(\forall
k\in \bbZ\). (For prior
appearance of such interpretation, cf. e.g. \cite{L1}). 

To do this, let \(\mathring{C}\) be a certain monopole chain complex
associated to the \(\Spin^c\) manifold \((M, \grs)\), and let \(\mathcal{Z}\subset \mathcal{B}:=\mathcal{C}^\sigma(M,
\grs)/C^\infty(M, U(1))\) denote the generating set of
\(\mathring{C}\). Denote the coefficient ring of \(\mathring{C}\) by
\(\bbK\). Recall that 
\(H^1(\mathcal{B};\bbZ)\simeq H^2(M; \bbZ)\), and therefore the class \(c_1(\grs)\in
H^2(M; \bbZ)\) defines a \(\bbZ\)-covering \(\pi\co
\tilde{\mathcal{B}}\to \mathcal{B}\). Let
\(\tilde{\mathcal{Z}}:=\pi^{-1}\mathcal{Z}\), and consider
the chain complex 
\((C, \partial_C)\) with \[C:=\bbK(\tilde{\mathcal{Z}})\quad
\text{and}\quad 
\langle\tilde{\grc}_1,  \partial_C\tilde{\grc}_2\rangle:=\sum_{\tilde{\grd}}\op{sign} (\tilde{\grd})\]
for any pair \(\tilde{\grc}_1, \tilde{\grc}_2\in
\tilde{\mathcal{Z}}\). Regarding elements in
\(\mathcal{M}_1(\grc_1, \grc_2)/\bbR\) as paths in \(\mathcal{B}\)
from \(\grc_1\) to \(\grc_2\), \(\tilde{\grd}\) in the sum above
stands for any lift of
some \(\grd\in \mathcal{M}_1(\grc_1, \grc_2)/\bbR\) to \(\tilde{\mathcal
  B}\) ending in \(\tilde{\grc}_1\) and \(\tilde{\grc}_2\), where
\(\grc_*:=\pi\tilde{\grc}_*\) and \(\op{sign}
(\tilde{\grd}):=\op{sign} (\grd)\). Since the spectral flow on
\(\mathcal{B}\) is controlled by \(c_1(\grs)\), the relative grading of this complex
\((C, \partial_C)\), defined by spectral flow along \(\tilde{\grd}\),
is \(\bbZ\)-valued. 
Fix a \(\tilde{\grc}_0\in
\tilde{\mathcal{Z}}\). \(\grc_0=\pi\tilde{\grc}_0\) is either even or
odd according to the canonical absolute \(\bbZ/2\)-grading. Set \(\gr
(\tilde{\grc}_0)=0\) if \(\grc_0\) is even, and set \(\gr
(\tilde{\grc}_0)=1\) if \(\grc_0\) is odd. Together with the
relative \(\bbZ\)-grading \(\gr (\cdot, \cdot)\) on
\((C, \partial_C)\), we have an absolute \(\bbZ\)-grading by setting
 \[\gr (\tilde{\grc}):=\gr (\tilde{\grc}_0)+\gr
 (\tilde{\grc}_0,\tilde{\grc})\] for any \(\tilde{\grc}\in
 \tilde{\mathcal Z}\).  (If
\(\tilde{\mathcal{Z}}=\emptyset\), \(\hat{C}\) is trivial). With this
definition,  we then have 
\(C_k=\mathring{C}_{k'}\) for any pair \(k\in \bbZ, k'\in \bbZ/c_\grs\)
with \(k=k'\mod c_\grs\). Moreover, given \(\tilde{\grc}_1\in
\tilde{\mathcal{Z}}\) with \(\gr (\tilde{\grc}_1)=k\) and a
\(\grd\in \mathcal{M}_1(\grc_1, \grc_2)/\bbR\), there is a unique lift \(\tilde{\grd}\)
of \(\grd\) starting from \(\tilde{\grc}_1\), whose end point is a lift \(\tilde{\grc}_2\in \tilde{\mathcal Z}\) of \(\grc_2\).  \(\gr (\tilde{\grc}_2)=k-1\). Thus, \(\partial_C|_{C_k}\co C_k\to C_{k-1}\) is
identical to \(\mathring{\partial}|_{\hat{C}_{k'}}\co \hat{C}_{k'}\to
\hat{C}_{k'-1}\) for the same pair of \(k, k'\) as before. 
\end{remarks}

\section{Balanced Floer homologies from monotone Floer chain
  complexes}
\setcounter{equation}{0}

This section re-introduces the fourth flavor of monopole Floer
homology denoted by \(\op{HM}^{tot}\) in \cite{L},  now
renamed \(\widetilde{\HM}\) in deference to Donaldson's notation. 
(Cf. p.187 of \cite{D}) This definition is a natural by-product of a
re-interpretation of \(\mathring{\HM}_*(M, \grs, c_b)\) in terms
{\em purely of the \(\bbK[U]\)-module \(\hat{C}_*(M, \grs,
  c_b)\)} (Corollary 5.3 in \cite{L}, re-stated as Proposition \ref{prop:KM-ES}
below). This result enables us to appeal to the third author's
``\(SW=Gr\)'' program, which in our context was carried out in part IV
of this series \cite{KLT4}.  The latter constructed an isomorphism from
an appropriate variant of \(ech\) to a {\em negative monotone} version
of monopole Floer homology, which is in turn related to the balanced
version via the following theorem of Kronheimer and Mrowka:

\begin{thm}[\cite{KM} Theorem 31.5.1]\label{thm:KM}
Suppose \(c_1(\grs)\) is non-torsion. 
Let \(\widehat{C}_*(M, \grs, c_b)\) and \(C_*(M, \grs,
c_-)=\widehat{C}_*(M, \grs, c_-)\) respectively denote the Seiberg-Witten-Floer chain
complexes with balanced and negative monotone
perturbations. 
Then there is a chain homotopy equivalence 
from the former to the latter. In particular, 
\(\widehat{\HM}_*(M, \grs, c_b)\simeq\HM_* (M, \grs, c_-)\).
\end{thm}
To be more precise, the statement of Theorem 31.5.1 in \cite{KM} concerns only the
Floer homologies. However, the chain homotopy equivalence referred to above
was constructed in its proof.

\begin{remarks}
The variant of \(ech\) relevant in this series of papers is related to
the negative monotone version of monopole Floer homology, and
therefore to \(\widehat{\HM}_*\) by the preceding theorem of
\cite{KM}. This is because the stable Hamiltonian structure used to
define the relevant \(ech\) is associated to an {\em nonexact} closed
2-form. Note in contrast that the ordinary embedded contact homology
associated to a contact structure is related to \(\widecheck{HM}_*\)
instead, since the relevant 2-form in this case is exact. As such, it 
belongs to the positive monotone situation, and the companion theorem
to the one just cited states that \(\widecheck{HM}_*(M, \grs, c_b)\simeq\HM_* (M, \grs, c_+)\).
\end{remarks}

\subsection{Some properties of the maps \(i\), \(j\), \(p\)}

In this section, unless otherwise specified, \(\mathring{C}_*=\mathring{C}_*(M, c_b)\) denotes 
the monopole Floer chain
complex associated to an oriented Riemannian, \(\Spin^c\) 3-manifold
with a {\em balanced perturbation}. Similarly, let
\(\mathring{HM}_*=\mathring{HM}_*(M, c_b)\).

Recall from Proposition 22.2.1 in \cite{KM} that \(\widehat{HM}_*\),
\(\widecheck{HM}_*\), \(\overline{HM}_*\) are related by a 
long exact sequence
\begin{equation}\label{eq:SW-fund-sq}
\cdots \to \overline{HM}_{*}\stackrel{i_*}{\to }\widecheck{HM}_*\stackrel{j_*}{\to }
\widehat{HM}_*\stackrel{p_*}{\to }\overline{HM}_{*-1}\stackrel{i_*}{\to }\cdots,
\end{equation}
which we shall call the {\em fundamental exact sequence of monopole
  Floer homologies}. The maps \(i_*\), \(j_*\), \(p_*\) in the sequence above are respectively 
induced by maps:
\[
i\co \bar{C}\to \check{C}, \quad j\co \check{C}\to  \hat{C}, \quad p\co
\hat{C}\to \bar{C}, 
\]
which, written in block form with respect to the decomposition 
\begin{equation}\label{eq:C-block}
\bar{C}=C^s\oplus C^u, \quad \check{C}=C^o\oplus C^s, \quad
\hat{C}=C^o\oplus C^u
\end{equation}
are given by 
\begin{equation}\label{eq:ijk}
i=\left[\begin{array}{cc} 0 &-\partial^u_o\\1 & -\partial^u_s
  \end{array}
\right], 
\quad
j= \left[\begin{array}{cc} 1 & 0\\0 & -\bar{\partial}^s_u
  \end{array}
\right],\quad 
p=\left[\begin{array}{cc} \partial ^o_s & \partial^u_s\\0 & 1
  \end{array}
\right].
\end{equation}
It is shown in \cite{KM} that they are respectively chain
maps of degree \(0\), degree \(0\), degree \(-1\). 
\begin{lemma}
The maps \(i\), \(j\), \(p\) are p-morphisms of \(H^*(BS^1)\)-modules.
\end{lemma}
\pf
It is verified in \cite{KM} that \([\mathring{\partial},
\mathring{U}]=0\) for the to-, from, and bar-versions of monopole
Floer chain complexes. 
A straightforward though tedious computation using (\ref{eq:ijk})
shows that:
\begin{equation}\label{eq:U-i}\begin{split}
i\bar{U}-\check{U}i+K_i\bar{\partial}+\check{\partial}K_i=0, \\
j\check{U}-\hat{U}j+K_j\check{\partial}+\hat{\partial}K_j=0, \\
p\hat{U}-\bar{U}p-K_p\hat{\partial}+\bar{\partial}K_p=0,
\end{split}
\end{equation}
where \(K_i\), \(K_j\), \(K_p\), written in block form with respect to
the same decompositions (\ref{eq:C-block}), are:
\[
K_i=\left[\begin{array}{cc} 0 &-U^u_o\\0 & -U^u_s
  \end{array}
\right], 
\quad
K_j= \left[\begin{array}{cc} 0 & 0\\0 & -\bar{U}^s_u
  \end{array}
\right],\quad 
K_p=\left[\begin{array}{cc} U^o_s & U^u_s\\0 & 0
  \end{array}
\right].
\]\epf

As was explained in the proof of Lemma \ref{lem:u-y}, the 
identities (\ref{eq:U-i}) can be rewritten as 
\[\begin{split}
S_U(i)S_U(\bar{\partial})-S_U(\check{\partial}) S_U(i) &=0, \quad\\
S_U(j)S_U(\check{\partial})-S_U(\hat{\partial}) S_U(j) &=0, \quad\\
S_U(p)S_U(\hat{\partial})+S_U(\bar{\partial}) S_U(p) &=0,
\end{split}
\] 
where
\(S_U(i)\), \(S_U(j)\), \(S_U(p)\), when written in block form with respect
to the same decomposition (\ref{eq:C-block}), as matrices with
coefficients in \(\bbK[y]\), are given as follows:
\[
S_U(i)=\left[\begin{array}{cc} 0 &-n^u_o\\1 & -n^u_s
  \end{array}
\right], 
\quad
S_U(j)= \left[\begin{array}{cc} 1 & 0\\0 & -\bar{n}^s_u
  \end{array}
\right],\quad 
S_U(p)=\left[\begin{array}{cc} N ^o_s & N^u_s\\0 & 1
  \end{array}
\right],
\]
where 
\[\begin{split}\bar{n}^s_u=\bar{\partial}^s_u+\bar{U}^s_uy, & \quad
n^*_*=\partial^*_*+U^*_*y, \\
\bar{N}^s_u=\bar{\partial}^s_u\jmath+\bar{U}^s_uy,  &\quad
N^*_*=\partial^*_*\jmath+U^*_*y,  
\end{split}
\]
these being homomorphisms of \(\bbK[y]\)-modules 
for any pair of super- and subscripts \(*\) among \(u, s, o\).

\begin{lemma}\label{lem:induced-KMa}
The induced maps from \(i\), \(j\), \(p\) fit into the
following long exact sequences:
\begin{gather}\label{eq:induced-KM1}
\cdots \longrightarrow H_*(S_U(\bar{C}))\stackrel{S_U(i)_*}{\longrightarrow }H_*(S_U(\check{C}))\stackrel{S_U(j)_*}{\longrightarrow }
H_*(S_U(\hat{C}))\stackrel{S_U(p)_*}{\longrightarrow
}H_{*-1}(S_U(\bar{C}))\stackrel{S_U(i)_*}{\longrightarrow }\cdots;\\
\begin{split}
\cdots \longrightarrow H_*(E_Y^\circ S_U(\bar{C})) &\stackrel{E_Y^\circ S_U(i)_*}{\longrightarrow }H_*(E_Y^\circ S_U(\check{C}))\stackrel{E_Y^\circ S_U(j)_*}{\longrightarrow }
H_*(E_Y^\circ S_U(\hat{C}))\\
& \quad \stackrel{E_Y^\circ S_U(p)_*}{\longrightarrow
}H_{*-1}(E_Y^\circ S_U(\bar{C}))\stackrel{E_Y^\circ S_U(i)_*}{\longrightarrow }H_{*-1}(E_Y^\circ S_U(\check{C}))\cdots.\label{eq:induced-KM}
\end{split}
\end{gather}
The first sequence is a sequence of \(H_*(S^1)\)-modules, and
the second one is a sequence of \(H^*(BS^1)\)-modules.
\end{lemma}
\pf {\em Part (a).} The proof is based on a modification of the proof of Proposition 22.2.1 in
\cite{KM}. Recall from \cite{KM} the definition of a ``mapping cone of \(-p\)'' \((\check{E}, \check{e})\):
\[
\check{E}=\hat{C}\oplus \bar{C}, \quad
\check{e}=\left[\begin{array}{cc}
\hat{\partial} &0\\
p &\bar{\partial}
\end{array}
\right] .
\]
The short exact sequence associated with \((\check{E}, \check{e})\),
\(0\to\bar{C}\to \check{E}\to \hat{C}\to 0\), 
induces a long exact sequence connecting the triple \(\overline{HM}\),
\(H(\check{E})\), \(\widehat{HM}\), with connecting map \(p_*\). \cite{KM}
 shows that \(\check{E}\) is chain-homotopic to \(\check{C}\). The
following diagram summarizes the construction: 
\[\xymatrix{
&&\check{C}\ar@{->}[d]^{k} &&\\
0\ar@{->}[r] &\bar{C}\ar@{->}[r]^{\bar{i}}\ar@{->}[ur]^{i}
&\check{E}\ar@{->}[r]^{\bar{j}}
\ar@{->}[d]^{l} &\hat{C}\ar@{->}[r]\ar@{->}[ul]_{j} & 0 \\
&&\check{C} &&
}
\]
where 
\[k=\left[\begin{array}{c} j\\ \Pi_s
   \end{array}
 \right], \quad l=\left[\begin{array}{cc} \Pi_o & i
   \end{array}
 \right],\]
with  \[
\begin{split}\Pi_s\co \check{C}=C^o\oplus C^s\to \bar{C}=C^s\oplus C^u,\\
\Pi_o\co \hat{C}=C^o\oplus C^u\to \check{C}=C^o\oplus C^s,\\
\Pi_u\co \bar{C}=C^s\oplus C^u\to
\hat{C}=C^o\oplus C^u
\end{split}
\] denoting projections to the \(s\), \(o\), \(u\)
components respectively.

In terms of these, the proof in \cite{KM}
reduces to the verification of the following identities:
\begin{gather}
lk=\op{Id}, \label{eq:1}\\
kl=\op{Id}+\check{e}K+K\check{e},\label{eq:2}\\
j=\bar{j}k, \label{eq:3}\\
ki-\bar{i}=\check{e}(K\bar{i})+(K\bar{i})\, \bar{\partial},\label{eq:4}
\end{gather}
where 
\[
K=\left[\begin{array}{cc}0 &-\Pi_u\\0&0
\end{array}\right].
\]
We now want to apply the preceeding constructions and identities to the
\(S_U\)-versions. To do so, first observe that the identities (\ref{eq:U-i}) imply
that \(\check{E}\) is an \(H^*(BS^1)\)-module, with the \(U\)-map
given by 
\[
U_{\check{e}}=\left[\begin{array}{cc}\hat{U} &0\\K_p &\bar{U}
  \end{array}
\right].
\]
With this defined, it is straightforward to check that \(\bar{i}\), \(\bar{j}\) are
\(H^*(BS^1)\)-morphisms. 
We can then use what was said in the
previous subsection to form the \(H_*(S^1)\)-modules and
morphisms \((S_U(\check{E}), S_U(\check{e}))\), \(S_U(\bar{i})\),
\(S_U(\bar{j})\). Lemma \ref{lem:u-y} ensures that 
\[
0\to S_U(\bar{C}) \stackrel{S_U(\bar{i})}{\longrightarrow}
S_U(\check{E}) \stackrel{S_U(\bar{j})}{\longrightarrow} S_U(\hat{C})\to  0
\]
is a short exact sequence of \(H_*(S^1)\)-modules. Meanwhile, the
identities (\ref{eq:U-i}) can be used again to verify that:
\begin{gather*}
k\check{U}-U_{\check{E}}k+\check{K}_j\check{\partial}+\check{e}\check{K}_j=0,\\
lU_{\check{E}}-\check{U}l+\check{K}_i\check{e}+\check{\partial}\check{K}_i=0,
\end{gather*}
where \[\check{K}_j=\left[\begin{array}{c} -K_j\\ 0
\end{array}\right], \quad \text{and}\quad 
\check{K}_i=\left[\begin{array}{cc} 0& K_i
\end{array}\right].
\]
This means that \(l\), \(k\) are both p-morphisms of
\(H^*(BS^1)\)-modules. By Lemma \ref{lem:u-y} we can then form the
\(H_*(S^1)\)-module morphisms \(S_U(l)\), \(S_U(k)\). 
The analogs of (\ref{eq:1}), (\ref{eq:2}) 
\begin{equation}\label{eq:S1}
S_U(l)\, S_U(k)=\op{Id}, \quad S_U(j)=S_U(\bar{j)}\, S_U(k) 
\end{equation}
now follow readily from the naturality property of \(S_U\) described in
Lemma \ref{lem:u-y}. Meanwhile, the analogs of (\ref{eq:2}),
(\ref{eq:4}) 
\begin{equation}\label{eq:S2}
\begin{split}
S_U(k)S_U(l) &=\op{Id}+S_U(\check{e})(K\otimes \, \jmath)+(K\otimes\,
\jmath)\,S _U(\check{e}), \\
S_U(k)S_U(i)-S_U(\bar{i})& =S_U(\check{e})\bbK+\bbK S_U(\bar{\partial}),
\quad \bbK:=(K\otimes \, \jmath)S_U(\bar{i})
\end{split}
\end{equation}
reduce to the following identities:
\[\begin{split}
&K_j\Pi_o-\Pi_uK_p=0\\
&\Pi_sK_i+K_p\Pi_u=0\\
&\hat{U}\, \Pi_u-\Pi_u\bar{U}-K_j\, i+jK_i=0;
\end{split}
\]
and these can be directly verified. This proves (\ref{eq:induced-KM1}).

To verify (\ref{eq:induced-KM}), we simply apply \(E_Y^\circ\) to the
\(S_U\)-version of \cite{KM}'s constructions and
identities obtained above. Since we have shown that
\((S_U(\check{E}), S_U(\check{e}))\), \(S_U(\bar{i})\),
\(S_U(\bar{j})\), \(S_U(l)\), \(S_U(k)\) are \(H_*(S^1)\)-morphisms,
Lemma \ref{lem:y-u} implies that \((E_Y^\circ S_U(\check{E}), E_Y^\circ S_U(\check{e}))\), \(E_Y^\circ S_U(\bar{i})\),
\(E_Y^\circ S_U(\bar{j})\), \(E_Y^\circ S_U(l)\), \(E_Y^\circ S_U(k)\)
are \(H^*(BS^1)\)-morphisms, and the analogs of the identities
(\ref{eq:S1}), (\ref{eq:S2}) follow without much ado by applying \(E_Y^\circ\) to them
and the naturality properties of \(E_Y^\circ\) described in Lemma \ref{lem:y-u}.
\epf

\subsection{The \(\bar{C}_*\) complex and localization}

\begin{lemma}\label{lem:induced-KMb}
\(H_*(S_U(\bar{C}))=0\).
\end{lemma}
\pf To compute \(H_*(S_U(\bar{C}))\), write
\[
S_U(\bar{C})=\bar{C}\otimes \bbK[y], \quad
S_U(\bar{\partial})=\bar{\partial}\otimes \, \jmath +\bar{U}\otimes y.
\]
Filtrate this complex by the degree in the factor \(\bar{C}\); this is
done just as in
the proof of Proposition \ref{prop:koszul}. The \(E_1\)-term is
\(\overline{HM}_*\); and \(d_1\) is the \(u\)-map on
\(\overline{HM}_*\). We claim that this map is invertible, and
therefore \(H_*(S_U(\bar{C}))\) vanishes. 

To see that this is indeed the case, 
write \begin{equation}\label{eq:bar-C}
\bar{C}=C_{\bbT}\otimes \bbK[x, x^{-1}],
\end{equation}
where \(C_\bbT\) is the
Morse complex of a Morse function on the torus of flat connections,
which is finitely generated. Recall that a generator \(a\otimes x^m\)
for \(C_{\bbT}\otimes \bbK[x, x^{-1}]\) corresponds to the \(m\)-th eigenvalue of
\(D_a\), the Dirac operator with the flat connection, where the
eigenvalues are ordered by their value in \(\bbR\), and \(1=x^0\)
corresponds to the minimal positive eigenvalue. 

The index of the \(C_\bbT\) factor
defines a finite length filtration on \(\bar{C}\), with respect to
which \(\bar{U}\) can be written as \(\sum_{k=0}^N\bar{U}_k\)
for some \(N\in \bbZ^{\geq 0}\). However, \(\bar{U}_0=x\) (understood as
multiplication), because the only
possible contribution to \(\bar{U}_0\) comes from the moduli space of
instantons from \(a\otimes x^m\) to \(a\otimes x^{m-1}\); and this 
consists of the space of gradient flows of the quadratic function
\(\sum_{m\in \bbZ}\lambda_m|\xi_m|^2\) on
\(\bbP\, (\op{Span}_\bbC\{\eta_m\}_m)\). Here, \(\eta_m\) denotes a
chosen unit-norm 
eigenvector of \(\lambda_m\). This moduli space is \(\bbC P^1\).
The fact that \(\bar{U}_0=x\) is an invertible operator on
\(C_{\bbT}\otimes \bbK[x, x^{-1}]\) then means that \(\bar{U}\) is 
invertible as well. \epf

It follows from the preceding Lemma and Lemma \ref{lem:induced-KMa} that \(S_U(j)\) induces an
\(H_*(S^1)\)-module isomorphism from \(H_*(S_U(\check{C}))\) to
\(H_*(S_U(\hat{C}))\). 
\begin{defn}\label{def:HM-tilde} (Cf. \cite{L} Equation (5.6))
We call the following group the ``total''-version of monopole Floer
homology:
\[\widetilde{HM}_*:=H_*(S_U(\hat{C}))\simeq H_*(S_U(\check{C})).\]
\end{defn}
The motivation for this definition comes from the theory of
\(S^1\)-equivariant theroy;
it is related to the equivariant versions of
Floer homologies \(\widehat{HM}\), \(\widecheck{HM}\),
\(\overline{HM}\) by properties expected of the
homology of their corresponding \(S^1\)-space. 
(The choice of the accent \(\sim\) in the notation reflects
  the fact that this is supposed to come from the 
space of framed configurations, in accordance with the
notation (5.1.1) in \cite{DK}). In particular, the
following lemma is a consequence of 
Proposition \ref{prop:koszul} (a) (ii) and Remarks \ref{rmk:ES} that 
\begin{lemma}\label{lem:1.11}
\(\widetilde{HM}_*\) is related to \(\widehat{HM}_*\) by the following
long exact sequence:
\begin{equation}\label{eq:HM-sq2}
\cdots\to \widehat{HM}_*\stackrel{U}{\to }\widehat{HM}_{*-2}\to \widetilde{HM}_*\to
\widehat{HM}_{*-1}\to \cdots.
\end{equation}
\end{lemma}

The following 
lemma is invoked in the next subsection: 
\begin{lemma}[Localization] 
Let \(\hat{C}\), \(\bar{C}\), \(\widehat{HM}\), \(\overline{HM}\)
denote the monopole Floer complexes or homologies for a balanced
perturbation. Then:

\textbf{ (a)} The map \(i_{V_*}\co H_*(E^-_YS_U(\bar{C}))\to
H_*(E^\infty _YS_U(\bar{C}))\) is an isomorphism.

\textbf{ (b)} 
The map  \(p_*\) induces an
isomorphism of \(\bbK[u, u^{-1}]\)-modules: 
\[p_*\co \widehat{HM}_*\otimes_{\bbK[u]}\bbK[u, u^{-1}]\to \overline{HM}_*\otimes_{\bbK[u]}\bbK[u, u^{-1}].\]
\end{lemma}
\pf
{\em Part (a).} By Proposition \ref{prop:koszul}, it is equivalent
to consider the localization map \(H_*(\bar{C})\to
H_*(\bar{C}\otimes_{\bbK[u]}\bbK[u, u^{-1}])\). However, we saw in the proof
of Lemma \ref{lem:induced-KMb} that the \(u\)-action is
invertible on \(H_*(\bar{C})\).

{\em Part (b).} \(\op{Tor}(\bbK[u], \bbK[u, u^{-1}])=0\); so we can work at 
the chain level. 
\[
H_*(\hat{C}\otimes  _{\bbK[u]}\bbK[u, u^{-1}])=H_*((\hat{U}^N\hat{C})\otimes_{\bbK[u]}\bbK[u, u^{-1}])
\]
for any \(N\in \bbZ^{\geq 0}\). 
There are
finitely many irreducible Seiberg-Witten solutions; and with a
balanced perturbation, the Seiberg-Witten actional functional is
real-valued. We can therefore order these finitely many irreducibles by their
values of action functional. A nonconstant Seiberg-Witten instanton always
decreases the actions unless it is reducible; so for sufficiently large \(N\),
\(\hat{U}^N\hat{C}\subset C^u\). 

Meanwhile, we saw in (\ref{eq:bar-C}) that \(C^u=C_\bbT\otimes
(x\bbK[x])\) and \(C^s=C_\bbT\otimes  \bbK[x^{-1}]\). We also saw in the proof
of Lemma \ref{lem:induced-KMb} that \(\bar{U}_0=x\). Therefore, \(C^u\) generates
\(\bar{C}\otimes _{\bbK[u]}\bbK[u, u^{-1}]\). 
This understood, the assertion follows because we can restrict our
attention to \(C^u\) and the \(u-u\) component of \(p\) is the identity.
\epf

\newpage

\subsection{Monopole Floer homologies from twisted tensor products}\label{sec:HM-ref}

The modules \(S_U(C)\) and \(E_Y(C)\) are  ``twisted
tensor products'' (in the sense of e.g. \cite{TT}, \cite{Le}),
 on which \(H^*(BS^1)\) and \(H_*(S^1)\) respectively
act by simple multiplications. On the other hand, the duality theorem 
Proposition \ref{prop:koszul} tells us the following: On the
homological level, we can
replace any \(H^*(BS^1)\) or \(H_*(S^1)\)-modules by such
twisted tensor products by applying \(E^-_YS_U\) or \(S_UE^-_Y\)
respectively. We shall reformulate the monopole Floer homologies 
\(\widehat{HM}_*\), \(\widecheck{HM}_*\), \(\overline{HM}_*\) defined
in \cite{KM} accordingly. In addition to these three flavors of
monopole Floer homologies, we will introduce a fourth flavor
\(\widetilde{HM}_*\)
from this point of view. These four flavors of
monopole Floer homologies will be regarded as a system and denoted
collectively as \(\mathring{HM}_*\) below. 
We call
\(\widehat{HM}_*\), \(\widecheck{HM}_*\), \(\overline{HM}_*\),
\(\widetilde{HM}_*\) the \,{\em from-, to, bar-, total-}versions of
monopole Floer homology respectively. 
Just as \(\widehat{HM}_*\),
\(\widecheck{HM}_*\), \(\overline{HM}_*\) are to be viewed as 
versions of equivariant homologies of the equivariant
Seiberg-Witten-Floer stable homotopy type (represented by a pointed
\(S^1\)-space) \(\op{SWF}\, (Y, \grc)\) that is introduced in
\cite{M}, what we denoted by 
\(\widetilde{HM}_*\) 
can be viewed as the (non-equivariant) homology of \(\op{SWF}\, (Y, \grc)\) itself.

We now state the main result of this subsection.
\begin{prop}\label{prop:KM-ES} (\cite{L} Corollary 5.3) 
Let \(\hat{C}\) denote \(\hat{C} (M, \grs, c_b)\) and
\(\mathring{\HM}\) denote \(\mathring{\HM}(M, \grs, c_b)\). 
There is a system of isomorphisms 
(as \(H^*(BS^1)\)-modules) from
\(H_*(E^\circ_YS_U(\hat{C}))\) to \(\mathring{HM}_{*}\), taking the fundamental exact
sequence of equivariant homologies for \(S_U(\hat{C})\) to 
the fundamental exact sequence of monopole Floer homologies. In
particular, 
we have the following commutative diagram of \(H^*(BS^1)\)-modules:
\begin{equation}\label{eq:KM}
\begin{split}
\xymatrix{
\llap{\(\cdots       \)}H_*(E^-_YS_U(\hat{C}))\ar@{->}[r]^{i_{V*}}\ar@{->}[d]
&H_*(E^\infty_YS_U(\hat{C}))\ar@{->}[r]\ar@{->}[d]
&H_*(E^+_YS_U(\hat{C})) \ar@{->}[r]^{\delta_{V*}}\ar@{->}[d]
&H_{*-1}(E^-_YS_U(\hat{C}))\rlap{\(\cdots\)}\ar@{->}[d]\\
\llap{\(\cdots       \)}\widehat{HM}_{*}\ar@{->}[r]^{p_*}
&\overline{HM}_{*-1}\ar@{->}[r]^{i_*}
&\widecheck{HM}_{*-1}\ar@{->}[r]^{j_*}
&\widehat{HM}_{*-1}\rlap{\(\cdots\)}}
\end{split}
\end{equation}
where the vertical arrows are \(H^*(BS^1)\)-module-isomorphisms.
\end{prop}

\paragraph{\it Proof of Proposition \ref{prop:KM-ES}.}

Regard \(\mathring{C}\) as  \(\bbZ\)-graded complexes as prescribed in
Remark \ref{rem:period_cpx} and 
consider the following diagram:

\begin{equation}\label{eq:large-diag}
\xymatrix{
\ar@{->}[d]^{E_YS_U(j)}&\ar@{->}[d]^{E_YS_U(j)}&\ar@{->}[d]^{E_YS_U(j)}
&\ar@{->}[d]^{E_YS_U(j)}\\
 \llap{\(\cdots
  \)}H_{*+2}(E^+_YS_U(\hat{C}))\ar@{->}[r]\ar@{->}[d]^{E_YS_U(p)} 
&H_{*+1}(E^-_YS_U(\hat{C}))\ar@{->}[r]^{i_{V*}}\ar@{->}[d]^{E_YS_U(p)}
&H_{*+1}(E^\infty_YS_U(\hat{C}))\ar@{->}[r]\ar@{->}[d]^{E_YS_U(p)}_{\simeq}\ar@{.>}[dl]_{\bar{h}}
&H_{*+1}(E^+_YS_U(\hat{C})) \rlap{\(\cdots\)}\ar@{->}[d]^{E_YS_U(p)}\\
\llap{\(\cdots \)}H_{*+1}(E^+_YS_U(\bar{C}))\ar@{->}[r]\ar@{->}[d]^{E_YS_U(i)}
&H_{*}(E^-_YS_U(\bar{C}))\ar@{->}[r]^{i_{V*}}_{\simeq}\ar@{->}[d]^{E_YS_U(i)}
&H_{*}(E^\infty_YS_U(\bar{C}))\ar@{->}[r]\ar@{->}[d]^{E_YS_U(i)}
&H_{*}(E^+_YS_U(\bar{C}))\rlap{\(\cdots\)}\ar@{->}[d]^{E_YS_U(i)}\\
\llap{\(\cdots
  \)}H_{*+1}(E^+_YS_U(\check{C}))\ar@{->}[r]_{\simeq}^{\delta_{V*}}\ar@{->}[d]^{E_YS_U(j)}_{\simeq}
& H_{*}(E^-_YS_U(\check{C}))\ar@{->}[r]^{i_{V*}}\ar@{->}[d]^{E_YS_U(j)}
&H_{*}(E^\infty_YS_U(\check{C}))\ar@{->}[r]\ar@{->}[d]^{E_YS_U(j)}
&H_{*}(E^+_YS_U(\check{C}))\rlap{\(\cdots\)}\ar@{->}[d]^{E_YS_U(j)}_{\simeq}\\
 \llap{\(\cdots
  \)}H_{*+1}(E^+_YS_U(\hat{C}))\ar@{->}[r]\ar@{->}[d]^{E_YS_U(p)}\ar@{.>}[ur]_{\check{h}}
 &H_{*}(E^-_YS_U(\hat{C}))\ar@{->}[r]^{i_{V*}}\ar@{->}[d]^{E_YS_U(p)}
 &H_{*}(E^\infty_YS_U(\hat{C}))\ar@{->}[r]\ar@{->}[d]^{E_YS_U(p)}_{\simeq}
 &H_{*}(E^+_YS_U(\hat{C}))\rlap{\(\cdots\)}\ar@{->}[d]^{E_YS_U(p)}\\
\llap{\(\cdots \)}H_{*+1}(E^+_YS_U(\bar{C}))\ar@{->}[r]\ar@{->}[d]^{E_YS_U(i)}
&H_{*}(E^-_YS_U(\bar{C}))\ar@{->}[r]^{i_{V*}}_{\simeq}\ar@{->}[d]^{E_YS_U(i)}
&H_{*}(E^\infty_YS_U(\bar{C}))\ar@{->}[r]\ar@{->}[d]^{E_YS_U(i)}
&H_{*}(E^+_YS_U(\bar{C}))\rlap{\(\cdots\)}\ar@{->}[d]^{E_YS_U(i)}\\
&&&
}
\end{equation}
All rows and columns above are exact sequences of
\(H^*(BS^1)\)-modules: the rows are fundamental exact sequences of
equivariant homologies of \(S_U(\hat{C})\), \(S_U(\bar{C})\),
\(S_U(\check{C})\), and the columns are the exact sequences from
(\ref{eq:induced-KM}). 

By Proposition \ref{prop:koszul}, the exact sequence in 
the second column is isomorphic to the first
fundamental exact sequence of the monopole Floer homologies
(\ref{eq:SW-fund-sq}), namely the second row in
(\ref{eq:KM}). Therefore we shall henceforth replace the second column
by 
\[
\cdots\stackrel{j_*}{\longrightarrow} \widehat{HM}_* \stackrel{p_*}{\longrightarrow} \overline{HM}_*\stackrel{i_*}{\longrightarrow} \widecheck{HM}_*\stackrel{j_*}{\longrightarrow} \widehat{HM}_{*-1}
\stackrel{p_*}{\longrightarrow}\cdots. \]
Our goal is therefore to construct an isomorphism from the exact
sequence in the first or fourth row to the exact sequence in the
second column. 
\begin{equation*}
\xymatrix{
\llap{\(\cdots       \)}H_{*+1}(E^-_YS_U(\hat{C}))\ar@{->}[r]^{i_{V*}}\ar@{->}[d]_{\hat{}}
&H_{*+1}(E^\infty_YS_U(\hat{C}))\ar@{->}[r]\ar@{->}[d]_{\bar{h}}
&H_{*+1}(E^+_YS_U(\hat{C})) \ar@{->}[d]_{\check{h}} \rlap{\(\cdots\)}\\
\llap{\(\cdots       \)}H_{*+1}(E^-_YS_U(\hat{C}))\ar@{->}[r]^{E_YS_U(p)}
&H_{*}(E^-_YS_U(\bar{C}))\ar@{->}[r]^{E_YS_U(i)}
&H_{*}(E^-_YS_U(\check{C}))\rlap{\(\cdots\)}}
\end{equation*}

To see this, note that in the third column of (\ref{eq:large-diag}), the map \[E_YS_u(p)\co
H_{*}(E^\infty_YS_U(\hat{C}))\to H_{*}(E^\infty_YS_U(\bar{C}))\]
is an isomorphism by 
the preceding Lemma and Proposition \ref{prop:koszul}. Therefore
\(H_{*}(E^\infty_YS_U(\check{C}))\) is trivial. 
This in turn implies that the map
\[\delta_{V*}\co H_{*+1}(E^+_YS_U(\check{C}))\to H_{*}(E^-_YS_U(\check{C}))\] in the
third row of (\ref{eq:large-diag}) is an isomorpism.
For the same reasons,
the map \[i_{V_*}\co H_{*}(E^-_YS_U(\bar{C}))\to H_{*}(E^\infty_YS_U(\bar{C}))\]
on the first and fifth rows of (\ref{eq:large-diag}) is an isomorphism as well, and therefore
\(H_{*}(E^+_YS_U(\bar{C}))\) is trivial too. This in turn implies
that the map \[E_YS_U(j)\co H_{*}(E^+_YS_U(\check{C}))\to H_{*}(E^+_YS_U(\hat{C}))\]
on the first and fourth columns is an isomorphism. We now take:
\[
\hat{h}=\op{Id}, \quad \bar{h}=i_{V*}^{-1}\circ E_YS_U(p), \quad
\check{h}= \delta_{V_*}\circ (E_YS_U(j))^{-1},
\]
and the Proposition follows.
\epf

\begin{remarks}\label{rmk:5.10}
The chain complex
\(\hat{C}\) in the statement of the preceding proposition may be
replaced by \(\CM (M, \grs, c_-)\) (yet \(\mathring{HM}\) still stands
for \(\mathring{\HM}(M, \grs, c_b)\)). When \(c_(\grs)\) is
non-torsion, this follows from 
Lemmas \ref{lem:S-htpy}, \ref{lem:E-htpy}, and Theorem
\ref{thm:KM}. When \(c_1(\grs)\) is torsion, this is simply because
\(\CM (M, \grs, c_-)=\CM (M, \grs, c_b)\). As we remarked previously in
Section \ref{sec:2.1}, in this case, monotone, balanced, and exact
perturbations are identical notions. 
\end{remarks}

\section{Monopole Floer homology under connected sum}\label{sec:6}
\setcounter{equation}{0}

We follow the (by now) traditional approach to connected sum formula
for Floer homologies that appeared in the instanton Floer homology
setting in \cite{F, D}. To proceed, some setting-up is required. 

\subsection{Preparations}\label{sec:6.1}
Let \(M_1, M_2\) be closed, oriented, connected 3-manifolds and \(\frak{s}_1\),
\(\frak{s}_2\) be \(\Spin^c\)-structures on \(M_1\), \(M_2\)
respectively. Denote by \(M_\sqcup =M_1\sqcup M_2\) the disjoint union
of \(M_1\), \(M_2\). Let \(\grs_\sqcup=(\grs_1, \grs_2)\) denote the \(\Spin^c\)
structure on \(M_\sqcup \) given by \(\grs_1\) and \(\grs_2\).

\paragraph{Part 1: Spin-c structures and gradings.}

Recall from \cite{KM} the interpretation of spin-c structures
and grading via oriented 2-plane fields on the 3-manifold
\(M\). Denote by \(\bbJ(M)\) the set of homotopy classes of oriented 2-plane fields on the 3-manifold
\(M\). According to Proposition 23.1.8 of \cite{KM}, this may be
identified with the set of gradings of the manifold \(M\), as defined
in \cite{KM} p.424.
There is a \(\bbZ\) action on \(\bbJ(M)\), defined by modifying a
representing plane field in a ball in \(M\)  (\cite{KM}, Definition 3.1.2).
Its quotient is the set of spin-c structures over \(M\),
\(\op{Spin}\,( M)=\bbJ(M)/\bbZ\). The orbit over \(\grs\in
\op{Spin}\,( M)\) is the set of gradings for the spin-c structure
\(\grs\), which we denote by
\(\bbJ (M, \grs)\). Let \(c_\grs\) be the
divisibility of \(c_1(\grs)\). The stabilizer of the orbit
\(\bbJ (M, \grs)\) is \(c_\grs\bbZ\); therefore \(\bbJ (M, \grs)\) is a torsor under \(\bbZ/c_\grs\bbZ\).  

Let \(B(p_1)\), \(B(p_2)\) be respectively open balls centered at
\(p_1\in M_1\), \(p_2\in M_2\), and \(\varphi\co
B(p_1)\backslash\{p_1\}\to  B(p_2)\backslash\{p_2\}\) be an
orientation reversing map such that  
\begin{equation}\label{def:M_sharp}
M_{\#}:=M_1\# M_2:=(M_1\backslash\{p_1\}\, \cup \,
 M_2\backslash\{p_2\})/\sim _\varphi.
\end{equation}
As described in \cite{KM}, the \(\bbZ\) action on \(\bbJ (M_i)\) is
induced from the \(C^0(Y_i, SO(3))\)-action on the space of plane fields on
\(M_i\). Each element in the group \(\bbZ\) is represented by
an even-degree element in \(C^0(M_i, SO(3))\) sending \(M_i\backslash B(p_i)\) to \(1\in
SO(3)\). Since this map has even degree, we may choose it to send
\(p_i\in B (p_i)\) to \(1\) as well. The orientation reversing map \(\varphi\) then defines an 
isomorphism: 
\[
\iota_{\bbJ}\co \big(\bbJ (M_1)\times \bbJ(M_2)\big)/\bar{\Delta}\to \bbJ(M_{\#}),
\]
where \(\bar{\Delta}\subset \bbZ\times \bbZ\) denotes the
anti-diagonal. This isomorphism is equivariant with respect to the
residual \(\bbZ\) action on \(\iota_{\bbJ}\co \big(\bbJ
(M_1)\times \bbJ(M_2)\big)/\bar{\Delta}\) and the \(\bbZ\) action of
\(\bbJ( M_{\#})\). Thus, by taking the quotient on
both sides above, one has an induced isomorphism
\[
\iota_S\co \op{Spin}(M_1)\times \op{Spin}(M_2)\to \op{Spin}(M_{\#}).
\]
Let \(\grs_\#=\iota_S (\grs_1, \grs_2)\).

Note that restrictions of \(\iota_\bbJ\) to orbits of the \(\bbZ\) actions
give rise to isomorphisms (also denoted \(\iota_\bbJ\))
 \[\iota_\bbJ\co (\bbJ (M_1, \grs_1)\times
\bbJ (M_2, \grs_2))/\bar{\Delta}\to \bbJ (M_\#, \grs_{\#})\] as
affine spaces under \(\bbZ/c_\#\bbZ\), where  
\(c_\#\) is the g.c.d of \(c_{\grs_1}\) and \(c_{\grs_2}\). 

Recall  from Part 4 of Section
\ref{sec:2.4} the definition of 
\((\hat{C}(M_\sqcup ), \hat{\partial}_{M_\sqcup })\) as a product
complex of \(\hat{C}(M_1)\) and \(\hat{C}(M_2)\). Through out this
section, we adopt
the same assumption that  the Floer
complex \(\mathring{C}(M_1)\) comes from a {\em non-balanced}
perturbation; in particular, \(\hat{C}(M_1)=CM (M_1)\) in our
notation. As observed there, this assumption implies that
\(\mathring{C}(M_\#)\) is also associated with a non-balance
perturbation, and \(\hat{C}(M_\#)=CM (M_\#)\) as well. 
Given an \(\grs_\sqcup=(\grs_1, \grs_2)\in
\op{Spin}(M_\sqcup)\simeq\op{Spin}(M_1)\times \op{Spin}(M_2)\), 
we use \(\bbJ(M_\sqcup,
\grs_{\sqcup})\) to denote the \(\bbZ_{c_\#}\)-grading \(\bbJ
(M_1)\times \bbJ(M_2)\big)/\bar{\Delta}\) on \((\hat{C}(M_\sqcup , \grs_\sqcup),
\hat{\partial}_{M_\sqcup })\) induced from the natural bi-grading of
the latter as a product complex. 

\begin{remarks}
Note that the canonical \(\bbZ/2\bbZ\)-grading (\cite{KM} \S 22.4) of \(\iota_\bbJ
(\xi_1, \xi_2)\) differs from the sum of the canonical \(\bbZ/2\bbZ\)
grading of \(\xi_1\) and \(\xi_2\) by 1. 
\end{remarks}

\paragraph{Part 2: \(\mathbf{A}_\dag\)-actions on
  \(\hat{C}(M_\sqcup)\) and \(\CM (M_\# )\).} 
Use the connected sum decomposition of \(M_\#\),
(\ref{def:M_sharp}), to define a splitting 
\[
H_1(M_\#;\bbZ)\simeq H_1(M_1;\bbZ)\oplus H_1(M_2;\bbZ),
\]
and correspondingly, a factorization of the algebra
\begin{equation}\label{A-sum}\begin{split}
\mathbf{A}_\dag (M_\#)& \simeq \mathbf{A}_\dag (M_1)\otimes
_{\bbK[u]}\mathbf{A}_\dag (M_2)\\ & \simeq \bbK[U]\otimes \textstyle{\bigwedge}
^*(H_1(M_1;\bbZ)/\text{Tors})\otimes
\textstyle{\bigwedge}^*(H_1(M_2;\bbZ)/\text{Tors}).
\end{split}
\end{equation}
This factorization is used to identify \(\mathbf{A}_\dag
(M_\#)\) with \[\begin{split}\mathbf{A}_\dag (M_\sqcup)& : =\bbK[U]\otimes \textstyle{\bigwedge}
^*(H_1(M_\sqcup;\bbZ)/\text{Tors})\\ & =\bbK[U]\otimes\textstyle{\bigwedge}^*(H_1(M_1;\bbZ)/\text{Tors})\otimes
\textstyle{\bigwedge}^*(H_1(M_2;\bbZ)/\text{Tors}).
\end{split}
\]
(Note that \(\mathbf{A}_\dag
(M_\sqcup)\neq H^*(\mathcal{B}^\sigma (M_\sqcup);\bbZ)=\mathbf{A}_\dag
(M_1)\otimes _\bbZ\mathbf{A}_\dag
(M_2)\) with  \(\mathbf{A}_\dag(M_\sqcup)\) so defined.) 

Recall the description of \(\mathbf{A}_\dag(M)\)-actions on monopole
Floer complex from Section \ref{sec:A-module}'s Part 2 and Remark
\ref{rmk:S_1-prod}. These depend
on the choice of a point \(p\in M\) for the \(U\)-map, and a circle
\(\gamma\subset M\) representing \(\grt\) for each element
\(\grt\in \{\grt_i\}_i\), the latter being a basis of \(H_1(M;\bbZ)/\text{Tors}\). We denote the associated maps on
the Floer complex by \(\mathring{\grm}_U=\mathring{U}_p=\mathring{U}\) or \(\mathring{\grm}_\gamma=\mathring{\grm}_\grt\),
depending on emphasis and context. The choices for the manifolds discussed in this section, \(M_*=M_1\),
\(M_2\), \(M_\#\), or \(M_\sqcup\), are given as follows. 

For the
\(H_1(M_*;\bbZ)/\text{Tors}\)-actions, we choose \(b_1(M_1)\)
mutually disjoint embedded circles \(\gamma^{[1]}_i\) in \(M_1\) so that
\(p_1\notin\bigcup_i\gamma^{[1]}_i\) and \(\{\grt_i^{[1]}=[\gamma^{[1]}_i]\}_i\)
forms a basis of \(H_1(M_1;\bbZ)/\text{Tors}\). Choose similarly \(b_1(M_2)\)
mutually disjoint embedded circles \(\gamma^{[2]}_j\) in \(M_2\) so that
\(p_2\notin\bigcup_j\gamma^{[2]}_j\) and \(\{\grt_j^{[2]}=[\gamma^{[2]}_j]\}_j\)
forms a basis of \(H_1(M_2;\bbZ)/\text{Tors}\).
Use \(\{\gamma^{[1]}_i\}_i\), \(\{\gamma^{[2]}_j\}_j\), and
\(\{\gamma^{[1]}_i\}_i\cup \{\gamma^{[2]}_j\}_j\) to define the
\(H_1(M_*;\bbZ)/\text{Tors}\) actions on \(\hat{C}(M_*)\)  
respectively for \(M_*=M_1, M_2, M_\sqcup\). For \(M_*=M_\#\), regard all the
\(\gamma^{[1]}_i\)'s and \(\gamma^{[2]}_j\)'s as embedded circles in \(M_\#\)
through (\ref{def:M_sharp}),
and use them to define the \(H_1(M_\#;\bbZ)/\text{Tors}\) action on \(CM(M_\#)\). 

The \(U\in \mathbf{A}_\dag(M_*)\)-actions for \(M_*=M_1\),
\(M_\sqcup\), or \(M_\#\) are given as follows. Choose a point \(p\in M_1\)
disjoint from \(\{p_1\}\cup \{\gamma^{[1]}_i\}_i\). This \(p\) can
also be viewed as a point in \(M_*=M_1\) or \(M_\sqcup\), or a point in \(M_*=M_\#\) via (\ref{def:M_sharp}). We use
\(\hat{U}_p=\hat{U}_p(M_*)\) to denote the associated \(U\)-action on the monopole
chain complex \(\hat{C}(M_*)\) for such \(M_*\). Note that 
\(\hat{U}_p(M_\sqcup)=U_p(M_1)\otimes 1\) on the product
complex \(\hat{C}(M_\sqcup)=CM(M_1)\otimes \hat{C}(M_2)\). 

\begin{remarks}
As noted previously, \(\hat{U}_p\) and \(\hat{U}_{p'}\) induce
different \(U\)-actions on \(\widehat{HM}(M_\sqcup)\). This is however
irrelevant for our purposes, namely deriving and applying the
connected sum formulas Propositions \ref{prop:conn-1} and
\ref{prop:conn-f}. See Lemma \ref{lem:U-diff} below. We choose \(p\) to be
on \(M_1\) because \(CM (M_1)\) is assumed to be associated with a
nonbalanced perturbation, and for our application this perturbation is
of the type discussed in Part 4 of Section \ref{sec:A-module}, where
\(U_p\) has a nice geometric interpretation. 
\end{remarks}
\paragraph{Part 3: \(\mathbf{A}_\dag (M_\sqcup)\)-actions on \(S_{U_\sqcup}\hat{C}_*(M_\sqcup
  )\). } Let  \(\hat{U}_\sqcup=\hat{U}_{p_2-p_1}\) be as in (\ref{eq:U-cup}). (To simplify 
notation, we shall frequently drop the hat from \(\hat{U}_\sqcup\);
that is,  \(U_\sqcup:=\hat{U}_\sqcup\) in what follows.)
The statement of the upcoming connected sum theorem relates
\(CM_*(M_\#)\) with \(S_{U_{\sqcup}}\big(\hat{C}_*(M_\sqcup
)\big)\),  the  Floer complex is obtained by
applying Section \ref{sec:S_U}'s \(S_U\) operation to \(\hat{C}_*(M_\sqcup )\),
with the latter regarded as an \(H^*(BS^1)\)-module generated by
the \(U=\hat{U}_{\sqcup}\) action above.


Abstractly, an \(\mathbf{A}_\dag (M_\sqcup)\)-action on
\(\hat{C}_*(M_\sqcup )\)  can be used to define 
a corresponding \(\mathbf{A}_\dag (M_\sqcup)\)-action on
\(S_{U_{\sqcup}}\hat{C}_*(M_\sqcup )\), due to the following observations: 
Given any generator \(Q=U\) or \(\grt_i\) of \(\mathbf{A}_\dag(M_\sqcup
 )\) and any map \(\hat{\grm}_Q\) on  \(\hat{C}(M_\sqcup ,
 \grs_{\sqcup })\) underlying the
 \(Q\)-action on the Floer homology, \(\hat{\grm}_Q\) is a p-morphism from
 \(\hat{C}(M_\sqcup , \grs_{\sqcup })\) to itself in the language of
 Section \ref{sec:4}. This in turn is because \(\hat{\grm}_Q\)
 and \(\hat{U}_\sqcup\) induce commutative maps on Floer
 homology. Thus by
 Lemma \ref{lem:u-y}, \(S_{U_{\sqcup}}(\hat{\grm}_Q)\) is defined
 (albeit non-uniquely), and is an 
 \(H_*(S^1)\)-morphism from \(S_{U_{\sqcup}}\big(\hat{C}_*(M_\sqcup ,
  \grs_{\sqcup })\big)\) to itself. 

With the \(\mathbf{A}_\dag (M_\sqcup)\)-action on
\(\hat{C}_*(M_\sqcup )\) fixed in the previous part, 
we define 
  the \(Q\)-action on \(S_{U_{\sqcup}}\hat{C}(M_\sqcup , \grs_{\sqcup })\) to be the
  map \(S_{U_{\sqcup}}(\hat{\grm}_Q)\). In the notation of
  (\ref{def:p-morphism}), these take the following form:  
\[
S_{U_\sqcup}(\hat{U}_p)= \hat{U}_p\otimes  1+\hat{K}_{\hat{U}_p}\otimes y; \quad S_{U_\sqcup}(\hat{\grm}_\gamma)= \hat{\grm}_\gamma\otimes  \jmath+\hat{K}_{\grm_\gamma}\otimes y,
\] 
where \(\hat{K}_{\grm_Q}\) are chain homotopy maps satisfying
(\ref{eq:S_UPhi}). In other words, they satisfy
\begin{equation}\label{eq:U=p-mor}
\begin{split}
[\hat{U}_{p}, \hat{U}_{\sqcup}] & =\hat{\partial}
(M_\sqcup)\, \hat{K}_{U_p}+\hat{K}_{U_p}\, \hat{\partial}(M_\sqcup),\\
[\hat{\grm}_\gamma, \hat{U}_{\sqcup}]&=\hat{\partial}
(M_\sqcup)\, \hat{K}_{\grm_\gamma}-\hat{K}_{\grm_\gamma}\, \hat{\partial}(M_\sqcup). 
\end{split}
\end{equation}

As previously noted, the choice  of the 
  chain homotopy maps \(\hat{K}_{\grm_Q}\) is not unique. In fact, it
  was observed in 
  Remarks \ref{rem:S_U-htpy} that even the homology
  \(H_*(S_{U_\sqcup}\grm_Q)\) depends on the choice of
  \(\hat{K}_{\grm_Q}\) (modulo homotopy).   In this
  article we adopt a 
  particular choice of these \(\hat{K}_{\grm_Q}\) that
  suits our purposes best and has certain nice properties. In
  particular, the homotopy class of \(\hat{K}_{U_p}\) or 
  \(\hat{K}_{\grm_\gamma }\) varies in a consistent manner with \(p\)
  or \(\gamma \), leading to the desired invariance result, Corollary
  \ref{cor:A-mod-inv} below.

To describe these particular choices of \(\hat{K}_{\grm_Q}\), first recall the notions 
  \(\mathring{n}[\op{u}]\), \(\ul{\op{h}}_{\hat{p}}\), \(\op{u}_\gamma \),
  \(\hat{\Theta }_p\) from Section \ref{sec:A-module}. Fixing a set of
  choices  for \(p\) and \(\gamma _i\)'s from Part 2, we set
\begin{equation}\label{def:K_U}
\begin{split}
\hat{K}_{U_p}& :=\hat{n}[d\ul{\op{h}}_\sqcup\wedge
d\ul{\op{h}}_{\hat{p}}]+[ \hat{n}[d\ul{\op{h}}_\sqcup], \hat{\Theta }_p]+[ \hat{\Theta}
_\sqcup , \hat{U}_\sqcup]\\
& =\hat{n}[\op{u}_\sqcup\wedge\op{u}_p],\\
\end{split}
\end{equation}
where \(\op{u}:=\op{u}_{p_2}-\op{u}_{p_1}\). 
For  each \(\gamma\in \{\gamma_i^{[1]}\}_i\cup \{\gamma_i^{[2]}\}_j\),
set 
\begin{equation}\label{def:K_t}
\begin{split}
\hat{K}_{\grm_\gamma } & :=\hat{n}[\op{u}_\gamma d\ul{\op{h}}_\sqcup]+[\hat{\grm}_\gamma, \hat{\Theta}
_\sqcup ]\\
& =\hat{n}[\op{u}_\gamma \op{u}_\sqcup ].
\end{split}
\end{equation}
We now verify that 
\begin{lemma}
The maps \(\hat{K}_{U_p}\), \(\hat{K}_{\grm_\gamma }\) given in
(\ref{def:K_U}) and (\ref{def:K_t}) above satisfy the identities (\ref{eq:U=p-mor}). 
\end{lemma}
\pf 
(i) To verify that (\ref{def:K_U}) satisfies the first identity in
(\ref{eq:U=p-mor}), let 
\(\N^+=(\N^+)_3\) be a 3-dimensional stratified submanifold of
\(\N^+_3(M_\sqcup)\). (Recall the notations \(\N_k (M)\), \(\N_k^+
(M)\) from Section \ref{sec:A-module}. Recall also that 
\((\N^+)_k\) stands for the \(k\)-th step in the stratification,
\(\emptyset\subset \cdots\subset (\N^+)_k\subset\cdots\subset \N^+\),  of \(\N^+\).) Since both \(d\ul{\op{h}}_p\) and 
\(d\ul{\op{h}}_\sqcup\) are closed forms, by \cite{KM}'s Lemma 21.3.1,
Equation (21.4), and Theorem
19.5.4, 
\[
0=\langle d (d\ul{\op{h}}_\sqcup\wedge
 d\ul{\op{h}}_{\hat{p}}), [\N^+]\rangle=\langle d\ul{\op{h}}_\sqcup\wedge
d\ul{\op{h}}_{\hat{p}}, \partial[\N^+]\rangle=\langle 
d\ul{\op{h}}_\sqcup \wedge d\ul{\op{h}}_{\hat{p}}, (\N^+)_2\rangle. 
\]
By construction, \( (\N^+)_2\) is a union of two types of product
spaces, 
the first being of the form \(\N_0^+(\grc_-,
\grc)\times\N_2^+(\grc, \grc_+)\) or \(\N_2^+(\grc_-,
\grc)\times\N_0^+(\grc, \grc_+)\), and the second of the form \(\N_1^+(\grc_-,
\grc)\times\N_1^+(\grc, \grc_+)\). Integrals of \(d\ul{\op{h}}_{\hat{p}}\wedge
d\ul{\op{h}}_\sqcup\) over these two subspaces of \( (\N^+)_2\) give
respectively the first and the second term of the right hand side of the following identity: 
\begin{equation}\label{eq:dn}
0=\big[\hat{\partial}, \hat{n}[
d\ul{\op{h}}_\sqcup \wedge d\ul{\op{h}}_{\hat{p}}]\big]-\big[\hat{n}[d\ul{\op{h}}_{\hat{p}}], \hat{n}[d\ul{\op{h}}_\sqcup]\big].
\end{equation}
By (\ref{e-hol}), \(\hat{n}[d\ul{\op{h}}_{\hat{p}}]=\hat{U}_p-[\hat{\partial},
\hat{\Theta}_p]\), and similarly for
\(\hat{n}[d\ul{\op{h}}_\sqcup]\). Together with the fact that
\(\hat{U}_p,\hat{U}_\sqcup\) are both chain maps, this implies: 
\[\begin{split}
\big[\hat{n}[d\ul{\op{h}}_{\hat{p}}], \hat{n}[d\ul{\op{h}}_\sqcup]\big] & =[\hat{U}_p,
\hat{U}_\sqcup]-\big[[\hat{\partial},\hat{\Theta}_p],  \hat{n}[d\ul{\op{h}}_\sqcup]\big]+\big[\hat{U}_p,
[\hat{\partial}_\sqcup, \hat{\Theta }_\sqcup]\big]\\
&=[\hat{U}_p,\hat{U}_\sqcup]-\big[\hat{\partial},[ \hat{\Theta}_p,
\hat{n}[d\ul{\op{h}}_\sqcup]]\big]
-\big[\hat{\partial}_\sqcup,[\hat{\Theta }_\sqcup, \hat{U}_p]\big].
\end{split}
\]
Inserting this back to (\ref{eq:dn}), the first line of
(\ref{eq:U=p-mor}) follows readily. 

(ii) The second identity in (\ref{eq:U=p-mor}) is verified using similar
arguments. Take now \(\N^+=(\N^+)_2\) to be a 2-dimensional stratified submanifold of
\(\N^+_2(M_\sqcup)\). The coefficients in
\(\big[\hat{\partial}, \hat{n}[\op{u}_\gamma  
d\ul{\op{h}}_\sqcup]\big]\) are given by terms of the form 
\[
\langle \op{u}_\gamma   
d\ul{\op{h}}_\sqcup, \partial[\N^+]\rangle=\langle \op{u}_\gamma  
d\ul{\op{h}}_\sqcup, (\N^+)_1\rangle. 
\]
\( (\N^+)_1\) is a union of product spaces of either the form 
\[\text{
\(\N_0^+(\grc_-,
\grc)\times\N_1^+(\grc, \grc_+)\) or \(\N_1^+(\grc_-,
\grc)\times\N_0^+(\grc, \grc_+)\).}
\] 
Since \(\op{u}_\gamma   \) is
locally constant, integrals of \(\op{u}_\gamma  
d\ul{\op{h}}_\sqcup\) over these spaces take the form of products 
\[
\text{\(\langle
\op{u}_\gamma  , \N_0^+(\grc_-,
\grc)\rangle \, \langle d\ul{\op{h}}_\sqcup, \N_1^+(\grc,
\grc_+)\rangle\) or \(\langle d\ul{\op{h}}_\sqcup, \N_1^+(\grc_-,
\grc)\rangle\, \langle\op{u}_\gamma  , \N_0^+(\grc, \grc_+)\rangle
\). }
\]
By (\ref{m-n}) and (\ref{e-hol}), this shows that 
\[\begin{split}
\big[\hat{\partial}, \hat{n}[\op{u}_\gamma  
d\ul{\op{h}}_\sqcup]\big]& =\big[\grm_\gamma , \hat{U}_\sqcup-[\hat{\partial},
\hat{\Theta}_\sqcup]\big]\\
& =[\grm_\gamma , \hat{U}_\sqcup]-[\hat{\partial},
[\grm_\gamma , \hat{\Theta }_\sqcup]],
\end{split}
\]
leading directly to the second identity of (\ref{eq:U=p-mor}). 
\epf

This understood, we may now justify the claim in Remark
\ref{rmk:S_1-prod} that the \(U\)-action on \(H_*(S_{U_\sqcup}\big(\hat{C}_*(M_\sqcup ,
  \grs_{\sqcup })\big)\) is independent of \(p\).

\begin{lemma}\label{lem:U-diff}
For any given \(p, p'\in M_\sqcup\), \(S_{U_\sqcup}(\hat{U}_{p'})\),
\(S_{U_\sqcup}(\hat{U}_p)\) are chain homotopic. 
\end{lemma}
\pf 
We wish to show that there is a map \(Z_*\co S_{U_\sqcup}(\hat{C})\to
S_{U_\sqcup}(\hat{C})\) such that 
\[
\begin{split}
& S_{U_\sqcup}(\hat{U}_{p'})-S_{U_\sqcup}(\hat{U}_p)
 =(\hat{U}_{p'}-\hat{U}_p)\otimes
 1+(\hat{K}_{U_p'}-\hat{K}_{U_p})\otimes y\\
& =[Z_*, D_\sqcup]. 
\end{split}
\]
We choose \(Z_*\) to be of the form   \(Z_*=Z_0\otimes \jmath +Z_1\otimes y\), with \(Z_0, Z_1\) being maps
from \(\hat{C}(M_\sqcup)\) to itself, respectively of degree \(-1\)
and \(-2\). The preceding identity now 
reads:
\begin{equation}\label{S_U-indep}
\begin{split}
(\hat{U}_{p'}-\hat{U}_p)& =[Z_0, \hat{\partial}]\\
(\hat{K}_{U_p'}-\hat{K}_{U_p})& =[Z_1, \hat{\partial}]+ [\hat{U}_\sqcup,
Z_0]. 
  \end{split}
\end{equation}

(i) In the case when \(p, p'\) belong to the same connected component
of \(M_\sqcup\), there is a path \(\lambda \) in \(M_\sqcup\) from
\(p\) to \(p'\) and an associated map \(\hat{\op{K}}_\lambda \),
defined in (\ref{def:opK}). Since 
by (\ref{def:opK}), \((\hat{U}_{p'}-\hat{U}_p)=[\hat{\op{K}}_\lambda
,\hat{\partial}]\), setting \(Z_0\) to be 
\[
Z_0^\lambda : = \hat{\op{K}}_\lambda
\]
suffices to validate the first line of (\ref{S_U-indep}). We claim
that with \(Z_0\) so chosen, and with \(K_{U_p}\) given by
(\ref{def:K_U}), the second line of (\ref{S_U-indep}) also holds if
\(Z_1\) is set to be 
\[
Z_1^\lambda :=\hat{m}[\uu_\sqcup\,\uu_\lambda ](\bbR\times M_\sqcup),
\]
where \(\uu_\sqcup=\uu_{\hat{p}_2}-\uu_{\hat{p}_1}\), and
\(\hat{p}_i\) denotes \(\bbR\times \{p_i\}\subset \bbR\times
M_\sqcup\). 
In other words, 
\[
\hat{n}[\op{u}_{p'}\op{u}_\sqcup]-\hat{n}[\op{u}_{p}\op{u}_\sqcup]=\big[\hat{m}[\uu_\sqcup\,\uu_\lambda
](\bbR\times M_\sqcup), \hat{\partial}_\sqcup\big]+ \big[\hat{n}[
\op{u}_\sqcup], \hat{m}[\uu_\lambda
](\bbR\times M_\sqcup)\big]. 
\]
This identity is essentially a higher-degree variant of
(\ref{def:opK}), and is proved by arguments similar to
(\ref{value1}). For more details the reader is referred to the proof
of (\ref{eq:A-action-intertwine0}) in the next subsection, which
differs by cosmetic changes from the proof for the preceding
identity. 

Let \(Z_{*}^\lambda =Z_0^\lambda \otimes \jmath +Z_1^\lambda \otimes y\)
denote the version of \(Z_*\) constructed using \(\lambda \). 


(ii) Now suppose that \(p, p'\) belong to different connected components of
\(M_\sqcup\). Let \(p\in M_1\) and \(p'\in M_2\); \(\lambda _1\subset M_1\) be a
path from \(p_1\) to  \(p\); \(\lambda _2\subset M_1\) be a
path from \(p_2\) to \(p'\). Then
by the discussion in Case (i) above, 
\[
S_{U_\sqcup}(U_{p'})-S_{U_\sqcup}(U_{p})=S_{U_\sqcup}(U_{\sqcup})+
[Z^{\lambda _2}_*-Z^{\lambda _1}_*, D_\sqcup].\] 
Meanwhile, according
our construction of \(K_{U_p}\), \(K_{U_\sqcup}=0\) and 
\[
S_{U_\sqcup}(U_{\sqcup})=\hat{U}_\sqcup\otimes 1=[1\otimes \partial_y,
D_\sqcup]. 
\]
So we have
\[
S_{U_\sqcup}(U_{p'})-S_{U_\sqcup}(U_{p})=[Z^{\lambda _2}_*-Z^{\lambda
  _1}_*+1\otimes \partial_y, D_\sqcup] 
\]
and we take \(Z_*=Z^{\lambda _2}_*-Z^{\lambda
  _1}_*+1\otimes \partial_y\) in this case. \epf


\begin{cor}\label{cor:A-mod-inv}
The \(\mathbf{A}_\dag (M_\sqcup)\) action on \(S_{U_\sqcup}(\hat{C}
(M_\sqcup))\), as defined above, induces an \(\mathbf{A}_\dag (M_\sqcup)\) action on \(H_*(S_{U_\sqcup}(\hat{C}
(M_\sqcup)))\) that is independent of choices of \(p\), \(\{\gamma
_i^{[1]}\}_i\), \(\{\gamma_j^{[2]}\}_j\).
\end{cor}
\pf The assertion regarding the \(U\)-action follows directly from the
previous lemma. To verify the assertion for
\(H_1(M_\sqcup;\bbZ)/\text{tors}\)-actions, 
take two embedded circles  \(\gamma \), \(\gamma '\) in \(M_\sqcup  \)
representing the same \([\gamma
]\in H_1(M_\sqcup; \bbZ)/\text{tors}\). Note that they must lie in the
same connected components of \(M_\sqcup\), and therefore there exists
an embedded surface in \(\bbR\times M_\sqcup\) from \(\gamma \) to
\(\gamma '\).  We wishes to show that there exists a map \(T_*\) from
\(S_{U_\sqcup}(\hat{C}(M_\sqcup))\) back to itself, satisfying: 
\[
S_{U_\sqcup}(\hat{\grm}_{\gamma '})-S_{U_\sqcup}(\hat{\grm}_\gamma
)=[T_*, D_\sqcup]. 
\]
Assume this time that \(T_*\) is of the form \(T_*=T_0\otimes
1+T_1\otimes y\), where \(T_0, T_1\) are maps
from \(\hat{C}(M_\sqcup)\) to itself, respectively of degree \(0\)
and \(-1\). Then the preceding identity now 
reads:
\begin{equation}\label{S_m-indep}
\begin{split}
(\hat{\grm}_{\gamma '}-\hat{\grm}_\gamma )& =[T_0, \hat{\partial}]\\
(\hat{K}_{\grm_{\gamma '}}-\hat{K}_{\grm_\gamma })& =[T_1, \hat{\partial}]+ [\hat{U}_\sqcup,
T_0]. 
  \end{split}
\end{equation}
By (\ref{eq:dF}), \(\hat{\grm}_{\gamma '}-\hat{\grm}_\gamma
=[\hat{m}[\textsc{f}_\Sigma ], \hat{\partial}_\sqcup ]\); so we set 
\[
  T_0=\hat{m}[\textsc{f}_\Sigma ](\bbR\times M_\sqcup). 
\]
so that the first line of (\ref{S_m-indep}) holds for this \(T_0\). 
To verify the second line, we choose 
\[
T_1=\hat{m}[\uu_\sqcup\textsc{f}_\Sigma ] (\bbR\times M_\sqcup). 
\]
With this choice and our construction of \(\hat{K}_{\grm_\gamma }\),  the second line of (\ref{S_m-indep}) says
\[
\hat{n}[\op{u}_{\gamma '}\op{u}_\sqcup]-\hat{n}[\op{u}_{\gamma }\op{u}_\sqcup]=\big[\hat{m}[\uu_\sqcup\,\textsc{f}_\Sigma 
](\bbR\times M_\sqcup), \hat{\partial}_\sqcup\big]+ \big[\hat{n}[
\op{u}_\sqcup], \hat{m}[\textsc{f}_\Sigma 
](\bbR\times M_\sqcup)\big]. 
\]
The proof of this identity is virtually identical to that for the
second line of (\ref{eq:A-action-intertwine0}), and the reader is
referred the reader to the next subsection for details. 
\epf


\paragraph{Part 4: The cobordisms \(\mathcal{V}\),
  \(\bar{\mathcal{V}}\), and the cobordism maps \(V_*\),
  \(V^\dag_*\).}  
Let \(\mathcal{V}:=(X, s)\) denote a cobordism described in
(\ref{(A.9a,11)}) and (\ref{eq:(A.12,15a)}), with
\(Y_-=M_\#\) and \(Y_+=M_\sqcup \). Assume that \(s\) has a unique
critical point of index 3 with critical value 0. 

There is a unique \(\Spin^c\)
structure \(\grs_X\) on such \(X\) with
\(c_1(\grs_X)|_{s^{-1}(-c)}=c_1(\grs_\#)\) and
\(c_1(\grs_X)|_{s^{-1}(c)}=c_1(\grs_\sqcup\)) for \(c\gg 0 \).
Meanwhile, given \([\varpi_i]\in H^2(M_i)\), there is a unique
\([\varpi_\#]\in H^2(M_{\#})\) and a \([\omega]\in H^2(X)\) that
restricts to \([\varpi_1]\), \([\varpi_2]\), \([\varpi_\#]\) respectively on
the \(M_1\), \(M_2\), \(M_{\#}\) ends of \(X\). Suppose as before that \([\varpi_1]\)
is non-balanced with respect to \(c_1(\grs_1)\) (and therefore
\([\varpi_\#]\) is also non-balanced with respect to
\(c_1(\grs_\#)\)). 
Let \(\omega\) be a
closed 2-form on \(X\) representing \([\omega]\) above, so that
\(\varpi_X =2\omega^+\) satisfies (\ref{eq:(A.13a)}). In particular,
\(\omega\) restricts respectively to (pull-backs) of closed 2-forms \(\varpi_1\),
\(\varpi_2\), \(\varpi_{\#}\) on the \(M_1\)-, \(M_2\)-, \(M_\#\)-ends
of \(X\). Let \(\varpi_\sqcup\) denote the two-form on
\(Y_+=M_\sqcup\) that restricts respectively to \(\varpi_1\),
\(\varpi_2\) on the \(M_1\) and \(M_2\) component of \(M_\sqcup\). 

Let
\(\bar{\mathcal{V}}:=(X, -s) \) denote the ``time-reversal'' of
\(\mathcal{V}\).
Given local systems \(\Gamma_i\) on \(\mathcal{B}^\sigma(Y_i)\),
\(i=1, 2\), let \(\Gamma_\sqcup=\Gamma_1\otimes \Gamma_2\) denote the
local system on \(\mathcal{B}^\sigma(M_1)\times
\mathcal{B}^\sigma(M_2)\simeq\mathcal{B}^\sigma(M_\sqcup )\). 
Note that \(X\) satisfies the condition that \(\delta_\pm\) are both
isomorphisms in the first bullet of Remark \ref{rem:2.1}, and thus
there is a unique \(\mathcal{V}\)-morphism \(\Gamma_{\cal V}\), which
together with its inverse \(\Gamma_{\bar{\cal V}}\), gives an 1-1
correspondence between local systems on \(\mathcal{B}^\sigma(M_\sqcup
)\) and \(\mathcal{B}^\sigma (M_\#)\). Let \(\Gamma_\#\) denote the
local system on \(\mathcal{B}^\sigma (M_\#)\) corresponding to
\(\Gamma_\sqcup\). Meanwhile, by the second bullet of Remark
\ref{rem:2.1}, \(\Gamma_\sqcup\) is (resp. strongly) \((\grs_\sqcup,
\varpi_\sqcup)\)-complete iff \(\Gamma_\#\) is (resp. strongly) \((\grs_\#,
\varpi_\#)\)-complete, and in this case 
\(\hat{m}[\mathpzc{u}](X, \grs_X, \varpi_X; \Gamma_X)\) is well-defined through
(\ref{eq:m-coeff-loc}). 


In what
follows, take 
\[\begin{split}
\hat{C}_*(M_\sqcup ) & =
\hat{C}_*(M_\sqcup ,
  \grs_{\sqcup },\varpi_{\sqcup }; \Gamma_{ \sqcup})\\
\CM_*( M_\#) & =\CM_*\big(M_{\#}, \grs_\#, \varpi_\# ;  \Gamma_{ \#}\big).
\end{split}
\]
The statement of the upcoming connected sum theorem involves certain
maps 
\[
\text{ \( V_*\co
\CM_*\big(M_{\#})\to  S_{U_\sqcup}\hat{C}_*(M_\sqcup )\) and  \(V_*^\dag\co
S_{U_\sqcup}\hat{C}_*(M_\sqcup )\to \CM_*\big(M_{\#})\). }
\]
These are constructed using
the moduli spaces \(\mathcal{M}_k(\mathcal{V}, \grc_\#,
\grc_{\sqcup})=\mathcal{M}_k(\bar{\mathcal{V}},  \grc_{\sqcup}, \grc_\#)\) of solutions to 
(\ref{eq:(A.14)}) associated to the \(\Spin^c\) 4-manifold \((X,
\grs_X)\) and the perturbation form \(\varpi_X\) described
above. 

Here is how they are defined. 
  Use the short-hand \((\hat{C}_\sqcup,
  \hat{\partial}_{\sqcup})=(\hat{C}(M_\sqcup  ),\hat{\partial}(M_\sqcup))\),
\((C_\#, \partial_\#)=(\CM (M_{\#}), \partial_{M_\#}))\) etc below. 
Write the chain module of \(S_{U_{\sqcup}}(\hat{C}_{\sqcup})\),
\(\hat{C}_{\sqcup}\otimes \bbZ[y]\), as the direct sum:
\begin{equation}\label{eq:S}
S_{U_{\sqcup}}(\hat{C}_{\sqcup})= 
\hat{C}_{\sqcup}\oplus y\,\hat{C}_{\sqcup}.
\end{equation}
With respect to this decomposition, it differential takes the
following block form: 
\begin{equation}\label{S}
D_{\sqcup}=\left[\begin{array}{cc}
\hat{\partial}_{\sqcup} & 0 \\
\hat{U}_{\sqcup} &- \hat{\partial}_{\sqcup}
\end{array}\right]. 
\end{equation}
Correspondingly, write the maps \(V_*\), 
\(V_*^\dag\) 
in  block form with respect to the decomposition (\ref{eq:S}) as
\begin{equation}\label{eq:V-m}
V_*=\left[ \begin{array}{c} V_0\\ V_1
  \end{array}\right], \quad
V^\dag_*= \left[ \begin{array}{cc}V^\dag_1 & V^\dag_0
  \end{array}\right],
\end{equation}
where \(V_i\co \hat{C}_\sqcup\to \hat{C}_\#\), \(V^\dag_i\co
\hat{C}_\#\to \hat{C}_\sqcup\), \(i=0,1\), are defined through
cobordism maps of the form  
\(\hat{m}[\uu](X, \grs_X, \varpi_X; \Gamma_X)\) for \(X=\mathcal{V}\)
or \(\bar{\mathcal{V}}\). These cobordism maps are defined as in Part 4 of
Section \ref{sec:2.4}, noting that \(\mathcal{V}\),
\(\bar{\mathcal{V}}\) satisfy the condition (\ref{eq:assumption-X}),
and assuming for the rest of this subsection  the same completeness condition for \(\Gamma _X\) 
alluded to in the end of Section \ref{sec:2.4}. 
Meanwhile, cochains \(u\) involved in the definition of these maps
are of the type 
introduced in Section \ref{sec:A-module}'s Part 3(a), with the relevant arc \(\lambda\)
chosen as follows. 
In the present section, let \(\lambda \) denote the ascending manifold of the unique critical
point of \(s\); it is a path in \(X\) asymptotic to 
 \((p_1, p_2)\in M_\sqcup =Y_+\). We orient it so that it begins from
 \(p_1\in M_1\) and ends at \(p_2\in M_2\).
Meanwhile, the descending manifold of this critical point will be
denoted \(\textsc{b}\); it is an embedded 3-ball in \(X\) that intersects
each \(s^{-1}(c)\simeq M_{\#}\) in a 2-sphere, \(\forall c\ll 0\). We
orient it so that it intersects with \(\lambda\) positively. Let
\(\bar{\lambda}\), \(\bar{\textsc{b}}\) respectively denote the descending and
ascending manifold from the unique critical point of \(-s\). These are
the same submanifolds in \(X\) as \(\lambda\), \(\textsc{b}\), but equipped with
the opposite orientation.

With the above said,  we are ready to write down the formulas for
where \(V_i\), \(V^\dag_i\), \(i=0,1\):
\begin{equation}\label{eq:V}
\begin{split}
V_0& =\hat{m}[1](\mathcal{V}, \grs_X, \varpi_X; \Gamma_{\mathcal{V}}), \\
V_1& = \hat{\op{K}}_\lambda (\mathcal{V}, , \grs_X, \varpi_X;
\Gamma_\V)=\hat{m}[\mathpzc{u}_\lambda](\mathcal{V}, , \grs_X, \varpi_X; \Gamma_\V)\\
& =\hat{m}[\theta_\lambda](\mathcal{V}, , \grs_X, \varpi_X;
\Gamma_\V)+\Theta_\sqcup\, V_0;\\
V_0^\dag&=\hat{m}[1](\bar{\mathcal{V}}, \grs_X, \varpi_X; \Gamma_\V),\\
V^\dag_1& =\hat{\op{K}}_{\bar{\lambda }}(\bar{\mathcal{V}}, \grs_X,
\varpi_X;
\Gamma_{\bar{\V}})=\hat{m}[\mathpzc{u}_{\bar{\lambda}}](\bar{\mathcal{V}},
\grs_X, \varpi_X; \Gamma_{\bar{\V}})\\
& =\hat{m}[\theta _{\bar{\lambda }}](\bar{\mathcal{V}}, \grs_X,
\varpi_X; \Gamma_{\bar{\V}})- V_0^\dag\, \Theta _\sqcup , 
\end{split}
\end{equation}
where \(\Theta _\sqcup \) denotes the map from
\(\hat{C}_\sqcup=\hat{C}(M_1)\otimes \hat{C}(M_2)\) to itself
\(1\otimes \Theta _{p_2}-\Theta _{p_1}\otimes 1\). 
 Cf. (\ref{def:K_lambda}) for the defintion of \(\hat{\op{K}}_\lambda \).

\begin{remarks}
With \cite{KM}'s notion of canonical \(\bbZ/2\)-gradings suitably
generalized, 
the maps \(V_*\), \(V_*^\dag\)  are of degree 0 with respect
this canonical \(\bbZ/2\)-grading. 
Recall the characteristic number \(\iota (X)\) for a cobordism \(X\)
from \(Y_-\) and \(Y_+\), \(Y_\pm \) both connected, from Definition 25.4.1
in \cite{KM}. When \(Y_\pm\) are allowed to be disconnected, we
generalize the formula in \cite{KM} as:
\[
\iota (X):=b^0(X)-b^1(X)+b^{2+}(X)+\frac{1}{2} \big( b_1(Y_+)-b_1(Y_-)-b_0(Y_+)+b_0(Y_-)\big), 
\] 
where \(b^{2+}(X)\) is the rank of maximal positive definite subspaces
in the image of \(H^2(X, \partial X)\to H^2(X)\) with respect to the
intersection form. 
With this generalized \(\iota (X)\), the statements of \cite{KM}'s Lemma 25.4.2 and Proposition 25.4.3
of \cite{KM} remain valid: \(\iota(X)\in \bbZ\) and is additive under
 composition of cobordisms, and a map  (if well-defined) of the form 
 \(\mathring{m}[\uu](X)\) 
 is of even or odd degree with respect to this canonical \(\bbZ/2\)
 grading depending on the parity of 
\[
\deg (\uu)-\iota(X).
\] 
For \(X=\V\), \(\bar{V}\), 
 \(\iota (\V)=0\) and \(\iota (\bar{\V})=1\); hence \(V_0\),
 \(V_1^\dag\) are of even degree, while \(V_1\) and \(V_0^\dag\) are
 of odd degree. Hence 
\(V_*=V_0+y\, V_1\co C_\#\to \hat{C}_\sqcup\otimes
\bbK[y]\), \(V_*^\dag\co V_1^\dag+V_0^\dag\partial_y\co
\hat{C}_\sqcup\otimes \bbK[y]\to C_\#\)  are both of even {\em degree} with respect to
the canonical grading.  (In fact, they are both of degree 0 when the canonical
 \(\bbZ/2\) grading lifts to an absolute grading (cf. Section 28.3 in
 \cite{KM}).)  This is not to be confused with the notion of
an {\em even or odd  map} in the sense of signs when it appears in
commutators. In the latter sense \(V_*\) is even, while \(V_*^\dag\)
is odd, since \(V_0\), \(V_0^\dag\) are even and \(V_1\),
\(V_1^\dag\) are odd. The parity of a map \(\mathring{m}[\uu](X)\) (in the sense of
commutators) is determined by {\em purely by \(\deg (\uu)\)}, 
independent of \(X\).   This is because only \(\deg (\uu)\)
contributes 
to the signs in gluing formulas. 
\end{remarks}

\subsection{A connected sum formula for non-balanced perturbations}

Adopt the notation and assumptions from the previous
subsection. 
\begin{prop}\label{prop:conn-1}
Under the above assumptions: 

{\bf (a)} Suppose that \([\varpi_\#]\)
is negative monotone, non-balanced with respect to
\(\grs_\#\). Let \(\Gamma_\#\) be arbitrary, and \(\Gamma_\sqcup\) be
determined by \(\Gamma_\#\) via \(\Gamma_X\). Then the maps \(V_*\)
\(V_*^\dag\) given in the previous subsection are well-defined chain
maps, and \(V_*\) defines a chain homotopy equivalence
\begin{equation}\label{eq:conn-sum}
V_*\co C_*\big(M_{\#}, \grs_\#, [\varpi_\#], \Gamma_{ \#}\big) \to 
S_{U_{\sqcup}}\big(\hat{C}_*(M_\sqcup ,
  \grs_{\sqcup }, [\varpi_{\sqcup }]; \Gamma_{ \sqcup})\big)
\end{equation}
respecting the 
 (relative) \(\bbZ/c_\#\)-grading on both sides. Moreover, the map
 \(V_*\) intertwines with 
the \[\mathbf{A}_\dag(M_{\#})\stackrel{(\ref{A-sum})}{\simeq}\textstyle{\bigwedge}
^*(H_1(M_1)/\text{Tors})\otimes
\textstyle{\bigwedge}^*(H_1(M_2)/\text{Tors})\otimes
\bbK[u]=\mathbf{A}_\dag(M_{\sqcup})\] actions on the two sides,
defined in the previous subsection's Parts 2 and 3 
using \(p\) and \(\{\gamma_i^{[1]}\}_i\cup \{\gamma_j^{[2]}\}_j\). 

{\bf (b)} Suppose that
\([\varpi_1]\) is nonbalanced with respect to
\(\grs_1\), and that \(\Gamma_i\) is strongly \((\grs_i, [\varpi_i])\)-complete
for \(i=1, 2\).  Then  the maps \(V_*\)
\(V_*^\dag\) are well-defined chain
maps, and \(V_*\) defines a chain homotopy equivalence
\[
V_*\co C_\bullet\big(M_{\#}, \grs_\#, 
[\varpi_\#], \Gamma_{ \#}\big) \to 
S_{U_{\sqcup}}\big(\hat{C}_\bullet(M_\sqcup ,
  \grs_{\sqcup }, [\varpi_{\sqcup }]; \Gamma_{ \sqcup})\big)
\]
respecting the 
 (relative) \(\bbZ/c_\#\)-grading on both sides. Moreover, the map
 \(V_*\) above intertwines with 
the \[\mathbf{A}_\dag(M_{\#})\stackrel{(\ref{A-sum})}{\simeq}\textstyle{\bigwedge}
^*(H_1(M_1)/\text{Tors})\otimes
\textstyle{\bigwedge}^*(H_1(M_2)/\text{Tors})\otimes
\bbK[u]=\mathbf{A}_\dag(M_{\sqcup})\] actions on the two sides defined
using \(p\) and \(\{\gamma_i^{[1]}\}_i\cup \{\gamma_j^{[2]}\}_j\). 
\end{prop}

\pf {\bf Part (a):} 
The proof has six steps.

\paragraph{\it Step 1.} In this part we show that 
the assumption on \(\Gamma _X\) of Part (a) ensures that \(\Gamma _X\)
satisfy the completeness conditions alluded to in Remark
\ref{rem:2.1}, so that the maps \(V_i\), \(V_i^\dag\), \(i=0,1\), are
well-defined. More precisely, we show that the sum
\begin{equation}\label{eq:m-V1}
\sum_{\grc_\#\in \grC(M_\#)}\sum_{(\grc_1, \grc_2)\in \grC(M_1)\times \grC(M_2)}
\sum_{z\in \pi_0(\mathcal{B}^\sigma(\mathcal{V}; \grc_\#, (\grc_1,
  \grc_2)))}\big\langle \mathpzc{u}, \mathcal{M}_{k, z}(\grc_\#, (\grc_1,
 \grc_2))\big\rangle
\end{equation}
has finitely many non-vanishing terms, and therefore 
\(V_*\) is a well-defined map between the ({\em pre-completed})
chain complexes \(\CM_*(M_\#)\), \(S_{U_\sqcup} (\hat{C}_*(M_\sqcup))\)
for any coefficients \(\Gamma_\#\) and its twin \(\Gamma_\sqcup\). 
To see this, observe that by the
well-known compactness property of spaces of 3-dimensional Seiberg-Witten
solutions, \(\CM\, (M_1)=C^o(M_1)\), \(C^o(M_2)\), \(\CM\,
(M_{\#})=C^o(M_{\#})\) are all finitely generated over \(\bbK\),
while \(C^u(M_2)\) is finitely generated over \(\bbK[u]\), with
\(u\) having degree \(-2\). Write
the generating sets of these free \(\bbK\)-modules respectively as
\(\grC(M_1)=\{\gra_i\}_i\), \(\grC^o(M_2)=\{\grb^o_j\}_j\),
\(\grC(M_{\#})=\{\grc_k\}_k\), and
\(\grC^u(M_2)=\{\grb^u_qu^n\}_{q,n}\), where there are finitely
many indices \(i, j, k, q\), and \(n\) runs through all non-negative
integers. 
Let \(\pi^\sigma\co \mathcal{B}^\sigma\to \mathcal{B}\)
denote the projection of the blown-up space. The index \(\imath_\grd\)
and the topological energy (cf. \cite{KM} Definition 4.5.4 and p.593
in the case of  non-exact perturbations) of an element
\(\grd\in \mathcal{M}(\mathcal{V})\) depends only its relative
homotopy class under \(\pi^\sigma\), and the former is controlled via
\(c_1(\grs_X)\), the latter through \([\omega]-2\pi[c_1(\grs_X)]\). The
monotonicity condition and the compactness property of
\(\mathcal{M}(\mathcal{V})\) under bounds on the topological energy
then ensures that only finitely \(\gra_i\), \(\grb^o_j\),
\(\grb^u_q\), \(z\) appear in the sum on the right hand side of
(\ref{eq:m-V1}). Meanwhile, since
\(\op{gr}(\grb^u_qu^n)-\op{gr}(\grb^u_qu^m)=-2(n-m)\), the index
bound \(\imath_\grd=k\) on the right hand side of (\ref{eq:m-V1})
implies that for each \(q\), only finitely many \(\grb^u_qu^n\) appears
on the right hand side of (\ref{eq:m-V1}).
(The aforementioned compactness result follows from a straightforward
generalization of Theorem 24.5.2 in
\cite{KM} to include nonexact perturbations).

The \(\hat{m}(\bar{\mathcal{V}})\) analog of
(\ref{eq:m-V1}) involves sum over \(\grC(M_{\#})\) instead, which
consists of finitely many elements. The finiteness of the relevant sum
then follows from the monotonicity assumption alone. 

\paragraph{\it Step 2.} In this step, we show that \(V_*\),
\(V^\dag_*\) are (respectively even and odd) chain
maps. 
This amounts to verifying the following identities:
\begin{equation}\label{eq:chain-maps}
\begin{split}
\hat{\partial}_\sqcup \, V_0-V_0\partial_\#=0\\
\hat{\partial}_\sqcup \, V_1+V_1\partial_\#-\hat{U}_\sqcup V_0=0\\
\partial_\#\, V^\dag_0-V^\dag_0\hat{\partial}_\sqcup=0,\\
V^\dag_1\hat{\partial}_\sqcup +\partial_\#\, V^\dag_1+V^\dag_0 \hat{U}_\sqcup=0.
\end{split}
\end{equation}
In view of Equation (\ref{eq:K-lambda}),  
these would have followed directly from \cite{KM}'s Proposition
25.3.4, if there latter's assumption on the connectedness of \(Y_\pm\)
could be removed. In the specific setting under discussion, such
generalization requires only simple modifications of what was in \cite{KM}.  
To do so, write the identities in full in terms
of \(m^o_{o\natural}\), \(m^{o\sharp}_o\), \(\partial^o_o(M_1)\),
\(\partial_o^o (M_{\#})\), \(\hat{\partial}(M_2)\),
\(\bar{\partial}^s_u(M_2)\) as given by (\ref{eq:hat-d}),
(\ref{eq:m-V}), (\ref{eq:V-m}). These can be reduced to 
the identities in Lemma 25.3.6 in \cite{KM} (with many vanishing terms),
with these substitutions:
\begin{itemize}
\item Drop the ``\(o\)'''s from the double superscript or subscripts
  \(o*\) of \(m\), 
\item Replace the entries of  \(\hat{\partial}(M_\sqcup )=(1\otimes \partial^\sharp_\natural(M_2)+\partial^o_o(M_1)\otimes
  1)\)) by \(\partial^\sharp_\natural\).
\end{itemize}
Theorem 24.7.2 in \cite{KM} conveniently supplies us with the general gluing theorem required for verifying these formulae. (We have at worst rank 1
boundary-obstruction). 


\paragraph{\it Step 3.} In the upcoming three steps, we show that the two chain 
complexes in (\ref{eq:conn-sum}) are chain-homotopy equivalent via
\(V_*\) and \(V_*^\dag\). More precisely, we shall show that their compositions satisfy the following
identities: 
\begin{equation}\label{eq:cob-comp1}\begin{split}
V^\dag_*\circ V_* -[H'_\#, \partial_\#]_{\text{even}}&=V^\dag_1\circ V_0+V^\dag_0\circ
V_1\\
&=\op{Id}_\#-[\smz_\#, \partial_\#]_{\text{odd}}, \\
\end{split}
\end{equation}
\begin{equation}\label{eq:cob-comp2}
\begin{split}
V_*\circ V^\dag_*-[H'_\sqcup, D_\sqcup]_{\text{even}}
&=\op{Id}_\sqcup\otimes
1-[\smz_\sqcup\otimes \jmath+\smx \otimes  y, D_\sqcup]_{\text{odd}} \quad \text{or in block form:}\\
  \left[  \begin{array}{cc}
V_0\circ V^\dag_1&
    V_0\circ V^\dag_0\\ V_1\circ V^\dag_1 &
    V_1\circ V^\dag_0
  \end{array}
\right]
 & -\left[\left[ 
\begin{array}{cc} A'&B'\\ C'& D'
\end{array}\right], \left[ 
\begin{array}{cc} \hat{\partial}_\sqcup &0\\ \hat{U}_\sqcup  &-\hat{\partial}_\sqcup
\end{array}\right] \right]_{\text{even}}\\
 & =\left[ \begin{array}{cc} \op{Id}_\sqcup & 0\\ 0 & \op{Id}_\sqcup
\end{array}\right]
-\left[\left[ 
\begin{array}{cc} \smz_\sqcup & 0\\ \smx& \smz_\sqcup
\end{array}\right], \left[ 
\begin{array}{cc} \hat{\partial}_\sqcup &0\\ \hat{U}_\sqcup  &-\hat{\partial}_\sqcup
\end{array}\right] \right]_{\odd}.
\end{split}
\end{equation}
for certain maps \(H'_\#\), \(\smz_\#\), from \(CM(M_\#)\) back to
itself, and maps  \(A'\), \(B'\), \(C'\), \(D'\), 
\(\smz_\sqcup, \smx\) from 
\(\hat{C}(M_\sqcup)\) back to itself. Here
\(\op{Id}_\#\), \(\op{Id}_\sqcup\) respectively denote the identity
maps from \(CM(M_\#)\) and \(\hat{C}(M_\sqcup)\) back to themselves.

The verification of the first identity (\ref{eq:cob-comp1}) involves
the cobordism \(W_\#\) obtained from composing  \(\mathcal{V}\) with
\(\bar{\mathcal{V}}\). This cobordism goes from \(M_{\#}\) to \(M_{\#}\),
and contains the circle \(\lambda_\#=\lambda\cup
\bar{\lambda}\) in its interior. A surgery along
\(\lambda_\#\) replacing a tubular neighborhood \(S^1\times B^3\) of
\(\lambda_\#\) with \(D^2\times S^2\) yields \(\bbR\times
M_{\#}\). On the other hand, to verify 
the second identity (\ref{eq:cob-comp2}), one composes in the
opposite order to get the cobordism \(W_\sqcup\) from \(M_\sqcup \) to
\(M_\sqcup \). There is an embedded 3-sphere \(S_\sqcup\subset W_\sqcup\)
obtained by joining the 3-balls \(\textsc{b}\) and
\(\bar{\textsc{b}}\) from Part 4 of the previous subsection. Doing a surgery along
\(S_\sqcup\), namely, replace a tubular neighborhood of it, \(I\times
S_\sqcup\), by a disjoint union of two 3-balls \(B_1\sqcup B_2\),
turns \(W_\sqcup\) into the product cobordism
 \(\bbR\times M_\sqcup \). One may find arcs \(\gamma _1\subset B_1\),
 \(\gamma _2\subset  B_2\) so that under this surgery they join with
 \((\lambda\cup \bar{\lambda})-I\times S_\sqcup\) to yield \(\bbR\times
 \{-p_1, p_2\}\subset\bbR\times M_\sqcup \). The cobordisms \(W_\#\)
 and \(W_\sqcup\) are equipped with metrics and \(\Spin^c\)
 structures \(\grs_{W_\#}\), \(\grs_{W_\sqcup}\) determined by the
 metric and \(\Spin^c\)-structure, \(\grs_X\), on 
 on \(X=\V\). The closed 2-form \(\omega\) on \(X\) likewise defines
 via concatenation 
 closed 2-forms \(\omega_\#\), \(\omega_\sqcup\) respectively on
 \(W_\#\) and \(W_\sqcup\). Let \(\varpi_{W_\#}:=2\omega_\#^+\),
 \(\varpi_{W_\sqcup}:=2\omega_\#^-\).

Note that like \(\mathcal{V}\) and \(\bar{\cal V}\), the composite
cobordisms \(W_\#\) and \(W_\sqcup\) also satisfy the assumption in
bullet 1 of Remark \ref{rem:2.1}. Therefore, given any
\(\Gamma_\#\), there is a unique \(\Gamma_{W_\#}\)-morphism 
which is an isomorphism from \(\Gamma_\#\) to itself. In
fact, \(\Gamma_{W_\#}=\Gamma_{\bar{\cal V}}\circ \Gamma_{\cal
  V}\). Similarly, given any
\(\Gamma_\sqcup\), there is a unique \(\Gamma_{W_\sqcup}\)-morphism from
\(\Gamma_\sqcup\) which is an isomorphism from \(\Gamma_\sqcup\) to
itself, and \(\Gamma_{W_\sqcup}=\Gamma_{\cal V}\circ \Gamma_{\bar{\cal
  V}}\). 

The proofs of (\ref{eq:cob-comp1}) and (\ref{eq:cob-comp2}) make use
of cobordism maps of the form \[
\hat{m}[\uu](W_\#, \grs_\#, \varpi_\#;
\Gamma _{W_\#}), \quad \hat{m}[\uu](W_\sqcup, \grs_\sqcup,
\varpi_{W_\sqcup}; \Gamma _{W_\sqcup}),\] as well as their
parametrized variants. (We often abbreviate these maps as
\(\hat{m}[\uu](W_\#)\), \(\hat{m}[\uu](W_\sqcup)\) below.) The
manifold \(W_\sqcup\) does not satisfy the condition
(\ref{eq:assumption-X}), but the formula for \(\hat{m}\) in
(\ref{eq:m-V}) has a straightforward adaptation in this 
context: simply replace terms of the form \(m^\#_\flat\) in
(\ref{def:m-map}) by \(m^{o\#}_{o\flat}\), and drop all the terms
\(\bar{m}^\#_\flat\). Replace \(\partial^\#_\flat\),
\(\bar{\partial}^\#_\flat\) respectively by \(\partial^o_0(M_1)\otimes
1+1\otimes \partial^\#_\flat(M_2)\) and \(1\otimes
\bar{\partial}^\#_\flat(M_2)\). 

We next describe the relevant cochains \(\uu\). 
Let \(\mathpzc{u}_{\lambda_\#}\in C^{1;\bbZ}_{\M (W_\#)}\) be the 1-cocycle 
associated to the circle
\(\lambda_\#\subset W_\#\), as defined in Section \ref{sec:A-module}'s Part
1(b).  Let \(\lambda_\sqcup\) denote the union
\(\lambda\cup\bar{\lambda}\subset W_\sqcup\), and use \(\lambda_{\sqcup -}\), \(\lambda_{\sqcup+}\) to denote
respectively the arcs \(\bar{\lambda }\), \(\lambda\) in
\(W_\sqcup\).
Let  \(\mathpzc{u}_{\lambda_{\sqcup+}}, \mathpzc{u}_{\lambda_{\sqcup-}}\in C^{1;\bbZ}_{\M (W_\sqcup)}\), be respectively the 
1-cochains 
defined in Section \ref{sec:A-module}'s Part
3(a). (The notation \(\mathpzc{u}_\lambda\),
\(\mathpzc{u}_{\bar{\lambda}}\) are usually reserved for the
1-cochains on 
\(\mathcal{B}^\sigma(\mathcal{V})\),
  \(\mathcal{B}^\sigma(\bar{\mathcal{V}})\) associated to the arcs \(\lambda\), \(\bar{\lambda}\)
in \(\mathcal{V}\), which 
appeared  previously in (\ref{eq:V}).)  
Define the 2-cochain
\(\mathpzc{u}_{\lambda_\sqcup}:=
\mathpzc{u}_{\lambda_{\sqcup-}}\mathpzc{u}_{\lambda_{\sqcup+}}\in
C^{2;\bbZ}_{\M (W_\sqcup)}\). Concretely, 
\begin{equation}\label{def:m-lambda1}
\begin{split}
 \hat{m}[\uu_{\lambda_{\sqcup-}}](W_\sqcup)& =\hat{\op{K}}_{\lambda_{\sqcup-}}(W_\sqcup)\\
& =\hat{m}[\theta_{\lambda_{\sqcup-}}](W_\sqcup)- \hat{m}[1](W_\sqcup)\, \Theta _\sqcup;\\
\hat{m}[\uu_{\lambda_{\sqcup+}}](W_\sqcup)& =\hat{\op{K}}_{\lambda_{\sqcup+}}(W_\sqcup)\\
& =\hat{m}[\theta_{\lambda_{\sqcup+}}](W_\sqcup)+\Theta _\sqcup\, \hat{m}[1](W_\sqcup);\\
\hat{m}[\uu_{\lambda_{\sqcup}}](W_\sqcup)& :=\hat{\op{K}}_{\lambda_{\sqcup}}(W_\sqcup)\\
& =\hat{m}[\theta_{\lambda_{\sqcup-}}\wedge\theta_{\lambda_{\sqcup+}}](W_\sqcup)-\hat{\op{K}}_{\lambda_{\sqcup+}}(W_\sqcup)\,
      \Theta _\sqcup +\Theta _\sqcup \, \hat{m}[\theta
  _{\lambda_{\sqcup-}}](W_\sqcup)\\
& =\hat{m}[\theta_{\lambda_{\sqcup-}}\wedge\theta_{\lambda_{\sqcup+}}](W_\sqcup)+ \Theta _\sqcup
      \,
      \hat{\op{K}}_{\lambda_{\sqcup-}}(W_\sqcup)-\hat{m}[\theta_{\lambda_{\sqcup+}}](W_\sqcup)\,
      \Theta _\sqcup , 
\end{split}
\end{equation}
where \(\Theta
_\sqcup:=1\otimes \Theta _{p_2}-1\otimes \Theta _{p_1}\). 
It will prove useful to denote the 0-cocycle \(1\in C^{0;\bbZ}_{\M (W_\sqcup)}\) on
\(\mathcal{B}^\sigma (W_\sqcup )\) by 
\(\mathpzc{u}_\emptyset \).


The proof of (\ref{eq:cob-comp1}),
(\ref{eq:cob-comp2}) involves two ingredients. The first is a set of gluing identities: 
\begin{equation}
\label{eq:ch-htpy1}
V^\dag_1V_0+V^\dag_0V_1=m[\mathpzc{u}_{\lambda_\#}](W_\#)+[
  H_\#, \partial_\# ]_{\text{even}};
\end{equation}
\BTitem\label{eq:ch-htpy}
\item[1)] The map \(V_0V^\dag_1=\hat{m}[\mathpzc{u}_{\lambda_{\sqcup-}}](
  W_\sqcup)+[A, \hat{\partial}_\sqcup]_{\text{even}}+B\hat{U}_\sqcup\).
\item[2)] The map \(V_0V^\dag_0=\hat{m}[\uu_\emptyset](
  W_\sqcup)-[\hat{\partial}_\sqcup, B]_{\text{odd}}\).
\item[3)] The map
  \(V_1V^\dag_1=\hat{m}[\mathpzc{u}_{\lambda_\sqcup}](
  W_\sqcup)+[\hat{\partial}_\sqcup, C]_{\text{odd}}-\hat{U}_\sqcup A+D\hat{U}_\sqcup \).
\item[4)] The map \(V_1V^\dag_0=\hat{m}[\mathpzc{u}_{\lambda_{\sqcup+}}](
  W_\sqcup)+[\hat{\partial}_\sqcup, D]_{\text{even}}-\hat{U}_\sqcup  B\).
\ETitem
The definition of the maps \(H_\#\), \(A\), \(B\), \(C\), \(D\), and
the verification of 
these identitites  occupy the remainder of this step and Step 4
below. 
In short, they all follow from an adaption of \cite{KM}'s Lemma 26.2.2,
together with a  parametrized variant of the identity (\ref{eq:K-lambda}). Rephrased in our
language, the composition identity in \cite{KM}, which was stated for
the check-version of monopole Floer homology, has the following
companion 
version in for the hat-version: Let \(W_1\)
be a connected cobordism from \(Y_-\) to \(Y_0\), and \(W_2\) a
connected cobordism from
\(Y_0\) to \(Y_+\). Let \(W=W_2\circ W_1\) denote the composite
cobordism of \(W_1\) and \(W_2\).  For $\uu_1\in \op{C}(\B^\sigma
((W_1)_c);\bbK)$ and $\uu_2\in \op{C}(\B^\sigma
((W_2)_c);\bbK)$, \cite{KM} defined an ``inner product''
of \(\uu_1\) and \(\uu_2\), denoted $\uu:=c (\uu_1\otimes
\uu_2)\in \op{C}(\B^\sigma (W_c);\bbK)$ (cf. \cite{KM}'s Equation
(26.9)). We have: 
\begin{equation}\label{KM-composition0}
\hat{m}[\uu_2](W_2)\, \hat{m}[\uu_1](W_1)=\hat{m}[\uu](W)+
[\hat{\textsc{k}}[\mathbf{u}](W), \hat{\partial}]+\hat{\textsc{k}}[\delta\mathbf{u}](W),
\end{equation}
where the maps  $\hat{\smk}$ are defined via integrations on a certain
parametrized
moduli space, and $\mathbf{u}$ is a parametrized version of
$\uu$. (Though cobordism maps $\mathring{m}[\uu](W)$ are previously
defined for cochains on \(\B^\sigma _{loc}(W)\) instead of those on
$\B^\sigma (W_c)$, there is a restriction map, $\ss$, from an open dense
subset of the former to the latter. As explained in \cite{KM},
because of unique continuation, it makes no practical difference to
work with either $\B^\sigma (W_c)$ or \(\B^\sigma _{loc}(W)\), or the
aforementioned open dense subset of \(\B^\sigma _{loc}(W)\). 
The identity (\ref{KM-composition0})  is the consequence of applying a 
Stokes' theorem to the compactification of the aforementioned
parametrized moduli space.  See (26.2), (26.3) in
\cite{KM} for the definition of the aforementioned parametrized moduli
space; Equations (26.11), (26.12) therein for the definition of  the
associated maps $\hat{\smk}$ (denoted \(\check{K}\) in \cite{KM});
Lemma 26.2.2 and its siblings in \cite{KM} for proofs of the key gluing
identities. 

Roughly speaking, the proof of (\ref{eq:ch-htpy1}) and
(\ref{eq:ch-htpy}) follow from applying variants of (\ref{KM-composition0}) to
$W=\bar{\V}\circ\V$, and $W=\V\circ \bar{\V}$ respectively, with $\uu$
taken to be $\uu_{\lambda _\#}$ in (\ref{eq:ch-htpy1}), and 
$\uu$ set to be 
$\uu_{\lambda_{\sqcup -}}, \uu_{\emptyset}, \uu_{\lambda_{\sqcup }},
\uu_{\lambda_{\sqcup +}}$ respectively in items 1)--4) of
(\ref{eq:ch-htpy}). The maps \(H_\#\), \(A\), \(B\), \(C\), \(D\) are
then 
given by 
\begin{equation}\label{def:H-D0}
\begin{split}
H_\#& =\hat{\smk}[\mathbf{u}_{\lambda _\#}](W_\#), \\
A& =\hat{\smk}[\mathbf{u}_{\lambda
  _{\sqcup -}}](W_\sqcup),\\ 
-B& =\hat{\smk}[\mathbf{u}_{\emptyset}](W_\sqcup)=\hat{\smk}[1](W_\sqcup),\\
 C& =\hat{\smk}[\mathbf{u}_{\lambda
  _{\sqcup }}](W_\sqcup),\\
 -D& =\hat{\smk}[\mathbf{u}_{\lambda _{\sqcup +}}](W_\sqcup).
\end{split}
\end{equation}
A couple of issues need to be addressed to be able to apply
(\ref{KM-composition0}) in our settting. 
Firstly, in \cite{KM}, \(Y_\pm\), \(Y_0\) are assumed to be
connected. As previously explained, there is no problem adapting
to the case when $Y_\pm$ is the disconnected manifold $M_\sqcup$. In
the case of (\ref{eq:ch-htpy1}), $W_\#=\bar{\V}\circ \V$ is glued
along $Y_0=M_\sqcup$. This creates no new troubles: The assumption that only the
\(M_2\) component of $M_\sqcup $ can be associated with balanced perturbations imply
that the straightforward sort of gluing argument applies with gluing 
along \(M_1\), leaving the more delicate analysis described in
\cite{KM} required for \(M_2\) alone. 
The second issue is related to the fact that, recalling the discussion
in Section \ref{sec:A-module}, the cochains \(\uu_\gamma
\) and their associated maps \(\hat{m}[\uu_\gamma ]\),
\(\hat{\smk}[\mathbf{u}_\gamma ]\) relevant to our discussion are of
a more general sort.  In particular, when \(\gamma \) is noncompact,
unlike those cochains on \(\B^\sigma (W_c)\) considered in \cite{KM},
our \(\uu_\gamma \in \op{C}(\B^\sigma _{loc}(W);\bbK)\) are sensitive
to the behavior of the Seiberg-Witten configurations over the ends
\(W\backslash W_c\). To explain this issue in more detail, as well
as to describe the modification to generalize
(\ref{KM-composition0}) to this context, some preliminary discussions
are required. 

Here are some key ingredients of \cite{KM}'s derivation of
(\ref{KM-composition0}). 
Let \(W(S)_c\) denote the variants of \cite{KM}'s
\(W(S)\) (cf. \cite{KM}'s (26.2) and thereabouts). We write it as: 
\begin{equation}\label{def:W(S)}
W(S)_c=(W_1)_c\cup ([-S/2,S/2]\times Y_0)\cup (W_2)_c. 
\end{equation}
Let \(W(S)\) be the (complete) manifold with
cylindrical ends containing \(W(S)_c\) as its ``compact piece''. (Recall the notation from
Section \ref{sec:2.2}; they were denoted \(W(S)^*\) in
\cite{KM}.) For example, the cobordisms \(W_\#(S)\) and \(W_\sqcup (S)\) are illustrated respectively in Figure
\ref{fig:Aa} and Figure \ref{fig:Ab} below, where the shaded regions
represent the ``necks'' of length \(S\). 
The parametrized moduli spaces involved in the proof of
(\ref{KM-composition0}) are of the following sort: 
\[\begin{split}
& \mathbf{M}_{k+1,z}(W, \grc_-, \grc_+) :=\bigcup_{S\in [0, \infty)}
\{S\}\times\mathcal{M}_{k,z}(W(S), \grc_-, \grc_+), \\
& \quad \mathbf{M}_{k+1}(W, \grc_-, \grc_+)
:=\bigcup_{z}\mathbf{M}_{k+1,z}(W, \grc_-, \grc_+);\\
& \mathbf{M}_{k+1,z}^+(W, \grc_-, \grc_+):=\bigcup_{S\in [0, \infty]}
\{S\}\times\mathcal{M}_{k,z}^+(W(S), \grc_-, \grc_+),\\
& \quad \mathbf{M}_{k+1}^+(W, \grc_-, \grc_+) :=\bigcup_{z}\mathbf{M}_{k+1,z}^+(W, \grc_-, \grc_+),\\
\end{split}
\]
with the ``fiber at $\infty$'', $\mathcal{M}^+_{k,z}(W(\infty), \grc_-,
\grc_+)=\mathcal{M}_{k,z}(W(\infty), \grc_-,
\grc_+)$, given in \cite{KM}'s (26.4). 
Their reducible variants are defined similarly. 
The compactified moduli space $\mathbf{M}_{k+1,z}^+(W,
\grc_-, \grc_+)$ maps to a smaller compactification,
$\bar{\mathbf{M}}_{k+1,z}(W, \grc_-, \grc_+)$ embedded in 
\begin{equation}\label{par-B}
[0, \infty]\times \B^\sigma ((W_1)_c)\times \B^\sigma ((W_2)_c),
\end{equation}
in a way similar to the map $\grr$ in Section
\ref{sec:2.4}. Cf. \cite{KM}'s (26.7). 
This map preserves the fibration (over \([0,\infty]\)) structure on
both spaces, and 
over the fiber $\{S\}\times \M_{k,z}^+(W(S),
\grc_-, \grc_+)\subset \mathbf{M}_{k+1,z}^+(W,
\grc_-, \grc_+)$, $S\in [0,\infty)$, this map factors through 
\begin{equation}\label{def:ss}
\begin{split}\M_{k,z}^+(W(S),
\grc_-, \grc_+)\stackrel{\grr}{\to } & \B^\sigma
_{loc}((W(S))^\circ\stackrel{\ss}{\to }\B^\sigma ((W(S))_c)^{\circ}\\
& \stackrel{\ss_1\times \ss_2}{\longrightarrow} \B^\sigma ((W_1)_c)\times \B^\sigma ((W_2)_c),
\end{split}
\end{equation}
where \(\grr\) is as in Section \ref{sec:2.4}, and for each \(i=1,2\),
\(\ss_i\co \B^\sigma((W(S))_c)^\circ\to \B^\sigma ((W_1)_c)\) denotes the map of restricting
to \((W_i)_c\subset W(S)\). Here, \(\B^\sigma
_{loc}((W(S))^\circ\) denotes a certain open dense subset of \(\B^\sigma
_{loc}((W(S))\); similarly for \(\B^\sigma (W_c)^{\circ}\). 
The cochains \(\uu_1\in \op{C}(\B^\sigma ((W_i)_c)\) from
(\ref{KM-composition0}) thus defines a cochain
\[
\uu= c(\uu_1\otimes \uu_2) :=(\ss_1\times \ss_2)^*(\uu_1\times \uu_2)
\]
in \(\op{C}(\B^\sigma (W_c))\), and \(\mathbf{u}\) in
\((\ref{KM-composition0})\) refers to the cochain on (\ref{par-B})
induced from \(\uu_1\) and \(\uu_2\).

Use \(\pmb{\grr}\) to denote the
aforementioned map from \(\mathbf{M}_{k+1}^+(W, \grc_-, \grc_+)\) to
(\ref{par-B}), and \(\pmb{\varsigma }\) for the embedding of \(
\bar{\mathbf{M}}_{k+1}(W, \grc_-, \grc_+)\) into
(\ref{par-B}). Use \(\pmb{\grr}|_S\), \(\pmb{\varsigma }|_S\) respectively to denote
the restriction of \(\pmb{\grr}\), \(\pmb{\varsigma }\) from the fiber
over \(S\) of 
\(\mathbf{M}_{k+1}^+(W, \grc_-, \grc_+)\to[0,\infty]\) or \(
\bar{\mathbf{M}}_{k+1}(W, \grc_-, \grc_+)\to[0,\infty]\) to
\(\B^\sigma ((W_1)_c)\times \B^\sigma ((W_2)_c)\).

The map \(\hat{\smk}[\mathbf{u}]\) is constructed from
\(\smk[\mathbf{u}]^\#_\flat\) and \(\bar{\smk}[\mathbf{u}]^\#_\flat\),
where \(\smk[\mathbf{u}]^\#_\flat\) is a sum
of terms with coefficients taking the form 
\begin{equation}\label{def:K[u]}
\langle\pmb{\grr}^* \mathbf{u},
\mathbf{M}_{k+1, z}^+(W, \grc_-, \grc_+)\rangle=\langle\pmb{\varsigma }^* \mathbf{u},
\mathbf{M}_{k+1, z}(W, \grc_-, \grc_+)\rangle,
\end{equation}
where \(k\) is the
degree of \(\uu\). Similarly for
\(\bar{\smk}[\mathbf{u}]^\#_\flat\). It is shown in \cite{KM}'s
Proposition 26.1.6 that for any \(k\), \(z\), \(\mathbf{M}_{k+1, z}^+(W, \grc_-, \grc_+)\)
is a stratified manifold where Stokes' theorem (in the sense of
\cite{KM}'s Lemma 21.3.1) is applicable. Equation (\ref{KM-composition0}) is
then the consequence of applying this Stokes' theorem to integrals of
the form 
\begin{equation}\label{par-stokes}
\langle \pmb{\grr}^*(\delta \mathbf{u}), \mathbf{M}_{k+1, z}^+(W,
\grc_-, \grc_+)\rangle=\langle \pmb{\grr}^*\mathbf{u}, \partial\, [\mathbf{M}_{k+1, z}^+(W,
\grc_-, \grc_+)]\rangle, 
\end{equation}
together with an analysis of the structure of \(\big(\mathbf{M}_{k+1, z}^+(W,
\grc_-, \grc_+)\big)_{k}\). That is, the codimension one stratified
submanifold, \(\big(\mathbf{M}_{k+1, z}^+(W,
\grc_-, \grc_+)\big)_{k}\subset \mathbf{M}_{k+1, z}^+(W,
\grc_-, \grc_+)\), is described as a union of the following form:
\begin{equation}\label{bdry-bfM}
\begin{split}
& (\mathbf{M}_{k+1,z}^+(W, \grc_-, \grc_+))_k
=\big( \{\infty\}\times \mathcal{M}_{k,z}^+(W(\infty), \grc_-,
\grc_+)\big) \\
& \qquad \cup \big( \{0\}\times \mathcal{M}_{k,z}^+(W(0), \grc_-, \grc_+)\big)
 \cup  \bigcup_{S\in (0, \infty)}
\{S\}\times\big(\mathcal{M}_{k,z}^+(W_\#(S), \grc_-, \grc_+)\big)_{k-1}.
\end{split}
\end{equation}
The first two terms on the right hand side of (\ref{bdry-bfM}) contribute
respectively
\begin{equation}\label{M-bdry-v}
\begin{split}
& \langle (\pmb{\grr}|_\infty)^*(\uu_1\times \uu_2), \mathcal{M}_{k,z}^+(W(\infty), \grc_-,
\grc_+)\rangle \quad \text{and} \\ 
& -\langle (\pmb{\grr}|_0)^*(\uu_1\times \uu_2), \mathcal{M}_{k,z}^+(W(0), \grc_-,
\grc_+)\rangle=-\langle\grr^* c(\uu_1\otimes\uu_2), \mathcal{M}_{k,z}^+(W, \grc_-,
\grc_+)\rangle 
\end{split}
\end{equation}
to the right hand side of (\ref{par-stokes}),
resulting respectively in the left hand side of
(\ref{KM-composition0}) and the first term on the right hand side of
(\ref{KM-composition0}). The last term in (\ref{KM-composition0})
arises from the left hand side of (\ref{par-stokes}). The contribution
from the last term of (\ref{bdry-bfM}) to the right hand side of
(\ref{par-stokes}) leads to the penultimate term in
(\ref{KM-composition0}), based on the straightforward adaptation of
\cite{KM}'s Proposition 25.3.4 to the parametrized context.

A simple reformulation of  \cite{KM}'s work suffices to make
(\ref{KM-composition0}) applicable to general \(\uu_i\in
\op{C}(\B^\sigma _{loc}(W_i);\bbK)\). 
Let \(W_1(S)\supset (W_1)_c\), \(W_2(S)\supset (W_2)_c\) be (the
closure of) the two
halves of \(W(S)\) when divided in the middle of ``the neck'', namely,
at the 3-manifold \(\{0\}\times Y_0\)  in (\ref{def:W(S)}). For \(i=1,2\), define
\(W_i(\infty)\) to be the previously introduced complete manifold, \(W_i\).  Instead of 
(\ref{par-B}), consider another space fibering
over \([0,\infty]\), whose fiber over \(S\in [0,\infty]\) is
\[
\B^\sigma (W_1(S))\times \B^\sigma (W_2(S))=: \B^\sigma (W_2\circ_S W_1),
\] 
where \(\B^\sigma
(W_i(S))\), \(i=1,2\), are both equipped with the topology inherited
from its embedding to \(\B^\sigma _{loc(W_i)}\). For all \(S\in
[0,\infty]\), \(\B^\sigma
(W_1(S))\) admits a well-defined (\(-\infty\))-limit map by construction, 
\(\Pi^{-\infty}:=\Pi^{-\infty}_{W_1}\co \B^\sigma
(W_1(S))\to \B^\sigma (Y_-)\); likewise, \(\B^\sigma
(W_2(S))\) have  a well-defined (\(+\infty\))-limit map,
\(\Pi^{\infty}:=\Pi^{\infty}_{W_2}\co \B^\sigma
(W_2(S))\to \B^\sigma (Y_+)\). (Recall that these maps  played in
important roles 
in the construction of the cochains in Section
\ref{sec:A-module}. These are not available with the spaces
\(\B^\sigma ((W_i)_c)\) used in \cite{KM}.) 
Denote the fibered space 
\[
 \pmb{\B}^\sigma (W_2\circ W_1):=\bigcup_{S\in [0, \infty]}\{S\}\times
 \B^\sigma (W_2\circ_S W_1). 
\]
(This space is homeomorphic to (\ref{par-B}) if endowed with the
stronger Banach topology.) 
When \(S\) is finite, let \(\Pi^{Y_0}_{W_*}\co \B^\sigma (W_*(S))^\circ\to
\B^\sigma (Y_0)\) denote the map of restricting to 3-manifold
\(\{0\}\times Y_0\subset W_*(S)\) for \(W_*=W, W_1, W_2\). (Again, the
superscript \(\circ\) is used to denote an appropriate open dense
subspace. It is sometimes dropped to make the notations less
cumbersome. As previously mentioned, this makes no practical
difference.)
As was done in Section \ref{sec:2}, when \(S=\infty\), let \(\Pi^{Y_0}_{W_i}\co \B^\sigma (W_i)^\circ\to
\B^\sigma (Y_0)\) denote the map of taking \((+\infty)\)-limits for
\(i=1\), and that of taking the \((-\infty)\)-limits for
\(i=2\). Slightly abusing notation, we now let \(\ss_i\co \B^\sigma
(W(S))^\circ\to\B^\sigma (W_i(S))\), \(i=1,2\), denote the map of
restricting to \(W_i(S)\subset W(S)\). Equation (\ref{def:ss}) has
straightforward analog here: For finite \(S\), the map
\(\ss_1\times \ss_2\) factors as: 
\begin{equation}\label{dec:B(W)}
\B^\sigma (W(S))^\circ\stackrel{\ss_1\times \ss_2}{\longrightarrow}\B^\sigma  (W_1(S))\times _{\B^\sigma (Y_0)}\B^\sigma
(W_2(S))\hookrightarrow \B^\sigma  (W_1(S))\times \B^\sigma(W_2(S)),
\end{equation}
In the above, the fiber product \(\B^\sigma  (W_1(S))\times _{\B^\sigma (Y_0)}\B^\sigma
(W_2(S))\) is regarded as a subspace of the product \(\B^\sigma
(W_1(S))\times \B^\sigma(W_2(S))\) where the maps \(\Pi^{Y_0}_{W_i}\co \B^\sigma  (W_i)\to
\B(Y_0)\), \(i=1,2\) take the same value. The previously introduced maps \(\pmb{\grr}\),
\(\pmb{\varsigma }\), \(\pmb{\grr}|_S\), \(\pmb{\varsigma }|_S\), as
well as the inner products \(c(\uu_1\otimes \uu_2)\), also admit
straightforward adaptions, which we denote by the same notation. For
finite \(S\),  the restriction of the maps \(\pmb{\grr}|_S\),
\(\pmb{\varsigma }|_S\) to to \(\mathcal{M}_{k,z}(W(S), \grc_-,
\grc_+)\) are
respectively the composition of \(\grr\) and \(\varsigma \) with
(\ref{dec:B(W)}). The arguments of \cite{KM} still apply with this
modification to establish (\ref{KM-composition0}) in the context of more general \(\uu_i\in
\op{C}(\B^\sigma _{loc}(W_i);\bbK)\). 

Even in this (slightly) generalized form, specific applications of
(\ref{KM-composition0}) in our context are often complicated by the
fact that the cochains \(\uu\) constructed in Section \ref{sec:A-module} do not
necessarily take the form of, or have no obvious interpretation as,
an inner product \(c(\uu_1\otimes \uu_2)\). 
In view of the roles played by bundles over \(\B^\sigma  (W)\) (such that \(\B^\sigma  _x(W)\) or \(\B^\sigma _\lambda
(W)\)) in the construction of the cochains from Section
\ref{sec:A-module}, we typically deal this problem
by going through various bundles over \(\pmb{\B}^\sigma (W_2\circ
W_1)\).  They are constructed in a manner similar to what was done in
Section \ref{sec:A-module}. For example, what will be called \(\tilde{\pmb{\B}}_x^\sigma (W_2\circ
W_1)\) is defined as follows. 

Take a point \(x\in W_c\) that lies in the 3-manifold
\(Y_0\subset W_c\) that separates \((W_1)_c\) and \((W_2)_c\). (The point \(x\in W\) is denoted \(\ul{x}\) when
regarded as a point in the 3-manifold \(Y_0\).) Recall the \(U(1)\)-bundles \(\pi_{\ul{x}}\co
\tilde{\B}^\sigma _{\ul{x}}(Y_0)\to \B^\sigma (Y_0)\),
\(\pi_x\co \tilde{\B}_x^\sigma (W)\to \B^\sigma  (W)\) from Section
\ref{sec:A-module}. 
For finite \(S\), the map
\(\Pi^{Y_0}_{W}\)  lifts to a map
\(\tilde{\Pi}^{Y_0}_{W}\co \tilde{\B}^\sigma  _x(W(S))\to
\tilde{\B}_{\ul{x}}(Y_0)^\sigma \) by construction. Meanwhile, for \(i=1,2\)
and any \(S\in [0,\infty]\), one may define 
\(\tilde{\Pi}^{Y_0}_{W_i}\), \(\pi^{W_i}_{x}\), \(\tilde{\B}^\sigma
_x(W_i)\) through the following commutative diagram: 
\begin{equation}\label{CD:Pi_i}
\begin{CD}
\tilde{\cal B}^\sigma _x(W_i(S)) @>\tilde{\Pi}^{Y_0}_{W_i}>>  \tilde{\cal B}^\sigma _{\ul{x}}(Y_0) \\
@V\pi^{W_i}_{x}VV @V\pi_{\ul{x}}VV\\
{\cal B}^\sigma (W_i(S)) @>\Pi^{Y_0}_{W_i}>>  {\cal B}^\sigma (Y_0). 
\end{CD}
\end{equation}
Now let 
\[
\begin{split}
\tilde{\B}^\sigma _x(W_2\circ_S W_1) & :=\tilde{\B}^\sigma
_x(W_1(S))\times \tilde{\B}^\sigma _x(W_2(S)); \quad\\
\tilde{\pmb{\B}}^\sigma  _x(W_2\circ W_1) &:=\bigcup_{S\in [0,
  \infty]}\{S\}\times \tilde{\B}^\sigma _x(W_2\circ_S W_1). 
\end{split}
\]
By construction, these are \(U(1)\times U(1)\)-bundles respectively over \(\B^\sigma
(W_2\circ_S W_1)\) and \(\pmb{\B}^\sigma  (W_2\circ W_1)\). The
fibered product \(\tilde{\B}^\sigma
_x(W_1(S))\times _{\tilde{\B}_x(Y_0)}\tilde{\B}^\sigma _x(W_2(S))\),
as a subspace of \(\tilde{\B}^\sigma
_x(W_1(S))\times \tilde{\B}^\sigma _x(W_2(S))\), is preserved under the
diagonal \(U(1)\)-action. The quotienting by this action is: 
\[\begin{split}
& \big( \tilde{\B}^\sigma
_x(W_1(S))\times _{\tilde{\B}_x(Y_0)}\tilde{\B}^\sigma _x(W_2(S))\big)
/U(1)_\Delta  \simeq \B^\sigma  (W_1(S))\times _{\B^\sigma
  (Y_0)}\B^\sigma(W_2(S))\\ 
& \qquad \hookrightarrow \B^\sigma  (W_1(S))\times \B^\sigma(W_2(S))=:\B^\sigma
(W_2\circ_S W_1),
\end{split}
\]
where \(U(1)_\Delta \subset U(1)\times U(1)\) denotes the
diagonal. 

The previously introduced 
restriction maps \(\ss_i \co \B^\sigma (W(S))^\circ\to
\B^\sigma ((W_i(S))\) lift to define maps \(\tilde{\ss}_i\co \tilde{\B}^\sigma _x(W(S))^\circ\to
\tilde{\B}^\sigma _x((W_i(S))\) for \(i=1,2\). 
With them we have the following variant of (\ref{dec:B(W)}) for {\em finite}
\(S\): 
\begin{equation}\label{dec:B_x(W)}
\begin{split}
& \tilde{\B}_x^\sigma (W(S))\stackrel{\tilde{\ss}_1\times \tilde{\ss}_2}{\longrightarrow}\tilde{\B}^\sigma  _x(W_1(S))\times _{\tilde{\B}_{\ul{x}}(Y_0)}\tilde{\B}^\sigma
_x(W_2(S))\\ &\qquad \quad \hookrightarrow \tilde{\B}^\sigma  _x(W_1(S))\times \tilde{\B}^\sigma
_x(W_2(S)):=\tilde{\B}^\sigma  _x(W_2\circ_S W_1);
\end{split}
\end{equation}
and together they form a commutative diagram:
\begin{equation}\label{CD:W-S}
\begin{CD}
 \tilde{\B}_x^\sigma (W(S))@>\tilde{\ss}_1\times \tilde{\ss}_2>> \tilde{\B}^\sigma  _x(W_1(S))\times _{\tilde{\B}_{\ul{x}}(Y_0)}\tilde{\B}^\sigma
_x(W_2(S)) @>\text{embeds}>> \tilde{\B}^\sigma  _x(W_2\circ_S W_1)\\
@V\pi^{W}_{x}VV @V\pi^{W_1}_{x}\times \pi^{W_2}_{x}VV @V\pi_x^{W_2\circ_S W_1}V:=\pi^{W_1}_{x}\times \pi^{W_2}_{x}V\\
\B^\sigma (W(S)) @>\ss_1\times \ss_2>> \B^\sigma  (W_1(S))\times _{\B^\sigma (Y_0)}\B^\sigma
(W_2(S))@>\text{embeds}>> \B^\sigma  (W_2\circ_S W_1);
\end{CD}
\end{equation}
the pair of horizontal maps \(\tilde{\ss}_1\times \tilde{\ss}_2\) and
\(\ss_1\times \ss_2\) in the left square above define a map between
\(U(1)\)-bundles (but not the right square), and the map
\((\tilde{\Pi}^{Y_0}_W, \Pi ^{Y_0}_W)\) between the \(U(1)\)-bundles
\(\pi_x^W\co \tilde{\B}_x^\sigma (W(S))\to \B^\sigma (W(S))\) and
\(\pi_{\ul{x}}\co \tilde{\cal B}^\sigma _{\ul{x}}(Y_0)\to {\cal
  B}^\sigma (Y_0)\) factors through this bundle map: 
\begin{equation}\label{CD:W-W_i}
\xymatrixcolsep{2.5pc}
\xymatrix{
\tilde{\cal B}^\sigma
_x(W(S))\ar@{->}[r]^{\pi^{W}_{x}}\ar@{->}[d]^{\tilde{\ss}_1\times \tilde{\ss}_2}
\ar@/^-6pc/[dd]_{\tilde{\Pi}^{Y_0}_{W}} 
& {\cal B}^\sigma (W(S))
\ar@{->}[d]^{\ss_1\times \ss_2}\ar@/^6pc/[dd]^{\Pi^{Y_0}_{W}}\\
\tilde{\B}^\sigma  _x(W_1(S))\times _{\tilde{\B}_{\ul{x}}(Y_0)}\tilde{\B}^\sigma
_x(W_2(S)) \ar@{->}[r]^{\pi^{W_1}_{x}\times \pi^{W_2}_{x}}
\ar@{->}[d]^{\tilde{\Pi}^{Y_0}_{W_1}=\tilde{\Pi}^{Y_0}_{W_2}} 
& \B^\sigma  (W_1(S))\times _{\B^\sigma (Y_0)}\B^\sigma
(W_2(S)) \ar@{->}[d]^{\Pi^{Y_0}_{W_1}=\Pi^{Y_0}_{W_2}}\\
\tilde{\cal B}^\sigma _{\ul{x}}(Y_0) 
\ar@{->}[r]^{\pi_{\ul{x}}}
& {\cal B}^\sigma (Y_0). 
}
\end{equation}

As a general rule, in what follows we adopt the convention of adding subscripts or
superscripts \(W\) in notations previously introduced in Section
\ref{sec:2} to emphasize the cobordism referred to. For example,
\(\Pi^{\pm \infty}_W\) denote the version of (\(\pm\infty\))-limit map
\(\Pi^{\pm \infty}\) for the cobordism \(W\), and \(\pi_x^W\) is the
\(W\)'s version of the projection map \(\pi_x\) in Section
\ref{sec:2}.
 
We shall also make use other variants of the bundle \(\pi_x^{W_2\circ W_1}\co \tilde{\pmb{\B}}^\sigma
_x(W_2\circ W_1)\to \pmb{\B}^\sigma(W_2\circ W_1)\). These are
constructed in a similar fashion, with the role of \(\tilde{\B}^\sigma
_x(Y_0)\) replaced by other bundles over \(\B^\sigma (Y_0)\), say \(\tilde{\B}^\sigma_\nu(Y_0)\), \(\nu \) being a
\(0\)-chain in \(Y_0\).  The composition formula (\ref{KM-composition0}) is
verified for a cochain \(\uu_\gamma \) from Section \ref{sec:A-module}
by applying the trick already used repeatedly in Section \ref{sec:2},
cf. e.g. the diagrams (\ref{CD-M+}), (\ref{CD-M+1}). That is, we
choose an appropriate lift of the embedding \(\pmb{\varsigma }\co
\bar{\bf M} \to \pmb{\B}^\sigma (W_2\circ W_1)\) to \(\tilde{\pmb{\varsigma }}\co
{\bf M}^+ \to \tilde{\pmb{\B}}^\sigma  (W_2\circ W_1)\) that fit in a
commutative diagram of the form
\begin{equation}\label{CD-bfM}\begin{CD}
\mathbf{M}^+@>\tilde{\pmb{\varsigma}} >>\tilde{\pmb{\B}}^\sigma(W_2\circ W_1)\\
@ V\pmb{\grr }VV @ V\pi VV\\
\bar{\bf M}  @>\pmb{\varsigma}>> \pmb{\B}^\sigma  (W_2\circ W_1)
\end{CD}
\end{equation} 
and apply Stokes' theorem over the top row of the preceding
diagram. The cochain \(\uu\in \op{C}(\B^\sigma (W);\bbK)\) is
typically interpreted in terms of inner products by considering a variant of the diagram (\ref{CD:W-S}) for
finite \(S\). Once in the inner product form, the cochain \(\uu\)
extends to be defined on the fiber at infinity, \(\M^+(W(\infty))\),
and consequently also a cochain \(\mathbf{u}\) suitable for applying the arguments
of (\ref{KM-composition0}). As outlined previously, the term on the
left hand side of (\ref{KM-composition0}), \(\hat{m}[\uu_2](W_2)\, \hat{m}[\uu_1](W_1)\), arises from integrals over strata in
\(\M^+(W(\infty))\). To put the integrals in a suitable product form, we must
factor the strata of \(\M^+(W(\infty), \grc_-, \grc_+)\) as products of two spaces,
provisionally written as \(\M^+_{W_1}(\grc_-, \grc)\times
\M^+_{W_2}(\grc, \grc_+)\), with \(\grc\in \B^\sigma (Y)\). Here, 
\(\M^+_{W_1}(\grc_-, \grc)\) consists of ``broken \(W_1\)-paths'' from
\(\grc_-\) to \(\grc\), and \(\M^+_{W_2}(\grc, \grc_+)\) consists of
``broken  \(W_2\)-paths'' from 
\(\grc\) to \(\grc_+\). Recall from \cite{KM}'s Equation (26.4) that a general
element of \(\M^+(W(\infty))\) is defined to be an element in a
product space, 
\begin{equation}\label{def:W-infty}
\begin{split}
& \N^+(Y_-, \grc_-, \grc'_-)\times \M (W_1, \grc'_-,
\grc_{0-})\times\N^+(Y_0, \grc_{0-}, \grc_{0+})\\
& \qquad \quad \qquad \times \M (W_2,
\grc_{0+}, \grc'_+)\times\N^+(Y_+, \grc'_+, \grc_+).
\end{split}
\end{equation}
There are different ways of organizing the preceding space in the form
of \(\M^+_{W_1}(\grc_-, \grc)\times
\M^+_{W_2}(\grc, \grc_+)\). When deriving the hat-version of
composition formula, we take the \(\M^+_{W_1}(\grc_-, \grc)\),
\(\M^+_{W_2}(\grc, \grc_+)\) respectively to be the first and the
second line of the preceding expression. (For the check-version, one
takes \(\M^+_{W_1}(\grc_-, \grc)\) to be the product of the first
two factors in (\ref{def:W-infty}), and \(\M^+_{W_2}(\grc, \grc_+)\)
the product of the remaining factors.) Applying \cite{KM}'s Proposition
26.1.6 to write out each entry of the identity from Stokes'
theorem as a sum in the manner of \cite{KM}'s (26.13) leads to a
variant of the composition formula (\ref{KM-composition0}). (But with our generalized definition of the cobordism maps).
 \begin{equation}\label{KM-composition}
\hat{m}[\uu](W(\infty))=\hat{m}[\uu](W)+
[\hat{\textsc{k}}[\mathbf{u}](W), \hat{\partial}]+\hat{\textsc{k}}[\delta\mathbf{u}](W).
\end{equation}
Depending on how \(\uu\) is expressed in terms of inner products, the
left hand side, \(\hat{m}[\uu](W(\infty))\), will expressed in terms of
products of the form \(\hat{m}[\uu_2](W_2)\, \hat{m}[\uu_1](W_1)\). 
 Note however that, compared with the long sum in the expression following \cite{KM}'s
(26.13), in our case there will be additional terms involving boundary
obstructed maps of the form \(\bar{n}^s_u[\op{u}](Y_0)\) (including the
boundary-obstructed differentials   \(\bar{\partial}^s_u(Y_0)\) that
appears in \cite{KM}'s formula). These additional terms are absorbed
in our generalized definition of \(\hat{m}[\uu]\) and
\(\hat{\smk}[\mathbf{u}]\). 
 \begin{equation}\label{KM-composition}
\hat{m}[\uu](W(\infty))=\hat{m}[\uu](W)+
[\hat{\textsc{k}}[\mathbf{u}](W), \hat{\partial}]+\hat{\textsc{k}}[\delta\mathbf{u}](W).
\end{equation}

\paragraph{\it Step 4.} 
We now apply the general discussion
above to derive the identities (\ref{eq:ch-htpy1}) and (\ref{eq:ch-htpy}). 
What follows describes the formulas in (\ref{def:H-D0}) in a more
explicit manner. The degree \(-1\) map 
 \(H_\#\co C_\#\to C_\#\) is given by: 
\[
H_\#=\sum_{\grc_1, \grc_2\in \grC_\#}\sum_{z\in
  \pi_0(\mathcal{B}^\sigma (M_\#, \grc_1, \grc_2))}
\langle\pmb{\varsigma }^*\mathbf{u}_{\lambda_\#}, \mathbf{M}_{1,z}(W_\#, \grc_1,
\grc_2)\rangle \Gamma_\# (z).\]
For practical purposes, it is often more convenient to work with the
more concrete variant of \(H_\#\): this is denoted by \(\dot{H}_\#=\hat{\smk}[\pmb{\theta }_{\lambda _\#}](W_\#)\)
below and is defined by replacing
\(\mathbf{u}_{\lambda_\#}\) in the preceding formula by \(\pmb{\theta
}_{\lambda _\#}\). 
The two maps \(H_\#\), \(\dot{H}_\#\) are related by 
\[
H_\#=\dot{H}_\#+\big[\hat{\partial}_\#, \hat{\smk}[\pmb{\varepsilon }_{\lambda _\#}]\big],
\]
where \(\pmb{\varepsilon }_{\lambda _\#}\) is the parametrized variant
of the 0-cochain \(\varepsilon _{\lambda _\#}\) defined in Section
\ref{sec:A-module}'s Part 1(b).

Likewise, the maps \(A, B, C, D\co \hat{C}_\sqcup\to \hat{C}_\sqcup\),
are assembled respectively from the constituents
\(A^{o\sharp}_{o\flat}, B^{o\sharp}_{o\flat}, C^{o\sharp}_{o\flat},
D^{o\sharp}_{o\flat}\co C_\sqcup ^{o\sharp}\to  C_\sqcup^{o\flat}\):
\begin{equation}\label{def:A-cons}
\begin{split}
A^{o\sharp}_{o\flat} & =\sum _{\grc_1\in \grC^{o\sharp}_\sqcup}\sum_{\grc_2\in \grC^{o\flat}_\sqcup}\sum_{z\in
  \pi_0(\mathcal{B}^\sigma (M_\sqcup, \grc_1, \grc_2))}
\langle\pmb{\varsigma }^*\mathbf{u}_{\lambda_\sqcup -}, \mathbf{M}_{1,z}(W_\sqcup, \grc_1,
\grc_2)\rangle \Gamma_\sqcup (z), \\
B^{o\sharp}_{o\flat} & =\sum _{\grc_1\in \grC^{o\sharp}_\sqcup}\sum_{\grc_2\in \grC^{o\flat}_\sqcup}\sum_{z\in
  \pi_0(\mathcal{B}^\sigma (M_\sqcup, \grc_1, \grc_2))}
\langle\pmb{\varsigma }^*\mathbf{u}_\emptyset, \mathbf{M}_{0,z}(W_\sqcup, \grc_1,
\grc_2)\rangle \Gamma_\sqcup (z), \\
C^{o\sharp}_{o\flat} & =\sum _{\grc_1\in \grC^{o\sharp}_\sqcup}\sum_{\grc_2\in \grC^{o\flat}_\sqcup}\sum_{z\in
  \pi_0(\mathcal{B}^\sigma (M_\sqcup, \grc_1, \grc_2))}
\langle\pmb{\varsigma }^*\mathbf{u}_{\lambda_\sqcup}, \mathbf{M}_{2,z}(W_\sqcup, \grc_1,
\grc_2)\rangle \Gamma_\sqcup (z), \\
D^{o\sharp}_{o\flat} & =\sum _{\grc_1\in \grC^{o\sharp}_\sqcup}\sum_{\grc_2\in \grC^{o\flat}_\sqcup}\sum_{z\in
  \pi_0(\mathcal{B}^\sigma (M_\sqcup, \grc_1, \grc_2))}
\langle\pmb{\varsigma }^*\mathbf{u}_{\lambda_{\sqcup+}}, \mathbf{M}_{1,z}(W_\sqcup, \grc_1,
\grc_2)\rangle \Gamma_\sqcup (z). \\
\end{split}
\end{equation}
The reducible variants of the above,
\(\bar{\smk}[\mathbf{u}](W_\sqcup)\), do not appear in the formulas
for \(A\), \(B\), \(C\), \(D\), as by assumption \(W_\sqcup\) is 
equipped with non-balanced perturbations. 
However, keep in mind that while \(B\), the simpler map among the four, is
assembled from the above according to the rule (\ref{def:m-map}) (substituting
\(m^\sharp_\flat[u]\) therein by \(B^{o\sharp}_{o\flat})\),  the more
general rule of Remark
\ref{rem:generalized-m} must be applied to construct  the more
complicated maps \(A\), \(C\), \(D\) from their constituents
above. Being a hat-version of a cobordism map, the formula for \(A\)
remains the same as that given in  (\ref{def:m-map}) since both end
points of \(\bar{\lambda }\) fall in \(Y_-\). The formulas for \(C\)
and \(D\) contain additional terms from the end points of \(\lambda \)
in \(Y_+\), similar to those in (\ref{def:m-lambda}) and
(\ref{def:m-lambda0}): for \(\smk=C\) or \(D\), the explicit formula for the
cobordism map \(C\) or \(D\) is 
 \begin{equation}\label{def:k-lambda}
\smk=\left[\begin{array}{cc}
\smk^{oo}_{oo} &\smk^{ou}_{oo}\\
\hat{\smk}^{oo}_{ou}
&\hat{\smk}^{ou}_{ou}
\end{array}\right],
\end{equation}
where (as the \(\bar{\smk}\)-terms vanish): 
\begin{equation}\label{def:k-lambda0}
\begin{split}
\hat{\smk}^{oo}_{ou}&
:=-(\bar{\partial}_\sqcup)^{os}_{ou}\, \smk^{oo}_{os}+(\bar{U}_\sqcup)^{os}_{ou}\, 
B^{oo}_{os},\\ 
\hat{\smk}^{ou}_{ou}&
:=-\bar{\smk}^{ou}_{ou}-(\bar{\partial}_\sqcup)^{os}_{ou}\, \smk^{ou}_{os}
+(\bar{U}_\sqcup)^{os}_{ou}\, B^{ou}_{os}.
\end{split}
\end{equation}

The maps \(A, C, D\) above also each has a companion version, denoted
respectively \(\dot{A}\), \(\dot{C}\), \(\dot{D}\). They are defined
from consitiuents given by the 
same formlas as in (\ref{def:A-cons}), with the cochains \(\uu_{\lambda _{\sqcup -}}\),
\(\uu_{\lambda _{\sqcup }}\), \(\uu_{\lambda _{\sqcup +}}\) therein replaced
by their more concrete variants: that is, 
respectively by \(\theta _{\lambda _{\sqcup -}}\), \(\theta
_{\lambda_{\sqcup -}}\wedge\theta _{\lambda _{\sqcup +}}\), \(\theta
_{\lambda _{\sqcup +}}\). The maps \(\dot{A}\), \(\dot{C}\),
\(\dot{D}\) are built form these constituents in manners parallel to
their sister versions above, but in the case of \(\dot{C}\),
\(\dot{D}\), all appearances of the boundary-obstructed map
\((\bar{U}_\sqcup)^{os}_{ou}\)'s  in (\ref{def:k-lambda0}) are 
replaced by its companion version,
 \[
1\otimes
\bar{n}^s_u[d\ul{\op{h}}_{\hat{p}_2}](M_2)-n^o_o[d\ul{\op{h}}_{\hat{p}_1}](M_1)\otimes
1=: \bar{n}^{os}_{ou}[d\ul{\op{h}}_{\hat{\sqcup}}] (M_\sqcup).
\] 
The aforementioned pairs of maps are related by the formulas:
\begin{equation}
\begin{split}\label{def:A-dot}
A& =\dot{A}-B\, \Theta _{\sqcup},\\
C& =\dot{C}+\Theta _{\sqcup}\, A-\dot{D}\, \Theta _\sqcup=\dot{C}+\Theta
_{\sqcup}\, \dot{A}-D\, \Theta _\sqcup,\\
D& =\dot{D}+\Theta _{\sqcup}\, B,
\end{split}
\end{equation}
where 
\[
\Theta _{\sqcup}=1\otimes\Theta _{p_2}-\Theta _{p_1}\otimes 1=\hat{n}[\ul{h}_{\hat{p}_1}]\otimes
1-1\otimes\hat{n}[\ul{h}_{\hat{p}_2}]=:-\hat{n}[\ul{h}_{\sqcup}].
\]
We shall make use of the following figure, which illustrates the construction of \(H_\#\) schematically. 
\begin{equation}\label{fig:Aa}
\text{\(-\infty\)-end at \(M_\#\)} \Opic{\Oblue\Ogrey\draw (3.1,2.1) node {\(\lambda_\#\)};
\draw (2,-0.8) node {$W_\#(S)\supset \lambda_\# $};}
\text{\(+\infty\)-end at \(M_\#\)} \quad \leadsto \quad \text{the map \(H_\#\)}.
\end{equation}

The construction of the maps \(A, B, C, D\co \hat{C}_\sqcup\to
\hat{C}_\sqcup\) (or their companions \(\dot{A}, \dot{B}:=B, \dot{C}, \dot{D}\)),
are illustrated in a similar fashion in the next set of pictures: 
\begin{equation}\label{fig:Ab}\begin{split}
\text{\(-\infty\)-end at \(M_\sqcup\)} \yjpic{\greyM\blueL \draw
  (2,-0.5) node {\(W_\sqcup (S)\supset \bar{\lambda }=\lambda_{\sqcup -}\)};
\draw (0.5,2.4) node {\(\bar\lambda\)};} \text{\(+\infty\)-end at
  \(M_\sqcup\)} & \leadsto \quad \text{the map \(A\) or \(\dot{A}\)}\\
\text{\(-\infty\)-end at \(M_\sqcup\)} \yjpic{\greyM\draw (2,-0.5)
  node {\(W_\sqcup (S)\supset \emptyset\)}; } \text{\(+\infty\)-end at
  \(M_\sqcup\)} & \leadsto \quad \text{the map \(B=\dot{B}\) }\\
\text{\(-\infty\)-end at \(M_\sqcup\)} \yjpic{\greyM\blueL\blueR \draw (2,-0.5) node {\(W_\sqcup (S)\supset \lambda \cup\bar{\lambda}=\lambda_\sqcup\)};
\draw (3.5,2.4) node {\(\lambda\)};
\draw (0.5,2.4) node {\(\bar\lambda\)};} \text{\(+\infty\)-end at
  \(M_\sqcup\)} & \leadsto \quad \text{the map \(C\) or \(\dot{C}\)}\\
 \text{\(-\infty\)-end at \(M_\sqcup\)} \yjpic{\blueR\greyM\draw
  (2,-0.5) node {\(W_\sqcup (S)\supset \lambda =\lambda_{\sqcup +}\)};\draw
  (3.5,2.4) node  {\(\lambda\)};} \text{\(+\infty\)-end at
  \(M_\sqcup\)} & \leadsto \quad \text{the map \(D\) or \(\dot{D}\)}\\
\end{split}
\end{equation}
The dotted 1-submanifold \(\gamma\) (possibly empty or disconnected)
in each of the cobordisms \(W\) in Figures
\ref{fig:Aa} and \ref{fig:Ab} is there to indicate that the map
on the right of the picture is constructed via coefficients given by
evaluating the cochain \(\mathbf{u}_\gamma \) (or \(\pmb{\theta
}_\gamma \)) associated to
\(\mathpzc{u}_\gamma \in \op{C}^*(\B^\sigma (W);\bbK)\) (or \(\theta
_\gamma \)) on relevant parametrized
moduli spaces \(\mathbf{M}\) associated to \(W\). (\(\gamma =\lambda
_\#\) in Figure \ref{fig:Aa}, \(\gamma =\lambda _{\sqcup -},
\emptyset, \lambda _\sqcup, \lambda _{\sqcup +}\) respectively in the
four lines in Figure \ref{fig:Ab}.) 

We now proceed with: 

\paragraph{\it (i) Verifying (\ref{eq:ch-htpy1}):} Re-expressed using
the more concrete companion, 
\(\theta _{\lambda _\#}\), of \(\uu_{\lambda _\#}\), this identity is 
equivalent to:  
\begin{equation}
\label{eq:ch-htpy1'}
\begin{split}
& \hat{m}[\theta _{\bar{\lambda }}](\bar{\V})\circ
\hat{m}[1](\V)+\hat{m}[1](\bar{\V})\circ \hat{m}[\theta _{\lambda
}](\V)\\
&\quad =m[\theta _{\lambda_\#}](W_\#)+\big[
  \hat{\smk}[\pmb{\theta }_{\lambda _\#}](W_\#), \partial_\#\big]\\
&\quad = m[\theta _{\lambda_\#}](W_\#)+[
  \dot{H}_\#, \partial_\#];
\end{split}
\end{equation}
To verify (\ref{eq:ch-htpy1'}), we shall apply 
(\ref{KM-composition}) to \(W=W_\#\), \(W_1=\V\),
\(W_2=\bar{\V}\), and \(\uu=\theta _{\lambda _\#}\). Note that
\(\delta\theta _{\lambda _\#}=0\), and thus the last term of
(\ref{KM-composition})  vanishes. The right hand side
of (\ref{eq:ch-htpy1'}) then coincides term by term with the first two
terms of the right hand side of (\ref{KM-composition}). To compute the
left hand side of (\ref{KM-composition}), namely
\(\hat{m}[\theta _{\lambda _\#}](W(\infty))\), we claim that 
\begin{equation}\label{theta-decomp0}
\theta _{\lambda _\#}=c(\theta _\lambda \otimes 1)+c(1\otimes \theta _{\bar{\lambda }}).
\end{equation}
This would then imply that \[
\hat{m}[\theta _{\lambda _\#}](W(\infty))=\hat{m}[\theta _{\bar{\lambda }}](\bar{\V})\circ
\hat{m}[1](\V)+\hat{m}[1](\bar{V})\circ \hat{m}[\theta _{\lambda }];
\]
thus establishing (\ref{eq:ch-htpy1'}). 

To verify (\ref{theta-decomp0}), recall the definitions \(\theta _{\lambda _\#}=d\op{hol}_{\lambda _\#}\),
\(\op{hol}_{\lambda _\#}\co \B^\sigma (W_\#)\to \bbR/\bbZ\) from Section \ref{sec:A-module}'s Part
1(b). Recall also the bundles and maps in (\ref{diag:B-lambda}) and
(\ref{CD-ttB}). Use \(\tilde{\tilde{\B}}_{p_1,p_2}^\sigma (M_\sqcup)\),
\(\tilde{\B}_{p_2-p_1}^\sigma (M_\sqcup)\), \(\B^\sigma (M_\sqcup)\),
\(\grC(M_\sqcup)\) to denote the terms, top to bottom, in the right
column of (\ref{diag:B-lambda}).  (\(\tilde{\B}_{p_2-p_1}^\sigma
(M_\sqcup)\) and \(\tilde{\B}_{p_1-p_2}^\sigma (M_\sqcup)\) denotes
the same space, but our convention is to use the notation 
\(\tilde{\B}_{p_2-p_1}^\sigma (M_\sqcup)\) when it is equipped with
the \(U(1)\)-action associated to the Thom form \(\vartheta
'_{p_2-p_1}\) given in (\ref{def:theta12}). Thus,
\(\tilde{\B}_{p_1-p_2}^\sigma (M_\sqcup)\) is endowed with the dual
\(U(1)\)-action.) 
 Recall the maps \(\op{hol}_{\lambda }\co \tilde{\B}^\sigma _\lambda
(\mathcal{V})\to \tilde{\tilde{\B}}_{p_1,p_2}(M_\sqcup)\), \(\op{hol}_{\bar{\lambda }}\co \tilde{\B}^\sigma _{\bar{\lambda}}
(\bar{\mathcal{V}})\to \tilde{\tilde{\B}}_{p_1,p_2}(M_\sqcup)\) from
  (\ref{CD-ttB}) and fix \(\vartheta '_{p_1}\), \(\vartheta
  '_{p_2}\), \(\rho _{\vartheta '_{p_1}}\), \(\rho _{\vartheta '_{p_2}}\) as was done there, 
  using them to define {\em both} \(\op{h}_\lambda \co \tilde{\B}^\sigma _\lambda
(\mathcal{V})\to \bbR/\bbZ\) and \(\op{h}_{\bar{\lambda }}\co \tilde{\B}^\sigma _{\bar{\lambda}}
(\bar{\mathcal{V}})\to \bbR/\bbZ\), as prescribed in Section
\ref{sec:A-module}'s Part 3(a). 

To interprete \(\uu_{\lambda _\#}\) in terms of inner products, consider now the analog of (\ref{CD:W-W_i}): 
Let \(\tilde{\Pi }^{M_\sqcup}_{W_\#}\co 
\tilde{\B}_{p_2-p_1}^\sigma (W_\#(S))\to\tilde{\B}^\sigma
_{p_2-p_1}(M_\sqcup)\) denote the pull-back of \(\Pi ^{M_\sqcup}_{W_\#}\co 
\B^\sigma (W_\#(S))\to\B^\sigma(M_\sqcup)\) under the the map \(\pi_{p_2-p_1}\co \tilde{\B}^\sigma
_{p_2-p_1}(M_\sqcup)\to \B^\sigma(M_\sqcup)\).  Let \(\tilde{\ss}_1\co
\tilde{\B}_{p_2-p_1}^\sigma (W_\#(S))\to \tilde{\B}_\lambda ^\sigma
(\V(S))\), \(\tilde{\ss}_2\co \tilde{\B}_{p_2-p_1}^\sigma (W_\#(S))\to
\tilde{\B}_{\bar{\lambda }}^\sigma  (\bar{\V}(S))\) be the direct 
analogs of their counterparts in 
(\ref{CD:W-W_i}). The analog of (\ref{CD:W-S}) in the present context
reads: 
 \begin{equation}\label{CD:W-S1}
\begin{CD}
 \tilde{\B}_{p_2-p_1}^\sigma (W_\#(S))@>\tilde{\ss}_1\times \tilde{\ss}_2>> \tilde{\B}^\sigma _\lambda
(\mathcal{V})\times  _{\tilde{\B}^\sigma _{p_2-p_1}(M_\sqcup)}\tilde{\B}^\sigma _{\bar{\lambda}}
(\bar{\V}) @>\text{embeds}>>\tilde{\B}^\sigma _\lambda
(\mathcal{V})\times   \tilde{\B}^\sigma _{\bar{\lambda}}
(\bar{\mathcal{V}}) \\
@V\pi^{W}_{\lambda _\#}VV  @V\pi^{\V}_{\lambda }\times
\pi^{\bar{\V}}_{\bar{\lambda }}VV @V
\pi^{\V}_{\lambda }\times \pi^{\bar{\V}}_{\bar{\lambda }}VV\\
\B^\sigma (W_\#(S)) @>\ss_1\times \ss_2>> \B^\sigma  (\V(S))\times _{\B^\sigma (M_\sqcup)}\B^\sigma
(\bar{\V}(S))@>\text{embeds}>> \B^\sigma  (\V(S))\times  \B^\sigma
(\bar{\V}(S).
\end{CD}
\end{equation}
Now observe that: 
\begin{itemize}
\item On the top row, the pull-back of the \(U(1)=\bbR/\bbZ\)-valued function
\begin{equation}\label{h-split}
\op{h}_\lambda \times 1 +1\times \op{h}_{\bar{\lambda }}
\end{equation}
on \(\tilde{\B}^\sigma _\lambda
(\mathcal{V})\times  \tilde{\B}^\sigma _{\bar{\lambda}}
(\bar{\mathcal{V}})\) to 
\(\tilde{\B}^\sigma _\lambda
(\mathcal{V})\times  _{\tilde{\B}_{p_2-p_1}(M_\sqcup)}\tilde{\B}^\sigma _{\bar{\lambda}}
(\bar{\mathcal{V}})\subset \tilde{\B}^\sigma _\lambda
(\mathcal{V})\times  \tilde{\B}^\sigma _{\bar{\lambda}}
(\bar{\mathcal{V}})\) 
does not depend on the choices of \(\vartheta '_{p_1}\), \(\vartheta
'_{p_2}\),  \(\rho _{\vartheta '_{p_1}}\), \(\rho _{\vartheta '_{p_2}}\). 
\item The preceding function is also invariant under the diagonal \(U(1)\) action on \(\tilde{\B}^\sigma _\lambda
(\mathcal{V})\times  _{\tilde{\B}_{p_2-p_1}(M_\sqcup)}\tilde{\B}^\sigma _{\bar{\lambda}}
(\bar{\mathcal{V}})\), and hence descends to define an \(\bbR/\bbZ\)-valued
function on the space in the middle of the bottom row of the diagram,
that is, \(\B^\sigma  (\V(S))\times _{\B^\sigma (M_\sqcup)}\B^\sigma
(\bar{\V}(S)).\)  
\item The \(\bbR/\bbZ\)-valued function \(\op{hol}_{\lambda _\#}\) on
  \(\B^\sigma  (W_\#)\) 
  agrees with the pull-back of the preceding function under the left
  arrow in the bottom row of the diagram, and we have 
\[
(\pi^{W_\#})^*\op{hol}_{\lambda
  _\#}=(\tilde{\ss}_1\times \tilde{\ss}_2)^*(\op{h}_\lambda \times 1
+1\times \op{h}_{\bar{\lambda }}).
\] 
Taking the differential on both sides, we have 
\[
(\pi^{W_\#})^*\theta _{\lambda
  _\#}=(\tilde{\ss}_1\times \tilde{\ss}_2)^*(\vartheta _\lambda \times 1
+1\times \vartheta _{\bar{\lambda }})
=(\tilde{\ss}_1)^*\vartheta_\lambda +(\tilde{\ss}_2)^*\vartheta
_{\bar{\lambda }}
\]
on \( \tilde{\B}_{p_2-p_1}^\sigma (W_\#(S))\). 
\item Recalling (\ref{vartheta-theta}), we then have 
\[
\theta _{\lambda _\#}=(\ss_1\times \ss_2)^*(\theta _\lambda   \times  1+1\times \theta
    _{\bar{\lambda }}), 
\] 
since \((\tilde{\ss}_1)^*(\tilde{\Pi}_{\lambda })^*\vartheta
  '_{p_2-p_1}=- (\tilde{\ss}_2)^*(\tilde{\Pi}_{\bar{\lambda }})^*\vartheta
  '_{p_1-p_2}\) on \(\tilde{\B}^\sigma _\lambda
(\mathcal{V})\times  _{\tilde{\B}^\sigma _{p_2-p_1}(M_\sqcup)}\tilde{\B}^\sigma _{\bar{\lambda}}
(\bar{\V})\). Meanwhile, the 1-form \(\theta _\lambda   \times  1+1\times \theta
    _{\bar{\lambda }}\) on \(\B^\sigma
(\V)\times \B^\sigma (\bar{\V})\) is nothing but \(c(\theta _\lambda
\otimes 1)+c(1\otimes \theta _{\bar{\lambda }}) \). 
\end{itemize}

\paragraph{\it (ii) Verifying (\ref{eq:ch-htpy}):} 
These identities follow directly from applying
(\ref{KM-composition0}) to \(W=W_\sqcup\), \(W_1=\bar{\V}\),
\(W_2=\V\), with \(\uu\) taken to be respectively to be
\(\uu_{\lambda _{\sqcup -}}\), \(\uu_\emptyset=1\), \(\uu_{\lambda
  _{\sqcup }}\), \(\uu_{\lambda _{\sqcup +}}\). These cochains have
natural interpretations as inner products: 
\[
\begin{split}
\uu_{\lambda _{\sqcup -}}& =c\, (\uu_{\bar{\lambda }}\otimes 1), \\
\uu_\emptyset&= c\, (1\otimes 1),\\
 \uu_{\lambda
  _{\sqcup }}& =c \, (\uu_{\bar{\lambda }}\otimes  \uu_\lambda ),\\
\uu_{\lambda_{\sqcup +}}& =c\, (1\otimes  \uu_{\lambda }).
\end{split}
\] 
In the case of item 2) of 
  (\ref{eq:ch-htpy}), \(\delta\uu_\emptyset =0\) and therefore the
  last term of (\ref{KM-composition0}) vanishes. Meanwhile, a
  straightforward adaptation of (\ref{eq:K-lambda}) to the
  parametrized setting identifies \(\hat{\smk}[\mathbf{u}]\) in the
  cases of \(\uu= \uu_{\lambda _{\sqcup -}}\), \(\uu_{\lambda_{\sqcup }}\), \(\uu_{\lambda _{\sqcup +}}\)
  respectively with the
  last terms of  (\ref{eq:ch-htpy})'s item 1), 2), 4). 


The figures below 
illustrate the identities  (\ref{eq:ch-htpy1}) in and
(\ref{eq:ch-htpy}), as well as hint on their origins. 

\(\bullet\) For  (\ref{eq:ch-htpy1}): 
\begin{equation}\label{fig:B1}\begin{split}
&\Opic{\OredM\OblueL}+
\Opic{\OredM\OblueR}=\\
&\qquad\qquad \Opic{\Oblue}-
\Opic{\OredL\Oblue\Ogrey}+
\Opic{\Oblue\Ogrey\OredR}
\end{split}
\end{equation}

\(\bullet\) The identity (\ref{eq:ch-htpy}) 1): 
\begin{equation}\label{fig:B5}\begin{split}
&\yjpic{\blueL\redM}=\yjpic{\blueL}-\\
&\qquad\qquad 
\yjpic{\blueLR\greyM\redL}+
\yjpic{\blueL\greyM\redR}+
\yjpic{\blueLL\greyM\redL}
  \end{split}
\end{equation}

\(\bullet\) The identity (\ref{eq:ch-htpy}) 2): 
\begin{equation}\label{fig:B3}
\yjpic{\redM}=
\yjpic{}+
\yjpic{\redL\greyM}-
\yjpic{\greyM\redR}
\end{equation}

\(\bullet\) The identity (\ref{eq:ch-htpy}) 3): 
\begin{equation}\label{fig:B4}\begin{split}
&\yjpic{\blueL\redM\blueR}=
\yjpic{\blueL\blueR}-
\yjpic{\redL\blueLR\greyM\blueR}\\
&\qquad\qquad+\yjpic{\blueL\greyM\blueRL\redR}+
\yjpic{\blueLL\redL\greyM\blueR}-
\yjpic{\blueL\greyM\redR\blueRR}
\end{split}
\end{equation}

\(\bullet\) The identity (\ref{eq:ch-htpy}) 4): 
\begin{equation}\label{fig:B2}\begin{split}
&\yjpic{\redM\blueR}=\yjpic{\blueR}\\
&\qquad\quad+
\yjpic{\redL\greyM\blueR}-
\yjpic{\greyM\redR\blueRL}-
\yjpic{\greyM\redR\blueRR}
\end{split}
\end{equation}
 

In each cobordism in the pictures, \(s\) increases from left to
right. They are read as follows: The dash lines (if present) in the cobordisms 
stand for 3-manifolds that split the cobordisms into a composition of
what we call ``factor cobordisms''. Each factor-cobordism (or the
cobordism itself, if it is not split) in the pictures is associated with a pair \((W,
\gamma )\), where \(W\) is a cobordism and \(\gamma \) 
is a 1-submanifold (possibly empty) of \(W\), the latter
being represented by 
dotted arcs or circles. This pair is associated with a cobordism map of the
form:
\begin{itemize} \item[(i)] 
\(\hat{m}[\uu_\gamma ]\) (resp. \(\hat{m}[\theta _\gamma
])\)) when \((W, \gamma )\) is not cylindrical;  
\item[(ii)] \(\hat{n}[\op{u}_\gamma ]\) (resp. \(\hat{m}[d \ul{\op{h}}_{\gamma
}])\)) when  \((W, \gamma
)\) is cylindrical, namely, it is of the form \(\bbR\times (Y, p)\),
\(p\) being a (possibly empty) 0-submanifold in \(Y\);
\item[(iii)] \(\hat{\smk}[\mathbf{u}_\gamma ]\)
(resp. \(\hat{\smk}[\pmb{\theta }_\gamma ]\)), when there is
a shaded region in the cobordism. 
\end{itemize}
Composition of cobordisms along the 
dashed lines correspond to compositions of maps associated to the
factor-cobordisms. For example, the dashed line in the first term of
(\ref{fig:B1}) splits the
composite cobordism \(W_\#\) into \(\mathcal{V}\supset \lambda\) on the left
and \(\bar{\cal V}\supset \emptyset\) on the right.  The left part
\(\mathcal{V}\supset \lambda\) corresponds to the map
\(V_1=\hat{m}[\mathpzc{u}_\lambda](\mathcal{V})\), and the right part
corresponds to the map
\(V_0^\dag=\hat{m}[\mathpzc{u}_\emptyset](\bar{\cal V})\); therefore
this term stands for \(V_0^\dag V_1\). The dashed line in the last
term splits \(W_\#(S)\) into \(W_\#(S)\supset \lambda_\#\) on the left
and the product cobordism \(\bbR\times M_\# \supset \emptyset\) on the
right. The former corresponds to the map \(H_\#\) according to Figure
\ref{fig:Aa}, and the latter corresponds to \(\partial_\#\). Thus this
term corresponds to the term \(\partial_\# \, H_\#\) in (\ref{eq:ch-htpy1}). 
With  \(\gamma \) again standing respectively for: \(\lambda
_\#\), \(\lambda _{\sqcup -}, \emptyset, \lambda _\sqcup, \lambda
_{\sqcup +}\) in Figures \ref{fig:B1}, \ref{fig:B5}, \ref{fig:B3},
\ref{fig:B4}, \ref{fig:B2},  the pictures suggest how each term of the identity arises from 
Stokes' theorem: that is, as the integral \(\mathbf{u}_\gamma \) over a constituent stratum of the
``boundary'' (to be more precise, cf. (\ref{bdry-bfM})) of the relevant compactifiied parametrized moduli space. 
Each such constituent stratum corresponds to the moduli space of  a 
particular type of  ``broken
\(W\)-paths''  (in
keeping with \cite{KM}'s terminology; cf.
Definition 23.3.2 therein).  The type for each term is specified by
the correponding picture, with dashed lines signifying  ``breaking
points'' of the broken \(W\)-path. 
The integrands in the identities, being defined from
differentials of holonomy maps along \(\gamma \), take the simple form
of an inner product under the decomposition when
the dotted arc/circle \(\gamma \) does not intersect the dashed
line. When they do intersect, the dashed
line splits \(\gamma \) into two arcs \(\gamma _1\) and \(\gamma
_2\), each lying in a factor-cobordism under the decomposition. 
The holonomy along \(\gamma \) being the product of the holonomy
along \(\gamma _1\) and that along \(\gamma _2\),
(cf. e.g. (\ref{h-split})), the integral of
\(\theta _\gamma \) over the spaces of such broken \(W\)-paths is thus
a sum of two terms, each involving integrating over one of the
\(\theta _{\gamma _i}\)'s. For example, this counts for the two terms on the left hand side of
(\ref{eq:ch-htpy1}), as well as the last two terms of (\ref{fig:B2}). 

\paragraph{\it Step 5.} The identities (\ref{eq:ch-htpy1}) and
(\ref{eq:ch-htpy}) reduce the proof of 
(\ref{eq:cob-comp1}) and (\ref{eq:cob-comp2}) to the next lemma, with
the maps \(H'\), \(A', B', C', D'\) from (\ref{eq:cob-comp1}) and
(\ref{eq:cob-comp2}) taken to be
\[\begin{split}
H' & =H-\textsc{h}\\
A' & =A-\textsc{a}\\
B' & =B-\textsc{b}\\
C' & =C-\textsc{c}\\
D' & =D-\textsc{d}, 
  \end{split}
\]
\(H\), \(A, B, C, D\) being the maps from (\ref{eq:ch-htpy1}) and
(\ref{eq:ch-htpy}), and \(\textsc{h}\), \(\textsc{a}, \textsc{b},
\textsc{c}, \textsc{d}\) being as in the lemma below. 
\begin{lemma}
There exists maps \(\textsc{h}\), \(\textsc{a}\), \(\textsc{b}\),
\(\textsc{c}\), \(\textsc{d}\) and \(\smz_\#, \smz_\sqcup\) such that 
\begin{equation}\label{eq:mW1}
 \op{Id}_\#-[\smz_\#, \partial_\#]_{\text{odd}}=\hat{m}[\mathpzc{u}_{\lambda_\#}](W_\#)+[\textsc{h}, \partial_\#]_{\text{even}};
\end{equation}
\begin{equation}\label{eq:mW}
\begin{split}
\op{1)}&  \quad \op{Id}_\sqcup -[\smz_\sqcup, \partial_\sqcup]_{\text{odd}}=\hat{m}[\mathpzc{u}_{\lambda_{\sqcup-}}](
  W_\sqcup)+[\textsc{a}, \hat{\partial}_\sqcup]_{\text{even}}+\textsc{b} \hat{U}_\sqcup,\\ 
 \op{2)} &  \quad 0=\hat{m}[1](
  W_\sqcup) -[\textsc{b}, \hat{\partial}_\sqcup ]_{\text{odd}},\\
\op{3)} &  \quad  [\hat{\partial}_\sqcup, \smx]_\ev -[\hat{U}_\sqcup, \smz_\sqcup ]_{\ev}=\hat{m}[\mathpzc{u}_{\lambda_\sqcup}](
  W_\sqcup) + [\textsc{c}, \hat{\partial}_\sqcup]_{\text{odd}}- \hat{U}_\sqcup \, \textsc{a} +\textsc{d} \hat{U}_\sqcup,, \\
\op{4)} & \quad \op{Id}_\sqcup-[\smz_\sqcup, \partial_\sqcup]_{\text{odd}}=\hat{m}[\mathpzc{u}_{\lambda_{\sqcup+}}](
  W_\sqcup) +[\hat{\partial}_\sqcup, \textsc{d}]_{\text{even}}-\hat{U}_\sqcup\,\textsc{b} .
\end{split}
\end{equation}
\end{lemma}
\pf These are also consequences of (\ref{KM-composition}), taking 
\(W=W_\#\), \(W_\sqcup\) respectively for (\ref{eq:mW1}) and
(\ref{eq:mW}), and with the same choices of \(\uu\) as in the previous
step. The splitting 3-manifolds \(Y_0\subset
W\)  however are chosen differently from those in the previous step. 

In the
case of \(W=W_\#\), we take \((W_1)_c\) to be a tubular neighborhood
of \(\lambda _\#\),  \(U(\lambda_\#)\), and so in this case 
\(Y_0=\partial (W_1)_c\simeq S^1\times S^2\); \((W_2)_c=W_\#\backslash
  U(\lambda_\#)\). There is a diffeomorphism taking the pair
  \(((W_1)_c, \lambda _\#)=(U(\lambda _\#), \lambda _\#)\) to \((S^1\times
  B^3, S^1\times \{0\})\), \(\{0\}\in B^3\) denoting the center of the 3-ball
  \(B^3\). We denote the embedded circle \( S^1\times \{0\}\subset S^1\times
  B^3\) by \(\gamma _0\). 
In the case when \(W=W_\sqcup\), we take \(Y_0\) to
  be the 3-sphere \(S_\sqcup\subset W_\sqcup\) described in Step
  3. This 3-sphere decomposes \(W_\sqcup\) as a connected sum of
  \(\bbR\times M_1\) and \(\bbR\times M_2\), and for both \(i=1,2\), \((W_i)_c\subset
  (W_\sqcup)_c\) is a manifold with boundary diffeomorphic to a
  product \([-1,1]\times M_i\) with an interior  4-ball
  removed. The rest of the proof is divided into several parts, 
  (i)--(viii) below. 

\paragraph{\it (i) Alternative metrics and perturbations.} 
A preliminary issue needs to be addressed before we are ready to apply
(\ref{KM-composition}). 
Recall that in the statement of the Lemma, the cobordism maps
\(\hat{m}[\uu](W_\#)\), \(\hat{m}[\uu](W_\sqcup)\) refer respectively
to \(\hat{m}[\uu](W_\#, \grs_{W_\#}, \varpi_{W_\#}; \Gamma _{W_\#}))\)
and \(\hat{m}[\uu](W_\sqcup, \grs_{W_\sqcup}, \varpi_{W_\sqcup};
\Gamma _{W_\sqcup})\), where the metrics and the closed 2-forms
\(\omega_{W_\#}\), \(\omega_{W_\sqcup}\) are defined via the
decompositions of \(W_\#\), \(W_\sqcup\) along \(M_\#\) and
\(M_\sqcup\). To apply the composition formula (\ref{KM-composition0})
or (\ref{KM-composition}) to the alternate decomposition described in
the preceding paragraph, we need to work with cobordism maps
associated to different choices of metrics and perturbation forms,
which are compatible with the aforementioned alternate decomposition
of \(W_\#\) and \(W_\sqcup\). However, we claim that  the identities (\ref{eq:mW1}),
(\ref{eq:mW}) are equivalent to identities of the same form for \(\hat{m}[\uu](W, \grs_W,
2\omega^+; \Gamma _W)\), \(W=W_\#\) or \(W_\sqcup\), with the latter
endowed with different metrics and perturbation forms \(\omega\), 
as along as 
\begin{equation}\label{def:var-metric}
\parbox{25em}{the differences are supported on compact regions in \(W\),
and in the case of the perturbation form \(\omega\), the difference is
exact. }
\end{equation}
(The maps \(\textsc{h}\), \(\textsc{a}\), \(\textsc{b}\)
\(\textsc{c}\), \(\textsc{d}\) will be altered, but that is
inconsequential.) This claim follows from \cite{KM}'s Proposition 25.3.8 (extended in the manner previously 
described, and with changes in perturbation forms incoporated). 


Slightly reformulated, the hat-version of the identity in \cite{KM}'s
Proposition 25.3.8 takes the following form: Suppose there is a path
of pairs consisting of a metric and a perturbation form on \(W\), so that
(\ref{def:var-metric}) holds for  the entire path.  Denote by
\(\hat{m}^+[\uu](W)\), \(\hat{m}[\uu](W)\) respectively for the
version of \(\hat{m}^-[\uu](W)\) associated to the pair at the end and
at the beginning of the path. Then: 
\begin{equation}\label{eq:var-metric}
\hat{m}^+[\uu](W)-\hat{m}^-[\uu](W)=[\hat{\smz}[\mathbf{u}], \hat{\partial}]+
\hat{\smz}[\delta \mathbf{u}], 
\end{equation}
where the \(\hat{\smz}\)-maps are defined using parametrized moduli spaces
associated to this path of metrics and perturbations parallel to the
definition of the previously introduced \(\hat{\smk}\)-maps, and
\(\mathbf{u}\) is the parametrized variant of \(\uu\) as before. (The
\(Z\)-maps are analogs of the \(\check{K}\)-maps in \cite{KM}'s Proposition 25.3.8.) Our
signs differ from those in \cite{KM}'s Proposition 25.3.8 because we
adopt the ``fiber-last'' convention of orienting the parametrized
moduli spaces, as opposed to \cite{KM}'s ``fiber-first''
convention. This is preferred as it is more consistent with the
orientation convention used for (\ref{KM-composition}). This identity in hand,
suppose identities of the form (\ref{eq:mW1}) and
(\ref{eq:mW}) are established for a particular pair of metric and
perturbation form. Use \(\hat{m}^-[\uu](W)\) for the version of cobordism maps
associated to this pair, and use \(\textsc{h}^-\), \(\textsc{a}^-\), \(\textsc{b}^-\)
\(\textsc{c}^-\), \(\textsc{d}^- \) to denote the version of maps \(\textsc{h}\), \(\textsc{a}\), \(\textsc{b}\)
\(\textsc{c}\), \(\textsc{d}\) in this version of (\ref{eq:mW1}) and
(\ref{eq:mW}). On the other hand, use \(\hat{m}^+[\uu](W)\) to denote the
the version of cobordism maps associated to the pair of metric and
perturbation appearing in the statement of the lemma. Then combining
the (\(-\))-versions of the identities  (\ref{eq:mW1}), 
(\ref{eq:mW}) with (\ref{eq:var-metric}), one would have a
(\(+\))-version of the identities  (\ref{eq:mW1}), 
(\ref{eq:mW}) with respect to a new set of maps \(\textsc{h}\), \(\textsc{a}\), \(\textsc{b}\)
\(\textsc{c}\), \(\textsc{d}\), if the latter is set to be: 
\[\begin{split}
\textsc{h}^+& = \textsc{h}^--\hat{\smz}[\pmb{1}](W_\#);\\
\textsc{a}^+ & =\textsc{a}^--\hat{\smz}[\mathbf{u}_{\lambda _{\sqcup -}}](W_\sqcup);\\
\textsc{b}^+& =\textsc{b}^-+\hat{\smz}[\pmb{1}](W_\sqcup);\\
\textsc{c}^+ &= \textsc{c}^--\hat{\smz}[\mathbf{u}_{\lambda_\sqcup}](W_\sqcup)\\
\textsc{d}^+& =\textsc{d}^-+\hat{\smz}[\mathbf{u}_{\lambda_{\sqcup
    +}}](W_\sqcup). 
\end{split}
\]
To reach the preceding conclusion, we made use of the identities: 
\[
\begin{split}
\hat{\smz}[\delta\mathbf{u}_{\lambda _{\sqcup
    -}}](W_\sqcup) & =\hat{U}_\sqcup\, \hat{\smz}[\pmb{1}](W_\sqcup);\\
\hat{\smz}[\delta\mathbf{u}_{\lambda _{\sqcup
    +}}](W_\sqcup) & =-\hat{\smz}[\pmb{1}](W_\sqcup) \, \hat{U}_\sqcup;\\
\hat{\smz}[\delta\mathbf{u}_{\lambda _{\sqcup}}](W_\sqcup) & =-\hat{\smz}[\mathbf{u}_{\lambda _{\sqcup
    +}}](W_\sqcup) \, \hat{U}_\sqcup  -\hat{U}_\sqcup\,  \hat{\smz}[\mathbf{u}_{\lambda _{\sqcup
    -}}](W_\sqcup). 
\end{split}
\]
In complete parallel to the \(\hat{\smk}\)-analogs mentioned in the
paragraph preceding (\ref{fig:B1}), these identities are also
straightforward adaptations of 
(\ref{eq:K-lambda}).  

Now permitted to work with alternative metrics and perturbation forms
by the preceding arguments, we endow \(W_\#\) and \(W_\sqcup \) with
the following sort of metrics and perturbations for the rest of this proof. For \(W=W_\sqcup\) or \(W_\#\), we require the
metric to:  
\begin{itemize}
\item agree with a product metric on
 a tubular neighborhood \(U(Y_0)\simeq[-1, 1]\times Y_0\) of
 \(Y_0\subset W\), where \(Y_0\), is \(S^1\times S^2\) in the case
 of \(W_\#\) and \(S^3\) in the case of \(W_\sqcup\). These are endowed with
 the standard metrics. (In particular, these metrics on \(Y_0\) have
 positive scalar curvature.)
\item agree with the original metric on the complements of
\((W_\#)_c\) and \((W_\sqcup)_c\). 
\item in the case of \(W_\#\), the restriction of the metric to
  \((W_1)_c\simeq S^1\times B^3\) has nonnegative scalar curvature and is
  invariant under rotation along the \(S^1\)-factor. 
\end{itemize}
The abbreviated notations
 \(\hat{m}[\uu](W_\#)\), \(\hat{m}[\uu](W_\sqcup)\) also take on
 different 
 meanings for the rest of this proof: they will stand respectively for \(\hat{m}[\uu](W_\#, \grs_{W_\#},
2\dot{\omega}^+_{\#}; \Gamma _{W_\#})\) and \(\hat{m}[\uu](W_\sqcup, \grs_{W_\sqcup},
2\dot{\omega}^+_{\sqcup}; \Gamma _{W_\sqcup})\), where the metrics are
as previously mentioned, and 
\(\dot{\omega}_{\#}\), \(\dot{\omega}_{\#}\) are closed 2-forms that:
\begin{itemize}
\item vanish over \(U(Y_0)\);
\item  are cohomologous respectively to 
\(\omega_{\#}\) and \(\omega_\sqcup\);
\item agree respectively with 
\(\omega_{\#}\) and \(\omega_\sqcup\) on the complements of
\((W_\#)_c\) and \((W_\sqcup)_c\);
\item \(\dot{\omega}_{\#}\) vanishes over \((W_1)_c\simeq
  S^1\times B^3\subset W_\#\). 
\end{itemize}
Such \(\dot{\omega}_{\#}\), \(\dot{\omega}_{\#}\) exist as  \(\omega_{W_\#}\) and \(\omega_{W_\sqcup}\) both restrict to
exact forms on \(Y_0\), and \(\omega_{\#}\) restricts to an exact form
on \((W_1)_c\simeq   S^1\times B^3\subset W_\#\). 

Now write \(W=W_2\circ W_1\) and define \(W(S)\) according to the
recipe (\ref{def:W(S)}). In the case \(W=W_\#\), \(W_1\) is regarded
as a cobordism from the empty set to \(S^1\times S^2\); the
\(-\infty\)-end of \(W_2\) consists of two connected components,
\(S^1\times S^2\) and \(M_\#\), but only the \(S^1\times S^2\)
component is ``glued to'' \(W_1\) to form \(W(S)\). See Figure
\ref{Figure (a)} (a) for an illustration. 
In the case \(W=W_\sqcup\), the
\(+\infty\)-end of \(W_1\) consists of two connected components,
\(M_1\) and \(S^3\), and the \(-\infty\)-end of \(W_2\) consists of
two connected components as well, \(S^3\) and \(M_2\), but only the
\(S^3\)-ends from both sides are ``glued'' to form \(W(S)\). See
Figure \ref{Figure (a)} (b) below. To indicate the 3-manifold where 
gluing take place in the composition, we write
\(W_\#=W_2\circ _{S^1\times S^2} W_1\) and \(W_\sqcup= W_2 \circ_{S^3} W_1\). 

\paragraph{\it  (ii) Surgered cobordisms.} Recall also from Step 3 above that surgery along \(\lambda
_\#\subset W_\#\) and \(S_\sqcup\subset W_\sqcup\) give respectively
the 
cobordisms \(W'_\#\simeq\bbR\times M_\#\) and \(W'_\sqcup\simeq
\bbR\times M_\sqcup\). We decompose these surgered cobordisms in a way
compatible with the decomposition of \(W_\#\), \(W_\sqcup\) above. 
In the case when \(W=W_\#\), the corresponding surgered manifold is decomposed as \(W_\#'=W_2\circ _{S^1\times S^2}W_1'\), with
\((W_1')_c\simeq D^2\times S^2\) and \(W_2\) being as in the
decomposition of \(W_\#\), \(W_\#=W_2\circ _{S^1\times S^2}W_1\). Like  \((W_1)_c\subset (W_\#)_c\),
we also equip 
\((W'_1)_c\) with a metric with non-negative scalar curvature. 
See Figure \ref{Figure (a)} (a\('\)) for an illustration. 
In the case when \(W=W_\sqcup\), the surgered manifold has two
connected components, 
\(W_\sqcup '=\hat{W}_1\sqcup \hat{W}_2\). Each connected component \(\hat{W}_i\simeq
\bbR\times M_i\), \(i=1,2\),  is obtained from \(W_i\subset W_\sqcup\) by filling in a 4-ball at
the boundary 3-sphere of \(W_i\).  Let \(\grB_c\) denote a closed
4-ball equipped with a metric which  has non-negative scalar curvature,
and that it is cylindrical on a collar of the boundary. Let \(\grB\) be
the corresponding manifold with one cylindrical end, regarded as a
cobordism from the empty set to \(S^3\), and let \(\bar{\grB}\) denote
the reverse cobordism. 
We decompose \(\hat{W}_1\), \(\hat{W}_2\) respectively as
\(\hat{W}_1=\bar{\grB }\circ _{S^3}W_1\), and \(\hat{W}_2=W_2\circ
_{S^3}\grB\),
where \(W_1\), \( W_2\) are as in the decomposition of \(W_\sqcup\). 
See Figure \ref{Figure (a)} (b\('\)) for an illustration.

We now apply the decomposition theorem (\ref{KM-composition0}),
(\ref{KM-composition}) to the composite cobordisms \(W=W_\#\),
\(W'_\#\), \(W_\sqcup\), \(W'_\sqcup\) described above,
and  illustrated
schematically in Figures \ref{Figure (a)} (a), (a\('\)), (b),
(b\('\)). 
In each  of the pictures, the 
dashed line represents \(Y_0\), the 3-manifold where composition
takes place. The shaded regions in each picture, \(W_2\subset
W_\#\), \(W_2\subset W'_\#\), and \(W_1\subset W_\sqcup\), \(W_1\subset
W'_\sqcup\), are associated with nonbalanced perturbation forms in the relevant
Seiberg-Witten equation, implying that the corresponding moduli spaces
of Seiberg-Wittten solutions contain no reducible elements.

\refstepcounter{equation}\phantomsection\label{Figure (a)}
\noindent\begin{tikzpicture}[scale=0.7]
\draw[fill=gray!40!white] (0,0) to [out=0,in=185] (3,0.1)
to [out=5,in=180] (8,0.5) -- (8,2) 
to [out=180,in=-10] (3,3.5) -- (3,2.3) 
to [out=-10,in=90] (4,1.7) 
to [out=-90,in=5] (3,1.1)
to [out=185,in=0] (0,1) -- (0,0);
\draw [thick] (0,0) to [out=0,in=185] (3,0.1)
to [out=5,in=180] (8,0.5);
\draw [thick] (8,2) to [out=180, in=-10] (3,3.5);
\draw [black, dotted, ultra thick] (3,3.5) -- (3,2.3);
\draw [thick] (3,2.3)
to [out=-10,in=90] (4,1.7) 
to [out=-90,in=5] (3,1.1)
to [out=185,in=0] (0,1);
\draw [thick] (3,3.5) to [out=170,in=90] (0.5,3) 
to [out=270,in=170] (3,2.3);
\draw [thick] (1.2,3.1) to [out=-30,in=160] (1.5,2.9)
to [out=-20,in=200] (2,2.9)
to [out=20,in=210] (2.3,3.1);
\draw [thick] (1.5,2.9) to [out=20,in=160] (2,2.9);
\draw [black, thick, dotted] (1.75,3) ellipse (0.8 and 0.4);
\node at (-0.5,0.5) {\(M_\#\)};
\node at (8.5,1.25) {\(M_\#\)};
\node at (5,1.5) {\(W_2\)};
\node at (1.7,4) {\(\overbrace{\qquad\qquad}^{W_1}\)};
\node [black, left] at (3.5,2.0) {\(\scriptstyle Y_0=S^1\times S^2\)};
\node [black, right] at (0.4,2.7) {\(\scriptstyle \lambda_\#\)};
\node [right] at (0,-2) {\parbox{0.4\hsize}{Figure \theequation\,
    (a):\\\hspace*{0.5cm}{\small \(W_\#=W_2\circ W_1\);\\
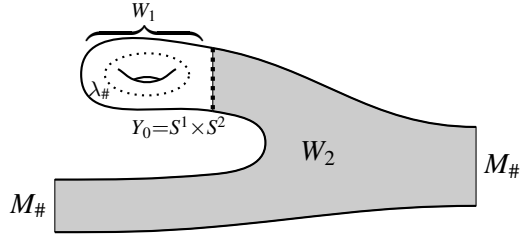
 \hspace*{0.5cm}\((W_1,\lambda_\#)\simeq (S^1\times B^3,\gamma_0)\);\\\hspace*{0.5cm}\(s\) increases
  from left to right.}}}; 
\end{tikzpicture}
\qquad\begin{tikzpicture}[scale=0.7]
\draw [thick] (0,0) to [out=0,in=185] (3,0.1)
to [out=5,in=180] (6,0.4);
\draw [thick] (6,2.4) to [out=160, in=-10] (3,3.5);
\draw [thick, dotted] (6,2.3) to [out=160, in=-10] (3,3.4);
\draw [thick] (3,2.3)
to [out=-10,in=90] (4,1.7) 
to [out=-90,in=5] (3,1.1)
to [out=185,in=0] (0,1);
\draw [thick, dotted] (3,2.4)
to [out=-10,in=90] (4.1,1.7) 
to [out=-90,in=5] (3,1.0)
to [out=185,in=0] (0,0.9);
\draw[fill=gray!40!white] (6,5.8) to [out=180,in=5] (3,5.7)
to [out=185,in=0] (0,5.4) -- (0,3.4) 
to [out=-20,in=170] (3,2.3) -- (3,3.5)
to [out=170,in=270] (2,4.1) 
to [out=90,in=185] (3,4.7)
to [out=5,in=180] (6,4.8) -- (6,5.8);
\draw [thick] (6,5.8) to [out=180,in=5] (3,5.7)
to [out=185,in=0] (0,5.4);
\draw [thick] (0,3.4) to [out=-20, in=170] (3,2.3);
\draw [thick, dotted] (0,3.5) to [out=-20, in=170] (3,2.4);
\draw [black, dotted, ultra thick] (3,2.3) -- (3,3.5);
\draw [thick] (3,3.5)
to [out=170,in=-90] (2,4.1) 
to [out=90,in=185] (3,4.7)
to [out=5,in=180] (6,4.8);
\draw [thick,dotted] (3,3.4)
to [out=170,in=-90] (1.9,4.1) 
to [out=90,in=185] (3,4.8)
to [out=5,in=180] (6,4.9);
\node at (-0.5,0.5) {\(M_2\)};
\node at (6.5,1.25) {\(M_2\)};
\node at (5,1.5) {\(W_2\)};
\node at (-0.5,4.5) {\(M_1\)};
\node at (6.5,5.25) {\(M_1\)};
\node at (1,4) {\(W_1\)};
\draw [black, fill] (3,2.4) circle [radius=0.05];
\draw [black, fill] (3,3.4) circle [radius=0.05];
\draw [black, fill] (0,0.9) circle [radius=0.05];
\draw [black, fill] (6,2.3) circle [radius=0.05];
\draw [black, fill] (0,3.5) circle [radius=0.05];
\draw [black, fill] (6,4.9) circle [radius=0.05];
\node [black, above] at (2.8,3.05) {\(\scriptscriptstyle x\)};
\node [black, below] at (2.8,2.8) {\(\scriptscriptstyle\overline{x}\)};
\node at (-0.2,0.95) {\(\scriptscriptstyle p_2\)};
\node at (6.3,2.3) {\(\scriptscriptstyle p_2\)};
\node at (-0.2,3.5) {\(\scriptscriptstyle p_1\)};
\node at (6.3,4.8) {\(\scriptscriptstyle p_1\)};
\node [black, right] at (3,2.7) {\(\scriptstyle Y_0=S^3\)};
\node at (4,4.5) {\(\scriptstyle \lambda_1\)};
\node at (4.8,3.15) {\(\scriptstyle \lambda_2\)};
\node at (1.5,2.5) {\(\scriptstyle \overline{\lambda}_1\)};
\node at (2,1.3) {\(\scriptstyle \overline{\lambda}_2\)};
\node [right] at (0,-2)
%
{\parbox{0.4\hsize}{Figure \theequation\,
    (b):\\\hspace*{0.5cm}{\small \(W_\sqcup=W_2\circ W_1\);\\
 \hspace*{0.5cm}\(\lambda=\lambda_1\cup_x\lambda_2\);\\\hspace*{0.5cm}\(\overline{\lambda}=\overline{\lambda}_1\cup_{\bar{x}}\overline{\lambda}_2\)}}};
\end{tikzpicture}

\noindent\begin{tikzpicture}[scale=0.7]
\draw[fill=gray!40!white] (0,0) to [out=0,in=185] (3,0.1)
to [out=5,in=180] (8,0.5) -- (8,2) 
to [out=180,in=-10] (3,3.5) -- (3,2.3) 
to [out=-10,in=90] (4,1.7) 
to [out=-90,in=5] (3,1.1)
to [out=185,in=0] (0,1) -- (0,0);
\draw [thick] (0,0) to [out=0,in=185] (3,0.1)
to [out=5,in=180] (8,0.5);
\draw [thick] (8,2) to [out=180, in=-10] (3,3.5);
\draw [black, dotted, ultra thick] (3,3.5) -- (3,2.3);
\draw [thick] (3,2.3)
to [out=-10,in=90] (4,1.7) 
to [out=-90,in=5] (3,1.1)
to [out=185,in=0] (0,1);
\draw [thick] (3,3.5) to [out=170,in=90] (1.5,3) 
to [out=270,in=170] (3,2.3);
\node at (-0.5,0.5) {\(M_\#\)};
\node at (8.5,1.25) {\(M_\#\)};
\node at (5,1.5) {\(W_2\)};
\node at (2.2,4) {\(\overbrace{\qquad\;\;\;}^{W_1'}\)};
\node [black, left] at (3.5,2.0) {\(\scriptstyle Y_0=S_1\times S_2\)};
\node [right] at (0,-2.5) {\parbox{0.4\hsize}{Figure
    \theequation\, (a\('\)):\\\hspace*{0.5cm}{\small \(W_\#'=W_2\circ
    W_1'\simeq \bbR\times M_\#\);\\
 \hspace*{0.5cm}\(W_1'\simeq D^2\times S^2\) }\\ \\     }};
\end{tikzpicture}\qquad
\begin{tikzpicture}[scale=0.7]
\draw [thick] (0,0) to [out=0,in=185] (3,0.1)
to [out=5,in=180] (6,0.4);
\draw [thick] (6,2.4) to [out=160, in=-10] (3,3.5);
\draw [thick, dotted] (6,2.3) to [out=160, in=-10] (3,3.4);
\draw [black, dotted, ultra thick] (3,3.5) -- (3,2.3);
\draw [thick] (3,2.3)
to [out=-10,in=90] (4,1.7) 
to [out=-90,in=5] (3,1.1)
to [out=185,in=0] (0,1);
\draw [thick, dotted] (3,2.4)
to [out=-10,in=90] (4.1,1.7) 
to [out=-90,in=5] (3,1)
to [out=185,in=0] (0,0.9);
\draw [thick] (3,3.5) to [out=170,in=90] (1.5,3) 
to [out=270,in=170] (3,2.3);
\draw[fill=gray!40!white] (6,7.5) to [out=180,in=5] (3,7.4)
to [out=185,in=0] (0,7.1) -- (0,5.1) 
to [out=-20,in=170] (3,4) -- (3,5.2)
to [out=170,in=270] (2,5.8) 
to [out=90,in=185] (3,6.4)
to [out=5,in=180] (6,6.5) -- (6,7.5);
\draw [thick] (6,7.5) to [out=180,in=5] (3,7.4)
to [out=185,in=0] (0,7.1);
\draw [thick] (0,5.1) to [out=-20, in=170] (3,4);
\draw [thick, dotted] (0,5.2) to [out=-20, in=170] (3,4.1);
\draw [black, dotted, ultra thick] (3,4) -- (3,5.2);
\draw [thick] (3,5.2)
to [out=170,in=-90] (2,5.8) 
to [out=90,in=185] (3,6.4)
to [out=5,in=180] (6,6.5);
\draw [thick, dotted] (3,5.1)
to [out=170,in=-90] (1.9,5.8) 
to [out=90,in=185] (3,6.5)
to [out=5,in=180] (6,6.6);

\draw [thick] (3,4) to [out=-10,in=-90] (5.0,4.5) 
to [out=90,in=-10] (3,5.2);
\draw [thick, dotted] (3,4.1) to [out=-10,in=-90] (4.9,4.5) 
to [out=90,in=-10] (3,5.1);
\node at (-0.5,0.5) {\(M_2\)};
\node at (6.5,1.25) {\(M_2\)};
\node at (5,1.5) {\(W_2\)};
\node at (-0.5,6) {\(M_1\)};
\node at (6.5,7) {\(M_1\)};
\node at (1,5.8) {\(W_1\)};
\node at (2.3,2.95) {\(\grB\)};
\node at (3.8,4.55) {\(\bar{\grB}\)};
\node at (1.3,2.95) {\(\scriptstyle \gamma_2\)};
\node at (4.7,4.55) {\(\scriptstyle \gamma_1\)};
\node at (2.8,5.0) {\(\scriptscriptstyle x\)};
\node at (2.8,4.2) {\(\scriptscriptstyle \overline x\)};
\node at (2.8,3.3) {\(\scriptscriptstyle x\)};
\node at (2.8,2.5) {\(\scriptscriptstyle \overline x\)};
\node at (-0.2,0.95) {\(\scriptscriptstyle p_2\)};
\node at (6.3,2.3) {\(\scriptscriptstyle p_2\)};
\node at (-0.2,5.2) {\(\scriptscriptstyle p_1\)};
\node at (6.3,6.6) {\(\scriptscriptstyle p_1\)};
\draw [black, fill] (3,2.4) circle [radius=0.05];
\draw [black, fill] (3,3.4) circle [radius=0.05];
\draw [black, fill] (3,4.1) circle [radius=0.05];
\draw [black, fill] (3,5.1) circle [radius=0.05];
\draw [black, fill] (0,0.9) circle [radius=0.05];
\draw [black, fill] (6,2.3) circle [radius=0.05];
\draw [black, fill] (0,5.2) circle [radius=0.05];
\draw [black, fill] (6,6.6) circle [radius=0.05];
\node at (4,6.2) {\(\scriptstyle \lambda_1\)};
\node at (4.8,3.15) {\(\scriptstyle \lambda_2\)};
\node at (1.5,4.2) {\(\scriptstyle \overline{\lambda}_1\)};
\node at (2,1.3) {\(\scriptstyle \overline{\lambda}_2\)};
\node at (3,8)
{\(\overbrace{\qquad\qquad\qquad\qquad\qquad\qquad}^{\hat W_1}\)};
\node at (3,-0.7)
{\(\underbrace{\qquad\qquad\qquad\qquad\qquad\qquad}_{\hat W_2}\)};
\node [right] at (0,-3.5) 
{\parbox{0.4\hsize}{Figure
    \theequation\, (b\('\)):\\\hspace*{0.5cm} {\small \(W_\sqcup '=\hat
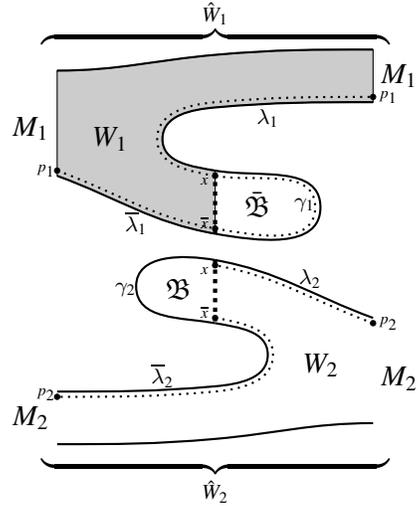
    W_1\sqcup \hat W_2\);\\
 \hspace*{0.5cm}\(\hat W_i=W_i \cup_{S^3} B^4\)\\
 \hspace*{0.5cm}\(\quad\;\; \simeq \bbR\times M_i\), \(i=1,2\).\\
 \hspace*{0.3cm} \( (-1)^i \hat{p}_i=\lambda _i\cup\gamma _i\cup\bar{\lambda }_i\)\\
\hspace*{0.3cm}\(\; \; \simeq\bbR\times \{ (-1)^i p_i\}\subset
 \bbR\times M_i\). }}};
\end{tikzpicture}

of Seiberg-Wittten solutions contain no reducible elements.

\paragraph{\it (iii) Some useful facts about \(\mathring{C}(Y_0)\). }
In all four pictures, the \(\Spin^c\) structure \(\grs_W\) on \(W\)
restricts to the trivial \(\Spin^c\) structure on \(Y_0\), denoted
\(\grs_0\) below. (\(\grs_0\) is characterized by the condition
\(c_1(\grs_0)=0\).) The 3-manifolds \(Y_0\) also all carry positive
scalar curvatures, so that \(\grC^o(Y_0)\) is trivial and the Floer
complex 
\(\hat{C}(Y_0, \grs_0)=C^u(Y_0, \grs_0)\). In fact, \(Y_0=S^1\times S^2, S^3\)
respectively in the cases of \(W_\#\), \(W_\sqcup\), and in both cases  \(\mathring{C}(Y_0,
\grs_0)\) are explicitly described in \cite{KM} (cf. e.g. Chapter 36 therein). 
We write 
\begin{equation}\label{C(Z)}
(\hat{C}(S^1\times S^2, \grs_0), \hat{\partial})=(\bbK[u,y]\hat{1}_Z,
0),
\end{equation}
where \(\bbK[u,y]\) is the (graded) polynomial algebra with
variables \(u\) and \(y\), 
\(\deg(u)=-2\), \(\deg (y)=1\), and \(\hat{1}_Z\) denotes an element
with degree \(-1\). Similarly,  
\[
(\hat{C}(S^3, \grs_0), \hat{\partial})=(\bbK[u]\hat{1},0), \quad
(\check{C}(S^3, \grs_0), \hat{\partial})=(\bbK[u^{-1}]\check{1},0), 
\] 
\(\hat{1}\in \hat{C}(S^3, \grs_0)\), \(\check{1}\in \check{C}(S^3,
\grs_0)\) being respectively the elements
of degree \(-1\) and \(0\) explicity described in \cite{KM}. (In our
convention, the plane field on \(S^3\) denoted \([\xi_-]\) and the
plane field on \(S^1\times S^2\) denoted \([\xi_0]\) in
\cite{KM} both have degree 0.) We also use the notation
\(u^n\hat{1}\), \(u^{-n}\check{1}\),  \(n\in \bbZ^{\geq 0}\) 
to denote respectively the element in \(\grC^u(S^3)\) with
\(\op{gr}\)-grading 
\([\xi_0]-1-2n\) (equivalently, \(\bar{\op{gr}}\)-grading
\([\xi_0]-2-2n\)), and the element in \(\grC^u(S^3)\) with
 \( \op{gr}\)-grading 
\([\xi_0]+2n\) (equivalently, \(\bar{\op{gr}}\)-grading
\([\xi_0]+2n\)). The \(\bar{U}\)-action on \(\bar{C}(S^3)\) is also
well-known: \(\langle\grc', \bar{U}(S^3)\grc\rangle= 1\) for any pair of
\(\grc, \grc'\in \grC(S^3)\) with
\(\bar{\op{gr}}(\grc)-\bar{\op{gr}}(\grc')=2\).

We are now ready to proceed with: 

\paragraph{\it \((iv)\) Verifying (\ref{eq:mW1}):} To compute
\(\hat{m}[\uu_{\lambda _\#}](W_\#)\), first note that the 1-cochain \(\uu_{\lambda _\#}\) has the simple
form of an inner product, \(\uu_{\lambda _\#}=c(\uu_{\gamma
  _0}\otimes  1)\) with respect to the decomposition of  \(W_\#\) shown in
Figure \ref{Figure (a)} (a). Thus, (\ref{KM-composition0}) is directly
applicable. Noting that \(\delta\uu_{\lambda_\#}=0\), this gives us: 
\begin{equation}\label{comp:W-Z}
 \hat{m}[1](W_2) \,  \hat{m}[\uu_{\gamma _0}](S^1\times
 B^3)=\hat{m}[\uu_{\lambda_\#}](W_\#)+\big[\hat{\textsc{k}}[\mathbf{u}_{\lambda
   _\#}](W_\#), \hat{\partial}_\#\big].
\end{equation}
This is compared to the formula obtained by applying 
(\ref{KM-composition0}) to \(W'_\#\),  decomposed as shown in
Figure \ref{Figure (a)} (a\('\)), and with the 1-cochain taken to be
\(\uu =1=c(1 \otimes 1)\). Here we have 
\begin{equation}\label{comp:W-Z1}
 \hat{m}[1](W_2) \,  \hat{m}[1](D^2\times
 S^2)=\hat{m}[1](W'_\#)+[\hat{\textsc{k}}[\pmb{1}](W_\#'), \partial_\# ]. 
\end{equation}
If the manifolds \(W_1\simeq S^1\times B^3\) and \(W_1'\simeq D^2\times
 S^2\) above are regarded cobordisms from the empty set to
 \((S^1\times S^2, \grs_0)\), then \(\hat{m}[\uu](W_1)\),
 \(\hat{m}[\uu](W'_1)\) are both elements of \(\hat{C}(S^1\times S^2,
 \grs_0)\). Alternatively, (in line with the definition of closed
 4-manifold invariants in \cite{KM}), for a cobordism \(W\) from the
 empty set to \(Y\), \(\hat{m}[\uu] (W)\in \hat{C}(Y)\) can be defined as
\[
\hat{m}[\uu] (W)=\hat{m}[\uu] (\dot{W})\hat{1}, 
\]
where \(\dot{W}\) is a cobordism from \(S^3\) to \(Y\) obtained by
removing a 3-ball from the interior of \(W\), and \(\uu\in \op{C}
(\B^\sigma  (\dot{W});\bbK)\) is used to
denote the cochain induced from that in \(\op{C}(\B^\sigma  (W);\bbK)\). 
With \(W_1\simeq S^1\times B^3\) and \(W_1'\simeq D^2\times
 S^2\) endowed with the metrics prescribed above, the values 
\(\hat{m}[\uu_{\gamma _0}](S^1\times B^3)\) and \(\hat{m}[1](D^2\times
 S^2)\) are also well-known (and follow from simple computations): In
 the notation of (\ref{C(Z)}), 
\[\begin{split}
\hat{m}[\uu_{\gamma _0}](S^1\times B^3)& =\hat{1}_Z\in \hat{C}(S^1\times
S^2),\quad \hat{m}[1](S^1\times B^3)=y\hat{1}_Z\in \hat{C}(S^1\times
S^2);\\
\hat{m}[1](D^2\times S^2) & =\hat{1}_Z\in \hat{C}(S^1\times
S^2). 
\end{split}
\]
Inserting these into (\ref{comp:W-Z}) and (\ref{comp:W-Z1}), we have: 
\begin{equation}\label{comp:W-Z2}
\hat{m}[\uu_{\lambda_\#}](W_\#)=\hat{m}[1](W'_\#)+\big[\hat{\textsc{k}}[\pmb{1}](W'_\#)-\hat{\textsc{k}}[\mathbf{u}_{\lambda
  _\#}](W_\#), \partial_\#\big].
\end{equation}
As observed before, \(W'_\#\simeq\bbR\times M_\#\). When the latter is
endowed with cylindrical metric and perturbation, \(\hat{m}[1](\bbR\times
M_\#)=\op{Id}\). Thus, since \(\delta (1)=0\),  by (\ref{eq:var-metric})
again, \[\hat{m}[1](W'_\#)=\op{Id}_\#-\big[\hat{\smz}[\pmb{1}](W'_\#), \partial_\#
 \big],\]
where \(\hat{\smz}[\pmb{1}](W'_\#)\) is defined using a path of
metrics/perturbations from the original version to the cylindrical
one.   Combining this with
(\ref{comp:W-Z2}), we arrive at (\ref{eq:mW1}). 




\paragraph{\it  \((v)\) Verifying (\ref{eq:mW}): Preparations.}
Consider Figure
\ref{Figure (a)} (b) and write \(\hat{m}[\uu](W_1)\) as a map from
\(\hat{C}(M_1)\) to \(\hat{C}(M_1)\otimes \hat{C}(S^3)\);
\(\hat{m}[\uu](W_2)\) as a map from \(\hat{C}(S^3)\otimes
\hat{C}(M_2)\)  to 
\(\hat{C}(M_2)\). 
 To simplify notation, we denote \(\lambda _{\sqcup +}, \lambda _{\sqcup
  -} \subset W_\sqcup\) respectively as \(\lambda , \bar{\lambda }\)
below. Recall that the point \(x=\lambda \cap Y_0\) (which is in \(W_\sqcup\)) separates \(\lambda \) 
into \(\lambda _1\cup \lambda _2\), with \(\lambda _i\subset W_i\) for
each \(i=1,2\). We denote the point of intersection of \(\bar{\lambda
}\) with \(Y_0\) as \(\bar{x}\), and  \(\bar{\lambda }=\bar{\lambda
}_1\cup \bar{\lambda }_2\). Recall also the arcs \(\gamma _1\subset
\bar{\grB}_c\), \(\gamma _2\subset \grB_c\) from Step 3. (In Step 3,
\(\bar{\grB}_c\), \(\grB_c\) were respectively denoted \(B_1\) and
\(B_2\). In the surgered manifold
\(W'_\sqcup=\hat{W}_1\sqcup\hat{W}_2\), for each \(i=1, 2\) \(\gamma
_i\) join with \(\lambda _i\cup\bar{\lambda }_i\) at \(\{x,
\bar{x}\}\subset S^3\) to form paths in \(\hat{W_i}\simeq \bbR\times
  M_i\). We denote the path in \(\hat{W}_1\) by \(-\hat{p}_1\) and
  that in \(\hat{W}_2\) by \(\hat{p}_2\), as they are diffeomorphic
  respectively to
  the paths \(\bbR\times \{-p_1\}\subset \bbR\times M_1\),
  \(\bbR\times \{p_2\}\subset \bbR\times M_2\) under suitable
  diffeomorphisms taking \(\hat{W}_i\) to \(\bbR\times M_i\). See
  Figure \ref{Figure (a)} (b\('\)). 

We begin with some computions of \(\hat{m}[\uu] (W_1)\). Express 
 \(\hat{m}[\uu](W_1)\) in block form as in  (\ref{def:m-lambda}). 
First, note the following facts:  \(\grC(M_1)=\grC^o(M_1)\) and
\(\grC^o(S^3)=\emptyset\); \(\bar{m}^\#_{\flat} (W_1)\) vanishes
for all \(\#, \flat\) while \(m^\#_{\flat} (W_1)\) vanishes except
when \(\#=o\) and \(\flat=os\). This means that, when 
 \(\hat{m}[\uu](W_1)\) is given by the simpler formula 
 (\ref{def:m-map}) (e.g. when \(\uu=1\)), only one 
term, \(-(\partial^o_o(M_1)\otimes \bar{\partial}^s_u(S^3)) \, 
m^o_{os}[\uu](W_1) \) on the lower left, can be non-vanishing. However,
\(\bar{\partial}^s_u(S^3)=0\). Thus, \(\hat{m}[\uu] (W_1)=0\) for such
\(\uu\). 


More generally, the lower row of  (\ref{def:m-lambda}) contains additional
 terms as explained in Remark \ref{rem:generalized-m}. These
 correspond to the last terms in both lines of (\ref{def:m-lambda0}). 
 According to the discussion following  (\ref{def:m-lambda0}),  in the
 cases of  \(\hat{m}[\theta _{\lambda _1}](W_1)\) and \(\hat{m}[\theta
 _{\bar{\lambda }_1}](W_1)\), these terms are of the same form as
 those in  (\ref{def:m-lambda0}), but with the map
 \(\bar{n}^s_u[d\ul{\op{h}}_{\hat{p}_2}]\) therein replaced
 respectively by: 
\[\begin{split}
 1\otimes
\bar{n}^s_u[d\ul{\op{h}}_{\hat{x}}](S^3)-\bar{n}^s_u[d\ul{\op{h}}_{\hat{p}_1}](M_1)\otimes
1=\bar{n}^s_u[d\ul{\op{h}}_{\hat{x}}](S^3) & \qquad
\text{for \(\hat{m}[\theta _{\lambda _1}](W_1)\)}, \\ 
1\otimes
\bar{n}^s_u[d\ul{\op{h}}_{\hat{\bar{x}}}](S^3) & \qquad
\text{for \(\hat{m}[\theta _{\bar{\lambda }_1}](W_1)\)}. 
\end{split}
\]
In the present context, the additonal term in the lower right corner of (\ref{def:m-lambda}),
being a product of one of the expressions above with \(
m^u_{os}[1](W_1)\), vanishes because the latter does. Thus, for both
\(\hat{m}[\theta _{\lambda _1}](W_1)\) and
\(\hat{m}[\theta_{\bar{\lambda }_1}](W_1)\), all entries in
(\ref{def:m-lambda0}) vanish except possibly for the lower left entry,
which is respectively 
\[ \begin{split}
\hat{m}^o_{ou}[\theta _{\lambda _1}](W_1)& =
\bar{n}^s_u[d\ul{\op{h}}_{\hat{x}}](S^3)\, m^o_{os}[1](W_1),\\
\hat{m}^o_{ou}[\theta_{\bar{\lambda }_1}](W_1) &
=\bar{n}^s_u[d\ul{\op{h}}_{\hat{\bar{x}}}](S^3)\, m^o_{os}[1](W_1). 
\end{split}
\]
Now, for any \(p\in S^3\), 
\[\begin{split}
\langle\grc', \bar{n}^s_u[d\ul{\op{h}}_{\hat{p}}](S^3)\,
\grc\rangle & =\langle\grc', (\bar{U}_p)^s_u(S^3)\, \grc\rangle\\
& =
\begin{cases}
1 & \text{when \(\grc=\check{1}\in \grC^s(S^3)\)
and \(\grc'=\hat{1}\in \grC^u(S^3)\),}\\
0 & \text{otherwise}.
\end{cases}
\end{split}
\]
Thus, \(\langle\grc_{+,1}\otimes \grc, \hat{m}[\theta _{\lambda
  _1}](W_1), \grc_{-,1}\rangle=\langle\grc_{+,1}\otimes \grc,
\hat{m}[\theta _{\bar{\lambda  }_1}](W_1), \grc_{-,1}\rangle\) vanishes
for all \(\grc\in \grC(S^3)\) and \(\grc_{\pm,1}\in \grC(M_1)=\grC^o(M_1)\)
except when \(\grc=\hat{1}\), in which
case 
\[\begin{split}
\langle\grc_{+,1}\otimes \hat{1}, \hat{m}[\theta _{\lambda_1}](W_1),
\grc_{-,1}\rangle& =\langle\grc_{+,1}\otimes \hat{1}, \hat{m}[\theta _{\bar{\lambda
  }_1}](W_1), \grc_{-,1}\rangle\\ 
& =\langle\grc_{+,1}\otimes \check{1},
\check{m}[1](W_1), \grc_{-,1}\rangle. 
\end{split}
\]
Note that \(\bar{n}[\op{u}](S^3)\) vanishes for all \(\op{u}\) of odd degrees,
since all pairs of \(\grc, \grc'\in \grC(S^3)\), the difference 
\(\bar{\op{gr}}(\grc')-\bar{\op{gr}}(\grc)\) is even. Together with
preceding arguments, this implies that in the block form
(\ref{def:m-lambda0}), all entries of \(\hat{m}[\theta _{\bar{\lambda}_1}\wedge
\theta _{\lambda _1}](W_1))\) also  vanish except possibly for the lower left entry,
which is 
\[ \begin{split}
\hat{m}^o_{ou}[\theta _{\bar{\lambda}_1}\wedge\theta _{\lambda
  _1}](W_1) & =\bar{n}^s_u[d\ul{\op{h}}_{\hat{\bar{x}}}](S^3)\,  m^o_{os}[\theta _{\lambda_1}](W_1)+\bar{n}^s_u[d\ul{\op{h}}_{\hat{x}}](S^3) \big) \,
m^o_{os}[\theta _{\bar{\lambda}_1}](W_1)\\
& = m^o_{os}[\theta
_{\lambda_1}](W_1)-+m^o_{os}[\theta _{\bar{\lambda}_1}](W_1), 
\end{split}
\]
and  \(\langle\grc_{+,1}\otimes \grc, \hat{m}[\theta _{\bar{\lambda  }_1}\wedge\theta _{\lambda
  _1}](W_1), \grc_{-,1}\rangle\) vanish
for all \(\grc\in \grC(S^3)\) and \(\grc_{\pm,1}\in \grC(M_1)=\grC^o(M_1)\)
except when \(\grc=\hat{1}\), in which
case
\[
\Big\langle\grc_{+,1}\otimes \hat{1}, \hat{m}[\theta _{\bar{\lambda}_1}\wedge
\theta _{\lambda _1}](W_1) \, \grc_{-,1}\Big\rangle =
\Big\langle\grc_{+,1}\otimes \check{1}, \check{m}[\theta _{\lambda
  _1}+\theta _{\bar{\lambda}_1}](W_1) \, \grc_{-,1}\Big\rangle  .
\]
Imitating physicists' notation, we use
\(\hat{m}[\uu](W_1)|\grc\rangle\), \(\grc\in \grC ^u(S^3)\) to denote
the map from \(\hat{C}(M_1)=\check{C}(M_1)\) to itself defined by 
\[
\Big\langle \grc_{+,1}, \hat{m}[\uu ](W_1) |\grc\rangle\, (\grc_{-,1})\Big\rangle
=\Big\langle \grc_{+,1}\otimes \grc, \hat{m}[\uu ](W_1)(\grc_{-,1})\Big\rangle. 
\]
Similarly, \(\check{m}[\uu](W_1)|\grc\rangle\), \(\grc\in \grC
 ^s(S^3)\) will denote
the map from \(\hat{C}(M_1)=\check{C}(M_1)\) to itself defined by 
\[
\Big\langle \grc_{+,1}, \check{m}[\uu ](W_1) |\grc\rangle\, (\grc_{-,1})\Big\rangle
=\Big\langle \grc_{+,1}\otimes \grc, \check{m}[\uu
](W_1)(\grc_{-,1})\Big\rangle. 
\]
Also, \(\langle\grc|\hat{m}[\uu ](W_2)\rangle\) will denote a map from
\(\hat{C}(M_2)\) to itself given by 
\[
\Big\langle\grc_{+,2},
 \langle\grc| \hat{m}[\uu ](W_2) \, (\grc_{ -,2})\Big\rangle=\Big\langle\grc_{+,2},
 \hat{m}[\uu ](W_2) \, (\grc\otimes \grc_{ -,2})\Big\rangle. 
\]

\paragraph{\it  (vi) Verifying 2) of (\ref{eq:mW}):}
In this case, the cochain \(\uu=1\in \op{C}(\B^\sigma
_{loc}(W_\sqcup);\bbK)\) can be written as an inner product,
\(1=c(1\otimes 1)\), and (\ref{KM-composition0}) applies.  
Noting that \(\delta (1)=0\), in this context
(\ref{KM-composition0}) gives, with respect to the factorization \(\hat{C}(M_\sqcup)=\hat{C}(M_1)\otimes
\hat{C}(M_2)\)
\[\begin{split}
 & \sum_{\grc\in \grC  ^u(S^3)}  \big(\hat{m}[1](W_1)
 |\grc\rangle\big) \otimes  \big( \langle \grc| \hat{m}[1](W_2)\big) \\
& \quad \quad =
 \hat{m}[1](W_\sqcup)+\big[\hat{\textsc{k}}[\mathbf{1}](W_\sqcup),
 \hat{\partial}_\sqcup, \big]. 
\end{split}
\]
The left hand side of the preceding formula vanishes, since we saw
that  \(\hat{m}[1] (W_1)=0\).  This directly leads to 2) of
(\ref{eq:mW}), with \(\textsc{b}\) set to be
\[
\textsc{b}=-\hat{\textsc{k}}[\mathbf{1}](W_\sqcup). 
\]

\paragraph{\it (vii) Verifying 1) and 4) of (\ref{eq:mW}):} 
As with in the proof of
(\ref{eq:ch-htpy}), in cases 1), 4) and 3), the relevant cochains do 
not take the simple form as an inner product under the decomposition,
and instead of (\ref{KM-composition0}),  
the more delicate formula (\ref{KM-composition}) is required. 
To begin, We again re-express the formulas 1) 3), and 4) of
(\ref{eq:mW}) in terms of the more concrete \(\theta _\lambda \) and
\(\theta _{\bar{\lambda }}\), making use of  the identities
(\ref{def:m-lambda1})  as well
as the previously established item 2) of (\ref{eq:mW}):  
\begin{equation}\label{eq:mW2}
\begin{split}
1')&  \quad \op{Id}_\sqcup- [\hat{\smz}_\sqcup, \hat{\partial}_{\sqcup}]_{\text{odd}}=\hat{m}[\theta_{\bar{\lambda}}](
  W_\sqcup)+[\dot{\textsc{a}}, \hat{\partial}_\sqcup]_{\text{even}}+\textsc{b}\, \hat{n}[d\ul{\op{h}}_{\sqcup}],\\ 
3') &  \quad  [\hat{\partial}_\sqcup, \smx ]_\ev
-[\hat{n}[d\ul{\op{h}}_{\sqcup}], \smz_\sqcup ]_{\ev}\\
& \qquad =\hat{m}[\theta_{\bar{\lambda}}\wedge\theta_\lambda](
  W_\sqcup) + [\dot{\textsc{c}}, \hat{\partial}_\sqcup]_{\text{odd}}-
  \hat{n}[d\ul{\op{h}}_{\sqcup}] \, \dot{\textsc{a}} +\dot{\textsc{d}}
  \, \hat{n}[d\ul{\op{h}}_{\sqcup}], \\
4') & \quad \op{Id}_\sqcup- [\hat{\smz}_\sqcup, \hat{\partial}_{\sqcup}]_{\text{odd}}=\hat{m}[\theta_{\lambda}](
  W_\sqcup) +[\hat{\partial}_\sqcup, \dot{\textsc{d}}]_{\text{even}}-\hat{n}[d\ul{\op{h}}_{\sqcup}]\,\textsc{b} ,
\end{split}
\end{equation}
where
\[
\hat{n}[d\ul{\op{h}}_{\sqcup}] (M_\sqcup):=1\otimes
\hat{n}[d\ul{\op{h}}_{\hat{p}_2}](M_2)-\hat{n}[d\ul{\op{h}}_{\hat{p}_1}](M_1)\otimes
1, 
\]
and \(\dot{\textsc{a}}, \dot{\textsc{c}},
\dot{\textsc{d}}\) are related to \(\textsc{a}\), \(\textsc{b}\),
\(\textsc{c}\), \(\textsc{d}\) via formulas parallel to those in
(\ref{def:A-dot}) relating \(\dot{A}\), \(\dot{C}\), \(\dot{D}\) to
\(A, B, C, D\).  

We omit the proof of \(1')\) above as it proof is entirely
  parallel to that for \(4')\). To
proceed with the proof of the latter, we first 
express \(\theta _{\lambda }\) as a sum of inner products in parallel
to (\ref{theta-decomp0}). Namely, we claim that in this context, 
\begin{equation}\label{theta-decomp1}
\theta _{\lambda }=c(\theta _{\lambda _1}\otimes 1)+c(1\otimes \theta _{\lambda _2}).
\end{equation}
 
Once this is established, applying 
(\ref{KM-composition}) in this case would yield:
\begin{equation}\label{decomp:W1-2}
\begin{split}
& \sum_{\grc\in \grC^u(S^3)}\Big( \hat{m}[\theta
_{\lambda_1}](W_1)|\grc\rangle\otimes \langle \grc|\hat{m}[1](W_2) 
+\hat{m}[1](W_1)|\grc\rangle\otimes \langle \grc|\hat{m}[\theta_{\lambda _2}](W_2)\Big) \\
& \qquad =\hat{m}[\theta _\lambda
](W_\sqcup)+\big[\hat{\smk}[\pmb{\theta }_\lambda ](W_\sqcup),
\hat{\partial}_\sqcup\big]_{\text{even}}+ \hat{n}[d\ul{\op{h}}_\sqcup]\,
\hat{\smk}[\pmb{1}](W_\sqcup)\\
& \qquad =\hat{m}[\theta _\lambda
](W_\sqcup)+\big[\hat{\smk}[\pmb{\theta }_\lambda ](W_\sqcup),
\hat{\partial}_\sqcup\big]_{\text{even}}-\textsc{b} \,
\hat{\smk}[\pmb{1}](W_\sqcup). 
\end{split}
\end{equation}
According to the computation of \(\hat{m}[\uu](W_1)\) in part (v) above, the
second term on the left hand side of the above vanishes, while the first term is given by 
\begin{equation*}
\begin{split}
& \hat{m}[\theta
_{\lambda_1}](W_1)|\hat{1}\rangle\otimes \langle \hat{1}|\hat{m}[1](W_2) 
=\check{m}[1](W_1)|\check{1}\rangle\otimes \langle
\hat{1}|\hat{m}[1](\hat{W}_2). 
\end{split}
\end{equation*}
Filling 3-balls at the \(S^3\)-end of \(W_1\) and \(W_2\) to get
\(\hat{W}_1\) and \(\hat{W}_2\) as shown in Figure \ref{Figure (a)}
(b\('\)), we see that  
\begin{equation}\label{claim:W-hatW}
\begin{split}
\check{m}[1](W_1)|\check{1}\rangle &
=\check{m}[1](\hat{W}_1)+[\check{\smk}[\pmb{1}](\hat{W}_1),
\check{\partial}(M_1)]_{\text{odd}}\\
& =\hat{m}[1](\hat{W}_1)+[\hat{\smk}[\pmb{1}](\hat{W}_1),
\hat{\partial}(M_1)]_{\text{odd}};\\
\langle \hat{1}|\hat{m}[1](W_2) & = \hat{m}[1](\hat{W}_2)+[\hat{\smk}[\pmb{1}](\hat{W}_2),
\hat{\partial}(M_2)]_{\text{odd}}. 
\end{split}
\end{equation}
Combining these with (\ref{decomp:W1-2}), we have: 
\begin{equation}\label{decomp:W1-2a}\begin{split}
& \hat{m}[1](W'_\sqcup)=\hat{m}[1](\hat{W}_1)\otimes
\hat{m}[1](\hat{W}_2)\\ 
& \quad =\hat{m}[\theta _\lambda
](W_\sqcup)+\big[\hat{\smk}[\pmb{\theta }_\lambda ](W_\sqcup),
\hat{\partial}_\sqcup\big]_{\text{even}}-\textsc{b} \,
\hat{\smk}[\pmb{1}](W_\sqcup)-[\hat{\smk}[\pmb{1}](W'_\sqcup),
\hat{\partial}_\sqcup]_{\text{odd}}. 
\end{split}
\end{equation}

Since for both \(i=1,2\) \(\hat{W}_i\simeq \bbR\times M_i\), and when equipped with
cylindrical metric and perturbation, \( \hat{m}[1](\bbR\times M_i)=\op{Id}\), by
(\ref{eq:var-metric}) we then 
have that 
\[
\hat{m}[1](W'_\sqcup)
=\op{Id}_\sqcup-\big[\hat{\smz}[\pmb{1}](W'_\sqcup), \hat{\partial}_\sqcup\big]_{\text{odd}},
\]
where \(\hat{\smz}[\pmb{1}](W'_\sqcup)\) is defined via a path of
metrics/perturbations from the original version on \(W'_\sqcup\simeq
\bbR\times M_\sqcup\) to the cylindrical version.  
Combining this with (\ref{decomp:W1-2a}), {\em modulo the proof of 
(\ref{theta-decomp1})}, we have verified item \(4'\)) of (\ref{eq:mW2}), with
\(\dot{\textsc{d}}\), \(\smz_\sqcup\) therein set respectively to be 
\begin{equation}\label{def:z-cup}
\dot{\textsc{d}}=-\hat{\smk}[\pmb{\theta }_\lambda
](W_\sqcup); \quad \smz_\sqcup= \hat{\smz}[\pmb{1}](W'_\sqcup)-\hat{\smk}[\pmb{1}](W'_\sqcup). 
\end{equation}
Item \(1'\)) of (\ref{eq:mW2}) is derived using the same arguments,
with \(\dot{\textsc{a}}\) set to be 
\[
\dot{\textsc{a}}=\hat{\smk}[\pmb{\theta }_{\bar{\lambda}}
](W_\sqcup). 
\]


We now return to the task of verifying (\ref{theta-decomp1}). This is
again done following the strategy outlined in the end of Step 3
above. Recall the definitions of the bundles
\(\tilde{\pi}_\lambda \co \tilde{\tilde{\B}}^\sigma _\lambda (W)\to
\B^\sigma (W)\), \(\pi_\lambda \co \td{\B}^\sigma _\lambda
(W)\to\B^\sigma (W)\) from the diagram (\ref{CD-ttB}). 
Let \(\tilde{\pi}_{\lambda _1}\co \tilde{\tilde{\B}}^\sigma _{\lambda
  _1}(W_1)\to \B^\sigma (W_1)\),
  \(\tilde{\pi}_{\lambda _1}\co \tilde{\tilde{\B}}^\sigma _{\lambda
    _2}(W_2)\to \B^\sigma (W_2)\) be \(U(1)\times U(1)\) bundles
  defined in a similar manner; namely, 
    by the commutative diagrams:
\[
\xymatrixcolsep{2pc}\xymatrix{
\tilde{\tilde{\mathcal{B}}}^\sigma
_{\lambda _1}(W_1(S))\ar@{->}[r]^(.4){\tilde{\tilde{\Pi}}^\partial_{\lambda _1}}
\ar@{->}[d]^{\td{\pi}_{\lambda _1}}
&\tilde{\mathcal{B}}^\sigma_{-p_1}(M_1)\times  \tilde{\B}_{\ul{x}}^\sigma (S^3)
\ar@{->}[d]^{\pi_{p_1}\times \pi_{\ul{x}}} \\
\mathcal{B}^\sigma(W_1(S))\ar@{->}[r]^(.4){\Pi^{\partial}_{\lambda _1}}
& \mathcal{B}^\sigma(M_1)\times \B^\sigma (S^3);
}
\quad 
\xymatrixcolsep{2pc}\xymatrix{
\tilde{\tilde{\mathcal{B}}}^\sigma _{\lambda _2}(W_2(S))\ar@{->}[r]^(.4){\td{\td{\Pi }}^\partial_{\lambda _2}}
\ar@{->}[d]^{\td{\pi}_{\lambda _2}}
&\tilde{\B}_{-\ul{x}}^\sigma (S^3)\times \tilde{\mathcal{B}}^\sigma_{p_2}(M_2)
\ar@{->}[d]^{\pi_{\ul{x}}\times \pi_{p_2}}\\
 \mathcal{B}^\sigma(W_2(S))\ar@{->}[r]^(.4){\Pi ^\partial_{\lambda _2}}  & \B ^\sigma (S^3)\times \mathcal{B}^\sigma(M_2),}
\]
where \(\Pi^\partial_{\lambda _i}\), \(i=1,2\), 
are defined similarly to their cousins in Section \ref{sec:2}:
\(\Pi^\partial_{\lambda _1}=\Pi ^{+M_1}_{W_1}\times\Pi ^{S^3}_{W_i}\), 
 \(\Pi^\partial_{\lambda _2}=\Pi^{S^3}_{W_2}\times
\Pi^{+M_2}_{W_2}\), with \(\Pi^{S^3}_{W_i}\) now  as defined in (\ref{CD:Pi_i}). In
the above, the sign
\(\pm\)  in the 
superscript of \(\Pi ^{\pm M_i}_{W_i}\) was introduced to distinguish between the two 
\(M_i\)-ends of \(W_i\): \(\Pi ^{\pm M_i}_{W_i}\) respectively denote the maps of
taking limits to the \(M_i\)-end at \(s\to\pm\infty\). Factor 
\(\td{\td{\Pi }}^\partial_{\lambda _1}\), \(\td{\td{\Pi
  }}^\partial_{\lambda _2}\) respectively as \(\td{\td{\Pi }}^\partial_{\lambda _1}=\td{\td{\Pi
}}_{W_1}^{+M_1}\times\td{\tilde{\Pi}}^{S^3}_{W_1}\),  \(\td{\td{\Pi
  }}^\partial_{\lambda _2}=\td{\tilde{\Pi}}^{S^3}_{W_2}\times
\td{\tilde{\Pi}}^{+ M_2}_{W_2}\). 

Let \(\td{\td{\tilde{\mathcal{B}}}}^\sigma_{x, \lambda}(W_\sqcup
(S))\) be the \(U(1)^{\times  3}\)-bundle over \(\mathcal{B}^\sigma(W_\sqcup(S))\)
  defined by the commutative diagram
\[
\xymatrixcolsep{4pc}\xymatrix{
\td{\td{\tilde{\mathcal{B}}}}^\sigma_{x, \lambda
}(W_\sqcup(S))\ar@{->}[r]^{\td{\td{\pi }}_\lambda } \ar@{->}[d]^{\td{\td{\pi
    }}_x}   \ar@{->}[rd]^{\td{\td{\td{\pi }}}_{x,\lambda }} 
& \tilde{\mathcal{B}}^\sigma_{x
}(W_\sqcup(S))
\ar@{->}[d]^{\td{\pi}_x} \\
\td{\tilde{\mathcal{B}}}^\sigma_{\lambda}(W_\sqcup(S))\ar@{->}[r]^{\td{\pi
  }_\lambda }\ar@{->}[d]^{\td{\pi }}  & \mathcal{B}^\sigma(W_\sqcup(S))\\
\tilde{\mathcal{B}}^\sigma_{\lambda}(W_\sqcup(S))\ar@{->}[ur]^{\pi_\lambda } &
}
\]
We have the following variant of (\ref{CD:W-S}) in the present context:
\begin{equation}\label{CD:Wcup-S}
\resizebox{.9 \textwidth}{!} 
{$
\minCDarrowwidth15pt
\begin{CD}
 \tilde{\tilde{\tilde{\B}}}_{x,\lambda }^\sigma (W_\sqcup(S))@>\tilde{\tilde{\ss}}_1\times \tilde{\tilde{\ss}}_2>> \tilde{\tilde{\B}}^\sigma  _{\lambda _1}(W_1(S))\times _{\tilde{\B}_{\ul{x}}(S^3)}\tilde{\tilde{\B}}^\sigma_{\lambda _2}(W_2(S)) @>
>> \tilde{\tilde{\B}}^\sigma  _{\lambda _1}(W_1(S))\times \tilde{\tilde{\B}}^\sigma_{\lambda _2}(W_2(S))\\
@V\td{\td{\tilde{\pi}}}_{x,\lambda  }VV  @V\td{\pi}_{\lambda _1}\times
\td{\pi}_{\lambda _2}VV   @V\td{\pi}_{\lambda _1}\times
\td{\pi}_{\lambda _2}VV\\
\B^\sigma (W_\sqcup(S)) @>\ss_1\times \ss_2>> \B^\sigma  (W_1(S))\times _{\B^\sigma (S^3)}\B^\sigma
(W_2(S))@>
>> \B^\sigma  (W_1(S))\times \B^\sigma
(W_2(S)).
\end{CD}
$}
\end{equation}
Fix now \(\vartheta '_{p_1}\), \(\vartheta '_{p_2}\), \(\vartheta
'_{\ul{x}}\), together with compatible trivializations \(\rho _{\vartheta '_{p_1}}\), \(\rho _{\vartheta '_{p_2}}\), \(\rho _{\vartheta
'_{\ul{x}}}\) and use them to define \(\bbR/\bbZ\)-valued functions
\(\op{h}_{\lambda }\), \(\op{h}_{\lambda _1}\), \(\op{h}_{\lambda
  _2}\), respectively on \(\td{\B}^\sigma _\lambda (W_\sqcup(S))\),
\(\td{\B}^\sigma _{\lambda _1} (W_1(S))\), \(\td{\B}^\sigma
_{\lambda _2}(W_2(S))\). Use \(\td{\td{\op{h}}}_\lambda \),
\(\td{\op{h}}_{\lambda _1}\), \(\td{\op{h}}_{\lambda _2}\) to
respectively denote their pull-backs to
\(\tilde{\tilde{\tilde{\B}}}_{x,\lambda }^\sigma (W_\sqcup(S))\),
\(\tilde{\tilde{\B}}^\sigma  _{\lambda _1}(W_1(S))\),
\(\tilde{\tilde{\B}}^\sigma_{\lambda _2}(W_2(S))\). Then arguing as
in the paragraph following (\ref{h-split}), we see that
\(\td{\td{\op{h}}}_\lambda \) agrees with the pull-back of the
function \((\td{\op{h}}_{\lambda _1}\times 1+1\times
\td{\op{h}}_{\lambda _2})\co \tilde{\tilde{\B}}^\sigma  _{\lambda _1}(W_1(S))\times
\tilde{\tilde{\B}}^\sigma_{\lambda _2}(W_2(S))\to \bbR/\bbZ\) via
\(\tilde{\tilde{\ss}}_1\times \tilde{\tilde{\ss}}_2\), and we have:
\[\begin{split}
\td{\td{\pi }}_x^*\td{\pi }^*\vartheta _\lambda & =d\td{\td{\op{h}}}_\lambda =(\tilde{\tilde{\ss}}_1\times \tilde{\tilde{\ss}}_2)^*(d\td{\op{h}}_{\lambda _1}\times 1+1\times
d\td{\op{h}}_{\lambda _2})\\
& =(\tilde{\tilde{\ss}}_1\times
\tilde{\tilde{\ss}}_2)^*\big(\tilde{\pi }^*_{W_1}\vartheta _{\lambda
  _1}\times 1 +1\times \tilde{\pi }^*_{W_2}\vartheta _{\lambda
  _2}\big). 
\end{split}
\]
Note that both sides of the preceding equation depends only on the
choices of \(\vartheta '_{p_1}\) and \(\vartheta '_{p_2}\),
independent of all other choices made to (simultaneously) define
\(\op{h}_{\lambda }\), \(\op{h}_{\lambda _1}\),
\(\op{h}_{\lambda_2}\). Meanwhile, note that 
\[\begin{split}
\td{\td{\td{\pi }}}_{x,\lambda }^*\theta _{\lambda } &=\td{\td{\pi
  }}_x^*\td{\pi }^*\vartheta _\lambda -\td{\td{\pi }}_x^*(\td{\Pi
  }_{W_\sqcup}^{+M_2})^*\vartheta '_{p_2}+\td{\td{\pi }}_x^*(\td{\Pi
  }_{W_\sqcup}^{+M_1})^*\vartheta '_{p_1};\\
\td{\pi }_{\lambda _1}^*\theta _{\lambda _1}& =\tilde{\pi }^*_{W_1}\vartheta _{\lambda _1}- (\td{\td{\Pi}}_{W_1}^{S^3})^*\vartheta '_{\ul{x}}+(\td{\td{\Pi}}_{W_1}^{+M_1})^*\vartheta '_{p_1};\\
\td{\pi }_{\lambda _2}^*\theta _{\lambda _2}& =\tilde{\pi }^*_{W_2}\vartheta _{\lambda _2}+ (\td{\td{\Pi}}_{W_2}^{S^3})^*\vartheta '_{\ul{x}}-(\td{\td{\Pi}}_{W_2}^{+M_2})^*\vartheta '_{p_2},
\end{split}
\]
and over \(\tilde{\tilde{\B}}^\sigma  _{\lambda _1}(W_1(S))\times
_{\tilde{\B}_{\ul{x}}(S^3)}\tilde{\tilde{\B}}^\sigma_{\lambda
  _2}(W_2(S))\hookrightarrow \tilde{\tilde{\B}}^\sigma  _{\lambda
  _1}(W_1(S))\times \tilde{\tilde{\B}}^\sigma_{\lambda _2}(W_2(S))\),
\[
\begin{split}
& \big(-(\td{\td{\Pi}}_{W_1}^{S^3})^*\vartheta
'_{\ul{x}}+(\td{\td{\Pi}}_{W_1}^{+M_1})^*\vartheta '_{p_1}\big)\times
1+ 1\times \big( (\td{\td{\Pi}}_{W_2}^{S^3})^*\vartheta
'_{\ul{x}}-(\td{\td{\Pi}}_{W_2}^{+M_2})^*\vartheta '_{p_2}\big)\\
& \qquad \qquad =(\td{\td{\Pi}}_{W_1}^{+M_1})^*\vartheta '_{p_1}\big)\times
1-1\times (\td{\td{\Pi}}_{W_2}^{+M_2})^*\vartheta '_{p_2}.
\end{split}
\]
Thus, 
\[
\td{\td{\td{\pi }}}_{x,\lambda }^*\theta _{\lambda }=(\tilde{\tilde{\ss}}_1\times
\tilde{\tilde{\ss}}_2)^*\big(\td{\pi }_{\lambda _1}^*\theta
_{\lambda _1}\times 1+1\times \td{\pi }_{\lambda _2}^*\theta
_{\lambda _2}\big), 
\]
and hence through (\ref{CD:Wcup-S}), 
\[
\theta _{\lambda }=(\ss_1\times\ss_2)^*\big(\theta_{\lambda _1}\times 1+1\times \theta_{\lambda _2}\big), 
\]
which means (\ref{theta-decomp1}).

\paragraph{\it  (viii) Verifying 3) of (\ref{eq:mW}):}
The composition formula (\ref{KM-composition}) in this case gives
\begin{equation}\label{decom:W-cup}
\begin{split}
 & \hat{m}[\theta _{\bar{\lambda }}\wedge\theta _{\lambda }]( W_\sqcup(\infty )) 
=\hat{m}[\theta _{\bar{\lambda }}\wedge\theta _{\lambda }](W_\sqcup)+\big[\hat{\textsc{k}}[\pmb{\theta}
_{\bar{\lambda}}\wedge \pmb{\theta }_\lambda ](W_\sqcup),
\hat{\partial}_\sqcup \big]_\odd \\
& \qquad \qquad \qquad -\hat{\textsc{k}}[\pmb{\theta}
_{\lambda}](W_\sqcup)\, \hat{n}[d\ul{\op{h}}_\sqcup](M_\sqcup)-
\hat{n}[d\ul{\op{h}}_\sqcup](M_\sqcup)\, \hat{\textsc{k}}[\pmb{\theta}
_{\bar{\lambda}}](W_\sqcup)\\
& \, \, = \hat{m}[\theta _{\bar{\lambda }}\wedge\theta _{\lambda }](W_\sqcup)+\big[\hat{\textsc{k}}[\pmb{\theta}
_{\bar{\lambda}}\wedge \pmb{\theta }_\lambda ](W_\sqcup),
\hat{\partial}_\sqcup \big]_\odd 
-\dot{\textsc{a}}\,
\hat{n}[d\ul{\op{h}}_\sqcup](M_\sqcup)+
\hat{n}[d\ul{\op{h}}_\sqcup](M_\sqcup)\,
\dot{\textsc{d}}.\\
\end{split}
\end{equation}
To compute the left hand side of the preceding formula, first use (\ref{theta-decomp1}) and its sister version for \(\theta
 _{\bar{\lambda }}\) to write: 
\[
\begin{split}
\theta _{\bar{\lambda }}\wedge\theta _{\lambda }=&\ss_1^*(\theta _{\bar{\lambda }_1}\wedge\theta _{\lambda_1
 })+ \ss_1^*\theta _{\lambda _1}\wedge\ss_2^*\theta _{\bar{\lambda
   }_2}\\
& \quad +\ss_1^*\theta _{\bar{\lambda }_1}\wedge \ss_2^*\theta _{\lambda _2}+
\ss_2^*(\theta _{\bar{\lambda }_2}\wedge\theta _{\lambda _2}).
\end{split}
\] 
Combining this with the computation of \(\hat{m}[\uu](W_1)\) in part
(v), 
we get: 
\begin{equation}\label{W-infty-cup}
\begin{split}
& \hat{m}[\theta _{\bar{\lambda }}\wedge\theta
_{\lambda}](W_\sqcup(\infty ))= 
\hat{m}[\theta _{\bar{\lambda}_1}\wedge\theta_{\lambda _1}](W_1)|\hat{1}\rangle\otimes \langle
\hat{1}| \hat{m}[1](W_2) \\
& \qquad \qquad \,\,  +\hat{m}[\theta _{\bar{\lambda}_1}](W_1)|\hat{1}\rangle\otimes \langle
\hat{1}| \hat{m}[\theta _{\lambda _2}](W_2) 
+\hat{m}[\theta _{\lambda_1}](W_1)|\hat{1}\rangle\otimes \langle
\hat{1}| \hat{m}[\theta _{\bar{\lambda }_2}](W_2)\\
& \qquad \qquad \, \, + \hat{m}[1](W_1)|\hat{1}\rangle\otimes \langle
\hat{1}| \hat{m}[\theta _{\bar{\lambda
  }_2}\wedge\theta_{\lambda _2}](W_2)\\
& \quad = \check{m}[\theta _{\bar{\lambda}_1}+\theta_{\lambda
  _1}](W_1)|\check{1}\rangle\otimes \langle
\hat{1}|\hat{m}[1](W_2)
+ \check{m}[1](\hat{W}_1)|\check{1}\rangle\otimes   \langle
\hat{1}|\hat{m}[\theta _{\bar{\lambda}_2}+\theta_{\lambda
  _2}](W_2)\\
\end{split}
\end{equation}

In comparison to the identities from decomposing \(W_\sqcup \),  we have the following identity obtained by applying 
 (\ref{KM-composition}) (and its check-version) to the decomposition
 of \(W'_\sqcup\) described in
 Figure \ref{Figure (a)} (b\('\)). Recall that for the underlying 
 decomposition of \(W'_\sqcup\), the paths \((-1)^i\hat{p}_i\) split as:
 \((-1)^i\hat{p}_i=(\bar{\lambda }_i\cup\lambda _i)\cup_{x, \bar{x}} (\gamma
 _i)\). The arguments leading to (\ref{theta-decomp0}) imply that in the
 present setting, 
\[
\theta _{(-1)^i\hat{p}_i}=c(\theta _{\lambda _i}\otimes 1)+ c(\theta_{\bar{\lambda }_i}\otimes 1)+c(1\otimes \theta _{\gamma _i}).
\]
Then 
\[\begin{split}
\hat{m}[\theta _{-\hat{p}_1}](\hat{W}_1(\infty))& = 
\check{m}[\theta _{-\hat{p}_1}](\hat{W}_1(\infty)) \\
&=\check{m}[\theta _{\lambda _1}+\theta _{\bar{\lambda
  }_1}](W_1)|\check{1}\rangle +\sum_{\grc\in\grC^s(S^3)}\check{m}[\theta _{\gamma _1}](\bar{\grB})(\grc)\,
\check{m}[1](W_1)|\grc\rangle;\\
\hat{m}[\theta _{\hat{p}_2}](\hat{W}_2(\infty)) &=\langle \hat{1}|\hat{m}[\theta _{\lambda _2}+\theta _{\bar{\lambda}_2}](W_2)+\langle \hat{m}[\theta _{\gamma _2}](\grB)|\hat{m}[1](W_2).
\end{split}
\]
Remember that \(\check{m}[\uu](\bar{\grB})\) is a map from
\(\check{C}(S^3)\) to \(\bbK\), while \(\hat{m}[\uu](\grB)\in
\hat{C}(S^3)\). However, when \(\deg(\uu)\) is odd, both
\(\check{m}[\uu](\bar{\grB})\) and \(\hat{m}[\uu](\grB)\) must vanish,
because all generators of \(\check{C}(S^3)\) (resp. \(\hat{C}(S^3)\))
are of even (resp. odd) degree. Thus, the last terms of both lines in the
preceding expression vanish. Combining these with (\ref{W-infty-cup}),
we have: 
\begin{equation}\label{comp:W-W'-infty}
\begin{split}
& \hat{m}[\theta _{\hat{p}_\sqcup}](W'_\sqcup(\infty))\\
&\quad =
\hat{m}[\theta _{-\hat{p}_1}](\hat{W}_1(\infty))\otimes
\hat{m}[1](\hat{W}_2(\infty))+\hat{m}[1](\hat{W}_1(\infty))\otimes
\hat{m}[\theta _{\hat{p}_2}](\hat{W}_2(\infty))\\
& \quad =
\check{m}[\theta _{\lambda _1}+\theta _{\bar{\lambda
  }_1}](W_1)|\check{1}\rangle \otimes \langle \hat{1}| \hat{m}[1](W_2)+\check{m}[1](W_1)|\check{1}\rangle \otimes \langle \hat{1}| \hat{m}[\theta _{\lambda _2}+\theta _{\bar{\lambda
  }_2}](W_2)\\
& \quad =\hat{m}[\theta _{\lambda _\sqcup}](W_\sqcup(\infty)),
\end{split}
\end{equation}
where \(\hat{p}_\sqcup\) denotes the 1-chain \(\hat{p}_2-\hat{p}_1\)
in \(W'_\sqcup=\hat{W}_1\sqcup \hat{W}_2\), and \(\theta
_{\hat{p}_\sqcup}:=\theta_{\hat{p}_2}-\theta_{\hat{p}_1}\). 
Now appy (\ref{KM-composition}) and (\ref{eq:var-metric}) to \(W'_\sqcup\) with \(\uu\) therein set to
be \(\theta_{\hat{p}_\sqcup}\); we get:
\[\begin{split}
& \hat{m}[\theta _{\hat{p}_\sqcup}](W'_\sqcup(\infty))\\
& \qquad =\hat{m}[\theta
_{\hat{p}_\sqcup}](W'_\sqcup)+ \big[\hat{\smk}[\pmb{\theta
}_{\hat{p}_\sqcup}](W'_\sqcup), \hat{\partial}_\sqcup\big]_\ev +\big[
\hat{\smk}[\pmb{1}](W'_\sqcup),
\hat{n}[d\ul{\op{h}}_\sqcup]\big]_\ev\\
&\qquad =\hat{m}[\theta
_{\hat{p}_\sqcup}](\bbR\times M_\sqcup)+ \big[\hat{\smk}[\pmb{\theta
}_{\hat{p}_\sqcup}](W'_\sqcup)-\hat{\smz}[\pmb{\theta
}_{\hat{p}_\sqcup}](W'_\sqcup), \hat{\partial}_\sqcup\big]_\ev \\
& \qquad \qquad +\big[
\hat{\smk}[\pmb{1}](W'_\sqcup)-\hat{\smk}[\pmb{1}](W'_\sqcup),
\hat{n}[d\ul{\op{h}}_\sqcup]\big]_\ev\\
&\qquad =\hat{m}[\theta
_{\hat{p}_\sqcup}](\bbR\times M_\sqcup)+ \big[\hat{\smk}[\pmb{\theta
}_{\hat{p}_\sqcup}](W'_\sqcup)-\hat{\smz}[\pmb{\theta
}_{\hat{p}_\sqcup}](W'_\sqcup), \hat{\partial}_\sqcup\big]_\ev -\big[
\smz_\sqcup, \hat{n}[d\ul{\op{h}}_\sqcup]\big]_\ev,
\end{split}
\]
where \(\smz_\sqcup\) is as in (\ref{def:z-cup}), and \(\hat{m}[\theta
_{\hat{p}_\sqcup}](\bbR\times M_\sqcup)=\hat{m}[\theta
_{\hat{p}_1}](\bbR\times M_1)\otimes \hat{m}[\theta
_{\hat{p}_2}](\bbR\times M_2)\) denotes the version of the 
cobordism map when the metric and peturbation form on \(\bbR\times
M_\sqcup\), as well as \(\hat{p}_i\subset \bbR\times M_\sqcup\), are
invariant under the \(\bbR\)-action. However, \(\hat{m}[\theta
_{\hat{p}_i}](\bbR\times M_i)=0\) by construction. (Recall
(\ref{vartheta-theta}) and (\ref{int-con}).)  Thus,  the
first term in the last line of the preceding formula vanishes. 
Putting all these together with (\ref{decom:W-cup}) and
(\ref{comp:W-W'-infty}), we have:
\[\begin{split}
 & \big[\hat{\smk}[\pmb{\theta
}_{\hat{p}_\sqcup}](W'_\sqcup)-\hat{\smz}[\pmb{\theta
}_{\hat{p}_\sqcup}](W'_\sqcup), \hat{\partial}_\sqcup\big]_\ev -\big[
\smz_\sqcup, \hat{n}[d\ul{\op{h}}_\sqcup]\big]_\ev\\
&\quad = \hat{m}[\theta _{\bar{\lambda }}\wedge\theta _{\lambda }](W_\sqcup)+\big[\hat{\textsc{k}}[\pmb{\theta}
_{\bar{\lambda}}\wedge \pmb{\theta }_\lambda ](W_\sqcup),
\hat{\partial}_\sqcup \big]_{\text{odd}}
-\dot{\textsc{a}}\,
\hat{n}[d\ul{\op{h}}_\sqcup](M_\sqcup)+
\hat{n}[d\ul{\op{h}}_\sqcup](M_\sqcup)\,
\dot{\textsc{d}} \\
\end{split}
\]
This implies (\ref{eq:mW2}) \(3')\) and hence also (\ref{eq:mW2}) 3),
if we set: 
\[
\dot{\textsc{c}}= \hat{\textsc{k}}[\pmb{\theta}
_{\bar{\lambda}}\wedge \pmb{\theta }_\lambda ](W_\sqcup), \quad \smx=\hat{\smz}[\pmb{\theta
}_{\hat{p}_\sqcup}](W'_\sqcup)-\hat{\smk}[\pmb{\theta
}_{\hat{p}_\sqcup}](W'_\sqcup)
\]
in these formulas. This finishes the proof of the lemma. 
\epf

\begin{remarks}
The preceding lemma has Yang-Mills analogs; see Theorem 7.16 and Corollary
7.21 of \cite{D}. 
A previous version of this article ({\tt arXiv:1204.0115v1}) contains
sketches of an alternative proof, where the underlying geometric
meanings of computations done here are clearer. 
\end{remarks}

We have now shown that \(V_*\) defines a chain-homotopy equivalence. 

\paragraph{\it Step 6.}
In this step we verify the claim that the maps \(V_*\), \(V_*^\dag\) intertwine
with the \(\mathbf{A}_\dag (M_\sqcup)\) action on \(S_{U_{\sqcup}}\big(\hat{C}_*(M_\sqcup ,
  \grs_{\sqcup }, r[w]_{\sqcup }; \Gamma_{ \sqcup})\big)\) and the 
  \(\mathbf{A}_\dag (M_\#)\) action on \(\CM_*(M_\# ,
  \grs_{\# }, r[w]_{\# }; \Gamma_{ \#})\) described in Parts 2 and
  3 of the previous subsection. More precisely, for each \(Q=U_p\),
  \(\grt_{\gamma}\), \(\gamma\in
  \{\gamma_i^{[1]}\}_i\cup\{\gamma_j^{[2]}\}_j\), we shall show that there exist
  homomorphisms \[\begin{split}Z_{Q*}\co \CM_*(M_\# ,
  \grs_{\# }, r[w_{\# }]; \Gamma_{ \#}) &\to S_{U_{\sqcup}}\big(\hat{C}_*(M_\sqcup ,
  \grs_{\sqcup }, r[w_{\sqcup }]; \Gamma_{ \sqcup})\big)\\
Z_{Q*}^\dag\co S_{U_{\sqcup}}\big(\hat{C}_*(M_\sqcup ,
  \grs_{\sqcup }, r[w_{\sqcup }]; \Gamma_{ \sqcup})\big) & \to \CM_*(M_\# ,
  \grs_{\# }, r[w_{\# }]; \Gamma_{ \#})
\end{split}
\] satisfying 
\[\begin{split}
V_*\grm_Q -S_{U_\sqcup}(\grm_Q)\, V_*& =Z_{Q*}\partial_\# +(-1)^{\deg
  Q} D_\sqcup Z_{Q*}.\\
V_*^{\dag}S_{U_\sqcup}(\grm_Q) -(-1)^{\deg Q}\grm_QV^\dag_*& =\partial_\# Z_{Q*}^\dag
  +(-1)^{\deg
  Q}Z_{Q*}^\dag D_\sqcup,
\end{split}
\]
We shall only verify the first line above, since the second line is
basically the adjoint of the first. Stated in terms of the
decomposition (\ref{eq:S}), this amounts to verifying the following
set of identities: 
\begin{equation}\label{eq:A-action-intertwine}
\begin{split}
\text{1)  } &  \quad V_0 U_p-\hat{U}_p V_0 =\hat{\partial}_\sqcup  Z_{U,0}+Z_{U,0}\partial_\#; \\
\text{2)  } &  \quad V_1U_p -\hat{U}_pV_1-\hat{K}_{U_p}V_0 =-\hat{\partial}_\sqcup  Z_{U,1}+Z_{U,1}\partial_\# +\hat{U}_\sqcup
Z_{U,0};\\
\text{3)  } &  \quad V_0\grm_{\gamma}-\grm_{\gamma} V_0 =-\hat{\partial}_\sqcup  Z_{\gamma,0}+Z_{\gamma,0}\partial_\#; \\
\text{4)  } &  \quad V_1\grm_{\gamma} +\grm_{\gamma} V_1-\hat{K}_\gamma V_0 =\hat{\partial}_\sqcup
Z_{\gamma ,1}+Z_{\gamma ,1}\partial_\# -\hat{U}_\sqcup Z_{\gamma,0}, 
\end{split}
\end{equation}
where 
\[Z_{U*}=\left[\begin{array}{c}Z_{U,0} \\Z_{U,1}
  \end{array}
\right]; \quad Z_{\gamma*}=\left[\begin{array}{c}Z_{\gamma,0} \\ Z_{\gamma,1}
  \end{array}
\right],
\]
and \(\hat{K}_{U_p}\), \(\hat{K}_\gamma\) are as defined in Part 3 of
last subsection. These are established by arguments similar to those
used to verify (\ref{eq:U=p-mor}), (\ref{eq:K-lambda}), and (\ref{eq:dF}). 

To proceed, we define
\(Z_{U*}\) and \(Z_{\gamma *}\) as follows.
Let \(\mathbf{p}\subset \mathcal{V}\) be a path so that on
\(\mathcal{V}-\mathcal{V}_c\simeq (\bbR^-\times M_\#)\cup
(\bbR^+\times M_\sqcup)\), 
\(\mathbf{p}\cap  (\mathcal{V}-\mathcal{V}_c)\) agrees with
\(\bbR^\pm\times \{p\}\) under the diffeomorphisms in (\ref{(A.9a,11)}). Suppose also that the path
\(\mathbf{p}\cup \bar{\bf p}\subset
\mathcal{V}\cup_{M_\sqcup}\bar{\cal V}=W_\#\) becomes the line
\(\bbR\times \{p\}\subset \bbR\times M_\#\simeq W_\#'\)  after the surgery of \(W_\#\)
along \(\lambda_\#\). (Equivalently, \(\bar{\bf p}\cup
\mathbf{p}\subset \bar{\cal V}\cup_{M_\sqcup}\mathcal{V}=W_\sqcup\) also 
becomes \(\bbR\times \{p\}\subset \bbR\times M_\sqcup\simeq W'_\sqcup\) after the surgery of \(W_\sqcup\)
along \(S_\sqcup\). For each \(\gamma\in
\{\gamma_i^{[1]}\}_i\cup\{\gamma_j^{[2]}\}_j\), define in a similar
fashion an embedded cylinder \(\Upsilon\subset
\mathcal{V}\) that ends at circles \(\gamma
\subset Y_\pm\) on both ends of \(\V\). 
Now set
\begin{equation}
\begin{split}
Z_{U,0} &:=
\hat{\op{K}}_{\mathbf{p}}(\mathcal{V};\Gamma_{\cal  V})
=\hat{m}[\theta_{\mathbf{p}}](\mathcal{V};\Gamma_{\cal V})+\hat{\Theta }_\sqcup V_0; \\
Z_{U,1}&:=\hat{m}[\uu_{\mathbf{p}}\mathpzc{u}_\lambda
](\mathcal{V};\Gamma_{\cal V})
=\hat{m}[\theta _{\mathbf{p}}\theta _\lambda
](\mathcal{V};\Gamma_{\cal V})+\big[\hat{m}[\theta _\lambda
](\mathcal{V};\Gamma_{\cal V}), \hat{\Theta }_\sqcup\big]-
\hat{\Theta}_\sqcup  Z_{U,0};\\
Z_{\gamma ,0} &:= \hat{m}[\textsc{f}_\Upsilon](\mathcal{V};\Gamma_{\cal  V});\\
Z_{\gamma ,1} &:=\hat{m}[\textsc{f}_\Upsilon \uu_\lambda
](\mathcal{V}; \Gamma_{\cal  V})= \hat{m}[\textsc{f}_\Upsilon \theta _\lambda
](\mathcal{V}; \Gamma_{\cal  V})+\Theta _\sqcup  Z_{\gamma ,0},
\end{split}
\end{equation} 
where \(\hat{\op{K}}_{\mathbf{p}}(\mathcal{V}; \Gamma_{\mathcal{V}})\) is as
defined in (\ref{def:K_lambda}) for \(X=\mathcal{V}\) and \(\lambda
=\mathbf{p}\) with respect to the \(X\)-morphism \(\Gamma
_{\mathcal{V}}\), and \(\textsc{f}_\Upsilon\) is as defined in
(\ref{scF_S}). 

With the preceding definitions, items 1) and 3) of
(\ref{eq:A-action-intertwine}) are direct consequences of
(\ref{eq:K-lambda}), and (\ref{eq:dF-X}).  
To derive items 2) and 4), first rewrite them in terms of the more
concrete cochains, \(\theta
_{\bf p}\), \(\theta _\lambda \), \(F_\Upsilon \), using the
now-established items 1), 3) and the identities (\ref{eq:chain-maps}):
\[\label{eq:A-action-intertwine0}
\begin{split}
2')  &  \quad \big[\hat{m}[\theta _\lambda ](\V), \hat{n}
[d\ul{\op{h}}_{\hat{p}}] \big]_{\ev}-\hat{n}[d\ulh_{\hat{p}} \wedge
d\ulh_\sqcup  ] \, \hat{m}[1](\V)\\
& \qquad \qquad =\big[\hat{m}[\theta _{\bf p}\wedge \theta _\lambda
](\V), \hat{\partial}\big]_\ev +\hat{n}[d\ulh_\sqcup  ]
(\V)\, \hat{m}[\theta _{\bf p}](\V);\\
4')  &  \quad \big[\hat{m}[\theta _\lambda ](\V), \hat{n}
[\op{u}_\gamma ]\big]_{\odd}-\hat{n}[\op{u}_\gamma \wedge
d\ulh_\sqcup  ]\, \hat{m}[1](\V)\\
& \qquad \qquad =\big[\hat{m}[\scf_\Upsilon \theta _\lambda ](\V),
\hat{\partial}\big]_\odd  -\hat{n}[d\ulh_\sqcup  ]
(\V)\, \hat{m}[\scf_\Upsilon ](\V).  
\end{split}
\]
where \(\op{u}_\gamma \) is the 0-cochain on \(\B^\sigma
_{loc}(Y_\pm)\) introduced in Part 2(b) of Section
\ref{sec:A-module}. Recall that in our case \(\gamma \) is used to denote both
embedded circles in \(Y_+\) and \(Y_-\); we use the same notation \(\op{u}_\gamma
\) for the corresponding cochains on \(\B^\sigma
_{loc}(Y_+)\) and \(\B^\sigma
_{loc}(Y_-)\). The same convention will be applied to the \(Y_+\)- and
\(Y_-\)'s versions of other cochains  associated to \(\gamma
\) that were constructed in Section \ref{sec:2}. 
 
According to Proposition 25.3.4 of \cite{KM}), verifying the identities
\(2')\), \(4')\) above is equivalent to verifying that: 
\begin{equation}\label{eq:A-action-intertwine0}
\begin{split}
\hat{m}[d(\theta _{\bf p}\wedge \theta _\lambda
) ](\V) & =\big[\hat{m}[\theta _\lambda ](\V), \hat{n}
[d\ul{\op{h}}_{\hat{p}}] \big]_{\ev}\\
& \quad -\hat{n}[d\ulh_{\hat{p}} \wedge
d\ulh_\sqcup  ] \, \hat{m}[1](\V)-\hat{n}[d\ulh_\sqcup  ]
(\V)\, \hat{m}[\theta _{\bf p}](\V);\\
\hat{m}[d (\scf_\Upsilon \theta _\lambda )](\V) &=\big[\hat{m}[\theta _\lambda ](\V), \hat{n}
[\op{u}_\gamma ]\big]_{\odd}\\
& \quad -\hat{n}[\op{u}_\gamma \wedge
d\ulh_\sqcup  ]\, \hat{m}[1](\V)+\hat{n}[d\ulh_\sqcup  ]
(\V)\, \hat{m}[\scf_\Upsilon ](\V).
\end{split}
\end{equation}

The rest of the this Step is devoted to verifying the preceding
identities. 

\paragraph{\it (i) Verifying the first line in (\ref{eq:A-action-intertwine0}):} We argue
similarly to (\ref{value1}). Let \(\bar{\M}\) denote a 3-dimensional moduli
space of the form \(\bar{\M}_{3,z}(\mathcal{V}; \grc_-,
\grc_+)\). Let \(\varsigma \co \bar{\M}\to\B^\sigma (\mathcal{V})
\subset \B^\sigma _{loc}(\V)  \) denote the embedding. Let \(\grM\), \(\M^+\) denote 
respectively the top dimensional stratum of \(\bar{\M}\) and \(\grr^{-1}\bar{\M}\) as
in Section \ref{sec:A-module}. The coefficients of the map on the left
hand side are given by integrals of the form 
\(\langle d(\theta _{\bf p}\wedge \theta _\lambda), \bar{\M}\rangle \). 
To compute them, let \(\td{\tilde{\mathcal{B}}}^\sigma _{{\bf p}\cup\lambda
}(\mathcal{V})\to \B^\sigma  (\V)\) be the \(U(1)\times U(1)\)-bundle
defined by the following commutative diagram: 
 \[\xymatrixcolsep{5pc}\xymatrix{
\td{\tilde{\mathcal{B}}}^\sigma _{{\bf p}\cup\lambda }(\mathcal{V})\ar@{->}[r]^{\pi'_{\mathbf{p}}}
\ar@{->}[d]^{\pi'_{\lambda }}
\ar@{->}[rd]^{\pi_{{\bf p}\cup\lambda }}&\tilde{\mathcal{B}}^\sigma_{\lambda }({\cal
  V})
\ar@{->}[d]^{\pi_{\lambda }} \\
 \tilde{\mathcal{B}}^\sigma_{\mathbf{p}}(\mathcal{V})\ar@{->}[r]^{\pi_{\mathbf{p}}}  & \mathcal{B}^\sigma({\cal V}).
}
\]
Similarly to (\ref{CD-M+}), we
shall 
choose a map \(\tilde{\tilde{\varsigma }}\co \M^+ \to  \tilde{\tilde{\mathcal{B}}}^\tau(\V)\) so that the diagram below commutes: 
\begin{equation}\label{CD:Z_U}
\xymatrix{
\M^+ \ar@{->}[rr]|{\;\tilde{\tilde{\varsigma }}\;}
\ar@{->}[dr]|{\;\tilde{\varsigma }_\lambda \;} \ar@{->}[drrr]|(0.25){\;\tilde{\varsigma }_{\mathbf{p}}\;} \ar@{->}[dd]^{\grr } & &
 \td{\tilde{\mathcal{B}}}^\sigma _{{\bf p}\cup\lambda }(\mathcal{V})
 \ar@{->}[dd]^{\pi  _{\mathbf{p}\cup \lambda }}|(0.325)\hole \ar@{->}[ld]^(0.55){\pi'_{\mathbf{p}}}|(0.475)\hole \ar@{->}[rd]^{\pi'_\lambda }\\
& \tilde{\mathcal{B}}_\lambda ^\sigma (\mathcal{V})
\ar@{->}[rd]^{\pi_\lambda } & & \tilde{\mathcal{B}}_{\bf p}^\sigma (\mathcal{V})
\ar@{->}[ld]^{\pi_{\bf p}}\\
\bar{\M} \ar@{->}[rr]^{\varsigma } && \B^\sigma (\mathcal{V}) 
}
\end{equation}

Consider the form \(((\pi _\lambda ')^*\vartheta _{\bf p})\wedge ((\pi
_{\bf p}')^*\vartheta _\lambda )\) on
\(\td{\tilde{\mathcal{B}}}^\sigma _{{\bf p}\cup\lambda
}(\mathcal{V})\). Use the identites \(\vartheta _{\bf p}=\pi _{\bf
  p}^* \theta _{\bf p}+(\td{\Pi }^{\infty})^*\vartheta '_p-(\td{\Pi}
^{-\infty})^*\vartheta '_p\) and \(\vartheta _{\lambda }=\pi _{\lambda
}^* \theta _\lambda +(\Pi ^\partial_\lambda )^*\vartheta '_{p_2-p_1}\)
to rewrite it as 
\[\begin{split}
& ((\pi _\lambda ')^*\vartheta _{\bf p})\wedge ((\pi_{\bf p}')^*\vartheta
_\lambda )
 =\pi _{\mathbf{p}\cup\lambda }^* (\theta _{\bf p}\wedge \theta
_\lambda )\\
& \quad +\big((\pi _\lambda ')^*(\td{\Pi }^{\infty})^*\vartheta '_p-(\pi _\lambda ')^*(\td{\Pi}
^{-\infty})^*\vartheta '_p\big) \wedge\pi _{\mathbf{p}\cup\lambda }^*
\theta _\lambda +(\pi _{\mathbf{p}\cup\lambda }^* \theta _{\bf
  p})\wedge(\pi _{\bf p}')^*(\td{\Pi }^\partial_\lambda )^*\vartheta '_{p_2-p_1}\\
&\quad 
+\big( (\pi _\lambda ')^*(\td{\Pi }^{\infty})^*\vartheta '_p-(\pi _\lambda ')^*(\td{\Pi}
^{-\infty})^*\vartheta '_p \big)\wedge(\pi _{\bf p}')^*(\td{\Pi }^\partial_\lambda
)^*\vartheta '_{p_2-p_1}. 
\end{split}
\]
Recall that the same notation \(p\) is used to denote corresponding
points in both \(Y_-\) and \(Y_+\). In the above, the same notation
\(\vartheta '_p\) is used to denote either the \(Y_-\) or the
\(Y_+\)'s version. 

Now pull back the preceding identity by \(\td{\td{\varsigma }}\) and
integrate over \([(\M^+)_2]=\partial [\M]\). (Here we again used \cite{KM}'s
Theorem 24.7.2 and Lemma 31.3.1. ) The integral over the left hand
side vanishes, because both \(\vartheta _\lambda \) and \(\vartheta
_{\bf p}\) are exact. Meanwhile, by way of the commutative diagram (\ref{CD:Z_U}) and Stokes' theorem, the integral over the first term on
the right hand sides is: 
\[
\langle \td{\td{\varsigma }}^*\pi _{\mathbf{p}\cup\lambda }^* (\theta _{\bf p}\wedge \theta
_\lambda ), \partial [\M^+]\rangle= \langle \td{\td{\varsigma }}^*\pi _{\mathbf{p}\cup\lambda }^* d(\theta _{\bf p}\wedge \theta
_\lambda ), [\M^+]\rangle=\langle d(\theta _{\bf p}\wedge \theta
_\lambda ), \bar{\M}\rangle. 
\]
This is exactly the coefficients of the map \(\hat{m}[\delta  (\theta _{\bf p}\wedge \theta
_\lambda )](\V)\) that we aim to compute. With a bit of diagram
chasing, the aforementioned integral identity then becomes 
\begin{equation}\label{int:Z_U}
\begin{split}
& \langle d(\theta _{\bf p}\wedge \theta
_\lambda ), \bar{\M}\rangle\\
& \quad =-\Big\langle \td{\varsigma }_{\bf p}^*\Big(\big((\td{\Pi }^{\infty})^*\vartheta '_p-(\td{\Pi}
^{-\infty})^*\vartheta '_p\big) \wedge\pi _{\mathbf{p}}^*
\theta _\lambda \Big), (\M^+)_2\Big\rangle\\
& \quad \quad -  \Big\langle  \td{\varsigma }_\lambda ^*\Big((\pi _\lambda ^* \theta _{\bf
  p})\wedge(\td{\Pi }^\partial_\lambda )^*\vartheta '_{p_2-p_1}\Big), (\M^+)_2\Big\rangle\\
&\quad \quad
-\Big\langle \td{\td{\varsigma }}^*\Big(\big( (\pi _\lambda ')^*(\td{\Pi }^{\infty})^*\vartheta '_p-(\pi _\lambda ')^*(\td{\Pi}
^{-\infty})^*\vartheta '_p \big)\wedge(\pi _{\bf p}')^*(\td{\Pi }^\partial_\lambda
)^*\vartheta '_{p_2-p_1}\Big), (\M^+)_2\Big\rangle. 
\end{split}
\end{equation}
According to \cite{KM}, \((\M^+)_2\) is a union of product
spaces of the forms:  
\begin{equation}\label{dec:M^+2}\begin{split}
& \N^+_0(Y_-\grc_-, \grc)\times \M_2 (\V; \grc_,
\grc_+), \quad \M_2 (\V; \grc_-, \grc)\times \N^+_0(Y_+; \grc, \grc_+), \\
&  \quad \quad \N^+_1(Y_-;\grc_-, \grc)\times \M_1 (\V; \grc_,
\grc_+), \quad \M_1 (\V; \grc_-, \grc)\times \N^+_1(Y_+; \grc,
\grc_+), \\
&  \quad \quad \N^+_2(Y_-;\grc_-, \grc)\times \M_0 (\V; \grc_,
\grc_+), \quad \M_0 (\V; \grc_-, \grc)\times \N^+_2(Y_+; \grc,
\grc_+),.
\end{split}
\end{equation}
The diagram requires that each \(\N^+\) factor of the preceding
product spaces must map to fibers of the bundles \(\td{\B}^\sigma _{\bf p}(\V)\),
\(\td{\B}^\sigma _\lambda (\V)\), or \(\td{\td{\B}}^\sigma
_{\mathbf{p}\cup \lambda } (\V)\), respectively under \(\td{\varsigma
}_{\bf p}\), \(\td{\varsigma }_\lambda \), and \(\td{\td{\varsigma
  }}\). 
Meanwhile, observe that on the right hand side of (\ref{int:Z_U}), the
first, second, and third term has respectively a factor of \((\td{\Pi }^{\pm\infty})^*\vartheta
'_p\), \((\td{\Pi }^\partial_\lambda)^*\vartheta '_{p_2-p_1}\), and  \((\pi _\lambda ')^*(\td{\Pi}
^{\pm\infty})^*\vartheta '_p \big)\wedge(\pi _{\bf p}')^*(\td{\Pi }^\partial_\lambda
)^*\vartheta '_{p_2-p_1}\). They restrict respectively to the volume
forms on the fibers of, respectively, \(\td{\B}^\sigma _{\bf p}(\V)\),
\(\td{\B}^\sigma _\lambda (\V)\), and \(\td{\td{\B}}^\sigma
_{\mathbf{p}\cup \lambda } (\V)\).  This means that 
product spaces of the types in the first line of (\ref{dec:M^+2}) never contribute to the integrals
on the right hand side of (\ref{int:Z_U}); those of the types in the second line of 
(\ref{dec:M^+2}) contribute only to the integrals in the first and
second terms 
on the right hand side of (\ref{int:Z_U}); those of the types in the third line of 
(\ref{dec:M^+2}) contribute only to the integrals in the last term 
on the right hand side of (\ref{int:Z_U}). Consequently, 
\begin{equation}\label{int:Z_U0}
\begin{split}
& \langle d(\theta _{\bf p}\wedge \theta
_\lambda ), \bar{\M}\rangle\\
& \quad =
-\sum_\grc\big\langle\theta _\lambda, \M_1 (\mathcal{V}; \grc-,
\grc)\big\rangle\, \big\langle d\ulh_p, \N^+_1(Y_+, \grc,
\grc_+)\big\rangle\\
& \quad \quad +
\sum_\grc\big\langle d\ulh_p, \N_1^+(Y_-; \grc_-,
\grc)\big\rangle\, \big\langle\theta _\lambda, \M_1 (\mathcal{V}); \grc,
\grc_+)\big\rangle \\
& \quad \quad -
\sum_\grc\big\langle\theta _{\bf p}, \M_1 (\mathcal{V}; \grc-,
\grc)\big\rangle\, \big\langle d\ulh_\sqcup, \N^+_1(Y_+, \grc,
\grc_+)\big\rangle\\
&\quad\quad  -\sum_\grc  \big\langle 1, \M_0 (\mathcal{V}; \grc-,
\grc_)\big\rangle\big\langle d\ulh_p\wedge d\ulh_\sqcup, \N_2^+(Y_+; \grc,
\grc_+)\big\rangle.
\end{split}
\end{equation}
This identity leads directly to the identity in the first line of
(\ref{eq:A-action-intertwine0}).

\paragraph{\it (ii) Verifying the second line in (\ref{eq:A-action-intertwine0}):}
 We proceed similarly, but now take  \(\bar{\M}=(\M)_2\) to be a 2-dimensional moduli
space of the form \(\bar{\M}_{2,z}(\mathcal{V}; \grc_-,
\grc_+)\). To compute \(\langle \delta (\scf_\Upsilon  \theta _\lambda
), \bar{\M}\rangle=\langle \scf_\Upsilon  \theta _\lambda,
(\bar{\M})_1 \rangle\), consider the bundles \(\pi _\gamma \co \hat{\mathcal{B}}^\sigma_{\gamma
 }(\mathcal{V})
\to \B^\sigma  (\V)\),  \(\hat{\td{\pi }}\co \hat{\td{\B}}^\sigma
(\V)\to \B^\sigma (V)\) defined by the following commutative diagram: 
\[
\xymatrixcolsep{5pc}\xymatrix{
\td{\tilde{\mathcal{B}}}^\sigma _{\gamma \cup\lambda }(\mathcal{V})\ar@{->}[r]^{\pi'_{\gamma }}
\ar@{->}[d]^{\pi'_{\lambda }}
\ar@{->}[rd]^{\pi_{\gamma \cup\lambda }}&\tilde{\mathcal{B}}^\sigma_{\lambda }({\cal
  V})
\ar@{->}[d]^{\pi_{\lambda }} \\
 \hat{\mathcal{B}}^\sigma_{\gamma
 }(\mathcal{V})\ar@{->}[r]^{\pi_{\gamma }} \ar@{->}[d]^{\hat{\Pi }^{\pm\infty}}  & \mathcal{B}^\sigma({\cal
   V})\ar@{->}[d]^{\Pi ^{\pm\infty}}\\ 
 \hat{\mathcal{B}}^\sigma_{\grt}(Y_\pm)\ar@{->}[r]^{\scp_\grt}  & \mathcal{B}^\sigma(Y_\pm),
}
\]
where \(\scp_\grt\co
\hat{\mathcal{B}}^\sigma_{\grt}(Y_\pm)\to \mathcal{B}^\sigma(Y_\pm)\)
was defined in Section \ref{sec:A-module}'s Part 2(b). Note that \(\hat{\mathcal{B}}^\sigma_{\gamma
 }(\mathcal{V})=(\Pi
 ^{+\infty})^*\hat{\mathcal{B}}^\sigma_{\grt}(Y_+)\simeq(\Pi
 ^{-\infty})^*\hat{\mathcal{B}}^\sigma_{\grt}(Y_-) \). Choose liftings
 \(\hat{\td{\varsigma }}\), \(\hat{\varsigma }_\gamma \), \(\td{\varsigma }_\lambda \)
 of the embedding \(\varsigma \co \bar{\M}\) that fit into the
 following commutative diagram: 
\begin{equation}\label{CD:B-hat-td}
\xymatrix{
\M^+ \ar@{->}[rr]|{\;\hat{\tilde{\varsigma }}\;}
\ar@{->}[dr]|{\;\tilde{\varsigma }_\lambda \;} \ar@{->}[drrr]|(0.25){\;\hat{\varsigma }_\gamma \;} \ar@{->}[dd]^{\grr } & &
 \td{\tilde{\mathcal{B}}}^\sigma _{\gamma \cup\lambda }(\mathcal{V})
 \ar@{->}[dd]^{\pi  _{\gamma \cup \lambda }}|(0.325)\hole \ar@{->}[ld]^(0.55){\pi'_{\gamma }}|(0.475)\hole \ar@{->}[rd]^{\pi'_\lambda }\\
& \tilde{\mathcal{B}}_\lambda ^\sigma (\mathcal{V})
\ar@{->}[rd]^{\pi_\lambda } & & \hat{\mathcal{B}}_\gamma ^\sigma (\mathcal{V})
\ar@{->}[ld]^{\pi_\gamma }\\
\bar{\M} \ar@{->}[rr]^{\varsigma } && \B^\sigma (\mathcal{V}). 
}
\end{equation}
Let \(\td{\scf}:=(\pi'_\lambda)^*\textsc{f}_\Upsilon  \). Noting that
\(\vartheta _\lambda \) is closed, we have by (\ref{eq:dF-X}) that
\[
d (\td{\scf}
\, \vartheta _\lambda )=((\pi_\lambda)^*(\Pi ^\infty)^*\mu _\gamma
)\wedge\vartheta _\lambda -((\pi_\lambda)^*(\Pi ^{-\infty})^*\mu _\gamma
)\wedge\vartheta _\lambda . 
\]
Pull back by \(\pi'_\gamma \) on both sides of the preceding
identity. Using the fact that \(\scp^*_\grt \mu _\gamma =d\op{x}_\gamma
\), a bit of diagram chasing then yields 
\[
d(\hat{\td{\scf}} \, \hat{\vartheta }_\lambda )=d\Big(\big((\pi '_\lambda )^*(\hat{\Pi }^{\infty})^*\op{x}_\gamma \hat{\vartheta }_\lambda \big)\Big)-d\Big(\big((\pi '_\lambda )^*(\hat{\Pi }^{-\infty})^*\op{x}_\gamma \hat{\vartheta }_\lambda \big)\Big),
\]
where \(\hat{\td{\scf}}:=(\pi _{\gamma \cup \lambda })^*\scf_\Upsilon
  \); \(\hat{\vartheta }_\lambda :=(\pi '_\gamma )^*\vartheta _\lambda
  \). Pull back both sides by \(\hat{\td{\varsigma }}\) and integrate
  over \(\M^+\). Then apply the Stokes' theorem (\cite{KM}'s Theorem
  24.7.2 and Lemma 31.3.1) to get: 
\[\begin{split}
& \big\langle \hat{\td{\varsigma }}^*(\hat{\td{\scf}} \, \hat{\vartheta}_\lambda ), (\M^+)_1\big\rangle\\
& \quad =\big\langle \hat{\td{\varsigma }}^*\big(\big((\pi '_\lambda
)^*(\hat{\Pi }^{\infty})^*\op{x}_\gamma  \hat{\vartheta }_\lambda
\big)\big),  (\M^+)_1\big\rangle 
- \big\langle \hat{\td{\varsigma }}^*\big(\big((\pi '_\lambda
)^*(\hat{\Pi }^{-\infty})^*\op{x}_\gamma  \hat{\vartheta }_\lambda
\big)\big),  (\M^+)_1\big\rangle .
\end{split}
\]
Recall that \(\vartheta _\lambda =\pi^*_\lambda \theta _\lambda
+(\tilde{\Pi }^\partial_\lambda )^*\vartheta '_{p_2-p_1}\). With a bit more
diagram chasing, the preceding formula
can be rewritten as: 
\begin{equation}\label{eq:m_gamma-int}
\begin{split}
& \big\langle \hat{\td{\varsigma }}^*\pi _{\gamma \cup \lambda}^*( \scf_\Upsilon\, \theta
_\lambda ),  (\M^+)_1\big\rangle +\big\langle \hat{\td{\varsigma
  }}^*\big( (\pi _{\gamma \cup \lambda}^* \scf_\Upsilon) \, (\pi '_\gamma
)^*(\tilde{\Pi }^\partial_\lambda )^*\vartheta '_{p_2-p_1}) ,
(\M^+)_1\big\rangle\\
& \quad =
\big\langle \hat{\td{\varsigma }}^*\big(\big((\pi '_\lambda
)^*(\hat{\Pi }^{\infty})^*\op{x}_\gamma \big)(\pi _{\gamma \cup \lambda
}^* \theta_\lambda )\big),  (\M^+)_1\big\rangle\\
&\quad \quad \quad -\big\langle \hat{\td{\varsigma }}^*\big(\big((\pi '_\lambda
)^*(\hat{\Pi }^{-\infty})^*\op{x}_\gamma \big)(\pi _{\gamma \cup \lambda
}^* \theta_\lambda )\big),  (\M^+)_1\big\rangle \\
&\quad \quad \quad +\big\langle \hat{\td{\varsigma }}^*\big(\big((\pi '_\lambda
)^*(\hat{\Pi }^{\infty})^*\op{x}_\gamma \big)\,  (\pi '_\gamma
)^*(\tilde{\Pi }^\partial_\lambda )^*\vartheta '_{p_2-p_1}) ,
(\M^+)_1\big\rangle\\
&\quad \quad \quad -\big\langle \hat{\td{\varsigma }}^*\Big(\big((\pi '_\lambda
)^*(\hat{\Pi }^{-\infty})^*\op{x}_\gamma \big)\,  (\pi '_\gamma
)^*(\tilde{\Pi }^\partial_\lambda )^*\vartheta '_{p_2-p_1}\Big) ,
(\M^+)_1\big\rangle. 
\end{split}
\end{equation}
 By the diagram (\ref{CD:B-hat-td}), the left most term in the preceding formula is 
\[
\big\langle \hat{\td{\varsigma }}^*\pi _{\gamma \cup \lambda}^*( \scf_\Upsilon\, \theta
_\lambda ),  (\M^+)_1\big\rangle =\big\langle \varsigma ^* (\scf_\Upsilon\, \theta
_\lambda ),  \partial[\M^+]\big\rangle =\big\langle d(\scf_\Upsilon\, \theta
_\lambda ),  \bar{\M}\big\rangle  , 
\]
namely, it is precisely the typical coefficient of \(\hat{m}[\delta  (\scf_\Upsilon\, \theta
_\lambda )](\V)\) that we seek to compute. 
To compute the other terms in (\ref{eq:m_gamma-int}),  
recall that according to \cite{KM}, \((\M^+)_1\) is a union of product
spaces of the forms:  
\begin{equation}\label{dec:M^+1}\begin{split}
& \N^+_0(Y_-\grc_-, \grc)\times \M_1 (\V; \grc_,
\grc_+), \quad \M_1 (\V; \grc_-, \grc)\times \N^+_0(Y_+; \grc, \grc_+), \\
&  \quad \quad \N^+_1(Y_-;\grc_-, \grc)\times \M_0 (\V; \grc_,
\grc_+), \quad \M_0 (\V; \grc_-, \grc)\times \N^+_1(Y_+; \grc, \grc_+).
\end{split}
\end{equation}
The map \(\grr\) is a diffeomorphism when restricted to spaces
described by the first line of the preceding expression; while
according to  (\ref{CD:B-hat-td}), spaces described by the second line
above lie in fibers of the \(U(1)\)-bundle \(\pi _{\gamma \cup \lambda
}\co \hat{\td{\B}}^\sigma  (\V)\to \B^\sigma  (\V)\).  Note that \((\pi '_\gamma
)^*(\tilde{\Pi }^\partial_\lambda )^*\vartheta '_{p_2-p_1}\) restricts
to Thom forms on fibers of \(\hat{\td{\B}}^\sigma  (\V)\), and both
\(p_1\) and \(p_2\) lie in the \(Y_+\)-end of \(\V\).  These imply
that only spaces of  the last type described in (\ref{dec:M^+1}) contribute
to the integrals in the second term on the left hand side of
(\ref{eq:m_gamma-int}), as well as to the integrals in the last two terms
on the right hand side.  Meanwhile, for the first two terms on the
right hand side of (\ref{eq:m_gamma-int}), only spaces described in
the first line of (\ref{dec:M^+1}) contribute.
Make use of these obersvations to rewrite (\ref{eq:m_gamma-int}) as
follows: 
\[\begin{split}
& \big\langle d(\scf_\Upsilon\, \theta
_\lambda ),  \bar{\M}\big\rangle  \\
& \quad =
\sum_\grc\big\langle\theta _\lambda, \M_1 (\mathcal{V}; \grc-,
\grc)\big\rangle\, \big\langle\op{u}_\gamma, \N^+_0(Y_+, \grc,
\grc_+)\big\rangle\\
& \quad \quad +
\sum_\grc\big\langle\op{u}_\gamma, \N_0^+(Y_-; \grc_-,
\grc)\big\rangle\, \big\langle\theta _\lambda, \M_1 (\mathcal{V}); \grc,
\grc_+)\big\rangle \\
&\quad\quad  -\sum_\grc  \big\langle 1, \M_0 (\mathcal{V}; \grc-, \grc_)\big\rangle\big\langle
\op{u}_\gamma \, d\ulh_\sqcup, \N_1^+(Y_+; \grc,
\grc_+)\big\rangle\\
& \quad\quad  +\sum_\grc\big\langle\textsc{f}_\Upsilon , \M_0
(\mathcal{V})(\grc-, \grc)\big\rangle \, \big\langle d \ul{\op{h}}_\sqcup, \N_1
(M_\sqcup)(\grc_, \grc_+)\big\rangle.\\
\end{split}
\]
Note that the last term in (\ref{eq:m_gamma-int}) is zero, because integrals of
the form \[\big\langle \hat{\varsigma }^*(\hat{\Pi
}^{-\infty})^*\op{x}_\gamma , \M_0 (\V; \grc_-, \grc)\big\rangle\]
vanish. Now, the preceding identity relates the coefficients in the
identity of maps in the second line of
(\ref{eq:A-action-intertwine0}), directly establishing the latter
identity.

\paragraph{Part (b):} The proof for Assertion (b) of this proposition differs
from part (a) only  in the mechanism to ensure that the right hand side of (\ref{eq:m-V}) and
its analog are well-defined. Instead of monotonicity, this is now
justified by the completeness condition on the local coefficients, and
by working with the grading-completed version of monopole Floer
complexes \(C_\bullet\). The relevant compactness theorem here is
Theorem 24.5.2 of \cite{KM}.
\epf
\begin{remarks}

(a) Recall from Section \ref{sec:2.1} that when \(c_1(\grs)\) is torsion, the following types of
perturbations are all equivalent: positive monotone, negative
monotone, balanced, exact. Thus, the assumption in part (a) implies
that \(c_1(\grs_\#)\) is nontorsion. On the other hand, the assumption that \([w_\#]\) is
monotone with respect to \(c_\#\) in part (a) implies that both
\([w_1]\), \([w_2]\) are respectively monotone with respect to
\(c_1(\grs_1)\), \(c_1(\grs_2)\)  with the same monotonicity constant. Combined
with the assumption that \([w_\#]\) is nonbalanced with respect to
\(c_1(\grs)\), this implies that \([w_i]\) is nonbalanced with
respect to \(c_1(\grs_i)\) for at least one of \(i=1\) or \(2\).  Keep
in mind that we always choose \(M_1\) to be the one endowed with a nonbalanced
perturbation.

(b) Our proof follows the ``standard'' cobordism argument that
appeared in \cite{F} and \cite{D}'s \S 7.4 in the Yang-Mills setting.
Bloom-Mrowka-Ozsvath \cite{BMO} proved a connected sum formula for the
case of exact
perturbations, using a different approach that 
involves surgery exact sequences. The use of the latter necessitates
the use of the completed version of monopole Floer homologies \(HM_\bullet\). 
 \end{remarks}

\subsection{Filtered monopole Floer homology and handle
  addition}\label{sec:filtered-conn}

Continue to work with the same settings and notation from earlier
parts of the section, but now specialize to the 3-manifolds and
cobordisms described in Sections \ref{sec:Ag)}, \ref{sec:Ah)}. More
specifically, the following two cases are considered: Fix an
\(r\gg\pi\).
\BTitem\label{eq:conditon-f}
\item[1)] Let \(M_1=Y_i\), for \(i=0,
  \ldots, \G -1\). Equip \(Y_i\) with the nontorsion \(\Spin^c\)
  structure and a metric from the set \(\Met\) in Proposition
  \ref{prop:A.6}. Let \(w_1\) be the corresponding harmonic 2-form
  \(w\) in Proposition \ref{prop:A.6}. Let \(M_2=S^1\times S^2\),
  \(\grs_2\) be the trivial \(\Spin^c\) structure, and \(w_2\equiv
  0\). Then \(M_\#\simeq Y_{i+1}\), and \([w_\#]=c_1(\grs_\#)\) is nontorsion. Choose the metric
  on \(Y_{i+1}\) to be from the set \(\Met\) from Proposition \ref{prop:A.6}.
\item[2)] Let \(M_1=S^1\times S^2\), with the nontorsion
  \(\Spin^c\)-structure \(\grs_1\), closed 2-form \(w_1\), and metric
  as described in Part 1 of Section \ref{sec:Ag)}. Let \((M_2, \grs_2)=(M, \grs)\)
  be a connected \(\Spin^c\) 3-manifold, with \(\varpi_2=rw_2\)
  for a closed 2-form \(w_2\) in the cohomlogy class
  \(c_1(\grs_2)\). Choose a metric on \(M\) with respect to which
  \(w_2\) is harmonic, and in the case when \(c_1(\grs)\) is
  non-torsion, having nondegenerate zeros. (When \(c_1(\grs)\) is
  torsion, \(w_2\) is necessarily 0). In other words, \(M_\sqcup\) is
  the \(Y_Z\) in Part 1 of Section \ref{sec:Ag)}. Thus \(M_\#\simeq
  Y_0\), and \([w_\#]=c_1(\grs_\#)\) is nontorsion. Choose the metric
  on \(Y_0\) to be from the set \(\Met\) from Proposition \ref{prop:A.6}.
 \ETitem
In both cases above, \(M_1\) is of the type \(Y_Z\) in Section
\ref{sec:Aa)}; and hence contains a special 1-cycle \(\gamma\). We
denote this by \(\gamma_1\). Consequently, assuming that \(p_1\) is
disjoint from \(\gamma_1\), both \(M_\sqcup\) and
\(M_\#\) inherit a 1-cycle from \(\gamma_1\subset M_1\).  They are
respectively denoted by \(\gamma_\sqcup\) and \(\gamma_\#\). According
to Section \ref{sec:3.8}, the filtered monopole Floer homologies
\(\HM^\circ (M_1, \langle w_1\rangle; \Lambda_{\gamma_1})\) and \(\HM^\circ(M_\#,
\langle w_\#\rangle; \Lambda_{\gamma_\#})\) are well-defined. In
parallel to what was done in Section \ref{sec:2.4}, define
\(\CM^\circ(M_\sqcup, \langle w_\sqcup\rangle;
\Lambda_{\gamma_\#})\) to be the product complex of \(\CM^\circ (M_1,
\langle w_1\rangle; \Lambda_{\gamma_1})\) and \(\hat{C}(M_2,
rw_2)\). The map \(U_\sqcup=\hat{U}_{M_\sqcup}\), as given in (\ref{eq:U-cup}),
acts on \(\CM^\infty(M_\sqcup)\) and maps \(\CM^-(M_\sqcup)\subset
\CM^\infty(M_\sqcup)\) into itself. The same notation is used to denote
its induced maps on \(\CM^+(M_\sqcup)\) and \(\CM^-(M_\sqcup)\). By
construction, the four flavors \(\CM^\circ(M_\sqcup)\) are related by
short exact sequences of the form (\ref{eq:fund-short}). Thus by Lemma
\ref{lem:u-y}, the \(H_*(S^1)\)-modules
\(S_{U_\sqcup}(\CM^\circ(M_\sqcup))\) are related by short exact
sequences of the same form. The long exact sequences induced are also
called the fundamental exact sequences. 

The remainder of this subsection consists of three parts. The first
part contains a filtered analog of Proposition \ref{prop:conn-1}. The
second part analyzes the filtered connected sum formula from Part
1. This last part derives Theorem \ref{thm:main} from this computation.

\paragraph{Part 1:} A filtered variant of Proposition \ref{prop:conn-1} states:
\begin{prop}\label{prop:conn-f}
Let \(M_\sqcup\), \(M_\#\) be as in either cases of
(\ref{eq:conditon-f}). Then there is a system of isomorphisms from \(\HM^\circ(M_\#,
\langle w_\#\rangle; \Lambda_{\gamma_\#})\) to
\(H_*(S_{U_\sqcup}(\CM^\circ(M_\sqcup))\), \(\circ=-, \infty, +,
\wedge\) as graded
\(\mathbf{A}_\dag(M_\#)\simeq\mathbf{A}_\dag(M_\sqcup)\)-modules,
which is natural with respect to the fundamental exact sequences on
both sides. 
\end{prop}
\pf Both cases in (\ref{eq:conditon-f}) satisfy the conditions of
Proposition \ref{prop:conn-1} (a). Take
\(\Gamma_1=\Lambda_{\gamma_1}\), \(\Gamma_2=\bbK\) (the constant local
coefficients). Then \(\Gamma_\sqcup=\Lambda_{\gamma_{\sqcup}}\) and
\(\Gamma_\#=\Lambda_{\gamma_\#}\). Repeat the proof of Proposition
\ref{prop:conn-1} using cobordisms \((X, w_X)\) constructed from Proposition \ref{prop:A.10a}
for Case 1) of (\ref{eq:conditon-f}), and Proposition \ref{prop:A.8}
for Case 2). Like in the previous section, we denote this by the 
shorthand \(\mathcal{V}\) when \(Y_-=M_\#\), and by \(\bar{\mathcal{V}}\) when
\(Y_-=M_\sqcup\). By construction, there is a cylinder \(C\subset X\)
ending at \(\gamma_\sqcup\subset Y_\sqcup\), and \(\gamma_\#\subset
Y_\#\) satisfying the constraints in Section \ref{sec:Ah)}. 
According to
Section \ref{sec:3.8}, this gives us chain maps \[\begin{split}
& m^\infty[u](X, \langle w_X\rangle; \Lambda_C)\co \CM^\infty(Y_-)\to
\CM^\infty(Y_+)\quad  \text{and} \\
& m^-[u](X, \langle w_X\rangle; \Lambda_C)\co \CM^-(Y_-)\to \CM^-(Y_+).
\end{split}
\] 
In parallel to (\ref{eq:V}), let 
\[\begin{split}
& V_0^\circ=m^\circ[1](\mathcal{V}; \Lambda_C), \quad V_1^\circ= m^\circ[\mathpzc{u}_\lambda](\mathcal{V}; \Lambda_C),\\
& V^{\dag, \circ}_1=m^\circ[\mathpzc{u}_{\bar{\lambda}}](\bar{\mathcal{V}}; \Lambda_C), \quad
V_0^{\dag, \circ}=m^\circ[1](\bar{\mathcal{V}}; \Lambda_C)
\end{split}\]
for \(\circ=-, \infty\), and use them to define \(V_*^\circ\),
\(V^{\dag, \circ}_*\) as in (\ref{eq:V-m}). Keeping in mind the
non-negativity of the integers \(n(\grd)\) entering the definitions of
\(\partial^\infty\) and \(m^\infty\), the rest of the proof of
Proposition \ref{prop:conn-1} may be repeated with only cosmetic
changes to see that \(V_*^\circ\) and 
\(V^{\dag, \circ}_*\) induce chain homotopy equivalences between \(\CM^\circ(M_\#,
\langle w_\#\rangle; \Lambda_{\gamma_\#})\) and 
\(S_{U_\sqcup}(\CM^\circ(M_\sqcup))\), for \(\circ=-, \infty\) These
fit into commutative diagrams with the fundamental exact sequences
(\ref{eq:fund-short}) on both sides of \(V_*^\circ\), \(V^{\dag, \circ}_*\).
This understood, the rest of the proposition follows from the Five Lemma.
\epf 

\paragraph{Part 2:} We next analyze the homologies
\(H_*(S_{U_\sqcup}(\CM^\circ(M_\sqcup))\) in the two cases of
(\ref{eq:conditon-f}) respectively. 

\paragraph{\it Case 1):}
Choose a product metric with constant
curvature on \(M_2=S^1\times S^2\). The moduli space of Seiberg-Witten
solutions over it is a circle of flat connections. Choose a real Morse
function on this circle with a pair of index 1 and index 0 critical
points, and two gradient flow lines between them. Perform a
perturbation to the Seiberg-Witten equations adapted to this Morse
function, as described in Chapter 33 of \cite{KM}. In this context
\(\hat{C}(M_2)=C^u=\bbK[u_2, y_2]\), \(\partial_{M_2}=0\), where: 
\begin{itemize}
\item the unit \(1\in \bbK[u_2, y_2]\) has grading \([\xi_+]\) in the
  notation of \cite{KM}, p.57.  
\item \(u_2\) has degree \(-2\), and the \(U_2\)-map acts by multiplication by \(u_2\). 
\item \(y_2\) has degree 1 and represents a generator of
  \(H_1(S^1;\bbZ)\) co-oriented with the moduli
  spaces. In particular, \(y_2^2=0\).
\end{itemize}
Thus, 
\begin{equation}\label{eq:conn-1}
\begin{split}
S_{U_\sqcup}(\CM^\circ  (M_\sqcup)) &=\CM^\circ (M_1)[u_2, y_2]\otimes \bbK [y], \\
D_\sqcup &=\partial_{M_1}\otimes\jmath +(U_1\otimes-1\otimes u_2)\otimes y.
\end{split}
\end{equation}
Write a generic element \(a\in S_{U_\sqcup}(\CM^\circ(M_\sqcup )\) as
\[a_0+a_1y, \quad \text{where \(a_0, a_1\in  \CM^\circ (M_1)[u_2, y_2]\).}\] Then \[D_\sqcup
a=\partial_{M_1}a_0-(\partial_{M_1}a_1)y+(U_1-u_2)(a_0)y.\]
Thus, \[\begin{split}
& H_*(S_{U_\sqcup}(\CM^\circ  (M_\sqcup))\\
& =\{\, a_0+a_1y\,|\, \partial_{M_1}a_0=0,
(U_1-u_2)\, a_0=\partial_{M_1}a_1 \}\otimes\bbK [y]\quad \mod\\
& \qquad \qquad \big(\partial_{M_1}b_1 y\sim 0,
u_2b_0 y\sim U_1b_0y-\partial_{M_1}b_0\big)\otimes \bbK[y_2]\\
&\simeq \HM^\circ  (M_1) \, y\otimes \bbK[y_2].
\end{split}\] 
Consequently,
\[H_*(S_{U_\sqcup}^\circ(M_\sqcup ))\simeq \HM^\circ( M_1)[y_2].\] 
(Alternatively, use a spectral sequence computation, filtrate
(\ref{eq:conn-1}) first by degree in \(y\), then by degree in \(u_2\)).

\paragraph{\it Case 2):} In this case, \(\grC(M_1)=\grC^o(M_1)\)
consists of a single irreducible point, \((A, (\alpha, \beta))=(0,
((2r)^{-1/2}, 0))\). (See e.g. \cite{Ds} for this well-known
fact). Thus, \(\CM^\circ (M_1)\) and the fundamental short exact
sequences relating them are simply the modules \(V^\circ\) and the
sequences in (\ref{defn:V}), (\ref{eq:V-seq1}),
(\ref{eq:V-seq2}). Write the variable \(u\) in (\ref{defn:V}) as
\(u_1\) below. As pointed out in Remark \ref{rmk:3.15} (a), \(u_1\)
stands both for the deck transformation and \(U\)-map on \(\CM^\circ (M_1)\).

This said, we have in this case
\begin{equation}\label{eq:S-f}\begin{split}
S_{U_\sqcup}(\CM^\circ (M_\sqcup ) )&=V^\circ(u_1)\otimes \CM\,( M, c_-)\otimes \bbK [y], \\
D_\sqcup &=1\otimes\partial_{M}\otimes\jmath +(u_1\otimes 1-1\otimes
U_2)\otimes y.
\end{split}
\end{equation}
This can alternatively be written as
\begin{equation}\label{eq:E}\begin{split}
E^\circ\big(\CM\,( M, c_-)\otimes \bbK [y], \partial_{M}\otimes\jmath -U_2\otimes y
\big)\\
=E^\circ\big(S_{U_2}(\CM\,(  M, c_-)\big).
\end{split}
\end{equation}
By Proposition \ref{prop:KM-ES} 
and Remark \ref{rmk:5.10}, the homology of the latter is
\(\mathring{HM}(M, c_b)\), and the isomorphisms from
\(H_*(S_{U_\sqcup}(\CM^\circ (M_\sqcup ) ))\) to the latter preserves
the \(\bbK[u]\)-module structure and are natural with respect to the
fundamental exact sequences. Since the \(U\)-map commutes with the
\(\bigwedge^*H_1(M;\bbZ)/\op{Tors}\)-actions on both sides, These are 
isomorphisms as \(\mathbf{A}_\dag(M)\)-modules.

To conclude, combining the above computation with Propostion
\ref{prop:conn-f}, we have:
\begin{cor}\label{cor:conn-f}
1) There is a system of isomorphisms of \(\mathbf{A}_\dag(M)\)-modules
\[\HM^\circ(Y_i , \langle w\rangle; \Lambda_\gamma)\simeq \HM^\circ
(Y_{i-1},  \langle w\rangle; \Lambda_\gamma)\otimes H_*(S^1)\quad \text{for
\(i=1, \ldots, \G\)}
\]
preserving the relative gradings and natural with respect to the
fundamental exact sequences.

2) There is a system of isomorphisms of \(\mathbf{A}_\dag(M)\)-modules
\[\HM^\circ(Y_0,  \langle w\rangle; \Lambda_\gamma)\simeq\mathring{\HM} (M, c_b)\]
preserving the relative gradings and natural with respect to the
fundamental exact sequences, respectively for \(\circ=-, \infty, +,
\wedge\) on the left hand side, and \(\circ=\wedge, -, \vee, \sim \)
on the right hand side. 
\end{cor}

\paragraph{\it Proof of Theorem \ref{thm:main}:}
(1): This follows from an iteration of  Corollary \ref{cor:conn-f} 1) and
Lemma \ref{lem:cor-0} below, in terms of the alternative notation (\ref{eq:nota-H}). 

\begin{lemma}\label{lem:cor-0}
There is a system of isomorphisms of \(\mathbf{A}_\dag(Y)\)-modules
\[\HM^\circ (Y, \langle w\rangle; \Lambda_\gamma)\stackrel{\simeq}{\longrightarrow}\HM^\circ(Y_{\G}, \langle w\rangle; \Lambda_\gamma)\]
preserving the relative gradings and natural with respect to the
fundamental exact sequences.
\end{lemma}
\pf \(Y\) and \(Y_{\G}\) stand for the same manifold with different
metrics and associated 2-form \(w\). As mentioned in Section
\ref{sec:2.4}, chain homotopies between the corresponding monopole
Floer complexes are provided by chain maps induced from cobordisms
\(X=\bbR\times Y\) equipped with metrics and self-dual 2-forms
interpolating those associated to \(Y_-\) and \(Y_+\). (See
e.g. Section IV.7.c for this type of argument.) In our setting, choose
\(X\) with the metrics and self-dual 2-forms over it to be those
constructed in  (\ref{prop:A.10b}). This construction also provides a
cylinder \(C\subset X\) ending at \(Y\)'s and \(Y_{\G}\)'s version of
\(\gamma\). which induces \(X\)-morphisms \(\Lambda_C\) between \(Y\)'s and \(Y_{\G}\)'s version of
\(\Gamma_\gamma\). The positivity result in Proposition \ref{prop:A.3}
guarantees that these chain maps are filtration-preserving, namely
they map \(Y_-\)'s version of \(\CM^-\subset \CM^\infty\) to
\(Y_+\)'s version of \(\CM^-\subset \CM^\infty\). As in the end of the
proof of Proposition \ref{prop:conn-f}, their induced maps
on homology together with the Five Lemma supply the isomorphisms
asserted in the lemma.
\epf

(2): This is a re-statement of Corollary \ref{cor:conn-f} 2) in
alternate notation, according to the second bullet of
(\ref{eq:nota-H}).  \epf

 \section{Properties of solutions to (\ref{eq:(A.4)})}\label{sec:B}
\setcounter{equation}{0}

This section supplies proofs for Lemma \ref{lem:A.1}
and Proposition \ref{prop:A.4}. Even so, much of what is done here is either
used in Section \ref{sec:C} or has analogs in Section \ref{sec:C}.   Section \ref{sec:Bc)} has the
proof of Lemma \ref{lem:A.1} and Section \ref{sec:Bh)} has the proof of Proposition \ref{prop:A.4}.

By way of a convention, the manifold \(Z\) is assumed implicitly to be
connected except in Section \ref{sec:Bh)}'s proof of Proposition
\ref{prop:A.4}.  

What follows is a brief outline of this section.

Section  \ref{sec:Ba)}: Lemmas \ref{lem:B.1}-\ref{lem:B.3} in this section establish
pointwise bounds on the norms of \(\psi\), \(\nabla_A\psi\), and
\(B_A\) when \((A,\psi) \) is a solution to some \((r,\mu)\)-version
of (\ref{eq:(A.4)}) in the case when \(r\) is large.

Section \ref{sec:Bb)}: Supposing that \(r\) is large, and
\((A,\psi)\) is a solution to an \((r,\mu)\)-version of (\ref{eq:(A.4)}), this section depicts on length scales that are \(O(r^{-1/2})\). This is the content of Lemma \ref{lem:B.4}.

Section  \ref{sec:Bc)}: This section introduces the notion of a
holomorphic domain. The principal examples are \(\mathcal{H}_0\) and
suitable neighborhoods of the special curve \(\gamma \) that is
described in Section \ref{sec:Aa)}3.2. Lemmas \ref{lem:B.4} and \ref{lem:B.5}  establish some very strong a priori bounds for solutions on holomorphic domains to \((r,\mu)\)-versions of (\ref{eq:(A.4)}) when \(r\) is large. This section has the proof of Lemma \ref{lem:A.1}.

Section \ref{sec:Bd)}: Lemma \ref{lem:B.7} in this section establishes
very strong a priori bounds on the 1-form \(B_A\) for a solution
\((A,\psi) \) to an \((r,\mu)\)-version of (\ref{eq:(A.4)}) where the \(w\) is harmonic.

Section \ref{sec:Be)}: Supposing that \((A,\psi)\) is a solution
to some \((r,\mu)\)-version of (\ref{eq:(A.4)}), there is a dichotomy between its
behavior where \(|\psi|\sim|w|\) and where \(|\psi|<|w|\). In the former
case, the \(\psi   \) is nearly \(A\)-covariantly constant and \(A\)
is nearly flat. This section and Lemma \ref{lem:B.8} in particular describes
\((A,\psi)\) where \(|\psi|\)  is significantly less than \(|w|\).

Section \ref{sec:Bf)}: This section gives a precise definition of the spectral flow function \(\grf_s\) (see (\ref{eq:(B.31)}) and summarizes some of its basic properties.

Section  \ref{sec:Bg)}: Lemma \ref{lem:B.9} in this section gives a priori, \(r\) and \(\grf_s\) dependent bounds for the functions \(\grc\grs\), \(\textsc{w}\) and \(\gra\) that appear in (\ref{eq:(A.5)}) and (\ref{eq:(A.6)}).

Section \ref{sec:Bh)}: This section has the proof of Proposition \ref{prop:A.4}.

\subsection{Pointwise bounds}\label{sec:Ba)}

Fix a Riemannian metric on \(Y_{Z}\) and a closed 2-form, 
denoted by \(w\), whose de Rham class is that of
\(c_{1}(\det (\bbS))\).  The four parts of this
subsection assume such data so as to supply a priori pointwise bounds
for the \(C^{\infty}(Y_{Z}; \bbS)\)-component of any given pair in \(\op{Conn} (E) \times
C^{\infty}(Y_{Z}; \bbS)\) that obeys (\ref{eq:(A.4)}).  
 
\paragraph{Part 1:}  The first lemma asserts relatively crude bounds which
are subsequently refined. 

\begin{lemma}\label{lem:B.1}  
There exists $\kappa >\pi $  with the following significance:  Fix \(r \geq
\kappa  \) and an element \(\mu \in \Omega \) with $\mathcal{P}$-norm less than 1.   Let
\((A, \psi)\)  denote a solution to the \((r, \mu)\)-version of (\ref{eq:(A.4)}).  Then \(|\psi| +
  r^{-1/2}|\nabla_{A}\psi|+r^{-1}|\nabla_{A}\nabla_{A}\psi|\leq \kappa
  \,( \sup_{Y_Z} |w|^{1/2} +r^{-1/2})\).
\end{lemma}

\pf If \(w\) is identically zero, write \(\psi = r^{-1/2} \lambda\).
 The pair \((A, \lambda)\) obeys the \(r = 1\) version of (\ref{eq:(A.4)}).  In this
case, the standard differntial equation techniques give the desired
bounds.  See for example what is said in Chapter 5 of \cite{KM}.  

 Granted that \(w\) is not identically zero, 
assume for what follows
that \(w \neq 0\) at points on \(Y_{Z}\).  The bound on
\(|\psi|\) follows by first using the Weitzenb\"ock
formula for the square of the Dirac operator to see that \(|\psi|^{2}\) obeys a differential
inequality that has the schematic form:
\begin{equation}\label{eq:(B.1)}
d^{\dag}d|\psi|^{2}+2|\nabla_{A}\psi|^{2}+ 2 r (|\psi|^{2} -|w| - c_{0} r^{-1})|\psi|^{2} \leq 0.  
\end{equation}
 
The maximum principle is now used with (\ref{eq:(B.1)}) to see that
\(|\psi|^{2} \leq c_{0} \sup_{Y_Z} |w|\) when \(r\geq c_0^{-1}\).
To say more about this, note that (\ref{eq:(B.1)}) in turn implies
that
\begin{equation}\label{:(7.2)}
d^\dag d|\psi|^2+ 2r\,( |\psi|^2-\sup_{Y_Z} |w| -c_0 r^{-1})\, |\psi|^2\leq 0.
\end{equation}
Now suppose that \(p\in Y_Z\) is a point where \(|\psi|^2\) achieves its maximum. The term \(d^\dag d|\psi|^2\)
in (\ref{:(7.2)}) is non-negative at \(p\) since it is \(-1\) times
the trace of the Hessian of \(|\psi|^2\) and the
Hessian of \(|\psi|^2\) at \(p\) is non-positive because \(p\) is a
point where \(|\psi|^2\) is maximal. It follows
as a consequence that the term \((|\psi|^2-\sup_{Y_Z} |w|
-c_0r^{-1})\) must be non-positive at \(p\), and this 
requires that \(|\psi|^2\) at \(p\) be less than \(\sup_{Y_Z} |w|+ c_0
r^{-1}\). The asserted bound follows from this. 

To see about the norm of $|\nabla_A\psi |$, digress for a moment and
fix a point \(p\in Y_Z\) and a number $\rho $ that is positive but
less than \(c_0^{-1}\). Use \(\chi \) to construct a function on \(Y_Z\) that is
equal to 1 on the ball of radius $\rho $ centered at \(p\)  and is equal to zero on the ball of radius $2\rho $
centered at \(p\). This function can and should be constructed so that the norm of its
differential is nowhere larger than $c_0\rho ^{-1}$ and so that the norm of the covariant derivative 0
of its differential is nowhere larger than \(c_0\rho ^{-2}\). Denote this function by \(\chi _\rho\). Now let \(B\)
denote the ball of radius $2\rho $ centered at \(p\). Multiply both
sides of (\ref{eq:(B.1)}) by \(\chi_\rho  \) and then integrate the
resulting inequality over \(B\). An integration by parts and an appeal
to the bounds for \(|d\chi_\rho  |\) and \(|\nabla d\chi_\rho  |\) and the bounds for \(|\psi|\) leads to the bound
\begin{equation}\label{:(7.3)}
\int _B \chi_\rho |\nabla_A\psi|^2 \leq c_0 \,( \sup_{Y_Z} |w| +r^{-1})
(\rho+\rho^3 r).
\end{equation}
To continue, let \(G_p\)  denote the Dirichlet Green's function for the
operator \(d^\dag d\) on \(B\) with pole at \(p\). This is a smooth,
non-negative function on \(B-p\) that vanishes on \(\partial B\) and obeys the pointwise bound
\begin{equation}\label{:(7.4)}
G_p(\cdot )\leq c_0 \dist \,( p, \cdot)^{-1} \quad \text{and}\quad  |dG_p|\leq c_0  \dist\, (p, \cdot)^{-2} 
 \end{equation}
at any given point \(q\in B-p\). Multiply both sides of (\ref{eq:(B.1)}) by \(\chi_\rho  G_p\) and then integrate the
resulting inequality over the ball \(B\). Use an integration by parts,
the bounds in (\ref{:(7.4)}), the a priori bounds on $d\chi_\rho $ and
$\nabla d\chi_\rho $, the a priori bound on \(|\psi|\) and (\ref{:(7.3)}) (with $\rho $ replaced by $2\rho $) to see that
\begin{equation}\label{:(7.5)}
\int _B \chi_\rho G_p |\nabla_A\psi|^2 \leq c_0 \, ( \sup_{Y_Z} |w| +r^{-1})
(1+\rho^2 r).
 \end{equation}

One last step is needed to obtain asserted the point wise bound for \(|\nabla_A\psi|^2\). To 
start this step, differentiate the equation \(D_A^2\psi =0\), commute covariant derivatives and 
then use the Bochner-Weitzenb\"ock formula again to obtain a
differential inequality for \(|\nabla_A\psi|^2\) that has the form: 
\begin{equation}\label{:(7.6)}
d^\dag d |\nabla_A\psi |^2 + 2\int_B|\nabla_A\nabla_A\psi |^2\leq c_0
\, r \,(  |\nabla_A\psi |^2 + 1) .\end{equation}
Multiply both sides of (\ref{:(7.5)}) by $\chi_\rho G_p$ and then integrate the resulting inequality over \(B\).
An integration by parts (for the left hand integral) leads to an
inequality that reads
\begin{equation}\label{:(7.7)}
|\nabla_A\psi |^2 (p)+ 2\int_B\chi_\rho G_p|\nabla_A\nabla_A\psi
 |^2\leq c_0 \rho ^3\int_B |\nabla_A\psi |^2+c_0\, r\int_B\chi_\rho G_p ( |\nabla_A\psi |^2 + 1).
\end{equation}
Granted (\ref{:(7.7)}), take $\rho =c_0^{-1}r^{-1/2}$ and then the
desired bound for \(|\nabla_A\psi |^2 (p)\) follows from
(\ref{:(7.7)}) with appeals to (\ref{:(7.3)}) and (\ref{:(7.5)}).
Much the same sort of argument using \(G_p\) and $\chi_\rho $ can be used to obtain the
asserted bounds for \(|\nabla_A\nabla_A\psi
 |^2\). Here is an outline of the argument: Multiplying the 
inequality in (\ref{:(7.6)}) by $ \chi_\rho $, integrating the result
over \(B\), then integrating by parts and using the now derived bounds
for \(|\psi|^2\) and \(|\nabla_A\psi |^2\)  leads to a \(c_0
(\sup_{Y_Z} |w| +r^{-1})(1+\rho^2 r)\) bound for 
the integral of \(\chi_\rho  |\nabla_A\nabla_A\psi|^2\). Multiplying
(\ref{:(7.6)}) by \(\chi_\rho G_p \) and integrating the result over \(B\)
leads to a \(c_0 (\sup_{Y_Z} |w| +r^{-1})(1+\rho^2 r)\) bound for the
integral of \(\chi_\rho G_p|\nabla_A\nabla_A\psi|^2\).
Meanwhile, differentiating the equation \(D_A^2\psi=0\) twice leads to
an inequality much like
(\ref{:(7.6)}) for \(d^\dag d |\nabla_A\psi |^2 \). Multplying the latter equation by $\chi_\rho G_p$ and integrating the result
over \(B\) leads to the desired bound on \(|\nabla_A\nabla_A\psi
 |^2(p)\) with the help of the previously derived bounds.
\epf

\paragraph{Part 2:}  This part of the subsection sets the notation for
what is to come in Part 3 and in the subsequent sections.  To start,
introduce \(K_{*}^{-1}\) to denote the
2-plane subbundle of the tangent bundle over the \(|w|
> 0\) part of \(Y_{Z}\) given by the kernel of \(*w\).  Orient \(K_{*}^{-1}\) by the
restriction of \(w\) and use the induced metric with this orientation to
view \(K_{*}^{-1}\) as a complex line
bundle.  Clifford multiplication by the 1-form \(*w\) on the
\(|w| > 0\) part of \(Y_{Z}\) writes \(\bbS\) as a direct sum of eigenbundles
\(E_{*} \oplus (E_{*}\otimes K_{*}^{-1})\) with \(E_{*}\) being the \(+i |w|\) eigenbundle.    
 
Use $\mathbb{I}_{\bbC}$ to denote the product
complex line bundle and \(\theta_{0}\) to denote the product
connection on $\mathbb{I}_{\bbC}$.  Let \(1_{\bbC}\) denote the
\(\theta_{0}\) -constant section of $\mathbb{I}_{\bbC}$ with value 1 at all points.
 Fix a unitary identification between \(E_{*}^{-1}\otimes_{\bbC} E_{*}\) and $\mathbb{I}_{\bbC}$ and
use the latter to write \(E_{*}^{-1}\otimes_{\bbC}\bbS \) as
$\mathbb{I}_{\bbC} \oplus K_{*}^{-1}$.  The bundle
\(K_{*}^{-1}\) has a canonical connection,
which we denote by \(A_{K_*}\), such that the section \((1_{\bbC}, 0)\) of the
bundle $\mathbb{I}_{\bbC} \oplus K_{*}^{-1}$ obeys the Dirac equation
as defined using the connection \(A_{K_*}+ 2\theta_{0}\) on its determinant line bundle.  The
norm of the curvature of \(A_{K_*}\) is bounded by \(c_{0}|w|^{-2}\) and the norm of the
\(k\)-th derivative of \(A_{K_*}\)'s curvature is bounded by \(c_{k}
|w|^{-2-k}\) with \(c_{k}\) being a constant.   

A section \(\psi\) of \(\bbS\) over \(U\) is written with respect to
this splitting as \(|w|^{1/2}(\alpha, \beta)\).  Meanwhile, the connection \(A\) on \(E\) defines a
corresponding connection on \(E_{*}\), that is, the connection
\(A_{*} = A -\frac{1}{2} (A_{K} -A_{K_*})\).  To keep the notation under control in what follows, the
\(A_{*}\)-covariant derivative on \(E_{*}\) is
also denoted by \(\nabla_{A}\), as is the
\(A_{*} + A_{K_*}\)-covariant derivative on \(E_{*} \otimes K_{*}^{-1}\).

 \paragraph{Part 3:} The next lemma refines Lemma \ref{lem:B.1}'s
bound on the \(|w| > 0\) part of \(Y_{Z}\).

\begin{lemma}\label{lem:B.2}   
There exists $\kappa >\pi $  with the following significance:  Fix \(r \geq\kappa  \) and an element \(\mu\in \Omega \) with $\mathcal{P}$-norm less than 1.  Let
\((A, \psi)\)  denote a solution to the \((r, \mu)\)-version of (\ref{eq:(A.4)}).  Fix \(m \in (\kappa,
\kappa r^{1/3} (\ln r)^{-\kappa})\)  and let \(U_{m}\)  denote the
\(|w| > m^{-1}\)  part of \(Y_{Z}\).  Write \(\psi\)  on \(U_{m}\)  as
\(|w|^{1/2} (\alpha,\beta)\).   Then the pair \((\alpha, \beta)\)
obeys the following on \(U_m\):
\begin{itemize}
\item  \(|\alpha|^{2} \leq 1 +\kappa m^{3}r^{-1}\).
\item  \(|\beta|^{2} \leq \kappa m^{3} r^{-1} (1 -|\alpha|^{2}) +\kappa^{3} m^{6}r^{-2}\).
\item  \(|\nabla_{A}\alpha|^{2}+ m^{-3}r |\nabla_{A}\beta|^{2}\leq \kappa m^{-1}r (1 -|\alpha|^{2}) +\kappa^{2}m^{2}\). 
\item  Denote by \(U_{*}\) the \(1 -|\alpha|^{2} \geq\kappa^{-1}\)
  part of \(U_{m}\).  Then \[|1-|\alpha|^{2} | \leq
(m^{2}e^{-\sqrt{r/m}\op{dist} (\cdot, U_*)/\kappa} +  \kappa m^{3}r^{-1}).\]
\end{itemize}
\end{lemma}
 
\pf Write \(\psi =|w|^{1/2} \eta \) on \(U_{2m}\).  The section
\(\eta\) on \(U_{2m}\) obeys an equation having the schematic form \(D_{A}\eta
+ \grR\cdot\eta  = 0\) with \(\grR\) being Clifford multiplication by the 1-form
\(\frac{1}{2}d(\ln  |w|)\).   Note in particular that \(|\grR| \leq
c_{0} \, m\) and the absolute
value of the covariant derivative of \(\grR\) is bounded by \(c_{0}\,
m^{2}\).  Use the Weitzenb\"ock formula for the operator
\(D_{A} + \grR\) to see that \(\eta\) obeys an equation that has the schematic form
\begin{equation}\label{eq:(B.2)}
\nabla_{A}^{\dag}\nabla_{A}\eta- \cl(B_{A})\cdot\eta +\grR_{1}\cdot\nabla_{A}\eta+ \grR_{0}\cdot\eta = 0,
\end{equation}
where \(\cl(\cdot)\) denotes the Clifford multiplication
endomorphism from \(T^*M\) to \(\End (\bbS)\) and where \(\grR_{1}\)
and \(\grR_{0}\) are linear and obey \(|\grR_{1}| \leq c_{0} m\)
and \(|\grR_{0}| \leq c_{0}m^{2}\).  Let \(\q\) denote the maximum of \(0\) and
\(|\eta|^{2} - 1 -c_{0}m^{3} r^{-1}\).  It
follows from (\ref{eq:(B.2)}) that \(\q\) on \(U_{2m}\) obeys
\begin{equation}\label{eq:(B.3)}
d^{\dag}d \q + 2r m^{-1} \q \leq 0.
\end{equation}

As Lemma \ref{lem:B.1} bounds \(\q\) by \(c_{0}m\) on  \(U_{2m}\), the comparison principle with the Green's function
for the operator \(d^{\dag}d + r m^{-1}\)
to see that \(\q \leq c_{0} m^{3}r^{-1}\) on \(U_{3m/2}\).  This implies the claim in the first bullet.  It also implies that
\(|\beta|^{2}\) is less than \(1 +c_{0} m^{3} r^{-1}\) on  \(U_{3m/2}\).  

 To see about the second bullet, project (\ref{eq:(B.2)}) onto the
\(E_{*} \otimes K_{*}^{-1}\) summand
of \(\bbS\) and take the fiberwise inner product of the resulting
equation with \(\beta\) to obtain a differential inequality that has
the form
\begin{equation}\label{eq:(B.4)}
d^{\dag}d|\beta|^{2} +  2r m^{-1}|\beta|^{2} \leq-|\nabla_A\beta|^{2} +
c_{0} r^{-1} m^{3}|\nabla_{A}\alpha|^{2}+ c_{0} r^{-1} m^{5}.
\end{equation}

Fix for the moment $\varepsilon$ > 0.  Project (\ref{eq:(B.2)})
next onto the \(E_{*}\) summand and take the pointwise inner
product with \(\alpha\) to obtain an equation for the function \(\w = 1- |\alpha|^{2}\) that has the form
\begin{equation}\label{eq:(B.5)}
d^{\dag}d \w + 2r m^{-1} \w =
2|\nabla_A\alpha|^{2 }+ r m^{-1} \w^{2} + \gre   , 
\end{equation}
where \(|\gre| \leq c_{0}\varepsilon|\nabla_A\beta|^{2} +c_{0}(1 + \varepsilon^{-1})m^{2}  + c_{0}m |\nabla_A\alpha|\).
 
It follows from (\ref{eq:(B.4)}) and (\ref{eq:(B.5)}) that there exist constants
\(z_{1}\) and \(z_{2}\) that are both bounded by
\(c_{0}\) and \(\varepsilon >c_{0}^{-1}\) such that the function \(q =
|\beta|^{2} - z_{1}r^{-1}m^{3} \w  - z_{2}r^{-2} m^{6}\) obeys the equation
\begin{equation}\label{eq:(B.6)}
d^{\dag}dq + 2r m^{-1} q \leq 0
\end{equation}
on \(U_{3m/2}\).  Granted this inequality, use the Green's function
for \(d^{\dag}d + r m^{-1}\) as before to see that \(|\beta|^{2}
\leq z_{1} m^{3}r^{-1} (1  -|\alpha|^{2}) + z_{2}m^{6}  r^{-2}\) on \(U_{m}\).
 
The proofs of the third and fourth bullets start by differentiating
(\ref{eq:(B.2)}) to obtain an equation for the components of
\(\nabla_{A}\eta\) and it then copies the
manipulations done in Step 2 of Section 4d in \cite{Ta1}
to obtain a differential inequality on \(U_{3m/2}\) for the function \(\grh :=
|\nabla_{A}\eta|^{2}\) that has the form
\begin{equation}\label{eq:(B.7)}
 d^{\dag}d\grh + 2r m^{-1} \grh \leq c_{0} (r m^{-1} \w \grh + m^{2} \grh + m^{4} +  r^{2}m^{-2}  \w^{2}).
\end{equation}
 
To prove the third bullet, use (\ref{eq:(B.4)}), (\ref{eq:(B.5)}) and (\ref{eq:(B.7)}) to find constants
\(z_{1}\), \(z_{2} > 0\) and \(z_{3}\), all with absolute value less
than \(c_{0}\),
such that the function \(q := \grh - z_{1} r  m^{-1} (1- |\alpha|^{2}) -z_{2} m^{2} + z_{3} r
m^{-1} |\beta|^{2}\) obeys (\ref{eq:(B.6)}) on \(U_{2m}\) when \(m < c_{0}^{-1}
r^{1/3}\).   Meanwhile, Lemma \ref{lem:B.1} implies that \(\grh\) is no
larger than \(c_{0} m r \) on \(U_{2m}\).  Given this last bound, the comparison argument that uses the
Green's function for \(d^{\dag}d + c_{0}^{-1}r m^{-1}\) says that
\(|\nabla_A\eta|^{2}\) is bounded by \(c_{0} m^{-1}r (1 -|\alpha|^{2}) +c_{0}^{2} m^{2}\) on
\(U_{m}\)  when \(m \leq c_{0}^{-1}r^{1/3}\).  This gives Lemma \ref{lem:B.2}'s bound for
\(|\nabla_{A}\alpha|^{2}\).
 The refinement that gives the asserted bound for \(|\nabla_A\beta|^{2}\) is obtained by
the same sort of argument after first projecting (\ref{eq:(B.2)}) onto the
\(E_{*}\otimes K_{*}^{-1}\)-summand of \(\bbS\) before differentiating so as to get an elliptic
equation for \(\nabla_{A}\beta\).  The details of
this part of the story are straightforward and omitted.

To prove the fourth bullet, use the first bullet of the lemma with (\ref{eq:(B.5)})
and (\ref{eq:(B.7)}) to see that \(q: = \grh + c_{0}^{-1}r m^{-1} \w - c_{0} m^{2}\)
obeys an equation on the \(\w \leq c_{0}^{-1 }\) part of  \(U_{2m}\)
that has the form \(d^{\dag}dq +c_{0}^{-1} r m^{-1} q\leq 0\) when \(m
\leq c_{0}^{-1}r^{1/3}\).
 Granted the latter and granted the a priori bound \(q \leq
c_{0} r m\) from Lemma \ref{lem:B.1}, then the comparison principle
using the Green's function for
\(d^{\dag}d + c_{0}^{-1} r m^{-1}\)  leads to the following:  If \(\c >
c_{0}\), then \(q \leq c_{0} r
m \, e^{-\sqrt{r/m}\dist (\cdot, U_\c)/c_0}\)
  where \(U_{\c}\) denotes the \(\w \geq \c^{-1}\) part of \(U_{2m}\).
  This last inequality implies Lemma \ref{lem:B.2}'s fourth bullet. \epf
 
\paragraph{Part 4:}  The final lemma of this subsection refines what
is said by Lemma \ref{lem:B.1} on the part of \(Y_{Z}\) where
\(|w|\) is positive but small.

\begin{lemma}\label{lem:B.3}   
There exists \(\kappa >1\) with the following property:  Fix \(m \in (\kappa,
\kappa^{-1} r^{1/3}  (\ln r)^{-\kappa })\).   Fix \(r \geq \kappa  \)
and fix \(\mu \in \Omega  \) with $\mathcal{P}$-norm less than
\(1\)  and let \((A, \psi)\) be a solution to the \((r,\mu)\)-version
of (\ref{eq:(A.4)}). Then 
\(|\psi|\leq \kappa m^{-1/2}\)  and \(|\nabla_{A}\psi| \leq\kappa
m^{-1} r^{1/2}\) on the \(|w| < m\)  part of \(Y_{Z}\).
\end{lemma}

\pf The maximum principle applied to (\ref{eq:(B.1)}) implies that \(|\psi|^{2}\)
can not have a local maximum where \(|\psi|^{2} >|w| + c_{0} r^{-1}\).
Indeed, if \(p\in Y_Z\) is a point where this
condition holds, then the left hand side of (\ref{eq:(B.1)}) at \(p\)
is strictly greater than \(d^\dag d|\psi|^2\) at \(p\). If
\(p\) is a local maximum of \(|\psi|^2\), then \(d^\dag d|\psi|^2\geq 0\)
at \(p\) and so the left hand side of (\ref{eq:(B.1)})  would be
positive which violates (\ref{eq:(B.1)}).

Since \(|\psi|^2\) can not have a local maximum where
\(|\psi|^2>|w|+c_0 r^{-1}\), it follows that \(|\psi|^2\)  can not
have a local maximum where \(|\psi|^2> m^{-1}+c_0 \, r^{-1}\) on the set
where \(|w|<  m^{-1}\). 
Meanwhile, Lemma \ref{lem:B.2} implies that \(|\psi|^{2}\leq c_{0}\, 
m^{-1}\) on the boundary of the set where  \(|w|<  m^{-1}\) (which is the boundary of 
\(U_{m}\)). Therefore, \(|\psi|^2\)   can not be greater than the
maximum of \(c_{0} \, m^{-1}\) and \(m^{-1}+c_0\,  r^{-1}\)  on the set
where \(|w|<  m^{-1}\). If \(m\in (c_0,
c_0^{-1} r^{1/3}  (\ln r)^{-c_0})\), 
then this maximum is \(c_{0} m^{-1}\).


To see about \(|\nabla_{A}\psi|\), let \(p
\in Y_{Z}\) denote a given point where \(w \leq 2m^{-1}\).  Fix Gaussian coordinates for a ball of
radius \(c_{0}^{-1}\) centered at \(p\) and then
rescale the coordinates so that the ball of radius
\(m^{-1/2} r^{1/2}\) about the origin in
$\mathbb{R}^{3}$ and radius 1.  Let \(\varphi\) denote
the corresponding map from the ball of radius \(1\) about the origin in
$\mathbb{R}^{3}$ to the original ball in
\(Y_{Z}\).  With this understood, the pull-back
\((\varphi^*A, m^{1/2}\varphi^*\psi)\) satisfies a version of (\ref{eq:(A.4)}) on the unit ball in \(\bbR^{3}\)
that is defined by the rescaled metric. It follows from the bound on
\(|\psi|\) that
\(|B_{A}| \leq c_{0}m^{-1} r\) and this implies that \(|\varphi^*B_{A}| \leq
c_{0}\).  This understood, standard elliptic regularity
techniques can be employed to see that the rescaled version of
\(m^{1/2}|\varphi^*(\nabla_{A}\psi)|\) has norm bounded by \(c_{0}\)
and so \(|\nabla_{A}\psi|\) has norm bounded by \(c_{0} m^{-1}r^{1/2}\).  \epf

\subsection{The micro-local structure of \((A,\psi)\)}\label{sec:Bb)}

 Part 3 of this section states and then proves Lemma \ref{lem:B.4}, this being
a lemma that describes solutions to (\ref{eq:(A.4)}) on the
\(|w| > 0\) part of \(Y_{Z}\)
when viewed with microscope that magnifies by a factor of the order of
\(r^{1/2}\).   Parts 1-2 of the subsection set the
notation that is used in particular for Lemma \ref{lem:B.4} but elsewhere as
well.
 
\paragraph{Part 1:}  This part of the subsection introduces the
 {\em vortex equations} on \(\bbC\).  This is a system of
equations that asks that a pair \((A_{0},\alpha_{0})\) of connection on a complex line bundle
over \(\bbC\) and section of this bundle obey
\begin{equation}\label{eq:(B.8)}
\begin{cases}
* F_{A_0}= -i (1 -|\alpha_{0}|^{2}), &\\
\bar{\partial}_{A_0}\bar{\alpha}_0= 0,&\\
|\alpha_{0}| \leq 1 .&
\end{cases}
\end{equation}
 
The notation here is such that \(*\) denotes the Euclidean Hodge dual on
\(\bbC\), while \(F_{A_0}\) and \(\bar{\partial}_{A_0}\) denote the respective curvature 2-form of \(A_{0}\) and the
d-bar operator defined by \(A_{0}\) on the space of sections
of the given complex line bundle.  Note that if \((A_{0},
\alpha_{0})\) is a solution to (\ref{eq:(B.8)}), then so is
\((A_{0}  - u^{-1}du, u\, \alpha_{0})\) with \(u\) being any smooth
map from \(\bbC\) to \(S^{1}\).

 Solutions with \(1 -|\alpha_{0}|^{2}\)
integrable are discussed at length in Sections 1 and 2 of \cite{T2},
Section IV.2b and Section IV.3a.  As noted in these references, if \(1 -
|\alpha_{0}|^{2}\) is
integrable then its integral is $2\pi $ times a non-negative integer.
 Fix \(\m \in \{0, 1, \ldots\}\).  The space of
\(C^{\infty}(\bbC; S^{1})\)
equivalence classes of solutions to (\ref{eq:(B.8)}) with the integral of \(1  -|\alpha_{0}|^{2}\)
equal to $2\pi \m$  has the structure of a smooth, \(2\m\)-dimensional
manifold.  This manifold is denoted in what follows by
\(\grC_{\m}\).  By way of a parenthetical remark, the space
\(\grC_{\m}\) has a natural complex structure that identifies it
with \(\bbC^{\m}\).  A solution with \(1 -|\alpha_{0}|^{2}\)
integrable is said here to be a  {\em finite energy solution} to the vortex equation. 

\paragraph{Part 2:}  Lemma \ref{lem:B.4} and some of the later subsections refer to
the notion of a transverse disk with a given radius through a
given \(|w| > 0\) point in
\(Y_{Z}\).  A {\em transverse disk} is the image via the
metric's exponential map of the centered disk of the
given radius in the 2-plane bundle \(\ker (*w)\) at the given point.
 There exists \(c_{0} > 100\) such that any
transverse disk with radius \(c_{0}^{-1}\) is
embedded with a priori bounds on the derivatives to any given order of
its extrinsic curvature.  If \(D \subset Y_{Z}\) is a
transverse disk centered at a point \(p\), and if \(\c \geq
c_{0}\), then \(|w|\) will be greater than  
\(\frac{1}{2}|w|(p)\) on the subdisk in \(D\) centered at \(p\) with radius
\(\c^{-1} |w|(p)\).  The constant \(\c\) can
be chosen so that the following is also true:  Let \(v\) denote the vector
field on the \(|w| > 0\) part of
\(Y_{Z}\) that generates the kernel of \(w\) and has pairing \(1\)
with \(*w\).  Then \(v\) is orthogonal to \(D\) at \(p\) and the length of the
projection to \(TD\) of \(v\) on the concentric disk in \(D\) of radius
\(\c^{-1}|w|(p)\) is no greater than
\(c_{0}\c^{-1}\).  Choose \(\c \geq
c_{0}\) with this property and use \(D_{p}\) to
denote the transverse disk through \(p\) of radius \(\c^{-1}
|w|(p)\).   
 
Reintroduce from Part 2 of Section \ref{sec:Ba)} the complex line bundle
\(K_{*}^{-1}\) defined over the \(|w| > 0\) part of \(Y_{Z}\).
 Recall that the underlying real bundle is the 2-plane bundle in
\(TY_{Z}\) annihilated by \(*w\).  Let \(p\) again denote a point
in the \(|w| > 0\) part of \(Y_{Z}\).  Fix an isometric isomorphism from
\(K_{*}^{-1}|_{p}\) to \(\bbC\).  Use \(\varphi\) in what follows to denote the map
from \(\bbC\) to \(Y_{Z}\) that is obtained by composing first the isomorphism with
\(K_{*}|_{p} =\ker (*w)|_{p}\) and then the metric's exponential map.  With \(r \geq 1\) given,
use \(\varphi_{r}\) to denote the composition of first
multiplication by \(r^{-1/2}|w(p)|^{-1/2}\) on \(\bbC\) and then
applying \(\varphi\).  

 To finish the notational preliminaries, suppose that \((A, \psi)\) is a
given pair in \(\op{Conn} (E) \times C^{\infty}(Y_{Z}; \bbS)\).
 Write \(\psi\) where \(|w| > 0\) as \(|w|^{1/2}(\alpha, \beta)\) to conform with Part 2 of Section
\ref{sec:Ba)}'s splitting of \(\bbS\) as \(E_{*}\oplus (E_{*}\otimes K_{*}^{-1})\).  Likewise reintroduce from
Part 2 of Section \ref{sec:Ba)} the connection \(A_{*}\) on the bundle
\(E_{*}\).  Given \(p \in Y_{Z}\) with \(|w(p)| > 0\), introduce \((A_{r}, \alpha_{r})\) to denote the
\(\varphi_{r}\)-pull back of the pair \((A_{*}, \alpha)\) to the
radius \(c^{-1}r^{1/2}|w(p)|^{1/2 }\) disk in \(\bbC\).

\paragraph{Part 3:}  Lemma \ref{lem:B.4} below characterizes the pair
\((A_{r}, \psi_{r})\).
 
\begin{lemma}\label{lem:B.4}  
There exists \(\kappa >10\)  and given \(\textsc{r} > \kappa \),
there exists \(\kappa_{\textsc{r}} >1\)  with the following property:  Fix \(r \geq
\kappa_{\textsc{r}}\)  and \(\mu \in\Omega  \) with $\mathcal{P}$-norm bounded by
\(1\) .  Suppose that \((A, \psi)\)  is a solution to the
\((r,  \mu)\)-version of (\ref{eq:(A.4)}).  Fix a point in
\(Y_{Z}\) where \(|w| >r^{-1/3}(\ln r)^\kappa \) and use the corresponding version of
\(\varphi_{r}\) to obtain the pair \((A_{r}, \alpha_{r})\) of
connection and section of a complex line bundle over
\(\bbC\).  There exists a solution to the vortex equation
on \(\bbC\) whose restriction to the radius \(\textsc{r}\) disk about the origin in
\(\bbC\) has \(C^{1}\)-distance less than \(\textsc{r}^{-4}\) from \((A_{r}, \alpha_{r})\) on this same
disk. Moreover, if \(1-|\alpha_{r}|^{2}< \frac{1}{2}\) at distances
between \(\textsc{r} + \kappa (\ln \textsc{r})^{2}\)  and  \(\textsc{r} -\kappa (\ln  \textsc{r})^{2}\) from the
origin, then \((A_{r}, \alpha_{r})\) has \(C^{1}\)-distance less than
\(\textsc{r}^{-4}\) in the radius \(\textsc{r}\)-disk about the origin in \(\bbC\) from a finite
energy solution to the vortex equations that defines
a point in some $\m \leq \pi \textsc{r}^{2}$ version of \(\grC_{\m}\).
\end{lemma}

\pf It follows from (\ref{eq:(A.4)}) and what
is said by the first three bullets of Lemma \ref{lem:B.2} that the curvature of
\(A_{r}\) and \(\alpha_{r}\) are such that
\begin{equation}\label{eq:(B.9)}
* F_{A_r} = - i (1-|\alpha_{r}|^{2}) +\gre_{0}  \quad\text{and}\quad 
 \bar{\partial}_{A_r}\alpha_r= \gre_{1}, 
\end{equation}
where \(|\gre_{0}| +|\gre_{1}| \leq c_{0} (\ln r)^{-c_0}\) on the disk in \(\bbC\) of radius less than \(\c^{-1}
r^{1/2}m^{-1/2}\).  The third bullet in
Lemma \ref{lem:B.2} also finds \(|\nabla  _{A_r}\alpha_r | \leq c_{0}\).  Granted (\ref{eq:(B.9)}), then the
argument used to prove Lemma 6.1 in \cite{TW1} can be used with only minor
modifications to prove the assertion with \(C^{1}\)-distance
replaced by the distance as measured by any \(\upsilon < 1-
\textsc{r}^{-1}\) H\"older norm.  The convergence in the \(C^{1}\)-topology follows using the arguments from
Section 6 in \cite{TW1} given also the second derivative bound from Lemma
\ref{lem:B.1}.   \epf

\subsection{Holomorphic domains}\label{sec:Bc)}

What follows directly sets the notation for what is to come in this
subsection.  An open set \(U \subset Y_{Z}\) is said to
be a {\em holomorphic domain} when the following criteria are met:
\begin{itemize}
\item  The metric has non-negative Ricci curvature on \(U\).
\item  The 2-form \(w\)  is non-zero on \(U\)  and covariantly constant.
\item  The curvature of \(A_{K}\)  on \(U\)  is a multiple of \(w\).
\item  The 1-form \(\mu\) on \(U\)  and its derivatives
to order 10 have norm less than \(e^{-r^2/2}\).
\end{itemize}
 
The following lemma strengthens the conclusions of Lemma \ref{lem:B.2} on a
holomorphic domain.

 \begin{lemma}\label{lem:B.5} 
Let \(U \subset Y_{Z}\) denote a holomorphic domain and let
\(U_{1} \subset U\)  denote an open set with
compact closure in \(U\) .  Use \(\textsc{d}\)  to denote the
function on \(U\)  that measures the distance to
\(Y_{Z}-U\) .  There exists $\kappa> \pi $  with the following significance:  Fix \(r
\geq \kappa  \) and a 1-form \(\mu \in\Omega  \) with $\mathcal{P}$-norm less than 1 whose
norm on \(U\)  and those of its first 10 derivatives is bounded by
\(e^{-r^2/2}\). Suppose that \((A, \psi)\)  is a solution to the \((r,  \mu)\)-version of (\ref{eq:(A.4)}).  Write
\(\psi\)  on \(U\)  as \(|w|^{1/2} (\alpha, \beta)\).  Then \(\beta\)
on \(U_{1}\) obeys:
\begin{itemize}
\item  \(|\beta| \leq \kappa e^{-\sqrt{r}\textsc{d}/\kappa}\).
\item  Given \(q \geq 1\), there exists \(\kappa_{q} \geq 1\)  such that
 \(|(\nabla_{A})^{q}\beta|\leq \kappa_{q} e^{-\sqrt{r}\textsc{d}/\kappa}\)
   with \(\kappa_{q}\)  depending only on the metric, \(A_{K}\) ,
\(U\)  and \(U_{1}\).
\end{itemize}
 \end{lemma}

\pf The proof that follows assumes
that \(\mu = 0\) on \(U\).  The proof in the general case differs little
from what is said below and is left to the reader.  

Keep in mind that the norm of \(|w|\) is constant on \(U\)
because \(w\) is covariantly constant.  Project the Weitzenb\"ock formula
for \(D_{A}^{2}\)  onto the \(E_{*}\otimes K_{*}^{-1}\) summand
of \(\bbS\) to obtain an equation for \(\beta\) on \(U\) that has the
schematic form:
\begin{equation}\label{eq:(B.10)}
\nabla_{A}^{\dag}\nabla_{A}\beta+  r |w| (1 +|\alpha|^{2} +|\beta|^{2}) \beta + \grR\beta = 0,
\end{equation}
with \(\grR\) determined solely by the metric and \(A_{K}\).  Granted
this, then by the conditions on the metric and \(A_K\) over \(U\), \(|\beta|\) obeys an equation of the form
\(d^{\dag}d |\beta| + r|w| |\beta|\leq 0\) on \(U\) when \(r\) is
larger than a constant that depends only on \(U\) and \(U_{1}\).  The bound in the first bullet of the lemma
follows from the latter equation using the comparison principle and the
Green's function for the operator
\(d^{\dag}d + r|w|\).  Given the bounds from Lemma \ref{lem:B.2}, very much the same strategy leads to the bounds in the subsequent
bullets after differentiating (\ref{eq:(B.1)}) to obtain an equation for
\((\nabla_{A})^{q} \beta\). \epf

Lemma \ref{lem:B.5} leads directly to the next lemma that describes \(\psi\)
on \(U_{\gamma}\) and $\mathcal{H}_0$.
 
\begin{lemma}\label{lem:B.6}
Given \(\varepsilon >0\), there exists $\kappa \geq \pi $  with the
following significance:  Introduce \(U\)  to denote
$U_{\gamma} \cup\mathcal{H}_0$ and let
\(\textsc{d}\)  denote the function on \(U\)  that measures the
distance to \(Y_{Z}-U\).   Introduce \(U_{\varepsilon} \subset U\)  to denote the
subset with \(\textsc{d} >  \varepsilon\).   Fix
\(r \geq \kappa  \) and a 1-form \(\mu \in \Omega  \) with $\mathcal{P}$-norm less than 1 whose
norm on \(U\)  and those of its first ten derivatives is bounded by \(e^{-r^2/2}\).  Let \((A, \psi)\) denote a solution to the \((r, \mu)\)-version
  of (\ref{eq:(A.4)}).  The following is true on \(U_{\varepsilon}\):
\begin{itemize}
\item The conclusions of Lemma \ref{lem:B.5} hold with \(U_{1}\)
  therein set to \(U_{\varepsilon}\).
\item  \( -\kappa e^{-\sqrt{r}\textsc{d}/\kappa}\leq 1  - |\alpha|^2  \leq \kappa e^{-\sqrt{r}\textsc{d}/\kappa}\).
\item Given \(q \geq 1\), there exists \(\kappa_{q} \geq 1\) such that
 \(|(\nabla_{A})^{q}\alpha| \leq \kappa_{q}  e^{-\sqrt{r}\textsc{d}/\kappa}\)
with \(\kappa_{q}\)  depending only on the metric, \(A_{K}\),
\(U\)  and \(\varepsilon\) .
\end{itemize}
 \end{lemma}

\pf The first bullet follows by virtue of the fact that 
$U_{\gamma} \cup\mathcal{H}_0$ is a holomorphic domain where the
constraints in (\ref{eq:(A.3a)}) and (\ref{eq:(A.3b)}) are obeyed.  To see about the other bullets of
the lemma, suppose for the moment that \(\delta > 0\),
that $p \in\mathcal{H}_0\cap U_{\varepsilon}$  and that \(1-|\alpha| >
\delta \) at \(p\).  As is proved in what follows, this assumption
leads to nonsense unless \(\delta\) is very small. So, supposing that \(1-|\alpha| >
\delta \) at \(p\), it
follows from the second bullet of Lemma \ref{lem:B.2} and from Lemma \ref{lem:B.5} that the
integral of \(*B_{A}\) on the radius
\(c_{0}^{-1} r^{-1/2}\delta \) disk in the constant \(u\) slice of
$\mathcal{H}_0$ through \(p\) is greater than \(c_{0}^{-1}\delta^{3}\).
 Lemma \ref{lem:B.5} implies that the pull-back of \(*B_{A}\) to the
constant \(u\) sphere through \(p\) can be written as   
 \(\frac{i}{4\pi}\textsc{f} \sin\theta d\theta \wedge d\phi \) and
 that \(\textsc{f} \geq -c_{0}  e^{-\sqrt{r}/c_0}\). This implies that 
 the integral of \(*B_{A}\) on this transverse sphere
in $\mathcal{H}_0$ will be positive if \(\delta
> c_{0}  e^{-\sqrt{r}/c_0}\).   
But the integral of \(*B_A\) on this transverse sphere is zero because
\(E\)'s first Chern class has zero pairing with the
\(H_2(\mathcal{H}_0; \bbZ)\)-summand in (\ref{eq:(A.2)}).  Therefore, it
must be the case that \(1-|\alpha|< c_{0}  e^{-\sqrt{r}/c_0}\) on \(\mathcal{H}_0\cap U_{\varepsilon}\).
 Now suppose that \(p \in U_{\gamma} \cap U_{\varepsilon}\).  The Dirac equation writes the  
\(\frac{\partial}{\partial t}\)-covariant derivative of \(\alpha\) as a linear combination of
covariant derivatives of \(\beta\).  This understood, Lemma \ref{lem:B.5}
implies that the absolute value of the
\(\frac{\partial}{\partial t}\)-covariant derivative of \(\alpha\) in \(U_{\gamma}\) is
bounded by \(c_{0}e^{-\sqrt{r}/c_0}\).  It follows as a consequence that if \(1-|\alpha|
> \delta \) at a point in \(U_{\gamma}\cap U_{\varepsilon}\), then
\(|\alpha| > \frac{1}{2}\delta \) at points in $\mathcal{H}_0\cap
U_{\varepsilon}$ if \(\delta > c_{0} e^{-\sqrt{r}/c_0}\), and as
explained previously, this is not
allowed if \(r \geq c_{0}\). 
Therefore, the conclusion is that
\(1-|\alpha|<c_{0}e^{-\sqrt{r}/c_0}\)  on the whole of
\((U_{\gamma}\cup  \mathcal{H}_0)\cup U_{\varepsilon}\). Much the same sort of argument proves that \(1-|\alpha|>-c_{0}e^{-\sqrt{r}/c_0}\)  on this same domain.


The assertion in the third bullet is proved by writing \(\psi =|w|^{1/2} \eta\)
on \(U\).  Keeping in mind that \(|w|\) is constant on \(U\),
project the Weitzenb\"ock formula for \(D_{A}^{2}\psi \) onto the \(E\)-summand of
\(\bbS\) and differentiating to obtain an equation for \((\nabla_{A})^{q} \alpha\).  Given
the first bullet of Lemma \ref{lem:B.6} and given Lemma \ref{lem:B.5}, the latter implies a
differential inequality for the function \(\sigma :=|(\nabla_{A})^{q}\alpha|\) of the form \(d^{\dag}d \sigma
+ r |w| \sigma \leq c_{q} e^{-\sqrt{r}/c_0}\)
when \(q = 1\), and it implies an equality of this same sort for \(q
> 1\) if the second bullet holds for all \(q'< q\).
 Here, \(c_{q}\) depends only on \(q\).  Use the
Green's function for \(d^{\dag}d + r|w| \)with this differential
inequality for \(\sigma\) to prove the third bullet's assertion. \epf

Lemma \ref{lem:B.6} in turn leads to the
 
\noindent\textit{Proof of Lemma \ref{lem:A.1}.}  If \(r \geq
c_{0}\), then Lemma \ref{lem:B.6} asserts that
\(|\alpha|\) is very close to \(1\) on a neighborhood of
$\gamma$ and so what is denoted in (\ref{eq:(A.7)}) as
\(\wp(|\alpha|)\) is equal to \(1\) on this
neighborhood.  With this in mind, note that \(\alpha|\alpha|^{-1}\) is \(\hat{A}\)-covariantly
constant where \(\wp = 1\).  This implies that \(\hat{A}\) has holonomy \(1\) along
$\gamma$.  Since \(A_{E}\) has holonomy \(1\) on $\gamma$, it
follows that \(\hat{A} - A_{E}\) on $\gamma$ can be written as \(i \,\hat{u} (t)
\, dt\) with \(\hat{u}\) being a function on \(\bbR/(\ell_\gamma\bbZ)\) whose integral is an integer multiple of $2\pi $.\epf

\subsection{The \(L^{1}\)-norm of \(B_A\) when \(w\) is harmonic}
\label{sec:Bd)}

 This section supplies a crucial bound for the integral of
\(|B_{A}|\) over \(Y_{Z}\) given an extra assumption about \(w\).  
 
\begin{lemma}\label{lem:B.7}   
Suppose that \(w\)  is a harmonic
2-form and that the zeros of \(w\)  are non-degenerate.  There
exists $\kappa \geq \pi $  with the following
significance: Fix \(r \geq \kappa  \) and a 1-form
\(\mu \in \Omega  \) with $\mathcal{P}$-norm less than 1.  Suppose that \((A,
\psi)\)  is a solution to the \((r,  \mu)\)-version of (\ref{eq:(A.4)}).  Then  
\(\int_{Y_Z}|w| |B_A| \leq \kappa    \) and \( \int_{Y_Z} |B_A |\leq \kappa r^{1/5}\).
\end{lemma}

By way of a look ahead, the lemma's bound of \(\kappa r^{1/5}\) for the \(L^{1}\)-norm of
\(B_{A}\) is replaced in Lemma \ref{lem:B.9} by the bound \((\ln r)^{c_0}\).
 
\pf The proof has three steps.  By
way of an overview, the plan is to compare the integrals of
\(|B_{A}|\) and \(|w|\, |B_{A}|\) with the integral of \(w\wedge iB_{A}\).  The point being that the absolute
value of the latter integral enjoys an \((A, \psi)\)-, \(r\)- and \(\mu\)-
independent bound by virtue of the fact that \(w\) is harmonic; it
computes the cup product pairing between the de Rham class of \(*w\) and
$2\pi $ times the first Chern class of the bundle \(E\).
 
\paragraph{\it Step 1:}  Fix \(m \in (c_{0}, c_{0}r^{1/3} (\ln r)^{-c_0})\) so as to invoke Lemmas \ref{lem:B.2} and \ref{lem:B.3}.  Use \(U_{m}\) to again
denote the part of \(Y_{Z}\) where \(|w|> m^{-1}\).  Since \(w\) has non-degenerate
zeros, the volume of \(Y_{Z}-U_{m}\) is less than \(c_{0}m^{-3}\).  Since
\(|B_{A}| \leq c_{0} r\, (|\psi|^{2} +|w|) + c_{0}\), it follows from Lemma \ref{lem:B.3}
that 
\begin{equation}\label{eq:(B.11)}
\int_{Y_Z-U_m} |B_A|\leq c_{0} r m^{-4}
\quad\text{and}\quad 
\int_{Y_Z-U_m} |w\wedge B_A |\leq c_{0} r m^{-5}.
\end{equation}
Save these bounds for the moment.  

\paragraph{\it Step 2:}  Fix \(m \in (c_{0}, c_{0}r^{1/3})\).  Use the equations in (\ref{eq:(A.4)})
and Lemma \ref{lem:B.2} to see that \(|B_{A}|\) on
\(U_{m}\) obeys \(|B_{A}|\leq r |w| (|1-|\alpha|^{2}| +|\beta|) + c_{0}\).  This
understood, the first and second bullets in Lemma \ref{lem:B.2} imply that
\begin{equation}\label{eq:(B.12)}
|B_{A}| \leq c_{0} r \, |w| (1-|\alpha|^{2}) + c_{0}|w| m^{3}   
\end{equation}
at all points in \(U_{m}\).  Meanwhile, use the equations in
(\ref{eq:(A.4)}) to see that
\begin{equation}\label{eq:(B.13)}
w \wedge iB_{A} \geq r\, |w|^{2} (1-|\alpha|^{2}) -c_{0}|w|
\end{equation}
on \(U_{m}\).  This lower bound and the upper bound in (\ref{eq:(B.12)})
imply that if \(q \in\{0, 1\}\), then
\begin{equation}\label{eq:(B.14)}
|w|^{q}|B_{A}| \leq c_{0}m^{1-q} (w \wedge iB_{A}) +c_{0} m^{2-q} 
\end{equation}
at all points in \(U_{m}\).  
 
\paragraph{\it Step 3:}  Fix for the moment \(m_{0} \geq
c_{0}\) and a positive integer \(N\) with an upper bound such that
\(2^{N}m_{0} < c_{0}^{-1}r^{1/3}\). For \(k \in \{1, 2,\ldots, N\}\),
set \(m_{k}:=2^{k}m_{0}\).  Noting that the volume of
\(U_{m_k}-U_{m_{k-1}}\) is bounded by \(c_{0} 2^{-3k}\), it follows
from (\ref{eq:(B.14)}) that 
\begin{equation}\label{eq:(B.15)}
\int_{U_{m_k}-U_{m_{k-1}}} |w|^q |B_A|
 \leq c_{0}m_{N}^{1-q} \int_{U_{m_k}-U_{m_{k-1}}} w\wedge iB_A + c_{0} 2^{-k}. 
\end{equation}
 
Sum the various \(k \in \{1, \ldots, N\}\) versions of (\ref{eq:(B.15)}) to see that
\begin{equation}\label{eq:(B.16)}
 \int_{U_{m_N}}|w|^q |B_A| \leq c_{0}m_{N}^{1-q}\int_{U_{m_N}}w\wedge iB_A + c_{0}.
\end{equation}
 
This last inequality and the \(m = m_{N}\) version of (\ref{eq:(B.11)})
imply that 
\begin{equation}\label{eq:(B.17)}
\int_{Y_Z} |B_A|\leq c_{0}m_{N}^{1-q}\int_{Y_Z}w\wedge iB_A + c_{0} (r m_{N}^{-4-q} +1).
\end{equation}

The integral on the right hand side of (\ref{eq:(B.17)}) is in any event
bounded by \(c_{0}\) and so what is written in (\ref{eq:(B.17)}) leads to
the bound
\begin{equation}\label{eq:(B.18)}
\int_{Y_Z} |w|^q |B_A|\leq c_{0} (m_{N}^{1-q} + r m_{N}^{-4-q}).
\end{equation}
 
This understood, take \(N\) so that \(r^{1/5} \leq m_{N} \leq c_{0} r^{1/5}\)
to obtain Lemma \ref{lem:B.7}'s assertion. \epf

\subsection{Where \(1-|\alpha|^{2}\) is not small}\label{sec:Be)}

Suppose that \((A, \psi)\) is a solution to a given \((r,
\mu)\)-version of (\ref{eq:(A.4)}).  Write \(\psi\) where
\(|w| > 0\) as \(|w|^{1/2}(\alpha, \beta)\) and denote the version of \(\kappa\) that appears in Lemma \ref{lem:B.4}
by \(\kappa_{\diamond}\). 
 
 The lemma that follows momentarily characterizes the
\(|w|  > r^{-1/3} (\ln r)^{\kappa_\diamond}\) part of \(Y_{Z}\) where \(1 -|\alpha|^{2}\) is not very small.
 To set the notation for the lemma, introduce \(v\) to denote the unit
length vector field on the part of \(Y_{Z}\) where
\(|w| >0\) that generates the kernel of \(w\)
and has positive pairing with \(*w\). A final bit of
notation concerns the version of \(\kappa\) that appears in Lemma
\ref{lem:B.2}. The latter is denoted in what follows by \(\kappa_\diamond\).
\begin{lemma}\label{lem:B.8}  
Assume that \(w\)  is a harmonic
2-form with non-degenerate zeros. There exists \(\kappa> \kappa_{\diamond}\) with the
following significance:  Fix \(r \geq \kappa  \) and \(\mu \in \Omega  \) with
$\mathcal{P}$-norm bounded by \(1\) and let \((A, \psi)\)  denote a
solution to the \((r, \mu)\)-version of (\ref{eq:(A.4)}).  
Fix a positive integer \(k\)  and set
\(m_{k} := (1 +\kappa^{-1})^{k}\kappa^{2}\).  If \(m_{k}< r^{1/3} (\ln
r)^{-\kappa}\), then there exists a set \(\Theta_{k}\), 
of at most \(\kappa\)  segments of integral curves of
\(v\)  with the following properties:
\begin{itemize}
\item Each segment from \(\Theta_{1}\)  is
properly embedded in the \(|w| \geq m_{2}^{-1}\)  part of \(Y_{Z}\)
and has length at most \(\kappa\). Moreover, the union of the radius \(\kappa
r^{-1/2}\)  tubular neighborhoods of the segments
in \(\Theta_{1}\)  contain all points in the \(|w| >\kappa^{-2}\)  part of
\(Y_{Z}\)  where \(1 -|\alpha|^{2} > \frac{1}{4}\kappa_{\diamond}^{-1}\).
\item If \(k> 1\), then each segment from \(\Theta_{k}\)  is properly embedded in the
\(|w| \in[m_{k+1}^{-1},m_{k-1}^{-1}]\)  part of
\(Y_{Z}\)  and the union of the radius \(\kappa m_{k}^{1/2}
r^{-1/2}\)  tubular neighborhoods of the segments
in \(\Theta_{k}\)  contain all \(1-|\alpha|^{2} >\frac{1}{4} \kappa_{\diamond}^{-1} \)
points in the \(|w| \in[m_{k+1}^{-1},m_{k-1}^{-1}]\) part of \(Y_{Z}\).  
\end{itemize}
\end{lemma}
 
\pf The proof has 8 steps.  By way
of a parenthetical remark, the proof follows a strategy like that used
in Section IV.2c to prove Proposition IV.2.4.
 
 {\it Step 1:}  This step states a fact about the finite energy solutions
to the vortex equations that plays a central role in the subsequent
arguments.  Keep in mind that a solution \((A_{0},
\alpha_{0})\) is a finite energy solution when \(1 -|\alpha_{0}|^{2}\) is
an \(L^{1}\)-function.  As noted in Part 1 of Section \ref{sec:Bb)},
if \((A_{0}, \alpha_{0})\) is a finite energy solution then the integral of \(1-
|\alpha_{0}|^{2}\) is $2\pi $ times a non-negative integer.  Use
\(\m\) to denote this integer.
The function \(\alpha_0\) vanishes at precisely \(\m\) points in \(\bbC\) (with repetitions allowed). This set of zeros of \(\alpha_0\)  is denoted by \(\vartheta\).
 As noted in Part 4 from Section 2a in \cite{T2}, 
\begin{equation}\label{eq:(B.19)}
1 -|\alpha_{0}|^{2}\leq c_{0}\sum_{z\in\vartheta}e^{-\dist(\cdot,z)}.
\end{equation}
with the number \(c_0\) in (\ref{eq:(B.19)}) being independent of \((A_{0},
\alpha_{0})\) and \(\m\). The bound in (\ref{eq:(B.19)}) 
with Lemma \ref{lem:B.4} has a number of consequences with regards to
the proof.  

To say more, return to the context of Lemma \ref{lem:B.4}.  Let
\(\kappa_{\diamond}\) denote the version of the constant \(\kappa\) that appears in this lemma.  Take
\(\textsc{r} >\kappa_{\diamond}\) so as to apply the Lemma
\ref{lem:B.4} when \(r\) is greater than the corresponding 
\(\kappa_{\textsc{r}}\).  With \(r \geq \kappa_{\textsc{r}}\) and \(\mu \in \Omega\)
with $\mathcal{P}$-norm bounded by 1, let \((A, \psi)\) denote a
solution to the \((r, \mu)\)-version of (\ref{eq:(A.4)}).  Fix \(p \in
Y_{Z}\) with |\(w(p)| \geq r^{-1/3}(\ln r)^{\kappa_\diamond}\) and use \(p\) to define the pair \((A_{r},
\alpha_{r})\) as instructed in Part 2 of Section \ref{sec:Bb)}.
 Assume for what follows that \(1 -|\alpha_{r}|^{2}< \frac{1}{2}\) at distances between \(\textsc{r} +\kappa_{\diamond} (\ln \textsc{r})^{2}\) and \(\textsc{r} -\kappa_{\diamond} (\ln \textsc{r})^{2}\) from the origin in \(\bbC\).   

Lemma \ref{lem:B.4} asserts that \((A_{r}, \alpha_{r})\) has \(C^{1}\)-distance at most
\(\textsc{r}^{-4}\) in the radius \(\textsc{r}\) disk about
the origin in \(\bbC\) from a finite energy vortex that defines a
point in some $\m \leq \pi \textsc{r}^{2}$  version of \(\grC_{\m}\). Let \((A_{0},
\alpha_{0})\) denote this solution.  It follows from
Lemma \ref{lem:B.4} that \(1-|\alpha_{0}|^{2}\) can be no greater than
 \(\frac{1}{2}+ 2\textsc{r}^{-4}\) at all points in \(\bbC\) with
distance between \(\textsc{r} -\kappa_{\diamond} (\ln
\textsc{r})^{2}\) and \(\textsc{r}\) from the origin in \(\bbC\)
(since otherwise, \((A_{0},
\alpha_{0})\) would have \(C^0\) distance greater than  \(\textsc{r}^{-4}\) in
the radius  \(\textsc{r}\) disk about the origin in \(\bbC\)). This implies that each zero of \(\alpha _0\) (which are the
points in the set \(\vartheta\) that appears in (\ref{eq:(B.19)}) has distance either less than \( \textsc{r}  -
\kappa_{\diamond} (\ln \textsc{r})^{2}\) from the 
origin in \(\bbC\) or it has distance greater than \(\textsc{r}\) from the origin in \(\bbC\). This understood, then it
follows as a consequence of (\ref{eq:(B.19)}) that
\(1-|\alpha_{0}|^{2}\leq\textsc{r}^{-4}\) on the annulus about the origin in \(\bbC\)
with inner radius  \(\textsc{r}  -
\kappa_{\diamond} (\ln \textsc{r})^{2} + c_{0} \ln  \textsc{r}\) and
outer radius  \(\textsc{r}-c_{0} \ln  \textsc{r}\). Indeed, at
distance \(\rho\) from the set \(\vartheta\), the sum on the right hand
side of (\ref{eq:(B.19)}) is at most \(c_0\m e^{-\rho}\). Since  $\m <
\pi \textsc{r}^{2}$, this is at most \(c_0\, \textsc{r}^2
e^{-\rho}\). Thus, if \(\rho > c_0\ln \textsc{r}\), then the sum on the
right hand side of (\ref{eq:(B.19)}) will be at most
\(\textsc{r}^{-4}\). Granted that
\(1-|\alpha_{0}|^{2}\leq\textsc{r}^{-4}\) on the annulus in \(\bbC\) centered at the 
origin with inner radius equal to  \(\textsc{r}  -
\kappa_{\diamond} (\ln \textsc{r})^{2} + c_{0} \ln  \textsc{r}\) and
outer radius  \(\textsc{r}-c_{0} \ln  \textsc{r}\), it then 
follows from Lemma  \ref{lem:B.4} that \(1-|\alpha_{r}|^{2}\leq
2\textsc{r}^{-4}\) on this same annulus.

 
If \(\textsc{r} > c_{0}\), then the preceding conclusion implies that \(1-|\alpha|^{2}\) is bounded by \(2\textsc{r}^{-4}\) on the annulus in
transverse disk centered at \(p\) with respective outer and inner radii
given by \((\textsc{r}-c_{0} \ln  \textsc{r})\,( r|w|(p))^{-1/2}\) and inner radius \((\textsc{r}-\kappa_{\diamond} (\ln 
\textsc{r})^{2} + c_{0} \ln  \textsc{r})\,( r|w|(p))^{-1/2}\). Since
\(\alpha\) is nowhere vanishing on this annulus, the connection
\(\hat{A}_*\) is defined on this annulus by the same formula
(\ref{eq:(A.7)}), and the last observation implies in particular
that the connection \(\hat{A}_{*}\) is flat and \(\alpha|\alpha|^{-1}\) is
\(\hat{A}_{*}\)-covariantly constant at points on this same annulus.   

In the applications to come, the integer \(\m\) will be bounded by
\(c_{0}\).  If this is the case, then (\ref{eq:(B.19)}) with Lemma \ref{lem:B.4}
implies that \(\hat{A}_{*}\) is flat and \(\alpha|\alpha|^{-1}\) is
\(\hat{A}_{*}\)-covariantly constant at all point on the radius
\((\textsc{r}  - c_{0} (\ln \textsc{r})^{2}) (r|w|(p))^{-1/2}\)
transverse disk centered at \(p\) except at distance less than
\(c_{0} (r|w|(p))^{-1/2}\) from a set of at most \(c_{0}\) points.

\paragraph{\it Step 2:}  Fix \(m_{0} > c_{0}\) so that the \(|w| \leq m_{0}^{-1}\) part of \(Y_{Z}\) is
a disjoint union of components with each component lying in the radius
\(c_{0} m_{0}^{-1}\) ball about
a zero of \(w\).  Require in addition that each such component lie in a
Gaussian coordinate chart centered on the nearby zero of \(w\) as the
embedded image of a closed ball in $\mathbb{R}^{3}$.

Fix \(z > m_{0}\) and let \(\kappa_{0}\) denote the sum of the versions of
\(\kappa\) that appear in Lemmas \ref{lem:B.1}, \ref{lem:B.2} and \ref{lem:B.7}; and let
\(\kappa_{z_0}\) denote the sum of \(\kappa_{0}\) and the \(\textsc{r} = z^{10}\)
version of the constant \(\kappa_{\textsc{r}}\) that
appears in Lemma \ref{lem:B.4}.  But for cosmetic changes, the arguments in
Section 6.4 of \cite{TW1} can be used with Lemmas \ref{lem:B.2}, \ref{lem:B.4} and \ref{lem:B.6} plus what
is said in Step 1 to find a \(z\)-independent \(\kappa_{1}\geq 100\kappa_{0}\) and a \(z\)-dependent
\(\kappa_{z} > \kappa_{z_0}\) such that the following is true:

{\em Fix \(r\geq \kappa_z\) and \(\mu\in \Omega\) with
  \(\mathcal{P}\)-norm bounded by 1. Suppose that \((A, \psi)\) is a
  solution to the \((r, \mu)\)-version of (\ref{eq:(A.4)}). There
  exists a positive integer \(n_0<\kappa_1\) and a set \(\Theta_0\),
  of at most \(n_0\) pairs of the form \((\gamma, \m)\) with
  \(\gamma\) being a properly embedded segment of an integral curve of
  \(v\) in the \(|w|\geq z^{-6}\) part of \(Y_Z\) with length less
  than \(\kappa_1\). Meanwhile, \(\m\) is a positive integer. The set
  \(\Theta_0\) has the following additional properties:}
\BTitem\label{eq:(B.20)}
 \item \(\sum_{(\gamma, \m)\in \Theta_0}\m \leq \kappa_{1}\).
\item Distinct curves from \(\Theta_{0}\)  are
separated by distance at least \(\kappa_{1}z^{4} r^{-1/2}\).
\item If \(p \in Y_{Z}\)  is such that \(|w\, (p)| \geq z^{-6}\)  and \(1-|\alpha|^{2} > \kappa_{\diamond}^{-1}\),
 then \(p\) has distance less than \(z^{4}r^{-1/2}\)  from a curve in \(\Theta_{0}\).  
\item If \((\gamma, \m) \in \Theta_{0}\), then the integral of
  \(\frac{i}{2\pi}F_{\hat{A}_*}\) over the radius \(z^{4}r^{-1/2}\)
  transverse disk centered at each point in $\gamma$  is equal to \(\m\).
\ETitem
 
What follows is a parenthetical remark concerning the fourth bullet.
 The condition in the third bullet of (\ref{eq:(B.20)}) implies that \(\alpha
|\alpha|^{-1}\) is  \(\hat{A}_{*}\)-covariantly constant near the
boundary of the radius \(z^{4} r^{-1/2}\)
transverse disk about each point in $\gamma$.  It follows as a
consequence that the integral of \(\frac{i}{2\pi}F_{\hat{A}_*}\)   
over this disk is an integer; and it follows from Lemma \ref{lem:B.2} that this
integer is non-negative.  This being the case, the fourth bullet adds
only that the integer is at least 1 and it is bounded a priori by a \(z\),
\((A, \psi)\)-, \(\mu\)- and \(r\)-independent number.
 
{\it Step 3:}  Fix a ball $B \subset Y_{Z}$  centered on a zero
of \(w\) that contains a component of the \(|w| \leq m_{0}^{-1}\) part of \(Y_{Z}\).
Suppose that $\varepsilon \in (0, 1)$ and that \(z >m_{0}\) have been specified.  With
\(\kappa_{z}\) as in Step 2, fix \(r \geq\kappa_{z}\), an element \(\mu \in \Omega\)
with $\mathcal{P}$-norm bounded by 1 and a solution, \((A, \psi)\),
to the \((r, \mu)\)-version of (\ref{eq:(A.4)}).  Let \(k\) denote the largest
integer with the properties listed below in (\ref{eq:(B.21)}).  By way of
notation, set \(m_{j}:=(1 +\varepsilon)^{j} z^{6}\).  
 
For each \(j \in \{1, \ldots, k\}\), there exists \(\c_{j} \in (100, (100)^{2^{\kappa_1}})\) and a set, \(\Theta_j\), that consists of data sets which have the form \((\gamma,\m, \textsc{d})\)  with \(\gamma\) being a properly
embedded segment of an integral curve of \(v\)  in the
\(|w| \in[m_{j+1}^{-1},m_{j-1}^{-1}]\)  part of \(B\),
with \(\m\)  being a positive integer and with \(\textsc{d}
\in (1, \c_{j})\). The set \(\Theta_j\) has the following additional properties:  
\BTitem\label{eq:(B.21)}
\item \(\sum _{\gamma, \m)\in \Theta_j}\m \leq \kappa_{1}\).
\item  Curves from distinct data sets in \(\Theta_j\) are separated by distance at least 
\(\frac{1}{2} \c_{j}^{2} zm_{j}^{1/2}r^{-1/2}\).
\item  If \(p \in Y_{Z}\)  is such that \(|w(p)| \in[m_{j+1}^{-1},m_{j-1}^{-1})\)  and \(1-
|\alpha|^{2} > \frac{1}{4} \kappa_{0}^{-1}\), then \(p\)
 has distance at most \(\textsc{d}zm_{j}^{1/2}r^{-1/2}\)  from a point on a curve from a datavset in \(\Theta_j\).  
\item  If \((\gamma, m, \textsc{d}) \in   \Theta_j\), then the integral of \(\frac{i}{2\pi}F_{\hat{A}_*}\)
over the radius \(\textsc{d}zm_{j}^{1/2}r^{-1/2}\)  transverse disk centered at each
point in \(\gamma\)  is equal to \(\m\).
\ETitem
 
The next steps find \((A, \psi)\), \(\mu\)- and \(r\)-independent choices
for \(\varepsilon\) and then \(z\), and an \((A, \psi)\), \(\mu\) and
\(r\)-independent \(\kappa_{*} \geq\kappa_{z}\) such that \(m_{k} \geq
r^{1/3}(\ln r)^{-\kappa_*}\) when \(r\) is greater than \(\kappa_{*}\).  Lemma \ref{lem:B.8}
follows if such \(\varepsilon\), \(z\) and \(\kappa_{*}\) exist.  
 
The upcoming steps find the desired conditions on \(\varepsilon\), \(z\) and
the lower bound for \(r\) such that the conditions of the integer \(k + 1\)
version of (\ref{eq:(B.21)}) are met if they are met for an integer \(k\) with
\(m_{k} <  r^{1/3}  (\ln r)^{-2\kappa_0}\).  This being the strategy, assume in what follows that \(k\) is such that
\(m_{k} < r^{1/3} (\ln r)^{-2\kappa_0}\) and (\ref{eq:(B.21)}) holds.

 {\it Step 4:}  The \(A_{*}\)-directional covariant derivative along the
vector field \(v\) is used momentarily to analyze the behavior of
\(\alpha\) at points along \(v\)'s integral curves.  This
directional derivative is denoted in what follows by
\((\nabla_{A}\alpha)_{v}\).  The
equations in (\ref{eq:(A.4)}) identify the latter with a linear combination of
\(A_{*}\) -covariant derivatives of \(\beta\).  This being
the case, Lemma \ref{lem:B.2} finds
\(|(\nabla_{A}\alpha)_{v}|\leq c_{0} m [(1-|\alpha|^{2}) + c_{0}r^{-1}
m^{3}]^{1/2}\) on the \(|w| > (2m)^{-1}\) part of \(Y_{Z}\) if \(m \leq
r^{1/3}(\ln r)^{-\kappa_0}\).  By way of a comparison, Lemma
\ref{lem:B.2} bounds the norm of the remaining components of \(\nabla_{A}\alpha \) by
\(c_{0} m^{-1/2} r^{-1/2} [(1  - |\alpha|^{2}) +c_{0}r^{-1}m^{3}]^{1/2}\).
 
What was said in the preceding paragraph about the norm of
\(|(\nabla_{A}\alpha)_{v}|\) has the following consequences for a point \(p \in Y_{Z}\)
where \(|w| \in[m_{k+2}^{-1},m_{k}^{-1}]\):  Let
\(\gamma_{p}\) denote the integral curve of \(v\) through \(p\)
and let \(p'\) denote a point on the segment of \(\gamma_{p}\)
where the distance to \(p\) is less than \(c_{0}^{-1}\kappa_{\diamond}^{-1}m_{k}^{-1}\).
\BTitem\label{eq:(B.22)}
\item If \(1  - |\alpha|^{2}> \frac{1}{4}\kappa_{\diamond}^{-1}\)  at \(p\), then \(1-|\alpha|^{2} > \frac{1}{8}\kappa_{\diamond}^{-1}\)  at \(p'\).
\item If \(1  - |\alpha|^{2}\leq \frac{1}{4}\kappa_{\diamond}^{-1}\)  at \(p\),
  then \(1-|\alpha|^{2} < \frac{1}{2}\kappa_{\diamond}^{-1}\)  at \(p'\).
\ETitem
 
This segment of \(\gamma_{p}\) is said in what follows to
be the  {\em short} segment of \(\gamma_{p}\). 

Note that if \(\varepsilon \leq c_{0}^{-1}\kappa_{\diamond}^{-2}\), then \(\gamma_{p}\)'s short segment has
points with \(|w| >m_{k-1}^{-1}\).  Assume in what follows
that \(\varepsilon \leq c_{0}^{-1}\kappa_{\diamond}^{-2}\) is
satisfied so as to invoke this fact about the short segment.

{\it Step 5:}  This step constitutes a digression to supply a coordinate
chart for any given \(|w| > 0\) point in
\(Y_{Z}\) that is used to exploit what is said in Step 4.  To
this end, suppose that \(m > 1\) has been specified.  Use
\(I_{m}\) to denote the interval \([-c_{0}^{-1}m^{-1}, c_{0}^{-1}m^{-1}]\) and
use \(D_{m}\) to denote the centered disk in \(\bbC\) with radius \(c_{0}^{-1}m^{-1}\).
Use \(t\) to denote the coordinate for the interval
\(I_{m}\) and use \(z\) for the complex coordinate on
\(D_{m}\). As will be explained momentarily, there is a
coordinate chart embedding from \(I_{m}\times D_{m}\) to \(Y_{Z}\) with the following properties:  
\BTitem\label{eq:(B.23)}
\item The point \((0, 0)\) is mapped to \(p\)  and \(I_{m} \times
  \{0\}\) is mapped to a segment of the integral curve of \(v\) through \(p\).
\item The image of any disk \(\{t\} \times D_{m}\) is a transverse disk centered at the image
of \((t, 0)\).
\item The function \(z \mapsto|z|\) on \(\{t\} \times D_{m}\)  is the pull-back of the
distance along the image of \(\{t\} \times D_{m}\)  to the image of \(\{t, 0\}\). 
\item The vector field \(v\) appears in these coordinates as 
\(\frac{\partial  }{\partial t}+ \gre\)  where \(|\gre| \leq c_{0}m |z|\).
\ETitem
 
To construct such a coordinate chart, fix an isometric isomorphism
between \(K_{*}^{-1}|_{p}\) and \(\bbC\).  By way of a reminder,
\(K_{*}^{-1}\) is used to denote the complex line bundle over the \(|w(p)| > 0\)
part of \(Y_{Z}\) whose underlying real bundle is the kernel
of \(*w\) with the complex structure defined using the metric and the
restriction of the form \(w\).  Let \(\gamma_{p}\) again
denote the integral curve of \(v\) through \(p\).  Parallel transport the
resulting frame for \(K_{*}^{-1}\) along \(\gamma_{p}\) to identify \(K_{*}^{-1}\) along
\(\gamma_{p}\) with \(\gamma_{p}\times\bbC\).  Fix a unit length affine parameter, \(t\), for
the segment of \(\gamma_{p}\) consisting of points with
distance \(c_{0}^{-1}m^{-1}\) or less from \(p\) with \(t = 0\)
corresponding to \(p\).  This identifies this segment with \(I_{m}\).
Granted this identification, compose the metric's exponential map from the
\(I_{m}\)  part of \(\gamma_{p}\) with the identification between \(K_{*}^{-1}\)
on this segment and the product \(\bbC\) bundle to define a map
from \(I_{m} \times \bbC \) into \(Y_{Z}\).  The restriction of this map to \(I_{m}
\times D_{m}\) gives the desired coordinate embedding.

{\it Step 6:}  Fix \(p \in Y_{Z}\) such that
\(|w(p)| \in[m_{k+2}^{-1},m_{k}^{-1}]\) and \(1-|\alpha|^{2} >\frac{1}{4}\kappa_{\diamond}^{-1}\).  Let \(p'\)
denote a chosen point on Step 4's short segment of
\(\gamma_{p}\) with\( |w(p')| =m_{k+1}^{-1}\).  It follows from (\ref{eq:(B.22)})
that \(1 - |\alpha|^{2}>\frac{1}{8} \kappa_{\diamond}^{-1}\) at \(p'\).
 This being the case, it follows from Lemma \ref{lem:B.4} and Lemma
 IV.2.8 that if \(z > c_{0}\) and if \(r >c_{0}\), then there is a
 point, with distance at most \(c_{0} m_{k+1}^{1/2}r^{-1/2}\) from \(p'\) where \(1  -
|\alpha|^{2} > \frac{1}{4} \kappa_{\diamond}^{-1}\).   It
then follows from the third bullet of (\ref{eq:(B.21)}) that there exists
\((\gamma, \m, \textsc{d}) \in \Theta_{k}\) such
that \(p'\) has distance at most \((\textsc{d}z  + c_{0}) m_{k+1}^{1/2} r^{-1/2}\)
from a point in \(\gamma\).  Let \(p_{*}\) denote the
latter point.  Use the coordinate chart in (\ref{eq:(B.23)}) to see that short
segment of \(\gamma_{p}\) intersects the transverse disk
through \(p_{*}\) at a point with distance at most \((1 +
c_{0}\varepsilon)(\textsc{d} z + c_{0})m_{k+1}^{1/2} r^{-1/2}\) from \(p_{*}\).  
 
Extend the curves from \(\Theta_{k}\) into the
\(|w| \geq m_{k+2}^{-1}\) part of \(Y_{Z}\)
by integrating the vector field \(v\).  Use \(\Gamma_{k+1}\)
to denote this set of extended curves.  Given \(\gamma \in
\Gamma_{k+1}\), fix a point \(p_{\gamma}\in \gamma \) where \(|w|=m_{k+1}^{-1}\).  The point
\(p_{\gamma}\) has its corresponding version of the
coordinate chart in (\ref{eq:(B.23)}) with \(\gamma\) appearing as an interval in
the \(z = 0\) locus that contains \((0, 0)\).  Let \(I_{\gamma}\) denote this interval.  
 
It follows from what was said in the preceding paragraph that the each
point in \(B\) where \(1-|\alpha|^{2} > \frac{1}{4} \kappa_{\diamond}^{-1}\) and
 \(|w| \in[m_{k+2}^{-1},m_{k}^{-1}]\) lies in the
\(|z| \leq (1 + c_{0}\varepsilon)(\textsc{d}z + c_{0})m_{k+1}^{1/2}
r^{-1/2}\) part of some \(\gamma \in \Gamma_{k+1}\) version
of \(I_{\gamma} \times D_{m_{k+1}}\).  In particular, if \(\varepsilon <c_{0}^{-1}\) and \(z >
c_{0}\), then this subset is contained in the subset where
\(|z| < \frac{3}{2} \textsc{d} z m_{k+1}^{1/2}r^{-1/2}\).  Assume that \(\varepsilon\) and \(z\) are such
that this is the case.  

 Note in this regard that if \((\gamma, \m, \textsc{d})\) and
\((\gamma', \m', \textsc{d}')\) are distinct elements in
\(\Theta_{k}\), then the respective subsets of \(B\) that are
parametrized via (\ref{eq:(B.23)}) by the \(|z| \leq 2\textsc{d} z m_{k+1}^{1/2}r^{-1/2}\) part of
\(I_{\gamma}\times D_{m_{k+1}}\) and the \(|z| \leq 2 \textsc{d}' z m_{k+1}^{1/2} r^{-1/2}\)
part of \(I_{\gamma'}\times D_{m_{k+1}}\)
are disjoint.  This is a consequence of the second bullet in (\ref{eq:(B.21)}).
  
\paragraph{\it Step 7:}  Fix \((\gamma,  \m, \textsc{d}) \in\Theta_{k}\).   It follows from what was said in Step
6 that \(\alpha |\alpha|^{-1}\)
is \(\hat{A}_{*}\)-covariantly constant in the solid annulus in
\(I_{\gamma} \times D_{m_{k+1}}\) that intersects any constant \(t\) slice as the annulus with inner radius \(\frac{3}{2} \textsc{d} z m_{k+1}^{1/2}r^{-1/2}\) and outer radius \(2\textsc{d} z
m_{k+1}^{1/2} r^{-1/2}\). Granted this, it then follows from the third bullet of (\ref{eq:(B.21)}) that
the integral of \(\frac{i}{2\pi}F_{\hat{A}_*}\) over the \(|z| < 2 \textsc{d} zm_{k+1}^{1/2} r^{-1/2}\)
part of any constant \(t\) disk in \(I_{\gamma} \times  D_{m_{k+1}}\) is the integer \(\m\).
 
To exploit the preceding observation, fix \(t \in I_{\gamma}\) and let \(p \in Y_{Z}\) denote
the point that corresponds to \((t, 0) \in I_{\gamma} \times
D_{m_{k+1}}\).  Associate to \(p\) the pair \((A_{r}, \alpha_{r})\) as
desribed in Part 2 of Section \ref{sec:Bb)}.  Use \(c_{z}\) in what follows to denote a constant that is greater
than 1 and depends only on \(z\).  It follows from Lemma \ref{lem:B.4} that if \(z
> c_{0}\) and if \(r >c_{z}\), then \((A_{r}, \alpha_{r})\) have \(C^{1}\)-distance less
than \(z^{-10}\) on the radius \(2\textsc{d} z\) disk in
\(\bbC\) from a finite energy solution to the vortex equations.
 Moreover, what is said by Lemma \ref{lem:B.4} implies that such a finite energy
solution must define a point in the space \(\grC_{\m}\).  Granted
this, then (\ref{eq:(B.19)}) and Lemma \ref{lem:B.4} imply the following when \(z
> c_{0}\) and \(r >c_{z}\):     
 
 {\em If \(z > c_{0}\)  and \(r \geq c_{z}\), then there is a set of at most
\(n_{0}\)  points in the \(|z|<\frac{3}{2} \textsc{d}
m_{k+1}^{1/2}r^{-1/2}\)  part of \(\{t\} \times  D_{m_{k+1}}\) such that}
\BTitem\label{eq:(B.24)}
\item Each point is a zero of \(\alpha\).
\item If \(1- |\alpha|^{2}\geq \frac{1}{8}\kappa_{\diamond}^{-1}\)  at
  \((t, z)\)  and \(|z| \leq 2 \textsc{d} \, m_{k+1} r^{-1/2}\), then \(z\)  has distance at
most \(c_{0} m_{k+1}^{1/2}r^{-1/2}\) from some point in this set.
\ETitem
 
Use \(\vartheta_{\gamma ,t}\) to denote this set of points and let \(\grU_{\gamma ,t}\) denote the set of
connected components of the union of the disks of radius
\(c_{0} m_{k+1}^{1/2}r^{-1/2}\) about the points in \(\vartheta_{\gamma ,t}\).  The next
assertion is a \(z > c_{0}\) and \(r >c_{z}\) consequence of (\ref{eq:(B.24)}) plus Lemma \ref{lem:B.4} and (\ref{eq:(B.19)}). 
\BTitem\label{eq:(B.25)}
\item The connection \(\hat{A}_{*}\)  is flat and \(\alpha|\alpha|^{-1}\)  is
\(\hat{A}_{*}\)-covariantly constant  on the complement of \(\bigcup_{U\in
  \grU_{\gamma ,t}}U \) in the  radius \(2 \textsc{d} m_{k+1}^{1/2} r^{-1/2}\)
 disk about the origin in \(\{t\} \times D_{m_{k+1}}\).
\item The integral of \(\frac{i}{2\pi}F_{\hat{A}_*}\) over any set \(U \in
\grU_{\gamma .t}\)  is a positive integer; and the sum of these integers is equal to \(\m\).
\ETitem
 
The next step constructs \(\Theta_{k+1}\) with the help of the various \((\gamma, \m, \textsc{d}) \in
\Theta_{k}\) versions of \(\vartheta_{\gamma ,0}\).

\paragraph{\it Step 8:}  To construct \(\Theta_{k+1}\), it is
necessary to cluster the points from the various \((\gamma,
\m,\textsc{d}) \in \Theta_{k}\) versions of \(\vartheta_{\gamma ,0}\) so that
points in the same cluster are pairwise much closer to each other than
they are to any point in another cluster.  This is necessary so as to
find the desired constant \(\c_{k+1}\) for the integer \(k+1\)
version of (\ref{eq:(B.21)}).   An appropriate clustering can be found by
invoking Lemma 2.12 in \cite{T4}.  In particular, an appeal to this
lemma finds \(\c_{k+1} \in (100, (100)^{2^{\kappa  _1}})\) and a set of at most \(\kappa_{1}\) pairs of
the form \((p, \textsc{d})\) where \(p \in B\) is such that
\(|w(p)| =m_{k+1}^{-1}\) and where \(\textsc{d} \in
(1, c_{k+1})\).  This set is denoted by \(\vartheta\) and
it has the properties in the list that follows.
\BTitem\label{eq:(B.26)}
\item If \((p, \textsc{d})\)  and \((p', \textsc{d}')\)
 are distinct elements in \(\vartheta\), then \(\dist\, (p,p') >
 \c_{k+1}^{2} z\, m_{k+1}^{1/2} r^{-1/2}\).
\item If \(p\) corresponds via (\ref{eq:(B.23)})  to a point in some
\((\gamma, \m, \textsc{d}) \in \Theta_{k}\) version of \(\vartheta_{k,0}\), then \(p\) has distance at most 
 \(\frac{1}{4}\textsc{d} z\, m_{k+1}^{1/2}r^{-1/2}\)  from a point of some pair from
\(\vartheta\).
\ETitem
 
Note for future reference that the bound in the first bullet of (\ref{eq:(B.26)})
has the following implication when \(z > c_{0}\) and \(r > c_{z}\):
\begin{equation}\label{eq:(B.27)}
\begin{split}
 \, &\textrm{If \((p, \textsc{d})\)  and \((p',\textsc{d}')\)  are distinct elements in \(\vartheta\),
then the distance between } \\ 
\, &\textrm{any two points on the respective short
segments \(\gamma_{p}\)  and \(\gamma_{p'}\)  is greater}\\ 
\, &\textrm{than \(\frac{1}{2}\c_{k+1}^{2} z
m_{k+1}^{1/2}r^{-1/2}\).}
\end{split}
\end{equation}

It follows from (\ref{eq:(B.25)}) and (\ref{eq:(B.26)}) that if \((\gamma, \m, \textsc{d})
\in \Theta_{k}\) and if \(U \in \grU_{\gamma ,0}\), then \(U\) is in the transverse disk of radius
 \(\frac{1}{2}\textsc{d} z m_{k+1}^{1/2}r^{-1/2}\) centered at a point of some pair in
\(\vartheta\). Granted this last conclusion, then the next
assertion is a direct consequence of what is said in Step 4 if \(z
> c_{0}\) and \(r >c_{z}\).  
\begin{equation}\label{eq:(B.28)}
\begin{split}
 \, &\textrm{If \((\gamma, m, \textsc{d}) \in
\Theta_{k} \) and \(t \in
I_{\gamma}\), then each \(U \in
\grU_{\gamma ,t}\)  is contained in the }\\
\, &\textrm{radius \(\textsc{d}\, z \,m_{k+1}^{1/2}
r^{-1/2}\)  tubular neighborhood of the integral
curve of \(v\)}\\  
\, &\textrm{through a point of some pair from
\(\vartheta\)}.  
\end{split}
\end{equation}

 Let \((p, \textsc{d}) \in \vartheta\).  What is said in
(\ref{eq:(B.27)}) and (\ref{eq:(B.28)}) has the following consequence:  
\begin{equation}\label{eq:(B.29)}
\begin{split}
 \, &\textrm{The integral of $\frac{i}{2\pi}F_{\hat{A}_*}$
  on the radius \(\textsc{d} \, z\, m_{k+1}^{1/2} r^{-1/2}\)
 transverse disk about}  \\
\, &\textrm{any point in the \(|w| \in[m_{k+2}^{-1},m_{k}^{-1}]\)  part of
\(\gamma_{p}\)  is a positive integer.}
\end{split}
\end{equation}
Let \(\m\) denote now this integer.  

 Define \(\Theta_{k+1}\) to be the set \(\{(p, m, \textsc{d})\, |\, 
(p, \textsc{d}) \in \vartheta)\}\).  It follows from (\ref{eq:(B.25)}) and
(\ref{eq:(B.27)})-(\ref{eq:(B.29)}) that the requirements for the integer \(k + 1\) version of
(\ref{eq:(B.21)}) are met using \(\c_{k+1}\) and the set
\(\Theta_{k+1}\) if \(\varepsilon < c_{0}\), \(z > c_{0}\) and \(r>
c_{z}\). \epf

 \subsection{The spectral flow function}\label{sec:Bf)}
  
This subsection constitutes a digression to say more about the
definition of \(\grf_{s}\).  Each pair \(\grc = (A, \psi)\) in
\(\op{Conn} (E) \times C^{\infty}(Y_{Z}; \bbS)\) and a given real number \(z\) determine an associated,
unbounded, self-adjoint operator on
\(L^{2}(Y_{Z}; iT^*Y_{Z}\oplus \bbS \oplus i\bbR)\).  This operator is
denoted by \(\grL_{\grc, z}\) and it is defined as follows:  A
given smooth section \(\grh = (b, \eta, \phi)\) of  \(iT^*Y_{Z} \oplus \bbS\oplus i\bbR \) is sent by \(\grL_{\grc, z}\)
to the section whose respective \(i T^*Y_{Z}\), \(\bbS\), and \(i\bbR\)-summands are
\begin{equation}\label{eq:(B.30)}
\begin{cases}
*db - d\phi - 2^{-1/2} z^{1/2}(\psi^{\dag}\tau\eta +\eta^{\dag}\tau\psi),&\\
 D_{A}\eta + 2^{1/2}z^{1/2}(\cl(b)\psi + \phi\psi),&\\
 *d*b - 2^{-1/2}z^{1/2}(\eta^{\dag}\psi -\psi^{\dag}\eta).&
\end{cases}
\end{equation}
 
The spectrum of this operator is discrete with no accumulation points
and has finite multiplicity.  The spectrum is also unbounded from above
and unbounded from below.

The section \(\psi_{E}\) of \(\bbS\) is chosen so that
the \((A_{E}, \psi_{E})\) and \(z = 1\) version
of (\ref{eq:(B.30)}) has trivial kernel.  If the \(z = r\) and \(\grc = (A, \psi)\)
version of (\ref{eq:(B.30)}) has trivial kernel, then the value of the spectral
flow function \(\grf_{s}(\grc)\) is a certain algebraic count of the
number of zero eigenvalues that appear along a continuous path
\(\grd\) of
operators that start at the \(z = 1\) and \((A_{E}, \psi_{E})\)
version of (\ref{eq:(B.30)}) and end at the \(z = r \) and
\((A, \psi)\) version and such that each member of the path differs
from \(\grL_{\grc ,r}\) by a bounded operator on
\(L^{2}(Y_{Z}; iT^*Y_{Z}\oplus \bbS \oplus i\bbR)\).  For the purposes
of the definition, it is sufficient to consider paths that are
parametrized by \([0, 1]\) such that the following conditions are met:
 Let \(\vartheta \subset [0, 1]\) denote the parameters that label
an operator with zero as an eigenvalue.  Then \(\vartheta\) is finite
and in each case, the zero eigenvalue has multiplicity 1 and the zero
eigenvalue crossing is transversal as the parameter varies in a small
neighborhood of the given point in \([0,1]\).  Having chosen such a
path, a given point in the corresponding version of \(\vartheta\)
contributes either \(+1\) or \(-1\) to \(\grf_{s}(\grc)\).  The
contribution is \(+1\) when the eigenvalue crosses zero from negative value
to positive value as the parameter in \([0, 1]\) varies near the given
point in \(\vartheta\); and it contributes \(-1\) to \(\grf_{s}(\grc)\)
if the eigenvalue crosses zero from a positive value to negative value
near the given point.

If \(\grL_{\grc ,r}\) has non-trivial kernel, then
\(\grf_{s}(\grc)\) is defined in the upcoming (\ref{eq:(B.31)}).  The
definition uses the following terminology:  Given \(\varepsilon> 0\), and \(\grc \in\op{Conn} (E) \times
C^{\infty}(Y_{Z}; \bbC)\), the
definition uses \(\grN_\varepsilon(\grc)\) to denote the
subset of pairs in \(\op{Conn} (E) \times
C^{\infty}(Y_{Z}; \bbC)\) with
the following two properties:  A pair \(\grc '\) is in
\(\grN_\varepsilon(\grc)\) if it has \(C^{1}\)-distance less than \(\varepsilon\) from \(\grc\), and if
\(\grL_{\grc ',r}\) has trivial kernel.  Standard
perturbation theory for ellipitic operators proves that
\(\grN_\varepsilon(\grc)\) is non-empty for any
\(\varepsilon > 0\).  With the notation set, define
\(\grf_{s}(\grc)\) by the rule
\begin{equation}\label{eq:(B.31)}
\grf_{s}(\grc) = \lim\sup_{\varepsilon\to 0}
\{\, \grf_{s}(\grc ')\, |\,  \grc ' \in \grN_\varepsilon(\grc)\}.
\end{equation}
 
Note by the way that the \(\lim \sup\) in (\ref{eq:(B.31)}) differs from the
corresponding \(\lim \inf\) by the dimension of the kernel of
\(\grL_{\grc ,r}\).

\subsection{The \(L^{1}\)-norm of \(B_{A}\), the spectral flow and the
functions \(\grc\grs^\grf\), \(\textsc{w}^\grf\), \(\gra^\grf\)} \label{sec:Bg)}
 
The functions 
\begin{equation}\label{eq:(B.32)}
\grc\grs^{\grf} = \grc\grs - 4\pi ^{2}\grf_{s} , \quad
\textsc{w}^{f} = \textsc{w} -2\pi  \grf_{s }  \quad\text{and}  \quad 
\gra^{\grf} = \gra +2\pi  (r  - \pi ) \grf_{s} 
\end{equation}
are invariant under \(C^{\infty}(Y_{Z};S^{1})\) action on \(\op{Conn}
(E) \times C^{\infty}(Y_{Z}; \bbS)\) that
has \(\hat{u} \in C^{\infty}(Y_{Z};S^{1})\) sending \((A, \psi)\) to \((A  -
\hat{u}^{-1}d\hat{u},\hat{u} \psi)\).   The upcoming Lemma \ref{lem:B.9}
supplies a priori bounds on the values of these functions when
evaluated on solutions to a given \((r, \mu)\)-version of (\ref{eq:(A.4)}).  It
also gives a better bound for the \(L^{1}\)-norm of the
curvature of the connection component of a solution than the bound in
Lemma \ref{lem:B.6}.
 
\begin{lemma}\label{lem:B.9}   
Suppose that \(w\)  is a harmonic
2-form with non-degenerate zeros.  There exists $\kappa
> \pi $  with the following significance:  Fix \(r\geq \kappa  \) and a 1-form \(\mu \in\Omega  \) with
$\mathcal{P}$-norm less than 1.
Suppose that \((A, \psi)\) is a solution to the
\((r, \mu)\)-version of (\ref{eq:(A.4)}).  Then:
\begin{itemize}\itemsep -.2pt
\item  The \(L^{1}\)-norm of \(B_{A}\)  is no greater than \(\kappa
\,( \ln r)^{4}\).
\item \(|\grc\grs^{\grf}| <r^{6/7}\),
\item  \(|\textsc{w}^{\grf}| <r^{6/7}\),
\item  \(|\gra^{\grf}| <r^{13/14}\).
\end{itemize}
\end{lemma}
As a parenthetical remark, the precise powers of \(r\) that appear in the
last three bullets are significant with regards to the applications to
come only to the extent that the power is less than 1 in the second and
third bullets and so less than 2 in the final bullet.

\pf By way of a look ahead, what is
said in Lemma \ref{lem:B.8} plays a vital role in the proof of all four bullets.
 The proof of Lemma \ref{lem:B.9} has 10 parts.  
 
 \paragraph{Part 1:} The proof of Lemma \ref{lem:B.9}'s first
bullet has four steps.  To set the notation for the proof, introduce
\(\kappa_{*}\) to denote the version of the constant
\(\kappa\) that appears in Lemma \ref{lem:B.8}.  As in Lemma \ref{lem:B.8}, set
\(m_{k} = (1 +\kappa_{*}^{-1})\kappa_{*}^{2}\) for \(k \in \{1,
2,\ldots\}\).  Assume in what follows that \(k\) is such that \(m_{k}
< r^{1/3}(\ln r)^{-\kappa_*}\).

{\it Step 1:}  Use the first bullet of Lemma \ref{lem:B.8} and the fourth bullet of
Lemma \ref{lem:B.2} to see the \(|B_{A}| <c_{0}\)  at points in the \(|w|>m_{1}^{-1}\) part of
\(Y_{Z}\) where the distance to all segments in \(\Theta_{1}\) is greater than \(c_{0} (\ln 
r)^{2} r^{-1/2}\).  This understood, this part of \(Y_{Z}\) contributes at most \(c_{0}\)
to the \(L^{1}\)-norm of \(B_{A}\).  Meanwhile,
the \(|w|>m_{1}^{-1}\) part of \(Y_{Z}\) of
the union of the radius \(c_{0} (\ln  r)^{2}r^{-1/2}\) tubular neighborhoods of the segments in
\(Y_{Z}\) contributes at most \(c_{0} (\ln r)^{4}\) to the
\(L^{1}\)-norm of \(B_{A}\).
 
{\it Step 2:}  Fix \(k >1\).  Use the integer \(k\) version of the
second bullet of Lemma \ref{lem:B.8} with the fourth bullet of Lemma \ref{lem:B.2} to see
that \(|B_{A}|\) is bounded by
\(c_{0} (1 + m_{k}^{2})\) at
points in the \(|w| \in[m_{k}^{-1},m_{k-1}^{-1}]\) part of 
\(Y_{Z-}\) where the distance to all segments in
\(\Theta_{k}\) is greater than \(c_{0}m_{k}^{1/2 } (\ln r)^{2} r^{-1/2}\).  Since this subset of
\(Y_{Z}\) has volume at most \(c_{0}m_{k}^{-3}\), so this portion of the
\(|w| \in[m_{k}^{-1},m_{k-1}^{-1}]\) subset in
\(Y_{Z}\) contributes at most \(c_{0}m_{k}^{-1}\) to the
\(L^{1}\)-norm of \(B_{A}\).  The volume of the remaining part of the \(|w| \in[m_{k}^{-1},m_{k-1}^{-1}]\) subset in 
\(Y_{Z-}\) is at most \(c_{0}r^{-1} (\ln  r)^{4}\).  Indeed, this can
be seen from (\ref{eq:(B.23)}) using the fact that each segment in
\(\Theta_{k}\) has length at most \(c_{0}m_{k}^{-1}\).  As \(|B_{A}|\) is no greater than
\(c_{0} m_{k}^{-1} r\) on this part of \(Y_{Z}\), so this part of \(Y_{Z}\)
contributes at most \(c_{0}m_{k}^{-1}(\ln r)^{4}\) to the \(L^{1}\)-norm of \(B_{A}\). 
 
{\it Step 3:}  Lemma \ref{lem:B.3} implies that \(|B_{A}|\)
is bounded by \(c_{0} r^{2/3}(\ln r)^{\kappa_*}\) on the subset of \(Y_{Z}\) where \(|w|
\leq c_{0} r^{-1/3}(\ln r)^{\kappa_*}\).  The volume of this subset is at most \(r^{-1}(\ln r)^{3\kappa_*}\)  and so the contribution from this part of \(Y_{Z}\) to the
\(L^{1}\)-norm of \(B_{A}\) is no greater than \(c_{0} r^{-1/4}\). 

 {\it Step 4:}  Sum the bounds in Steps 1-3 to see that the
\(L^{1}\)-norm of \(B_{A}\) is no greater than
\(c_{0} (\ln  r)^{4}\sum_{k=0,1,\ldots}  (1 +1/\kappa_{*})^{-k}\).  This sum is
bounded by \(c_{0}\kappa_{*} (\ln r)^{4}\).

\paragraph{Part 2:}  The proof of the last three bullets of the lemma
starts with the following observation:  There is a smooth map,
\(\hat{u}\co 
Y_{Z} \to S^{1}\), such that the connection \(A' = A -\hat{u} ^{-1}d\hat{u}\) can
be written as \(A' = A_{E} + \hata _{A'}\) where
\(\hata _{A'}\) is a coclosed, \(i \bbR \)-valued 1-form whose
\(L^{2}\) orthogonal projection to the space of harmonic
1-forms on \(Y_{Z}\) is bounded by \(c_{0}\).  The
upcoming Lemma \ref{lem:B.10} asserts the pointwise bound
\(|\hata _{A'}| \leq c_{0}r^{1/3} (\ln r)^{c_0} \). 
Assume this bound for the time being.  
 
Introduce \(\grc'\) to denote \((A  -\hat{u} ^{-1}d\hat{u},\hat{u} \psi)\).
The supremum bound for \(|\hata _{A'}|\) and the \(L^{1}\)-bound for \(B_{A}\) from Lemma
\ref{lem:B.9}'s first bullet imply directly that
\(|\grc\grs(\grc')| \leq c_{0}r^{1/2} (\ln r)^{c_0}\).  The \(L^{1}\)-bound for \(B_{A}\) also implies
that  \(|\textsc{w}(\grc')| \leq c_{0}(\ln r)^{c_0}\).  Thus, \(|\gra(\grc')| \leq c_{0} r (\ln r)^{c_0}\).
Granted these bounds, then the last three bullets of Lemma \ref{lem:B.9}
follow if
\begin{equation}\label{eq:(B.33)}
|\grf_{s}(\grc')| \leq r^{6/7}.
\end{equation}
 
The fact that (\ref{eq:(B.33)}) holds given the assumptions of the lemma is proved
in the remaining parts of this subsection.    
 
\paragraph{ Part 3:}  The proof of the last three bullets of Lemma \ref{lem:B.9}
invoked a pointwise bound for \(|\hata _{A'}|\).
 The lemma that follows supplies the asserted bound.
 
\begin{lemma}\label{lem:B.10}   
There exists $\kappa >\pi $  with the following significance:  Fix \(r \geq
\kappa  \) and an element \(\mu \in \Omega \) with $\mathcal{P}$-norm less than \(1\). Let
\((A, \psi)\) denote a solution to the \((r, \mu)\)-version of (\ref{eq:(A.4)}).  Write \(A\)  as
\(A_{E} + \hata _{A}\)  and assume that
\(\hata _{A}\)  is a coclosed 1-form. Use \(\c\) to
denote the \(L^{2}\)-norm of the \(L^{2}\)-orthogonal projection of
\(\hata _{A}\)  to the space harmonic 1-forms.  Then
\(|\hata _{A}| \leq r^{1/2}(\ln r)^\kappa+ \kappa \c \).
 \end{lemma}
 
\pf The proof that follows has three steps.
 
\paragraph{\it Step 1: } Write \(\hata _{A}\) as \(\hata ^{\perp} + \grp\) where  \(\hata ^{\perp}\) is
\(L^{2}\)-orthogonal to the space of harmonic 1-forms and
where \(\grp\) is a harmonic 1-form.  The norm of \(\grp\) is bounded by
\(c_{0}\c\).  To bound \(\hata ^{\perp}\), let
\(\mathcal{C}^{\perp} \subset C^{\infty}(Y_{Z}; T^*Y_{Z})\)
denote the subspace of coclosed 1-forms that are \(L^{2}\)-orthogonal to the space of harmonic 1-forms.  The operator \(*d\) maps
\(\mathcal{C}^{\perp}\) to itself and Hodge theory gives a
Green's function inverse.  Given \(p \in M\), the
corresponding Green's function with pole at \(p\) is denoted by \(G^{\perp}_{p}(\cdot)\).
 This function is smooth on the complement of \(p\) and it obeys the pointwise
bound \(|G^{\perp}_{p}(\cdot)|\leq c_{0} \dist\, (\cdot,p)^{-2}\).
 
\paragraph{\it Step 2:}  Introduce \(\kappa_{*}\) to denote Lemma
\ref{lem:B.8}'s version of \(\kappa\).  Reintroduce from Lemma
\ref{lem:B.8} the sequence \(\{m_{k} = (1 +\kappa_{*}^{-1})^{k}\kappa^{2}\}_{k=1,2\ldots, N}\) with \(N\)
being the greatest integer such that \(m_{k} <
r^{1/3} (\ln r)^{-\kappa  _*}\).  Let $\mathcal{U}_{1}$ denote the
\(|w|> m_{2}^{-1 }\) part of \(Y_{Z}\). For \(k \in \{1, \ldots, N-1\}\), use
$\mathcal{U}_{k}$ to denote the \(|w|\in [m_{k+1}^{-1},m_{k-1}^{-1}]\) part of \(Y_{Z}\),
and use $\mathcal{U}_{N}$ to denote the part of
\(Y_{Z}\) where \(|w| \leq m_{N-1}^{-1}\).   Given \(k \in \{1,\ldots,
N-1\}\), let $\Gamma_{k}$ denote the set of curves from \(\Theta_{k}\)'s data sets.
 By way of a reminder, there are at most $\kappa_{*}$
curves in $\Gamma_{k}$ and each is a properly embedded
segment of an integral curve of \(v\) in $\mathcal{U}_{k}$.
 
Lemmas \ref{lem:B.2} and \ref{lem:B.8} supply \(c_{*} \in (1,c_{0})\) with the following property:  If \(p \in\mathcal{U}_{k}\) has distance greater than
\(c_{*} m_{k} r^{-1/2} (\ln r)^{2}\) to any curve from \(\Gamma_{k}\),
then \(1  - |\alpha|^{2} \leq c_{0} m_{k}^{3}r^{-1}\).  Denote by $\mathcal{T}_{k1}$
the union of the radius \(\c_{*}m_{k}r^{-1/2} (\ln  r)^{2}\) tubular
neighborhoods of the curves from $\Gamma_{k}$.  Since
\(*d\hata ^{\perp} = B_{A}\), it follows from
Lemmas \ref{lem:B.1} and \ref{lem:B.2} that \(|B_{A}|\leq c_{0} m_{k}^{2}\) on
$\mathcal{U}_{k}-\mathcal{T}_{k1}$, and it follows
from Lemma \ref{lem:B.2} and Lemma \ref{lem:B.3} that
\(|B_{A}| \leq c_{0}m_{k}^{-1} r\) on \(\mathcal{T}_{k1}\).
Note also that the volume of $\mathcal{U}_{k}$ is at most \(c_{0} m_{k}^{-3}\) and that
of \(\mathcal{T}_{k1}\) at most \(c_{0}m_{k}^{-1}r^{-1} (\ln  r)^{4}\). 
 
\paragraph{\it Step 3:}  Suppose that \(k \in \{1, \ldots, N-1\}\) and that 
$p \in\mathcal{U}_{k}$.  Keeping in mind that the volume of
$\mathcal{U}_{k}$ is bounded by \(c_{0}m_{k}^{-3}\), it follows from what is said
about \(G^{\perp}_{p}\) in Step 1 and what is said about \(|B_{A}|\) in Step 2 that
\begin{equation}\label{eq:(B.34)}
|\hata ^{\perp}|(p) \leq c_{0 }\int_{\mathcal{T}_k}\dist\, (\cdot, p)^{-2}
  |B_A| +  c_{0} (m_{k} +  (\ln r)^{c_0}).
\end{equation}
Use the various \(\gamma \in \Gamma_{k}\) versions
of (\ref{eq:(B.23)}) to see that the integral on the right hand side
of (\ref{eq:(B.34)}) is no greater than \(c_{0}m_{k}^{-1/2}r^{1/2}(\ln r)^{c_0}\).  

Suppose that $p \in \mathcal{U}_{N}$.  In this case,
what is said about \(G^{\perp}_{p}\) in Step 1 and what is said in Step 2 about
\(|B_{A}|\) imply that \(|\hata ^{\perp}|(p) \leq c_{0} r^{1/3} (\ln r)^{c_0}\). \epf

\paragraph{Part 4:}  Fix \(\c > 1\) and suppose that \(\grc = (A,
\psi)\) solves (\ref{eq:(A.4)}) and is such that the \(i \bbR \)-valued
1-form \(\hata _{A} = A - A_{E}\) is coclosed and that
the \(L^{2}\)-norm of its \(L^{2}\)-orthogonal
projection to the space of harmonic 1-forms on \(Y_{Z}\) is
less than \(\c\).  The value of \(\grf_{s}\) will be computed by
choosing a convenient, piecewise continuous path of self-adjoint
operator from the \((A_{E}, \psi_{E})\) and \(z=1\) version of (\ref{eq:(B.30)}) to \(\grL_{c,r}\).  This
path is the concatentation of the three real analytic segments that are
described below.  The absolute value of \(\grf_{s}(\grc)\) is no
greater than the absolute value of the sum of the absolute values of
the spectral flow along the three segments.
 
By way of notation, each segment is parametrized by \([0, 1]\) and the
operator labeled by a given \(s \in[0, 1]\) in the \(k\)'th
segment is denoted by $\mathcal{L}_{k,s}$.  The first
segment's operator $\mathcal{L}_{1,s}$ for
\(s \in [0, 1]\) is the \((A_{E}, \psi_{E})\)
and \(z = 1- s\) version of (\ref{eq:(B.30)}).  This path has no dependence on \((A,
\psi)\) or \(r\), and so the absolute value of the spectral flow along
this path is no greater than \(c_{0}\).  The remaining two segments are: 
\BTitem\label{eq:(B.35)}
\item The second segment's operator $\mathcal{L}_{2,s}$  for \(s \in [0, 1]\) 
is the \((A_{E} + s \, \hata _{A}, 0)\) version of (\ref{eq:(B.30)}).
\item The third segment's operator $\mathcal{L}_{3,s}$  for \(s \in [0, 1]\) 
is the \((A, \psi)\), \(z = s^{2} r\)-version of (\ref{eq:(B.30)}).
\ETitem
 
The strategy for bounding the absolute value of the spectral flow along
(\ref{eq:(B.35)})'s two segments borrows heavily from Section 3 of
\cite{TW2}.  To say more about this, suppose that \(\mathcal{L}\) is an unbounded,
self-adjoint operator on a given separable Hilbert space with discrete
spectrum with no accumulation points and finite multiplicities.  Let
\(\{e_{s}\}_{s\in[0,1]}\)
denote a real analytic family of bounded, self-adjoint operators on
this same Hilbert space.  Of interest is the spectral flow between the
\(s = 0\) and \(s = 1\) members of the family \(\{\mathcal{L}_{s}= \mathcal{L} +e_{s}\}_{s\in[0,1]}\).
 To obtain a bound, fix for the moment \(T > 0\) and let
\(\grn_{T,s}\) denote the number of linearly independent
eigenvectors of \(\mathcal{L}_{s}\) whose eigenvalue has
absolute value no greater than \(T\).  Set \(\grn_{T} =\sup\, \{\grn_{T,s}\}_{s\in[0,1]}\).
 As explained in \cite{Tasf}, the spectral flow for the family
\(\{\mathcal{L}_{s}\}_{s\in[0,1]}\)
has absolute value no greater than
\begin{equation}\label{eq:(B.36)}
  \frac{1}{2T}\grn_{T}\sup \, \Big\{\, \Big\|  \frac{d}{ds}e_{s}\Big\|_{op}\Big\}_{s\in[0,1]},
\end{equation}
where the norm \(\| \cdot\|_{op}\) here denotes the operator norm.
 
The supremum in (\ref{eq:(B.36)}) for the family \(\{\mathcal{L}_{2,s}\}_{s\in[0,1]}\)
is bounded by \(c_{0}|\hata _{A'}|\), and thus by \(c_{0}r^{1/2 }(\ln r)^{c_0}\).
It follows from Lemma \ref{lem:B.3} that the supremum that appears in
(\ref{eq:(B.36)}) for the family \(\{\mathcal{L}_{3,s}\}_{s\in[0,1]}\)
is \(c_{0} r^{1/2}\).  This understood, then (\ref{eq:(B.36)}) in either case leads to
\begin{equation}\label{eq:(B.37)}\begin{split}
 &\textrm{The absolute value of the spectral flow along the families}\\
&\textrm{\(\{\mathcal{L}_{2,s}\}_{s\in[0,1]}\) and
\(\{\mathcal{L}_{3,s}\}_{s\in[0,1]} \) is no greater than \(c_{0} r^{1/2}  (\ln r)^{c_0}\frac{1}{T} \grn_{T}\).}
\end{split}
\end{equation}
The next part of the subsection describes the strategy that is used to
bound \(\grn_{T}\) for a suitable choice of \(T\).

 \paragraph{Part 5:}  A bound for \(\grn_{T}\) is obtained with
the help of the Weitzenb\"ock formula in (IV.A.12) for a given \(z \geq 0\) version of  
\(\grL_{\grc, z}^{2}\).   This formula writes \(\grL_{\grc, z}^{2}\) as
$\nabla_{A}^{\dag}\nabla_{A}+ \mathcal{Q}$ where $\mathcal{Q}$ denotes an endomorphism of \(i
T^*Y_{Z} \oplus \bbS \oplus i\bbR\) and \(\nabla_{A}\) denotes here the connection on the
bundle \(i T^*Y_{Z} \oplus \bbS \oplus i\bbR \) given by the Levi-Civita connection on the \(i
T^*Y_{Z}\)-summand, the Levi-Civita connection and \(A\) on the
\(\bbS\) summand, and the product connection on the \(i\bbR\) 
summand.  This rewriting of  \(\grL_{\grc, z}^{2}\) is used to write the square of the
\(L^{2}\)-norm of \(\grL_{\grc, z}\grq\) as
\begin{equation}\label{eq:(B.38)}
\int_{Y_Z} |\grL_{\grc, z}\grq\, |^2 = \int_{Y_Z} \big(|\nabla_A\grq|^2+\langle \grq, \mathcal{Q}\grq\rangle \big),
\end{equation}
with \(\langle \cdot,\cdot \rangle\) denoting here the Hermitian inner product
on \(i T*Y_{Z} \oplus \bbS \oplus i\bbR\).  If \(\grq\) is a linear combination of eigenvectors of  \(\grL_{\grc, z}\)
with the norm of the eigenvalue bounded by \(T\), then what is written in
(\ref{eq:(B.38)}) is no greater than \(T^{2}\) times the square of the
\(L^{2}\)-norm of \(\grq\).  
 
The formula in (\ref{eq:(B.38)}) is exploited to bound \(\grn_{T}\) using
the following observation:  Suppose that \(\grU\) is an open cover of
\(Y_{Z}\) such that no point is contained in more than
\(c_{0}\) sets from \(\grU\).  Let \(h\) denote for the moment a given
function on \(Y_{Z}\).  Then
\begin{equation}\label{eq:(B.39)}
\int_{Y_Z} h^2 \leq \sum_{U\in\grU}  \int_U h^2\leq  c_{0}  \int_{Y_Z} h^2.
\end{equation}
 Hold onto this last observation for the moment.  Use
\(c_{\diamond}\) to denote the version of \(c_{0}\) that appears in this last inequality.

The endomorphism $\mathcal{Q}$ is self-adjoint, so it can be written
at any given point as a sum $\mathcal{Q}^{+} +\mathcal{Q}^{-}$  with
$\mathcal{Q}^{+}$ being positive semi-definite and $\mathcal{Q}^{-}$  being negative
definite.  With this fact in mind, suppose now that each set \(U \in
\grU\) has an assigned, finite dimensional vector subspace
\(V_{U} \in C^{\infty}(U; i T^*M\oplus \bbS \oplus i\bbR)\) with the
following property: 

\begin{equation}\label{eq:(B.40)}\begin{split}
 &\textrm{If \(q \in C^{\infty}(U; i T^*M \oplus\bbS \oplus i\bbR)\) is
\(L^{2}\)-orthogonal to \(V_{U}\), then}\\
&\int_U \big( |\nabla_A\grq |^2+\langle\grq, \mathcal{Q}^+\grq\rangle \big)
 > 2c_{\diamond} (T^{2} + \sup_{U}|\mathcal{Q}^{-}|) \int_U |\grq|^2.
\end{split}
\end{equation}
 
Given \(V_{U}\), define \(\Phi_{U}\co C^{\infty}(Y_{Z}; i T^*M \oplus\bbS \oplus i\bbR) \to
V_{U}\)  to be the composition of first restriction to \(U\)
and then the \(L^{2}\)-orthogonal projection.  Set
$\mathcal{V}=\bigoplus_{U\in\grU}V_{U}$  and denote by \(\Phi\) the linear map from
\(C^{\infty}(Y_{Z}; i T^*M \oplus\bbS \oplus i\bbR)\) to \(\mathcal{V}\) given by \(\bigoplus_{U\in\grU}\Phi_{U}\).  

 The inequalities in (\ref{eq:(B.39)}) and (\ref{eq:(B.40)}) have the
 following immediate consequence:  If \(\grq \in \ker (\Phi)\), then the
\(L^{2}\)-norm of \(\grL_{\grc, z}\) is greater than \(T\).  Given that such is the case, it then follows
directly that \(\grn_{T} \leq\sum_{U\in\grU}\dim (V_{U})\).  
 
The subsequent parts of the proof define a version of \(\grU\) for suitable \(T\)
with associated vector spaces \(\{V_{U}\}_{U\in\grU}\)
such that (\ref{eq:(B.40)}) holds.  The resulting bound for \(\grn_{T}\)
leads via (\ref{eq:(B.37)}) to the bound in (\ref{eq:(B.33)}) for \(|\grf_{s}|\).

\paragraph{Part 6:}  Part 5 alludes to a certain open cover of
\(Y_{Z}\).  This part of the subsection defines this cover.
 To this end, reintroduce from Step 2 of the proof of Lemma \ref{lem:B.10} the
sets \(\{\mathcal{U}_{k}\}_{1\leq k\leq N}\). The cover in question is given as \(\grU =\bigcup_{k=1,2,\ldots N} \grU_{k}\) where all \(U\in \grU_{k}\) are subsets of
\(\mathcal{U}_{k-1} \cup \mathcal{U}_{k} \cup\mathcal{U}_{k+1}\).  The
definition requires the choice of a constant \(\c > 1\).  Part 10 of the proof gives a lower
bound for \(\c\) by \(c_{0}\).  Any choice above this bound
suffices.

 To define a given \(k \in \{1, \ldots, N-1\}\) version of
\(\grU_{k}\), reintroduce from Step 2 of the proof of Lemma \ref{lem:B.10}
the set \(\Gamma_{k}\), this being the set of curves from \(\Theta_{k}\)'s data sets.  By way of a
reminder, there are at most \(\kappa_{*}\) curves in
$\Gamma_{k}$ and each is a properly embedded segment of
an integral curve of \(v\) in \(\mathcal{U}_{k}\).  This same
step in the proof of Lemma \ref{lem:B.8} introduced a constant
\(c_{*}\) such that \(1-|\alpha|^{2} <c_{0} m_{k}^{3}r^{-1}\) at points with distance
\(c_{*}m_{k}^{1/2}r^{-1/2} (\ln  r)^{2}\) or more to all
curves from $\Gamma_{k}$.  The discussion that follows uses \(\textsc{r}_{k}\) to denote
\(c_{*}m_{k}^{1/2}r^{-1/2} (\ln  r)^{2}\) and
\(\rho_{k}\) to denote \(\c^{-1} \min (T, m_{k}^{-1})\).  
 
The collection \(\grU_{k}\) for \(k \in \{1, \ldots, N-1\}\)
is written as \(\grU_{k-} \cup \grU_{k0} \cup \grU_{k+}\).  The sets from \(\grU_{k-}\) are balls of
radius \(\rho_{k}\) whose centers have distance at least \(\rho_{k}\) to all curves from
\(\Gamma_{k}\).  These balls cover the complement in $\mathcal{U}_{k}$
of the union of the radius \(2\rho_{k}\) tubular neighborhoods of the curves from
\(\Gamma_{k}\).  A cover as just described can be found
with less than \(c_{0}\rho_{k}^{-3}m_{k}^{-3}\) balls, and such is the case
with the cover \(\grU_{k-}\).  
 
The sets from \(\grU_{k0}\) are balls with distance between
\(2\rho_{k}\) and \(\textsc{r}_{k}\) to at least one curve from \(\Gamma_{k}\).  Let \(U\) denote a give
ball from \(\grU_{k0}\) and let \(\textsc{d}\) denote its distance to
the union of the curves from $\Gamma_{k}$.  The radius
of \(U\) is equal to \(\frac{1}{8}\textsc{d}\).  The various \(\gamma \in
\Gamma_{k}\) versions of (\ref{eq:(B.23)}) can be used to see that
a collection of \(c_{0} \ln (\rho_{k}/\textsc{r}_{k})(\textsc{r}_{k}m_{k})^{-1}\)
balls of this sort can be found whose union contains every point in
$\mathcal{U}_{k}$ with distance between \(\rho_{k}\) and \(2\textsc{r}_{k}\) to at
least one curve from \(\Gamma_{k}\).  The set \(\grU_{k0}\) is such a collection of balls.
 
The set \(\grU_{k+}\) consists of balls of radius \(\c^{-1} m_{k}^{1/2}
r^{-1/2}\) whose centers have distance at most \(\textsc{r}_{k}\) to some curve from
\(\Gamma_{k}\).  The balls from \(\grU_{k+}\) cover
the set of points with distance \(\textsc{r}_{k}\) or less to
some curve from \(\Gamma_{k}\).  The collection
\(\grU_{k+}\) has at most \(c_{0} \c^{3}(\ln  r)^{4}m_{k}^{-3/2}r^{1/2}\) balls.
 
The sets that comprise \(\grU_{N}\) are balls of radius 
\(r^{-1/3}(\ln  r)^{-\c}\) with centers in $\mathcal{U}_{N}$.  These sets define an
open cover of $\mathcal{U}_{N}$.  A cover of this sort
can be found with less than \(c_{0}(\ln  r)^{c_0\c}\) elements, and such is the case for \(\grU_{N}\).

 \paragraph{Part 7:}  This part of the subsection defines the vector spaces
\(\{V_{U}\}_{U\in\grU}\).
 The next lemma is needed for the definition.  
 
\begin{lemma}\label{lem:B.11}   
There exists \(\kappa \geq 1\)  with the following significance:  Let \(U \subset
Y_{Z}\)  denote a ball of radius \(\rho \in(0, \kappa^{-1})\).  Fix an
isometric isomorphism between \(E|_{U}\)  and
\(U\times \bbC \).  Use the latter to view the product connection on \(U\times \bbC  \) as a
connection on \(E|_{U}\) .  Use \(\nabla_{0}\)  to denote the corresponding
covariant derivative on \(C^{\infty}(U; iT^*M\oplus \bbS \oplus i\bbR)\).  There
exists a \(\kappa\)-dimensional vector space \(W_{U}\in C^{\infty}(U; iT^*M \oplus \bbS
\oplus i\bbR)\)  such that if \(\grq\)  is a section over \(U\)  of \(iT^*M \oplus \bbS \oplus i\bbR  \) which is
\(L^{2}\)-orthogonal to \(W_{U}\), then \(\int_U |\nabla_0\grq |^2
\geq \kappa^{-1}\rho^{-2}\int_U |\grq|^2\).
\end{lemma}
This lemma will be proved momentarily; so assume it to be true for now.
 
Fix \(U \subset \grU\).  If \(\c \geq c_{0}\) then the radius
of each ball from \(\grU\) will be smaller than Lemma \ref{lem:B.12}'s
version of \(\kappa^{-1}\) and each ball from \(\grU\) will
sit in the Gaussian coordinate chart about its center point.  With
this understood, fix \(U \in \grU\) and let \(p\) denote \(U\)'s
center point.  Fix an isometric isomorphism between
\(E|_{p}\) and \(\bbC\) and use \(A\)'s parallel transport along the radial geodesics from
\(p\) to extend this identification to one between
\(E|_{U}\) and the product bundle \(U \times\bbC\).  Define \(V_{U}\) to be Lemma \ref{lem:B.11}'s vector space \(W_{U}.\)
 
\paragraph{\it Proof of Lemma \ref{lem:B.11}.}  If \(\rho <c_{0}^{-1}\), then \(U\) has a Gaussian
coordinate chart centered at its center point.  Fix an isometric
identification between \(K^{-1}\) at the center point of \(U\)
with \(\bbC\) and use the \(A_{K}\) parallel transport
along the radial geodesics through the center point to extend this
isomorphism to one between
\(K^{-1}|_{U}\) and \(U \times\bbC\).  Use the coordinate basis with the identification
\(K^{-1}|_{U} = U \times\bbC \) and the chosen identification \(E|_{U}
= U \times \bbC \) to give a product structure to \(T^*M\) and
\(\bbS\) over \(U\).  Having done so, rescale the coordinates so the
ball of radius \(\rho\) becomes the ball of radius 1; then invoke
the next lemma. \epf
 
\begin{lemma}\label{lem:B.12}   
Let \(U \subset \bbR^{3}\)  denote the ball of radius \(1\)
 centered on the origin.  If \(\grh \in C^{\infty}(U; \bbC)\)  is such that \(\int_U\grh= 0 \), then 
\(\int_U |d\grh |^2\geq\frac{1}{4}\int_U\grh^2\).
\end{lemma}
 
\pf It is sufficient to prove the
bound for functions  that depend only on \(z\) through its absolute value.
 This understood, use \(\rho\) to denote \(|z|\) and let
\(h\) denote a function that depends only on \(\rho\) and has integral
zero over the unit ball.  Let \(h_{*} = h  - h(1)\). Use
integration by parts to see that
\begin{equation}\label{eq:(B.41)}
\int_0^1\grh_*^2\rho^{-1}d\rho \leq 2 \int_0^1
|\frac{\partial}{\partial\rho}\grh|\, |\grh_*|\, \rho d\rho.
\end{equation}
 
What is written in (\ref{eq:(B.41)}) implies that
\begin{equation}\label{eq:(B.42)}
\int_0^1\grh_*^2d\rho\leq 4  \int_0^1|d\grh|^2\rho^2d\rho.
\end{equation}
 
Meanwhile, \(\int_0^1\grh_*^2d\rho\geq \int_0^1\grh_*^2\rho^2d\rho\),
the latter being the integral of
\(\grh_{*}^{2}\) over the unit ball.  This last
integral is \(\frac{1}{3}\grh(1)^{2}\) plus the integral of \(\grh^{2}\)
because the integral of \(\grh\) is zero. \epf

\paragraph{Part 8:}  This step sets the stage for the specification of \(\c\)
and \(\{\rho_{k}\}_{1\leq k\leq N-1}\)
so as to guarantee (\ref{eq:(B.40)}).  To start, let \(U \subset \grU\) denote a
given ball and let \(p\) denote the center point of \(U\).  Fix an isometric
isomorphism between \(E|_{p}\) and \(\bbC\) and
then use \(A\)'s parallel transport along the radial
geodesics from \(p\) to extend this isomorphism to give an isomorphism
between \(E|_{U}\) and \(U \times \bbC\).
 Let \(\theta_{0}\) denote the product connection on \(U
\times \bbC\).  Use the isomorphism just defined to view
\(\theta_{0}\) as a connection on
\(E|_{U}\).  Having done so, write \(A\) on \(U\) as
\(\theta_{0} + \hata _{A,U}\) with
\(\hata _{A,U}\) being an \(i \bbR \)-valued 1-form on \(U\).  Let
\(\textsc{d}_{U}\) denote the radius of \(U\).  The norm of
\(\hata _{A,U}\) is bounded by \(c_{0}\textsc{d}_{U} \sup_{U}|B_{A}|\). 
 
Fix \(k \in \{1, \ldots, N-1\}\); let \(U\) denote a ball from either
\(\grU_{k-}\) or \(\grU_{k0}\).  It follows from what
Lemma \ref{lem:B.2} that \(|B_{A}| \leq c_{0} m_{k}^{2}\) on \(U\) and so
\(|\hata _{A,U}| \leq c_{0}\c^{-1} \rho_{k}m_{k}^{2}\).  If \(U \in
\mathcal{U}_{k+}\), then it follows from Lemma \ref{lem:B.2} that
\(|B_{A}| \leq c_{0}m_{k}^{-1} r\) on \(U\) and so
\(|\hata _{A,U}| \leq c_{0}\c^{-1} m_{k}^{-1/2}r^{1/2}\) on \(U\).  If \(U\) is from \(\grU_{N}\), then
Lemma \ref{lem:B.3} finds \(|B_{A}| \leq c_{0} r^{2/3} (\ln r)^{c_0}\) on \(U\) and so \(|\hata _{A,U}| \leq
c_{0} c^{-1} r^{1/3} (\ln r)^{c_0}\).  
 
Given \(U \subset \grU\), use the isomorphism defined above between
\(E|_{U}\) and \(U\times \bbC \) to again view \(\theta_{0}\) as a connection on
\(E|_{U}\).  Use \(\nabla_{0}\) to denote the corresponding covariant derivative on
\(C^{\infty}(U; iT^*M \oplus \bbS\oplus V)\).  Since
\(|\hata _{A,U}|^{2} \leq c_{0 }\sup_{U}|B_{A}|\) in all cases, so 
\begin{equation}\label{eq:(B.43)}
|\nabla_{A}\grq|^{2}
\geq \frac{1}{2} |\nabla_{0}\grq|^{2} -c_{0}\, (\sup_{U}|B_{A}|)  |\grq|^{2}
\end{equation}
for all \(\grq \in C^{\infty}(U; iT^*M \oplus\bbS \oplus V)\).  
 
Consider next the endomorphism $\mathcal{Q}$ that appears in (\ref{eq:(B.40)}).
 A look at the formula in (IV.A.12) finds
\begin{equation}\label{eq:(B.44)}
|\mathcal{Q}^{-}| \leq c_{0} (1 +|B_{A}| +  z^{1/2}|\nabla_{A}\psi|)\quad
   \text{and}\quad
|\mathcal{Q}^{+}| \geq c_{0}^{-1} z  |\psi|^{2}.  
\end{equation}
 
To say more about the bounds in (\ref{eq:(B.44)}) on the sets from \(\grU\), fix first \(k
\in \{1, \ldots, N-1\}\) and let \(U\) denote a ball from
\(\grU_{k-}\) or \(\grU_{k0}\).  Lemma \ref{lem:B.2} finds
\(|\nabla_{A}\psi| \leq c_{0} m_{k}^{1/2}\) and
\(|\psi|^{2} \geq c_{0} m_{k}^{-1}\) on \(U\).
 Since \(|B_{A}|\) on \(U\) is bounded by
\(c_{0}m_{k}^{2}\), the inequalities in (\ref{eq:(B.43)}) and (\ref{eq:(B.44)}) imply that
\begin{equation}\label{eq:(B.45)}
|\nabla_{A}\grq|^{2} +\langle \grq, \mathcal{Q}^{+}\grq \rangle  \geq  
|\nabla_{0}\grq|^{2} +
2c_{\diamond} \sup_{U}
|\mathcal{Q}^{-}|\, |\grq|^{2} -c_{0}m_{k}^{2}|\grq|^{2} 
\end{equation}
for all \(\grq \in C^{\infty}(U; iT^*M \oplus\bbS \oplus V)\).  Meanwhile, if \(U\) is a ball from
\(\grU_{k+}\), then Lemma \ref{lem:B.3} finds
\(|\nabla_{A}\psi| \leq c_{0} m_{k}^{-1}r^{1/2}\) and \(|B_{A}|\leq c_{0}m_{k}^{-1} r\).
 This being the case, then (\ref{eq:(B.43)}) and (\ref{eq:(B.44)}) find
\begin{equation}\label{eq:(B.46)}
|\nabla_{A}\grq|^{2} +\langle \grq, \mathcal{Q}^{+}\grq \rangle  \geq  
\frac{1}{2}|\nabla_{0}\grq|^{2} +2c_{\diamond} \sup_{U}
|\mathcal{Q}^{-}|\, |\grq|^{2} - c_{0}m_{k}^{-1} r |\grq|^{2} 
\end{equation}
for all \(\grq \in C^{\infty}(U; iT^*M \oplus\bbS \oplus V)\).
 
Suppose next that \(U\) is a ball from \(\grU_{N}\).  What is
said in Lemma \ref{lem:B.3} implies that \(|B_{A}|\leq c_{0} r^{2/3}
(\ln r)^{c_0}\) and \(|\nabla_{A}\psi| \leq c_{0} r^{1/6}\)  on \(U\), so (\ref{eq:(B.43)}) and (\ref{eq:(B.44)})
lead to the inequality
\begin{equation}\label{eq:(B.47)}
|\nabla_{A}\grq|^{2} + \langle \grq, \mathcal{Q}^{+}\grq \rangle  \geq  \frac{1}{2} |\nabla_{0}\grq|^{2} +
2c_{\diamond} \sup_{U}
|\mathcal{Q}^{-}|\, |\grq|^{2} - c_{0}r^{2/3}(\ln r)^{c_0} |\grq|^{2} 
\end{equation}
for all \(\grq \in C^{\infty}(U; iT^*M \oplus\bbS \oplus V)\).

 \paragraph{Part 9:}  This part of the subsection specifies \(\c\) and
\(\{\rho_{k}\}_{1\leq k\leq N-1}\)
so as to satisfy (\ref{eq:(B.40)}).  To this end, suppose that \(k \in \{1,
\ldots, N-1\}\).  Suppose that \(U\) is from \(\grU_{k-}\) or
\(\grU_{k0}\).  If \(\grq \in C^{\infty}(U;iT^*M \oplus \bbS \oplus V)\) is
\(L^{2}\)-orthogonal to the subspace \(V_{U}\),
then Lemma \ref{lem:B.11} and (\ref{eq:(B.45)}) find
\begin{equation}\label{eq:(B.48)}
 \int_U\big( |\nabla_A\grq|^2+\langle \grq, \mathcal{Q}^+\grq\rangle\big)\geq (c_{0}^{-1}
\rho_{k}^{-2} -c_{0}m_{k}^{2 }+2c_{\diamond} \sup_{U}
|\mathcal{Q}^{-}|)\int_U |\grq|^2.
\end{equation}
 
It follows as a consequence that (\ref{eq:(B.40)}) holds if
\(\rho_{k}^{-2} \geq c_{0}  (T^{2} +m_{k}^{2})\) and this is so if
\(\c > c_{0}\).  Suppose next that \(U\) is from
\(\grU_{k+}\) and that \(\grq \in C^{\infty}(U;iT^*M \oplus \bbS \oplus V)\) is
\(L^{2}\)-orthogonal to \(V_{U}\).  Lemma \ref{lem:B.11}
and (\ref{eq:(B.46)}) imply that
\begin{equation}\label{eq:(B.49)}
\int_U\big( |\nabla_A\grq|^2+\langle \grq, \mathcal{Q}^+\grq\rangle\big) \geq
\big((c_{0}^{-1}\c^{2}   -c_{0}) m_{k}^{-1} r +2c_{\diamond} \sup_{U}
|\mathcal{Q}^{-}|\big)\int_U |\grq|^2,
\end{equation}
if \(q\) is \(L^{2}\) -orthogonal to \(V_{U}\).  Thus
(\ref{eq:(B.40)}) holds if \(c \geq c_{0} (1 + m_{k}r^{-1} T^{2})\); and in particular,
(\ref{eq:(B.40)}) holds for \(c > c_{0}\) if the eigenvalue
bound \(T\) is less than \(r^{1/6}(\ln r)^{-c_0}\).
 
The last case to consider is that where \(U\) comes from
\(\grU_{N}\).  Lemma \ref{lem:B.11} and (\ref{eq:(B.47)}) imply for such \(U\) that
\begin{equation}\label{eq:(B.50)}
\int_U\big( |\nabla_A\grq|^2+\langle \grq, \mathcal{Q}^+\grq\rangle\big)
 \geq \big((c_{0}^{-1}  (\ln r)^{2\c}-(\ln r)^{c_0}r^{2/3} + 2c_{\diamond}
\sup_{U} |\mathcal{Q}^{-}|\big) \int_U |\grq|^2
\end{equation}
if \(\grq\) is \(L^{2}\)-orthogonal to \(V\).  It follows as a
consequence that (\ref{eq:(B.40)}) holds for such \(U\) if both \(\c >
c_{0}\) and the eigenvalue bound \(T\) is less than \(r^{1/3}\).   
 
Granted all of the above, and given that \(T <
r^{1/6} (\ln r)^{-\c}  \), then (\ref{eq:(B.40)}) holds for all sets from \(\grU\) if \(\c >
c_{0}\).  This understood, choose \(\c\) to be twice this lower bound.

\paragraph{Part 10:}  The dimension of each \(U \in \grU\) version of
\(V_{U}\) is bounded by \(c_{0}\), and so it follows
from what is said at the end of Part 5 that \(\grn_{T}\) is no
greater than \(c_{0}\) times the number of sets in the
collection \(\grU\).  
 
An upper bound for size of \(\grU\) is obtained by summing upper bounds for the
sizes of the various \(k \in \{1, \ldots, N\}\) versions of
\(\grU_{k}\).  Let \(N_{T}\) denote the largest value
of \(k\) such that \(T > m_{k}\) and suppose first
that \(k \in \{1, \ldots, N_{T}\}\).  It follows from
what is said in Part 6 that \(\grU_{k-}\) contains no more than
\(c_{0} T^{3}m_{k}^{-3}\) sets.  Meanwhile,
\(\grU_{k0}\) and \(\grU_{k+}\) together contain at most
\(c_{0} m_{k}^{-3/2}r^{1/2} (\ln r)^{c_0}\) balls.  Thus 
 \(\bigcup_{1\leq k\leq N_T}\grU_{k}\) contains at most \(c_{0}
(T^{3} + r^{1/2} (\ln r)^{c_0})\) balls.  Suppose next that \(k \in \{N_{T}+1, \ldots,
N  - 1\}\).  In this case, \(\grU_{k-}\) has at most \(c_{0}\) balls while \(\grU_{k0}\) and
\(\grU_{k+}\) again have at most \(c_{0}m_{k}^{-3/2} r^{1/2} (\ln r)^{c_0} \)
  balls.  Thus, \(\bigcup_{N_T<k\leq N-1}\grU_{k}\) contains at most \(c_{0}
T^{-3/2}r^{1/2} (\ln r)^{c_0}\) balls.  As noted in Part 6, the set \(\grU_{N}\) has at most
\(c_{0}(\ln r)^{c_0}\) balls.  

 Given that \(T \leq c_{0} r^{1/6}(\ln r)^{c_0}\), the bounds just stated imply that \(\grn_{T} \leq
c_{0} r^{1/2}(\ln r)^{c_0}\).  Therefore, (\ref{eq:(B.37)}) bounds the spectral flow along the
families \(\{\mathcal{L}_{2,s}\}_{s\in[0,1]}\) and \(\{\mathcal{L}_{3,s}\}_{s\in[0,1]}\)
by \(c_0 T^{-1} r (\ln r)^{c_0}\).  This understood, take \(T = r^{1/7} (\ln r)^{c_0}\)
  to obtain the bound in (\ref{eq:(B.33)}). \epf

 \subsection{The proof of Proposition \ref{prop:A.4}}\label{sec:Bh)}

 If \(Y_{Z}\) has a single component, then the function
\(\grf_{s}\) is defined in Section \ref{sec:Bf)}.  Proposition
\ref{prop:A.4}'s assertion in this case is implied directly by
Lemma \ref{lem:B.9}'s fourth bullet.

Suppose now that \(Y_{Z}\) has more than 1 component.  To
define \(\grf_{s}\) in this case, introduce \(\mathcal{Y }\) to
denote the set of components of \(Y_{Z}\).  The space \(\op{Conn} (E)
\times C^{\infty}(Y_{Z};\bbS)\) can be written as
\(\prod_{Y'\in\mathcal{Y} }(\op{Conn} (E|_{Y'})\times
C^{\infty}(Y';\bbS|_{Y'}))\).  Section \ref{sec:Bf)} defines any given \(Y' \in \mathcal{Y} \) version of \(\grf_{s}\) on
\(\op{Conn} (E|_{Y'})\times C^{\infty}(Y';\bbS|_{Y'})\).  Denote the latter by
\(\grf_{s ;Y'}\).  Set \[\grf_{s} =\sum_{Y'\in\mathcal{Y} }\grf_{s ;Y'}.\]  

 Each \(Y'\in \mathcal{Y}\)  has its version of the function \(\gra\) on
\(\op{Conn} (E|_{\mathcal{Y}'})\times C^{\infty}(Y';\bbS|_{Y'})\).  Use \(\gra_{Y'}\) to
denote the latter.  Then \(\gra^{\grf} =
\sum_{Y'} (\gra_{Y'} + 2\pi   (r  - \pi) \grf_{s ;Y'})\).  This understood, it is
enough to bound \(|\gra_{Y'} + 2\pi  (r  - \pi) \grf_{s ;Y'}|\) for each \(Y' \in
\mathcal{Y}\) .   Lemma \ref{lem:B.9} supplies a suitable bound when
\(c_{1}(\det (\bbS|_{Y'}))\) is not torsion.  This understood, suppose \(Y'\in \mathcal{Y}\)  and
\(c_{1}(\det (\bbS)|_{Y'})\) is torsion.  Thus, \(w = 0\) on \(Y'\).  

Write \(\psi\) on \(Y'\) as \(r^{ -1/2 }\lambda \) to see
that the set of solutions to (\ref{eq:(A.4)}) on \(Y'\) is \(r\)-independent.  It follows
as a consequence of what is said in Chapter 5 of \cite{KM} that the space of
\(C^{\infty}(Y'; S^{1})\)-orbits of
solutions to (\ref{eq:(A.4)}) on \(Y'\)  is compact.  Hold on to this fact for the
moment.  Write \(\psi\) in the \(Y'\) version of (\ref{eq:(B.30)}) as
\(r^{1/2} \lambda \) and write the sections \(b\) and
\(\phi\) as \((r z)^{1/2} b'\) and \((rz)^{1/2} \phi'\) to see that the spectrum of the
operator in (\ref{eq:(B.30)}) depends neither on \(r\) nor \(z\).  What was just said
about compactness and what was just said about the spectrum implies
directly that \(|\gra_{Y'} + 2\pi   (r  - \pi )\grf_{s ;Y'}| \leq c_{0}\).

 \section{Cobordisms and the Seiberg-Witten equations}\label{sec:C}
\setcounter{equation}{0}
 
This section proves Propositions \ref{prop:A.3} and \ref{prop:A.5}.  
Here is a brief of what is to come.

Section  \ref{sec:Ca)}: This section states three key lemmas (Lemmas  \ref{lem:C.1}-\ref{lem:C.3})
that are used in Section \ref{sec:Cb)} to prove Proposition
\ref{prop:A.3}. These are used
subsequently also. These lemmas establish a priori estimates on the
norms of \(\psi\), \(\nabla_A\psi\) and the curvature \(F_A\) when \((A,\psi)\) is an instanton
solution to (\ref{eq:(A.14)}) and \(r\) is large. Lemmas \ref{lem:C.1}-\ref{lem:C.3}  are proved in
subsequent subsections of Section \ref{sec:C}.

Section \ref{sec:Cb)}: This section uses the lemmas in Section \ref{sec:Ca)} to prove
Proposition \ref{prop:A.3}.

Section  \ref{sec:Cc)}: This section ties up a loose end by giving the
proof of Lemma \ref{lem:C.1} from Section \ref{sec:Ca)}.

Section \ref{sec:Cd)}: This section ties up a loose end by giving the proof of
Lemma \ref{lem:C.2}  from Section \ref{sec:Ca)}.

Section  \ref{sec:Ce)}: This section gives the proof of Lemma \ref{lem:C.3}  
from Section \ref{sec:Ca)} modulo Lemma \ref{lem:C.5}  which is an assertion about the behavior of
\(\psi \) on
certain domains in a cobordism.

Section \ref{sec:Cf)}: This section proves Lemma \ref{lem:C.5}.

Section  \ref{sec:Cg)}: This section uses the results in the
previous sections of Section \ref{sec:C} and the results from Section \ref{sec:B} to prove
Proposition \ref{prop:A.5}.


 \subsection{The three key lemmas}\label{sec:Ca)}

The three parts of this subsection supply three lemmas that assert
pointwise bounds for \(\psi\), the curvature of \(A\) and for the covariant
derivative of \(\psi\). These bounds are used in the next
subsection to prove Proposition \ref{prop:A.3}.   All three lemmas assume
implicitly that the conditions in Section \ref{sec:Ac)} are satisfied.
 Additional assumptions are stated when needed.
 
 \paragraph{Part 1:}  The first lemma starts the story with a pointwise
bound for \(|\psi|\) and \(L^{2}\)-bounds
on \(F_{A}\) and the covariant derivatives of \(\psi\).  With
regards to notation, this lemma uses \((\nabla_{\bbA}\psi)_{s}\) to
denote the section of \(\bbS^{+}\) over the \(|s | > 1\) part of \(X\) that gives the
pairing between \(\nabla_{\bbA}\psi \) and the vector field
\(\frac{\partial}{\partial s}\).  
 
\begin{lemma}\label{lem:C.1}  
There exists \(\kappa >1\) such that given any \(\c \geq \kappa \), there exists
\(\kappa_{\c}\) with the following significance: Fix \(r \geq
\kappa_{\c}\).  If \(X\) is not the product cobordism, assume that
the metric obeys (\ref{eq:(A.12,15a)}) with \(L \leq \c\), that the norm of
the Riemann curvature is bounded by \(r^{1/\c}\) and that the norm of \(w_{X}\)  is bounded by
\(\c\).   Fix \(\mu_{-}\)  and \(\mu_{+}\) from the \(Y_{-}\)  and \(Y_{+}\)  versions
of \(\Omega\) with \(\mathcal{P}\)-norm bounded by 1
 and use this data to define the equations in (\ref{eq:(A.14)}).  Suppose
that \(\grd = (\bbA, \psi)\)  is an instanton solution to these equations.  Then \(|\psi| \leq
\kappa_{\c}\).   If \(X\) is not the product cobordism, assume in addition that the
volume of the \(s\)-inverse image of any length 1 interval is
bounded by \(\c\) and that the metric's injectivity radius is greater
than \(r^{-1/\c}\).  Also assume in this case that \(L_{tor}\leq \c
r\) and that \(w_X\) obeys (\ref{eq:(A.13b)})
plus Item c) of the fourth bullet of (\ref{eq:(A.16)}). Let \(\grc_{-}\) and \(\grc_{+}\)  denote the respective \(s \to-\infty  \) and \(s \to \infty  \) limits of \(\grd\)  and suppose that \(\gra(\grc_{-}) -\gra(\grc_{+}) \leq \c r^{2}\).  Then
\begin{itemize}
\item The \(L^{2}\)-norms of
 \(|F_{\bbA}( \frac{\partial}{\partial s}, \cdot) |\)  and  \(r^{1/2}|(\nabla_{\bbA}\psi)_{s}| \) on the \(|s | \geq L\)  part of
\(X\)  are less than \(\kappa_{\c}r\). 
\item The \(L^{2}\)-norms of \(F_{\bbA}\)  and \(r^{1/2}\nabla_{\bbA}\psi  \) on the
\(s\)-inverse image of any length 1 interval in \(\bbR\)
 are no greater than \(\kappa_{\c}r\). 
\end{itemize}
\end{lemma}
This lemma is proved in Section \ref{sec:Cc)}.
 
\paragraph{Part 2:}  The next lemma supplies a refined set of bounds for
\(|\alpha|\) and its covariant derivatives on \(U_{C}\) and \(U_{0}\).  This lemma and the
subsequent lemma implicitly write \(\bbS^{+}\) on
\(U_{C}\) and \(U_{0}\) as \(E \oplus (E\otimes K^{-1})\).  Having
done so, they then write \(\psi\) with respect to this splitting as
\((\alpha, \beta)\); and they write the connection \(\bbA\) as
\(\bbA = A_{K} + 2A\) with \(A\) being a connection on \(E\).

The notation in these upcoming lemmas refers to the complex structure on
\(U_{C}\) and \(U_{0}\) that is defined using the
metric and the compatible symplectic form ds \(\wedge *w + w\).
  The \((1, 0)\)-part of the complexified cotangent space for this
complex structure is the direct sum of the span of \(ds + i*w\) and \(dz\) on
\(U_{C}\) and it is the direct sum of the span of \(ds + i*w\)
and the \((1, 0)\)-part of the tangent space to the constant-\((s,
u)\) spheres in \(U_{0}\) with the complex structure on \(S^{2}\)
being the standard one.  These lemmas write
\(\nabla_{A}\alpha\) with respect to the \((1, 0)\)- and
\((0, 1)\)-splitting of the complexified cotangent bundle as
 \(\partial _{A}\alpha +\bar{\partial }_{A}\alpha \) with  \(\partial _{A}\alpha \) denoting \((1, 0)\)-part of
\(\nabla_{A}\alpha \) and with \(\bar{\partial }_{A}\alpha \) denoting the \((0, 1)\)-part.  The lemma and the subsequent also
introduce \(\rho_{D}\) to denote the diameter of the
cross-sectional disk \(D\) that is used to define \(U_{C}\).
 
\begin{lemma}\label{lem:C.2}   
There exists \(\kappa  >100 \,( 1 + \rho_{D}^{-1})\) such that
given any \(\c \geq \kappa \), there exists
\(\kappa_{\c} \geq \kappa  \) with the following significance:  Fix \(r \geq
\kappa_{\c}\)  and assume that the metric obeys
(\ref{eq:(A.12,15a)}), (\ref{eq:(A.15b)}), and the \((\c, \r = r)\)-versions of the conditions in the first two bullets of (\ref{eq:(A.16)}).
 Assume that \(|w_{X}| \leq\c\)  and that \(w_{X}\)  obeys (\ref{eq:(A.13c)}).  Fix elements \(\mu_{-}\)  and
\(\mu_{+}\)  from the \(Y_{-}\)  and \(Y_{+}\)  versions
of \(\Omega\) with  \(\mathcal{P}\)-norm bounded by 1.
 Assume in addition that their norms and those of their derivatives to
order 10 on \(U_{\gamma}\)  and \(\mathcal{H}_0\) are bounded by
\(e^{-r^2}\).   Use this data to define the equations in
(\ref{eq:(A.14)}).  Let \(\grc_{-}\)  and \(\grc_{+}\)  denote
respective solutions to the \((r,  \mu_{-})\)-version of (\ref{eq:(A.4)}) on \(Y_{-}\)  and the \((r, \mu_{+})\)-version of (\ref{eq:(A.4)}) with
\(\gra(\grc_{-}) - \gra(\grc_{+}) \leq \c r^{2}\), and suppose that \(\grd = (\bbA,
\psi)\)  is an instanton solution to (\ref{eq:(A.14)}) with \(s\to -\infty  \) limit equal to
\(\grc_{-}\)  and \(s \to \infty\) limit equal to \(\grc_{+}\).  If \(p\)  is
a point in one of the domains \(U_{C}\)  or \(U_{0}\) with distance greater than
\(\kappa^{2} r^{-1/2}(\ln r)^{2}\)  from the domain's boundary, then the following holds at \(p\):
\begin{itemize}
\item \(|\beta|^{2} \leq e^{-\sqrt{r}/\kappa^2}\)
   and    \(|\alpha|^{2}\leq 1 +e^{-\sqrt{r}/\kappa^2}\).
\item \(|\nabla_{A}\beta|  + |\nabla_{A}\nabla_{A}\beta|\leq   e^{-\sqrt{r}/\kappa^2}\).         
\item \(|\bar{\partial  }_A\alpha | \leq  e^{-\sqrt{r}/\kappa^2}\).
\item If \(|\alpha|^{2} \in(\kappa^{-1}, 1 -\kappa^{-1})\)  at \(p\),
  then either \(|\nabla_{A}\alpha|^{2}\geq \kappa^{-3} r\) at \(p\)
  or the Hessian \(\nabla d|\alpha|^{2} \)
at \(p\)  has an eigenvalue with absolute value greater than \(\kappa^{-3} r\).
\item \(|\nabla_{A}\alpha| +r^{-1/2}|\nabla_{A}(\nabla_{A}\alpha)|\leq
  \kappa r^{1/2}\) if \(|F_{A}| \leq \c r\) on the radius \(\kappa  r^{-1/2}\)-ball centered at \(p\).
\end{itemize}
\end{lemma}
This lemma is proved in Section \ref{sec:Cd)}.

\paragraph{Part 3:}  The final lemma here writes \(F_{A}\) on
\(U_{C}\) and \(U_{0}\) as \(F_{A} = ds \wedge \mathcal{E}_{A} + *\mathcal{B}_{A}\) with
\(\mathcal{E}_{A}\) and \(\mathcal{B}_{A}\) denoting \(s\)-dependent, 
\(i\mathbb{R}\)  valued 1-forms on either \(\bbR/(\ell _\gamma\bbZ) \times D\) or $\mathcal{H}_0$ as the
case may be.  These 1-forms are written as 
\begin{equation}\label{eq:(C.1)}\begin{cases}
\mathcal{E}_{A} = -i (1  - \sigma) (r (1-|\alpha|^{2}) +\grz_{A})  \, dt+ \grr + \grX  &    \text{and}\\
\mathcal{B}_{A} = -i \sigma (r (1 -|\alpha|^{2}) +\grz_{B})  \, dt + \grr - \grX, &
\end{cases}
\end{equation}
where \(\sigma \), \(\grz_{A}\) and \(\grz_{B}\) are functions, and where
both \(\grr\) and \(\grX\) annihilate the vector field
\(\frac{\partial}{\partial t}\). Note that \(\mathcal{E}_{A} +
\mathcal{B}_{A}= -i\, (r (1-|\alpha|^{2}) +\grz_{A}+\grz_{B})  \, ds+ 2 \grr \) which means that \(\grr \) and the combination 
\(\grz_{A}+\grz_{B}\) contain the terms with \(\beta\) that appear in the left most equation of (\ref{eq:(A.14)}).

\begin{lemma}\label{lem:C.3}   
There exists $\kappa >\pi $  such that given any  \(\c \geq \kappa \), there
exists \(\kappa_{\c}  > 200 (1 +\rho_{D}^{-1})\) with the
following significance:  Fix \(r \geq\kappa_{\c}\)  and assume that the metric and
\(w_{X}\)  are  \((\c, \r = r)\)-compatible.  Fix
elements \(\mu_{-}\)  and \(\mu_{+}\)  from the \(Y_{-}\)-  and \(Y_{+}\)-versions
of \(\Omega\) with \(\mathcal{P}\)-norm bounded by 1.
 Assume in addition that their norms and those of their derivatives up
 to order 10 on \(U_{\gamma}\) and \(\mathcal{H}_0 \) are bounded by
 \(e^{-r^2}\).   Use all of these data to define the equations in (\ref{eq:(A.14)}).   Let \(\grc_{-}\)  and \(\grc_{+}\)
 denote the respective solutions to the \((r, \mu_{-})\)-version of (\ref{eq:(A.4)}) on
\(Y_{-}\)  and the \((r, \mu_{+})\)  version of (\ref{eq:(A.4)}) on \(Y_{+}\)  with \(\gra(\grc_{-})  -
\gra(\grc_{+}) \leq r ^{2-1/\c}\).  Suppose that \(\grd = (\bbA, \psi)\) is an instanton solution to (\ref{eq:(A.14)}) with \(s \to
-\infty  \) limit equal to \(\grc_{-}\)  and \(s\to \infty  \) limit
equal to \(\grc_{+}\).  Let \(p\)  denote a point in either one of the
domains \(U_{C}\) or \(U_{0}\) with distance \(\kappa^{-1}\) or more from the
domain's boundary.  Then the following are
true at \(p\): 
\begin{itemize}
\item \(- r ^{-100} < 1 - \sigma < 1+ r ^{-100}\).
\item \(|\grz_{A}| +|\grz_{B}| \leq r^{-100} \).
\item \(|\grr| \leq \kappa r^{-100}\) .
\item \(|\grX|^{2} \leq 2r^{2}\sigma (1  - \sigma) (1 -|\alpha|^{2}) + \kappa r^{-100}\).
\item \(|\nabla \mathcal{E}_{A}| +|\nabla\mathcal{B}_{A}| \leq \kappa r^{3/2}\).
\end{itemize}
\end{lemma}

Lemma \ref{lem:C.3} is proved in Section \ref{sec:Ce)} modulo a key lemma which is proved in
Section \ref{sec:Cf)}.

 \subsection{Proof of Proposition \ref{prop:A.3}}\label{sec:Cb)}
 
This part of the subsection uses what is said in Lemmas \ref{lem:C.1}-\ref{lem:C.3} to
prove Proposition \ref{prop:A.3}.  The argument assumes that the integral of \(i
F_{\hat{A}}\) over \(C\) is negative so as to derive nonsense.  This
is done in the eight parts that follow.  Before starting, note that
the assumptions in this proposition allow Lemmas \ref{lem:C.1} and \ref{lem:C.3} to be
invoked, and the conclusions of Lemma \ref{lem:C.3} imply in particular that
Lemma \ref{lem:C.2} can be invoked as well.
 
 \paragraph{Part 1:}  This first part of the proof sets the stage for what
is to come by supplying two observations about the pull-back of
\(iF_{\hat{A}}\) to \(C\).  What follows is the first observation:
\begin{equation}\label{eq:(C.2)}
 \text{The integral of \(\frac{i}{2\pi}F_{\hat{A}}\) over \(C\)  is an integer.}
\end{equation}
 
This follows from Lemma \ref{lem:B.6} since the latter implies that \(\hat{A}\) is flat and
 \(\alpha/|\alpha|\) is \(\hat{A}\)-covariantly constant
where \(|s | \gg 1\) on \(C\). 
 
The second observation concerns the function \(\textsc{f}\) on \(C\) that is
defined by writing the pull-back to \(C\) of \(i F_{\hat{A}}\) as \(\textsc{f }ds \wedge dt\):  
\begin{equation}\label{eq:(C.3)}
\text{ The function \(\textsc{f}\) is nearly non-negative in the
sense that \(\textsc{f} \geq -c_{0}r^{-100 }\).}
\end{equation}
 
This follows directly from the formula given below for \(\textsc{f}\) using
the second bullet of Lemma \ref{lem:C.2} and the first and second bullets of
Lemma \ref{lem:C.3}.  The upcoming formula for \(\textsc{f}\) uses
\((\partial _{A}\alpha)_{0}\) to denote
the \(ds + i*w\) component of  \(\partial _{A}\alpha \) and
use \((\bar{\partial }_{A}\alpha)_{0}\) to denote the \(ds - i *w\) component of \(\bar{\partial }_{A}\alpha\). 
Here is the promised formula for \(\textsc{f}\):
\begin{equation}\label{eq:(C.4)}
\textsc{f} = (1  - \wp) \,( 1  -  \sigma) \,( r (1 -|\alpha|^{2})  +
\grz_{A}) + \wp'(|( \partial _{A}\alpha)_{0}|^{2}-  |(\bar{\partial }_{A}\alpha)_{0}|^{2}).
\end{equation}

This formula follows directly from (\ref{(3.9)v2}) 
and (\ref{eq:(C.1)}).

 \paragraph{Part 2:}  Let \(\mathbb{I}  \subset \bbR \) denote
the set characterized as follows:  A point \(s\) is in \(\mathbb{I }\) if the
the integral of \(\textsc{f}\) over the slice \(\{s \} \times \gamma\)
in \(C\) is negative.  The following assertion is a direct consequence of
(\ref{eq:(C.2)}) and (\ref{eq:(C.3)}):  
\begin{equation}\label{eq:(C.5)}
 \text{If  \(\int_CiF_{\hat{A}}< 0\) then the measure of the set  \(\mathbb{I}\)  
is greater than \(c_{0}^{-1}r^{100}\).  }
\end{equation}
 
Granted (\ref{eq:(C.5)}), there are at least \(c_{0}^{-1} r^{100}\)
disjoint open intervals of length 1 in \(\bbR\) with center point in
\(\mathbb{I}\) .  This understood, use the first bullet of Lemma
\ref{lem:C.1} to
find an interval \(I \subset \bbR \) of length 1 with center
point in \(\mathbb{I }\), with \(|s | > L + 2\) and such that
\begin{equation}\label{eq:(C.6)}
\int_{I\times Y}\big( |F_A (\frac{\partial}{\partial s},
\cdot)|^2+r|(\nabla_A\psi)_s|^2\big) <  r ^{-97}.  
\end{equation}
 This inequality enters the story in Parts 3 and 7.

\paragraph{Part 3:}  
Supposing that \(I \subset \bbR \) is given by Part 2, let \(s\) denote its center point, this being a
point for which the integral of \(\textsc{f}\) over \(\{s\}\times  \gamma\)
is negative. This part proves that \(|\alpha|^{2} \leq
\frac{5}{8}\)  
on \(\{s\}\times \gamma\). To see why this is, suppose for the sake of argument that this condition is
violated at \(p \in \{s\} \times \gamma \). Since the integral of
\(\textsc{f}\) on \(\{s\}\times \gamma \) is negative, there must be some
point where the function \(\wp\) is less than 1 and thus
\(|\alpha|^2\leq \frac{9}{16}\). As a consequence, the 
variation of \(|\alpha |\) on \(\{s\}\times\gamma \) must be greater than \(c_0^{-1}\). As explained next, this variation is in 
fact no greater than \(c_0 r^{-15}\) if \(r\geq c_0^{-1}\). To start
the explanation, suppose that \(\varepsilon >0\) and that 
there are points on \(\{s\}\times\gamma \) with their respective
values of \(|\alpha |\) differing by more that \(\varepsilon\). Let
\(\frac{\partial}{\partial t}\) denote the unit length tangent vector
to \(\{s\}\times\gamma \)  and let \((\nabla_A\alpha )_t\)  denote the directional
covariant derivative of \(\alpha \) along \(\frac{\partial}{\partial t}\). It follows as a consequence of the fundamental 
theorem of calculus that there is a point  \(q\in \{s\}\times\gamma \)
where \(|(\nabla_A\alpha)_t|> c_0^{-1}\varepsilon \). Let \((\nabla_A\alpha)_s\)
denote the directional covariant derivative of \(\alpha \) along the vector field \(\frac{\partial}{\partial s}\). Granted this 
lower bound for \(|(\nabla_A\alpha)_t|\) at \(q\), then the
inequality in the third bullet of Lemma \ref{lem:C.2} requires that
\(|(\nabla_A\alpha)_s|\geq c_0^{-1}\varepsilon \) at \(q\) also if
\(\varepsilon \) is greater than \(c_0 e^{-\sqrt{r}/c_0}\). Assuming
\(r\geq c_0\), then this will
be the case when \(\varepsilon > r^{-15}\). The \(c_0^{-1}\varepsilon
\) lower bound for \(|(\nabla_A\alpha)_s|\) at \(q\), what is said by
the fifth bullet of Lemma  \ref{lem:C.2} and what is said by Lemma
\ref{lem:C.3} imply that \(|(\nabla_A\alpha)_s|\geq
c_0^{-1}\varepsilon\) in the ball in \(U_C\)
of radius \(c_0^{-1}\varepsilon  r^{-1}\) centered at \(q\). The
latter bound implies in turn that the integral of
\(|(\nabla_A\alpha)_s|^2\) on this same ball is greater than
\(c_0^{-1}\varepsilon  ^6r^{-4}\), which is nonsense if \(\varepsilon
> r^{-15}\)  because it
runs afoul of what is said in (\ref{eq:(C.6)}).

\paragraph{Part 4:} Let \(I\) and \(s\in I\) be as in Part 3. Keep in
mind that the integral of \(\textsc{f}\) over \(\{s\}\times\gamma \) is
negative. As will be explained momentarily, the lower bound in
(\ref{eq:(C.3)}) for \(\textsc{f}\) leads to the following observation:
\begin{equation}\label{eq:(C.7)}
  \text{The variation of \(\wp\) over \(\{s \} \times\gamma  \) is no greater than \(c_{0}r^{-50}\).}
\end{equation}
 To prove this, first use the fundamental theorem of
calculus to see that
\begin{equation}\label{eq:(C.8)}
\sup_{\{s\}\times \gamma  }  \wp - \inf_{\{s\}\times \gamma  } \wp
\leq c_{0} \Big( \int_{\{s\}\times \gamma  }\wp' |(\nabla_A\alpha)_t|^2\Big)^{1/2}.
\end{equation}
The bound in (\ref{eq:(C.7)}) follows from (\ref{eq:(C.8)}) using the lower bound for
\(\textsc{f}\) and the third bullet of Lemma \ref{lem:C.2}.  
 
 
 

\paragraph{Part 5:} This part uses the conclusions of Part 3 to
deduce the following:  
\begin{equation}\label{eq:(C.9)}
\text{The function \(\sigma\) on the \(|u|
< 1\)  part of  \(\{s \} \times (\gamma \cap\mathcal{H}_0) \) obeys \(\sigma <
c_{0} r^{-33}\).}
\end{equation}
 
To see why this is the case, let \((s, p)\) denote a given point in the
\(|u| < 1\) part of \(\{s \} \times(\gamma \cap \mathcal{H}_0)\) where \(\sigma
> 0\).  Let \(S\) denote the cross-sectional sphere in
$\mathcal{H}_0$ that contains \(p\).  Use (\ref{(3.9)v2}) to write
the pull-back of \(F_{A}\) to \(S\) as \(\frac{1}{2} \textsc{b} dz \wedge d \bar{z}\) 
with \(\textsc{b} = \sigma (r (1-|\alpha|^{2}) +\grz_{B})\).  Use \(\varepsilon\) to denote the value of
\(\sigma\) at \((s, p)\).  Invoke the first and second bullets of Lemma
\ref{lem:C.3} to conclude (using what is said in Part 3 to the effect
that \(|\alpha |^2\leq \frac{5}{8}\)) that value of \(\textsc{b}\) at \((s, p)\) is greater than  
\(\frac{3}{8} r \varepsilon - c_{0} r^{-100}\).  The
fifth bullet of Lemma \ref{lem:C.3} finds that \(\textsc{b} >
c_{0}^{-1} r \varepsilon \) on the radius \(c_{0}^{-1} r^{1/2}\varepsilon  \) disk in the cross-sectional sphere \(\{s \} \times S\)
with center at \((s, p)\).  Meanwhile, the first bullets of Lemma \ref{lem:C.3} and
Lemma \ref{lem:C.2} imply that \(\textsc{b} > -c_{0}r^{-99}\) on the whole of \(\{s \} \times S\), and so
the integral of \(\textsc{b}\) over \(\{s \} \times S\) is no less than
\(c_{0}^{-1}\varepsilon^{3} - c_{0}r^{-99}\).  This integral must be zero because the
first Chern class of \(E\) has zero pairing with the cross-sectional
spheres in \(\mathcal{H }_0\) .  Thus \(\varepsilon\leq c_{0} r^{-33}\).

 \paragraph{Part 6:} What is said in Part 5 implies that \((1-\wp) < c_{0} r^{-50}\) on \(\{s\} \times \gamma\).  Indeed, if this bound is violated, then it
follows from (\ref{eq:(C.7)}) and the formula for \(\textsc{f}\) in (\ref{eq:(C.4)}) that the
integral of \(\textsc{f}\) over the \(|u| < 1\) part of \(\{s \} \times \gamma \) is greater than \(c_{0}
r^{-49}\).  Given the lower bound in (\ref{eq:(C.3)}), this last
lower bound runs afoul of the assumption that
\(\textsc{f}\)'s integral over \(\{s \} \times \gamma\) is negative.  The small size of \(1 -\wp \) implies in particular that
\(|\alpha|^{2} > \frac{3}{8} \) on \(\{s \} \times \gamma\).  
 
 \paragraph{Part 7:}  Granted the conclusions of Parts 5 and 6, then
the fourth bullet of Lemma \ref{lem:C.2} asserts that one or the other of the
following are true at each point in the \(|u|< 1\) part of \(\{s \} \times (\gamma \cap\mathcal{H}_0)\):  Either
\(|\nabla_{A}\alpha|^{2}
\geq c_{0}^{-1}r\) or the Hessian
matrix \(\nabla d|\alpha|^{2}\) has
an eigenvalue with absolute value greater than
\(c_{0}^{-1} r\).  As explained next, this has the following consequence:
\begin{equation}\label{eq:(C.10)}\begin{split}
 &\textrm{Let \(( \partial _{A}\alpha)_{1}\)
 denote the component of
 \(\partial _{A}\alpha  \) that annihilates both
 $\frac{\partial}{\partial s}$}\\
&\textrm{and the kernel of \(w\). Then
\(|( \partial _{A}\alpha)_{1}|^{2}\)
 is greater than \(c_{0}^{-1} r^{1/2}\)  at all points }\\ 
&\textrm{in a radius \(c_{0}^{-1}r^{-1/2}\)
 ball with center at distance less than \(c_{0}r^{-1/2}\)  }\\
& \textrm{from each point in
the \(|u| < 1\)  part of \(\{s \}\times (\gamma  \cap   \mathcal{H}_0)\). }
\end{split}
\end{equation}

To prove this, suppose first that \(|\nabla_{A}\alpha|^{2}\geq c_{0}^{-1}r\) at a given point.
 Use the third bullet of Lemma \ref{lem:C.2} to see that one or both of
\(|( \partial _{A}\alpha)_{1}|^{2}\) and \(|( \partial _{A}\alpha)_{0}|^{2}\)
are greater than \(c_{0}^{-1}r\).  In the
latter case, the third bullet of Lemma \ref{lem:C.2} implies that
\(|(\nabla_{A}\alpha)_{s}|^{2}\) is greater than \(c_{0}^{-1}r\) at the point,
and the second derivative bound from the fifth bullet of Lemma \ref{lem:C.2}
implies that \(|(\nabla_{A}\alpha)_{s}|^{2}\geq c_{0}^{-1}r\) at all points in a
radius \(c_{0}^{-1}r^{-1/2}\) ball
centered on this point.  This being the case, the integral of
\(|(\nabla_{A}\alpha)_{s}|^{2}\) over this ball is greater than \(c_{0}^{-1}r^{-1}\) and this violates (\ref{eq:(C.6)}).  Granted that
\(|(\nabla_{A}\alpha)_{1}|^{2}\geq c_{0}^{-1} r\) at the given point,
then the second derivative bound from the fifth bullet of Lemma \ref{lem:C.2}
implies what is asserted by (\ref{eq:(C.10)}). 

 Now suppose that the Hessian matrix
\(\nabla d|\alpha|^{2}\) at the given point has an eigenvalue that is greater than
\(c_{0}^{-1} r\).  Let \(v\) denote a
unit length eigenvector at the point with such an eigenvalue.  As will
be 
explained directly, this vector must be such that
\(|ds(v)| + |dt(v)| < \frac{1}{100} \).  To see why this is the case, suppose to the contrary that the latter bound
is violated at a given point.  It then follows from the first and
fifth bullets of Lemma \ref{lem:C.2} that
\(|(\nabla_{A}\alpha)_{s}|\geq c_{0}^{-1}r^{1/2}\) at all points in some ball of radius
\(c_{0}^{-1}r^{-1/2}\) whose center has distance at most \(c_{0} r^{1/2}\)
from the given point.  This implies in particular that the integral of
\(|(\nabla_{A}\alpha)_{s}|^{2}\) over this same ball is no less than \(c_{0}
r^{-1}\).  But this is nonsense as it runs afoul of (\ref{eq:(C.6)}).  

   The fact that \(v\) is a unit length vector implies that
\(|dz(v)| >  \frac{1}{2} \).  Use this lower bound for \(|dz(v)|\) with the third
bullet of Lemma \ref{lem:C.2} and the second derivative bounds from the fifth
bullet of Lemma \ref{lem:C.2} to see that
\(|( \partial _{A}\alpha)_{1}|^{2}\geq c_{0}^{-1}r^{1/2}\) at all points in a ball of radius
\(c_{0}^{-1} r^{-1/2}\) whose center point has distance at most
\(c_{0} r^{-1/2}\) from the given point.

\paragraph{Part 8:}  Introduce the connection \(\hat{A}_{1}\) on
\(E\)'s restriction to 
$I \times\mathcal{H}_0$ that is obtained from \((A, \alpha)\) by
the formula \(\hat{A}_{1} = A - \frac{1}{2}(\bar{\alpha}\nabla_{A}\alpha -\alpha\nabla _{A} \bar{\alpha})\).  The curvature 2-form of \(\hat{A}_{1}\) is 
\begin{equation}\label{eq:(C.11)}
  F_{\hat{A}_{1}} = (1  - |\alpha|^{2})F_{A} + \nabla_{A}\alpha \wedge\nabla_{A}  \bar{\alpha}.
\end{equation}
Let \((s', p')\) denote the center point of a ball that is described by
(\ref{eq:(C.10)}).  Introduce $S \subset \mathcal{H}_0$ to
denote the cross-sectional sphere that contains the point \(p'\).  Use
(\ref{(3.9)v2}) to write the pull-back of the curvature of \(\hat{A}_{1}\) to
\(\{s \} \times S\) as \(\textsc{b}_ 1 dz \wedge d \bar{z}\)  with \(\textsc{b}_1\) given by
\begin{equation}\label{eq:(C.12)}
\textsc{b}_ 1 = \sigma (1 -|\alpha|^{2}) (r (1 -|\alpha|^{2})  +
\grz_{A}) +|( \partial _{A}\alpha)_{1}|^{2}-  |( \bar{\partial }_{A}\alpha )_{1}|^{2},
\end{equation}
with \(( \bar{\partial }_{A}\alpha )_{1}\) denoting here the
\(d\bar{z}\) component of \(\bar{\partial }_{A}\alpha \). The function \(\textsc{b}_1\) is also very
nearly non-negative in the sense that \(\textsc{b}_1\geq -c_{0} r^{-100}\), this
being a consequence of what is said in the first and third bullets of
Lemma \ref{lem:C.2} and the first and second bullets of Lemma \ref{lem:C.3}.  This
understood, then it follows from (\ref{eq:(C.10)}) and this lower bound for
\(\textsc{b}_1\) that the integral of
\(\textsc{b}_1\) over \(\{s' \} \times S\) is positive.
 But this is nonsense because the latter integral computes $2\pi $
times the pairing of the first Chern class of \(E\) with the homology class
defined by \(S\), and this pairing is equal to zero. \epf

\subsection{Proof of Lemma \ref{lem:C.1}}\label{sec:Cc)}

The bounds in the lemma constitute a particular case of bounds that
are used in Chapter 24 of \cite{KM}.  As most of the machinery behind what
is done in \cite{KM} is not needed for the proofs, the argument for Lemma
\ref{lem:C.1} is presented momentarily.  What follows directly lays a convention
that is invoked implicitly in the arguments for Lemma \ref{lem:C.1}
and in some of the subsequent lemmas.
 
If \(X\) is the product \(\bbR \times Y_{Z}\), the
bundles \(E\) and \(K^{-1}\) over \(Y_{Z}\) pull back
via the projection to define bundles over \(X\); their connections
\(A_{E}\) and \(A_{K}\) likewise pull back to define
connections on these bundles.  The bundle
\(\det (\bbS^{+})\) is isomorphic to \(E^2\otimes
K^{-1}\) and thus to the pull-back of \(\det (\bbS)\).
Fix once and for all an isometric isomorphism.  

Suppose now that \(X\) is not a product.  Use the embedding in the second
bullet of (\ref{(A.9a,11)}) to identify the \(s \leq -1\) part of \(X\) with
\((-\infty, -1] \times Y_{-}\), and then use the
projection to \(Y_{-}\) to view the \(Y_{-}\)
version of the bundle \(\bbS\) as bundles over the \(s \leq -1\)
part of \(X\).  The bundles \(\bbS^{+}\) and
\(\bbS^{-}\) are isometrically isomorphic to
\(\bbS\) via an isomorphism that covers the isomorphisms between
both \(\Lambda^{+}\) and
\(\Lambda^{-}\) and \(T^*Y\) given by the interior product
with \(\frac{\partial}{\partial s}\).   Fix such an isomorphism once and for all.  This induces a
Hermitian isomorphism between the bundle
\(\det (\bbS^{+})\) over the \(s < -1\) part of
\(X\) and the \(Y_{-}\) version \(\det (\bbS)\).  Fix once
and for all an isometric isomorphism between these bundles.  Use this
isomorphism with the pull back via the composition of the embedding
from (\ref{(A.9a,11)})'s second bullet and the projection to
\(Y_{-}\) to view \(A_{K} + 2A_{E}\) as a Hermitian connection on the \(s \leq -1\) part of
\(\det (\bbS^{+})\).  The analogous constructions can
be made on the \(s > 1\) part of \(X\) using the
\(Y_{+}\) version of \(\bbS\) and so define an incarnation
of the \(Y_{+}\) version of \(A_{K} +
2A_{E}\) as a Hermitian connection on
\(\det (\bbS^{+})\).

  Suppose for the moment that \(\bbA\) is a given Hermitian
connection on \(\det (\bbS^{+}) \to X\).
 If \(X\) is the product \(\bbR \times Y_{Z}\), then
\(\bbA\) can be written as \(A_{K} + 2A\) with \(A\) now a
connection on the bundle \(E \to X\).  There is a map \(\hat{u}: X
\to S^{1}\) such that \(A  -
\hat{u}^{-1}d\hat{u} = A_{E} + \hata _{A}\)
where \(\hata _{A}\) annihilates the vector field \(\frac{\partial}{\partial s}\).
If \(X\) is not the product, then \(\bbA\) can be written as
\(A_{K} + 2A\) on the \(s \leq -1\) and \(s \geq 1\) parts of
\(X\) with \(A\) being a connection on the incarnation of \(E\) over the relevant
part of \(X\)   In this case, there exists a map \(\hat{u}\) as just described but
with domain the \(s \leq -1\) part of \(X\), and likewise there exists such
a map with domain the \(s \geq 1\) part of \(X\).   

 The map \(\hat{u}\) in the case when \(X = \bbR \times
Y_{Z}\) is unique up to multiplication by an \(s\)-independent
map from \(Y_{Z}\) to \(S^{1}\), and in the other
cases, it is unique up to a map from the either the \(s \leq -1\) or \(s
\geq 1\) part of \(X\) whose differential annihilates \(\frac{\partial}{\partial s}\).  The convention in each case is to take a map \(\hat{u}\) whose restrictions
to the constant s slices of its domain are homotopic to the constant
map to \(S^{1}\). 

 The connection \(A_{*} = A  - \hat{u}^{-1}d\hat{u}\)
can be viewed as a map from \(\bbR\) or \((-\infty, -1]\) or \([1,
\infty)\) to \(\op{Conn} (E|_{Y_{*}})\) with \(Y_*\) either \(Y_{Z}\) or
\(Y_{-}\) or \(Y_{+}\) as the case may be.  If
\(\psi\) is a given section over \(X\) of \(\bbS\), then
\(\psi_{*} = \hat{u}\psi \) can likewise be viewed as a map
from \(\bbR\) or \((-\infty, -1]\) or \([1, \infty)\) to
\(C^{\infty}(Y_{*};\bbS|_{Y_{*}})\).  When viewed in this light, the
equations in (\ref{eq:(A.14)}) can be written as equations for
\((A_*, \psi_*)\) on the whole of \(X\) when \(X\) is the product cobordism, and on the
\(s  \leq   -L\) and \(s \geq L\) parts of \(X\) when \(X\) is not
the product cobordism. These equations are:
\begin{equation}\begin{cases}\label{eq:(C.13)}
 \frac{\partial}{\partial s}A_{*} + B_{A_*} -r \,( \psi_{*}^\dag \tau\psi_{*} - i*w_{X_*}) -\frac{1}{2}B_{A_K} 
 -i d\mu_{*} =0 & \text{and}\\
\frac{\partial}{\partial s} \psi_{*} + D_{A_*}\psi_* =   0. &
\end{cases}
\end{equation}
The notation here uses \(w_{X*}\)
to denote the 2-form \(w\) when \(X\) is the
product cobordism. When \(X\) is not the product cobordism,
\(w_{X*}\) denotes the \(s\)-dependent 2-form that is defined
on the relevant constant \(s\) slices of \(X\) by the pull-back of
\(w_{X}\). In particular, \(w_{X*}= w\) on the components of the
\(s \leq -L\) and \(s >  L\) parts of \(X\) where \(c_{1}(\det(\mathbb{S}))\)
is not torsion. What is denoted in (\ref{eq:(C.13)}) by \(\mu_{*}\) is
either \(\mu\), \(\mu_{-}\), or \(\mu_{+}\) as the case may be.

\pf The proof has four steps.

\paragraph{\it Step 1:}  The assertion that \(|\psi| \leq
\kappa \) is proved by using the Weitzenb\"ock formula to write
\(\mathcal{D}_{\bbA}^-\mathcal{D}_{\bbA}^+\psi\) as
\(\nabla_{\bbA}^{\dag}\nabla_{\bbA}\psi+ \op{cl}(\, F_{\bbA}^{+}) \psi +
\frac{1}{4}\textsc{r} \psi \) where \(\textsc{r}\) denotes the scalar curvature of
the Riemannian metric.  Granted this rewriting, it then follows from
(\ref{eq:(A.14)}) and from the assumed bound on the norm of Riemann
curvature that the function \(|\psi|\) obeys the differential
inequality \(d^{\dag}d|\psi| + r(|\psi|^{2} -|w_{X}| -c_{\c } )|\psi| \leq 0\).   Use the maximum
principle with this last inequality and the large \(|s|\)
bounds on \(|\psi|\) that follow from Lemma \ref{lem:B.1} to
see that \(|\psi| \leq \c + c_{0}\).
 
\paragraph{\it Step 2: } Let \(L_*\) denote either \(L\) or \(L_{tor}\). Then
use \(\mathbb{I}\subset \bbR \) to denote either \(\bbR\), \((-\infty, -L]\) or \([L, \infty)\). Define \(Y_{*}\) to be
\(Y_{Z}\) in the case when \(\mathbb{I} =\mathbb{R}\). When
\(\mathbb{I} = (-{\infty}, L_{*}]\) or \([L_{*}, {\infty})\) and \(L_{*} = L\), define
\(Y_{*}\) to be the union
of the components of the constant \(s\in\mathbb{I}\)  slices
of \(X\) where \(c_{1}(\det(\mathbb{S}))\) is not torsion. In the case when
\(L_{*} =L_{tor}\), define \(Y_{*}\) to be the union
of the components of the constant \(s\in \mathbb{I}\) slices of \(X\) where
\(c_{1}(\det(\mathbb{S}))\) is torsion.   
Write \(\bbA\) on \(\mathbb{I }\times Y_{*}\) as \(A_{K} + 2A \) and
introduce by way of notation \(\grd|_{s}\) to denote the
pull-back to \(\{s \} \times Y_{*}\) of \((A, \psi)\). Also introduce
\(\grB_{(A,\psi)}\) to denote
\begin{equation}\label{eq:(C.14)}
\grB_{(A,\psi)} =B_{A} - r (\psi^{\dag}\tau\psi
- i *w) + i *d\mu_{*} +  \frac{1}{2}B_{A_K},
\end{equation}
with \(\mu_{*}\) denoting either
\(\mu_{-}\) or \(\mu_{+}\) as the case
may be.  Use \(D_{A}\) in what follows to denote the Dirac
operator on \(Y_{*}\)  as defined using the connection
\(A_{K} + 2A\) for the \(\Spin^c\)-structure with spinor bundle \(\bbS =\bbS^{+}\).  Suppose that \(s' > s\) are
two points in \(\mathbb{I}\) .   Take the \(L^{2}\)-norm of
the left hand expressions in both equations of (\ref{eq:(C.13)}) over \([s, s']
\times Y_{*}\).  The square of these norms are zero.
 This being the case, integration by parts in the expressions the
square of these \(L^{2}\)-norms results in an identity of
the form
\begin{equation}\label{eq:(C.15)}
\frac{1}{2}\int_{[s, s']\times     Y_*}\big(|\frac{\partial}{\partial
  s}A_*|^2+|\grB_{(A, \psi)}|^2+2r\big(|\frac{\partial}{\partial s}\psi_*|^2+|D_A\psi|^2\big)\big)=\gra(\grd|_{s}) -
\gra(\grd|_{s '}) .
\end{equation}
Taking limits in (\ref{eq:(C.15)}) as \(s \to -\infty\) or as \(s' \to
\infty \) as the case may be leads to the identities

\begin{equation}\begin{cases}\label{eq:(C.16)}
\frac{1}{2}\int_{\mathbb{I}\times Y_*}(|\frac{\partial}{\partial
  s}A_*|^2+|\grB_{(A, \psi)}|^2+2r(|\frac{\partial}{\partial s}\psi_*|^2+|D_A\psi|^2))= \gra(\grc_{-}) - \gra(\grd|_{s=-L_*}).&\\
\frac{1}{2}\int_{\mathbb{I}\times Y_*}(|\frac{\partial}{\partial
  s}A_*|^2+|\grB_{(A, \psi)}|^2+2r(|\frac{\partial}{\partial s}\psi_*|^2+|D_A\psi|^2))
   = \gra(\grd|_{s=L_*} ) - \gra(\grc_{+}). &
\end{cases}
\end{equation}
 Note that the identities in (\ref{eq:(C.15)}) and (\ref{eq:(C.16)}) hold with \(\grd =
(A_{K} + 2A, \psi)\) on the right hand side. By way of an explanation, the integration by
parts proves the analogs that have \(\grd_{*} = (A_{K}+ 2A_{*}, \psi_{*})\) used on the right
hand side, and if they hold using \(\grd_{*}\), then they hold using \(\grd\)
because the restriction of the map \(\hat{u}\) to any slice \(\{s \}
\times Y_{*}\) in \(\bbI\times Y_*\) is homotopic to the constant map to \(S^{1}\). 
 
\paragraph{\it Step 3:}  The assertion made by the first bullet of Lemma \ref{lem:C.1}
follows directly from (\ref{eq:(C.16)}) when the data is such that \(X\) is the
product cobordism.  The proof in the general case and the proof of the
second bullet of Lemma \ref{lem:C.1} use an integral version of the
Weitzenb\"ock formula for the operator
\(\mathcal{D}_{\bbA}^-\mathcal{D}_{\bbA}^+\).  The details follow directly.

Integrate \(|F_{\bbA}^{+} - r(\psi^{\dag}\tau\psi  - \frac{i}{2}w_{X}) + i\grw_{\mu}^{+}|^{2}+ 2 r |\mathcal{D}_{\bbA}^+\psi|^{2}\)
over \(s^{-1}([-L  - 3, L + 3])\) and denote the result
by \(\mathcal{I}\).  Integrate this same expression over the respective
\(|s| \in[L,L_{*}- 4]\) and \([L_{*}-5, L_{*}+ 1]\) parts of \(X_{tor}\). Denote these
integrals as \(\mathcal{I}_{tor0}\) and \(\mathcal{I}_{tor1}\). In each case, let \(X_{*}\) denote the
region of integration and let \(\partial_{-}X_{*}\) and \(\partial_{+}X_{*}\)
denote the two boundaries of the relevant region of integration with
\(\partial_{-}X_{*}\) at the smaller value of \(s\) and \(\partial_{+}X_{*}\)
at the larger value. Use the Weitzenb\"ock formula for \(\mathcal{D}_{\bbA}^-\mathcal{D}_{\bbA}^+\psi\)
from Step 1 with Stokes' theorem to rewrite the identities \(\mathcal{I} = 0\),
\(\mathcal{I}_{tor0} = 0\) and \(\mathcal{I}_{tor1}= 0\) respectively as
\begin{equation}\label{eq:(C.17)}
 \frac{1}{2}\int_{X_*}( |F_{\mathbb{A}}|^2+r^2|\psi ^\dag\tau\psi-iw_X|^2+2r|\nabla_{\mathbb{A}}\psi|^2)+ \gri _*= \gra(\grd|_{\p_-X_*}) - \gra(\grd  |_{\p_+X_*}),
\end{equation}
with \(\gri_*\) in the case of \(\mathcal{I}\) and
\(\mathcal{I}_{tor1}\) denoting a term with absolute value no greater
than \(c_{0} \c r ( \int_{s^{-1}([-L  - 3, L + 3])}
(|F_{\mathbb{A}}|^2)^{1/2} + c_{0} \r^{1+c_0/\c}\).  In the case of
\(\mathcal{I}_{tor0}\), the absolute value of \(\gri_*\) is no greater
than \(c_0\c L_*r\). This bound on \(|\gri_*|\) in the case of \(\mathcal{I}\) and
\(\mathcal{I}_{tor0}\) is a direct consequence of
the bounds on the norms of the Riemannian curvature tensor and
\(w_{X}\), the size of \(L\), the volume of the \(s\)-inverse image
of intervals, and the bound \(|\psi|^{2} \leq 2 \c\) from Step 1.  
In the case of \(\mathcal{I}_{tor0}\), the bound for \(|\gri_{*}|\) is a consequence of the fact that
\(dw_{X}=0\) on the integration domain, this being the assumption made by Item
c) of the fourth bullet of (\ref{eq:(A.16)}). By way of an explanation, \(\gri_{*}\)
in this case can be written as sum of three terms, these denoted by
\(\gri_{\grg}\), \(\gri_{\mu}\) and \(\gri_{K}\). The term that is denoted by \(\gri_{\grg}\)
gives the contribution of the scalar curvature term in the Weitzenb\"ock
formula for \(\mathcal{D}_{\bbA}^-\mathcal{D}_{\mathbb{A}}^+\). As such, it is bounded by the integral of
\(c_{0} r |\psi|^{2}\) over the \(|s| \in[L, L_{*}-4]\) part of \(X_{tor}\). The
bound \(|\psi|^{2}\leq c_{0}\c\) leads to a bound on \(|\gri_{\grg}|\) by \(c_{0}\c rL_{*}\).

The term that is denoted by \(\gri_{\mu}\) comes by writing
\(|F_{\mathbb{A}}^{+}-r (\psi^{\dag}\tau\psi-\frac{i}{2}w_{X}) +i\grw_{\mu}^{+}|^{2}\) as the sum of
\(|F_{\mathbb{A}}^{+}-r (\psi^{\dag}\tau\psi-\frac{i}{2}w_{X})|^{2}\) with terms that
involve \(w_{\mu}^{+}\). One of these terms has the inner product between
\(F_{\mathbb{A}}^{+}\) and \(\grw_{\mu}^{+}\). Stokes' theorem identifies the integral of the latter
with the contributions to the boundary terms on the right hand side of
(\ref{eq:(C.17)}) from the \(\gre_{\mu}\) part of the functional \(\gra\). The other \(\grw_{\mu}^{+}\)
terms are bounded by the integral over \(X_{*}\) of \(c_{0}\big(r \Big|
|\psi|^{2}- |\frac{i}{2}w_{X}|\Big| |\grw_{\mu}^{+}|+|\grw_{\mu}^{+}|^{2}\big)\).
This understood, the bounds on \(|\psi|^{2}\) and \(|w_{X}|\) lead to a bound
on \(|\gri_{\mu}|\) by \(c_{0}\c rL_*\).

What follows explains how the term \(\gri_{K}\) in \(\gri_{*}\) arises. The
\(dw_{X}= 0\) assumption is used to derive a suitable bound on \(|i_{K}|\).
As noted above, the derivation starts by writing
\(|F_{\mathbb{A}}^{+}-r (\psi^{\dag}\tau\psi-\frac{i}{2}  w_{X}) +iw_{\mu}^{+}|^{2}\)
as \(|F_{\mathbb{A}}^{+}-r (\psi^{\dag}\tau\psi- \frac{i}{2}w_{X})|^{2}\) plus terms
that involve \(\grw_{\mu}^{+}\). The norm \(|F_{\mathbb{A}}^{+}-r (\psi^{\dag}\tau\psi-w_{X})|^{2}\)
is then written as a sum of \(|F_{\bbA}^{+}|^{2}\), \(r^{2}|\psi^{\dag}\tau\psi-\frac{i}{2}w_{X}|^{2}\)
and twice the inner product between \(F_{\mathbb{A}}^{+}\) and \(r(\psi^{\dag}\tau\psi- \frac{i}{2}w_{X})\).
The integral over \(X_{*}\) of the term with the inner product between
\(F_{\mathbb{A}}^{+}\) and \(r\psi^{\dag}\tau\psi\)  is cancelled by the contribution from the
\(F_{\mathbb{A}}^{+}\) term in the Weitzenb\"ock formula for \(\mathcal{D}_{\bbA}^-\mathcal{D}_{\bbA}^+\psi\).
The inner product between \(F_{\mathbb{A}}^{+}\) and \(- \frac{i}{2}rw_{X}\)
is equal to that of \(F_{\mathbb{A}}\) with \(-\frac{i}{2} rw_{X}\)
and thus its integral is that of \(r F_{\mathbb{A}}\wedge w_{X}\).
Stokes' theorem identifies most of the latter with the contributions
to the boundary terms on the right hand side of (\ref{eq:(C.17)}) from the
\(r \textsc{w}\) term in \(\gra\). The term designated by \(\gri_{K}\) is
what remains after the application of Stokes' theorem. To say more
about \(\gri_{K}\), note that the application here of Stokes' theorem
requires writing \(\mathbb{A}\) as \(A_{K}+2A_{E}+\hata _{\mathbb{A}}\) with \(\hata _{\mathbb{A}}\) being
an \(i\)-valued 1-form on \(X_{*}\). Stokes' theorem involves only \(\hata _{A}\). The \(\gri_{K}\)
term is the integral of \(\frac{i}{2}rF_{A_K+2A_E}\wedge w_{X}\). This understood, the bound \(|\gri_{K}|\leq c_{0}\c rL_{*}\)
follows from the \(|w_{X}|\leq \c\) assumption.

There is one other subtle point with regards to the derivation of
(\ref{eq:(C.17)}) in the case when \(X_{*}\) is the \(|s| \leq L+3\) part of \(X\), this being that
the application of Stokes' theorem
requires a Hermitian connection on the bundle
\(\det (\bbS^{+})\) whose curvature has norm bounded
by $r^{c_1/\c}$ with \(c_{1}\) being a constant that is independent of \(\grd\), \(r\),
\(\c\), the metric and \(w_{X}\).  The pull back of this
connection from the \(s \leq -L\)  and \(s \geq L\)  part of \(X\) via the
embeddings from the second and third bullets should also be the
respective \(Y_{-}\) and \(Y_{+}\) versions of
\(A_{K } + 2A_{E}\).  Such a connection can be constructed using the
isomorphism between de-Rham cohomology and the \v{C}ech cohomology that is defined by a cover of the
\(|s | \leq L + 1\) part of \(X\) by Gaussian coordinate
charts with the property that the any given number of charts have
either empty or convex intersection (see Chapter 8 in \cite{BT}).
 The $r^{1/\c}$-bound on the norm of Riemannian curvature and the 
$r^{-1/\c}$ lower bound on the injectivity radius can be used to obtain such a
cover by sets of radius greater than $r^{-c_0/\c}$.   As the connection is constructed from the de-Rham isomorphism
using a subbordinate partition of unity, this lower bound on the
minimum chart radius can be used to construct a connection on
\(\det (\bbS^{+})\) with an $r^{c_0/\c}$ bound on the norm of its curvature.  
 
Section \ref{sec:Cf)} says more about \(\gri_*\) when the \((\c, \r = r)\)
version of (\ref{eq:(A.16)}) is assumed.
 
\paragraph{\it Step 4:}  Define \(X_*\), \(\partial_-X_*\), and \(\partial_+X_*\) as in
Step 3. Granted Step 3's bound for the norm of the \(\gri_*\) term in (\ref{eq:(C.17)}), then
(\ref{eq:(C.15)}) and (\ref{eq:(C.17)}) imply that
\begin{equation}\label{eq:(C.18)}
\gra(\grd|_{\partial_+X_*}) \leq \gra(\grd|_{\partial_-X_*} ) + c_{0} \c^{2} r^{2} .
\end{equation}
This inequality with the top identity in (\ref{eq:(C.16)}) imply that
\(\gra(\grc_+)\leq \gra(\grd |_{s}) \leq \gra(\grc_{-}) + c_{0}\c^{2}
r^{2}\) when \(s\geq L\); and the identity in the
bottom bullet of (\ref{eq:(C.16)}) and (\ref{eq:(C.18)}) imply the inequalities
\(\gra (\grc_-)\geq \gra(\grd|_{s} ) \geq \gra(\grc_{+}) - c_{0}\c^{2}
r^{2}\) when \(s\leq -L\).
Given these inequalities, then (\ref{eq:(C.17)}) implies that
\begin{equation}\label{eq:(C.19)}
 \frac{1}{2}\int_{X_*}
 (|F_{\mathbb{A}}|^2+r^2|\psi
 ^\dag\tau\psi-iw_X|^2+2r|\nabla_{\mathbb{A}}\psi|^2) \leq \gra(\grc_{-}) - \gra(\grc_{+}) +c_{0} \c^{2} r^{2}.
\end{equation}

This last inequality with the identities in (\ref{eq:(C.15)}) and (\ref{eq:(C.16)}) imply
directly the assertion made by the first bullet of Lemma \ref{lem:C.1} and it
implies the second bullet when the length one interval is part of
\([-L-3, L+3]\) or \([-L_*-1, -L_*+5]\) or \([L_*-5, L_*+1]\).  

Granted what was just said, the second bullet of Lemma \ref{lem:C.1} holds if
its assertion is true when the length one interval is disjoint from
\([-L, L ]\), \([-L_*, -L_*+4]\) and \([L_*-4, L_*]\).  To prove the
assertion for these cases, use (\ref{eq:(C.18)}) with (\ref{eq:(C.15)}) and (\ref{eq:(C.16)})
to see that \(\gra(\grd|_{s}) -\gra(\grd|_{s '}) <c_{0} \c^{2} r^{2}\)
if \(s> s'\) and if both are in the same component of the complement
in \(\bbR\) of any of these three intervals.  This fact is exploited for the case \(s'= s +1\) using an integration by parts argument to rewrite the integrand on the
left hand side of the \(s' = s + 1\) version of (\ref{eq:(C.15)}) so as to have the
same form as the integrand on the left hand side of (\ref{eq:(C.17)}).  The
resulting inequality with the bound \(\gra(\grd|_{s}) -\gra(\grd|_{s+1}) < c_{0} \c^{2}
r^{2}\) leads directly to what is asserted by Lemma
\ref{lem:C.1}'s second bullet. \epf

 \subsection{Proof of Lemma \ref{lem:C.2}}
\label{sec:Cd)}
 
The proof of Lemma \ref{lem:C.2} has five steps.  By way of a look ahead, the
arguments depend crucially on the fact that the metric with the 2-form
\(ds \wedge *w + w\) define a K\"ahler structure on \(U_{C}\)
and on \(U_{0}\).  The proof that follows considers only the
special case where both \(\mu_{-}\) and
\(\mu_{+}\) vanish on the respective
\(Y_{-}\) and \(Y_{+}\) versions of
\(U_{\gamma}\) and $\mathcal{H}_0$.  The
argument in the general case is little different and so not given. 
 
\paragraph{\it Step 1:}  Let \(V_{*}\) denote either \(U_{C}\)
or \(U_{0}\).  The fact that the metric with \(ds \wedge
*w + w\) defines an integrable complex structure on
\(V_{*}\) has following consequence:  View \(\beta\) as a
section of the \((0, 2)\)-part of
\(\wedge^{2}T^*V_{*} \otimes \bbC\).  Then the right most equation in
(\ref{eq:(A.14)}) can be written on either \(U_{C}\) or \(U_{0}\) as 
\begin{equation}\label{eq:(C.20)}
\bar{\partial}_A\alpha +\bar{\partial}_A^\dag\beta = 0.
\end{equation}
This last equation implies that \(\beta\) obeys
\begin{equation}\label{eq:(C.21)}
\nabla_{A}^{\dag}\nabla_{A}\beta +  r (1 + |\alpha|^{2} +|\beta|^{2}) \beta + \grr\beta = 0,
\end{equation}
where \(\grr\) is determined solely by the metric.  In particular, the absolute
value of \(\grr\) and its derivatives to any specified order are also bounded
by \(c_{0}\).  The equation just written implies that
\(|\beta|^{2}\) obeys the differential inequality
\begin{equation}\label{eq:(C.22)}
d^{\dag}d|\beta|^{2}+r |\beta|^{2} +|\nabla_{A}\beta|^{2}\leq 0 .
\end{equation}
This last inequality is exploited momentarily with the help of the
Green's function for the operator \(d^{\dag}d + r\).   

Let \(x \in V_{*}\) denote a given point and let \(G_{x}(\cdot)\) denote the Dirichlet
Green's function for \(d^{\dag}d + r\)
with pole at \(x\).  Keep in mind for what follows the following fact
about \(G_{x}(\cdot)\):  It is non-negative and it obeys:
\begin{equation}\label{eq:(C.23)}
G_{x}(\cdot) \leq c_{0}  \frac{1}{\dist (x, \cdot)^2}e^{-\sqrt{r}\dist
  (x, \cdot)}.
  \end{equation}

Introduce \(\textsc{d}\co V_{*} \to [0, c_{0})\) to denote the function that measure the distance to
the boundary of \(V_{*}\).  Fix \(x\) in the interior of
\(\textsc{d}_{*}\), multiply both sides of (\ref{eq:(C.22)}) by
\(G_{x}( \cdot )\) and integrate the resulting
inequality over \(V_{*}\).  An integration by parts in the
left hand integral using the bound
\(|\beta|^{2} \leq c_{0}\c\) from Lemma \ref{lem:C.1} leads directly to the following
inequalities:
\begin{equation}\label{eq:(C.24)}
\begin{cases}
& |\beta|^{2} \leq
c_{0} \c e^{-\sqrt{r}\, \textsc{d}},\\
& \int_BG_x|\nabla_A\beta|^2 
 \leq c_{0} \c\frac{1}{\textsc{d}^2}e^{-\sqrt{r}\, \textsc{d}}.
\end{cases}
\end{equation}
The second inequality is used in Step 3 to derive bounds on the higher
order derivatives of \(\beta\).
 
\paragraph{\it Step 2:}  This step constitutes a digression to state some very
crude bounds for the norms of \(F_{\bbA}\), \(\nabla_{\bbA}\psi\) and their covariant derivatives.
 The following lemma states these bounds.  
\begin{lemma}\label{lem:C.4}
There exists $\kappa>\pi $ such that given any \(\c >  \kappa\),
 there exists \(\kappa_{\c}\) with the
following significance:  Fix \(r \geq\kappa_{\c}\)  and assume the \((\c, \r = r)\) 
version of the first two bullets of (\ref{eq:(A.16)}).  Assume in addition that
 \(|w_{X}| \leq \c\)  and that the norms of its derivatives to order 10 are bounded by 
\(r^{1/\c}\).  Fix respective elements \(\mu_{-}\)  and \(\mu_{+}\)
 from the \(Y_{-}\)  and \(Y_{+}\)  versions of \(\Omega\) with
 \(\mathcal{P}\)-norm bounded by 1.   Use this data to define
the equations in (\ref{eq:(A.14)}).  Let \(\grd = (\bbA, \psi)\) 
denote an instanton solution to (\ref{eq:(A.14)}) with
\(F_{\bbA}\)  and \(r^{1/2}|\nabla_{\bbA}\psi|\) having \(L^{2}\)-norm less than \(\c
r\)  on the \(s\)-inverse image of any
length 1 interval in \(\bbR\).   Then the norm of
\(F_{\bbA}\)  and \(|\nabla_{\bbA}\psi|\),
and those of their derivatives up through order 4  are bounded
everywhere by \(\kappa_{\c} r^{\kappa  _{\c}}\).  
\end{lemma}
 
\pf This follows using a standard
elliptic boot-strapping argument since the equations in (\ref{eq:(A.14)}) can be
viewed as elliptic equations on any given ball in \(X\) for a suitable pair
on the \(C^{\infty}(X; S^{1})\)-orbit of
\((\bbA, \psi)\).  Except for one remark, the details of this
bootstrapping are completely straightforward and so they will not be
presented.  The remark concerns the fact that the assumed lower bound
for the injectivity radius is needed for the proof so as to invoke
various Sobolev embedding theorems using embedding constants that are
bounded by powers of \(r\). \epf
 
The bounds supplied by Lemma \ref{lem:C.4} are used in the next step.

\paragraph{\it Step 3:}  To obtain the asserted bound for the covariant derivative
of \(\beta\), differentiate (\ref{eq:(C.21)}) and commute covariant derivatives to
obtain an equation for \(\nabla_{A}\beta\) that has
the schematic form
\begin{equation}\label{eq:(C.25)}\begin{split}
& \nabla_{A}^{\dag}\nabla_{A}(\nabla_{A}\beta)+ r (1 + |\alpha|^{2}+|\beta|^{2})\nabla_{A}\beta + \\
& \quad \quad  \grR_{0}(F_{A})\nabla_{A}\beta+ \grR_{1}(\nabla F_{A})\beta + r \grR_{2}( \nabla_A\psi)\nabla_A\beta +\grr_{1}\nabla_{A}\beta = 0,
\end{split}
\end{equation}
where \(\grR_{0}\), \(\grR_{1}\) and \(\grR_{2}\) are
endomorphisms that are linear functions of their entries and are such
that \(|\grR_{*}(b)| \leq c_{0}|b|\).  Meanwhile, \(\grr_{1}\) is such that \(|\grr_{1}|
\leq c_{0}\) also.  Take the inner product of both
sides of (\ref{eq:(C.25)}) with \(\nabla_{A}\beta\) and invoke
Lemma \ref{lem:C.4} to see that
\begin{equation}\label{eq:(C.26)}
d^{\dag}d(|\nabla_{A}\beta|^{2})+ r|\nabla_{A}\beta|^{2}+|\nabla_{A}\nabla_{A}\beta|^{2}
\leq c_{\c} r^{c_{\c}} (|\nabla_{A}\beta|^2+ |\beta|^{2}),
\end{equation}
where \(c_{\c}\) here and in what follows denotes a constant
that is greater than 1 and depends only on \(\c\).  The value of
\(c_{\c}\) can be assumed to increase between consecutive
appearances. 

 Fix a point \(x \in V_{*}\) with distance greater than \(c_{0} r^{-1/2} (\ln  r)^{2}\)
from the boundary of \(V_{*}\).  Having done so, multiply
both sides of (\ref{eq:(C.26)}) by \(G_{x}\) and integrate both sides of
\(V_{*}\).  Use the second bullet in (\ref{eq:(C.24)}) to bound
integral on the right hand side of the resulting inequality by
\(c_{0}e^{-\sqrt{r}/c_0}\) when \(r \geq c_{\c}\) .  An integration by parts on the
left hand side using Lemma \ref{lem:C.4} to bound \(|\nabla_{A}\beta|\) on the
boundary of \(V_{*}\) and the bound just stated implies that
\begin{equation}\label{eq:(C.27)}
|\nabla\beta_{A}|^{2}(x)
+ \int_B G_x |\nabla_A\nabla_A\beta |^2 \leq c_{0} e^{-\sqrt{r}/c_0}
\end{equation}
when \(r \geq c_{\c}\).  This gives the desired bound for \(|\nabla_{A}\beta|\).

To obtain the bound for
\(|\nabla_{A}\nabla_{A}\beta|\), differentiate (\ref{eq:(C.25)}) twice
and take the inner product of both sides with \(\nabla_{A}\nabla_{A}\beta\)
after commuting covariant derivatives.  The result is an equation that
looks much like (\ref{eq:(C.26)}) with \(\nabla_{A}\beta\) replaced by
\(\nabla_{A}\nabla_{A}\beta \) on the left hand side and with the addition of the term \( r^{c_{\c}}  
|\nabla_{A}\nabla_{A}\beta|^{2}\)
on the right hand side.  Granted that this is the case, then the same
Green's function argument that led to (\ref{eq:(C.27)}) leads to
an analogous bound for \(|\nabla_{A}\nabla_{A}\beta|^{2}\).
 
\paragraph{\it Step 4:}  This step and Step 5 addresses the assertions of Lemma \ref{lem:C.2} that
concern \(\alpha\).  To start, act by \(\bar{\partial}_A^\dag\) on both sides of (\ref{eq:(C.20)}), commute covariant derivatives and use the
bounds from Lemma \ref{lem:C.2} for \(|\beta|\) to see that
\(\alpha\) obeys an equation that has the form
\begin{equation}\label{eq:(C.28)}
\nabla_{A}^{\dag}\nabla_{A}\alpha- r (1 - |\alpha|^{2}) \alpha= \gre ,
\end{equation}
where \(|\gre| \leq e^{-\sqrt{r}/c_0}\) when \(r \geq c_{\c}\).  This
equation implies that \(\w = 1- |\alpha|^{2}\) obeys a differential inequality of the form
\begin{equation}\label{eq:(C.29)}
d^{\dag}d\w + r \w \geq|\nabla_{A}\alpha|^{2}+ r \w^{2} - e^{-\sqrt{r}/c_0}.
\end{equation}
Use of the Green's function \(G_{x}\) with the
fact that \(|\w| \leq c_{0}\c\) on the boundary of \(V_{*}\) along the
same lines as in Parts 1 and 3 finds \(\w \geq e^{-\sqrt{r}/c_0}\) at distances greater than \(c_{0} r^{-1/2}(\ln r)^{2}\) from the boundary of \(V_{*}\) when \(r\geq c_{\c}\).  This is the
\(|\alpha|^{2}\) assertion in the first bullet of Lemma \ref{lem:C.2}.
 
The assertion in the third bullet follows directly from (\ref{eq:(C.20)}) given
Lemma \ref{lem:C.2}'s bounds for \(|\beta|\)
and \(|\nabla_{A}\beta|\).  The
assertion in the fourth bullet follows directly from (\ref{eq:(C.29)}) given that 
\(\w (1  - \w) = |\alpha|^{2} (1-|\alpha|^{2})\) and that this is
greater than \(\frac{1}{2}\delta^{2}\) at points where
\(|\alpha|^{2}\) is between \(\delta\)
and \(1  - \delta\).  The assertions in the fifth bullet about the
covariant derivatives of \(\alpha\) are proved in Step 5.  
 
\paragraph{\it Step 5:}  This step derives the asserted bounds in the fifth bullet
for the norms of the covariant derivatives of \(\alpha\).  To do this, 
suppose that \(x \in V_{*}\) is such that
\(|F_{A}| \leq \c r\) on the ball of radius \(c_{0} r^{-1/2}\)
centered at \(x\).   Use \(\r_{r}\) in what follows to denote the
rescaling map from \(\bbC^{2}\) to \(\bbC^{2}\) that is given by the rule \(x
\mapsto \r_{r}(x) = r^{-1/2} x\).
 The pull-back of \((A, \psi)\) by this map is denoted by
\((A_{r}, \psi_{r})\).  The bound \(|F_{A}| \leq \c r\) implies that the
absolute value of the curvature of \(A_{r}\) is bounded in the
radius 1 ball about the origin in \(\mathbb{C}^{2}\)  is
bounded by \(\c\).  Meanwhile, the pull-back of the equations in (\ref{eq:(A.14)}) by
this map constitutes a uniformly elliptic system of equations (modulo
the action of
\(C^{\infty}(\bbC^{2};S^{1})\) in the radius 1 ball about the origin in
\(\mathbb{C}^{2}\)  with coefficients that have
\(r\)-independent bounds for their absolute values and for those of their
derivatives to any a priori chosen order.  This understood, the fact
that \(|\psi_{r}|\leq 2\) in this
ball and the afore-mentioned bound by \(\c\) for the norm of the curvature
of \(A_{r}\) imply via standard elliptic bootstrapping
arguments that the \(A_{r}\)-covariant derivatives of
\(\psi_{r}\) through order 2 are bounded by \(c_{0} \c\) in the radius \(c_{0}^{-1}\) ball about the
origin in \(\bbC^{2}\).  Granted these bounds, use
the chain rule of calculus to obtain the bounds asserted by the fifth
bullet of Lemma \ref{lem:C.2} for the covariant derivative of
\(\alpha\). \epf

\subsection{Proof of Lemma \ref{lem:C.3}}\label{sec:Ce)}
Use \(V_{*}\) again to denote either \(U_{C}\)
or \(U_{0}\). The functions \(\grz_{A}\) and
\(\grz_{B}\) are both equal to
\(r|\beta|^{2}\) on
\(V_{*}\) and so what is asserted by the second bullet of
Lemma \ref{lem:C.3} follows from the first bullet of Lemma \ref{lem:C.2}.  The absolute
value of \(r\) is bounded by \(c_{0}r|\alpha| |\beta|\) on
\(V_{*}\) and so the third bullet of Lemma \ref{lem:C.3} also follows
from the first bullet of Lemma \ref{lem:C.2}.  The bounds in the first bullet of
Lemma \ref{lem:C.3} follow from the bound in the fourth bullet and that for
\(|\alpha|^{2}\) in the first bullet
of Lemma \ref{lem:C.2}.  If the bounds in first through fourth bullets of Lemma
\ref{lem:C.3} hold, then \(|F_{A}|\) is bounded by
\(c_{0} r\) at the points in \(V_{*}\) with distance \(\frac{1}{200}\rho_{D}\) from the boundary of \(V_{*}\).
 Granted that this is the case, then the rescaling argument in Step 5
of the proof of Lemma \ref{lem:C.2} can be used to derive the bound given in the
fifth bullet of Lemma \ref{lem:C.3}.  
 
The upcoming Lemma \ref{lem:C.5} is the critical ingredient for the proof of the
fourth bullet of Lemma \ref{lem:C.3}.   The \(\gra(\grc_{-}) -\gra(\grc_{+}) \leq r^{2-1/\c} \) assumption in Lemma \ref{lem:C.3} and the
final  three bullets of (\ref{eq:(A.16)}) are needed only to invoke Lemma \ref{lem:C.5}.

\begin{lemma}\label{lem:C.5}   
There exists \(\kappa >100 (1 + \rho_{D}^{-1})\) such that 
given any \(\c \geq \kappa \), there exists \(\kappa_{\c} > \kappa  \) with the
following significance:  Fix \(r \geq\kappa_{\c}\) and assume that the metric and 
\(w_{X}\) are \((\c, \r = r)\)-compatible.  Fix elements \(\mu_{-}\)  and
\(\mu_{+}\) from the \(Y_{-}\)  and \(Y_{+}\) versions
of \(\Omega\) with  \(\mathcal{P}\)-norm bounded by 1
and use this data to define the equations in (\ref{eq:(A.14)}).   Let
\(\grc_{-}\)  and \(\grc_{+}\)  denote
solutions to the \((r,  \mu_{-})\)-version of
(\ref{eq:(A.4)}) on \(Y_{-}\)  and the \((r, \mu_{+})\)  version of (\ref{eq:(A.4)}) on
\(Y_{+}\)  with \(\gra(\grc_{-}) -\gra(\grc_{+}) \leq  r^{2-1/\c}  \).  Let \(\grd = (\bbA, \psi)\)  denote an
instanton solution to (\ref{eq:(A.14)}) with \(s \to
-\infty  \) limit equal to \(\grc_{-}\)  and \(s\to \infty  \) limit equal to
\(\grc_{+}\).  Use \(B\)  to denote a ball of radius \(\kappa^{-2}\)
in the domain \(U_{C}\)  or in the domain \(U_{0}\)  with center at distance
\(\kappa^{-1}\)  or more from the domain's boundary.  Then 
\(r\int_B |1-|\psi|^2|\leq  \kappa_{\c} r^{1-1/\kappa  _{\c}}\).
\end{lemma}
 
Lemma \ref{lem:C.5} is proved in Section \ref{sec:Cf)}.  Granted Lemma
\ref{lem:C.5}, then the six steps that follow prove the fourth bullet
of Lemma \ref{lem:C.3} in the case when \(\mu_{-}\) and \(\mu_{+}\) are zero on
the \(Y_{-}\) and \(Y_{+}\) version of \(U_{\gamma}\) and \(\mathcal{H}_ 0\) .  The
proof when they are not zero but bounded by \(e^{-r^2}\)
is little different and so not given. 
 
\paragraph{\it Step 1: }  Let \(V_{*}\) denote either
\(U_{C}\) or \(U_{0}\).  Keep in mind that metric on
\(V_{*}\) has non-negative Ricci curvature tensor, that the
2-form \(w_{X} = w\) is covariantly constant on
\(V_{*}\), that both \(\grw_{\mu} = 0\) and that \(B_{A_K}\)  
 is covariantly constant on \(V_{*}\).  These facts with the
bounds from Lemma \ref{lem:C.2} for \(|\beta|\) and
\(|\nabla_{A}\beta|\) have the following implication:  Let \(\ss\) denote \(|\mathcal{E}_{A} -
\mathcal{B}_{A}|\).  Granted that \(r \geq
c_{\c}\), then the equations in (\ref{eq:(A.14)}) imply that \(\ss\) obeys the
differential inequality
\begin{equation}\label{eq:(C.30)}
d^{\dag}d\ss + r|\alpha|^{2} \ss \leq r|\nabla_{A}\alpha|^{2}+ e^{-\sqrt{r}/c_0}
\end{equation}
at the points in \(V_{*}\) with distance greater than
\(c_{0}r^{-1/2} (\ln r)^{2}\)
from the boundary of \(V_{*}\).  Let \(\w\) again denote \(1-|\alpha|^{2}\) and let
\(q_{0}\) denote \(\ss - r \w\).  It follows from (\ref{eq:(C.29)})
and (\ref{eq:(C.30)}) that 
\begin{equation}\label{eq:(C.31)}
d^{\dag}d q_{0} + r|\alpha|^{2} q_{0}\leq    e^{-\sqrt{r}/c_0}
\end{equation}
at the points in \(V_{*}\) with distance \(c_{0}
r^{-1/2} (\ln r)^{2}\) or more from
\(V_{*}\)'s boundary if \(r \geq c_{\c}\).

\paragraph{\it Step 2:}  Fix \(\rho_{*} > 0\) but such
that \(\rho_{*} <10^{-8} \rho_{D}\).  Fix
\(s_{0} \in \bbR\).  Let \(V' \subset V_{*}\) denote the set of points in the
\((s_{0}- 1 - \rho_{*}, s_{0} + 1 + \rho_{*})\) part of
\(V_{*}\) with distance \(\rho_{*}\) or
more from the boundary of \(V_{*}\), and let \(V \subset
V'\) denote the set of points in \(V_{*}\) with distance
greater than \(2 \rho_{*}\)  from the boundary of
\(V_{*}\).  Thus, each point in \(V\) has distance
\(\rho_{*}\) or more from the boundary of \(V'\).  
 
Fix a sequence \(\{ \varsigma _{n}\}_{n=1,\ldots}\) of
smooth, non-negative functions on \(V'\) with the following properties:
 Each function in this series is bounded by 1 and is equal to 1 on \(V\).
 Second,  \(\varsigma_{1}\)  has compact support and for
each \(n \geq 1\), the function  \(\varsigma_{n+1}\)  has
compact support where  \(\varsigma _{n} = 1\).  Finally,
the absolute values of the first and second derivatives of the
functions in this series enjoy \(s_{0}\)-independent upper
bounds.
 
\paragraph{\it Step 3:}  For each integer \(n \geq 1\), set \(q_{n} =
\max (\varsigma _{n} q_{n-1}, 0)\).  Use \(\q_{0}\) to denote the maximum of \(q_{0}\); and for
\(n \geq 1\), use \(\q_{n}\) to denote the maximum of
\(q_{n}\).  Note that \(\q_{n} \leq\q_{n-1}\).  It follows from (\ref{eq:(C.31)}) that if \(r \geq
c_{\c}\), then any given \(n \geq 1\) version of
\(q_{n}\) obeys
\begin{equation}\label{eq:(C.32)}
d^{\dag}dq_{n} +  r|\alpha|^{2} q_{n}\leq (dd^{\dag} \varsigma _{n})
\, q_{n-1} +  2 \langle d \varsigma _{n},dq_{n-1} \rangle  + c_{0}e^{-\sqrt{r}/c_0},
\end{equation}
where  \(\langle  \cdot , \cdot \rangle  \) denotes the metric inner product.  Fix
a constant \(z_{n} \geq 1\) to be determined shortly, and
let \(q_{n*}\)  denote the maximum of 0 and
\(q_{n}  - r^{-1}z_{n}\q_{n-1}\).  The function \(q_{n*}\)  obeys
\begin{equation}\label{eq:(C.33)}\begin{split}
d^{\dag}d q_{n *} & + r|\alpha|^{2}q_{n *} \\
 \leq  & z_{n}\q_{n-1} \w + \big(-z_{n} \q_{n-1} +
(dd^{\dag}\varsigma _{n})\, q_{n-1} +  2 \langle d \varsigma _{n},dq_{n-1} \rangle \big) + c_{0} e^{-\sqrt{r}/c_0}.  
\end{split}\end{equation}
Note also that \(q_{n*}\) has compact support
in \(V'\) since \(q_{n} -r^{-1}z_{n} \q_{n-1} =
-r^{-1} z_{n} \q_{n-1}\) on the complement of the support of  \(\varsigma_{n}\) .
 
 \paragraph{\it Step 4:}  Fix \(x\) in the interior of \(V'\) and let \(G_{x}\) now
denote the Dirichlet Green's function for the operator
\(d^{\dag}d\) on \(V'\) with pole at \(x\).  The function
\(G_{x}\) is non-negative, \(|G_{x}(\cdot)| \leq c_{0} \dist (x, \cdot)^{-2}\)
and \(|dG_{x}(\cdot)|\leq c_{0} \dist(x,\cdot)^{-3}\).  Multiply both sides of
(\ref{eq:(C.33)}) by \(G_{x}\) and integrate the two sides of the
resulting inequality over \(V'\).  Integrate by parts on both sides to
remove derivatives from \(q_{n*}\)  and \(q_{n-1}\) to obtain the inequality
\begin{equation}\label{eq:(C.34)}
q_{n *}(x) \leq z_{n}\q_{n-1}\int_{V'}\big(\frac{1}{\dist(x,\cdot)^{2}}\w\big)
+  ( - c_{0}^{-1}z_{n }+ \mathpzc{e}_{n}) \q_{n-1} +  e^{-\sqrt{r}/c_0}.
\end{equation}
where \(\mathpzc{e}_{n} \leq c_{0}\sup_{x\in V'}(|d^{\dag}d \varsigma _{n}|+ |d \varsigma _{n}|)\).  Granted this
bound, a purely \(n\)-dependent choice for \(z_{n}\) leads from
(\ref{eq:(C.34)}) to the inequality
\begin{equation}\label{eq:(C.35)}
q_{n *}(x) \leq z_{n}\q_{n-1} \int_{V'}\big(\frac{1}{\dist(x,\cdot)^{2}}\w\big)+  e^{-\sqrt{r}/c_0},
\end{equation}
Lemma \ref{lem:C.5} is used to exploit this inequality.

\textit{  Step 5:}  Fix \(\rho > 0\) and break up the integral
in (\ref{eq:(C.35)}) into the part where \(\dist (x, \cdot)\) is greater
than \(\rho\) and the part where \(\dist (x, \cdot )\) is
less than \(\rho\).  Having done so, appeal to Lemma \ref{lem:C.5} and the
first bullet of Lemma \ref{lem:C.2} to see that
\begin{equation}\label{eq:(C.36)}
q_{n *}(x) \leq z_{n}(\rho^{-2} r^{-1/c_0}+ \rho^{2}) \, \q_{n-1}  +e^{-\sqrt{r}/c_0}
\end{equation}
when \(r \geq c_{\c }\).  Let
\(c_{*}\) denote the value of \(c_{0}\) that
appears in (\ref{eq:(C.36)}).  Take \(\rho =r^{-1/4c_*}\)
 in (\ref{eq:(C.36)}).  The resulting right hand side is independent of \(x\); and
this leads directly to the inequality
\begin{equation}\label{eq:(C.37)}
q_{n} \leq z_{n}  r^{-1/2c_*}q_{n-1} + e^{-\sqrt{r}/c_0}
\end{equation}
when \(r \geq c_{\c}\).  As Lemma \ref{lem:C.4} finds
\(q_{0} <r^{c_{\c}}\), what is written in (\ref{eq:(C.37)}) implies that an \(n = c_{\c}\)
version of \(\q_{n}\) is bounded by \(r^{-200}\).

 \paragraph{\it Step 6:}  Since  \(\varsigma _{n} = 1\) on \(V\), the
conclusion from Step 5 implies that
\begin{equation}\label{eq:(C.38)}
|\mathcal{E}_{A} - \mathcal{B}_{A}| <r (1  - |\alpha|^{2})  +r^{-200 }
\end{equation}
at all points in \(V\).  Square both sides of (\ref{eq:(C.38)}).  What with the
bounds for \(|\grz_{A}|\) and \(|\grz_{B}|\) from Lemma \ref{lem:C.3}'s
second bullet, the resulting inequality implies that 
\begin{equation}\label{eq:(C.39)}
(1  - 2\sigma)^{2} r^{2} (1-|\alpha|^{2}) +|\grX|^{2} \leq r^{2}(1  - |\alpha|^{2}) +c_{0} r^{-198};
\end{equation}
and rearranging terms writes this as 
\begin{equation}\label{eq:(C.40)}
|\grX|^{2} \leq 2r^{2} \sigma (1  - \sigma) (1  -|\alpha|^{2}) + c_{0}r^{-198}.
\end{equation}
This is gives the bound stated in the fourth bullet of Lemma \ref{lem:C.3}.  \epf

 \subsection{Proof of Lemma \ref{lem:C.5}}\label{sec:Cf)}

 The proof has six parts.  Parts 1 and 2 revisit the formula in
(\ref{eq:(C.15)}) and Part 3 revisits the formula in (\ref{eq:(C.17)}).
These steps present the proof in the case when
\(c_{1}(\det(\mathbb{S}))\) is non-torsion on all components of the
\(|s| > 1\) part of \(X\). But for the two remarks that follow, the proof when \(X_{tor} \neq \emptyset\)
differs only cosmetically. 

The first remark concerns the formula in (\ref{eq:(C.17)}) in the case when \(X_{*}\) is the
respective \(|s| \in[L,L_{*}- 4]\) part of \(X_{tor}\), the remark being that the absolute value of \(\gri_{*}\) in this case is bounded by
\(c_{0} \c^{2}r \ln r\). The reason is as follows: As noted subsequent to (\ref{eq:(C.17)}), the absolute value of
the relevant version of \(\gri_{*}\) is bounded in any event by \(c_{0}\c r L_{*}\).
Meanwhile, the first bullet of (\ref{eq:(A.16)}) bounds \(L_{*}\) by \(\c \ln r\). 

The second remark concerns (\ref{eq:(C.17)}) in the case when
\(X_{*}\) is the \(|s| \in[L_{*}-4, L_{*}]\) part of \(X_{tor}\), this
being that the absolute value of the corresponding version of
\(\gri_{*}\) is at most \(c_{0}\) when \(r\) is larger than a purely
\(\c\)-dependent constant. Given Item
d) of the fourth bullet of (\ref{eq:(A.16)}), the proof that this
is so differs only in notation from what is said below in Part 2 to
prove the analogous bound for the version of
\(\gri_{*}\) that appears in (\ref{eq:(C.17)}) when \(X_{*}\) is the \(|s| \in[L -4, L]\) part of \(X\).

\paragraph{Part 1:}  Write \(\grd = (\bbA, \psi)\).  When \(X\), the
metric and \(w_{X}\) are described by the first bullet of
(\ref{eq:(A.17)}), use this pair as instructed in the proof of Lemma
\ref{lem:C.1} to define the map \((A_{*}, \psi_{*})\) from
\(\bbR\) to \(\op{Conn} (E) \times C^{\infty}(Y_{Z}; \bbS)\).  When
the second bullet of (\ref{eq:(A.17)}) is relevant, then \((A_{*}, \psi_{*})\) as defined in the proof of Lemma \ref{lem:C.1}
denotes a map from \((-\infty, -1]\) to \(\op{Conn} (E|_{Y_-}) \times C^{\infty}(Y_{-};
\bbS  |_{Y_-})\) and also a map from \([1, \infty)\) to \(\op{Conn} (E|_{Y_+}) \times C^{\infty}(Y_{+};\bbS   |_{Y_+})\). 
 
Set \(I_{L} = [-L, L]\) when \(X\), the metric and
\(w_{X}\) are described by the first bullet of (\ref{eq:(A.17)}), and set
\(I_{L}\) to be either \([-L, -L + 4]\) or \([L  - 4, L]\) otherwise.
 Use \(Y_{*}\) to denote the constant \(s \in
I_{L}\) slice of \(X\), this being either \(Y_{Z}\),
\(Y_{-}\) or \(Y_+\).  Write the metric on
\(I_{L} \times Y_{*}\) as
\(ds^{2} + \grg\) with \(\grg\) denoting an \(s\)-dependent metric on
\(Y_{*}\).   Define the \(s\)-dependent 1-form
\(w_{*}\) on \(Y_{*}\) by writing
\(w_{X}\) as \(ds \wedge *w_{*} +
w_{*}\) with the Hodge dual defined here by \(\grg\).  The two
equations in (\ref{eq:(A.14)}) on the \(s \in I_{L}\) part of \(X\) are
equivalent to equations for \((A_{*}, \psi_{*})\) that can be written as
\begin{equation}\label{eq:(C.41)}
\begin{cases}
\frac{\partial   }{\partial s}A_{*} + \grB_{\grd}  = 0 &   \text{and}\\ 
 \frac{\partial   }{\partial s}\psi_{*} +  D_{A_*}\psi_* = 0, &  
 \end{cases}
\end{equation}
with \(\grB_{\grd}\) denoting the following \(s \in I_{L}\) dependent 1-form on \(Y_{*}\):
\begin{equation}\label{eq:(C.42)}
\grB_{\grd} = B_{\bbA} - r(\psi^{\dag}\tau\psi - iw_{*}) +
i\grw_{\mu}^{+}(\frac{\partial  }{\partial s}, \cdot ) +  \frac{1}{2}B_{A_K}.
\end{equation}
 
By way of notation, \(D_{A_*}\) in (\ref{eq:(C.41)}) denotes the Dirac
operator defined by the metric \(\grg\), its Levi-Civita connection and the connection \(A_{K} +2A_{*}\) on the \(\{s \} \times Y_{*}\) version of \(\det (\bbS)\).  
 
 \paragraph{Part 2:}  If \(X\), the metric and \(w_{X}\) are described
by the first bullet of (\ref{eq:(A.17)}), then the integration and use of
Stokes' theorem that leads to (\ref{eq:(C.15)}) can be repeated
with the domain of integration being \(s^{-1}([-L, L])\)
to find that
\begin{equation}\label{eq:(C.43)}
 \begin{split}\frac{1}{2}\int_{\mathbb{R}\times Y_Z} & \big(
 |\frac{\partial}{\partial s} A_*|^2 +|\grB_{(A, \psi)}|^2+2r(|\frac{\partial}{\partial s}\psi_*|^2+|D_A\psi|^2)\big)
 + \gri_{\mu} \\
& = \gra(\grd|_{s=-L}) - \gra(\grd|_{s=L}),
\end{split}
\end{equation}
where \(\gri_{\mu} = 0\) when \(\grw_{\mu}\) is such that \(X\), the metric, \(w_{X}\) and
\(\grw_{\mu}\) define the product metric, and where
\(|\gri_{\mu}| \leq c_{0} \big(\int_{s^{-1}([-L, L])} |\frac{\partial}{\partial s} A_*|^2  \big)^{1/2}\) in any event.  This
being the case, the second bullet of Lemma \ref{lem:C.1} implies that
\(|\gri_{\mu}| \leq c_{0} r\).   
 
Assume now that \(X\), the metric and \(w_{X}\) are described by
the second bullet in (\ref{eq:(A.17)}).  The derivation of
(\ref{eq:(C.15)}) and (\ref{eq:(C.43)}) can be repeated with the
domain of integration being \(s^{-1}([-L,  -L +4])\) and also
\(s^{-1}([L  - 4, L])\) to obtain the following identities:
\begin{equation}\label{eq:(C.44)}
\begin{cases}
& \frac{1}{2}\int_{[-L, -L+4]\times Y_Z} \big(
 |\frac{\partial}{\partial s} A_*|^2+|\grB_{(A, \psi)}|^2+2r(|\frac{\partial}{\partial s}\psi_*|^2+|D_A\psi|^2)\big)
 + \gri    \\
&\qquad = \gra(\grd|_{s=-L}) - \gra(\grd|_{s=-L+4}),\\
& \frac{1}{2}\int_{[L-4, L]\times Y_Z} \big(
 |\frac{\partial}{\partial s} A_*|^2+|\grB_{(A, \psi)}|^2+2r(|\frac{\partial}{\partial s}\psi_*|^2+|D_A\psi|^2)\big)
 + \gri  \\
& \qquad = \gra(\grd|_{s=L-4}) - \gra(\grd|_{s=L}),\\
\end{cases}
\end{equation}
where \(\gri\) in this case is such that \(|\gri | \leq c_{0}
r^{2-1/\c}\) when \(\c > c_{0}\) and \(r > c_{\c}\) with \(c_{c}\) denoting a constant that
depends only on \(\c\).  The paragraphs that follow explain how this bound
comes about. 
 
The term denoted by \(\gri\) can be written as the sum of three integrals, \(\gri =
\gri_{\grg }+ \gri_{w} +  \gri_{\mu}\).
 What is denoted by \(\gri_{\mu}\)  appears here for the
same reason it appears in (\ref{eq:(C.43)}) and it has the analgous bound,
\(|\gri_{\mu}| \leq c_{0} r\).   The integral denoted by \(\gri_{\grg}\)
accounts for the \(s\)-dependence of the metric \(\grg\) on \(Y_{*}\)
when commuting the operators \(\frac{\partial}{\partial s}\) and \(D_{A_*}\).  In particular, the integrand that defines \(\gri_{\grg}\) is
bounded by \(c_{0} r \big( | \frac{\partial}{\partial s} \grg|  |\psi|
|\nabla_{A}\psi | +|\grR_{\grg}( \frac{\partial}{\partial s}, \cdot)|, |\psi|^{2}\big)\) with
\(\grR_{g}\) denoting the Riemannian  curvature tensor of the metric
\(ds^{2} + \grg\).  This understood, (\ref{eq:(A.16)}) with Lemma
\ref{lem:C.1}'s bounds for \(|\psi|^{2}\) and the \(L^{2}\)-norm of
\(|\nabla_{A}\psi|\) imply that \(|\gri_{\grg}| \leq c_{0} r^{3/2+1/\c}  \).

The integral \(\gri_{w}\) arises from the contribution to the
integral of \(|\frac{\partial}{\partial s}A _{*} +\grB_{\grd}|^{2}\) of the metric inner
product of \(\frac{\partial}{\partial s}A_{*}\) with \(-i r *w_{*}\).  The
integral of this inner product is written as \(\int_{I_L} h(s)ds\)
with \(I_{L}\) denoting \([-L, -L + 4]\) or \([L  - 4, L]\) as the
case may be; and with \(h(s)\) denoting the integral of the 3-form \(-i r   
 \frac{\partial}{\partial s}A_{*} \wedge w_{*}\) over \(\{s \}
\times Y_{*}\).  Only a portion of the integral of \(-i
r \frac{\partial}{\partial s}A_{*} \wedge  w_{*}\) contributes to
\(\gri_{w}\).  To say more, write \(A_{*}\) as
\(A_{E} + \hata _{A}\) with \(\hata _{A}\)
denoting an \(s\)-dependent 1-form on \(Y_{*}\).  The
integral of the 3-form \(-i r \frac{\partial}{\partial s}A_{*}\wedge  w_{*}\) over \(\{s\} \times  Y_{*}\) is written using \(\hata _{A}\)
as
\begin{equation}\label{eq:(C.45)}
 -ir  \frac{\partial}{\partial s}\big(\int_{\{s\}\times Y_*}\hata_A\wedge w_{*}\big)
+ i r \big(\int_{\{s\}\times Y_*}\hata_A\wedge \frac{\partial}{\partial s}w_{*}\big).
\end{equation}
 
The contributions of the function \(\textsc{w}\) in (\ref{eq:(A.6)}) to
the right hand side of
(\ref{eq:(C.44)}) are given by the integral over \(I_{L}\) of the left
most term in (\ref{eq:(C.45)}), this being a consequence of the fundamental
theorem of calculus.  What is denoted by \(\gri_{w}\) is the
integral over \(I_{L}\) of the right most term in (\ref{eq:(C.45)}).  A
bound for the absolute value of the latter is obtained by using the the
assumption in Item b) of the fourth bullet of (\ref{eq:(A.16)}) to write  
\(\frac{\partial}{\partial s}w_{*}\) as \(d\mathpzc{b}\) with \(\mathpzc{b}\) as described by this same part of
(\ref{eq:(A.16)}).   Stokes' theorem equates the the right most
integral in (\ref{eq:(C.45)}) with the integral of \(i r \, d\hata _{A}
\wedge \mathpzc{b}\).  This being the case, it follows from (\ref{eq:(A.16)}) that this
second contribution to \(\gri_{w}\) has absolute value less than
\(c_{0}r^{2-1/\c}\).
 
\paragraph{Part 3:}  Integrate
\(|F_{\bbA}^{+} - r(\psi^{\dag}\tau\psi - iw_{X}^{+})   - i\grw_{\mu}^{+}|^{2}+ r|D_{\bbA}\psi|^{2}\)
over \(s^{-1}([-L+ 4, L  - 4])\).  Integrate by parts
using the fact this integral is zero to derive an identity that can be
written as
\begin{equation}\label{eq:(C.46)}
\begin{split}
\frac{1}{2}\int    _{s^{-1}([-L+4, L-4])}
&\big(|F_\bbA|^2+r^2|\psi^\dag\tau\psi-iw_X^+|^2+2r|\nabla_\bbA\psi|^2\big)+
\gri_{L} \\
& = \gra(\grd|_{s=-L+4}) - \gra(\grd  |_{s=L-4})
\end{split}
\end{equation}
with \(\gri_{L}\) such that \(|\gri_{L}|
\leq c_{0}r^{1+c_0/\c}\).    The paragraphs that follow momentarily
derive the latter bound. By way of comparison, the absolute value of
the term \(\gri\) in (\ref{eq:(C.17)}) has the bound \(c_{0} \c  r \big(\int    _{s^{-1}([-L+4, L-4])}|F_\bbA  |^2\big)^{1/2} + c_{0} r^{1+c_0/\c}\).  The difference can be traced to the assumption that
 \(w_{X}\) is a closed 2-form on \(s^{-1}([-L+ 4, L  -4])\).     

 The bound on \(|\gri_{L}|\) can be seen by
writing \(\gri_{L}\) as a sum of four integrals, these denoted by
\(\gri_\psi\), \(\gri_{cs}\), \(\gri_{w}\) and \(\gri_{\mu}\).  The
integrand of \(\gri_\psi\) is \(\frac{1}{4}r |\psi|^{2} \textsc{r}\) with
\(\textsc{r}\) denoting the scalar curvature of \(X\).  By way of an
explanation, this term comes from the integration by parts and
subsequent commuting of covariant derivatives that rewrites the
integral of
\(r|D_{\bbA}\psi|^{2}\) as an integral over the \(s^{-1}(-L + 4)\) and
\(s^{-1}(L  - 4)\) boundaries of the integration domain
plus an integral over \(s^{-1}([-L + 4, L  - 4])\) whose
integrand is the sum of
\(r|\nabla_{\bbA}\psi|^{2}\),
a curvature term involving
\(F_{\bbA}^{+}\) and the product of  
 \(\frac{1}{4}r|\psi|^{2} \textsc{r}\) with
\(\textsc{r}\) denoting the scalar curvature of the metric on \(X\).  The
boundary terms account for the right most integral in
(\ref{eq:(A.5)})'s formula for \(\gra\).  Use the bounds from the first
two bullets of (\ref{eq:(A.16)}) with the bound
\(|\psi|^{2} \leq c_{0}\c\) from Lemma \ref{lem:C.1} to see that \(|i_\psi| \leq c_{0}r^{1+2/\c}\)
  if \(r > c_{\c}\) with \(c_{\c}\) again
denoting a constant that depends only on \(\c\).  
 
The integrals \(i_{cs}\) and \(\gri_{w}\) involve a
chosen Hermitian connection on \(\det (\bbS^{+})\)
whose curvature has norm bounded by \(\c r^{c_0/\c}\) 
 and whose pull back from the \(s \leq -L + 8\) and \(s \geq L  - 8\)
part of \(X\) via the embeddings from the second and third bullets is the
respective \(Y_{-}\) and \(Y_{+}\) versions of \(A_{K } + 2A_{E}\).  Step 3 of
the proof of Lemma \ref{lem:C.1} explains why such connections exist.  Let
\(\bbA_{\bbS}\) denote a chosen connection with this property.  
 
The integral \(\gri_{cs}\) comes by first writing
\(|F_{\bbA}^{+}|^{2}\) as \( \frac{1}{2}|F_{\bbA}|^{2}\) plus the
Hodge star of \(\frac{1}{2}F_{\bbA}\wedge F_{\bbA}\).
 The latter is rewritten using an integration by parts after writing
\(\bbA\) as \(\bbA_{\bbS} +
\hata _{A}\) with \(\hata _{A}\) being an \(i \bbR\)-valued 1-form on \(X\).  Writing \(\bbA\) in this way yields 
\begin{equation}\label{eq:(C.47)}
 \frac{1}{2}F_{\bbA}\wedge F_{\bbA }=\frac{1}{2} d\hata _{A} \wedge d\hata _{A} +
d\hata _{A} \wedge    F_{\bbA_\bbS}
 +\frac{1}{2}F_{\bbA_\bbS} \wedge     F_{\bbA_\bbS}.
\end{equation}
 
An integration by parts writes the integrals of the first two terms on
the right side of (\ref{eq:(C.47)}) as boundary integrals, these giving the
respective cs contributions to \(\gra(\grd|_{s=L+4})\) and \(\gra(\grd |_{s=L-4})\).  The integral of the right most term in (\ref{eq:(C.43)}) is
\(\gri_{cs}\).  Thus \(|\gri_{cs}|\leq c_{0}  r^{c_0/\c}\). 
 
The integral \(\gri_{w}\) is obtained by invoking
Stokes' theorem to rewrite the term from the inner
product between \(F_{\bbA}^{+}\) and \(\frac{i}{2}r w_{X}\) that arises when
\(|F_{\bbA}^{+} - r(\psi^{\dag}\tau\psi  - \frac{i}{2}w_{X}) + i\grw_{\mu}^{+}|^{2}\)
is written as \(|F_{\bbA}^{+}|^{2}+ r |\psi^{\dag}\tau\psi   - w_{X}|^{2}\) plus remainder terms.
 One of these remainder terms is twice the inner product of
\(F_{\bbA}^{+}\) with \(\frac{i}{2} r w_{X}\).  The integral of the latter is the integral of
the 4-form \(-i r F_{\bbA}^{+}\wedge w_{X}\).  Write \(- i rF_{\bbA} \wedge w_{X}\) as the sum
of    \(-i r  d\hata _{A} \wedge w_{X}\) and  
\(- 2 ir F_{\bbA}\wedge  w_{X}\).  Because \(w_{X}\) is closed,
an integration by parts writes the integral of the first of these as an
integral over the boundary of the integration domain.  The latter
accounts for the respective \(\textsc{w}\) contributions \(\gra(\grd|_{s=-L+4})\) and \(\gra(\grd  |_{s=L-4})\).  The integral of \(-2i r F_{\bbA}\wedge  w_{X}^{+}\) is
\(\gri_{w}\).  This being the case, the bound
\( |\gri_{w}| \leq c_{0} \c r^{1+c_0/\c} \) follows directly from the (\ref{eq:(A.16)}) and what is said in Step 3 of the
proof of Lemma \ref{lem:C.1} about \(| F_{\bbA_\bbS}|\).

  The integral denoted by \(\gri_{\mu}\)  has two
contributions.  The first accounts for the terms with
\(\grw_{\mu}\) that arise in the aformentioned rewriting of
\(|F_{\bbA}^{+} - r(\psi^{\dag}\tau\psi  - iw_{X}) + i\grw_{\mu}^{+}|^{2}\).
 It follows from the left hand equation in (\ref{eq:(A.14)}) that the integrand
for this part of \(\gri_{\mu}\)  is bounded by
\(c_{0}\).  The second contribution is proportional to the
integral of \(d\hata _{A} \wedge \grw_{\mu}\);
it appears when Stokes' theorem is used to write the
respective \(\gre_{\mu}\) parts of \(\gra(\grd|_{s=-L+4})\) and \(\gra(\grd  |_{s=L-4})\) as a term that has norm bounded by \(c_{0}\) and another
whose integrand is proportional to \(d\hata _{A}\wedge\grw_{\mu}\).  The norm of the latter is bounded by
\(c_{0 }\c (|F_{\bbA}| +\c^{2})\).  Granted this, it follows that
\(|\gri_{\mu}| \leq
c_{0} \c \Big( \big(\int_{s^{-1}([-L+4, L-4])}|F_\bbA  |^2\big)^{1/2} + \c^{2}\Big)\) and this is guaranteed
by Lemma \ref{lem:C.1} to be less than \(c_{0} \c (r +\c^{2})\). 
 
  \paragraph{Part 4:}  If the first bullet of (\ref{eq:(A.17)}) holds, the
assumption \(\gra(\grc_{-})  - \gra(\grc_{+}) < r^{2-1/\c} \) with (\ref{eq:(C.16)}) and (\ref{eq:(C.43)}) imply that  
\begin{equation}\label{eq:(C.48)}\begin{split}
&\frac{1}{2}\int_{\bbR\times Y_Z} \big(|\frac{\partial  }{\partial s}
A_*|^2+|\grB_\grd|^2+2r(|\frac{\partial}{\partial s}\psi_*|^2+|D_A\psi|^2)\big)\\
  &\quad\leq \gra(\grc_{-}) - \gra(\grc_{+}) +c_{0} r \\
&\quad \leq c_{0}  r^{2-1/\c}
\end{split}
\end{equation}
when \(\c > c_{0}\) and \(r\) is greater than a
constant that depends only on \(\c\).  If the second bullet of (\ref{eq:(A.17)})
holds, the assumption \(\gra(\grc_{-}) - \gra(\grc_{+})< r^{2-1/\c}\) 
with (\ref{eq:(C.16)}), (\ref{eq:(C.44)}) and (\ref{eq:(C.46)}) imply
the bounds that follow when \(\c > c_{0}\) and \(r\) is greater than a constant that
depends only on \(\c\):
\begin{equation}\label{eq:(C.49)}
\begin{cases}
& \int_{(-\infty, -L+4]\times Y_-} \big(|\frac{\partial  }{\partial s}
A_*|^2+|\grB_\grd|^2+2r(|\frac{\partial}{\partial s}\psi_*|^2+|D_A\psi|^2)\big)
  \leq c_{0} r^{2-1/\c}.\\
 & \int_{[L-4, \infty)\times Y_+} \big(|\frac{\partial  }{\partial s}
A_*|^2+|\grB_\grd|^2+2r(|\frac{\partial}{\partial s}\psi_*|^2+|D_A\psi|^2)\big)
  \leq c_{0} r^{2-1/\c}.\\
 & \int_{s^{-1}([-L+4, L-4])} \big(| F_\bbA|^2+ r^2|\psi^\dag\tau\psi-iw_X^+|^2+2r|\nabla_\bbA\psi|^2\big)
  \leq c_{0}r^{2-1/\c}.
\end{cases}
\end{equation}
 Put away for now the bounds in (\ref{eq:(C.48)}) and those in the first two
bullets of (\ref{eq:(C.49)}). Assuming that the second bullet of
(\ref{eq:(A.17)}) holds, the bound in the third bullet of (\ref{eq:(C.49)}) implies the bound
\begin{equation}\label{eq:(C.50)}
r \int_{s^{-1}([-L+4, L-4])}|\psi^\dag\tau\psi-iw_X^+|\leq c_{0} r^{1-1/\c}
\end{equation}
when \(\c> c_{0}\) and \(r\) is greater than a
constant that depends only on \(\c\).  Let \(B\) denote the given ball from
Lemma \ref{lem:C.5}.  Use the second and third bullets of (\ref{eq:(A.12,15a)}) and (\ref{eq:(A.15b)}), the
first bullet of Lemma \ref{lem:C.2}, and (\ref{eq:(C.50)}) to see that
\begin{equation}\label{eq:(C.51)}
r\int_{B\cap s^{-1}([-L+4, L-4])} |1-|\alpha|^2|
 \leq c_{0} r^{1-1/\c},
\end{equation}
when \(r\) is greater than a purely \(\c\)-dependent constant.

\paragraph{Part 5:} If the first bullet of (\ref{eq:(A.17)}) holds, then \(I\) denotes
in what follows any given length 1 interval in \(\bbR\).  If the
second bullet of (\ref{eq:(A.17)}) holds, then \(I\) denotes a length 1 interval in
either \((-\infty, -L + 4]\) or in \((L  - 4, \infty)\).  In either
case, reintroduce the 1-form \(\upsilon_{X}\) from the
fifth bullet of (\ref{eq:(A.17)}).  Take the inner product of both sides of
(\ref{eq:(C.41)}) with \(i \upsilon_{X}\); then integrate the
resulting identity over \(s^{-1}(I)\).  The left hand
side of the result can be written as a sum of four integrals; and the
assertion that this sum is zero can be rewritten as the identity
\begin{equation}\label{eq:(C.52)}\begin{split}
\int_I & \big(\int_{Y_*}\upsilon_X\wedge r(w_*+*i\psi^\dag\tau\psi)\big)
ds  = \int_I\big(\int_{Y_*}\upsilon_X\wedge id\hata_\bbA\big) ds\\
& 
+  \int_I\big(\int_{Y_*}\upsilon_X\wedge *\frac{\partial}{\partial s} A_*\big) ds
+ \int_I\big(\int_{Y_*}\upsilon_X\wedge
*(-\grw_\mu^+(\frac{\partial}{\partial s}, \cdot)+\frac{1}{2}iB_{\bbA_\bbS})\big) ds. 
\end{split}
\end{equation}
 
Use what is said by either the first bullet in (\ref{eq:(A.17)}) or the second and
fifth bullets of (\ref{eq:(A.16)}) to bound the absolute value of the right most
integral in (\ref{eq:(C.52)}) by a purely \(\c\)-dependent constant.  Meanwhile,
Stokes' theorem finds the middle integral on the right
hand side of (\ref{eq:(C.52)}) equal to zero.  The absolute value of the left
most integral on the right hand side of (\ref{eq:(C.52)}) is bounded by
\(c_{0}\c\) times the \(L^{2}\)-norm over
\(s^{-1}(I)\) of \(\ps A_{*}\).  This being the case, use either (\ref{eq:(C.48)}) or the
first two bullets in (\ref{eq:(C.49)}) to bound the absolute value of the left
most integral on the right side of (\ref{eq:(C.52)}) by \(r^{1-1/(2\c)}\)
when \(r\) is greater than a purely \(\c\) -dependent constant.  
 
It follows as a consequence of what was just said in the preceding
paragraph that the absolute value of the integral on the left hand side
of (\ref{eq:(C.52)}) is no greater than \(r^{1-1/(2\c)}\)
when \(r\) is large.  The plan for what follows is to rewrite this
integral as the sum of two terms, one being the integral of
\(r|\upsilon_{X}| \big| |w_{*}| -|\psi|^{2}\big|\) and the other bounded by \(r^{1-1/c_{\c}}\).  This is done in Part 7.  Part 6 supplies the necessary tools.  A
bound of this sort with the second and third bullets of (\ref{eq:(A.12,15a)}) and
(\ref{eq:(A.15b)}) plus the first bullet of Lemma \ref{lem:C.2} leads directly to the bound
 \begin{equation}\label{eq:(C.53)}
  r\int_{B\cap s^{-1}(I)}|1-|\alpha|^2|
 \leq c_{0} r^{1-1/\c}
\end{equation}
when \(B\) is any given ball from Lemma \ref{lem:C.5}.  This bound implies Lemma
\ref{lem:C.5}'s assertion if the first bullet of (\ref{eq:(A.17)}) holds.
 This bound with (\ref{eq:(C.52)}) imply Lemma \ref{lem:C.5}'s bound when
the second bullet of (\ref{eq:(A.17)}) holds.

\paragraph{Part 6:}  The two lemmas that are stated momentarily and
then proved supply what is needed for Part 7.  To set the stage for
the first lemma, note that Clifford multiplication by
\(w_{X}\) splits \(\bbS^{+}\) where
\(w_{X} \neq 0\) as a direct sum of eigenbundles for the
endomorphism given by Clifford multiplication by \(w_{X}\).
 Write this direct splitting as \(\bbS^{+} =E_{X} \oplus (E_{X}\otimes K_{X}^{-1})\) with it
understood that the left most summand is the \(i
|w_{X}|\)-eigenspace.  The upcoming lemma
writes a section \(\psi\) of \(\bbS^{+}\) where
\(w_{X} \neq 0\) as
\(|w_{X}|^{1/2} \eta\)
and it writes \(\eta\) with respect to the direct sum decomposition of
\(\bbS^{+}\) as \((\alpha, \beta)\).   The
lemma that follows asserts bounds for \(|\alpha|\)
and \(|\beta|\) that are the analogs of those
asserted by the first two bullets of Lemma \ref{lem:B.2}.
 
\begin{lemma}\label{lem:C.6}   
There exists \(\kappa  > 100\), and given \(\c \geq \kappa \), there exists
\(\kappa_{c}\)  with the following significance:
 Fix \(r \geq \kappa_{c}\)  and assume that
the metric obey the \((\c, \r = r)\)  version of the constraints in
the first three bullets of (\ref{eq:(A.16)}) and
\(|w_{X}| \leq \c\) , or that the first bullet of (\ref{eq:(A.17)}) holds.  Fix elements
\(\mu_{-}\)  and \(\mu_{+}\) from the respective \(Y_{-}\)  and
\(Y_{+}\)  versions of \(\Omega\) with \(\mathcal{P}\)-norm bounded by 1 and use all of this
data to define the equations in (\ref{eq:(A.14)}).  Let \(\grd = (\bbA,
\psi)\)  denote an instanton solution to these equations. Fix \(m > 1\).  Then
\[
|\alpha|^{2} \leq 1 +\kappa_{c }m^{3} r^{-1+\kappa/\c}\quad 
\text{and} \quad |\beta|^{2}< \kappa m^{3} r^{-1+\kappa/\c} (1  - |\alpha|^{2}) +
\kappa^{3} m^{6} r^{-2+\kappa  /\c}
\]
at the points in \(X\) where \(|w_{X}| >m^{-1}\).
\end{lemma}
\pf The proof is much like that of
the first two bullets in Lemma \ref{lem:B.2} with the only salient difference
being the \(r\)-dependent bounds for the norms of the Riemannian curvature
and the covariant derivatives of \(w_{X}\).  The paragraphs
that follow briefly explain how this \(r\)-dependence is dealt with.
 
The section \(\eta = (\alpha, \beta)\) of \(\bbS^{+}\) obeys an equation of the form
\(\mathcal{D}_{\bbA}\eta +\grR\cdot\eta = 0\) with \(\grR\) being an endomorphism that is
bounded by \(c_{\c} m^{-1}r^{1/\c}\) on $U_{2m}$ .  The Weitzenb\"ock formula for the operator
\((\mathcal{D}_{\bbA} + \grR)^{2}\) leads to an equation for \(\eta\) that has the schematic form
\begin{equation}\label{eq:(C.54)}
\nabla_{\bbA}^{\dag}\nabla_{\bbA}\eta- \frac{1}{2}\op{cl}(F_{\bbA}^{+})\eta +\grR_{1}\cdot\nabla_{\bbA}\eta+ \grR_{0}\cdot\eta = 0,
\end{equation}
where \(|R_{1}| \leq c_{c}m^{-1}\) and \(|R_{0}| \leq c_{c}m^{-2}\).  As in the proof of Lemma \ref{lem:B.2}, introduce \(\q\) to denote the maximum of
\(0\) and \(|\eta|^{2} - 1\).  It follows
from (\ref{eq:(C.54)}) that \(\q\) obeys the inequality
\begin{equation}\label{eq:(C.55)}
d^{\dag}d \q +  rm^{-1}\q \leq c_{\c}m^{-2}   r^{2/\c}
\end{equation}
on \(U_{2m}\) when \(r \geq c_{c}\).  It follows from Lemma \ref{lem:C.1} that \(\q
\leq c_{\c} m\) on the boundary of \(U_{2m}\). This understood, the
comparison principle using the Green's function for
\(d^{\dag}d + rm^{-1}\) can be used to see that \(\q - c_{c}m^{3}r^{-1+2/\c}\)
is no greater than \(c_{\c} m e^{-\sqrt{r/(2m)}}\) on \(U_{2m}\).  
This bound on \(\q\) implies what is said by Lemma \ref{lem:C.6} about
\(|\alpha|^{2}\).
 
 To see about the bound for \(|\beta|^{2}\),  project (\ref{eq:(C.54)}) to
the \(E_{X}\otimes K_{X}^{-1}\) summand of \(\bbS^{+}\) to see that
\(|\beta|^{2}\) obeys a differential inequality on \(U_{2m}\)
 that has the schematic form
\begin{equation}\label{eq:(C.56)}
d^{\dag}d |\beta|^{2} +  rm^{-1}|\beta|^{2} \leq
-2|\nabla\beta|^{2} +c_{k}   r^{-1+c_0/\c}m^{3}|\nabla_{A}\alpha|^{2}+ c_{0} m^{2}r^{c_0/\c}
\end{equation}
when \(r \geq c_{\c}\).  Meanwhile, the projection of
(\ref{eq:(C.54)}) to the \(E_{X}\) summand can be used to see that
\(\w = 1- |\alpha|^{2}\) on \(U_{2m}\) obeys the following analog of any given \(\varepsilon >0\) version of (\ref{eq:(B.5)})
\begin{equation}\label{eq:(C.57)}
d^{\dag}d \w + r m^{-1} \w \geq|\nabla\alpha |^{2} -c_{0} \varepsilon|\nabla\beta|^{2} -c_{0} (1 + \varepsilon^{-1})m^{2}r^{c_0/\c}.
\end{equation}
It follows from (\ref{eq:(C.56)}) and (\ref{eq:(C.57)}) that there are constants
\(z_{1}\) and \(z_{2}\) that are both bounded by
\(c_{\c}\), and there exists an \(\varepsilon >
c_{\c}^{-1}\) such that \(q :=|\beta|^{2} - z_{1} r^{-1+c_0/\c}m^3\w-
z_{2}r^{-2+c_0/\c}m^{6}\) obeys the equation \(d^{\dag}dq + rm^{-1}q \leq 0\) on \(U_{2m}\).  This being the case, a comparison principle argument much like that
used in the preceding paragraph bounds \(q\) by \(c_{c}m
e^{-\sqrt{r/2m}} \) on \(U_{2m}\).  This bound implies Lemma
\ref{lem:C.6}'s assertion about \(|\beta|^{2}\). \epf
 
The next lemma supplies an analog for \(X\) of Lemma \ref{lem:B.3}.

\begin{lemma}\label{lem:C.7}
There exists \(\kappa  >100\) , and given \(\c \geq \kappa \), there exists
\(\kappa_{\c}\) with the following significance:
 Fix \(r \geq \kappa_{\c}\)  and assume that the metric obeys the
 \((\c, \r = r)\)  version of the constraints in the first three
 bullets of (\ref{eq:(A.16)}) and \(|w_{X}| \leq \c\), or
that the first bullet of (\ref{eq:(A.17)}) holds.   Fix elements
\(\mu_{-}\)  and \(\mu_{+}\) from the respective \(Y_{-}\) and
\(Y_{+}\) versions of \(\Omega\) with \(\mathcal{P}\)-norm bounded by 1  and use this data
to define the equations in (\ref{eq:(A.14)}).  Let \(\grd =
(\bbA,\psi)\)  denote an instanton solution to these equations. Fix \(m > 1\).  Then
\(|\psi|^{2} \leq\kappa_{\c} (m^{-1} +c_{\c}r^{-1+\kappa   /\c})\)  at
points in \(X\)  where \(|w_{X}| \leq m^{-1}\).
\end{lemma}

\pf The Weitzenb\"ock formula for \(\mathcal{D}_{\bbA}^{2}\)  was used
in Step 1 of the proof of Lemma \ref{lem:C.1} to write the differential
inequality \(d^{\dag}d|\psi| + r(|\psi|^{2} -|w_{X}| -c_{c }r^{-1+1/\c})|\psi| \leq 0\).  The maximum principle
precludes a local maximum for \(|\psi|^{2}  -m^{-1}  - c_{c}r^{-1+1/\c} \)
on \(X-U_{m}\) and Lemma \ref{lem:C.6} implies that
\(|\psi|^{2} \leq 2(m^{-1} + c_{c} m^{2} r^{-1+c_0/\c})\) on the boundary of \(X-U_{m}\).   \epf

 \paragraph{Part 7:}  Fix \(m > 1\) for the moment and write \((ds \wedge \upsilon_{X})^{+}\) on
\(U_{m}\) as \(q_{X} w_{X} +\b_{X}\) with \(\b_{X}\) being a self-dual 2-form
that obeys \(\b_{X} \wedge w_{X} = 0\).  Note
in this regard that 
\begin{equation}\label{eq:(C.58)}
q_{X}|w_{X}|^{2}  = *(ds \wedge \upsilon_{X} \wedge w_{X})
\end{equation}
with the \(*\) here denoting the Hodge star that is defined by the metric
\(ds^{2} + \grg\) on \(I \times Y_{*}\).
  Granted (\ref{eq:(C.58)}), it follows either from the first bullet of (\ref{eq:(A.17)})
or from the fourth bullet and Item c) of the fifth bullet of (\ref{eq:(A.16)})
that 
\begin{equation}\label{eq:(C.59)}
q_{X}|w_{X}|^{2}\geq -c_{\c}r^{-1/\c}.
\end{equation}
  Noting that \(*(ds \wedge \upsilon_{X}\wedge w_{X})\) is also the \(\grg\)-Hodge star on
\(Y_{*}\) of \(\upsilon_{X} \wedge w_{*}\), the integrand of the \(U_{m}\) part of
the integral on the left hand side of (\ref{eq:(C.52)}) is 
\begin{equation}\label{eq:(C.60)}
r q_{X}|w_{*}|^{2}(1 -|\alpha|^{2} +|\beta|^{2}) +  \grr\quad 
\text{where  \(|r| \leq c_{\c}r |\b| |w_{X}|  |\alpha| |\beta|\).}
\end{equation}
 
Use the bound in (\ref{eq:(C.59)}) and the bounds supplied by Lemma
\ref{lem:C.6} to see that the \(U_{m}\) part of the integral on the
left side of (\ref{eq:(C.52)}) can be written as
\begin{equation}\label{eq:(C.61)}
r\int_{U_m}|q_X| |w_*| \big| |w_*|-|\psi|^2\big| + \gre     \quad 
\text{where \(|\gre| \leq c_{\c} (r^{1-c_0/\c} +  m^{3}r^{c_0/\c})\).}
\end{equation}
 
Meanwhile, it follows from Lemma \ref{lem:C.7} that the contribution to the
integral on the left side of (\ref{eq:(C.52)}) from \(X-U_{m}\) is no
greater than \(c_{\c} (r m^{-1} +m^{2}r^{c_0/\c})\).  Lemma
\ref{lem:C.7} also gives such a bound for the integral of \(|q_{X}|
|w_{*}|  \big| |w_{*} | - |\psi|^{2}\big|\) over the part of \(I \times Y_{*}\) in \(X -U_{m}\).
 Granted these bounds, fix for the moment \(\varepsilon> 0\) but with \(\varepsilon <
c_{0} \c^{-1}\) and take \(m = r^{\varepsilon /\c}\).  Use the just stated bounds and (\ref{eq:(C.61)}) to see that 
\begin{equation}\label{eq:(C.62)}
\int_{I\times Y_*} |q_X|\, |w_*|\, \big| |w_*|-|\psi|^2\big| \leq 
\int_{I\times Y_*} ds \wedge \upsilon_X\wedge r
(w_*+*i\psi^\dag\tau\psi )+  c_{\c} r^{1-\varepsilon  /\c}.
\end{equation}
 
This last bound with what is said at the end of Part 5 implies Lemma
\ref{lem:C.5}. \epf

 \subsection{Proof of Proposition \ref{prop:A.5}}\label{sec:Cg)}

  Fix a smooth, \(r\)-independent metric on \(X\) whose pull-back via the
embeddings from the second and third bullets of (\ref{(A.9a,11)}) restricts to the
\(s < -2\) and \(s > 2\) parts of \(X\) as the product metric \(ds^{2 }+ \grg_{*}\), where
\(\grg_{*}\) denotes the given metric on \(Y_{-}\)
and \(Y_{+}\) as the case may be.  Use \(\grm_{X}\) to
denote this metric.  Use this metric to define the bundles
\(\bbS^{+}\) and \(\bbS^{-}\) over \(X\).  The constructions at the
beginning of Section \ref{sec:Cc)} can be repeated to view the \(Y_{-}\) and \(Y_{+}\)
version of \(\bbS\) as the restrictions to the respective \(s < -1\)
and \(s > 1\) parts of \(X\) of the \(\grm_{X}\) versions of \(\bbS^{+}\) and
\(\bbS^{-}\).  Use this view of these versions of \(\bbS\) to view the \(Y_{-}\) and \(Y_{+}\)
versions \(A_{K} + 2A_{E}\) as a Hermitian
connection on the restriction of the \(\grm_{X}\) version of the
bundle \(\det (\bbS^{+})\) to the \(|s| > 1\) part of \(X\).  This connection has smooth,
\(r\)-independent extensions to the whole of \(X\) as a Hermitian connection on
the \(\grm_{X}\) version of \(\det (\bbS^{+})\).
 Fix such an extension and denote it by \(\bbA_{\bbS}\).  
 
Use the \(s < -1\) and \(s > 1\) isomorphisms between the
\(Y_{-}\) and \(Y_{+}\) versions of \(\bbS\) to
view the corresponding versions of \(\psi_{E}\) as a
section of the \(\grm_{X}\) version of
\(\bbS^{+}\) over the \(|s |> 1\) part of \(X\).  Fix a smooth extension of the latter to
the whole of \(X\) and denote it by \(\psi_{\bbS}\).
 
The metric \(\grm_{X}\) and the pair \(\grd_{\bbS}=(\bbA_{\bbS},\psi_{\bbS})\) defines a version of the operator that appears in (2.61) of \cite{T3}. This operator defines a map from
\(C^{\infty}(X; i T^*X \oplus\bbS^{+})\) to \(C^{\infty}(X; \Lambda^{+} \oplus \bbS^{-} \oplus i \bbR)\).  The
latter defines an unbounded, Fredholm operator between the
\(L^{2}\)-versions of these spaces, and so it has a
corresponding Fredholm index, this denoted in what follows as
\(\imath_{\bbS}\).

Fix \(\c > c_{0}\) so that Proposition \ref{prop:A.4} can be
invoked using \(Y_{-}\) and \(Y_{+}\).  Fix \(r\gg 1\) and pairs \(\mu_{-}\)
and \(\mu_{+}\) from the respective \(Y_{-}\) and \(Y_{+}\) versions
of \(\Omega\) with \(\mathcal{P}\)-norm less than 1; and suppose that \(\grc_{-}\) and
\(\grc_{+}\) are the corresponding solutions to the
\(Y_{-}\) and \(Y_{+}\) versions of (\ref{eq:(A.4)}).  
Let \(\grm\) denote a metric on \(X\)
that obeys (\ref{eq:(A.12,15a)}) and (\ref{eq:(A.15b)}).  Suppose that
\(\grd = (\bbA, \psi)\) is a pair of connection on \(\det (\bbS^{+})\) over \(X\) and section
over \(X\) of \(\bbS^{+}\) with \(s  \to-\infty \) limit \(\grc_{-}\) and \(s \to \infty\)
limit \(\grc_{+}\).  This metric \(\grm\) and \(\grd\) together define a
corresponding version of the operator that appears in (2.61) of
\cite{T3}.  If both \(\grc_{-}\) and \(\grc_{+}\) are
non-degenerate then this operator has an unbounded, Fredholm extension
whose domain and range are the respective spaces of square integrable
sections of \(i T^*X \oplus \bbS^{+}\) and \(i\Lambda^{+} \oplus\bbS^{-} \oplus i \bbR\).  Assume
this to be the case for the moment, and let
\(\imath_{\grd+}\)  denote the corresponding
Fredholm index.  It follows using the excision theorem for the index
(or from what is said in \cite{APS}) that \(\imath_{\bbS}  =\imath_{\grd +} +\grf_{s}(\grc_{-}) -\grf_{s}(\grc_{+})\). 
 
With the preceding understood, write \(\gra(\grc_{-})  -\gra(\grc_{+})\) as 
\begin{equation}\label{eq:(C.63)}
\gra^{\grf}(\grc_{-})-\gra^{\grf}(\grc_{+})  - 2\pi (r  - \pi )(\grf_{s}(\grc_{-}) -\grf_{s}(\grc_{+}))
\end{equation}
and then use the formula in the last paragraph to write
\begin{equation}\label{eq:(C.64)}
\gra(\grc_{-}) -\gra(\grc_{+}) =\gra^{\grf}(\grc_{-})-\gra^{\grf}(\grc_{+}) + 2\pi  (r  - \pi )(\imath_{\grd +} - \imath_{\bbS}).
\end{equation}
 
Since \(\imath_{\bbS}\) is independent of \(r\) and \(\c\),
this last formula proves Proposition \ref{prop:A.5} when both \(\grc_{-}\)
and \(\grc_{+}\) are non-degenerate.

If one or neither is non-degenerate, fix \(\varepsilon >0\) and fix \(\grc_{-}'\) in the set
\(\grN_\varepsilon(\grc_{-})\) from Section \ref{sec:Bf)} that takes
on the supremum in the \(\grc_{-}\) version of (\ref{eq:(B.31)}).  Fix \(c_{+}'\) in
\(\grN_\varepsilon(c_{+})\) with the analogous property.  With \(\grc_{-}'\) and \(\grc_{+}'\) as just
described, choose a pair \(\grd'\) of connection on \(\det (\bbS^{+})\) and section of
\(\bbS^{+}\) with \( s \to -\infty\) limit \(\grc_{-}'\) and \(s \to \infty \) limit
\(\grc_{+}'\).  The metric \(\grm\) with \(\grd'\) define an unbounded, but
now Fredholm version of the operator from (2.62) in \cite{T3} with domain
and range being the respective spaces of square integrable sections of
\(i T^*X \oplus \bbS^{+}\) and \(i\Lambda^{+} \oplus\bbS^{-} \oplus i\bbR\).  Let \(\imath_{\grd'}\) denote the Fredholm index
of this operator.  Define \(\imath_{\grd+}\)  to be \(\imath_{\grd'+}\).  Note that this
definition does not depend on \(\grc_{-}'\), \(\grc_{+}'\) or \(\grd'\).   

The arguments that lead to (\ref{eq:(C.64)}) can be repeated verbatim
to obtain the modified version that has \(\grc_{-}\) replaced by
\(\grc_{-}'\) and \(\grc_{+}\) replaced by \(\grc_{+}'\).   Keeping this in mind, choose \(\grc_{-}'\) so that
\(|\gra(\grc_{-}') -\gra(\grc_{-})| < 1\), and choose \(\grc_{+}'\) so
that  \(|\gra(\grc_{+}') -\gra(\grc_{+})| < 1\).  It follows using (\ref{eq:(B.31)}) that \(|\gra^{\grf}(\grc_{-}') -
\gra^{\grf}(\grc_{-})| < 1\) and \(|\gra^\grf(\grc_{+}')  - \gra^\grf(\grc_{+})|< 1\).  The latter bound with the \((\grc_{-}',
\grc_{+}')\) analog of (\ref{eq:(C.64)}) implies what is asserted by
Proposition \ref{prop:A.5} when the non-degeneracy condition does not
hold for one or both of \(\grc_{-}\) and \(\grc_{+}\). \epf

\section{Constructing 2-forms on cobordisms}\label{sec:D}
\setcounter{equation}{0}

 This section mainly supplies proofs for Propositions \ref{prop:A.6},
\ref{prop:A.8}, \ref{prop:A.10a} and \ref{prop:A.10b}.  The proof of Proposition \ref{prop:A.6} is in Section \ref{sec:Db)}, that of
Proposition \ref{prop:A.8} is in Section \ref{sec:Dd)}, that of
Proposition \ref{prop:A.10a} is in Section
\ref{sec:De)}, and Section \ref{sec:Df)} contains the proof of
Proposition \ref{prop:A.10b}.  
The basic issue in each proof is to construct metrics and closed 2-forms on cobordisms with certain prescribed properties. These constructions occupy most of theses subsections. By way of a look ahead, these constructions are, on the whole, quite intricate. Note that there is little by way of the Seiberg-Witten equations in this section.

A proof of Proposition \ref{prop:ech-compute} is given in Section
\ref{pf:prop1.5}, using notions introduced in Section \ref{sec:De)}.

\subsection{\(\mathrm{Met}_T\) metrics on \(\{Y_k\}_{k\in\{0, \ldots, \G\}}\)}\label{sec:Da)}
 
The eight parts of this section describe a set of preferred metrics
on each \(k\in\{0, \ldots, \G\}\) version of \(Y_k\). These parts also
describe the associated harmonic 2-forms with de
Rham cohomology class that of \(c_{1}(\det (\bbS))\). Let \(Y_*\)
denote \(Y_k\) for any \(k\in \{0, \cdots, \G\}\). As the
\(M_\delta\cup \mathcal{H}_0\) part of \(Y_*\) and \(Y\) are
canonically isomorphic, notions defined on any of them are defined for
others and are denoted by the same notation.
  
\paragraph{Part 1:}   This part of the subsection summarizes various
properties of \(Y_*\) that concern \(\mathcal{H }_0\) and the curve
\(\gamma^{(z_0)}\).   Most of what is said below can be found in Section II.1.  

The handle \(\mathcal{H}_0\)  in \(Y_*\) has coordinates \((u, \theta, \phi)\) with \((\theta, \phi)\) being
the standard spherical coordinates on the 2-sphere and with \(u \in
[-R  - \ln (7\delta_{*}), R +\ln (7\delta_{*})]\). As can be seen in (IV.1.5), the
2-form \(w\) and the 1-form \(\upsilon_{\diamond}\) restrict to this handle as  
\begin{equation}\label{eq:(D.1)}
w  = \sin \theta   \, d\theta \wedge d\phi \quad 
\text{and} \quad \upsilon_{\diamond} = 2\big(\chi_{+}e^{2(|u|-R)} +
\chi_{-}e^{-2(|u|+R)}\big) \, du,
\end{equation}
where \(\chi_{+} = \chi(-u  -\frac{1}{4}R)\) and \(\chi_{-} = \chi(u -\frac{1}{4}R)\).  The curve \(\gamma^{(z_0)}\) intersects $\mathcal{H}_0$ as the \(\theta = 0\) line.
Meanwhile, the \(M_{\delta}\) part of  \(\gamma^{(z_0)}\) has a
tubular neighborhood with coordinates \((t, (\theta, \phi))\) with \(t \in [\delta^2, 3  - \delta^2]\), with
\(\theta \in [0, \theta_{*})\) and with \(\phi\) the affine coordinate on \(\bbR/(2\pi\bbZ)\) .  Here,
\(\theta_{*}\) is positive, smaller than \(\frac{1}{100}\delta_{*}\) but greater \(100
\delta^{3}\).  The 2-form \(w\) here appears as in (\ref{eq:(D.1)})
and \(\upsilon_{\diamond}\) appears as \(dt\).  The
coordinate transition function identifies \(t\) with
\(e^{-2(R-u)}\) near the index 0 critical point and with
\(e^{-2(R+u)}\) near the index 3 critical point.  

Recall the function \(\ff\) on \(M\) that plays a central role in much of
\cite{KLT1}-\cite{KLT4}. This is described in detail in Section II.1.
Recall also the vector field \(v\) in \cite{KLT2} p.19.
 Set \(\varepsilon_{*} =\delta_{*} \sin (\frac{1}{2}\theta_{*})\).  The coordinates just described can be
used to construct a piecewise smooth embedded 2-sphere in the \(\ff
\in [\varepsilon_{*}^{2}, 3 -\varepsilon_{*}^{2}]\) part of \(M_{\delta}\) as follows:  
\BTitem\label{eq:(D.2)}
\item The 2-sphere intersects the complement of the
  radius-\(\delta_{*}\)  coordinate balls about the index 0 and 3
  critical points of $\ff$ as the
cylinder where \(\theta = \frac{1}{2}\theta_{*}\).  
\item The 2-sphere intersects the \(\ir \in(\varepsilon_{*},
  \delta_{*}]\) part of the radius-\(\delta_{*}\) coordinate ball
  centered on the the index 0  and index 3 critical points of $\ff$ as
  the locus where \((\ir, \theta,\phi)\) are such that \(\cos \theta > 0\) 
and \(\ir \sin \theta =  \delta_{*}\sin \, (\frac{1}{2}\theta_{*})\).  
\item The 2-sphere intersects the \(\ir =\varepsilon_{*}\) spheres
centered about the index 0 and index 3 critical points as the locus where
\(\cos \theta \leq 0\).
\item The 2-sphere is tangent to \(v\) on the rest of \(M_\delta\).
\ETitem
 
As can be seen, this embedding is smooth except along the following
loci:  It is \(C^{1}\) on the \(\cos \theta = 0\) circle in the
boundary of the respective radius \(\varepsilon_{*}\) coordinate balls about the index 0
and index 3 critical points of \(f\).  It is only
\(C^{0}\) on the \(\theta = \frac{1}{2}\theta_{*}\) circle in the boundary of the respective
radius \(\delta_{*}\) coordinate balls about the index 0 and index 3 critical points.  

The piecewise smooth embedding just described can be smoothed to any
desired accuracy so that the vector field \(\frac{\partial}{\partial \phi}\)  
along the resulting 2-sphere is everywhere tangent, the vector field \(v\)
along the 2-sphere is tangent everywhere on the \(\ff=\frac{3}{2}\)  
circle but nowhere else, and so that the restriction of \(\ff\) to
this sphere has
just two critical points (both nondegenerate), these at the points with \(\theta = 0\) and
$\theta  = \pi $ on the boundary of the radius
\(\varepsilon_{*}\) coordinate balls about its
respective index 0 and index 3 critical points.

\paragraph{Part 2:}  It proves useful for what follows to be somewhat more
precise about the smoothing of the surface from (\ref{eq:(D.2)}) near
the \(\ff=\frac{3}{2}\) circle.  To this end, introduce first
\(\rho_{*}\) to denote \(\frac{3}{2}- \varepsilon_{*}^{2}\) and
\(\rho_{1 *}=\rho_{*} +  \sqrt {2} \big(1- \cos
(\frac{1}{2}\theta_{*})\big)\).  Return to the \(\ff\in [\delta^2, 3- \delta^2]\) tubular
neighborhood of \(\gamma^{(z_0)}\) with the coordinates \((t,(\theta ,
\phi ))\) as described above.  Replace the coordinate \(\theta\) on a neighborhood of the \(\theta = \frac{1}{2}\theta_{*}\) locus by the function 
\(\hat{\rho}=  \sqrt{2}(1   - \cos \theta)^{1/2}\).  Fix
\(\varepsilon_{1} \in (0,c_{0}^{-1}\varepsilon_{*}^{2})\) and use the
coordinate \(\hat{\rho}\) to define the smoothing of the \(\ff \in
(\frac{3}{2}-\varepsilon_1, \frac{3}{2} + \varepsilon_{1})\) part of the surface defined by
(\ref{eq:(D.2)}) to be the locus where 
\begin{equation}\label{eq:(D.3)}
\hat{\rho  }= \rho_{1*} -\big(\rho_{*}^{2}  - (t -\frac{3}{2})^{2}\big)^{1/2}.
\end{equation}
 
Note that the vector field \(v\) is tangent to the locus defined by (\ref{eq:(D.3)})
only along the \(t = 3/2\) circle, and note that the corresponding lines are
tangent from the inside.   Introduce by way of notation \(S\) to
denote a smoothing as just described of the original piecewise smooth
embedding given by (\ref{eq:(D.2)}). (This is the sphere denoted by
\(S_z\) in \cite{L}, about Equation (6.2).)
 
\paragraph{Part 3:}  Use \((x_{1}, x_{2},x_{3})\) for the Euclidean coordinates on
\(\mathbb{R}^3\).  The function \(\ff\) and the \(\bbR/(2\pi
\bbZ)\)-valued coordinate function \(\phi\) can be used to embed a
neighborhood of \(S\) into \(\bbR^{3}\) 
as the sphere of radius \(\rho_{*}-\frac{3}{2}\varepsilon _*^2\) about
the origin by taking \(x_3=\ff-\frac{3}{2}\) and by 
setting the pair  \((x_{1}, x_{2})\) to equal \(\big((\rho_{*}^2-x_3^2)^{1/2} \cos \phi, (\rho_{*}^2-x_3^2)^{1/2} \sin \phi\big)\). 
Note in this regard 
that the values of \(x_3\) on the image of \(S\) range from \(-\rho_{*}\) to \(\rho_{*}\) because
the values of \(\ff\) on \(S\) range from \(\varepsilon _*^2\) to \(3-\varepsilon _*^2\).

 
This embedding is extended to a neighborhood of \(S\) by exploiting the fact
that the \(|\ff -\frac{3}{2} | > \frac{1}{2}\varepsilon_{1}\) part of \(S\) has a neighborhood with the
following property:  Let \(p\) denote a point in this neighborhood.  Then
\(p\) sits on an integral curve of \(v\) that intersects \(S\), and there is
precisely one such intersection point with distance
\(c_\varepsilon^{-1}\varepsilon_{1}^{3}\) or less from \(p\).  Here, \(c_\varepsilon > 1\) is a constant that
depends on \(\varepsilon_{1}\).  Such a neighborhood
exists because \(v\) is tangent to \(S\) only on the \(\ff=\frac{3}{2}\)
circle in \(S\).  Let \(\mathcal{N}_{1}\) denote this neighborhood.  Given \(p \in \mathcal{N}_{1}\), let
\(\eta(p) \in S\) denote the unique point on the integral curve of
\(v\) through \(p\) with distance less than
\(c_\varepsilon^{-1}\varepsilon_{1}^{3}\) from \(p\).  Associate to \(p\) the point in \(\bbR^{3}\) with the coordinates
\begin{equation}\label{eq:(D.4)}
x_{1}(p) = x_{1}(\eta(p)), \quad x_{2}(p) = x_{2}(\eta(p)),\quad x_{3}(p) = \ff(p)-\frac{3}{2}.
\end{equation}
 
To complete the definition of the embedding, suppose next that \(p\) is a
point near the \(\ff \in (\frac{3}{2}-\varepsilon_1, \frac{3}{2}+
\varepsilon_{1})\) part of \(S\) where the coordinates \((t,
\hat{\rho}, \phi)\) are defined.  Associate to \(p\) the point in \(\mathbb{R}^{3}\)  with the coordinates
\begin{equation}\label{eq:(D.5)}\begin{split}
x_{1}(p) & = |\hat{\rho}(p) - \rho_{1*}| \cos \phi(p),\\
x_{2}(p) & = |\hat{\rho}(p)  -  \rho_{1*}|\sin \phi(p),  \\ 
x_{3}(p) & = t(p) -\frac{3}{2}.
\end{split}\end{equation}
 
Note in particular that if \(p\) is also in \(\mathcal{N}_{1}\),
then it follows from the definition of the function \(\hat{\rho}\)
and the definition of \(\rho_{1}\) that the points
given by (\ref{eq:(D.4)}) and (\ref{eq:(D.5)}) are the same.
 
What is said at the end of the preceding paragraph has the following
implication:  The map from \(\mathcal{N}_{1}\) to
\(\bbR^{3}\) and the map described in the preceding
paragraph together define a smooth, \(\phi\)-equivariant embedding of
a neighborhood of \(S\) into \(\bbR^{3}\) that maps \(S\) to
the radius \(\rho_{*}\) sphere and maps \(v\) to \(\frac{\partial}{\partial x_3}\).  
 
Fix \(\varepsilon > 0\) so that the region in \(\bbR^{3}\) with
\((x_{1}^{2} +x_{2}^{2} +x_{3}^{2})^{1/2} \in(\rho_{*}  -
\varepsilon, \rho_{*} + \varepsilon)\) is in the image of
the embedding of \(\mathcal{N}_{1}\).  By way of notation,
\(\mathcal{N}_\varepsilon\) is used in the subsequent
discussion to denote both this region in
\(\bbR^{3}\) and its inverse image in
\(M_{\delta}\).  It is worth keeping in mind for what
follows that the points in the \(\bbR^{3}\)
incarnation of \(\mathcal{N}_\varepsilon\)  with
distance  {\em greater} than \(\rho_{*}\) from the
origin are in the \(\mathcal{H}_{0}\) component of \(Y-S\).
 
By construction, the 1-form \(\upsilon_{\diamond}\) appears on the \(\bbR^{3}\) version of
\(\mathcal{N}_\varepsilon\)  as \(dx_{3}\).
 Meanwhile, the 2-form \(w\) must appear here as \(\textsc{k}
dx_{1} \wedge dx_{2}\) with \(\textsc{k}\) being a strictly positive function of \(x_{1}^{2} +
x_{2}^{2}\).  This is because \(w\) is closed,
it annihilates \(v\) and \(v\) appears on the \(\bbR^{3}\)
version of \(\mathcal{N}_\varepsilon\) as \(\frac{\partial}{\partial x_3}\). 
 
Use \(\rho\) to denote the function \((x_{1}^{2} +x_{2}^{2} )^{1/2}\) on
\(\bbR^{3}\) and introduce the \(\bbR/(2\pi\bbZ)\)-valued function \(\phi\) by writing \(x_{1}\)
and \(x_{2}\) as \(\rho \cos \phi \) and \(\rho\sin \phi\).  The
observations from the preceding paragraph, the fact that \(w\) is
harmonic and the fact that its metric Hodge dual is \(\upsilon_{\diamond}\) have the following
implication:  The metric from \(M_{\delta}\) appears on the \(\bbR^{3}\) incarnation of
\(\mathcal{N}_\varepsilon  \) as 
\begin{equation}\label{eq:(D.6)}
\grg = \textsc{k}\,  (h^{-2}d\rho^{2} + h^{2}\rho^{2} d\phi^{2}) +dx_{3}^{2}
\end{equation}
with \(h\) denoting a strictly positive function of \(\rho^2\).

 \paragraph{Part 4:}  This part of the subsection says something of the
topological significance of \(S\) and Part 3's embedding of
\(S\) and its neighborhood \(\mathcal{N}_\varepsilon\) in
\(\mathbb{R}^3\) .  To set the stage, recall that
\(Y_{0}\) was obtained from \(M\) by attaching the 1-handle
\(\mathcal{H}_0\) .  This was done by first deleting the
radius \(7\delta_{*}\) coordinate balls about the index
\(0\) and index \(3\) critical points of \(\ff\) to obtain a manifold with
boundary.  The resulting boundary spheres were then glued to the \(u = R
+ \ln (7\delta_{*})\) and \(u = -R -\ln (7\delta_{*})\) boundary spheres of \([-R -
\ln (7\delta_{*}), R +\ln (7\delta_{*})] \times S^{2}\). 
 
The sphere \(S\) enters a second description of \(Y_{0}\) as the
connected sum of \(M\) with the manifold \(S^{1} \times
S^{2}\).  (Cf. \cite{L}, (6.2)). The connected sum description constructs
\(Y_{0}\) by deleting the respective 3-balls from \(M\) and
\(S^{1} \times S^{2}\) and gluing the
resulting two boundary spheres to the boundary spheres of the product
of an interval with \(S^{2}\).   Denote this product as \(I\times S^{2}\) with \(I \subset \bbR \) being
an interval.  As explained below, the surface \(S\) can be viewed as a
cross-sectional sphere of \(I \times S^{2}\).
 
To see directly this connected sum depiction of \(Y_{0}\), first
view \(S\) and \(\mathcal{N}_\varepsilon \) as subsets in
\(\bbR^{3}\).  Let \(\ir =(\rho^{2} +x_{3}^{2})^{1/2}\) denote the
radial coordinate on \(\bbR^{3}\).  The connected sum
picture of \(Y_{0}\) results in an embedding of \(I \times S^{2}\) into \(\bbR^{3}\) whose image
is the \(\ir \in [\rho_{*} -\frac{1}{16}\varepsilon,  \rho_{*} +  \frac{1}{16}\varepsilon]\) part of \(\mathcal{N}_\varepsilon\).
 This depiction of \(I \times S^{2}\) in \(Y_{0}\) identifies the \(\ir
 = \rho_{*} + \frac{1}{16}\varepsilon \) sphere in \(\mathcal{N}_\varepsilon \) with the
boundary of the complement of a ball in \(S^{1} \times
S^{2}\).  This missing ball can be identified with the \(\ir
< \rho_{*} + \frac{1}{16}\varepsilon \) part of \(\bbR^{3}\).  Indeed, the
\(Y_{0}\) incarnation of the \(\ir =\rho_{*} +\frac{1}{16}\varepsilon \) sphere in \(\bbR^{3}\) splits
\(Y_{0}\) into two components.  The component that contains
the \(\ir > \rho_{*} +\frac{1}{16}\varepsilon \) part of \(\mathcal{N}_\varepsilon\)  is
the complement of a ball in \(S^{1} \times S^{2}\); and \(S^{1} \times
S^{2}\) is reconstituted in full when this complement is
filled in by adding the \(\ir \leq \rho_{*} + \frac{1}{16}\varepsilon \) part of \(\bbR^{3}\) to the \(\ir
> \rho_{*} + \frac{1}{16}\varepsilon \) incarnation of \(\mathcal{N}_\varepsilon\).
 
The \(Y_{0}\) incarnation of the \(\ir =
\rho_{*} -\frac{1}{16}\varepsilon \) sphere in \(\bbR^{3}\) also
separates \(Y_{0}\) into two components.  The component that
has the \(\ir < \rho_{*} -\frac{1}{16} \varepsilon \) part of  \(\mathcal{N}_\varepsilon\)
is the complement of a ball in \(M\).  This ball is attached to give back
\(M\) by viewing the complement of its center point as the \(\ir >
\rho_{*}- \frac{1}{16}\varepsilon_{*}\) part of \(\bbR^{3}\).  To see this, take a second copy of
\(\bbR^{3}\) and use \(\ir'\) to denote the distance to
the origin in the latter.  Use \((\theta', \varphi')\) to denote
the associated spherical coordinates.  The manifold \(M\) is obtained by
attaching the \(\ir ' \leq (\rho_{*} -\frac{1}{16}\varepsilon)^{-1}\) ball in this second copy of
\(\bbR^{3}\) to the \(\ir =\rho_{*}-\frac{1}{16}\varepsilon \) sphere in the original copy of
\(\bbR^{3}\) via the identifications \(\ir' =
\ir^{-1}\) and \((\theta '= \pi   - \theta, \phi '= \phi)\). 
 
Since \(S\) splits \(Y_{0}\) into two parts, it likewise splits \(Y_*\)
into two parts.  The component of \(Y_*-S\) that contains \(\gamma^{(z_0)}\)
has its canonical identification with the \(\gamma^{(z_0)}\) component of
\(Y_{0}-S\). The other component of \(Y_* - S\) is
obtained from the complementary component of
\(Y_{0}-S\) by attaching the  \(\grp  \in \Lambda  \) labeled 1-handles.
 
Both \(Y_* - \mathcal{N}_\varepsilon \) and \(Y_{0} -\mathcal{N}_\varepsilon \)
likewise have two components because \(\mathcal{N}_\varepsilon\) 
is a tubular neighborhood of \(S\).  
A given \(k\in\{0, \ldots, \G\}\) version of \(Y_k\) is obtained from \(Y_0\) by attaching \(k\) 1-handles with attaching regions that are disjoint from the component of \(Y_0-\mathcal{N}_\varepsilon\) that contains \(\gamma^{(z_0)}\). This understood, \(\mathcal{N}_\varepsilon\) can be viewed as a subset of \(Y_k\) and \(Y_k-\mathcal{N}_\varepsilon\) also has two components.
By way of notation, the component of \(Y_*-\mathcal{N}_\varepsilon\) or
any given \(k\in\{0, \ldots, \G\}\) version of \(Y_k-\mathcal{N}_\varepsilon\)
that contains \(\gamma^{(z_0)}\) is denoted in what follows by \(\mathcal{Y}_{0}\) 
and the other component is denoted by \(\mathcal{Y}_{M}\).
(\(\mathcal{Y}_M\) has a natural interpretation as a sutured manifold,
which is denoted by \(M(1)\) in Remark \ref{rmk:suture}).

\paragraph{Part 5:}   This part of the subsection introduces a
family of distinguished metrics on the \(k\in\{0, \ldots, \G\}\) version
of \(Y_k\) that play
central roles in the subsequent discussions.  Parts 6 and 8 say more
about this set.
 
This distinguished set of metrics is parametrized by a
parameter \(T\) which is in all cases greater than 1.  With \(T\) chosen, the
corresponding set of metric is denoted in what follows by
\(\op{Met}_{T}\) .  The metrics from \(\op{Met}_{T}\) 
are constructed momentarily from the set of metrics on
 \(\mathcal{Y}_{M} \cup \mathcal{N}_\varepsilon \) that are given by (\ref{eq:(D.6)}) on
 \(\mathcal{N}_\varepsilon \). This set of metrics on \(\mathcal{Y}_{M} \cup \mathcal{N}_\varepsilon \)
is denoted by \(\op{Met} ^{\mathcal{N}}\).
 Note with regards to (\ref{eq:(D.6)}) that its formula depicts a 1-parameter
family of metrics with the parameter being the length of the curve
\(\gamma^{(z_0)}\).  The length of \(\gamma^{(z_0)}\) plays no role of
significance. In
any event, the length is assumed to be the same for all metrics in
\(\op{Met} ^{\mathcal{N}}\) whether defined on \(Y\) or on a
\(k\in\{0, \ldots, \G\}\) version of \(Y_k\).

The criteria for membership in \(\op{Met}_{T}\) follow directly:  All metrics in
\(\op{Met}_{T}\)   agree on \(\mathcal{Y} _{0} \cup \mathcal{N}_\varepsilon \);
the metric they define on this set is denoted in what follows by \(\grg_{T}\).
The metric \(\grg_{T}\) on \(\mathcal{Y}_0\) is the metric from
(\ref{eq:(A.3b)}). 
Meanwhile, the metric \(\grg_{T}\)  
on \(\mathcal{N}_\varepsilon\) is defined in the three steps that follow.

 \paragraph{\it Step 1:} Introduce \(\chi_{\r}\)  to denote the function on  \(\mathbb{R}^3\)  given by
\(\chi(64\varepsilon^{-1}(\ir  -\rho_{*})-1)\).  This function equals \(1\) where
\(\ir <  \rho_{*} + \frac{1}{64}\varepsilon \) and equals \(0\) where \(\ir >
\rho_{*} + \frac{1}{32}\varepsilon\). Fix \(T> 1\)  and introduce \(\r_{T}\) 
to denote \((1  - \chi_{\r} +\frac{1}{T}\chi_{\r}) \r\) .  The \(\r\) derivative of \(\r_{T}\)
is strictly positive because that of \(\chi_{\r}\) is non-positive.  Set \(\rho_{T} = \r_{T}
\sin \theta \) and \(x_{3T} = \r_{T} \cos \theta\).  Noting that \(d\rho_{T}\) and
\(dx_{3T}\) are linearly independent, the quadratic form
\begin{equation}\label{eq:(D.7)}
\textsc{k}(\rho_{ T })\, (h^{-2}(\rho_{T})\, d\rho_{T}^{2} +
h^{2}(\rho_{T})\rho_{T}^{2}d\phi^{2}) + dx_{3T}^{2}
\end{equation}
defines a smooth metric on \(\bbR^{3}\).  The metric
\(\grg_{T}\) on the \(\ir >\rho_{*} -\frac{1}{4} \varepsilon \) part
of  \(\mathcal{N}_\varepsilon\) is given by (\ref{eq:(D.7)}).  
 
\paragraph{\it Step 2:}  The definition of \(\grg_{T}\) on the \(\ir  \in
[\rho_{*}-\frac{1}{2}\varepsilon, \rho_{*}-\frac{1}{4}\varepsilon]\) part of  \(\mathcal{N}_\varepsilon\)
requires yet another function of \(r\).  This one is defined by the rule \(\r
\mapsto\chi(4\varepsilon^{-1}(\ir  -\rho_{*}) + 2)\) and it is
denoted by \(\chi_{\ir *}\).  The function
\(\chi_{\ir *}\)  is equal to 1 where \(\r
< \rho_{*}-\frac{1}{2}\varepsilon \) and it is equal to 0 where \(\ir >
\rho_{*}-\frac{1}{4}\varepsilon\).  Set \(x_{3T*}\) to denote
the function \(\big(1  - \chi_{\ir  *} +\frac{1}{T} \chi_{\ir  *}\big)
\, x_{3}\). Introduce by way of notation
\(\textsc{k}_T\) and \(h_{T}\) to denote the functions \(\textsc{k}(
\rho/T)\) and \(h(\rho/T)\).  Noting that \(dx_{1}\), \(dx_{2}\)
and \(dx_{3T}\) are linearly independent, the quadratic form
\begin{equation}\label{eq:(D.8)}
\frac{1}{T^2}   \textsc{k}_T(h_{T}^{-2}d\rho^{2} +h_{T}^2\rho ^2d\phi^{2})  +  \frac{1}{T^2}dx_{3T *}^{2} 
\end{equation}
defines a smooth metric on the \(\ir \in [\rho_{*}
 - \frac{1}{2}\varepsilon, \rho_{*}  - \frac{1}{4}\varepsilon]\) part of  \(\mathcal{N}_\varepsilon\).
 The latter extends the metric given in (\ref{eq:(D.7)}) because
\(\rho_{T} = \frac{1}{T}\rho \) and \(x_{3T} =\frac{1}{T}x_{3}\) where \(\rho <\rho_{*} + \frac{1}{64}\varepsilon\).

\paragraph{\it Step 3:}  The definition of \(\grg_{T}\) on the \(\ir <\rho_{*}-\frac{1}{2}\)
 part of  \(\mathcal{N}_\varepsilon\)  requires one more function of
 \(r\).  This one is denoted by \(\chi_{\ir **}\)   and it is defined by the
rule \(\ir \mapsto\chi(4\varepsilon^{-1}(\ir  -\rho_{*}) + 3)\).  This function is equal to 0
where \(\ir > \rho_{*} -\frac{1}{2}\varepsilon \) and it is equal to 1 where \(\ir <
\rho_{*}  -\frac{3}{4}\varepsilon\).  With this function in hand, define the function
\(T_{*}\) to be \(T (1-\chi_{\ir  ** }) +\chi_{\ir  ** }\).  The function
\(T_{*}\) is equal to \(T\) where \(\ir >
\rho_{*} -\frac{1}{2} \varepsilon \) and it is equal to 1 where \(\ir <
\rho_{*}  - \frac{3}{4}\varepsilon\).  The metric \(\grg_{T}\) is defined on the \(\ir
\leq \rho_{*}-\frac{1}{2} \varepsilon \) part of
\(\mathcal{N}_\varepsilon\) to be the quadratic form 
\begin{equation}\label{eq:(D.9)}
\frac{1}{T_*^2} \textsc{k}_{T_*}  (h_{T_*}^{-2}d\rho^{2} +h_{T_*}^2 d\phi^{2})  + 
\frac{1}{T_*^4} dx_{3}^{2}.  
\end{equation}
This definition of \(\grg_{T}\) smoothly extends the metric
defined in (\ref{eq:(D.8)}).  Moreover, the metric \(\grg_{T}\) as just
defined is the metric in (\ref{eq:(D.6)}) where \(\ir < \rho_{*}  - \frac{3}{4} \varepsilon\).

\paragraph{Part 6:} This part of the subsection and Part 8 point out some
key properties of the \(\op{Met}_{T}\) metrics.  This part focuses
on the metric \(\grg_{T}\), this being the restriction of each
\(\op{Met}_{T}\) metric to \(\mathcal{Y}_{0} \cup \mathcal{N}_\varepsilon \).  As explained in the
subsequent two paragraphs, each \(T > 1\) version of
\(\grg_{T}\) on the complement in \(\mathcal{Y}_{0}\cup  \mathcal{N}_\varepsilon \) of the \(\ir \leq
\rho_{*}\) part of \(\mathcal{N}_\varepsilon\) can be viewed as the
pull-back of a \(T\)-independent metric on \(S^{1} \times
S^{2}\) by a \(T\)-dependent embedding of the \(\gamma^{(z_0)}\)
component of \(Y_*-S\) or \(Y-S\)  as the case may be.
 The embedding is denoted by \(\Phi_{T}\).
 
To define this \(T\)-independent metric on \(S^{1}\times S^{2}\), view
\(S^{1}\times S^{2}\)  as in Part 4.  By way of a reminder, this
view comes with a distinguished ball with a distinguished
diffeomorphism onto the \(\ir <\rho_{*} + \frac{1}{16}\varepsilon  \) ball in
 \(\bbR^{3}\) centered on the origin.  There is in addition, a distinguished
identification between the complement of the concentric \(\ir  \leq\rho_{*} \)
ball in \(S^{1}\times S^{2}\)  and the union of \(\mathcal{Y}_{0}\) and the \(\ir \geq\rho_{*}\)
part of \(\mathcal{N}_\varepsilon \). The latter identifies the metric from Section 1 on
 \(\mathcal{Y}_{0}\) with a metric on \(S^{1} \times S^{2}\) whose
restriction to the  \(\ir  \leq\rho_{*}+ \frac{1}{16}\varepsilon  \) ball in the
distinguished coordinate chart appears as
\(\textsc{k} (\rho ) (h ^{-2} ( \rho ) d \rho ^{2} +h ^{2} ( \rho ) \rho ^{2} d \phi ^{2})+dx_{3} ^{2}\).
 This is the desired \(T\)-independent metric on
\(S^{1} \times S^{2}\).  This \(S^{1}\times S^{2}\)  metric is denoted by \(\grg_{*}\).
 
 Fix \(T \geq 1\).  The promised embedding of the \(\mathcal{Y}_0\) component of
\(Y_* - S\) into \(S^{1}\times S^{2}\)  is defined
as follows:  This embedding agrees with the embedding from the
preceding paragraph on \(\mathcal{Y}_{0}\) 
and on the  \(\ir  >\rho_{*}+\frac{1}{32}\varepsilon  \) part of \(\mathcal{N}_\varepsilon\).
 Meanwhile, the promised embedding on the
 \(\ir \in (\rho_{*}, \rho_{*} +\frac{1}{16} \varepsilon )\) part
 of \(\mathcal{N}_\varepsilon\) maps the latter onto the \(\ir 
 \in (T ^{-1} \rho_{*}, \rho_{*} +\frac{1}{16}\varepsilon )\) ball in the
distinguished coordinate chart.  The map here sends the point with
spherical coordinates \((\ir , \theta , \phi )\) to that with the
spherical coordinates \(( \ir_{T}, \theta , \phi )\).

\paragraph{Part 7:}  This part of the subsection describes a certain
closed 2-form on a given \(k\in \{0, \ldots, \G\}\) version of \(Y_k\) with compact support in
\(\mathcal{Y}_M\)  and with the following additional
property:  The de Rham class of this 2-form annihilates all but the
\(H_{2}(M; \bbZ)\) summand in the Mayer-Vietoris
direct sum decomposition for \(H_{2}(Y; \bbZ)\) in
(IV.14) or in the analogous direct sum decomposition for
\(H_{2}(Y_k; \bbZ)\).  Meanwhile, it acts on the \(H_{2}(M; \bbZ)\) summand as
\(c_{1}(\det (\bbS))\).  A version of this 2-form is also defined on
\(M\). In all cases, the 2-form is denoted by \(\p\).
 It is used in the upcoming Lemma \ref{Lemma D.1} and in later
subsections.   The construction of \(\p\) follows directly.
 
View \(M_\delta\) as being a subset of each \(k\in\{0, \ldots, \G\}\)
version of \(Y_k\). As such, it sits in the \(\mathcal{Y}_M\) part of \(Y_k\).
It follows from the description of \(H_{2} (Y; \bbZ )\) in Part 4 of
Section II.1c that there exists a finite set of the form \(\Theta\) whose elements
are pairs of the form \(( \gamma ,\textsc{z}_\gamma )\), with
\(\gamma\)  being a loop in a level set of \(M_{\delta}\)  of
the function \(\ff\) on \(M\). 
Meanwhile,
\(\textsc{z}_\gamma\) is an integer.  The loops from \(\Theta\) generate the
image in any given \(k\in\{0, \ldots, \G\}\) version of \(H_{1} (Y_k; \bbZ )/\op{tors}\) of \(H_{1} (M;\bbZ )/\op{tors}\) via the
Mayer-Vietoris homomorphism for the \(Y_k\) analog
of the direct sum decomposition in (IV.1.4).  Meanwhile, the paired integers are such
that \(\sum_{ \gamma\in\Theta }\textsc{z}_\gamma \gamma  \) represents the image of the Poincar\'e
dual of the restriction of \(c_{1} (\det ( \bbS ))\) to the \(H_{2}
(M; \bbZ )\) summand in this same direct sum decomposition.
  Let \((\gamma , \textsc{z}_\gamma )\) denote a pair from
\(\Theta\). The loop \(\gamma\) has a tubular neighborhood in
 \(M_{\delta}\) which is the image via an embedding of
\(S^{1} \times  D\) where \(D\subset  \bbR ^{2}\) is a small radius disk about the origin and where
 \(\gamma\) corresponds to the image of \(S^{1} \times  \{0\}\).   Use \(\mathcal{T}_{\gamma}\) 
in what follows to denote a tubular neighborhood of this sort.  These are to be chosen so that the pairwise distinct versions have disjoint closure that is disjoint from the boundary of the closure of the \(M_\delta\) part of \(\mathcal{N}_\varepsilon\).

 Note that there exists such a tubular neighborhood with an
embedding that has the following property:  The pull back of \(d\ff\) via
the embedding is a constant 1-form from the \(D\) factor of
\(S^{1} \times  D\) and the kernel of the pull back via the embedding of the 2-form
 \(w\)  is a constant vector field that
is tangent to this \(D\) factor.  The existence of such an embedding
follows from two facts, the first being that \(\gamma\) is in an \(\ff\)-level
set.  The second fact follows from the definition in the first bullet
of (IV.1.3) of \(w\) on \(\mathcal{T}_{\gamma}\) 
as the area form for the \(\ff\)-level sets.  An embedding of this sort
is 
used in Part 7 of the upcoming Section \ref{sec:De)}.

Fix a compactly supported 2-form on \(D\) whose integral is equal to 1.  View
this 2-form first as an \(S^{1}\)-independent form on \(S^{1}\times
D\) and then as a 2-form on \(M\) and on each \(k\in \{0, \ldots,
\G\}\) versions of \(Y_k\) with compact support in \(\mathcal{T}_{\gamma}\).
 Use \(\p_{\gamma}\) to denote the latter incarnation; then set
 \(\p=\sum_{(\gamma, \textsc{z}_\gamma )\in \Theta}\textsc{z}_\gamma \p_{\gamma}\).  By construction, the de Rham class of
 \(\p\)  agrees with \(c_{1} (\det ( \bbS ))\) on the \(H_{2} (M; \bbZ )\) summand of the
Mayer-Vietoris direct sum decomposition of \(H_{2} (Y; \bbZ )\) in (IV.1.4) or its
analog for \(H_{2} (Y_{0};\bbZ )\) as the case may be. The de Rham class of  \(\p\)  also
annihilates the \(H_{2} (\mathcal{H}_{0}; \bbZ )\)-summand in these
direct sum decompositions.  In the case of \(H_{2} (Y; \bbZ )\), the de Rham class
of \(\p\)  also annihilates the \(\bigoplus_{\grp  \in{\Lambda }} H_{2} (\mathcal{H}_{\grp}; \bbZ )\)-summand in
(IV.1.4). 
 
\paragraph{Part 8:} Fix \(k\in\{0, \ldots, \G\}\). Given \(T> 1\) and
a metric from \(\op{Met}_{T}\) on \(Y_k\), the next lemma uses 
\(w_{T}\) to denote the associated harmonic 2-form on
\(Y_k\) whose de Rham cohomology class is that of \(c_{1}(\det (\bbS))\).  
 
\begin{lemma}\label{Lemma D.1}   
There exists \(\kappa>1\)  with the following significance:  Fix a
metric from \(Y_k\)'s version of \(\op{Met}_{T}\)  so as to define \(w_{T}\).  Let \(\|
\p\|_{2}\)  denote the metric \(L^{2}\)-norm of \(\p\), and let \(w\)
be the closed 2-form from (\ref{eq:(A.3a)}). Then the \(L^{2}\)-norm of
\(w_{T}\)  is at most \(\kappa (1 + \| \p \|_{2})\) and the
\(C^{1}\)-norm of \(w_{T}-w\)  on \(\mathcal{Y}_{0}\) and on the \(\ir
> \rho_{*} + \frac{1}{2}\varepsilon   \) part of
\(\mathcal{N}_\varepsilon\)   is at most \(\kappa T^{-1/2}\).
\end{lemma}
 
\pf The proof has four steps. 
 
\paragraph{\it Step 1:} The \(L^{2}\)-norm of \(w_{T}\) as defined by the metric from
\(\op{Met}_{T}\)  on \(Y_k\) is greater than \(c_{0} ^{-1} \)
because the integral of \(w_{T}\) over \(\mathcal{H}_{0}\) must be
greater than \(c_{0}^{-1}\) 
so as to have integral 2 on each cross-sectional 2-sphere.  As
explained directly, the \(L^{2}\)-norm of \(w_{T}\) is also less than \(c_{0 }(1 + \| \p \|_{2})\).
The proof that this is so uses the fact that a
given harmonic form minimizes the
\(L^{2}\)-norm amongst all closed forms in its de Rham cohomology
class.  To obtain such a form, reintroduce the coordinates \((t,z)\) for
\(U_{\gamma}\) and let \(\textsc{b}\)  denote a smooth
function with compact support centered on the origin in \(\bbC \) and
with integral 2. Choose a \(T\)-independent version of \(\textsc{b}\)  so that its
incarnation as a function on \(U_{\gamma}\)  has support in \(U_{\gamma}  \cap \mathcal{H}_{0}\).
 With \(\textsc{b}\) chosen, set \(\p_{0}\) to denote
 \(\frac{i}{2}\textsc{b}\, dz  \wedge d\bar{z}\).  This is a closed, compactly supported 2-form in \(\mathcal{Y}_{0}\) 
whose de Rham cohomology class when viewed in either
\(H^{2} (Y_k; \bbZ )\) has pairing zero with all but the
\(H_{2} (\mathcal{H}_{0};\bbZ )\)-summand in the \(Y_k\) version of (IV.1.4).
By construction, the de Rham
cohomology class of \(\p_{\bbS} =\p_{0} + \p\)  is that of \(c_{1} (\det ( \bbS ))\).
 The metric \(L^{2}\)-norm of \(\p_{\bbS}\) is less than \(c_{0} \,( 1 +\| \p\|_{2})\).

\paragraph{\it Step 2:}  Use \(\sigma\) to denote the function on \(\gamma^{(z_0)}\)'s component of
\(Y_k - S\) that equals 1 on \(\gamma^{(z_0)}\)'s component of
\(Y_k - \mathcal{N}_\varepsilon \) and is given near \(S\) by the function on the
\(\ir  \geq \rho_{*}\) part of \(\bbR^{3}\) by the radial function \(\ir \mapsto\chi (2 -128 \varepsilon ^{-1} ( \ir - \rho_{*} ))\).
The function \(\sigma\)  is equal to 1 where \(\r
>\rho_{*}+\frac{1}{64}\varepsilon  \) and it is equal to 0 where \(\ir  <\rho_{*}+ \frac{1}{128}\varepsilon \).  
 
Use \(\e_{T}\) to denote the \(\Phi_{T}^{-1}\)  -pull-back to \(S^{1}
\times S^{2}\)  of the 2-form \(\sigma w_{T} \).
 This 2-form is supported on the complement in
\(S^{1}\times S^{2}\)  of the \(\ir  < \frac{1}{T}( \rho_{*} +\frac{1}{128}\varepsilon )\) part of the
distinguished coordinate ball.   It follows from what is said in Step
1 that the \(L^{2}\)-norm of \(\e_{T}\) is bounded from below by \(c_{0}^{-1}\) 
and bounded from above by \(c_{0}\). 

Use  \(*\) to denote the \(\grg_{*}\)-Hodge dual on \(S^{1}\times S^{2}\).  Note that
\(d \e_{T} \) and \(d * \e_{T}\) are equal to zero on the complement of the
\(\ir  \leq \frac{1}{T}(\rho_{*} + \frac{1}{64}\varepsilon )\) part of the
distinguished coordinate chart.  Meanwhile, the norms of both are
bounded by \(c_{0}T | ( \Phi_{T} ^{-1} )*w_{T} |_{\grg_*}\) on this same ball.  This observation, the fact that the \(\grg_{T}\)-metric is the \(\Phi_{T}\)-pullback of \(\grg_{*}\) 
and the fact that the \(\grg_{*}\)-volume of the \(\ir  \leq
\frac{1}{T}( \rho_{*}+ \frac{1}{64}\varepsilon )\) coordinate ball is bounded by \(c_{0}T ^{-3} \)
implies that the \(L^{1}\)-norm of both \(d \e_{T}\) and \(d * \e_{T}\) 
is bounded by \(c_{0} (1 + \| \p\|_{2})T ^{-1/2}\).

 \paragraph{\it Step 3:}  The 2-form \(w\)  appears in the \(\r
 \geq\rho_{*}\) part of the \(\bbR^{3}\) incarnation of \(\mathcal{N}_\varepsilon\) 
as \(\textsc{k} ( \rho ) \rho  d \rho \wedge d \phi \). The latter form extends smoothly to the
 \(\ir  \leq\rho_{*}\) part of \(\bbR^{3}\) as a
 \(\grg_{*}\)-harmonic 2-form.  It follows as a consequence that
 \(w\)'s restriction to \(\mathcal{Y}_{0}\) 
and to the  \(\ir  \geq \rho_{*} + \frac{1}{32}\varepsilon  \) part of \(\mathcal{N}_\varepsilon\) 
is the pull-back by all \(\Phi_{T}\) of the \(\grg_{*}\)-harmonic
2-form on \(S^{1}\times S^{2}\) whose de Rham class has pairing equal
to 2 with the generator of \(H_{2} (S^{1}\times S^{2} ; \bbZ )\).
This corresponding form on
\(S^{1}\times S^{2}\)  is \(\frac{1}{2\pi}\sin \theta d \theta  \wedge d \phi  \) and also denoted by \(w\). 
 
\paragraph{\it Step 4:}  Introduce the operator
 \(\grD_{*}=  * d + d * \) on \(S^{1}\times S^{2}\) and use it to
 write the 2-form \(\e_{T}\)  as \((1+\grz_{T} ) w  + u_{T} \) with \(\grz_{T}\) 
denoting a constant with norm bounded by \(c_{0} T ^{-3/2}\) and with
 \(u_{T}\) denoting a 2-form which is \(L^{2}\)-orthogonal to \(w\)
 and such that \(\grD u_{T} = \grD \e_{T}\).
  As the Green's function kernel for  \(\grD\)  is smooth on the complement of
the diagonal in \(\times_{2} (S^{1}\times S^{2})\), the fact that \(\grD \e_{T} \)
has support where \( \ir  < \frac{1}{T}( \rho_{*} +\frac{1}{64}\varepsilon )\) and the
\(c_{0} (1 + \| \p \|_{2})T ^{-1/2}\) bound on its \(L^{1}\)-norm implies that
\( | {u_{T} |+ |\nabla u_{T}} | \leq c_{0} (1+ \| \p\|_{2})T ^{-1/2}\) on \(\mathcal{Y}_0\)   
and also on the  \(\ir  > \rho_{*} +\frac{1}{2}\varepsilon  \) part of \(\mathcal{N}_\varepsilon\). \epf

\subsection{Proof of Proposition \ref{prop:A.6}}\label{sec:Db)}
The three parts of
this subsection prove the assertion made by Proposition
\ref{prop:A.6}. 
 
\paragraph{Part 1:}  Let \(Y_{Z}\) denote a given compact,
oriented 3-manifold and let \(\textsc{z}\) denote a non-zero class in
\(H^{2}(Y_{Z}; \bbZ)/\op{tors}\).  Hodge
theory associates to each metric on \(Y_{Z}\) a harmonic
2-form whose de Rham cohomology class is \(\textsc{z}\).  Of specific
interest in what follows are metrics whose associated harmonic 2-form
has transverse zeros.  There is a residual set of metrics on
\(Y_{Z}\) with this property, see for example \cite{Ho} for a proof. 

Fix \(k\in \{0, \ldots, \G\}\). Let \(\grg_{\mathcal{N}}\)
denote a metric in the \(Y_k\) version of \(\op{Met}^{ \mathcal{N}}\).  Fix \(T > 1\) and use
\(\grg_{\mathcal{N}}\) to define a metric in \(\op{Met}_{T}\), this denoted by \(\grg_{1}\).  Let
\(w_{1}\) denote the associated harmonic 2-form with de Rham
cohomology class \(c_{1}(\det (\mathbb{S}))\).  If
\(w_{1}\) has degenerate zeros, fix a second metric,
\(\grg_{2}\), on \(Y_k\) with the following
properties:  Let \(w_{2}\) denote the corresponding
\(\grg_{2}\) harmonic 2-form.  Then \(w_{2}\) has
non-degenerate zeros, and the \(\grg_{1}\) norms of
\(w_{2}  - w_{1}\) and \(\grg_{2}-\grg_{1}\), and those of their \(\grg_{1}\)-covariant
derivatives to order 100 are less than \(T^{-1}\). If \(w_1\) has
nondegenerate zeros, take \(\grg_2=\grg_1\).

\paragraph{Part 2:}  Write \(w_{2}\) on \(\mathcal{Y}_{0}\) and on the \(\ir
>\rho_{*} +\frac{1}{2} \varepsilon  \) part of
 \(\mathcal{N}_\varepsilon\) as \(w  +u_{2}\). By Lemma \ref{Lemma D.1}, 
 the 2-form is such that \(|u_{2} |\leq c_{0}T ^{-1/2}\).
 This 2-form is also exact; but more to the point, \(u_{2}\) 
can be written as \(dz_{2}\) where \(z_{2}\) is a 1-form with \(| z_{2}| < c_{0}T ^{-1/2}\) 
on the \(\ir \geq\rho_{*}+ \frac{5}{8} \varepsilon  \) part of
\(\mathcal{N}_\varepsilon\). Hold on to \(z_{2}\) for the moment.  Let \(\sigma_{\perp}\) 
denote the function of  \(\r\)  on \(\mathcal{N}_\varepsilon\) given by
 \(\sigma_{\perp} =  \chi (8 \varepsilon ^{-1} (\ir - \rho_{*} ) -   5)\). This function is equal to 1 where
 \(\ir  <\rho_{*}  +\frac{5}{8} \varepsilon  \) and it is equal to 0 where  \(\ir  > \rho_{*} + \frac{3}{4}\varepsilon \).  Use
 \(w_{3}\) to denote the closed 2-form on \(Y_*\) that is given by
 \(w_{2}\) on \(\mathcal{Y}_{M}\), given by  \(w\)  on \(\mathcal{Y}_{0}\) 
and given by  \(w + d( \sigma_{\perp} z_{2})\) on \(\mathcal{N}_\varepsilon\) .
 The 2-form \(w_{3}\) has the same de Rham class as \(w_{2}\),
the same zero locus as it agrees with \(w_{2}\) 
where both are zero, and \(|w_{2}- w_{3} | \leq c_{0} T ^{-1/2}\).
 
Use \(\upsilon_{\diamond}\) to denote the \(\grg_{*}\)-Hodge dual on \(S^{1} \times S^{2}\)  of the
2-form  \(w  = \sin \theta  \, d \theta \wedge d \phi \).
 Write the \(\grg_{2}\)-Hodge star of \(w_{2}\) as \(\upsilon_{\diamond} +\q_{2}\) on \(\mathcal{Y}_{0}\) 
and on the \(\ir > \rho_{*} +\frac{1}{2} \varepsilon  \) part of
\(\mathcal{N}_\varepsilon\). As both the \(\grg_{2}\)-Hodge star of \(w_{2}\) and \(\upsilon_{\diamond}\) 
are exact on \(\mathcal{N}_\varepsilon\), it follows that \(\q_{2} =d
\o_{2}\) on \(\mathcal{N}_\varepsilon\). Moreover, such a function \(\o_{2}\) 
can be found with \( |\o_2 | \leq c_{0}T ^{-1/2}\) on the  \(\ir  >\rho_{*} + \frac{1}{2}\varepsilon  \) part of \(\mathcal{N}_\varepsilon\).
 This is so because \(|w -w_{2} | < c_{0} T ^{-1}\) and \(|\grg_{2} -\grg_{*} | < c_{0}T ^{-1}\) 
on this part of \(\mathcal{N}_\varepsilon\). Fix a version of \(\o_{2}\) that obeys this bound.  Let \(\upsilon_{3}\) denote the closed 1-form on \(Y_*\) given by \(\upsilon_{\diamond}\) on \(\mathcal{Y}_{0}\),
by the \(\grg_{2}\)-Hodge star of \(w_{2}\) on \(\mathcal{Y}_{M}\) and given on \(\mathcal{N}_\varepsilon\) 
by \(\upsilon_{\diamond} +d( \sigma_{\perp}\o_2 )\). This closed 1-form is such that
 \(w_{3} \wedge \upsilon_{3} \geq 0\) when \(T >c_{0}\)  with equality
 only at the zeros of \(w_{3}\).

 With \(T > c_{0}\) chosen, the
upcoming Lemma \ref{lem:D.2} uses what was just said about \(w_{3}\) and
 \(\upsilon_{3}\) as input to supply a metric on \(Y_{*}\)  with the
properties in the list that follows.  This new metric is denoted by \(\grg_{3T}\).
 The  \(\grg_{3T}\)-Hodge star sends \(w_{3}\) to \(\upsilon_{3}\); thus \(w_{3}\) 
is \(\grg_{3T}\)-harmonic. The metric \(\grg_{3T}\) on \(\mathcal{Y}_{0}\) 
and on the  \(\ir  > \rho_{*} + \frac{3}{4} \varepsilon   \) part of \(\mathcal{N}_\varepsilon\) is the metric \(\grg_{*}\).
 The metric \(\grg_{3T}\) on the \(\ir \in[\rho_{*} +
 \frac{1}{2}\varepsilon , \rho_{*} +\frac{3}{4}\varepsilon ]\) part of \(\mathcal{N}_\varepsilon\) can be written as
 \(\grg_{2} +\grh\) with \(\grh\) and its \(\grg_{2}\)-covariant
 derivaties to order 20 having \(\grg_{2}\)-norm less than \(c_{0}T ^{-1}\). Finally, the metrics \(\grg_{3T}\) and \(\grg_{2}\) are identical except on the rest of \(Y\).
 
 Any sufficiently large \(T\) version of the metric \(\grg_{3T}\) 
meets the requirements of Proposition \ref{prop:A.6}'s space \(\Met\). 
 Conversely, each metric in \(\Met\) is a sufficiently large \(T\) version of a
metric \(\grg_{3T}\) that is constructed as described above from
some metric in \(\op{Met}^{\mathcal{N}}\).  The lower bound on
\(T\) depends on various properties of the chosen
\(\op{Met}^{\mathcal{N}}\) metric, these being an upper bound
on the norm of the metric's Riemann curvature, the
metric volume of \(\mathcal{Y}_{M}\), and a lower bound on
the metric's injectivity radius.

\paragraph{Part 3:} The existence of the metric \(\grg_{3T}\) 
follows from the first lemma below.
 
\begin{lemma}\label{lem:D.2}   
Let \(Y_{Z}\)  denote an
oriented 3-manifold and let \(\grg\)  denote a given Riemannian metric
on \(Y_{Z}\).  Use \(*\)  in what follows to
denote the Hodge star defined by \(\grg\).  Suppose that \(U\)
 and \(V\)  are open sets in \(Y_{Z}\) 
with the closure of \(V\)  being a compact subset of \(U\).
 Let \(\omega\) and \(\upsilon\) denote respectively
a 2-form and a 1-form on \(Y_{Z}\)  such that
\(\omega  \wedge \upsilon  > 0\)  on \(U\) and such that  \(*\omega  =
\upsilon  \) on \(Y_{Z}-V\). 
\begin{itemize}
\item There are smooth metrics on \(Y_{Z}\)  which
equal \(\grg\)  on \(Y_{Z}-U\)  and have Hodge star sending \(\omega\) to \(\upsilon\).
  Moreover, there exists metric of this sort whose volume
3-form is the same as the \(\grg \)-volume 3-form.
\item  Fix \(k \in \{0, 1, \ldots\}\)  and \(\textsc{d}> 1\).   There
  exists \(\kappa> 1\)  with the following significance:  Suppose that the
  \(C^{k}\)-norms on \(U\)  of \(\omega\), \(\upsilon\) and the Riemann curvature
tensor of \(\grg\) are less than \(\textsc{d}\).   Then
\(Y_{Z}\)  has a metric that obeys the conclusions of
the first bullet and differs from \(\grg\)  by a tensor whose
\(\grg\)-norm and those of its \(\grg\)-covariant derivatives
to order \(k\)  are bounded by \(\kappa\) times the \(C^{k}\)-norm of  \(*\omega  -\upsilon \).
\end{itemize}
\end{lemma}

Lemma \ref{lem:D.2} has a generalization that holds for 1-parameter families of
data sets.  This parametrized version is given below but used in the
next subsection.

\begin{lemma}\label{lem:D.3}  
Let \(\{(\grg_{ \tau }, \omega _{ \tau },\upsilon_{ \tau })\}_{\tau\in
  [0,1]}\) denote a smoothly parametrized family of metrics, 2-forms and
1-forms on \(Y_{Z}\)  with \(\omega _{ \tau } \wedge\upsilon_{ \tau } > 0\) on
\(U\)  and such that the \(\grg_{\tau }\)-Hodge dual of \(\omega_{\tau
}\) is \(\upsilon_{\tau}\) on \(Y_Z-V\) .  There is
a corresponding smooth, 1-parameter family of metrics such that each
 \(\tau\in [0, 1]\) member obeys the conclusion of first
bullet of Lemma \ref{lem:D.2}.  Moreover, this new  family of metrics
can be chosen to obey the properties listed below. 
\begin{itemize}
\item  Let \(I \subset [0, 1]\) denote an open neighborhood of
one or both of the end points.  Suppose that the conclusions of the
first bullet of Lemma \ref{lem:D.2} hold for \((\grg_{\tau }, \omega _{ \tau }, \upsilon_{ \tau })\)  when  \(\tau  \in
I\).  There is a neighborhood \(I' \subset I\)  of the
endpoints such each  \(\tau  \in I'\)  member of the new
family is the corresponding \(\grg_{\tau }\).
\item Given a non-negative integer \(k\)  and \(\textsc{d} \geq 1\), there exists \(\kappa  > 1\) 
with the following significance:  Suppose that the conditions of the
second bullet of Lemma \ref{lem:D.2} are satisfied for each  \(\tau 
\in[0,1]\)  and that the \(C^{k}\)-norms
of the \(\tau\)-derivatives to order \(k\)  of \(\{(g_{ \tau },  \omega _{ \tau },\upsilon_{ \tau })\}_{ \tau\in [0,1]}\)
are also bounded by \(\textsc{d}\).  There is a
1-parameter family of metrics that obeys the preceding bullet and the
first and second bullets of Lemma \ref{lem:D.2}.  In addition, each  \(\tau 
\in [0, 1]\)  member of the family differs from the
corresponding metric \(\grg_{\tau}\) by a tensor whose \(\tau\)-derivatives to order \(k\)  have
\(C^{k}\)-norm bounded by \(\kappa\) times the \(C^{k}\)-norm of the sum
of the \(\tau\)-derivatives to order \(k\)  of the difference between
\(\upsilon_{\tau  }\) and the \(\grg_{\tau}\)-Hodge star of \(\omega_{\tau }\).
\end{itemize}
\end{lemma}

\textit{Proof of Lemmas \ref{lem:D.2} and \ref{lem:D.3}.}  Let \(\Omega\) denote \(\grg\)'s volume
3-form.  Write \(\omega  \wedge \upsilon \) as \(\grq \Omega \) with
\(\grq\) being a non-negative function on \(U\).  Let \(v\) denote the vector field
on \(U\) that is annihilated by \(\omega\) and has pairing \(\grq\) with
\(\upsilon\).  Let \(\ker (\upsilon) \subset TU\) denote the 2-plane
bundle that is annihilated by \(\upsilon\).  The 2-form  \(\omega\)  is
symplectic on \(\ker (\upsilon)\) and so orients \(\ker (\upsilon)\).
 Choose an \(\omega\)-compatible almost complex structure on
\(\ker (\upsilon)\), denoted by \(J\) below.  Note in this regard that there
are no obstructions to finding such an almost complex structure.  This
is so because the space of almost complex structures that are
compatible with a constant symplectic form on
\(\bbR^{2}\) is contractible.  The construction
just given yields a new metric with volume 3-form \(\Omega\).

With \(J\) chosen, a metric on \(U\) is defined as follows:  The vector field \(v\)
has norm \(\grq^{1/2}\) and is orthogonal to \(\ker (\upsilon)\).
The inner product between vectors \(v\) and \(v'\) in a given fiber of
\(\ker (\upsilon)\) is \(\grq^{-1/2}  \omega (v, J v')\).  A metric of this sort has  \(*\omega = \upsilon\)
and is such that both \(\omega\) and \(\upsilon\) have norm
\(\grq^{1/2}\).  Moreover, any metric with these two
properties is of the form just described.  In particular, any two
differ only with respect to the choice of the almost complex structure
on the \(\ker (\upsilon)\).  
 
Let \(J_{1}\) denote a chosen,  \(\omega\)-compatible almost
complex structure on \(\ker (\upsilon)|_{U}\) and
let \(\grg_{1}\) denote the corresponding metric. The metric \(\grg\) on
\(U-V\) is by necessity of the sort just described, thus it differs from
\(\grg_{1}\) only on \(\ker (\upsilon)\).  In particular, the
metric \(\grg\) on \(\ker (\upsilon\)) is given by \(\grq^{-1/2}
 \omega (v, J_{\grg}v')\) with \(J_{\grg}\)  being an
 \(\omega\)-compatible almost complex structure on \(\ker (\upsilon)|_{U- V}\) .
 As noted above, if point \(p \in U\), then the space of \(\omega|_{p}\)-compatible
almost complex structures on \(\ker (\upsilon)|_{p}\) is contractible.  This
understood, there are no obstructions to choosing an
\(\omega\)-compatible almost complex structure on \(\ker (\upsilon)|_{U}\) that agrees with
\(J_{\grg}\) near \(Y_{Z}-U\) and agrees with \(J_{1}\) on \(V\).  Let \(J_{2}\) denote an almost
complex structure of this sort.  The metric defined as instructed
above by \(J_{2}\) has the properties that are asserted by the
first bullet of Lemma \ref{lem:D.2}.
 
The assertions of the second bullet of Lemma \ref{lem:D.2} and those
of Lemma \ref{lem:D.3} are proved by taking care with the choice of \(J_{2}\) and its
 \(\tau  \in [0, 1]\) counterparts.  As the details are
straightforward and rather tedious, they are omitted. \epf
 
\subsection{\(\mathrm{Met}_T\) metrics on cobordisms}\label{sec:Dc)}
 
Lemma \ref{Lemma D.1} has an analog given below that concerns
self-dual forms on cobordisms.  The cobordism manifold is denoted
below by \(X\) and it is assumed to be of the sort that is described in
Section \ref{sec:Ac)} with its constant  \(s\) slices where  \(s  < -1\) and
 \(s  > 1\) given as follows:  Either one is \(Y\) and the other is \(Y_{\G}\);
or one is some \(k\in \{1, \ldots, \G\}\) version of 
\(Y_k\) and the other is \(Y_{k-1}\sqcup (S^1\times S^2)\), or one is
\(Y_0\) and the other is \(M\sqcup S^{1} \times S^{2}\).  The case
when both are \(Y\) or both some \(k\in \{1, \ldots, \G\}\) version of 
\(Y_k\) is also
allowed, but only the case where both are \(Y_{\G}\) are needed in
what is to come.
The topology of
\(X\) is further constrained by the requirement that  \(s\) have
1 critical point when it is not diffeomorphic to a product with
\(\bbR\). If one of these slices is \(Y\) and the other \(Y_{\G}\), or
if both are \(Y\) or both \(Y_k\) for \(k\in \{1, \ldots, \G\}\), then
\(s\) has no critical points and the cobordism manifold \(X\) is
\(\bbR\times Y\) or \(\bbR\times Y_k\) as the case may be, with the projection to \(\bbR\) given by the function \(s\).

One more constraint on \(X\) is needed.
By way of background, what is said in Part 4 of Section \ref{sec:Da)}
identifies \(\mathcal{Y}_{0} \cup  \mathcal{N}_\varepsilon \) as a
subset of \(Y\) and \(Y_{k}\), and also \(S^{1} \times S^{2}\). This extra constraints uses \(\mathcal{Y}_{0\varepsilon}\) to denote the union of \(\mathcal{Y}_{0}\) 
and the \(\ir  > \rho_{*} +\frac{1}{128}\varepsilon  \) part of
\(\mathcal{N}_\varepsilon\). Here is the extra constraint:
\begin{equation}\label{eq:(D.10)}\begin{split}
&\text{There is a distinguished embedding of \(\bbR \times  \mathcal{Y} _{0\varepsilon }\) into
\(X\) with the following property:}\\  
& \text{The respective \(s <0\) and \(s >0\) slices of the image of
  this embedding,}\\
&\text{when
written using the diffeomorphisms from the second and third bullets of (\ref{(A.9a,11)}),}\\
&\text{appear as the incarnation of \(\mathcal{Y}_{0\varepsilon}\) in
  either \(Y\), \(Y_k\), or \(S^{1}\times S^{2}\) as the case may be.}
\end{split}
\end{equation}

The metric for \(X\) is assumed to obey a constraint that requires membership in an analog
for \(X\) of the various \(T >1\) versions of the space \(\op{Met}_{T}\). The definition of this \(X\) version of
\(\op{Met}_{T}\) requires the a priori selection of metrics \(\grg_{-}\) and \(\grg_{+}\) from the respective \(Y_{-}\)  and
\(Y_{+}\)  versions of \(\op{Met}_{T}\) with it
understood that \(\op{Met}_{T}\) in the case of \(M   \sqcup (S^{1} \times S^{2})\) is the
space consisting of the metric \(\grg_{*}\)
on \(S^{1} \times S^{2}\) and a metric on \(M\) of the following sort: If
\(c_1(\det\,( \bbS))\) is torsion on \(M\), then any metric on \(M\) is allowed. If this class is not torsion,
then the metric's associated harmonic 2-form with de Rham coholomogy
class \(c_1(\det\,( \bbS|_M))\) has non-degenerate zeros. Meanwhile,
\(\op{Met}_T\) for any given \(k\in \{1, \ldots, \G\}\) version of
\(Y_k\sqcup (S^1\times S^2)\) consists of a \(\op{Met}_T\) metric for
\(Y_{k-1}\) and any metric for \(S^1\times S^2\).
Reintroduce from Part 5 of Section \ref{sec:Da)} the
metric \(\grg_{T}\) on \(\mathcal{Y}_{0}\cup  \mathcal{N}_\varepsilon
\). Of immediate interest in what follows is \(\grg_{T}\)'s restriction to \(\mathcal{Y} _{0\varepsilon }\).
 By way of a reminder, \(\grg_{T}\) on \(\mathcal{Y}_{0\varepsilon}\) is
 the metric \(\grg_{*}\) on \(\mathcal{Y}_{0}\) 
and it is the metric in (\ref{eq:(D.7)}) on the
\(\ir > \rho_{*} + \frac{1}{128}\varepsilon  \) part of \(\mathcal{N}_\varepsilon\).
 
The analog of \(\op{Met}_{T}\) for \(X\)
consists of the space of metrics with the following three properties:  
\BTitem\label{eq:(D.11)}
\item The metric obeys the \(L =100\) version of (\ref{eq:(A.12,15a)}).
\item The metric pulls back via the embedding in (\ref{eq:(D.10)}) as the metric
\(d s^{2}+  \grg_{T}\).
\item The metric pulls back from the \(s  \leq -104\) part of \(X\) via the embedding in
the second bullet of (\ref{eq:(A.8)}) as \(d s^{2} + \grg_{-}\), and it pulls back from the \(s 
> 104\) part of \(X\) via the embedding from the third bullet of (\ref{(A.9a,11)}) as
\(d s^{2} + \grg_{+}\).
\ETitem
 
This analog for \(X\) of \(\op{Met}_{T}\)  is denoted in what follows by
\(\op{Met}_{T}\) also, its dependence on \(\grg_{-}\) and \(\grg_{+}\) being implicit.  

 Lemma \ref{lem:D.4} given momentarily supplies the promised
analog to Lemma \ref{Lemma D.1}.  To set the notation, suppose that a
metric on \(X\) has been specified and that \(\p_{X}\) 
is a differential form on \(X\).  The lemma uses \(\langle  \p_{X}  \rangle _{2}\) to denote the
\(L^{2}\)-norm of \(\p_{X}\) over the \(|s |  < 104\) part of \(X\).  Lemma \ref{lem:D.4}
uses \(w_{-}\) and \(w_{+}\) to denote the respective \(\grg_{-}\) and \(\grg_{+}\) 
harmonic 2-forms with de-Rham cohomology class that of
\(c_{1} (\det\, (\bbS ))\); and it uses the embeddings from the second
and third bullets of (\ref{(A.9a,11)}) to view \(w_{-}\) and \(w_{+}\) 
as 2-forms on the  \(s  \leq -1\) and  \(s  > 1\) parts of \(X\). 
 
 \begin{lemma}\label{lem:D.4} 
Let \(X\) denote a cobordism manifold of the sort described above.
Given metrics \(\grg_{-}\) and \(\grg_{+}\) in the respective \(Y_{-}\) and
\(Y_{+}\) versions of \(\op{Met}_{T}\), there exists \(\kappa > 1\) with the following
significance:  Fix \(T >1\), and fix a Riemannian metric on
\(X\) from the corresponding set \(\op{Met}_{T}\).
 There is a self-dual, harmonic 2-form on \(X\) whose  pull-back to the
constant \(s\)-slices of \(X\) converges as \(s \to -\infty \) to
\(w_{-}\) and as \(s \to \infty \) to \(w_{+}\). Let \(\p_{X}\) denote a closed 2-form on
\(X\) that equals \(w_{-}\) where \(s  <-102\), that equals \(w_{+}\)
where \(s>102\), and with de Rham cohomology class that of \(c_{1}(\det\, (\bbS ))\).
\begin{itemize}
\item The \(L^{2}\)-norm of this harmonic self-dual 2-form on the
\(s\)-inverse image of any length 1 interval in \(\bbR\) is bounded by \(\kappa
\langle  \p_{X} \rangle _{2}\).
\item The pull-back of this harmonic self-dual
2-form to the constant \(s > 1\) and \(s <-1\) slices differs in the \(C^{1}\)-topology 
from \(w_{-}\) and \(w_{+}\) by at most \(\kappa \langle
\p_{X}\rangle _{2}\, e^{-|s|/z}\) with \(z \geq 1\) depending on the corresponding limit metric.
\item  The pull-back of this harmonic self-dual 2-form to \(\bbR\times \mathcal{Y}_{0 \varepsilon}\) 
via the embedding from (\ref{eq:(D.10)}) differs from
\(d s \wedge\upsilon_{\diamond} + w\) by a 2-form whose \(C^{1}\)-norm
on \(\bbR \times  \mathcal{Y}_{0}\) and on the \(\ir  >\rho_{*}+\frac{1}{2}\varepsilon  \) part of \(\bbR
\times\mathcal{N}_\varepsilon \) is less than
\(\kappa\langle\p_X\rangle_2
T ^{-1/2}\).
\end{itemize}
\end{lemma}

\pf The existence of a closed, self-dual harmonic 2-form with the desired
 \(s \to  -\infty \) and \(s  \to  \infty \) limits follows from the
 index theorem in \cite{APS}.  This 2-form is denoted in what
follows by   \(\omega\).  Given
the first bullet, then the assertion in the second
 bullet follows from the eigenfunction expansion that
is depicted below in (\ref{eq:(D.13)}).  As explained next, the third bullet also
follows from the second bullet. 
 
To prove the third bullet, fix \(s_{0} \subset  \bbR \) and introduce \(\sigma_{0}\) 
to denote the function on \(\bbR\)  given by the rule \(s  \mapsto \chi (|s-  s_{0}|-   1)\).  This
function equals 1 where \(|s -  s_{0}|\) is less than 1 and it equals zero where
\(|s - s_{0}|\) is greater than 2.  Let \(\sigma\) denote the function
from Step 2 of the proof of Lemma \ref{Lemma D.1} and let \(\Phi_{T}\) 
denote the embedding from Part 5 of Section \ref{sec:Da)}.  View the
 \(\Phi_{T}^{-1}\) pull-back of \(\sigma_{0} \sigma  \omega   \) as a 
 2-form on \(\bbR \times (S^{1}\times S^{2})\) with support where \(|s - s_{0}| < 2\).  The assumed
\(L^{2}\)-bound for \(\omega\) with a Green's function argument much
like that used in Step 4 of the proof of Lemma \ref{Lemma D.1} can be used to derive the pointwise bound
that is asserted by Lemma \ref{lem:D.4}.  The derivation differs little from
that in Step 4 of the proof of Lemma \ref{Lemma D.1} save for the fact that the
Green's function in question is that for the elliptic operator 
\begin{equation}\label{eq:(D.12)}
 \mathcal{D}\co C ^{\infty} (\bbR \times(S^{1} \times S^{2}); \Lambda ^{+} \oplus  \bbR )\to 
C ^{\infty} (\bbR \times(S^{1}\times S^{2}); T^*( \bbR  \times(S^{1} \times S^{2} ))
\end{equation}
given by the formula \(\mathcal{D}  =*_{X} d_{X} + d_{X}\)  where \(d_{X}\) denotes the
4-dimensional exterior derivative \(d s   \wedge \ps (\cdot) + d\) and where \(*_{X}\) 
denotes the Hodge star for the metric \(d s ^{2} + \grg_{*}\).
  
The lemma's first bullet is proved in the four steps
that follow.
 
\paragraph{\it Step 1:}  Let \(\omega\) denote the relevant
closed, self-dual harmonic form.  Fix an integer \(\n \in  \{106, 107, \ldots\}\) and
introduce by way of notation \(I_{\n} \subset  \bbR  \) to denote a
closed interval of length \(2 \n\)  whose
endpoints have distance 106 or more from the origin.  Let
 \(\textup{C}\)  denote the space of closed two forms on the domain
 \(s ^{-1} (I_{\n})\) that agree with \(\omega\) on
some neighborhood of the \(s\)-inverse images of the boundary points of
\(I_{\n}\).  The 2-form \(\omega\) is the minimizer in \(\textup{C}\)  of the functional that is
defined by the rule \(\grw  \mapsto\int_{s^{-1}(I_{\n})}|\grw^+|^2\).
 
\paragraph{\it Step 2:}  Use the embedding from the
second bullet of (\ref{(A.9a,11)}) to write the \(s  \geq 100\) part of \(X\) as
\([100, \infty)\times Y_{+}\) and likewise write the \(s < -100\) part
of \(X\) as \((-\infty, -100] \times Y_{-}\).  Let \(Y_{*}\) for the moment denote either
\(Y_{+}\) or \(Y_{-}\). Let \(*\)  denote either the \(\grg_{-}\)- or
\(\grg_{+}\)-version of the Hodge star on \(Y_{*}\).  The corresponding operator \(d*\) 
defines an unbounded, self-adjoint operator on the space of closed
2-forms on \(Y_{*}\). Let \(\Xi^{-}\) denote an \(L^{2}\)-orthonormal
basis of eigenvectors of \(d*\) on the space of closed 2-forms with
negative eigenvalue and let \(\Xi^{+}\) denote an \(L^{2}\)-orthonormal
basis of eigenvectors of \(d*\) with positive eigenvalue.  The eigenvalue of
\(d*\) on a given eigenvector, \(\op{a}\), is denoted by \(\lambda_{\op{a}}\).
 
 The 2-form \(\omega\) on \((-\infty, -1]\times Y_{-}\) and on \([1, \infty)  \times Y_{+}\)  can be written as 
\begin{equation}\label{eq:(D.13)}
\begin{cases}
 \omega   =d s   \wedge* w_{-} + w_{-}
 +\sum_{\op{a}\in\Xi^+}\textsc{z}_{\op{a}} \, 
  e^{\lambda_{\op{a}}(s+1)} (d s   \wedge* \op{a} + \op{a})
  &\text{where \(s  \leq -104\).} \\
 \omega   =d s   \wedge* w_{+} +w_{+} +
 \sum_{\op{a}\in\Xi^-}\textsc{z}_{\op{a}} \, 
  e^{\lambda_{\op{a}}(s-1)} (d s \wedge * \op{a} + \op{a}) &
\text{where \(s  \geq 104\).}
\end{cases}
\end{equation}
 
What is denoted by \(\textsc{z}_{(\cdot)}\)
in (\ref{eq:(D.13)}) is a real number.  Keep in mind for what follows that any
given version of \(e^{\lambda_{\op{a}}s}(d s  \wedge * \op{a} + \op{a})\) is the
exterior derivative on its domain of definition of the 1-form
\(\q_{\op{a}} = \lambda_{\op{a}}^{-1}e^{\lambda_{\op{a}}s}* \op{a}\).  
 
 \paragraph{\it Step 3:}  Fix \(m  > 1\).  Let \(\op{a}\) denote an eigenvector in the \(Y_{-}\) version of
  \(\Xi^{+}\). Introduce \(\sigma_{\op{a}}\) to denote the function on
 \(\bbR\) given by the rule \(s  \mapsto \sigma_{\op{a}} ( s ) = 1 - \chi (- m ^{-1} \lambda_{\op{a}} (s + 102)-1)\).  This
function equals 0 where \(s > -102 - m \lambda_{\op{a}}^{-1}\) 
and it equals 1 where \(s < -102 -2 m \lambda_{\op{a}}^{-1}\).
 If \(\op{a}\) is in the \(Y_{+}\)  version of \(\Xi^{-}\), then \(\sigma_{\op{a}}\) 
is given by the rule \(s  \mapsto\sigma_{\op{a}} (s )= 1 -\chi (- m ^{-1} \lambda_{\op{a}} (s - 102) - 1)\).  This
version of \(\sigma_{\op{a}}\) is 0 where \(s <102 + m |\lambda_{\op{a}} |^{-1}\) and it is 1 where
\(|s |> 102 +2 m |\lambda_{\op{a}} |^{-1}\). Meanwhile, use \(\chi_{*}\) to denote the function
\(\chi (102-|s |)\).  This function is 1 where \(|s |  > 102\) and 0 where \(|s |  < 101\).
 
Use \(\p_{X}\) and these functions to define the 2-form \(\grw\) on \(X\) by the rule
\begin{equation}\label{eq:(D.14)}
 \begin{split}
\grw  & = \chi_{*}  d s  \wedge* \p_{X} + \p_{X} + \sum_{\op{a}\in\Xi^+}\textsc{z}_{\op{a}} \,
 d (\sigma_{\op{a}} \lambda_{\op{a}}^{-1}e^{\lambda_{\op{a}}(s+102)}* \op{a})\\
&  \qquad + \sum_{\op{a}\in\Xi^-}\textsc{z}_{\op{a}}\,
 d( \sigma_{\op{a}} \lambda_{\op{a}}^{-1}e^{\lambda_{\op{a}}(s-102)}* \op{a}).
\end{split}
\end{equation}
  
This is a closed 2-form whose de Rham cohomology class
is the same as \(\omega\).  Let \(\textsc{e}\) denote the smallest of the numbers from the set
\(\{\lambda_{\op{a}} \, |\, \op{a}  \in \Xi  ^{+} \} \cup \{|\lambda_{\op{a}}| \, |\, \op{a}  \in   \Xi  ^{-} \}\)
with it understood that \(\Xi^{+}\) refers to the \(Y_{-}\)  version and \(\Xi  ^{-}\) 
refers to the \(Y_{+}\) version.  The 2-form \(\grw\)  equals \(\omega\) where
\(|s | \geq 1 + 2 m \, \textsc{e}^{-1}\).

\paragraph{\it Step 4:}  The square of the \(L^{2}\)-norm of \(\grw^{+}\) over the
\(|s |  \leq 102 +  2 m \textsc{e}^{-1}\)  part of \(X\) is no greater than 
\begin{equation}\label{eq:(D.15)}\begin{split}
 \int_{s^{-1}([-102, 102])}|\p_X|^2 & +c_{0} m ^{-2}
 e^{-2m}\sum_{\op{a}\in \Xi^+\cup\Xi^-} |\lambda_{\op{a}} |^{-1} 
|\textsc{z}_{\op{a}}|^{2}+ 4 m \textsc{e}^{-1} ( \|w_{-} \|_{2} ^{2} +  \| w_{+} \|_{2} ^{2} )\\
 & 
+ \sum_{\op{a}\in \Xi^+\cup\Xi^-} |\lambda_{\op{a}} |^{-1} 
|\textsc{z}_{\op{a}}|^{2}(e^{-2m}-e^{-4  |\lambda_{\op{a}} |m/\textsc{e}}).
\end{split}
\end{equation}
Meanwhile, the integral of \(\omega\) over this same part of \(X\) is equal to 
\begin{equation}\label{eq:(D.16)}\begin{split}
\int_{s^{-1}([-102, 102])}|\omega|^2 & +\sum_{\op{a}\in \Xi^+\cup\Xi^-}  |\lambda_{\op{a}} |^{-1} 
|\textsc{z}_{\op{a}}|^{2}(1-e^{-2m})+ 4 m \textsc{e}^{-1} (\|w_{-} \|_{2}^{2} +  \|w_{+} \|_{2}^{2} ) \\
 & + \sum_{\op{a}\in \Xi^+\cup\Xi^-}  |\lambda_{\op{a}} |^{-1} 
|\textsc{z}_{\op{a}}|^{2}(e^{-2m}-e^{-4 |\lambda_{\op{a}} | m/\textsc{e}}). 
\end{split}
\end{equation}

 As noted in Step 1, the expression in (\ref{eq:(D.16)}) can
not be greater than what is written in (\ref{eq:(D.15)}).  This being the case,
the \(m>c_{0}\)  versions of (\ref{eq:(D.15)}) and (\ref{eq:(D.16)}) imply  the bound
\begin{equation}\label{eq:(D.17)}
\int_{s^{-1}([-102, 102])}|\omega|^2 +\sum_{\op{a}\in \Xi^+\cup\Xi^-}  |\lambda_{\op{a}} |^{-1} 
|\textsc{z}_{\op{a}}|^{2}  \leq c_{0} (1 + \langle  \p_{X}  \rangle _{2}^{2} ).
\end{equation}

This last bound has the following corollary:  Let \(I\subset \bbR  \) denote any interval
of length 1.  Then \(\int_{s^{-1}(I)}|\omega|^2 \leq c_0\,( 1 +  \langle
\p_{X}  \rangle _{2}^{2} )\). \epf

 \subsection{Proof of Proposition \ref{prop:A.8}}\label{sec:Dd)}
 To explain the first bullet, identify a neighborhood of the critical point of the function \(s\)  with a ball about the origin in \(\bbR^{4}\)  
using coordinates \((y_{1}, y_{2}, y_{3}, y_{4})\) and write
 \(s\)  in terms of these coordinates as  \(s  =y_{4}^{2} -y_{1}^{2} -y_{2}^{2} -y_{3}^{2}\) 
when the constant,  \(s  <-1\) slices of \(X\) are \(Y_{0}\)  and the constant,  \(s  > 1\)
slices are \(M \sqcup  (S^{1}  \times S^{2})\).  With the ends reversed, the function \(s\) 
appears as \(s  =-  y_{4} ^{2} +y_{1} ^{2} +y_{2} ^{2} +y_{3} ^{2}\).
 The embeddings given in the second and third bullets of (\ref{(A.9a,11)}) are
defined using a pseudogradient vector field for \(s\).  This pseudogradient vector
field in the \(Y_{-} = Y_{0}\)  and \(Y_{+}  =  M\sqcup  (S^{1} \times S^{2})\) case can be
chosen so as to have the following properties:  The inverse image of
the descending 3-ball from the critical point via the embedding given
by the second bullet of (\ref{(A.9a,11)}) appears as the locus \((-\infty, 0)\times  S\) with \(S\) being the
2-sphere that is described in Part 4 of Section \ref{sec:Da)}.  Meanwhile, the
inverse image via the embedding given by the third bullet of (\ref{(A.9a,11)}) of
one of the ascending arc from this critical point intersects the
\((0,\infty)\times (S^{1}\times S^{2})\) component of \((0, \infty)\times  (M\sqcup (S^{1} \times S^{2}))\) as the
locus \((0, \infty) \times p_{*}\)  with \(p_{*}  \in S^{1}\times S^{2}\)  being the \(\ir  = 0\) point in the ball that is
described in the third paragraph of Part 4 in Section \ref{sec:Da)}.  The other
ascending arc intersects the \((0, \infty) \times  M\) component as the
 \(\ir ' = 0\) point in the ball that is described in the fourth
 paragraph of Section \ref{sec:Da)}.   There is a completely analogous
 picture of \(X\) when \(Y_{0}\)  is the constant  \(s> 0\) slice of
 \(X\) and \(S^{1} \times S^{2}\) is the constant \(s < 0\) slice.
 
 What is said above about the descending and
ascending submanifolds from the critical point has the following
consequence:  The pseudogradient vector field that defines the
embeddings from the second and third bullets of (\ref{(A.9a,11)}) can be chosen so
that (\ref{eq:(A.9b)}) are obeyed and likewise the
condition in (\ref{eq:(D.10)}).   These properties are assumed in what follows.
 The fact that \(S\) carries no homology implies that the fourth
 bullet of (\ref{(A.9a,11)}) holds for \(X\).
 
Parts 1-10 of this subsection construct large \(L\) versions of the form \(w_{X}\) 
and the metric that are used in Part 11 to satisfy the requirements of the
second bullet of Proposition \ref{prop:A.8}.  These constructions require the
choice of parameters \(T \gg 1\), \(L_{0} \gg 1\) and \(L_{1} >L_{0} + 1\).  Granted
large choices, Parts 1-10 construct a closed 2-form denoted by \(\omega_{T *}\)  
and a metric denoted by \(\grm_{T *}\) that makes \(\omega _{T *}\) self-dual.  Any \(L >L_{1}+ 20\) version of \(\omega _{T *}\) can serve for Proposition \ref{prop:A.8}'s desired 2-form \(w_{X}\) 
and the corresponding version of \(\grm_{T *}\) can serve for the desired metric. 

Proposition \ref{prop:A.8} requires as input a metric on \(M \sqcup
(S^1\times S^2)\) and asserts that such a metric determines a certain subset of the set \(\Met\) on \(Y_0\). To say more about this
subset, recall from Part 2 of Section \ref{sec:Db)} that each metric in
 \(\Met\)  is determined in part by a metric from Section \ref{sec:Da)}'s set
\(\op{Met}^{\mathcal{N}}\) 
and a large choice for a number denoted by \(T\).  A metric of this sort
was denoted by \(\grg_{3T}\) in Section \ref{sec:Db)}.  As noted at the
end of Part 2 of Section \ref{sec:Db)}, a lower bound on \(T\) is
determined by certain properties of the metric from \(\op{Met}^{\mathcal{N}}\).
 A metric of this sort is in Proposition \ref{prop:A.8}'s subset
if and only if \(T\) is greater than a new lower bound that is determined
by the afore-mentioned properties of the \(\op{Met}^{\mathcal{N}}\) 
metric.   Suffice it to say for the purposes of the proof that this
new lower bound is defined implicitly by the constructions in the
subsequent eleven parts of this subsection.    

The upcoming Parts 1-10 are written so as to
simultaneously supply a metric and a closed, self-dual 2-form for Section
\ref{sec:De)}'s proof of Proposition \ref{prop:A.10a} and Section
\ref{sec:Df)}'s proof of Proposition \ref{prop:A.10b}.  This is done by considering a cobordism space \(X\) as described in the
previous section whose limit manifolds \(Y_-\) and \(Y_+\) are as follows: Either one is \(Y_0\) and the
other is \(M  \sqcup( S^{1}  \times S^{2})\); or one is some
\(k\in\{1, \ldots, \G\}\) version of \(Y_k\) and the other is \(Y_{k-1}  \sqcup( S^{1}  \times S^{2})\), or both are \(Y_{\G}\). Although not needed for what follows, the constructions in 
Parts 1-10 can be done when both limit manifolds are \(Y\) or both are some
\(k\in\{1, \ldots, \G\}\) version of \(Y_k\). 

\paragraph{Part 1:}  When \(Y_-\) or \(Y_+\) is not \(M  \sqcup  (
S^1\times S^2 )\), choose metrics, \(\grg_{1-}\) and \(\grg_{1+}\) 
in the respective \(Y_{-}\)  and \(Y_{+}\)-versions of
\(\op{Met}^{\mathcal{N}}\) as the case may be. In the case when one of
\(Y_-\) or \(Y_+\) is some \(k\in\{1, \ldots, \G\}\) version of \(Y_k\)
and the other is \(Y_{k-1}\sqcup (  S^1\times S^2)\), what is denoted \(\op{Met}^{\mathcal{N}}\)
allows any metric for the \(S^1\times S^2\) component. Fix a \(T\gg
1\); in particular so that Lemma \ref{Lemma D.1} can be invoked for the metric in \(\op{Met}_T\) defined using \(\grg_{1-}\) in the case of \(Y_-\) and \(\grg_{1+}\) in the case of \(Y_+\). Use \(\grg_{1-}\) to choose a metric \(\grg_2\) as directed in Part 2 of Section \ref{sec:Db)} on \(Y_-\). Then set \(\grg_-=\grg_2\). Meanwhile, use \(\grg_2\) to construct a version of the metric \(\grg_{3T}\) and denote it by \(\grg_{-T}\). Do the same using \(\grg_{1+}\); denote the chosen \(\grg_2\) metric on \(Y_+\) by \(\grg_+\) and use \(\grg_{+T}\) to denote the resulting \(\grg_{3T}\) metric.
If either of \(Y_{-}\)  or \(Y_{+}\)  is \(M   \sqcup(S^{1} \times
 S^{2} )\), take the metric of the sort described in Part 1 
of Section \ref{sec:Ag)} for \(M\) and the metric \(\grg_*\) on
\(S^1\times S^2\). Denote the resulting
metric on \(M   \sqcup(S^{1} \times S^{2} )\) as \(\grg_-\) in the \(Y_-\)
case and \(\grg_+\)  in the \(Y_+\) case. With \(T\geq 1\)
chosen, this same metric is also denoted at times by \(\grg_{-T}\) and \(\grg_{+T}\) as the case may be.

By way of notation, the constant \(c_0\) in what follows depends
implicitly on the various properties of the metrics \(\grg_{1-}\) and
\(\grg_{1+}\). In particular, \(c_0\) depends on an upper bound for
the norm of the metric’s curvature, upper and lower bounds on the
metric’s volume and a lower bound on the injectivity radius.

Let \(\grm\) denote a chosen metric in the \(\grg_-\) and \(\grg_+\) version of \(\op{Met}_T\) on \(X\). Certain constraints on \(\grm\) are imposed later in this subsection. Note that some of the latter impose constraints on \(\grg_{1-}\) and \(\grg_{1+}\).

\paragraph{Part 2:}   Use \(w_{-}\) and \(w_{+}\) to denote the
respective \(\grg_{-}\) and \(\grg_{+}\) harmonic 2-forms on \(Y_{-}\)
and \(Y_{+}\) with de Rham cohomology class that of
\(c_{1} (\det( \bbS ))\). Fix for the moment a closed 2-form \(\p_{X}\) 
on \(X\) as described in Lemma \ref{lem:D.4}.  Use \(\omega\) to denote the
self-dual 2-form on \(X\) given by Lemma \ref{lem:D.4} for the case when the metric
on \(X\) is  \(\grm\).  The distinguished
embedding from (\ref{eq:(D.11)}) pulls \(\omega\) back to \(\bbR \times \mathcal{Y}_{0\varepsilon }\)
as a 2-form that can be written as 
\begin{equation}\label{eq:(D.18)}
 \omega   =d s   \wedge \upsilon_{\diamond} +  w  +d s   \wedge\frac{\partial}{\partial s}\q   +   d \q ,
\end{equation}
 with  \(\q\)  being an \(s\)-dependent 1-form on \(\mathcal{Y}_{0\varepsilon}\). 
 Lemma \ref{lem:D.4} says that the \(C^{1}\)-norms of 
  \(\frac{\partial}{\partial s}\q\)  and \(d \q\) on \(\bbR \times  \mathcal{Y} _{0}\) and on the  \(\ir  >
 \rho_{*} + \frac{1}{2} \varepsilon  \) part of \(\bbR \times
 \mathcal{N}_\varepsilon \) are less than \(c_{0} \langle  \p_{X}  \rangle _2 T ^{-1/2}\).

 An  \(s\)-  and \(T\)-independent open cover of \(\mathcal{Y}_{0\varepsilon }\) 
by balls of radius \(c_{0}^{-1} \varepsilon \) can be used to write \(\q\)  on \(\mathcal{Y}_{0}\) 
and on the \(\ir  >\rho_{*}+ \frac{17}{32}\varepsilon  \) part of \(\mathcal{N}_\varepsilon\) 
as \(\q_{0}  + d \k\)  with \(\q_{0}\) obeying
\( | \q_{0} | \leq c_{0}  \langle  \p_{X}  \rangle _{2} T ^{-1/2} \)
and \( | \frac{\partial}{\partial s}\q_{0}| \leq c_{0}\langle  \p_{X}
\rangle _{2 }T^{-1/2}\).  Meanwhile, \(\k\) is a smooth function with
\(|d( \k)| \leq c_{0} \langle  \p_{X}  \rangle _{2}T^{-1/2}\).  Both
\(\q_{0}\) and \(\k\) can be constructed so as to depend smoothly on \(s\).  It follows as a
consequence of the bound \(|d( \frac{\partial}{\partial s}\k)| \leq
c_{0} \langle  \p_{X} \rangle _{2}T^{-1/2}\) that an \(s\)-dependent
constant can be added to \(\k\) so that the resulting function, \(\k_{0}\), depends
smoothly on \(s\) and obeys \(|\frac{\partial}{\partial s}\k_{0}| \leq c_{0} \langle  \p_{X}  \rangle _{2}T^{-1/2}\).

 Reintroduce \(\sigma_{\perp}\) from Part 2 of Section \ref{sec:Db)}.  The 2-form \(w +
d(\sigma_{\perp} \q_{0})\) is equal to \(w\)
on \(\bbR \times\mathcal{Y}_{0} 
\) on the \(\ir > \rho_{*} +\frac{3}{4}  \varepsilon \) part of \(\bbR \times
 \mathcal{N}_\varepsilon\).  Meanwhile, it is equal
to \(w + d\q\) on the \(\ir < \rho_{*} +\frac{5}{8} \varepsilon \) part of \(\bbR \times
 \mathcal{N}_\varepsilon\).  Moreover, the norm of
the difference between this 2-form and \(w\) on the \(\ir >
\rho_{*} +\frac{1}{2}  \varepsilon \) part of  \(\mathcal{N}_\varepsilon\) is
bounded by \(c_{0} \langle  \p_{X}  \rangle _{2 }T^{-1/2}\), this being a consequence of the bounds in
the preceding paragraph for \(\q_{0}\).

Of interest in what follows is the 2-form on \(\bbR \times\mathcal{Y}_{0\varepsilon }\) given by
\begin{equation}\label{eq:(D.19)}
ds \wedge \b +  w + d(\sigma_{\perp}\q_{0}) \quad    \text{with \(\b =
\upsilon_{\diamond} +\sigma_{\perp}  \frac{\partial}{\partial s}\q_{0} + d(\sigma_{\perp}  \frac{\partial}{\partial s}\k_{0})\)}.
\end{equation}
 
This is a closed 2-form on  \(\mathbb{R}\times\mathcal{Y}_{0\varepsilon}\) which is
\(ds \wedge *w + w\) on \(\bbR  \times\mathcal{Y}_{0}\) and on the
\(\ir >\rho_{*} + \frac{3}{4}\varepsilon \) part of \(\bbR \times \mathcal{N}_\varepsilon\).   The
bounds given above on the norms of \(\k_{0}\), its
\(s\)-derivative, and on the norms of \(\q_{0}\),
\(d\q_{0}\) and \(\frac{\partial}{\partial s}\q_{0}\) imply the following:  There exists
\(c_{\diamond} > 1\) such that each \(s \in\bbR \) version of the
3-form \(\b \wedge (w + d(\sigma\q_{0}))\) is strictly positive on \(\mathcal{Y}_{0\varepsilon}\) if 
\begin{equation}\label{eq:(D.20)}
 \langle  \p_{X}  \rangle _{2} T ^{-1/2} \leq c_{\diamond}^{-1}.  
\end{equation}
 Assume in what follows that this bound holds.
 Granted (\ref{eq:(D.20)}), then Lemma \ref{lem:D.3} supplies a smooth, \(s\)-dependent
metric on \(\mathcal{Y}_{0\varepsilon}\) with the properties listed below.  The notation uses
\(\grg_{X}\) to denote the metric at any given \(s \in\bbR\).
\BTitem\label{eq:(D.21)}
\item The Hodge star of \(\grg_{X}\)  sends \(w +d(\sigma \q_{0})\)  to \(\b\). 
\item The metric \(\grg_{X}\) is
\(\grg_{*}\)  on \(\bbR \times  \mathcal{Y}_{0}\)  and on the \(r >
\rho_{*} + \frac{3}{4}\varepsilon  \) part of \(\bbR \times \mathcal{N}_\varepsilon\).
\item The metric \(\grg_{X}\)  is the metric in
(\ref{eq:(D.7)}) on the \(r < \rho_{*} +\frac{5}{8}\varepsilon  \) part of \(\bbR \times  \mathcal{N}_\varepsilon  \).
\item Given \(k\in \{1, 2, \ldots\}\), there exists
\(c_{k} > 1\)  such that the \(s< -104\)  and \(s > 104\) versions of
\(\grg_{X}\)  and their derivatives to order \(k \geq 1\)  differ by
at most \(c_{k}e^{-|s|/c_0}\) from the metric \(\grg_{-T}\) on \(Y_{-}\)  incarnation of
 \(\mathcal{N}_\varepsilon\) or \(\grg_{+T}\)  on the \(Y_{+}\) 
incarnation as the case may be.
\ETitem

By way of an explanation for the fourth bullet, this follows from (\ref{eq:(D.19)})
and the third bullet of Lemma \ref{lem:D.4} given the following fact:  The
derivatives to order \(k\) of any given coclosed eigenvector of \(*d\) on
\(Y_{-}\) or \(Y_{+}\) with \(L^{2}\)-norm 1 is bounded by a
polynomial function of the norm of the eigenvalue with coefficients
that are determined solely by the given metric.

\paragraph{Part 3:}  Let \(\grm_{T}\) denote the metric on \(X\)
that is equal to \(\grm\) on complement of the image of
(\ref{eq:(D.10)})'s embedding and whose pull-back to \(\bbR
\times\mathcal{Y} _{0\varepsilon }\) via this embedding is the metric \(ds^{2} + \grg_{X}\).
 This is a smooth metric on \(X\) whose pull-back by the embeddings from
the second and third bullets of (\ref{(A.9a,11)}) converge as \(s \to
-\infty \) to the metric \(ds^{2} + \grg_{-T}\) and converge as \(s \to \infty \) to the metric
\(ds^{2} + \grg_{+T}\).  These pull backs are also independent of \(s\) for \(|s | > 104\) at
points of the form \((s, p)\) if \(p\) is in either
\(\mathcal{Y}_{0}\), the \(\ir >\rho_{*} + \frac{3}{4}\varepsilon \) part of  \(\mathcal{N}_\varepsilon\),
or \(\mathcal{Y}_M\). 
 
Let  \(\omega_{T}\) denote the closed 2-form on \(X\) given by
 \(\omega\) on the complement of the image of (\ref{eq:(D.10)})'s
embedding and whose pull-back to \(\bbR\times\mathcal{Y}_{0\varepsilon }\) via this
embedding is the 2-form in (\ref{eq:(D.19)}).  The 2-form
 \(\omega_{T}\)  is closed.  This 2-form is also self-dual
when self-duality is defined by the metric \(\grm_{T}\), this
being a consequence of the first bullet in (\ref{eq:(D.21)}).  Let
\(w_{-T}\) and \(w_{+T}\) denote the
\(\grg_{-T}\) and \(\grg_{+T}\) harmonic 2-forms with de
Rham cohomology class that of \(c_{1}(\det(\bbS))\).
 Use \(*\) in what follows to denote either the \(\grg_{-T}\) or
\(\grg_{+T}\) Hodge dual.  The pull-backs of
 \(\omega_{T}\) via the embedding from the second bullet of (\ref{(A.9a,11)})
differs from \(ds \wedge *w_{-T} +w_{-T}\) in the \(C^{1}\) topology by at most
\(c_{T}e^{-|s|/c_T}\) with \(c_{T}\) > 1 being a constant.  The
pull-back via the embedding from the third bullet of (\ref{(A.9a,11)}) differs from
\(ds \wedge *w_{+T} + w_{+T}\) in the \(C^{1}\)-topology by at most
 \(c_{T}\) .  By way of an explanation, these bounds follow from the
 second and third bullet of Lemma \ref{lem:D.4}.  Keep in mind that  \(\omega_T\)  
obeys the second and third bullets of (\ref{(A.9a,11)}).  

Neither  \(\omega_{T}\) nor \(\grm_{T}\) are
likely to be \(s\)-independent where \(|s |\) is sufficiently
large.   This is a defect that is remedied in Parts 4-7 below.

\paragraph{Part 4:}  Both \(w_{-T}\) and
\(w_{+T}\) have nondegenerate zeros
on the components of \(Y_{-}\) and \(Y_{+}\) where they are
not identically zero, these being the components where
\(c_{1}(\det(\mathbb{S}))\) is not torsion. Let \(Y_{*} \subset Y_{-}\) denote such a
component and let \(p \in Y_{*}\) denote a zero of \(w_{-T}\).
Let \(B \subset Y_{*}\) denote a small radius ball centered on
\(p\) with the following properties:  The point \(p\) is the only zero of
\(w_{-T}\) in the closure of \(B\); and \(B\) is disjoint from
\(\mathcal{Y}_{0}\) and from the \(\ir >\rho_{*} + \frac{3}{4}\varepsilon  \)
  part of  \(\mathcal{N}_\varepsilon\).  Since
\(w_{-T}\) vanishes transversely at \(p\), there exists
\(L_{0} > 1\) such that each \(s <-L_{0}\) version of \(w_{-T} + d\q\) vanishes
transversely in the closure of \(B\) at a single point.  Let
\(p_{s}\) denote this point.  Note in particular that \(\dist \, (p,
p_{s}) \leq c_{-}\).   Granted that \(\dist \, (p, p_{s}) \ll 1\)
for \(s \ll -1\), there exists \(s_{0}> 1\) such that \(\dist \, (p, p_{s})\) is less than  
\(\frac{1}{8}\) times the radius of \(B\) when \(s \leq - s_{0}\).  This
being the case, there exists \(L_{0} >
s_{0}\), \(c_{-} > 1\) and a family
of diffeomorphisms of \(Y_{*}\) parametrized by
\((-\infty, -L_{0}]\) with the properties in the list that
follows.  The list uses \(\Psi_{s}\) to denote the
diffeomorphism labeled by a given \(s \in (-\infty,-L_{0}]\).
\BTitem\label{eq:(D.22)}
\item If \(s > - L_{0}-1\),  then \(\Psi_{s}\) is the identity map.
\item Every \(s \in (-\infty, -L_{0}]\)
version of \(\Psi_{s}\)  is the identity where
\(\dist \, (\cdot, p) > 2 \dist \, (p, p_{s})\).
\item Every \(s \in (-\infty, -L_{0}]\) version of \(\Psi_{s}\)  differs
from the identity in the \(C^{10}\)-topology by at most \(e^{-|s|/c}\).
\item \(\Psi_{s}(p) =p_{s}\)  when \(s < - L_{0} -2\).
\ETitem

This family of diffeomorphisms defines a diffeomorphism of \(X\) which is
the identity on the \(s > -L_{0} -1\) 
part of \(X\), and on the image in \(X\) of \((-\infty,-L_{0}]\times(Y_--Y_{*})\)
via the diffeomorphism in the second bullet of (\ref{(A.9a,11)}). This
diffeomorphism is defined on the image of \((-\infty,-L_{0}]\times Y_{*}\) via the second
bullet of (\ref{(A.9a,11)}) by that of \((-\infty, -L_{0}] \times
Y_{*}\) that sends a given point \((s, q)\) to \((s,
\Psi_{s}(q))\).  Use \(\Psi_{p}\) to
denote this diffeomorphism of \(X\).  Various versions of this
diffeomorphism are defined by the zeros of \(w_{-T}\) on the
components of \(Y_-\) where \(c_1(\det\, (\bbS))\) is not torsion. These
diffeomorphisms pairwise commute. Use \(\Psi\) to denote their composition.

Introduce \(\grm_{T0}\) to denote \(\Psi^*\grm_{T}\) and \(\omega_{T0}\) to denote
 \(\Psi^*\omega_T\).  The 2-form \(\omega_{T0}\)  is closed and it is self-dual if the
notion of self-duality is defined using \(\grm_{T0}\).  The form
 \(\omega _{T0}\) can be written on \((-\infty, -L_0]\times Y_*\) as \(ds \wedge(*w_{-T} + \n) + (w_{-T} + \mm)\) where \(\n\) 
and \(\mm\) have \(C^{1}\)-norm less than \(c_{-}e^{-|s|/c_-}\) 
 and both vanish on \((-\infty, -1] \times
\mathcal{Y} _{0}\) and on the \(\ir >\rho_{*} + \frac{3}{4}\varepsilon \) part of \((-\infty, -1] \times
 \mathcal{N}_\varepsilon\).  By way of notation,
\(c_{-}\) denotes here and in what follows a constant that
is greater than 1.  Its value can increase between successive
appearances.   Note that the fact that  \(\omega_{T0}\) 
is closed requires that \(d\n\) equals \(\frac{\partial}{\partial s}\mm\).

The pull-back of \(w_{-T} + \mm\) to each constant \(s\) slice of
\((-\infty, -1] \times Y_{*}\) defines the same
cohomology class as \(w_{-T}\).  This implies in particular
that \(\mm = d\uu\) with \(\uu\) being an \(s\)-dependent 1-form on
\(Y_*\).  Any \(s\)-dependent, closed 1-form can be added
to \(\uu\) without changing \(d\uu\), and this fact is used to choose \(\uu\) so that the
conditions that follow hold.      

\BTitem\label{eq:(D.23)}
\item The 1-form \(\uu\) is zero on \(\mathcal{Y}_0\) and on the   \(\ir >\rho _{*} + \frac{3}{4}\varepsilon  \) part of
 \(\mathcal{N}_\varepsilon\).
\item The \(C^{2}\)  norm of \(u\)  is less
than \(c_{-}e^{-|s|/c_-}\).
\item  Let \(p\)  denote a zero of \(w_{-T}\) in \((-\infty,
  -L_0]\times Y_*\). Then \(|u|\) and the norm of \(u\)'s covariant
    derivative  along \(\frac{\partial}{\partial s}\) at any \(s\in (-\infty,
  -L_0]\times Y_*\)  is bounded by \(c_{-} \dist \, (\cdot,p)^{2}e^{-|s|/c_-}\).
\ETitem

To explain how the third bullet can be satisfied, let \(p\) again denote
a zero of \(w_{-T}\) .  Use the metric \(\grg_{-T}\)
to construct a Gaussian coordinate chart centered at \(p\) so as to
identify \(B\) with a small radius ball in  \(\mathbb{R}^3\).
 The corresponding coordinate map to \(\mathbb{R}^3\) 
is denoted by \(x\) or \((x_{1}, x_{2},x_{3})\).  Write the two form  \(\omega_{T0}\) 
as \(ds \wedge (*w_{-T} + \n) +(w_{-T} + \mm)\).  The 2-form \(\mm\) appears in these
coordinates as
\begin{equation}\label{eq:(D.24)}
\mm = \frac{1}{2}\gro^{ij} x^{i}   \varepsilon^{jnm} dx^{n} dx^{m}+\cdot\cdot\cdot
\end{equation}
where the summation convention over repeated indices is used.  The
unwritten terms in (\ref{eq:(D.24)}) are \(\mathcal{O} (|x|^{2})\).  What is
denoted by \(\{\varepsilon^{jnm}\}_{1\leq j,n,m\leq 3}\) is
anti-symmetric with respect to interchanging indices and so defined by
the rule \(\varepsilon^{123} = 1\).  Meanwhile, \(\{\gro^{ij}\}_{i, j =  1,2,3}\) are the
entries of a traceless, \(s\)-dependent matrix whose norm and that of its
\(s\)-derivative are at most \(c_{-}e^{-|s|/c_1}\).  The matrix is traceless because \(\mm\) is closed.  The fact that this
matrix \(\gro\) is traceless implies that \(\mm\) on \(B\) can be written as
\(d\uu_{B}\) with \(\uu_{B} = \frac{1}{6}\gro^{ij}x^{i} x^{n}   \varepsilon^{jnm} dx^{m} +\cdot\cdot\cdot \) where
the unwritten terms are \(\mathcal{O} (|x|^{3})\).  Since \(\uu -
\uu_{B} = d\p\) on \(B\), it follows that \(\uu\) can be modified with no
change near the boundary of \(B\) so that \(\uu = \uu_{B}\) on a small
radius ball in \(B\) centered at \(p\).

 \paragraph{Part 5:}  Fix \(L_{1 }>L_{0} + 1\) and let \(\chi_{1}\) denote the
function on \(\bbR\) given by \(\chi(- L_{1}  - s)\).
 This function equals zero where \(s < - L_{1}
 - 1\) and it equals 1 when \(s >   - L_{1}\).
 Use  \(\chi_{1}'\) to denote the derivative of
\(\chi_{1}\).  The function \(\chi_{1}\)
and the 2-form  \(\omega_{T0}\)  are used next to define the
2-form on \((-\infty, -L_0]\times Y_*\) to be denoted by  \(\omega_{T1}\) .  This
2-form is  \(\omega_{T0}\) on the \(s > -
L_{1}\) part of \((-\infty, -L_0]\times Y_*\), and it is given where \(s \leq
-L_{1}\) by the formula that follows for its pull-back via
the embedding from (\ref{(A.9a,11)})'s second bullet:
\begin{equation}\label{eq:(D.25)}
 \omega _{T1} = ds \wedge (*w_{-T} +\chi_{1} \n + \chi_{1}'\uu ) +
w_{-T} + \chi_{1} d\uu .
\end{equation}
 
The 2-form  \(\omega_{T1}\) is a closed 2-form on \((-\infty, -L_0]\times Y_*\).  The
remainder of this part of the subsection and Part 6 describe a metric
on the \(s \in (-\infty, - L_{0}]\times Y_*\) that
makes  \(\omega _{T1}\) self-dual.  This new metric is
equal to \(\grm_{T0}\) where \(s \geq - L_{1} + 1\) and it is equal to \(ds^{2} + \grg_{-T}\) where
\(s < - L_{1} -2\).  This new metric is denoted
below by \(\grm_{T1}\). 
 The five steps that follow describe the metric
\(\grm_{T1}\) at points in \((-\infty, -L_{0}]\times Y_{*}\)
 that project to \(Y_{*}\) near the zero locus of \(w_{-T}\).   
 
\paragraph{\it Step 1:}  The 2-form \(w_{-T}\) and the 1-form
\(*w_{-T}\) on \(B\) can be written using the Gaussian
coordinates \((x_{1}, x_{2}, x_{3})\) on \(B\) as
\begin{equation}\label{eq:(D.26)}
  w_{-T} =\frac{1}{2}
\textsc{a}^{ij}  x^{i} \varepsilon^{jnm} dx^{n}
dx^{m}+\cdot\cdot\cdot\quad 
\text{and}\quad *w_{-T} =\textsc{a}^{ij}x^{i}dx^{j} +\cdot\cdot\cdot 
\end{equation}
with summations over repeated indices implicit.  The various \(i\), \(j
\in \{1, 2, 3\}\) versions of \(\textsc{a}^{ij}\) in  (\ref{eq:(D.26)})
are the entries of an invertible matrix, this denoted by \(\textsc{a}\).
 The unwritten terms in (\ref{eq:(D.26)}) vanish to order \(|x|^{2}\).  The fact that
\(w_{-T}\) is closed implies that \(\textsc{a}\) is traceless
and the fact that \(*w_{-T}\) is self-dual implies that
\(\textsc{a}\) is symmetric. 

  The unwritten terms in (\ref{eq:(D.26)}) are incorporated using the notation
whereby \(w_{-T}\) and \(*w_{-T}\) on \((-\infty, - L_{0}] \times B\) are written as 
\begin{equation}\label{eq:(D.27)}
ds \wedge (f_{1} \hat{e}^{1} +f_{2} \hat{e}^{2} + f_{3}\hat{e}^{3}) + f_{1} \hat{e}^{2} 
\wedge \hat{e}^{3} + f_{2} \hat{e}^{3 }\wedge \hat{e}^{1} +f_{3}  \hat{e}^{1 }\wedge \hat{e}^{2},
\end{equation}
where \(\{\hat{e}^{k}\}_{1\leq k\leq 3}\)
denotes a \(\grg_{-T}\)-orthonormal set for \(T^*B\) with
\(\{\hat{e}^{k} = dx^{k} +\mathcal{O} (|x|^{2})\}_{1\leq k\leq 3}\); and where \(\{f_{k}\}_{1\leq k\leq 3}\)
are functions with \(\{f_{k} =\sum _{1\leq i\leq 3}\textsc{a}^{ik}x^{i} + \mathcal{O} (|x|^{2})\}_{ 1\leq k\leq 3}\).
  Note in particular that these are such that \(df_{1}\wedge df_{2} \wedge df_{3}
> \frac{1}{2}\det(\textsc{a})\) on a concentric ball in \(B\) centered at the origin.
 This ball is denoted by \(B'\).  It is assumed in what follows that
\(L_{0}\) is chosen so that  \(\omega _{T0} =\omega _{T}\) on the complement of a concentric ball in
\(B'\) with radius one fourth that of \(B'\).  In particular, it is assumed
that (\ref{eq:(D.22)})'s diffeomorphism \(\Psi_{s}\)
is the identity for all \(s\) on a neighborhood in \(B\) of \(B-B'\).

\paragraph{\it Step 2:}  The \(\Psi\)-pull back of \(\{ds, \hat{e}^{1}, \hat{e}^{2}, \hat{e}^{3}\}\) is
\(\grm_{T}\)-orthonomormal.  The \(\Psi\)-pull back of \(ds\) is
\(ds\).  Meanwhile, \(\Psi\) can be chosen so that 
\begin{equation}\label{eq:(D.28)}
\Psi^*\hat{e}^{k} = \hat{e}^{k} +\sum_{1\leq k\leq 3} \grp^{k} ds +\sum_{1\leq j\leq 3} \grp^{kj}\hat{e}^{j},
\end{equation}
where \(\sum_{1\leq k\leq 3}|\grp^{k}| \leq c_{-}e^{-|s|/c_-}\)   and \(\sum_{1\leq k,j\leq 3}|\grp^{kj}| \leq c_{-} |x|e^{-|s|/c_-}\)  
   when \(s < -L_{0}  -1\).  This is done by
defining (\ref{eq:(D.22)})'s diffeomorphism
\(\Psi_{s}\) using the Gaussian coordinates in (\ref{eq:(D.25)}) by the rule \(x \mapsto
\Psi_{s}(x) = x + p_{s}\) at points \((s, x)\)
with \(|x| < \frac{3}{2}|p_{s}|\) and \(s < -L_{0}  - 1\).  Use \(\{\hat{e}^k_s\}_{1\leq k\leq 3}\)
to denote \(\{\Psi^*\hat{e}^{k}\}_{1\leq k\leq 3}\).
 Granted this notation, the 2-form  \(\omega '_{T0}\) near \(p\) can be written as
\begin{equation}\label{eq:(D.29)}
 \omega _{T0}' = ds \wedge(f_{s_1}\hat{e}^{1}_s +f_{s_2}\hat{e}^{2}_s +f_{s_3}\hat{e}^{3}_s) +
f_{s_1}\hat{e}^2_s \wedge \hat{e}^{3}_s +f_{s_2} \hat{e}^{3}_s \wedge \hat{e}^{1}_s +
f_{s_3} \hat{e}^{1}_s \wedge \hat{e}^{2}_s,
\end{equation}
 where \(\{f_{s_k} = \Psi_{s}^*f_{k}\}_{1\leq k\leq 3}\).

\paragraph{\it Step 3:}  Introduce
\(\{e^{k}_{s_\chi} =\hat{e}^{k} + \chi_{0}\sum_{1\leq j\leq 3}\grp^{kj}\hat{e}^{j}\}_{1\leq k\leq 3}\).  Use this
\(s\)-dependent basis to write the (\ref{eq:(D.25)})'s 2-form
\(w_{-T} + \chi_{1} d\uu\) on \(B'\) as 
\begin{equation}\label{eq:(D.30)}
w_{-T} + \chi_{1} d\uu =f_{s_{\chi 1}}e^{2}_{s_ \chi}\wedge e^{3}_{s_ \chi} +
f_{s_\chi 2}e^{3}_{s_ \chi}\wedge e^{1}_{s_ \chi} +f_{s_\chi 3}e^{1}_{s_ \chi}\wedge e^{2}_{s_ \chi},
\end{equation}
where \(\{f_{s_\chi k}\}_{1\leq k\leq 3}\)
are smoothly varying functions of \(s\) and the coordinate \(x\) with the
property that \(f_{s_\chi (\cdot)}= f_{(\cdot)}\) when \(s < -L_{1}  -1\) and
\(f_{s_\chi (\cdot)}= f_{s (\cdot)}\) when \(s> -L_{1}\).  This depiction can be derived
from the fact that \(\{f_{k}\}_{1\leq k\leq 3}\) generate \(C^{\infty}(B')\).  Note that
\(f_{s_\chi k} =f_{k} +\cdots\)
with the unwritten terms such that their norms are bounded by
\(c_{-} e^{-|s|/c_-}|x|\) and such that their first derivatives have norms
bounded by \(c_{-}e^{-|s|/c_-}\).   This implies in particular that the functions
\(\{f_{s_\chi k}\}_{1\leq k\leq 3}\) also generate \(C^{\infty}(B')\) and that \(df_{s_\chi 1} \wedge
df_{s_\chi 2} \wedge df_{s_\chi 3}>    \det\, (\textsc{a})\) on \(B'\) when \(L_{0} >c_{-}\).

 The 1-form \(*w_{-T} + \chi_{0} \n +\chi_{0}'\uu\) can be written schematically on
\((-\infty, - L_{0}] \times B'\) using the basis \(\{e^{k}_{s_\chi}\}_{1\leq k\leq 3}\) as
\begin{equation}\label{eq:(D.31)}
*w_{-T} + \chi_{0} \n +\chi_{0}'\uu =
\sum_{1\leq k,i\leq 3}f_{s_\chi k} \textsc{c}_{ki} \, e^{i}_{s_\chi},
\end{equation}
with \(\{\textsc{c}_{ki}\}_{1\leq i,k\leq 3}\)
denoting a matrix of smooth functions of \(s\) and the coordinate \(x\).
 Given that the functions \(\{f_{s_\chi k}\}_{1\leq k\leq 3}\)
also generate \(C^{\infty}(B')\), such a depiction
follows because both \(\n\) and \(\uu\) vanish at \(p\).   Keep in mind for what
follows that the matrix with coefficients \(\{\textsc{c}_{ki}\}_{1\leq i,k\leq 3}\) 
differs from the identity matrix by at most \(c_{-}e^{-|s|/c_-}\).

\paragraph{\it Step 4:}  A particular set of three smooth functions of \(s \in
(-\infty,  -L_{0}]\) and the coordinate \(x\) is specified
momentarily.  Let \(\{\grq^{k}\}_{1\leq k\leq 3}\) denote any given set
of such functions.  Use this set to define 1-forms \(\{\hat{e}^{k}_{s_\chi}\}_{1\leq k\leq 3}\)
on \((-\infty,  -L_{0}] \times B'\) by the rule
\begin{equation}\label{eq:(D.32)}
\hat{e}^{k}_{s_\chi} =e^{k}_{s_\chi} -\grq^{k} ds. 
\end{equation}
 
Given the formula in (\ref{eq:(D.31)}) and (\ref{eq:(D.32)}), it follows that
 \(\omega_{T1}\)  on \((-\infty, -L_{0}]\times B'\) can be written using \(\{\hat{e}^{k}_{s_\chi}\}_{1\leq k\leq 3}\) as
\begin{equation}\label{(eq:(D.33)}
ds \wedge\big(f_{s_\chi k}(\textsc{c}_{ki} + \varepsilon^{kni}
\grq^{n}) \big) \hat{e}^{i}_{s_\chi} +
 \frac{1}{2}\, f_{s_\chi k} \varepsilon^{knm} 
\hat{e}^{n}_{s_\chi} \wedge \hat{e}^{m}_{s_\chi}.
\end{equation}
This equation uses the summation convention over repeated indices.

\paragraph{\it Step 5:}  The set \(\{\grq^{k}\}_{1\leq k\leq 3}\) is
introduced for the following reason:  There is a unique choice for
\(\{\grq^{k}\}_{1\leq k\leq 3}\)  that makes the matrix with entries \(\{\textsc{c}_{ki} +
\varepsilon^{kni}\grq^{n}\}_{1\leq i,k\leq 3}\) a symmetric
matrix, this being \(\{\grq^{k} = \frac{1}{2} \varepsilon^{kin} \textsc{c}_{ni}\}\).
  This choice is used in what follows.  With this choice understood,
a metric is defined on \((-\infty, -L_{0}] \times B'\)
by the following rules:
\BTitem\label{eq:(D.34)}
\item \(ds\)  has norm 1  and it is orthogonal to
\(\{\hat{e}^{k s_\chi}\}_{1\leq k\leq 3}\).
\item Given \((i, k) \in \{1, 2, 3\}\), then the inner product between
\(\hat{e}^{k}_{s_\chi }\) and \(\hat{e}^{i}_{s_\chi}\) is \(\textsc{c}_{ki} +\varepsilon^{kni} \grq^{n}\).  
\ETitem
The inner product defined by the second bullet is positive definite if
\(L_{0} > c_{0}c_{-}\) because of the afore-mentioned fact that the
matrix defined by \(\{\textsc{c}_{ki}\}_{1\leq i,k\leq 3} \)
differs by at most \(c_{-}  e^{-|s|/c_-} \) from the identity matrix.  
 
The metric just defined is the metric \(\grm_{T0}\) when \(s
> -L_{1}\); and it is the metric
\(ds^{2} + \grg_{-T}\) when \(s < -L_{1}  - 1\).  Moreover, the 2-form
 \(\omega_{T1}\)  is self-dual on \((-\infty,-L_{0}] \times B'\) when self-duality is defined by this
metric.   Denote this metric by \(\grm_{T1p}\).  
 
Let \(B'' \subset B'\) denote the concentric ball whose radius is one
half that of \(B'\).  The desired metric \(\grm_{T1}\) is defined to
equal \(\grm_{T1p}\)  on \((-\infty, - L_{0}]\times B''\).

\paragraph{Part 6:}  Use \(U\) to denote the union of the various versions
of the ball \(B''\).  The two steps that follow directly describe the
metric \(\grm_{T1}\) on \((-\infty, - L_{0}]\times (Y_{*}-U)\).  

\paragraph{\it Step 1:}  This step describes a metric on \((-\infty, -
L_{0}] \times (Y_*-U)\) to be denoted by \(\grm_{T1 \diamond}\).  The metrics \(\grm_{T1}\) and \(\grm_{T1 \diamond}\)
agree on the product of \((-\infty, - L_{0}]\) with the
complement in \(Y_*\) of the union of the various versions
of the ball \(B'\).  The definition of this metric
\(\grm_{T1 \diamond}\) assumes that \(L_{0} > c_{\diamond}\) with
\(c_{\diamond}\) such that  \(\omega  =w_{-T} + \chi_{1} d\uu\) and \(\upsilon =
*w_{-T} + \chi_{1} \n +\chi_{1}' \uu\) from (\ref{eq:(D.25)}) obey \(\upsilon \wedge\omega  > 1/c_{\diamond}\) on
\((-\infty, - L_{0}] \times (Y_*-U)\).   The existence of \(c_{\diamond}\) follows from (\ref{eq:(D.23)}).
   Let \(p\) denote a zero of \(w_{-T}\) and let \(B_{\diamond} \subset B'\) denote the concentric ball
whose radius is three quarters that of \(B'\).  Use \(V\) to denote the union
of the various versions of \(B_{\diamond}\). Invoke Lemma
\ref{lem:D.3} on \((-\infty, -L_{0}] \times (Y_*  - U)\) using
 \(\omega\) and \(\upsilon\) to obtain a smooth family of metrics on
\(Y_*-V\) parametrized by \((-\infty, -L_{0}]\) with the properties listed in the upcoming (\ref{eq:(D.34)}).
 The notation uses \(\grg_{\diamond}\) to denote any given \(s
\in (-\infty, - L_{0}]\) member of the family.  To
explain more of the notation,  note first that pull-backs of \(\grm\) and
Part 4's metric \(\grm_{T0}\) via the embedding
from the second bullet of (\ref{(A.9a,11)}) agree on \((-\infty,
-L_{0}] \times (Y_* -U)\).  In particular, the pull-back of \(\grm_{T0}\) to this part of
\((-\infty, - L_{0}] \times Y_*\) can be written as \(ds^{2} + \grg_{X}\) with
\(\grg_{X}\) denoting here a smooth, \(s\)-dependent metric on
\(Y_*-U\).  This metric \(\grg_{X}\) is the metric \(\grg_{-T}\) on \(\mathcal{Y}_{M}-U\) and
it is the metric from (\ref{eq:(D.19)}) on \(\mathcal{Y}_{0\varepsilon}\).  
\BTitem\label{eq:(D.35)}
\item Each \(s \in (-\infty, - L_{1} -1]\)  version of \(\grg_{\diamond}\)  is
\(\grg_{-T}\) and each \(s \in [-L_{1},- L_{0}]\)  version is the corresponding version of \(\grg_{X}\).    
\item The \(\grg_{X}\)-Hodge dual of the 2-form \(w_{-T} +\chi_{1} d\uu\)  on
\(Y_*-V\)  is the 1-form \(*w_{-T} + \chi_{1} \n +\chi_{1}'\uu\).
\ETitem
 
 The metric \(\grm_{T1\diamond}\) on \((-\infty, - L_{0}] \times(Y_*-U)\) is defined to be \(ds^{2}
+\grg_{\diamond}\).  It follows directly from the second bullet in (\ref{eq:(D.35)}) that the 2-form  \(\omega _{T1}\) is self-dual on \((-\infty, - L_{0}] \times(Y_*-V)\) when the notion of self duality is defined
using the metric \(\grm_{T1\diamond}\).

  \paragraph{\it Step 2:}  Let \(p\) denote a zero of \(w_{-T}\).  The metrics \(\grm_{T1 \diamond}\) and
\(\grm_{T1p}\) are both metrics on \((-\infty, -L_{0}] \times (B'-B_{\diamond})\).
 The 2-form  \(\omega _{T1}\) is self-dual on \((-\infty,- L_{0}] \times (B'-B_{\diamond})\)
when the latter notion is defined by either metric.  Use
\(z_{\diamond}\) and \(z_{p}\) to denote the respective \(\grm_{T1 \diamond}\) and
\(\grm_{T1p}\) norms of  \(\omega _{T1}\).  Since \(\omega _{T1} \wedge  \omega _{T1}> c_{-}^{-1}\) here, there
is a  \(\omega _{T1}\)-compatible almost complex structures for
\((-\infty, - L_{0}] \times(B'-B_{\diamond})\), these denoted by \(J_{\diamond}\) and \(J_{p}\), such that 
\begin{equation}\label{eq:(D.36)}
\grm_{T1 \diamond} =z_{\diamond}^{-1}\omega _{T1}( \cdot,J_{\diamond}(\cdot)) \quad   \text{and}\quad
 \grm_{T1p}  = z_{p}^{-1}
 \omega _{T1}( \cdot,J_{p}(\cdot)).
\end{equation}
  As the space of  \(\omega _{T1}\)-compatible almost
complex structures on \((-\infty, - L_{0}] \times(B'-B_{\diamond})\) is contractible, there exists such
an almost complex structure with two properties, the first of which is
as follows:  The almost complex structure is \(J_{p}\) at
points with \(B'-B_{\diamond}\) component in a neighborhood of the boundary of the closure of
\(B_{\diamond}\); and it is \(J_{\diamond}\) at points with \(B'-B_{\diamond}\) component in the \(B'\)
part of a neighborhood of the boundary of the closure of \(B'\) in \(B\).  To
state the second property, keep in mind that \(J_{\diamond} = J_{p}\) in some neighborhood
of \((-\infty, - L_{1 } - 1] \times(B'-B_{\diamond})\) and also in some neighborhood of
\([- L_{1}, -L_{0}] \times(B'-B_{\diamond})\).  What follows is the second
property:  The new almost complex structure is \(J_{\diamond}\) and thus \(J_{p}\) in slightly
smaller neighborhood of \((-\infty, - L_{1 }- 1]\times
(B'-B_{\diamond})\) and \([-L_{1}, -L_{0}] \times(B'-B_{\diamond})\).
Use \(J_{*}\) to denote an almost complex structure of the sort just described.

   Fix a smooth, strictly positive function on \((-\infty, -L_{1 } - 1]
  \times(B'-B_{\diamond})\) that is equal to \(z_{\diamond}\) where \(J_{*} =J_{\diamond}\) and equal to \(z_{p}\) where
\(J_{*} = J_{p}\).  Denote this function by
\(z_{*}\).  Use \(J_{*}\) and \(z_{*}\) to define the metric \(\grm_{ T 1}\) on
\((-\infty, - L_{1 } - 1] \times(B'-B_{\diamond})\) by the rule \(\grm_{T1} =z_{*}^{-1}\omega _{T1}( \cdot ,J_{*}(\cdot))\). 
This metric smoothly extends the metrics defined in Step 1 and in
Step 5 of Part 5 and it has all of the desired properties.

\paragraph{Part 7:}  Let \(Y_{*} \subset Y_{-}\) now denote a component
where \(w_{-T}\) is identically zero, thus a component where
\(c_{1}(\det(\mathbb{S}))\) is torsion. Suppose that \(L> 1\)  has been chosen. Let
\(\omega_{T0}\) now denote the pull-back of \(\omega_{T}\) to \((-\infty, -L] \times
Y_{*}\) via the embedding from the second bullet of (\ref{(A.9a,11)}). It follows from Lemma \ref{lem:D.4}
that the \(C^{1}\) norm of \(\omega_{T0}\) is bounded by \(c_{0}
\langle\p_{X}\rangle_{2}e^{-|s|/c_0}\). The 2-form \(\omega_{T0}\) is exact on \((-\infty,
-L]\times Y_{*}\), it can be written as \(ds\wedge
\frac{\partial}{\partial s}u +du\) with \(d\) denoting here the exterior
derivative along the constant \(s\) slices of \((-\infty, -L]\times Y_{*}\) and with
\(u\) denoting a smooth, \(s\)-dependent 1-form on \(Y_{*}\) with \(|u|\), \(|du|\) and |\(u|\) bounded
by \(c_{0}\langle \p_{X}\rangle_{2}\, e^{-|s|/c_0}\).

With the preceding understood, fix \(L_{tor} >L +4\) and let \(\chi_{*}\)
denote the function on \(\mathbb{R}\) defined by the rule \(s
\mapsto\chi(-L_{tor}+ 3 -s)\). This function equals 1 where \(s >-L_{tor}+ 3\) and it equals 0
where \(s <-L_{tor} + 2\). Use \(\chi_{*}\) to define a self-dual form on
\((-\infty, -L]\times Y_{*}\) by the following rules: This form is equal
to \(\omega_{T0}\) on \([-\infty,-L_{tor}+ 4, -L] \times Y_{*}\), it is
identically 0 on \([-\infty,-L_{tor}]\times Y_{*}\) and it is equal to
\(\chi_{*}(ds\wedge \frac{\partial}{\partial s}u +du)\) on \([-L_{tor}, -L_{tor}+ 4]\times
Y_{*}\). Denote this 2-form by \(\omega_{T1}\).

The 2-form \(\omega_{T1}\) can be written as \(ds \wedge*w_{*}+w_{*}\) with
\(w_{*}=d(\chi_{*}u)\) with it understood again that \(d\) here denotes the exterior derivative along
\(Y_{*}\). Let \(\chi_{*}'\) denote the derivative of the function \(s\mapsto\chi_{*}(s)\).
The 2-form \(w_{*}\) on \([-L_{tor}, -L_{tor}+ 4]\times Y_{*}\) can be written
as \(d\b\)  with \(\b =\chi_{*}'u+\chi_{*}u\). Note in particular that
\(|\b|\leq c_{0}\c\langle \p_{X}\rangle_{2}e^{-|s|/c_0}\).

Fix \(\c >1\). The bound just given for \(|\b|\) leads to the following conclusion: Fix
\(r  > 1\). Then \(|\b|\leq r^{-10}\) if \(L_{tor} >c_{0}(|\ln(\langle\p_{X}\rangle_{2}|+ \ln r)\). 

\paragraph{Part 8:} Define the 2-form \(\omega_{T*}\) on the \(s \leq 0\) part of \(X\) as
follows: This 2-form is equal to \(\omega_{T}\) where \(s\in [-L, 0]\).
Meanwhile, its pull-back to each component of the \(s < -1\) part of \(X\) via the
embedding from the second bullet of (\ref{(A.9a,11)}) is the corresponding version
of the 2-form \(\omega_{T1}\).
Modulo notation, what is said in Parts 4-7 can be repeated for the
\(s>0\) part of \(X\) to extend the definition of
\(\omega_{T*}\) and the metric \(\grm_{T*}\) to the whole of \(X\). The form
\(\omega_{T*}\) is self-dual if the latter notion is defined by \(\grm_{T*}\).
This construction has the following additional property: Suppose that
\(\p_{X}\) obeys (\ref{eq:(D.20)}). Fix \(\c >c_{0}\). If \(\ir > 1\) has been chosen to be greater than a purely
\(\c\)-dependent constant, then the \((L =\c, L_{tor} =\c \ln r)\) version of
\(\grm_{T*}\) and \(\omega_{T*}\) obey the constraints given by
(\ref{eq:(A.12,15a)}), (\ref{eq:(A.13b)}), (\ref{eq:(A.13c)}),
(\ref{eq:(A.15b)}) and the \((\c,\r)\) version of (\ref{eq:(A.16)}). 
Here, the closed 1-form \(\upsilon_X\) can be chosen so that it is
\(s\)-independent and \(\upsilon_X=*w_{\pm T}\) over
constant \(s\)-slices of \(X\) where \(|s|>L-4\). The bounds in items
4b), 4d), and 5c) of (\ref{eq:(A.16)}) follows from the bounds on
\(\uu\) in (\ref{eq:(D.23)}) and those for \(\b\) in Part 7 above.

\paragraph{Part 9:} The happy conclusions of Part 8 are contingent on the existence of a closed
2-form, \(\p_{X}\), on \(X\) with the following properties: The de
Rham cohomology class of \(\p_{X}\) is \(c_{1}(\det(\mathbb{S}))\),
  it equals \(w_{-}\) where \(s < -102\), it equals \(w_{+}\) where \(s
> 102\), and it obeys the bound in (\ref{eq:(D.20)}).  
 
  The subsequent four steps in this part of the subsection construct
\(\p_{X}\) on various parts of \(X\).  These constructions are
used in Part 11 and they are also used in the proofs of Proposition
\ref{prop:A.10a} and \ref{prop:A.10b}. 
 
  \paragraph{\it Step 1:}   This step first states and then proves a
lemma that supplies a crucial tool for what is to come.  

\begin{lemma}\label{lem:D.5} 
Let \(U\) denote a 3-manifold and let \(V\subset U\) denote an open
set with compact closure in \(V\).  Given the data set
consisting of \(U\), \(V\), and a Riemannian metric on \(U\), there
exists \(\kappa  >1\) with the following significance:  Let \(u\) 
denote a closed, exact 2-form on \(U\).  There is a 1-form on \(U\), this denoted
by \(\q\), with \(\int_V|\q |^2\leq\kappa  \int_U|\uu|^2\) and such that \(d \q  =\uu\).
\end{lemma}
To set the notation used below, the \(L^{2}\)-norm of a
function or differential form over a given set \(W \subset  U\) is denoted by
 \(\| \cdot \|_{W}\).

\pf 
The set \(V\) has a finite cover by Gaussian
coordinate balls with centers in \(U\) with the property that the mutual
intersection of balls from this cover is either empty or convex.  This
cover can also be chosen so that each ball has the same radius and such
that no ball intersects more than \(c_{0}\)  others.  The
minimal number of balls in such a cover, their common radius and the
combinatorical properties of the mutual intersections is determined a
priori by \(U\), \(V\) and the metric.  Let \(\grU\) denote such a cover and let
 \(\sigma\) denote the radius of its constituent balls. 

Let \(B  \in \grU\).  The fact that \(B\) is convex
can be used to write  \(\uu\)  on \(B\) as
 \(\uu  =d \q_{B} \) where \(\| \q_{B} \|_{B} \leq c_{0} \sigma \|  \uu \|_{B}\).
  Let \(B\) and \(B'\) denote two sets from \(\grU\) .  Then
 \(d \q_{B} -d \q_{B'} = 0\) on their intersection, and so
\( \q_{B} -\q_{B'} =d \k_{BB'} \) with \(\k_{BB'} \)
being a function on \(B' \cap B\).  It follows that  \( \|d \k_{BB'} \|_{B' \cap B} 
\leq c_{0} \sigma  (  \|\uu \|_{B} +  \| \uu \|_{B'} )\).
 Changing  \(\k_{BB'} \) by a constant if needed produces a version with
\( \| \k_{BB'} \|_{B \cap B'} \leq c_{0} \sigma  \| d \k_{BB'} \|_{B' \cap B} \)
and thus   \(\|d \k_{BB'} \|_{B' \cap B} \leq c_{0} \sigma ^{2} (  \|\uu \|_{B} +  \|\uu \|_{B'} )\).

 Now suppose that \(B\), \(B'\), and \(B''\) are from
 \(\grU\)  with a point in common.  Let \(\c_{BB'B''}\) 
denote \(\k_{BB'} + \k_{B'B''} + \k_{B''B}\) . This \(\c_{BB'B''}\) 
is constant and the collection of such numbers is a \v{C}ech cohomology
cocycle whose cohomology class gives the class of
 \(\uu\)  via the de Rham isomorphism.
 It follows that this cocycle is zero, and so
 \(\c_{BB'B''}=\c_{BB'} +\c_{B'B''} +\c_{B''B}\) 
with each term being constant.  Noting that
\(| {\c_{BB'B''}} |\leq c_{0} \sigma ^{-1} (  \|\uu \|_{B} +  \|\uu \|_{B'} +  \|\uu \|_{B''} )\),
it follows that \(| {\c_{BB'}} | \leq c_{*} \sigma ^{-1} \sup_{B''\in U :B''\cap B' \cap B \neq \emptyset} (  \|\uu \|_{B} +  \|\uu \|_{B'} +  \|\uu \|_{B''} )\)
with \(c_{*} \geq 1\) determined a priori by the combinatorics of the cover \(\grU\). 
 
 Let \(\{\chi_{B} \}_{ B\in \grU} \)
denote a partition of unity subbordinate to the cover
 \(\grU\) .  Note that these functions can be chosen so that
\( | d \chi_{B} | \leq c_{0} \sigma ^{-1} \).
 Define now a 1-form  \(\q\)  on \(B\) by
the rule \(\q|_{B} = \q_{B} + d(\sum_{B'} \chi_{B'} ( \k_{BB'}  - \c_{BB'} ))\).
 This defines a smooth 1-form on \(V\) with
\(d \q  =\uu\)  and with \(\|\q \|_{V} \leq c_{*} \sigma \|  \uu \|_{U}\).
 
\paragraph{\it Step 2:}  This lemma that is stated
and then proved in this step makes the first application of Lemma \ref{lem:D.5}.
  
\begin{lemma}\label{lem:D.6} 
There exists \(\kappa  >0\) with the following significance:
Fix \(k\in\{0, \ldots, \G\}\) and then \(T >1\) so as to define \(\op{Met}_T\) on \(Y_k\).
 Let \(\grg\) denote a \(\op{Met}_{T}\) metric on \(Y_k\) and let \(w_{\grg}\)
denote the corresponding harmonic 2-form whose de Rham cohomology class
is that of \(c_{1} (\det\, ( \bbS ))\).
 The 2-form \(w_{\grg}\) on the \(\ir \in [ \rho_{*} -\frac{1}{16}  \varepsilon ,
 \rho_{*} +\frac{1}{128}\varepsilon ]\) part of \(\mathcal{N}_\varepsilon\)
can be written as \(d \q\) with \(\q\) being a 1-form
whose \(L^{2}\)-norm on this part of \(\mathcal{N}_\varepsilon\)
is bounded by \(\kappa /T\) times that of \(w_{\grg}\).
\end{lemma}

\pf
The metric on the \(\ir \in[\rho_{*}  - \frac{1}{8} \varepsilon, \rho_{*} + \frac{1}{64} \varepsilon]\) part of  \(\mathcal{N}_\varepsilon\)
is the metric given by (\ref{eq:(D.7)}) with \(\rho_{T} =  
 \rho /T\) and with \(x_{3T} =
 x_3/T\).  The functions \(\textsc{k}\) and \(h\) are smooth around
\(\rho = 0\) with \(h(0\)) and \(\textsc{k}(0) = 1\).  It follows as a
consequence that the metric in the region of interest when written
using \(\rho_{ T }\) and \(x_{T}\) is uniformly close for \(T > c_{0}\) to the
Euclidean metric on the part of the radius  
 \((\rho_{*} +\frac{1}{64}\varepsilon)/T\) ball about the origin in
\(\bbR^{3}\) that lies outside the concentric ball of radius  
 \((\rho_{*}  -\frac{1}{8}\varepsilon)/T\).  Take this to be the region \(U\) for Lemma \ref{lem:D.5} and take
\(V\) to be the part of this same ball where the radius is between  
 \( (\rho_{*}  -\frac{1}{16}\varepsilon)/T\) and 
 \((\rho_{*} +\frac{1}{128} \varepsilon)/T\).  A cover \(\grU\) can be found as in the proof of Lemma \ref{lem:D.5}
with a \(T\)-independent bound on the number of sets, a \(T\)-independent
combinatorical structure to the intersections between them, and a
common radius for the balls, \(c_{0}\).   
This can be done because the \(T\)-dependence is just given by
scaling the coordinates.   Granted all of this, then the claim by the
lemma follows by appeal to Lemma \ref{lem:D.5}. \epf

  \paragraph{\it Step 3:}   This step supplies a part of what will be
\(\p_{X}\) on the \(s \in [-102, -101]\) part of \(X\) when
\(Y_{-} \) is a \(k\in\{0, \ldots, \G\}\) version of \(Y_k\), and on the \(s  \in
[100, 102]\) part of \(X\) when \(Y_{+}\) is a \(k\in\{0, \ldots, \G\}\) version of \(Y_k\).
 The constructions that follow use the embeddings from the second and
third bullets of (\ref{(A.9a,11)}) to view the \(s < 0\) and \(s >
0\) parts of \(X\) as \((-\infty, 0) \times Y_{-}\) and as \((0, \infty) \times Y_{+}\). 

 Let \(\chi_{\diamond 1}\)  denote the
function on \(\bbR\) given by the rule \(\chi(|s| - 101)\).  Denote
its derivative by \(\chi_{\diamond 1}'\). This function
is equal to 0 where \(|s | \geq 102\) and it is equal
to 1 where \(|s | \leq 101\).  Use \(\chi\) to
construct a smooth function on
 \(\mathcal{N}_\varepsilon\) that equals 0 where
\(|\ir - \rho_{*}|  > \frac{1}{128} \varepsilon \) and equals 1 where \(|\ir -\rho_{*}| < \frac{1}{256}
 \varepsilon\).  Construct this function of \(\r\) so that its derivative
is bounded by \(c_{0}\).  Use \(\sigma_{1}\) to
denote this new function of \(\r\).

If \(Y_{-}\) is a \(k\in\{0, \ldots, \G\}\) version of \(Y_k\), let \(\q_{1-}\) denote the \(w_{\grg} =
w_{-}\) version of \(\q\) that is given by Lemma \ref{lem:D.6}.  Define
\(\p_{ \mathcal{N} 1}\) where \(s \in [-102,-101]\) to be 
\begin{equation}\label{eq:(D.37)}
\p_{ \mathcal{N} 1} = - ds \wedge\chi_{\diamond 1}'\sigma_{1} \q_{1-}  + w_{-}  - \chi_{\diamond 1 }d(\sigma_{1} \q_{1-}).
\end{equation}
 
This is a closed form with de Rham cohomology class that of
\(c_{1}(\det(\bbS))\) and it equals
\(w_{-}\) where \(s \leq -102\).  Of particular note is the
fact that \(\p_{ \mathcal{N} 1} = 0\) on the
\(|\ir - \rho_{*}| < \frac{1}{256}\varepsilon \) part of  \(\mathcal{N}_\varepsilon\)
where \(s > -101\) and that it equals \(w_{-}\) on
the complement of the \(|\ir  -\rho_{*}| <\frac{1}{128}\varepsilon \) part of  \(\mathcal{N}_\varepsilon\).
 It follows from Lemma \ref{lem:D.6} that the \(L^{2}\) norm of
\(\p_{ \mathcal{N} 1}\) at any given \(s \in
[-102, -101]\) is bounded by \(c_{0}\) times that of \(w_{-}\).  

 If \(Y_{+}\) is a \(k\in\{0, \ldots, \G\}\) version of \(Y_k\), then very much the
same formula defines an \(s \in [101, 102]\) analog to
\(\p_{ \mathcal{N} 1}\).  The latter is
obtained by using Lemma \ref{lem:D.6} with \(w_{\grg} =
w_{+}\).  Lemma \ref{lem:D.6} supplies a 1-form \(\q_{1+}\)
with \(d\q_{1+} = w_{+}\) on the \(|\ir -\rho_{*}| < \frac{1}{128} \varepsilon \)
part of  \(\mathcal{N}_\varepsilon\).
 Use \(w_{+}\) and \(\q_{1+}\) in (\ref{eq:(D.37)}) in lieu of
\(w_{-}\) and \(\q_{-}\) to define \(\p_{ \mathcal{N} 1}\) where \(s \in [101,102]\).

\paragraph{\it Step 4:}  This step extends the definition of
\(\p_{ \mathcal{N} 1}\) to the \(s \in [-101,
-100]\) part of \(X\) when \(Y_{-}\) is a \(k\in\{0, \ldots, \G\}\)
version of \(Y_k\),
and to the \(s \in [100, 101]\) part of \(X\) when \(Y_{+} \) is a \(k\in\{0, \ldots, \G\}\)
version of \(Y_k\).  The embeddings from the second and third
bullets of (\ref{(A.9a,11)}) are again used to view the \(s < 0\) and \(s
> 0\) parts of \(X\) as \((-\infty, 0) \times Y_{-}\) and as \((0,\infty) \times Y_{+}\).  

  Thie extension of \(\p_{ \mathcal{N} 1}\) uses the function \(\chi_{\diamond 2}\)
on \(\bbR\) that is given by \(\chi(|s | - 100)\). The latter function is 0 where \(|s | \geq 101\)
and it is equal to 1 where \(|s | \leq 100\).  The
derivative of \(\chi_{\diamond 2}\) is denoted by \(\chi_{\diamond 2}'\).
 Reintroduce the closed 2-form \(\p_{0}\) from Step 1 of the
proof of Lemma \ref{Lemma D.1}.  By way of a reminder, this 2-form has compact
support on \(\mathcal{Y}_0\) ; and it has integral 2 over
each cross sectional 2-sphere in  \(\mathcal{H}_0\).  

 Suppose that \(Y_{-}\) is a \(k \in\{0, \ldots, \G\}\)
version of \(Y_k\).  The
extension of \(\p_{ \mathcal{N} 1}\) will equal
\(\p_{ \mathcal{N} 1}\) on the complement in
\(Y_{-}\) of the union of \(\mathcal{Y}_{0}\) and
the \(\ir \geq \rho_{*} +\frac{1}{512} \varepsilon \) part of  \(\mathcal{N}_\varepsilon\).
 Lemma \ref{lem:D.5} is used momentarily to obtain a 1-form to be denoted by
\(\q_{2-}\) with the following properties:  The 1-form
\(\q_{2-}\) has compact support on \(\mathcal{Y}_{0}\) and the \(\ir \geq
\rho_{*} + \frac{1}{512} \varepsilon  \) part of  \(\mathcal{N}_\varepsilon\),
its exterior derivative is \(w_{\grg} = w_{-} -\p_{0}  +
d(\sigma_{1}\q_{1-})\), and its \(L^{2}\)-norm is bounded by \(c_{0}\)
times that of \(w_{-}\).  Granted such a 1-form, the extension of
\(\p_{ \mathcal{N} 1}\) is given by
\begin{equation}\label{eq:(D.38)}
\p_{ \mathcal{N} 2} = - ds \wedge\chi_{\diamond 2}'\q_{2-}  +  w_{-} -
d(\sigma_{1} \q_{1-}) + \chi_{\diamond 2}\, d\q_{2- }.  
\end{equation}
 
This is a closed 2-form that equals
\(\p_{ \mathcal{N} 1}\) where \(s \leq -101\)
and for all \(s \in [-101, -100]\) on the complement of
\(\mathcal{Y}_{0}\) and the \(\ir \geq\rho_{*} +\frac{1}{512} \varepsilon \) part of  \(\mathcal{N}_\varepsilon\).
 This 2-form for \(s \geq -100\) is equal to \(\p_{0}\) on
\(\mathcal{Y}_{0}\) and the \(\ir \geq\rho_{*}\)  part of \(\mathcal{N}_\varepsilon\).
 
 The application of Lemma \ref{lem:D.5} takes \(U = V = S^{1}
\times S^{2}\).  The diffeomorphism \(\Phi_{T}\) in Part 6 of Section \ref{sec:Da)} is used to view
\(\p_{0} - (w_{-} -d(\sigma_{1} \q_{1-}))\) as a smooth 2-form on \(S^{1} \times S^{2}\), and viewed as
such, Lemma \ref{lem:D.5} is applied using this 2-form for \(w_{\grg}\).
 Lemma \ref{lem:D.5} then finds a 1-form, \(\q\), on \(S^{1} \times
S^{2}\) with \(d\q = \p_{0} -(w_{-} - d(\sigma_{1}\q_{1-}))\) and with \(L^{2}\)-norm bounded by
\(c_{0}\) times the \(L^{2}\)-norm of \(w_{-}\) on \(Y_{-}\).  The next two paragraphs
explains how to obtain \(\q_{2-}\) from \(\q\).

 View \(\p_{0} - (w_{-} -d(\sigma_{1} \q_{1-}))\) as a 2-form on
\(S^{1} \times S^{2}\) as done in the
preceding paragraph.  As explaind in Part 4 of Section \ref{sec:Da)}, the
coordinates \((\rho, \phi, x_{3})\) for \(\mathcal{N}_\varepsilon\) can be viewed where \(\r
\leq \rho_{*} + \frac{1}{16}\varepsilon \) as coordinates for a ball of this same radius in
  \(S^{1}\times S^2\).   The 2-form \(\p_{0} - (w_{-} -d(\sigma_{1} \q_{1-}))\) vanishes on the
concentric ball of radius 
 \((\rho_{*} + \frac{1}{256}  \varepsilon)/T\).  It follows as a consequence that \(\q\) can be written
as \(d\k\) with \(\k\) being a smooth function on this ball.  Since the
\(L^{2}\) norm of \(d\k\) on this ball is bounded by \(c_{0}\) times the \(L^{2}\)-norm of
\(w_{-}\) over \(Y_{-}\), it follows that \(\k\) can
be modified by adding a constant if nessecary so that its
\(L^{2}\)-norm over this ball is bounded by
\(c_{0}/T\) times the \(L^{2}\)-norm of \(w_{-}\) over \(Y_{-}\).  
 
Use \(\chi\) to construct a smooth function of the radial coordinate on
this ball with compact support that equals 1 on the concentric ball of
radius \((\rho_{*} +\frac{1}{512}\varepsilon)/T\) ball.  In particular, such a function can be
constructed so that the absolute value of its derivative is bounded by
\(c_{0} T\).  Let \(\sigma_{2}\) denote such a function and define \(\q_{*}\) to be \(\q - d(\sigma \k)\).
 This 1-form has the same properties as \(\q\) but it is zero on the
complement of the image of the embedding \(\Phi_{T}\) from
Part 6 of Section \ref{sec:Da)}.  The desired 1-form \(\q_{2-}\) is
\(\Phi_{T}^*\q_{*}\).
 
If \(Y_{+}\) is either a \(k\in\{0, \ldots, \G\}\)
version of \(Y_k\), then there is an analogous construction that defines
\(\p_{ \mathcal{N} 2}\) on the \(s \in [100,
101]\) part of \(X\).  The formula for the latter is given by replacing
\(w_{-}\), \(\q_{1-}\) and \(\q_{2+}\) by \(w_{+}\), \(\q_{1+}\) and a 1-form,
\(\q_{2+}\), that is defined by the rules given in the
preceding paragraph with \(w_{+}\) and \(\q_{1+}\) used in lieu of \(w_{-}\) and \(\q_{1-}\).

\paragraph{Part 10:}  Constructions in Part 11 and in the proof of
Proposition \ref{prop:A.10a} require a particular choice for the metric
\(\grm\) on certain 
parts of \(X\).   The constraint given momentarily holds on the \(s \in
[-100, -96]\) part of \(X\) when \(Y_{-}\) is a \(k\in\{0, \ldots, \G\}\)
version of \(Y_k\), and it holds on the \(s \in [96, 100]\) part of \(X\) when
\(Y_{+}\) is a \(k\in\{0, \ldots, \G\}\)
version of \(Y_k\).  

 The statement of the constraint uses the embeddings from the second
and third bullets of (\ref{(A.9a,11)}) to view the \(s < 0\) and \(s
> 0\) part of \(X\) as \((-\infty, 0] \times Y_{-}\) and as \((0, \infty) \times 
Y_{+}\).  Viewed this way, the constraint on the metric \(\grm\)
involves only the \(\ir \in [\rho_{*} - \frac{15}{16}  \varepsilon, \rho_{*} )\) parts of  \([-100,
-96] \times  \mathcal{N}_\varepsilon\) and \([96,100] \times  \mathcal{N}_\varepsilon\).  To
define \(\grm\) on these parts of \(X\), construct a smooth, non-decreasing
function on \(\bbR\) to be denoted by \(T_{\diamond}\):  This function equals \(T\) where
\(|s | \geq 99\) and it equals \(1\) where \(|s| \leq 98\).   The ubiquitous function \(\chi\) can be
used to define this function \(T_{\diamond}\).  Reintroduce the metric \(\grg_{T}\) on
 \(\mathcal{N}_\varepsilon\) that is defined in Part 5 of Section
 \ref{sec:Da)}.  The assignment \(s \mapsto\grg_{T_{\diamond}}\)
 defines a 1-parameter family of metrics on
\(\mathcal{N}_\varepsilon \) with parameter space
either \([-100, -96]\) or \([96, 100]\) as the case may be.  The \(|s
| = 100\) end member of this family is \(\grg_{T}\) and
the \(|s | = 96\) member is the metric in (\ref{eq:(D.6)}).  
 
Use \(\chi\) to construct a smooth function of the coordinate \(\ir\) on
 \(\mathcal{N}_\varepsilon\) that is equal to 1 where \(\ir
< \rho_{*}  -\frac{1}{1024}  \varepsilon \) and equal to 0 where
\(\ir >\rho_{*}  - \frac{1}{2048}\varepsilon\).  Use \(\sigma_{*}\)
to denote this function. 

The metric \(\grm\) is constrained by requiring that its pull-back to \([-100,
-96] \times  \mathcal{N}_\varepsilon \) via the
embedding from the second bullet of (\ref{(A.9a,11)}) or to \([96, 100] \times
 \mathcal{N}_\varepsilon \) via the embedding from the
third bullet of (\ref{(A.9a,11)}) to be the metric
\begin{equation}\label{eq:(D.39)}
ds^{2} + \sigma_{*}\, \grg_{T_\diamond}+ (1  - \sigma_{*}) \, \grg_{T}.
\end{equation}

Note in particular that this metric smoothly extends
\(ds^{2} + \grg_{T}\) near \(|s |= 100\) and it smoothly extends \(ds^{2} +
\grg_{T}\) from the \(\ir \leq \rho_{*}- \frac{15}{16} \varepsilon \) part of  \(\mathcal{N}_\varepsilon\)
for all \(s\) in the relevant interval.

An important observation is given momentarily about the   
 versions of the \(L^{2}\)-norm of \(w_{-} -
d(\sigma_{1} \q_{1-})\) on the \(\ir \leq\rho_{*} - \frac{1}{512}\varepsilon \) part of  \(\mathcal{N}_\varepsilon\).
 Keep in mind in what follows that this 2-form is zero on the \(\ir
> \rho_{*}  - \frac{1}{256}\varepsilon \) part of  \(\mathcal{N}_\varepsilon\).
 Given \(s  \in [-100, -96]\), the notation uses
 \(\| w_{-} -d(\sigma_{1} \q_{1-})\|_{s}\) to denote   
  version of the \(L^{2}\)-norm of \(w_{-} -
d(\sigma_{1} \q_{1-})\) on the \(\ir <\rho_{*}  - \frac{1}{512}\varepsilon \) part of  \(\mathcal{N}_\varepsilon\).
 There is the analogous definition for \(\|w_{+}  - d(\sigma_{1}\q_{1+}) \|_{s}\) when \(s
\in [96, 100]\). Here is the key observation:  
\BTitem\label{eq:(D.40)}
\item Each \(s \in [-100, -96]\) version of
\(\| w_{-} -d(\sigma_{1} \q_{1-})\|_{s}\)  is bounded by
\(c_{0}\) times the \(L^{2}\)-norm of \(w_{-}\) on \(Y_{-}\).
\item Each \(s   \in [96, 100]\) version of
\(\| w_{+} - d(\sigma_{1}\q_{1+}) \|_{s}\)  is bounded by \(c_{0}\)  times the \(L^{2}\)-norm of
\(w_{+}\)  on \(Y_{+}\).
\ETitem

To see about (\ref{eq:(D.40)}), write any \(s \in[-100, -96]\) or \(s \in
[96, 100]\) version of \(\grg_{T_\diamond}\) at any given point in the \(\ir < \rho_{*}- \frac{1}{512}\varepsilon \) part of
 \(\mathcal{N}_\varepsilon\) as
\begin{equation}\label{eq:(D.41)}
\grg_{T_\diamond}   = \lambda_{1}\hat{e}^{1} \otimes \hat{e}^{1} +
\lambda_{2}\hat{e}^{2}\otimes \hat{e}^{2} + \lambda_{3} \hat{e}^{3} \otimes \hat{e}^{3} 
\end{equation}
with \(\lambda_{1}\), \(\lambda_{2}\) and \(\lambda_{3}\) being positive numbers and with
\(\{\hat{e}^{k}\}_{k=1,2,3}\) being a
\(\grg_{T}\)-orthonormal frame.  It follows from (\ref{eq:(D.7)})-(\ref{eq:(D.9)})
that each \(\lambda_{k}\) can be written as \((T/T_\diamond   )^{2 }\e_{k}\) where the numbers
\(\e_{1}\), \(\e_{2}\) and \(\e_{3}\) are such
that \(c_{0}^{-1} \leq \e_{1}, \e_{2} \leq c_{0}\) and
\(c_{0}^{-1} \leq \e_{3}\leq c_{0}(T/T_\diamond     )^{2}\).  It follows from this that the volume form of  
the metric  is less than \(c_{0} (T/T_\diamond      )^{4}\) times that of \(\grg_{T}\).  It also
follows from this that the square of the \(\grg_{T_\diamond}\)-norm of \(w_{-} - d(\sigma_{1}\q_{1-})\) is less than \(c_{0}(T/T_\diamond      )^{4}\) times the square of its \(\grg_{T}\) norm.
 These last two observations imply that the integrand whose integral
gives \(\| w_{-} -d(\sigma_{1} \q_{1-})\|_{s}^{2}\) is no greater
than \(c_{0}\) times the integrand whose integral computes the
square the \(\grg_{T}\) version of the \(L^{2}\)-norm
of \(w_{-} - d(\sigma_{1}\q_{1-})\).  This last fact implies directly the first bullet
of (\ref{eq:(D.40)}).  But for replacing \(-\) subscripts with \(+\) subscripts, the same
argument proves the second bullet of (\ref{eq:(D.40)}).

\paragraph{Part 11:}  This part of the subsection completes the proof
of Proposition \ref{prop:A.8}.   According to Part 8, it is sufficient to find
the closed 2-form \(\p_{X}\) with certain special properties.
 This is done given two more constraints on \(\grm\). The
first contraint affects \(\grm\) only on the
\(|s |  \in  [96, 100]\) part of \(X\).  The
statement of this uses the embeddings from the second and third bullets
of (\ref{(A.9a,11)}) to view the \(s < 0\) and \(s > 0\) parts of \(X\)
as \((-\infty, 0] \times Y_{-}\) and as \((0,\infty) \times  Y_{+}\):
\begin{equation}\label{eq:(D.42)}\begin{split}
 & \text{The metric \(\grm\)  on \([-100, -96] \times\mathcal{Y}_{M}\)
is the product metric \(d s ^{2} +\grg_{-}\)}\\
& \text{when \(Y_{-}= Y_{0}\); and  when \(Y_{+} = Y_{0}\), the metric \(\grm\) on \([96, 100]\times \mathcal{Y}_{M}\)}\\
&\text{ is the product metric
\(d s ^{2}+\grg_{+}\).}
\end{split}
\end{equation}
 
To state the second constraint, re-introduce from
Part 7 of Section \ref{sec:Da)} the set \(\Theta\) and the associated
collection \(\{\mathcal{T}_{\gamma}\}_{(\gamma, Z_\gamma)\in \Theta}\) of subsets of \(M_{\delta}\).
 The following observation views them as subsets of
\(Y_{0}\)  and \(M\):

There exists an embedding of \(\bbR \times 
( \bigcup_{(\gamma, Z_\gamma)\in \Theta}\mathcal{T}_{\gamma} )\) into  \(X\)
with the following two properties:
\BTitem\label{eq:(D.43)}
\item The function \(s\) on \(X\) pulls back via the embedding to its namesake
on the \(\bbR  \)-factor of \(\bbR \times 
( \bigcup_{(\gamma, Z_\gamma)\in \Theta}\mathcal{T}_{\gamma} )\).
\item  The composition of this embedding of the
\(|s |  > 1\) part of \(\bbR \times ( \bigcup_{(\gamma, Z_\gamma)\in
  \Theta}\mathcal{T}_{\gamma} )\) with the inverse of the embeddings
from the second and third bullets of (\ref{(A.9a,11)}) is the tautological inclusion
map.
\ETitem
 
 The existence of such an embedding is implied by what
is said in the first paragraph of this section about the ascending and
descending manifolds from the critical point of \(s\).  The second constraint uses
 \(\grm_{-}\) and \(\grm_{+}\) to denote the metrics \(d s ^{2}+\grg_{-} \) and
\(d s ^{2}+\grg_{+} \) on the product  \(\bbR \times ( \bigcup_{(\gamma, Z_\gamma)\in \Theta}\mathcal{T}_{\gamma} )\).

\begin{equation}\label{eq:(D.44)}
 \begin{split}
&\text{There exists a \(T\)-independent constant,
\(c_{*} > 1\), with  the following }\\ 
& \text{significance: The pull-back of \(\grm\) via the embedding in
(\ref{eq:(D.43)}) obeys}\\
&  \text{\(c_{*} ^{-1} \grm_{-} \leq  \grm  \leq c_{*} \grm_{-}\)
and \(c_{*} ^{-1}\grm_{+} \leq  \grm  \leq c_{*} \grm_{+}\).}
\end{split}
\end{equation}

Granted these constraints, the three steps that follow construct
 \(\p_{X}\) when \(Y_{-}  =Y_{0}\).  The construction when \(Y_{+}  =Y_{0}\)  is not given
as it has the identical description but for changes of \(-\) to \(+\) in
various places.

\paragraph{\it Step 1:} Define \(\p_{X}\) on the \(s \in  [-102, -101]\) part of \(X\) to
be \(\p_{ \mathcal{N} 1}\) and define \(\p_{X}\) on the \(s \in [-101, -100]\) part of \(X\) to
be \(\p_{ \mathcal{N} 2}\).
 The rest of this step extends the definition of \(\p_{X}\) 
to the  \(s \in  [-100, -98]\) part of \(X\).
 To this end, use the embedding from the second bullet of (\ref{(A.9a,11)}) to
view this part of \(X\) as \([-100, -98]\times Y_{0}\). 
 
The 2-form \(\p_{ \mathcal{N} 2}\) near \(s = -100\)
is the \(s\)-independent 2-form on \(Y_{0}\) given by \(\p_{0}\) on \(\mathcal{Y}_{0}\)  and
\(w_{-} - d(\sigma_{1}q_{1-})\) on the rest of \(Y_{0}\).  This
understood, \(\p_{X}\) is extended to the \(s \in [-100, -98]\)
part of \(X\) as this \(s\)-independent 2-form on \(Y_{0}\).  
 
Write the metric \(\grm\) appearing on \([-100, -98] \times Y_{0}\)
as \(ds^{2} + \grg\) with \(\grg\) denoting an \(s\)-dependent metric on
\(Y_{0}\).  The constraint in (\ref{eq:(D.42)}) asserts that \(\grg =
\grg_{-}\) on \(\mathcal{Y}_{M}\).  Meanwhile, \(\grg\) is Part 9's metric   
 on the \(\ir < \rho_{*} -\frac{1}{512}\varepsilon \) part of  \(\mathcal{N}_\varepsilon\).
  It therefore follows from (\ref{eq:(D.40)}) that the
 \(L^{2}\)-norm of \(\p_{X}\) on \(Y\) as defined by any \(s \in
 [-100, -98]\) version of \(\grg\) is bounded by \(c_{0}\).   
 
\paragraph{\it Step 2:}  This step extends the definition of \(\p_{X}\) to the
\(s \in [-98, -96]\) part of \(X\).  To do this, view the \(s \in [-98,
-96]\) part of \(X\) as \([-98, -96] \times~Y_{0}\) as in Step
1.  Keep in mind for what follows that the metric \(\grm \) here has the form
\(ds^{2} + \grg_{M}\) with \(\grg_{M}\) being an \(s\)-independent metric on \(Y_{0}\).  Note in
particular that \(\grg_{M} = \grg_{-}\) on \(\mathcal{Y}_{M}\) and it is the metric that is
depicted in (\ref{eq:(D.6)}) on the \(\ir \leq \rho_{*} - \frac{1}{1024} \varepsilon \) part of  \(\mathcal{N}_\varepsilon\).

  Lemma \ref{lem:D.5} is invoked momentarily to construct a 1-form on the union of
\(\mathcal{Y}_{M}\) and the \(\ir <\rho_{*}\) part of
 \(\mathcal{N}_\varepsilon\), with the three properties
listed momentarily.  The list of properties denotes the 1-form by
\(\q_{3-}\) and it reintroduces the 2-form \(\p\) from Part 7 of
Section \ref{sec:Da)}.  Here are the three properties:  The 1-form
\(\q_{3-}\) obeys \(d\q_{3-} =\p -w_{-} + d(\sigma_{1}\q_{1-})\), it vanishes on the \(r \geq
\rho_{*}  -\frac{1}{512} \varepsilon \) part of  \(\mathcal{N}_\varepsilon\),
and its \(L^{2}\)-norm as defined by the \(\grg_{M}\)
is bounded by \(c_{0}\) times the \(L^{2}\)-norm of \(w_{-}\) on \(Y\).

 Let \(\chi_{\diamond 3}\) denote the function on \(\bbR\) given by
\(\chi(|s |  - 97)\).  The function \(\chi_{\diamond 3}\) equals 0 where \(|s | \geq 98\) and it equals 1
where \(|s | \leq 97\).  Introduce \(\chi_{\diamond 3}'\) to denote its
derivative.  The 2-form \(\p_{X}\) on the \(s \in [-98, -96]\)
part of \(X\) is \(\p_{0}\) on \(\mathcal{Y}_0\) and
it is given on the rest of \(Y_{0}\) by 
\begin{equation}\label{eq:(D.45)}
ds \wedge \chi_{\diamond 3}'\q_{3-} + w_{-} -d(\sigma_{1} \q_{1-})
+\chi_{\diamond 3} \, d\q_{3-} .
\end{equation}
 
Of particular note is that the \(\grm\) version of the \(L^{2}\)-norm of the 2-form \(\p_{X}\) on \([-98, -96] \times Y_{0}\) is bounded by \(c_{0}\).  What follows is
a key point to keep in mind for Step 3:  The 2-form \(\p_{X}\)
on \([-97, -96] \times Y_{0}\) is the 2-form \(\p_{0} + \p\) from \(Y_{0}\). 
 
Lemma \ref{lem:D.5} is invoked using for the set \(U\) the union of 
\(\mathcal{Y}_{M}\) and the \(\ir <\rho_{*} - \frac{1}{1024}\varepsilon \) part of  \(\mathcal{N}_\varepsilon\).
 Lemma \ref{lem:D.5}'s set \(V\) is the union of
\(\mathcal{Y}_{M}\) and the \(\ir <\rho_{*}  - \frac{1}{512} \varepsilon \) part of  \(\mathcal{N}_\varepsilon\).
 The 2-form \(w_{\grg}\) is \(\p - w_{-} +d(\sigma_{1} \q_{1-})\).  Note that this
2-form is zero on the \(\ir > \rho_{*} -\frac{1}{256} \varepsilon \)
part of \(U\).  The metric used for the lemma is the metric \(\grg_{M}\).  It follows from (\ref{eq:(D.40)}) and (\ref{eq:(D.42)}) that the
\(L^{2}\)-norm of  \(\p - w_{-} +d(\sigma_{1} \q_{1-})\) as
defined by \(\grg_{M}\) is bounded by \(c_{0}\).  As
neither \(U\), \(V\) nor \(\grg_{M}\) depend on \(T\), Lemma \ref{lem:D.5} finds a
1-form \(\q\) on \(U\) with \(d\q = \p - w_{-} +d(\sigma_{1} \q_{1-})\) whose
\(L^{2}\)-norm on \(V\) is bounded by \(c_{0}\).  To
obtain \(\q_{3-}\) from \(\q\), note that \(\q\) on the \(\ir >
\rho_{*} -\frac{1}{256} \varepsilon \) part of  \(\mathcal{N}_\varepsilon\) is
given by \(d\k\) with \(\k\) denoting a smooth function.   Changing \(\k\) by a
constant if necessary supplies a version whose \(L^{2}\)-norm is
bounded by \(c_{0}\) times that of \(d\k\); thus by \(c_{0}\).  Take
such a version.  Meanwhile, use \(\chi\) to construct a smooth function of \(\r\) on
 \(\mathcal{N}_\varepsilon\) that equals 0 where \(\r
\geq \rho_{*}  -\frac{1}{512}  \varepsilon\), equals 1 where \(r \leq\rho_{*} -\frac{1}{256}\varepsilon   \)
  and whose derivative has norm bounded by \(c_{0}\).  Denote
this function by \(\sigma_{3}\) and set \(\q_{3-} = \q - d(\sigma_{3}\k)\).

\paragraph{\it Step 3:}  This step extends the definition of \(\p_{X}\) to the
\(s \in [-96, 102]\) part of \(X\).  To this end, consider first the
definition of \(\p_{X}\) on the \(s \in [-96, 100]\) part of \(X\).
 As \(\p_{0}\) is supported in the image of the embedding from
(\ref{eq:(D.10)}) and as the 2-form \(\p\) is supported in the image of the embedding
from (\ref{eq:(D.43)}), these embeddings can be used to view \(\p_{0} +
\p\) as a 2-form on the \(s \in [-96, 100]\) part of \(X\).  View them in
this light and define \(\p_{X}\) on this same part of \(X\) to be
\(\p_{0} + \p\).  The constraint in (\ref{eq:(D.44)}) has the following
implication:  The \(L^{2}\)-norm of \(\p_X\) on
the \(s \in [-96, 100]\) part of \(X\) is bounded by \(c_{0}\). 

 The definition of \(\p_{X}\) on the \(s \in [100, 102]\) part of
\(X\) views this part of \(X\) via the embedding from the third bullet of (\ref{(A.9a,11)})
as \([100, 102] \times (M\sqcup (S^{1} \times S^{2}))\).  The 2-form \(\p_{0}\) on
\(S^{1} \times S^{2}\) can be written as \(w + d\q_{0}\) with \(\q_{0}\) being a smooth 1-form.
 Likewise, the 2-form \(\p\) on \(M\) can be written as
\(w_{+}|_{M} + d\q_{M}\)
with \(\q_{M}\) denoting a smooth 1-form.  Set
\(\q_{4+} = \q_{0} + \q_{M}\). Let \(\chi_{\diamond 4}\) denote the
function on \(\bbR\) given by \(\chi(s   - 100)\).  This function
\(\chi_{\diamond 4}\) is equal to 1
where \(s < 100\) and it is equal to 0 where \(s >101\).  Use \(\chi_{\diamond 4}'\) to
denote its derivative.

 Define \(\p_{X}\) on the  \(s \in[100, 102]\) part of \(X\) to be
the 2-form
\begin{equation}\label{eq:(D.46)}
ds \wedge \chi_{\diamond 4}'\q_{4+} + \p_{0} + \p -\chi_{\diamond 4}d\q_{4+}
\end{equation}
 This form is closed, and it extends \(\p_{X}\) as a 2-form that
equals \(w_{+}\) where \(s > 10\)1.  Of particular
note is that the \(L^{2}\)-norm of \(\p_{X}\) on the
\(s \in [100, 102]\) part of \(X\) is bounded by \(c_{0}\). \epf

\subsection{Proof of Proposition \ref{prop:A.10a}}\label{sec:De)}

The proof of this proposition has two parts.  Of the two possible
cases, only that where 
\(Y_{-} = Y_k\) and \(Y_{+} = Y_{k-1}\sqcup(S^1\times S^2)\)
is discussed as the case when the roles are switched is proved
with the same argument but for changing the direction of various
inequalities and signs that involve \(s\).  

 Part 1 of what follows proves the first bullet of Proposition \ref{prop:A.10a}.
 The subsequent parts of this subsection address the assertion in the
second bullet and in doing so, they define implicitly the required
subset \(\op{Met}(Y_k)\).
To make the definition only slightly
less implicit, return momentarily to what is said about
\(\Met\) just prior to Part 1 of Section \ref{sec:Dd)}.
 By way of a reminder, each metric in \(\Met\)
is determined in part by a metric from the \(Y_{0}\) version
of Section \ref{sec:Da)}'s set \(\op{Met}^{ \mathcal{N}}\)
and a sufficiently large choice for a number denoted by \(T\).  A lower
bound on \(T\) is determined by certain properties of the chosen
\(\op{Met}^{\mathcal{N}}\)   metric.  This said, a metric from
\(\Met\) is in Proposition \ref{prop:A.10a}'s subset \(\op{Met}(Y_k)\)
if and
only if the chosen value for \(T\) is larger than a new lower bound.  This new lower bound is determined in
part by the same properties of the chosen \(\op{Met}^{\mathcal{N}}\) 
metric that determine the \(\op{Met}(Y_0)\) lower bound. The
chosen metrics on the \(S^{1} \times  S ^{2}\)  components also
determine in part the lower bound for \(T\). By the way, no generality is
lost by taking the metrics on these components to be the product of the
standard Euclidean \(S^{1}\)  and the standard round metric on \(S^{2}\).
 In any event, this new lower bound is determined implicitly by the
 constructions in Parts 2-10.

\paragraph{Part 1:} This part discusses the first bullet of the
proposition. The notation used below is that used to describe \(Y\) and
its geometry in \cite{KLT1}-\cite{KLT4}. In particular, the manifold
\(Y\) and its 2-form \(w\) are described in
Section II.1. A summary of the salient features can be found in
Section IV.1a. The notation used below is the same as that used in
Section II.1 and Section IV.1a. 

To set the stage, label the \(\G\) pairs in the set \(\Lambda\) as
\(\{\grp_1, \ldots, \grp_{\G}\}\). A \(k\in \{1, \ldots, \G\}\)
version of the manifold \(Y_k\) is obtained from 
\(Y_0\) by attaching \(k\) 1-handles,
these being the handles from from the set
\(\{\mathcal{H}_\grp\}_{\grp\in\{\grp_1, \ldots, \grp_k\}}\). Thus,
\(Y_k\) is obtained from \(Y_{k-1}\) by attaching just the handle \(\mathcal{H}_{\grp_k}\).
By way of a short review, \(Y\) is
obtained from \(Y_{0}\) by a surgery that attaches \(\textsc{g}\) 1-handles to
\(Y_{0}-\mathcal{H}_{0}\).
The attaching region of each handle are disjoint coordinate balls
centered around a pair of points in \(Y_{0}-\mathcal{H}_{0}\).
The set of such pairs is denoted by \(\Lambda\). The 1-handle that
corresponds to a given pair \(\grp \in \Lambda \) is denoted by
\(\mathcal{H}_{\grp}\). The geometry of \(Y_k\) near
\(\mathcal{H}_{\grp_k}\) is as follows: The handle \(\mathcal{H}_{\grp_k}\) is diffeomorphic to
\([-R -7\ln\delta_{*}, R + 7\ln\delta_{*}]\times S^{2}\) given by the 
preferred coordinates
\((u,(\theta, \phi))\) with \(u\) denoting the Euclidean coordinate for
interval factor and with \((\theta, \phi)\) denoting spherical
coordinates on the constant \(u\) cross-sectional spheres of
\(\mathcal{H}_{\grp_k}\). The handle is attached to \(Y_{k-1}\) using the
identifications given in (\ref{eq:(A.1)}) with it understood that
\((\ir_{+}, (\theta_{+}, \phi_{+}))\) and \((\ir_{-},(\theta_{-}, \phi_{-}))\) are
certain preferred spherical coordinates for respective balls about the
two points that comprise the pair \(\grp_k\). 


The definition of \(X\) requires choosing 
a properly embedded arc in the \(\mathcal{Y}_{M}\) part of \(Y_{k-1}\). The arc has one end
point at one of the points in \(\grp_k\) and the other end point at the
other. This arc intersects a neighorhood of the boundary of the radius
\(7\delta_{*}\) coordinate ball centered at the points from \(\grp_k\) as a ray from the
origin when viewed using the coordinate system that is specified in
Section II.1.a. Part 7 of Section 9.1 introduces a finite set of pairs
\(\Theta\) in \(M_\delta\) with one partner in each pair being an embedded loop in
\(M_\delta\). Part 7 of Section \ref{sec:Da)} associates each such loop a small radius
tubular neighorhood, this being \(\mathcal{T}_{\gamma}\) when \(\gamma\) is the loop in question. The
arc must be chosen so as to lie in the complement of the closure of
all such tubular neighborhoods. The arc can and should be chosen to
intersect that \(\ff = \frac{3}{2}\) Heegaard surface in \(M_\delta\) transversely in a single
point. Denote this arc by \(\lambda_{\grp_k}\).

Let \(S_{\grp_k}\subset \mathcal{Y}_M\) denote an embedded 2-sphere boundary of
neighborhood of the arc \(\lambda_{\grp_k}\) with each point having distance between
\(2\delta\) and \(4\delta\) from the arc. The neighborhood in question
and \(S=S_{\grp_k}\) should be disjoint from the closures of the tubular
neighborhoods of the loops from \(\Theta\). The sphere \(S\) appears
in \(Y_k\) as an embedded 2-sphere that separates \(Y_k\) into two
components. One of these contains \(\mathcal{H}_{\grp_k}\) and is
diffeomorphic to the complement in \(S^1\times S^2\) of an embedded ball.

The following is a consequence of what is said
above about the descending and ascending submanifolds from the critical
points of \(s\): The pseudogradient
vector field that defines the embeddings from the
second and third bullets of (\ref{(A.9a,11)}) can be chosen so that
(\ref{eq:(A.9b)}) 
are obeyed and likewise (\ref{eq:(A.10)}) and the
conditions in (\ref{eq:(D.10)}) and (\ref{eq:(D.43)}). These
properties are assumed in what follows. The condition for the first
Chern class 
is satisfied if and only it 
has zero pairing with the cross-sectional 2-spheres in each \(\grp\in
\{\grp_1, \ldots, \grp_{k-1}\}\) version of the \(Y_{k-1}\) version of
\(\mathcal{H}_{\grp}\) and annihilates the generator of
\(H_2(S^1\times S^2;\bbZ)\).


\paragraph{Part 2:} Proposition  \ref{prop:A.10a} requires as input a
metric from a certain subset of a set of metrics on \(Y_{k-1}\) that
is denoted by \(\op{Met}(Y_{k-1})\) and a metric from a set of metrics
on \(Y_k\), this denoted by \(\op{Met}(Y_k)\). 
These subsets of metrics are in the respective
\(Y_{k-1}\) and \(Y_k\) versions of \(\Met\). They are defined roughly as
follows: Let \(Y_*\) for the moment denote either \(Y_{k-1}\) or
\(Y_k\). Each metric in the \(Y_*\) version of \(\Met\) is determined in
part by a metric from the corresponding version of \(\op{Met}_N\), this
defined in Section \ref{sec:Da)}. The second input for the definition is a large
choice for the parameter \(T\). A metric in \(\Met\) of this sort is denoted in
Section \ref{sec:Db)} by \(\grg_{3T}\). A \(Y_*\) metric \(\grg_{3T}\) is in \(\op{Met}(Y_*)\) if
\(T\) is greater than a certain lower bound that is determined
implicitly by the chosen \(\op{Met}_N\) metric. As in the case of
Proposition \ref{prop:A.10a}'s implicit definition of \(\op{Met}(Y_0)\),
this lower bound is determined implicitly by the requirements of
subsequent constructions. In any event, it is determined by certain
curvature norms, injectivity radius lower bounds and volume.

The construction of a suitable metric on \(X\) starts by choosing metrics \(\grg_{1-}\) and \(\grg_{1+}\) from the respective \(Y_-\) and \(Y_+\) versions of \(\op{Met}_N\). This done, use what is said in Parts 1-10 of Section \ref{sec:Dd)} to define a metric \(\grm_{T*}\) and self-dual 2-form \(\omega_{T*}\) on \(X\). It then follows from what is said in Part 8 and at the start of Part 9 of Section \ref{sec:Dd)} that the pair \(\grm_{T*}\) and \(\omega_{T*}\) satisfy the requirements of Proposition \ref{prop:A.10a} if there exists a suitable closed 2-form \(\p_X\) on
\(X\) with the following properties: The de Rham cohomology class of
\(\p_X\) is that of \(c_1(\det(\bbS))\). In addition, \(\p_X\) must
equal \(w_-\) where \(s < -102\) and \(w_+\) where \(s > 102\) with
\(w_-\) and \(w_+\) being the respective \(\grg_-\) and \(\grg_+\)
harmonic 2 forms on \(Y_-\) and \(Y_+\) with de Rham cohomology class
that of \(c_1(\det(\bbS))\).

The construction of \(\p_X\) is this case differs in only one respect from the construction
described in Parts 9-11 of Section \ref{sec:Dd)}, this involving Step
3 in Part 11 of Section \ref{sec:Dd)}. To say more about this
difference, require as in Part 11 of Section \ref{sec:Dd)} that the
metric \(\grm\) obey (\ref{eq:(D.42)}). Require in addition that
(\ref{eq:(D.43)}) is obeyed; as noted in Part 1 above, such a
requirement can be met. With (\ref{eq:(D.43)}) understood, the metric
\(\grm\) is chosen so as to obey the constraints in
(\ref{eq:(D.44)}). Proceed with the constructions in Steps 1 and 2 of
Part 11 in Section \ref{sec:Dd)}. Step 3 in Part 11 of Section
\ref{sec:Dd)} is replaced with the following Step \(3'\):

\paragraph{\it Step \(3'\):} Define \(\p_X\) on the \(s\in [-96, 96]\)
part of \(X\) by viewing \(\p_0 +\p\) as a 2-form on this
part of \(X\) via the embeddings in (\ref{eq:(D.10)}) and (\ref{eq:(D.43)}). The constraint in (\ref{eq:(D.44)}) implies that
such a definition yields a version of \(\p\) with \(L^2\)-norm bounded
by \(c_0\) on the \(s\in [-96, 96]\) part of \(X\). Extend \(\p_X\) to the \([96, 102]\) part of \(X\) by copying almost
verbatim what is done in Steps 1 and 2 with the direction of \(s\)
reversed and with the metric \(\grg_+\) in (\ref{eq:(D.42)}) used in
lieu of \(\grg_-\).

\subsection{Proof of Proposition \ref{prop:ech-compute}}\label{pf:prop1.5}
This subsection provides a proof of Proposition
\ref{prop:ech-compute} and thus completes the proof of Theorem \ref{cor:main}.


Explicit formulae for the differentials and \(\mathbf{A}_\dag(Y)\)
actions on the chain complex used to define \(ech^\circ\) are
given in Theorem 1.1 of \cite{KLT3}. These formulae were also written in terms of a
factorization of  \(\mathbf{A}_\dag(Y)\) into a tensor product \(\mathbf{A}_\dag(M)\otimes
H_{-*}(S^1)\otimes H_{-*}(S^1)^{\otimes\G}\), which is however {\em
different} from that in (\ref{eq:A-inclusion}), the
factorization used in the statement of the proposition.
The difference originates from a different choice of splitting for 
\(H_1(Y;\bbZ)/\op{Tors}\) from 
that in (\ref{eq:H_Y-split1}). 

In \cite{KLT3}, an
``\(M\)-adapted 1-cycle basis'' is assigned to
\(H_1(Y;\bbZ)/\op{Tors}\), whose basis elements are represented by
``\(M\)-adapted 1-cycles'' in \(Y\). Each ``\(M\)-adapted 1-cycle''
is of one of the following three types: \begin{itemize}
\item \(\hat{i}^{(z)}\) for every \(z\in \yen
-z_0\); \item \(\gamma^{(z_0)}\); and \item \(\hat{i}_\grp\) for each \(\grp\in
\Lambda\). \end{itemize}
Decompose \(H_1(Y;\bbZ)/\op{Tors}\) accordingly into 
\begin{equation}\label{def:H_split2}
H_1(Y;\bbZ)/\op{Tors}\simeq H_1(M;\bbZ)/\op{Tors}\oplus H_1(S^1\times S^2;\bbZ)\oplus\bigoplus _{\grp\in \Lambda}H_1((S^1\times S^2)_\grp;\bbZ),
\end{equation}
with the first, second, and third summand generated by the ordered
sets \(\{[\hat{i}^{(z)}]\}_{z\in \yen
-z_0}\); \(\{[\gamma^{(z_0)}]\}\); and \(\{[\hat{i}_\grp]\}_{\grp\in
\Lambda}\) respectively. 

On the other hand, in Section \ref{sec:1.1} we split
\(H_1(Y;\bbZ)/\op{Tors}\) differently using a connected-sum
decomposition of \(Y\). Namely, by 
combining (\ref{eq:H_Y-split0})
and (\ref{eq:H_Y-split1}) we get another splitting
\begin{equation}\label{def:H_split3}\begin{split}
& H_1(Y;\bbZ)/\op{Tors}\\
& \simeq H_1(\ul{M};\bbZ)/\op{Tors}\bigoplus _{\grp\in
\Lambda}H_1((S^1\times S^2)_\grp;\bbZ)\\
&\simeq H_1(M;\bbZ)/\op{Tors}\oplus H_1(S^1\times S^2;\bbZ)\oplus\bigoplus _{\grp\in \Lambda}H_1((S^1\times S^2)_\grp;\bbZ).\end{split}
\end{equation} 
Note that the preceding 
splitting (\ref{def:H_split3}) depends on the relative homlogy class
of the chosen arcs \(\gamma_z\) and \(\lambda_\grp\)'s. The summands from this 
splitting are generated by elements in
\(H_1(Y;\bbZ)/\op{Tors}\) represented by the following sets of
1-cycles in \(Y\). 
\begin{itemize}
\item For the first summand \(H_1(M;\bbZ)/\op{Tors}\): \(b_1(M)\)
1-cycles from the \(M\)-summand of the connected sum decomosition \(Y\simeq
M\#_{\G +1} (S^1\times S^2)\), so that as 1-cycles in \(M\) they avoid all the arcs \(\gamma_z\)
and \(\lambda_\grp\) and their homology classes together
form an (arbitrary) basis for
\(H_1(M;\bbZ)/\op{Tors}\). For example, the set of 1-cycles \(\{[\hat{i}^{(z)}]\}_{z\in \yen
-z_0}\) is a possible choice. 
\item For the second summand \(H_1(S^1\times
S^2;\bbZ)\): The cycle coming from the
1-cycle \(\ul{\gamma }\) in the \(\ul{M}\)-summand of the connected
sum \(Y\simeq \ul{M}\#_{\G} (S^1\times S^2)\). (This cycle in \(Y\) was
called  \(\gamma^{(z_0)}\) in \cite{KLT3}). 
\item For each \(H_1((S^1\times
S^2)_\grp;\bbZ)\)-summand in \(\bigoplus _{\grp\in
\Lambda}H_1((S^1\times S^2)_\grp;\bbZ)\): a 1-cycle
\(\ul{\lambda}_\grp\subset Y\) constructed
from the arc \(\lambda_\grp\subset M\) in Part 1, in a way parellel to the construction of
\(\gamma^{(z_0)}\) or \(\ul{\gamma}\) from \(\gamma_z\subset M\). 
\end{itemize}
Here is a more precise description of the cycles \(\ul{\lambda}_\grp\). For \(\grp=\grp_k\), \(k=1, \ldots,
\G\), let \(\ul{\lambda}_{\grp_k}\) be a 1-cycle in \(Y_k\) characterized
by the properties listed below. Recall the sphere \(S_{\grp_k}\)
from Part 1. Let
\(\mathcal{N}_{k, \varepsilon}\) denote the version of
\(\mathcal{N}_\varepsilon\) from Part 4 of Section \ref{sec:Da)} when
the sphere \(S\) therein is set to be \(S_{\grp_k}\). Then 
\begin{itemize}
\item on
\(Y_k-\mathcal{N}_{k, \varepsilon}\), \(\ul{\lambda}_{\grp_k}\) agrees
with \(\lambda_{\grp_k}\cup \hat{i}_{\grp_k}\);
\item on \(\mathcal{N}_{k, \varepsilon}\), \(\ul{\lambda}_{\grp_k}\) is
transverse to the spheres \(\{\rho\}\times S^2\) \(\forall \rho\in I\)
under the identification \(\mathcal{N}_{k, \varepsilon}\simeq
I\times S^2\) in Part 4 of Section \ref{sec:Da)}.
\end{itemize}
Recall that \(Y_{k+1}, \ldots, Y_{\G}=Y\) are constructed from \(Y_k\) by
iteratively connected summing with \(S^1\times S^2\),  and thus they
all contain a 1-cycle inherited from the \(\ul{\lambda}_{\grp_k}\subset
Y_k\) described above. We use the same
notation \(\ul{\lambda}_{\grp_k}\) for all such cycles in
\(Y_{k+1}, \ldots, Y_{\G}=Y\). 

With the above understood, the splitting 
(\ref{def:H_split2}) adopted in \cite{KLT3} is related to the
splitting (\ref{def:H_split3}) used in this article's Section
\ref{sec:1.1} via a transformation matrix of the following block form:
\begin{equation}\label{transf-matrix}\left[\begin{array}{ccc}
\bbX & 0 & 0\\
0 & \op{Id} & 0\\
\bbY & \vec{1} & \op{Id}
\end{array}\right],
\end{equation}
where \(\vec{1}\) denotes a row vector of all entries 1, \(\bbX\) is an automorphism of \(H_1(M;\bbZ)/\op{Tors}\), and \(\bbY\)
depends on the relative homology classes of \(\lambda_{\grp}\)'s. One may
 choose the arcs \(\lambda_\grp\)'s so that the entry \(\bbY\)
vanishes. {\em Such a choice of \(\lambda_\grp\)'s is adopted in this article.}

Use \((C_{ech}^\circ, \partial_{ech})\) to denote the underlying chain
complex of \(ech^\circ\), and let 
\((CF^\circ,\partial_{HF})\) be the Heegaard Floer
complex.  In \cite{KLT3}'s notation, the chain module \(C_{ech}^\infty\) is generated by the set 
\(\hat{\mathcal{Z}}_{ech, M}\), which is a \(\bbZ\)-bundle over
the set \(\mathcal{Z}_{ech, M}\).  The latter is written in
\cite{KLT3}'s (1.10) as a product of \(\mathcal{Z}_{HF}\), the
generating set for the Heegaard Floer chain module \(\widehat{CF}\),
and for each \(\grp\in
\Lambda\), a copy of \(\bbZ\times \textsc{o}\). This can be used to
write the \(ech\) chain module \(C_{ech}^\circ\) as a tensor product of \(CF^\circ\)
and, for each \(\grp\in
\Lambda\), a polynomial algebra \(C_\grp=\bbZ[\tau_\grp, \tau^{-1}_\grp, y_\grp^+,
y_\grp^-]\). Here, \(\tau_\grp\) is an even variable and corresponds to the
generator \(1\in \bbZ\) of the first factor in \(\bbZ\times \textsc{o}\),
and \(y_\grp^+, y_\grp^-\) are odd variables so that the polynomials
\(1, y_\grp^+, y_\grp^-, y_\grp^+y_\grp^-=-y_\grp^-y_\grp^+\)
correspond respectively to the elements \(0\), \(1\), \(-1\),
\(\{1,-1\}\) of \(\textsc{o}\) in \cite{KLT3}'s notation. 

Recall that a \(U\)-map on the \(ech\)-chain complex \(C_{ech}^\circ\), and for each \(M\)-adapted 1-cycle \(\hat{\gamma}\), a
map we shall denote by \(\grt_{\hat{\gamma }}\), were defined in
\cite{KLT2}'s Appendix and \cite{KLT3}'s Section 1. Together they
define the \(\mathbf{A}_\dag(Y)\)-action on \(ech^\circ\). 

Stated in the language of this article, Theorem 1.1 of \cite{KLT3}
asserts the following: With respect to the
aforementioned decomposition of the chain module 
\begin{equation}\label{ech-module}C_{ech}^\circ\simeq CF^\circ\otimes
\textstyle{\bigotimes}_{\grp\in \Lambda}C_\grp,\end{equation} 
\begin{itemize}
\item \((C_{ech}^\circ , \partial_{ech})\) is the product complex of
the Heegaard Floer chain complex \((CF^\circ, \partial_{HF})\) and 
for each \(\grp\in
\Lambda\), the chain complex 
\[(C_\grp, \partial_\grp): =\big.(\bbZ[\tau_\grp, \tau^{-1}_\grp, y_\grp^+,
y_\grp^-], (1+\tau_\grp)(\partial_{y_\grp^+}+\tau^{-1}_\grp \partial_{y_\grp^-})\big),\]
where \(\tau _\grp\) has degree 0 and \(y_\grp^+,
y_\grp^-\)  both have degree 1. Note that the homology
\(H(C_\grp, \partial_\grp)\) has two generators, one of degree 0 and
the other of degree 1, and they are respectively 
represented by the elements \(1\) and \(y_\grp^+-\tau_\grp y_\grp^-\)
in the polynomial algebra \(\bbZ[\tau_\grp, \tau^{-1}_\grp, y_\grp^+,
y_\grp^-]\). 
\item The \(U\)-map on \(C_{ech}^\circ\) acts only on the
\(CF^\circ\)-factor; namely,
\(U_{ech}=U_{HF}\otimes\bigotimes_{\grp\in \Lambda}\op{Id}\), and the
map \(U_{HF}\) on \(CF^\circ\) induces the \(U\)-action on
\(HF^\circ\). 
\item The \(\grt_{\gamma^{(z_0)}}\)-action on \(C_{ech}^\circ\) likewise has the form
\(\grt_{HF}^{(z_0)}\otimes\bigotimes_{\grp\in \Lambda}\op{Id}\) under the
decomposition (\ref{ech-module}).
\item The \(\grt_{\hat{i}^{(z)}}\)-action on \(C_{ech}^\circ\) has the form
\(\grt_{HF}^{(z)}\otimes\bigotimes_{\grp\in \Lambda}\op{Id}\) under the
decomposition (\ref{ech-module}), and the map \(\grt_{HF}^{(z)}\)
induces the action of \([\hat{i}^{(z)}]\in H_1(Y;\bbZ)/\op{Tors}\) on
\(HF^\circ\). 
\item For \(\grp=\grp_k\), \(k=1, \ldots, \G\), the
\(\grt_{\hat{i}_\grp}\)-action on on \(C_{ech}^\circ\) is the tensor
product of \(\partial_{y_\grp^+}\) on the \(C_\grp\)-factor, and 
\(\op{Id}\) on all other factors of \(C_{ech}^\circ\). Note that
with \(\partial_{y_\grp^+}\) identified as the generator of the algebra
\(H_{-*} (S^1)\simeq \bbZ[\partial_{y_\grp^+}]\simeq \bigwedge^*H_1((S^1\times S^2)_\grp;\bbZ)\), the homology
\(H(C_\grp, \partial_\grp)\) is identified with the module \(H_*(S^1)\) with the
standard \(H_{-*} (S^1)\) action. 
\end{itemize}
View \(HF^\circ\) as a module over \(\mathbf{A}_\dag(M)\otimes
\bigwedge^*H_1(S^1\times S^2;\bbZ)\), with the
\(\bigwedge^*H_1(S^1\times S^2;\bbZ)\)-factor part of the action
generated by the induced map from \(\grt_{HF}^{(z_0)}\) in the fourth
bullet above. It follows that there is an isomorphism between 
\(ech^\circ\) and  \(HF^\circ\boxtimes H_*(S^1)^{\boxtimes\G}\) as modules
over \(\big(\mathbf{A}_\dag(M)\otimes
H_{-*}(S^1)\big)\otimes H_{-*}(S^1)^{\otimes\G}\). Here, the  \(\big(\mathbf{A}_\dag(M)\otimes
H_{-*}(S^1)\big)\otimes H_{-*}(S^1)^{\otimes\G}\)-module structure
on \(ech^\circ\) comes from the decomposition
(\ref{def:H_split2}) to identify \(\mathbf{A}_\dag (Y)\) with 
\[\begin{split}
 \mathbf{A}_\dag (Y)
& \stackrel{i_{cycle}}{\simeq}\big(\mathbf{A}_\dag(M)\otimes
\textstyle{\bigwedge}^*H_1(S^1\times S^2;\bbZ)\big)\otimes\textstyle{\bigotimes }_{\grp\in
\Lambda}\textstyle{\bigwedge}^*H_1((S^1\times
S^2)_\grp;\bbZ)\\
& \simeq\big(\mathbf{A}_\dag(M)\otimes
H_{-*}(S^1)\big)\otimes H_{-*}(S^1)^{\otimes\G},\end{split}\]   
is isomorphic to the external tensor
product \(HF^\circ\boxtimes H_*(S^1)^{\boxtimes\G}\) as modules
over \(\big(\mathbf{A}_\dag(M)\otimes
H_{-*}(S^1)\big)\otimes H_{-*}(S^1)^{\otimes\G}.\) The two
factorizations of \(\mathbf{A}_\dag(Y)\), \(i_{sum}\)
in (\ref{eq:A-inclusion}) and above \(i_{cycle}\), are related via
(\ref{transf-matrix}) (where \(\bbY=0\)). 
According to Theorem \ref{thm:main} (2), the middle \(
H_{-*}(S^1)\)-factor in (\ref{eq:A-inclusion})'s factorization of
\(\mathbf{A}_\dag (Y)\) acts trivially on \(H^\circ(Y)\). Recalling from 
\cite{KLT4} that 
\(H^\circ(Y)\) and \(ech^\circ\) are canonically isomorphic as
\(\mathbb{A}_\dag(Y)\)-modules, this means that the \(\mathbf{A}_\dag(M)\otimes
H_{-*}(S^1)^{\otimes\G}\)-action on \(ech^\circ\) in the statement of Proposition
\ref{prop:ech-compute} is the same as the one arising from composing
the inclusion \(\mathbf{A}_\dag(M)\otimes
1 \otimes H_{-*}(S^1)^{\otimes\G}\hookrightarrow\mathbf{A}_\dag(M)\otimes
H_{-*}(S^1)\otimes H_{-*}(S^1)^{\otimes\G}\) with \(i_{cycle}\). The
assertion of the Proposition then follows from the isomorphism
\(ech^\circ\simeq HF^\circ\boxtimes H_*(S^1)^{\boxtimes\G}\)
described above. \epf

\subsection{Proof of Proposition \ref{prop:A.10b}}
\label{sec:Df)}

The construction of the cobordism manifold \(X\), its metric and self-dual 2-form has nine parts.

\paragraph{Part 1:} This part sets some of the notation for the
construction in the subsequent parts of the subsection of the desired
metric on \(X\) and the 2-form \(w_{X}\). Fix a
metric on \(Y\) of the sort that is described in Part 2 of Section \ref{sec:Ag)} and
denote the latter by \(\grg_{Y}\). The 2-form \(w\) on \(Y\) has
\(\grg_{Y}\) norm equal to 1 and its Hodge dual is the 1-form \(\hat{a}\)
that is described in Section II.3a; see also (IV.1.6). The constant
\(L\) for use in (\ref{eq:(A.12,15a)}) is specified at the end of the proof. Assume
until then that \(L > 100\) has been chosen. 

The description of the metric for \(X\) and the 2-form
\(w_{X}\) on the \(s \in [-L, -L + 8]\) part of \(X\) requires the
formula for \(w\) on a given \(\grp \in \Lambda \) version of
\(\mathcal{H}_\grp\) from (IV.1.3):

\begin{equation}\label{eq:(D.47)}\begin{split}
w & = 6 x  \cos\theta \sin\theta d\theta du - \sqrt{6}\,f'
\cos\theta \sin^{2}\theta \, du\,  d\phi \\
&\quad\quad  +\sqrt{6}\, f (1 - 3\cos^{2}\theta)
\sin\theta \, d\theta  \, d\phi .
\end{split}\end{equation}

The notation here uses \(x\) and \(f\) to denote a pair of non-negative
functions on \(\mathcal{H}_\grp\), these given in (IV.1.2),
with \(f'\) denoting the derivative of \(f\). Both \(x\) and \(f\) are invariant
under the reflection \(u \mapsto -u\). The function \(x\) has compact
support and is a non-zero constant where \(|u|
< 2\). This constant is denoted by \(x_{0}\). The function \(f\) on the \(|u| < 4\) part of
\(\mathcal{H}_\grp\) is given by the rule \(u \mapsto f(u) = x_{0} + 4 e^{-2R}\cosh(2u)\). 

The 1-form \(\upsilon_{\diamond}\) given in
(IV.1.5) plays a central role in what follows. This 1-form near on
the \(|u| < 4\) part of \(\mathcal{H}_\grp\) can be written as
\begin{equation}\label{eq:(D.48)}
\upsilon_{\diamond} = 4 e^{-2R}\cosh (2u) (1 - 3 \cos^{2} \theta)  \, du
+ 12\, e^{-2R} \sinh(2u) \cos\theta \sin\theta \, d\theta .
\end{equation}

The 1-form \(\upsilon_{\diamond}\) is a closed form on
\(Y\), and its zero locus are the loci in each \(\grp \in \Lambda\)
version of \(\mathcal{H}_\grp\) where both \(u\) and the function
\(1 - 3 \cos^{2} \theta \) are zero. Note also that \(*w = \upsilon_{\diamond}\) on the complement in \(Y\)
of the \(|u| \geq R + \ln\delta  -  9\) parts of each \(\grp \in \Lambda \) handle \(\mathcal{H}_\grp\). A
second point to note is that \(*(w \wedge\upsilon_{\diamond}) \geq c_{0}^{-1}|\upsilon_\diamond|^{2}\)
on the whole of \(Y\).

\paragraph{Part 2:} Let \(*\) denote for the moment the Hodge star of the
metric \(\grg_{Y}\) on \(Y\). The desired metric for \(X\) must pull
back to \((-\infty, -L] \times Y\) via the embedding from the
second bullet of (\ref{(A.9a,11)}) as the metric \(ds^{2 }+
\grg_{Y}\). Meanwhile, the corresponding pull-back of
\(w_{X}\) must equal \(ds \wedge *w + w\). This 2-form
is self-dual but it is not closed; this is because \(d*w \neq 0\) on
the \(|u| \leq R + \ln\delta - 9\) part of each \(\grp\in \Lambda \) version of \(\mathcal{H}_\grp\). This
last fact follows from the formula in (IV.1.6).

This rest of this part of the subsection describes \(w_{X}\)
for \(s \in [-L, -L + 3]\). The metric on this part of \(X\) still pulls
back as \(ds^{2} + \grg_{Y}\) via the second bullet of (\ref{(A.9a,11)}).

Let \(\chi_{\diamond 1}\) denote the function on \(\mathbb{R}\) given by the rule \(s \mapsto
\chi(-s - L + 2)\). This function is equal to 0 where \(s
< -L + 1\) and it is equal to 1 where \(s > -L + 2\). The derivative of \(\chi_{\diamond 1}\) is denoted in
subsequent equations by \(\chi_{\diamond 1}'\). Fix \(m> 1\) and introduce \(\chi_{m}\) to denote the
function of the coordinate \(s\) given by the rule \(s \mapsto
\chi(m |u| - 1)\). This function equals 0 where \(|u| > 2m^{-1}\) and it
equals 1 where \(|u| <m^{-1}\). By way of a look ahead, \(m\) will be set equal
to \(\r^{1/c_0\c}\) when the time comes to verify the requirements of Proposition \ref{prop:A.10a}.

Use \(w_{1}\) to denote the \(s\)-dependent 2-form on \(Y\) that is
equal to \(w\) on the \(M_{\delta} \cup\mathcal{H}_{0}\) part of \(Y\), and equal to the following
2-form below on each \(\grp \in \Lambda \) version of \(\mathcal{H}_\grp\): 
\begin{equation}\label{eq:(D.49)}\begin{split}
w_{1} & = d\big(x  (1 -\chi_{\diamond 1}\chi_{m})\,( 1 -
3\cos^{2}\theta) \, du\big) - \sqrt{6}\,f'
\cos\theta \sin^{2}\theta \, du \, d\phi \\
&\quad + \sqrt{6}\, f \, (1 - 3\cos^{2}\theta) 
\sin\theta \, d\theta \, d\phi .
\end{split}
\end{equation}

Note that \(|w_1| \leq c_{0}\). Meanwhile, \(\frac{\partial }{\partial s}w_{1} = d\b\) with \(\b = -x \chi_{\diamond 1}'
\chi_{m} (1 - 3\cos^{2}\theta) du\). As \(\chi_{m}= 0\)  where \(|u|
> 2m^{-1}\), the \(L^{2}\)-norm of \(\b\) on \([-L, -L + 3] \times Y\) is no greater than \(c_{0}
m^{-1}\). The appearance of \(\chi_{\diamond 1}\) in the definition guarantees that
\(w_{1} = w\) where \(s \leq -L\). Note that
\(w_{1}\) is a closed 2-form on \(Y\) for each \(s\). A key point to
note is that zero set of the \(s > -L + 1\) versions of
\(w_{1}\) consists of two circles in each \(\grp \in\Lambda \) version of \(\mathcal{H}_\grp\), these being the
circles where \(u\) and \(1 - 3 \cos^{2}\theta \) are both zero. 

The desired 2-form \(w_{X}\) pulls back to \([-L, - L + 3]
\times Y\) via the embedding from the second bullet of (\ref{(A.9a,11)}) as \(ds 
\wedge *w_{1} + w_{1}\).

\paragraph{Part 3:} What follows directly describes the desired metric
and the 2-form \(w_{X}\) on the \(s \in [-L + 3, -L + 4]\)
part of \(X\). To this end, let \(\chi_{\diamond 2}\) denote the
function on \(\mathbb{R}\) that is given by the rule \(s \mapsto
\chi(s  +  L - 3)\). This function is equal to 1 where \(s \leq
- L + 3\) and it is equal to 0 where \(s \geq - L + 4\). A smooth
metric on \(Y\) will be constructed momentarily whose Hodge star sends the \(s
\geq - L + 3\) versions of \(w_{1}\) to \(\upsilon_{\diamond}\), thus making \(w_{1}\)
harmonic. Let \(\grg_{1}\) denote this metric. Use \(\grg\) to
denote the \(s\)-dependent metric
\(\chi_{\diamond 2}\grg_{Y}+ (1 -\chi_{\diamond 2})\grg_{1}\) and let \(*\) now denote its Hodge dual. The metric
on \(X\) pulls back \([- L + 3, - L + 4] \times Y\) via the embedding
from the second bullet of (\ref{(A.9a,11)}) as \(ds^{2} + \grg\). The
pull back of \(w_{X}\) to \([- L + 3, - L + 4] \times Y\) is
the 2-form \(ds \wedge *w_{1} + w_{1}\).
This 2-form is self-dual when \(s\) is near \(- L + 4\). The two steps that
follow construct the metric \(\grg_{1}\).

\paragraph{\it Step 1: } The 2-form \(w_{1}\) is equal to \(w\) on the
\(M_{\delta} \cup \mathcal{H}_{0}\) part of \(Y\) and its \(\grg_{Y}\) Hodge star here is
\(\upsilon_{\diamond}\). This understood, the metric
\(\grg_{1}\) on \(M_{\delta} \cup\mathcal{H}_{0}\) is set equal to \(\grg_{Y}\). To
define \(\grg_{1}\) on a given \(\grp \in \Lambda \) version of
\(\mathcal{H}_\grp\), note first that the function
\(\chi_{\diamond1}\) in (\ref{eq:(D.49)}) is
equal to 1 when \(s \in [- L + 3, - L + 4]\). This implies that
\(w_{1}\) is \(s\)-independent when \(s \in [- L + 3, - L +
4]\). More to the point, it also implies that the \(s \in [- L + 3,
- L + 4]\) version of \(w_{1}\) shares the same zero locus
with the closed 1-form \(\upsilon_{\diamond}\), this being the circles in each \(\grp \in \Lambda \) version of
\(\mathcal{H}_\grp\) where \(u\) and \(1 -3\cos^{2}\theta \) are both zero. Meanwhile
\(w_{1} \wedge \upsilon_{\diamond}> 0\) on the complement of their common zero locus. This
last observation can be used with Lemma \ref{lem:D.2} to construct the desired
metric \(\grg_{1}\) on any part of the complement in
\(\mathcal{H}_\grp\) of the \(u = 0\) and \(1 -
3\cos^{2}\theta = 0\) locus as a smooth extension of
the metric \(\grg_{Y}\) from \(M_{\delta} \cup\mathcal{H}_{0}\). 

\paragraph{\it Step 2:} Let \(\mathcal{T} \subset \mathcal{H}_\grp\)
denote the \(|u| <   m^{-1}\) part of \(\mathcal{H}_\grp\). The function
\(\chi_{m}\) in (\ref{eq:(D.36)}) is equal to 1 on \(\mathcal{T}\) and \(f = x_{0}  + 4 e^{-2R}
\cosh (2u)\) on \(\mathcal{T}\). This being the case, it follows from
(\ref{eq:(D.47)}) and (\ref{eq:(D.48)}) that the metric on
\(\mathcal{T}\) with volume 3-form \(\Omega= \sin\theta \, du\,
d\theta \, d\phi \) and Hodge star
defined by the rules:
\begin{equation}\label{eq:(D.50)}
\begin{cases}
& * \sin\theta d\theta d\phi = \frac{1}{\sqrt{6}}\frac{4e^{-2R}\cosh
    (2u)}{x_0+4r^{-2R}\cosh (2u)} du,\\
&* \sin\theta d\phi\, du   =\frac{3}{2\sqrt{2}}d\theta ,\\
&* du \, d\theta =  \frac{3}{2\sqrt{2}}\sin\theta \, d\phi 
\end{cases}
\end{equation}
sends \(w_{1}\) to \(\upsilon_{\diamond}\).
Note that a suitable change of coordinates near the \(\theta = 0\)
and \(\theta = \pi  \) loci can be used to prove that the metric
defined by (\ref{eq:(D.50)}) is smooth on the whole of \(\mathcal{T}\). 

As noted previously, Lemma \ref{lem:D.2} can be used to extend the metric defined
in (\ref{eq:(D.50)}) to the whole of \(\mathcal{H}_\grp\) so as to agree
with \(\grg_{Y}\) on \(\mathcal{H}_\grp \cap M_{\delta}\). This must be done with some care so as
to obtain an \(m =r^{1/c_0\c}\)
  extension that can be used to satisfy the second item of (\ref{eq:(A.16)}).
With this goal in mind, note that Lemma \ref{lem:D.2} can be used to find an
extension with the following three properties: 
\BTitem\label{eq:(D.51)}
\item The norm of the Riemannian curvature tensor and those of its
covariant derivatives to order \(20\) are bounded by \(c_{0}\). 
\item The injectivity radius is bounded from below by \(c_{0}^{-1}\).
\item The metric volume of \(Y\) is at most \(c_{0}\).
\ETitem

The first bullet of Lemma \ref{lem:D.2} gives metrics that obey the third bullet
of (\ref{eq:(D.51)}) and the second bullet of Lemma \ref{lem:D.2}
supplies metrics that obey all three bullets.

\paragraph{Part 4:} The desired metric for \(X\) and the 2-form
\(w_{X}\) on the \(s \in [-L + 4, - \frac{3}{4}L + 2]\) portion of \(X\) are described below. This is done by specifying
their pull-backs via the embedding from the second bullet of (\ref{(A.9a,11)}) to
\([-L + 4, - \frac{3}{4}L + 2] \times Y\). 
In this part, we use
\(\chi_{\diamond 2}\) to denote the function on \(\mathbb{R}\) given by the rule \(s \mapsto \chi\big( \frac{4}{L-20}
 (s +  L - 5)\big)\). This function is equal to 1 where \(s < -L +
5\) and it is equal to zero where \(s > -\frac{3}{4} L\). Use \(\chi_{\diamond 2}'\) to
denote the derivative of \(\chi_{\diamond 2}\). Note in particular that
\(|\chi_{\diamond 2}'|\leq c_{0} L^{-1}\).

Let \(w_{2}\) denote the \(s\)-dependent 2-form on \(Y\) given by
\(w_{1}\) for \(s < -L + 4\), given by \(w\) on
\(M_{\delta} \cup \mathcal{H}_{0}\), and given on each \(\grp \in \Lambda \) version of
\(\mathcal{H}_\grp\) for \(s \geq -L + 4\) by
\begin{equation}\label{eq:(D.52)}\begin{split}
w_{2} &=\chi_{\diamond 2} \, d\big(x (1 -\chi_{m}) (1 -
3\cos^{2}\theta)\, du\big) - \sqrt{6}f' \cos\theta \sin^{2}\theta \, du\,
d\phi \\
& \quad + \sqrt{6} f (1 - 3\cos^{2}\theta) \sin\theta \, d\theta \, d\phi .
\end{split}
\end{equation}

The 2-form \(w_{2}\) is a closed 2-form on \(Y\) for each \(s\), it has
the same zero locus as \(w_{1}\) and it has the property that
\(w_{2} {\wedge\upsilon}_{\diamond} =w_{1} \wedge \upsilon_{\diamond}\).

An \(s\)-dependent metric on \(Y\) is described momentarily for the cases
when \(L > c_{0}\). This metric is denoted by
\(\grg\). Let \(*\) denote the corresponding Hodge dual. By way of a look
ahead, \(\grg\) is chosen so that \(d*w_{2}=\frac{\partial}{\partial s}
 w_{2}\). The pull back of the desired metric on \(X\) to \([-L +
4, -\frac{3}{4}L +2] \times Y\) via the embedding from the second bullet of (\ref{(A.9a,11)}) is
the quadratic form \(ds^{2} + \grg\), and the corresponding
pull back of \(w_{X}\) is \(ds \wedge *w_{2}+ w_{2}\). Note in particular that \(w_{X}\) is
self-dual and closed if self-duality is defined by the metric \(ds^{2} + \grg\). 

The metric \(\grg_{1}\) from Part 3 is \(s\)-independent and so it
is defined where \(s > -L + 4\). This understood, the metric
\(\grg\) is set equal to \(\grg_{1}\) where \(s < -L + 5\). It
is also set equal to \(\grg_{1}\) for all \(s \in [-L + 4, - \frac{3}{4}L +2]\) on \(M_{\delta} \cup
\mathcal{H}_{0}\). This is to say that it equals
\(\grg_{Y}\) for all such \(s\) on \(M_{\delta}\cup \mathcal{H}_{0}\).  The metric \(\grg\) is chosen
where \(s \geq -L + 5\) on each \(\grp \in \Lambda \) version of
\(\mathcal{H}_\grp\) so that its Hodge star on each \(\grp \in\Lambda \) version of \(\mathcal{H}_\grp\) acts on
\(w_{2}\) as 
\begin{equation}\label{eq:(D.53)}
* w_{2} =\chi_{\diamond 2}' x \,( 1 -\chi_{m}) (1 - 3\cos^{2}\theta)
\, du+ \upsilon_{\diamond}.
\end{equation}

As will be explained directly, if \(L > c_{0}\), there are
metrics of the sort just described that obey the \(c_0 =1 \) version of
(\ref{eq:(D.51)}) where \(s > - \frac{3}{4} L + 1\).

To see about these requirements, consider first constructing a
metric of the desired sort where \(s > -\frac{3}{4} L\). The metric
that is defined by (\ref{eq:(D.50)}) with volume form \(\sin\theta \, du\,
d\theta \, d\phi \) satisfies the requirements where \(|u| < 2\). Since \(w_{2}\wedge \upsilon_{\diamond} > 0\) on
the rest of \(\mathcal{H}_\grp\) and the \(\grg_{Y}\) Hodge star of \(w_{2}\) is
\(\upsilon_{\diamond}\) on \(M_{\delta}\cup \mathcal{H}_{0}\), Lemma \ref{lem:D.2} finds an extension
of the latter metric from the \(|u| < 1\) part
of each \(\mathcal{H}_\grp\) that has the desired properties.
Use \(\grg_{2}\) to denote this \(s\)-independent metric.

Consider next the story where \(s <  -\frac{3}{4}L + 1\). The metric on any given
\(\grp \in \Lambda \) version of \(\mathcal{H}_\grp\) that is defined by (\ref{eq:(D.50)}) with volume
form   \(\sin\theta \, du \, d\theta \, d\phi \) has Hodge star sending \(w_{2}\) to \(\upsilon_{\diamond}\) where
\(|u| < m^{-1}\). Let \(\upsilon\) denote the 1-form on the right hand side of (\ref{eq:(D.53)}). The
3-form \(\upsilon  \wedge w_{2}\) can be written where \(|u| \geq   \frac{1}{2}m^{-1}\) as \(\grq \upsilon_{\diamond}
\wedge w_{2}\) and it follows from the fact that
\(|\chi_{\diamond 2}'|< c_{0} L^{-1}\) that \(\grq> c_{0}^{-1} -c_{0} L^{-1}\). Thus, \(\upsilon \wedge w_{2} >  0\) where
\(|u| > \frac{1}{2}m^{-1}\). Given this positivity and given what was said in
the preceding paragraphs, Lemmas \ref{lem:D.2} and \ref{lem:D.3} can be used to construct
an \(s\)-dependent metric where \(s < -\frac{3}{4} L + 1\) that equals \(\grg_{2}\) where \(s > - 
\frac{3}{4}L +\frac{1}{2}\) , that equals \(\grg_{1}\) where \(s < -L + 5\) and equals
\(\grg_{Y}\) on \(M_{\delta} \cup \mathcal{H}_{0}\).

\paragraph{Part 5:} This part and Part 6 construct the desired metric for
\(X\) and the 2-form \(w_{X}\) where \(s \in[-\frac{3}{4} L + 1,
-\frac{1}{2} L + 2]\). By way of a look ahead, the metric pulls back from this part
of \(X\) via the embedding from the second bullet of (\ref{(A.9a,11)}) as
\(ds^{2} + \grg_{3}\) with \(\grg_{3}\)
being an \(s\)-dependent metric on \(Y\) that equals the metric
\(\grg_{*}\) for all \(s\) on the set \(\mathcal{Y}_{0\varepsilon }\) from (\ref{eq:(D.10)}). 

The metric \(\grg_{3}\) is independent of \(s\) on the whole of \(Y\) when
\(s \in [-\frac{1}{2}L + 1, -\frac{1}{2}L + 2]\). This \(s\)-independent version of \(\grg_{3}\) is in a
large \(T\) version of the space \(\op{Met}_{T}\) that is defined in
Part 5 of Section \ref{sec:Da)}. For the purposes to come, the choice of \(T\)
requires choosing \(L > c_{T }\) with
\(c_{T}\) denoting here and in what follows a constant that
depends on \(T\) and is greater than \(c_{0}T^{2}\)
in any event. The value of \(c_{T}\) may increase between
appearances. 

Use \(*\) now to denote the \(\grg_{3}\) Hodge star on \(Y\). The
2-form \(w_{X}\) pulls back via the embedding from the second
bullet of (\ref{(A.9a,11)}) to \([-\frac{3}{4}L + 1, -\frac{1}{2} L + 2] \times Y\) as \(ds \wedge *w_{3} +
w_{3}\), with \(w_{3}\) denoting an \(s\)-dependent,
closed 2-form on \(Y\). The 2-form \(w_3\) is also independent
of \(s\) where \(s \in [- \frac{1}{2} L + 1, -\frac{1}{2}L + 2]\) and it is independent of \(s\) on
\(\mathcal{Y}_{0\varepsilon }\) for all \(s\). With regards to the
motivation for what follows below and in Part 6, keep in mind that \(ds \wedge *w_{3} +
w_{3}\) is closed if and only if both \(dw_{3} =
0\) and \(d(*w_{3}) = \frac{\partial}{\partial s}w_{3}\) for all \(s\).

This part of the subsection makes the assumption that
\(c_1(\det(\mathbb{S}))\) annihilates the \(H_{2}(M; \mathbb{Z})\) summand of the direct sum
decomposition for \(H_{2}(Y; \mathbb{Z})\) given in
(IV.1.4). This assumption makes for a simpler construction. Even
so, much of what is done here is used again for Part 6's construction for the general case. 

The construction that follows has six steps. Note that some of these
steps use notation from Section \ref{sec:Da)}.

\paragraph{\it Step 1:} Let \(\chi_{\diamond 3}\)
denote the function of \(s\) given by \(\chi\big(\frac{3}{L-8}  (s
+\frac{3}{4}  L  - 2)\big)\). This function equals 1 for \(s < -
\frac{3}{4}L + 2\) and it equals 0 for \(s \geq -
\frac{1}{2}L\). Reintroduce the notation from Section \ref{sec:Da)}
and let \(\chi_{\ir}\) denote the function on \(\mathbb{R}^{3}\) given by
\(\chi(64\, \varepsilon_{*}^{-1}(\ir- \rho_{*}) - 1)\). This function equals 1
where \(\ir < \rho_{*} + \frac{1}{64}\varepsilon\) and it equals 0 where \( \ir >
\rho_{*} +\frac{1}{32} \varepsilon\). Given \(T \geq 1\) and use \(\chi_{\r}\)
with \(\chi_{\diamond 3}\) to define the \(s\)-dependent function on \(\mathbb{R}^{3}\) given by
\begin{equation}\label{eq:(D.54)}
\ir_{s_T} =\chi_{\diamond 3}\ir + (1 - \chi_{\diamond 3}) \big(1- \chi_{\ir} + \frac{1}{T}\chi_{\r} \big) \ir.
\end{equation}

Note in particular that \(\frac{\partial}{\partial s}\r_{s_T} > 0\) because
\(\chi_{\ir}\) is a non-increasing function of \(\ir\). Use \(\rho_{s_T}\) and
\(x_{s_T3 }\) to denote the respective \(s\)-dependent functions on \(\mathbb{R}^{3}\) given by
\(\ir_{s_T}\sin\theta \) and \(\ir_{s_T} \cos\theta\). 

Define the \(s\)-dependent 2-form \(w_{3}\) on \(Y\) by setting
\(w_{3} = w_{2}\) for \(s \leq -\frac{3}{4}L + 2\) and setting it equal to \(w\) on the \(\mathcal{Y}_{0}\)
component of \(Y-\mathcal{N}_{\varepsilon}\). The 2-form \(w_{3}\) is defined on
\(\mathcal{N}_{\varepsilon}\) by specifying it on the \(\mathbb{R}^{3}\) incarnation of
\(\mathcal{N}_{\varepsilon}\) to be \(\textsc{k}(\rho_{s_T})\,
\rho_{s_T}\, d\rho_{s_T} \, d\phi\). The
definition of \(w_{3}\) on the rest of \(Y\) uses \(\tau\) to
denote the function of \(s\) given by
\((\chi_{\diamond 3} + (1 - \chi_{\diamond 3})/T)^{2}\). The latter function equals 1 where \(s
< - \frac{3}{4} L + 2\) and it is equal to \(\frac{1}{T^2}\)
 where \(s > - \frac{1}{2}L\). The 2-form \(w_{3}\) is defined on
\(\mathcal{Y}_{M} \cap M_{\delta}\) to be \(\tau w_{2}\); and it is defined on each \(\grp \in
\Lambda \) version of \(\mathcal{H}_\grp\) in the upcoming
(\ref{eq:(D.55)}). This upcoming definition uses
\(\chi_\Delta\) to denote the function of \(u\) and
\(\theta\) given by \(\chi(|u|^{2}- 1) \, \chi(4 (1 - 3\cos^{2}\theta) - 1)\). The
function \(\chi_\Delta\)   is equal to 1 where both \(|u| < 1\) and \(|1 -
3\cos^{2}\theta| <\frac{1}{4}\) ; and it is equal to 0 where either \(|u|
> 2\) or \(|1 -3\cos^{2}\theta|  >\frac{1}{2}\). Note in particular that the support of
\(\chi_{\Delta}\) consists of two open sets. These
are mirror images under the involution \(\theta \mapsto \pi- \theta\),
with one being a neighborhood of the \(u = 0\) and \(\cos\theta = \frac{1}{\sqrt{3}}\)
  circle with \(0 < \theta <\frac{\pi}{2}\) on its closure. Define 
\begin{equation}\label{eq:(D.55)}
w_{3} = -\sqrt{6} \tau  \, d\big(f \cos\theta\sin^{2}\theta d\phi - (x_{0} + 4
e^{-2R})   \op{sign} (\cos\theta) \chi_\Delta d\phi\big)
\end{equation}
on \(\mathcal{H}_\grp\) for \(s>\frac{3}{4}L+2\).

By way of comparison, the 2-form \(w_{2}\) on
\(\mathcal{H}_\grp\) can be written as \(\sqrt{6} \, d(f\cos\theta \sin^{2}\theta d\phi)\). What is
written in (\ref{eq:(D.55)}) adds a 2-form with support on
\(\mathcal{H}_\grp\) to \(\tau w_{2}\).

The 2-form \(w_{3}\) on \(Y\) is closed for each \(s\). Moreover, it
defines the \(s\)-independent de Rham cohomology class
\(c_{1}(\det(\mathbb{S}))\) because the latter class is
assumed to annihilate the \(H_{2}(M; \mathbb{Z})\) summand
in (IV.1.4). 

\paragraph{\it Step 2:} The \(s\)-dependent metric \(\grg_{3}\) is defined when \(s
\in [-\frac{3}{4}  L + 1, - \frac{1}{2}L + 2]\) with the help of a
certain \(s\)-dependent 1-form, \(\b\). The 1-form \(\b\) should obey
\(d\b =  \frac{\partial}{\partial s}w_{3}\). There are four additional constraints on \(\b\). The
first is that \(\b\) should vanish on \(\mathcal{Y}_{0}\) and on
the part of \(\mathcal{N}_{\varepsilon}\) where \(\r
> \rho_{*} +\frac{1}{16} \varepsilon\). The second constraint specifies \(\b\) on the
\(|u| < 4\) part of \(\mathcal{H}_\grp\): 
\begin{equation}\label{eq:(D.56)}
\b = -\sqrt{6} \tau' \big(f \cos\theta\sin^{2}\theta  - (x_{0} + 4
e^{-2R})   \op{sign}(\cos\theta) \chi_\Delta\big) d\phi ,
\end{equation}
where \(\tau'\) denotes \(\frac{\partial}{\partial s}\tau\). The third constraint asks that \(\b\)'s norm at
\(s \in [-\frac{3}{4}  L + 1, - \frac{1}{2}L + 2]\) when measured by the metric \(\grg_{Y}\) obeys \(|\b|_{\grg_-} \leq c_{T} L^{-1}\). The fourth
contraint requires the following: Fix \(k \in \{0, \cdots, 20\}\).
Then the \(\grg_{Y}\)-covariant derivatives up to order \(20\) of
\((\frac{\partial}{\partial s})^{k}\b\) are bounded by \(c_{T}L^{-k-1}\).

To see about satifying these constraints, note first that \(\b\) can be
chosen to vanish on \(\mathcal{Y}_{0}\) and on the \(\ir
> \rho_{*} +\frac{1}{16}\varepsilon \) part of \(\mathcal{N}_{\varepsilon}\)
because \(w_{3}\) is constant on these parts of \(Y\), and because
the first cohomology of the \(\ir \in [\rho_{*} +\frac{1}{32}
  \varepsilon, \rho_{*} +\frac{1}{64}\varepsilon]\) part of \(\mathcal{N}_{\varepsilon}\)
is zero. The \(c_{0} L^{-1}\) bound on \(|\chi_{\diamond 3}'|\)
implies that \(\b\) can be chosen to vanish on
\(\mathcal{Y}_{0}\) and so that its norm elsewhere when
measured by the metric \(\grg_{Y}\) is bounded by
\(c_{0} L^{-1}\). A 1-form of this sort can be chosen so that the
\(\grg_{Y}\) norms of its derivatives also have the required norm
bound. Let \(\b_{*}\) denote such a choice, and let \(\b_{\Lambda}\) denote the 1-form on any
given \(\grp \in \Lambda \) version of \(\mathcal{H}_\grp\) given by (\ref{eq:(D.47)}). Their difference, \(\b_{*} -
\b_{\Lambda}\), is a closed 1-form on \(\mathcal{H}_\grp\). As
\(H^{1}(\mathcal{H}_\grp \cap M_{\delta};  \mathbb{R}) = 0\), this difference can be
written as \(d\k\) with \(\k\) denoting a function on
\(\mathcal{H}_\grp\). The function \(\k\) can be taken so that
\(|\k | \leq c_{0}L^{-1} \)since the \(\grg_{Y}\) norms of both
\(\b_{*}\) and \(\b_{\Lambda}\) obey a similar \(c_{0} L^{-1}\) bound. Granted this bound
on \(\k\), then \(\b = \b_{*} - d(\chi(|u| - 4)\k)\) has all of the requisite
properties.

\paragraph{\it Step 3:} The definition of the upcoming Steps 4 and 6 use
observations made below about \(w_{3}\) and \(\b\) on the
\(|u| \leq 4\) part of each \(\grp \in \Lambda\) version of \(\mathcal{H}_\grp\). The first series of
observations concern \(w_{3}\). To start, note that the zero
locus of the 2-form in (\ref{eq:(D.56)}) is the same as that of
\(\upsilon_{\diamond}\), this being the locus where
both \(u = 0\) and \(1 - 3\cos^{2}\theta = 0\). The reason being that \(f'\)  and  
 \(\chi_\Delta\) have the same sign where \(\chi_\Delta \neq 0\), and likewise the
functions \((1 - 3\cos^{2}\theta)\) and
\(\op{sign}(\cos\theta) \, \chi_\Delta\) have the same sign where 
 \(\chi_\Delta \neq 0\). In fact, these comments about the derivatives of \(\chi_\Delta\)
imply that \(w_{3}\) on \(\mathcal{H}_\grp\)  can be written schematically as
\begin{equation}\label{eq:(D.57)}\begin{split}
w_{3} &=  -(1 + \textsc{a}_{1})\tau \sqrt{6}\, f'\cos\theta
\sin^{2}\theta \, du\, d\phi \\
&\quad + (1 + \textsc{a}_{2}) \tau
\sqrt{6} \, f (1 - 3\cos^{2}\theta) \sin\theta \, d\theta \, d\phi  
\end{split}
\end{equation}
where \(\textsc{a}_{1}\) and \(\textsc{a}_{2}\) are
smooth, non-negative functions of \(u\) and \(\theta\) that equal zero
where both \(|u| < 1\) and \(|1 - 3\cos^{2}\theta| <\frac{1}{4}\) and where either \(|u| > 2\) or \(|1
- 3\cos^{2}\theta|  >\frac{1}{2}  \). Given that \(w_{2}\) on \(\mathcal{H}_\grp\)  is
\(-\sqrt{6}\, d(f \cos\theta \sin^{2}\theta d\phi)\), these last remarks imply that 
\begin{equation}\label{eq:(D.58)}
w_{3} \wedge \upsilon_{\diamond}\geq \tau w_{2}
\wedge\upsilon_{\diamond}\quad  \text{on \(\mathcal{H}_\grp\)}
\end{equation}
with the inequality being a strict one only where \(d\chi_\Delta \neq 0\).

The next series of remarks concern the 1-form \(\b\) on the
\(|u| \leq 4\) part of
\(\mathcal{H}_\grp\). The first point of note being that \(f(u) \cos\theta \sin^{2}\theta \) is equal to
\((x_{0} + 4 e^{-2R}) \frac{2}{3\sqrt{3}}\op{sign}(\cos\theta)\) on the zero locus of
\(\upsilon_{\diamond}\). It follows as a consequence that \(\b\) can be written as
\begin{equation}\label{eq:(D.59)}
\b = -\textsc{b}_{1 }\tau' f' \cos\theta
\sin^{2} \theta d\phi + \textsc{b}_{2 }\tau' f (1 -3\cos^{2}\theta) \sin\theta d\phi ,
\end{equation}
where \(\textsc{b}_{1}\) and \(\textsc{b}_{2}\) are
smooth functions of \(u\) and \(\theta\).

\paragraph{\it Step 4:} The metric \(\grg_{3}\) on each \(\grp \in \Lambda\)
version of \(\mathcal{H}_\grp\) is defined to be the metric
from Part 5 for \(s < -\frac{3}{4} L + 2\). The metric \(\grg_{3}\) on
\(\mathcal{H}_\grp\) at other values of \(s\) is defined in part
so that its Hodge star obeys
\begin{equation}\label{eq:(D.60)}
*w_{3} = \tau\upsilon_{\diamond} + \b.
\end{equation}

There is one other constraint. To explain it, note first that the
metric \(\grg_{2}\) does not depend on \(s\) when if \(s \in [-\frac{3}{4}
L + 1, -\frac{3}{4} L + 2]\). Use \(\grg_{2+}\) to denote this \(s\)-independent metric.
Look at (\ref{eq:(D.45)}) to see that the \(s > -\frac{1}{2}  L + 1\) version of \(w_{3}\) on the \(|u|
> 4\) part of each \(\mathcal{H}_\grp\) is  \(\frac{1}{T^2}w_{2}\). Since \(\b\) is zero when \(s > -\frac{1}{2}
L + 1\), the constraint in (\ref{eq:(D.60)}) is satisfied by taking the Hodge star
to be that defined by \(\grg_{2+}\). This understood, the final
constraint is as follows:
\begin{equation}\label{eq:(D.61)}\begin{split}
& \text{The metric \(\grg_{3}\) on each \(\grp \in\Lambda  \) version of
\(\mathcal{H}_\grp\)  when \(s > - \frac{1}{2}L + 1\) must be both}\\ 
 & \text{\(s\)-independent and \(T\)-independent; and it must equal
\(\grg_{2+}\) where \(|u|> 4\).}
\end{split}
\end{equation}

As explained in what follows, an \(s\)-dependent metric with all of these
requisite properties exists if \(L\) is greater than a \(T\)-dependent
constant.

Consider first the existence of a metric with the desired properties
where \(|u| < 1\) and \(|1 -3\cos^{2}\theta| <\frac{1}{4}\), this
being a neighborhood of the common zero locus of \(w_{3}\) and \(\upsilon_{\diamond}\). The
metric \(\grg\) is defined on this part of \(\mathcal{H}_\grp\) by
its volume 3-form \(\Omega= \sin\theta \, du\, d\theta \, d\phi\)
and the Hodge duals
\begin{equation}\label{eq:(D.62)}
\begin{cases}
& * \sin\theta \, d\theta \, d\phi = \frac{1}{\sqrt{6}}\frac{4e^{-2R}\cosh
  (2u)}{x_0+4e^{-2R}\cosh (2u)}du + \tau^{-1}\tau'\textsc{b}_{2}
\sin\theta \, d\phi,\\
& *\sin\theta  \, d\phi \, du = \frac{\sqrt{ 3}}{2\sqrt{2}}d\theta - \frac{1}{\sqrt{6}}\tau^{-1}\tau' \textsc{b}_{1}
\sin\theta \, d\phi ,\\
& *du\, d\theta =  \frac{\sqrt{ 3}}{2\sqrt{2}} \sin\theta \, d\phi  +\frac{1}{\sqrt{6}}\tau^{-1}\tau' \textsc{b}_{2} du
- \frac{1}{\sqrt{6}}\tau^{-1}\tau'\textsc{b}_{1}d\theta .
\end{cases}
\end{equation}

These formulas for the Hodge dual define a symmetric, bilinear form on
the cotangent bundle of this part of \(\mathcal{H}_\grp\).
This bilinear form is positive definite if
\(\tau^{-1}|\tau'| <c_{0}^{-1}\), which is guaranteed if
\(T^{2}L^{-1} < c_{0}^{-1}\) since \(\tau^{-1} < T^{2}\) and \(|\tau'| < c_{0}L^{-1}\). 

To see about defining \(\grg_{3}\) on the rest of
\(\mathcal{H}_\grp\), use the fact that \(|\b| \leq c_{0} L^{-1}\) to draw
the following conclusion: If \(L >
c_{0}T^{2}\), then \(w_{3}\wedge (\tau \upsilon_{\diamond} + \b)> 0\) on the complement in \(Y\) of the \(|u|
< \frac{1}{2}\)    and \(|1 - 3\cos^{2}\theta|< \frac{1}{8} \) part of each \(\grp \in \Lambda \) version of
\(\mathcal{H}_\grp\). This being the case, then Lemma \ref{lem:D.3}
can be used directly to obtain a family of metrics on
\(\mathcal{H}_\grp\) parametrized by the set \([-\frac{3}{4}L + 1, -\frac{1}{2} L + 2]\) so as to obey (\ref{eq:(D.60)}) and (\ref{eq:(D.61)}). Use
\(\grg_{3\Lambda}\) to denote this family of metrics on \(\bigcup_{\grp\in\Lambda}\mathcal{H}_\grp\).

\paragraph{\it Step 5:} The 1-form \(\upsilon_{\diamond}\) is used
here to construct another closed, \(s\)-dependent 1-form that plays a
central role in the upcoming definition of the \(s \in [-\frac{3}{4} L + 1, -\frac{1}{2}L + 2]\) versions of \(\grg_{3}\) on \(M_{\delta}\cup \mathcal{H}_{0}\). This new 1-form is denoted
by \(\upsilon_{\diamond 3}\) and its definition is given in the subsequent paragraph. 

The 1-form \(\upsilon_{\diamond 3}\) on \(\mathcal{Y}_{0}\) is
\(\upsilon_{\diamond}\) and it is defined on the \(\ir 
> \rho_{*} - \frac{1}{4}  \varepsilon \) part of \(\mathcal{N}_{\varepsilon}\) to
be \(dx_{s_T3}\) with the latter defined in
Step 1. Since \(\upsilon_{\diamond} =dx_{3}\) on \(\mathcal{N}_{\varepsilon}\), it
follows from the definition of \(x_{s_T3}\)
that \(\upsilon_{\diamond 3}\)  as defined so far is a 1-form on the union of
\(\mathcal{Y}_{0}\) and the \(\ir >\rho_{*} - \frac{1}{4}  \varepsilon \) part of \(\mathcal{N}_{\varepsilon}\).
The definition of \(\upsilon_{\diamond 3}\) on the \(\ir \in [\rho_{*}
-\frac{1}{2} \varepsilon, \rho_{*}-\frac{1}{4}  \varepsilon]\) part
of \(\mathcal{N}_{\varepsilon}\) requires the reintroduction of the function
\(\chi_{\ir *}\) from Step 2 in Part 5 of
Section \ref{sec:Da)}. This function is used here to define
\(x_{s_T3*} =\big(\chi_{\diamond 3}+ (1 - \chi_{\diamond 3}) (1
-\chi_{\ir *} +  \frac{1}{T}\chi_{r*}) \big)\, x_{3}\).
Define \(\upsilon_{\diamond 3}\) on
the \(\ir \in [\rho_{*} - \frac{1}{2} \varepsilon, \rho_{*} -\frac{1}{4} \varepsilon]\) part of \(\mathcal{N}_{\varepsilon}\)
to be  \(\tau^{1/2}dx_{s_T3*}\). It follows
from the definitions of \(x_{s_T3}\) and
\(x_{s_T3*}\) that the definition just given
defines a smooth 1-form on the union of \(\mathcal{Y}_{0}\)
with the \(\ir > \rho_{*} -\frac{1}{2}\varepsilon \) part of \(\mathcal{N}_{\varepsilon}\).
As the latter's restriction near the \(\ir =
\rho_{*} -\frac{1}{2}\varepsilon \) is \(\tau dx_{3}\), a smooth 1-form on
\(\mathcal{Y}_{0} \cup\mathcal{N}_{\varepsilon}\) is defined by setting
\(\upsilon_{\diamond3} = \tau dx_{3}\) on the \(\ir \leq \rho_{*}- \frac{1}{2}\varepsilon \) part of \(\mathcal{N}_{\varepsilon}\).
Noting that \(\tau dx_{3} = \tau\upsilon_{\diamond}\), defining \(\upsilon_{\diamond 3}\) on \(\mathcal{Y}_{M}\) to be \(\tau
\upsilon_{\diamond}\) defines a smooth, closed 1-form on \(Y\). 

The 1-form \(\upsilon_{\diamond 3}\) has the four properties that are listed below.

\paragraph{\sc Property 1:} {\em The 1-form \(\upsilon_{\diamond 3}\) is
equal to \(\upsilon_{\diamond}\) where \(s\in [- \frac{3}{4}L + 1, - \frac{3}{4} L + 2].\)}

This follows because \(\chi_{\diamond 3}= 1\) at these values of \(s\).

\paragraph{\sc Property 2:} {\em The zero locus of each \(s
\in [- \frac{3}{4}L +1, -\frac{1}{2}L+ 2]\) version of \(\upsilon_{\diamond 3}\) is identical to that of \(\upsilon_{\diamond}\). }

This is because \(\upsilon_{\diamond 3}\)  has no zeros
on \(\mathcal{Y}_{0} \cup\mathcal{N}_{\varepsilon}\) and it is equal to
\(\tau \upsilon_{\diamond}\) on \(\mathcal{Y}_{M}\).

\paragraph{\sc Property 3:} {\em Each \(s \in [- \frac{3}{4}L + 1, - \frac{1}{2}L + 2]\) version of \(w_{3} \wedge
\upsilon_{\diamond 3}\) is positive on the 
complement of the common zero locus of \(w_{3}\) and \(\upsilon_{\diamond 3}\).}

This property follows directly from the definitions on
\(Y-(\bigcup_{\grp\in\Lambda}\mathcal{H}_\grp)\) and from (\ref{eq:(D.57)}) on each \(\grp \in
\Lambda \) version of \(\mathcal{H}_\grp\).

To set the stage for the fourth property, note that \(w_{3}\) and \(\upsilon_{\diamond 3}\) do not
depend on \(s\) when \(s \in [- \frac{3}{4} L + 1, -\frac{3}{4}L +
2]\). Use \(w_{3+}\) and \(\upsilon_{\diamond 3+}\)  to denote these \(s\)-independent differential forms. To continute the stage setting,
let \(\grg_{3\Lambda_+}\) denote the \(s\)-independent metric on
\(\bigcup_{\grp\in\Lambda}\mathcal{H}_\grp\) given by the \(s \in [- \frac{3}{4} L + 1, -\frac{3}{4}L + 2]\) version of Part 5's metric
\(\grg_{3\Lambda}\). What with (\ref{eq:(D.51)}), this metric on \(\bigcup_{\grp\in\Lambda}\mathcal{H}_\grp\) with \(\grg_{2+}\) on
\(Y-\{\bigcup_{\grp\in\Lambda}\mathcal{H}_\grp\}\) define a smooth, \(s\)- and \(T\)- independent
metric on \(Y\). Denote the latter by \(\grg_{\diamond}\). The restriction of \(\grg_{\diamond}\) to
\(\mathcal{Y}_{M} \cup\mathcal{N}_{\varepsilon}\) is in the space
\(\op{Met}^{\mathcal{N}}\) from Part 5 of Section \ref{sec:Da)}. This understood, let \(\grg_{\diamond T}\) denote the
\(\op{Met}_T\) metric that is constructed in Part 5 of Section \ref{sec:Da)}
from \(T\) and \(\mathcal{Y}_{M} \cup \mathcal{N}_{\varepsilon}\) part
of \(\grg_{\diamond}\).

\paragraph{\sc Property 4:}{\em The \(\grg_{\diamond T}\) Hodge star of \(w_{3+}\) is
\(\upsilon_{\diamond 3+}\). }

The definitions in Part 5 of Section \ref{sec:Da)} with those given above for
\(w_{3+}\) and \(\upsilon_{\diamond 3+}\) imply this on
\(Y-(\bigcup_{\grp\in\Lambda} \mathcal{H}_\grp)\) and
(\ref{eq:(D.60)}), (\ref{eq:(D.61)}) imply this on \(\bigcup_{\grp\in\Lambda}\mathcal{H}_\grp\). 

\paragraph{\it Step 6:} This step completes the definition of \(\grg_{3}\) on \(Y\)
so as to satisfy five constraints, the first being that
\(*w_{3} =\upsilon_{\diamond 3} + \b\) at each \(s\in [-\frac{3}{4} L
+ 1, -\frac{1}{2} L + 2]\). The second contraint asks that the \(s
\in[-\frac{3}{4}L + 1, -\frac{3}{4} L + 2]\) versions are independent
of \(s\); and the third asks that the
\(s \in [- \frac{1}{2} L + 1, -\frac{1}{2}  L + 2]\) versions are also independent of \(s\) and
that this \(s\)-independent metric is \(\grg_{\diamond T}\). The fourth
constraint asks that \(\grg_{3} =\grg_{3\Lambda}\) on the
\(|u| < R + \ln\delta \) part of each \(\grp\in \Lambda \) version of \(\mathcal{H}_\grp\). The
fifth and final constraint asks that \(\grg_{3} =
\grg_{*}\) on \(\mathcal{Y}_{0}\) and on the\( \r
>  \rho_{* } + \frac{1}{16}  \varepsilon \) part of \(\mathcal{N}_{\varepsilon}\).

Use \textsc{Property 3} and what is said in Step 4 with the bound
\(|\b|_{\grg_-}< c_{0} L^{-1}\) to see that
\(w_{3} \wedge(\upsilon_{\diamond 3} + \b)>  0\) on the complement in \(Y\) of the common zeros of
\(w_{3}\) and \(\upsilon_{\diamond 3}\)  if \(L \geq c_{T}\). Given this bound, Lemma \ref{lem:D.3} with the input from
Step 4 and \textsc{Property 4} of Step 5 find a metric with all of the
desired properties. Take such a metric for \(\grg_{3}\). Note
for future reference that the \(s\)-independent, \(s > - \frac{1}{2}
 L + 1\) version of \(\grg_{3}\) is equal to \(\grg_{Y}\) on \(\mathcal{Y}_{M} \cap M_{\delta}\).

\paragraph{Part 6:} This part of the subsection puts no constraints on
the restriction of \(c_{1}(\det(\mathbb{S}))\) to the
\(H_{2}(M; \mathbb{Z})\) summand in \(H_{2}(Y;
\mathbb{Z})\). The \(s\)-dependent metric \(\grg_{3}\) and the
2-form \(w_{3}\) in this case are identical to their namesakes
in Part 5 on \(Y-(\bigcup_{(\gamma, \textsc{z}_\gamma)\in\Theta}\mathcal{T}_{\gamma})\).
The three steps that follow define \(\grg_{3}\) and
\(w_{3}\) on \(\bigcup_{(\gamma, \textsc{z}_\gamma)\in\Theta}\mathcal{T}_{\gamma}\).

\paragraph{\it Step 1:} Reintroduce from Part 7 of Section \ref{sec:Da)} the closed 2-form \(\p\)
on \(Y\). By way of a reminder, the de Rham class of \(\p\) has pairing 0 with
the \(H_{2}(\mathcal{H}_{0}; \mathbb{Z})\oplus (\bigoplus_{\grp\in\Lambda}\mathcal{H}_\grp)\) summand in (IV.1.4)'s
decomposition of \(H_{2}(Y; \mathbb{Z})\) and its pairing
with the \(H_{2}(M; \mathbb{Z})\) summand is the same as
that of \(c_{1}(\det(\mathbb{S}))\). Since \(\p\)'s support lies in \(\bigcup_{(\gamma, \textsc{z}_\gamma)\in\Theta}\mathcal{T}_{\gamma}\)
and thus in
\(\mathcal{Y}_{M}-(\bigcup_{\grp\in\Lambda}\mathcal{H}_\grp)\), setting \(w_{3}\) on \(\mathcal{Y}_{M} \cup
(\bigcup_{\grp\in\Lambda}\mathcal{H}_\grp)\) to be \(w_{3} =  \tau w_{2} + (1 - \tau) \p\) defines a closed 2-form on \(Y\)
for each \(s \in[-\frac{3}{4}  L + 1, - \frac{1}{2} L + 2]\) with de Rham cohomology class
\(c_{1}(\det(\mathbb{S}))\). 

The metric \(\grg_{3}\) is defined on \(\bigcup_{(\gamma, \textsc{z}_\gamma)\in\Theta}\mathcal{T}_{\gamma}\)
so that its Hodge star maps \(w_{3}\) to \(\tau\upsilon_{\diamond} + \b\) with \(\b\) denoting a certain
1-form with \(d\b = \frac{\partial}{\partial   s}w_{3}\). As done previously, Lemma \ref{lem:D.3} will be used to
construct a metric with this property that meets all of the other requirements.

\paragraph{\it Step 2:} The definition of \(\grg_{3}\) and \(\b\) on \(\bigcup_{(\gamma, \textsc{z}_\gamma)\in\Theta}\mathcal{T}_{\gamma}\)
requires what is said here about the \(w_{2}\) and \(\p\) in the
support of \(\p\). To start, reintroduce from Part 7 of Section \ref{sec:Da)} the set
\(\Theta\) and write \(\p\) as  \(\sum_{(\gamma, \textsc{z}_\gamma)\in\Theta}\textsc{z}_{\gamma}\p_{\gamma}\) with each \((\gamma,\textsc{z}_{\gamma})\) version of
\(\p_{\gamma}\) being a closed 2-form with support in the
tubular neighborhood \(\mathcal{T}_{\gamma}\) that is
described in Part 7 in Section \ref{sec:Da)}. Part 7 of Section \ref{sec:Da)}
describes a diffeomorphism from \(S^{1} \times D\) to
\(\mathcal{T}_{\gamma}\) with \(D\) denoting a small radius
disk about the origin in \(\mathbb{R}^{2}\). The
diffeomorphism identifies \(\gamma\) with \(S^{1}\times \{0\}\) and it
has two important property that concern the 2-form \(w\)  on \(Y\) and the function \(\ff\) from Section II.1. As noted in
Part 7 of Section \ref{sec:Da)}, the 1-form \(d\ff\) pulls back via the embedding
of \(S^{1} \times D\) as a constant 1-form on the \(D\)
factor and the kernel of the pull back via the
embedding of the 2-form \(w\) is a
constant vector field that is tangent to this \(D\) factor. These last
properties are exploited in the next paragraph.

As can be seen in (IV.1.5), the 1-form
\(\upsilon_{\diamond}\) on \(\mathcal{T}_{\gamma}\) is \(d\ff\). Meanwhile, the
2-form \(w_{2}\) on \(\mathcal{T}_{\gamma}\) is
still the original 2-form \(w\) on \(Y\) as described in (IV.1.3). This
understood, what was said above about \(d\ff\) and the kernel of \(w\) imply
that \(S^{1} \times D\) has coordinates \((t, (x, y))\) with
\(t\) denoting an affine coordinate for the \(S^{1}\) factor and
\((x, y)\) coordinates for \(D\) with the following two properties: The
1-form \(\upsilon_{\diamond}\) pulls back as \(dx\) and the 2-form \(w_{2}\) pulls back as
\(\textsc{h}_{\gamma}(y, t)\, dy \, dt\) with \(\textsc{h}_{\gamma}\) denoting a positive function.
Granted these coordinates, the 2-form \(\p_{\gamma}\) has
the form \(\grh_{\gamma}(x, y)\, dx\, dy\) with
\(\grh_{\gamma}\) denoting a function with compact support in
a small radius disk about the origin in the \((x, y)\)-plane and with total
integral equal to 1.

\paragraph{\it Step 3:} An almost verbatim repeat of
what is said in Step 2 of Part 6 supplies a version of the 1-form
\(\b\) which obeys the four properties
listed in the first paragraph of Step 2 in Part 6 with it understood
that \(w_{3}\) is now defined as in Step 1. 

It follows as a consequence of what is said in Step 2 that
\begin{equation}\label{eq:(D.63)}
(\tau w_{2}+(1-\tau)\p) \wedge \upsilon_{\diamond}= \tau w_{2}\wedge\upsilon_{2},
\end{equation}
and thus the \(\grg_Y\)-norm of \((\tau w_2+(1-\tau)\p)\wedge(\tau\upsilon_{\diamond}+\b)\)
is no less than \(\tau^{2}(c_{0}^{-1}-c_{T}T^{2}L^{-1})\).
This being the case, Lemma \ref{lem:D.3} supplies an
\(s\)-dependent metric on \(Y\) with all
of the desired properties if \(L\) is larger than a purely \(T\)-dependent
constant.

Let \(\grg_{3+}\) denote the \(s\)-independent metric
on \(Y\) given by the \(s\in [- \frac{1}{2} L + 1, -\frac{1}{2}L + 2]\) versions of \(\grg_{3}\). This is the metric
\(\grg_{Y}\) on \((\mathcal{Y}_{M} \cap M_{\delta})-(\bigcup_{(\gamma, \textsc{z}_\gamma)\in\Theta}\mathcal{T}_{\gamma})\).
It proves necessary for what follows to take some care with regards to
the choice of \(\grg_{3+}\) on \(\bigcup_{(\gamma, \textsc{z}_\gamma)\in\Theta}
\mathcal{T}_{\gamma}\). In particular, Lemmas \ref{lem:D.2} and \ref{lem:D.3} will construct a version of
\(\grg_{3}\) with \(\grg_{3+}\) on each \(\mathcal{T}_{\gamma}\) by \(\grg_{Y}\)-volume 3-form
\(\textsc{h}_{\gamma}dx\, dy\, dt\) and the Hodge star rules: 

\begin{equation}\label{eq:(D.64)}
\begin{cases}
&*dx \, dy  =\textsc{a}_{\gamma}dt -\textsc{a}_{\gamma}\tau^{-1}(1-\tau)\,\textsc{h}_{\gamma}^{-1}\textsc{z}_{\gamma}\grh_{\gamma}dx +
\textsc{b}_{\gamma}dy, \\
&*dy\, dt  =\textsc{h}_{\gamma}^{-1}(1+ \tau^{-2}\textsc{h}_{\gamma}^{-1}\textsc{a}_{\gamma}(1-\tau)^{2}\textsc{z}_{\gamma}
\grh_{\gamma})dx \\
&\qquad \qquad \quad -\textsc{a}_{\gamma}\tau^{-1}(1- \tau)\textsc{h}_{\gamma}^{-1}\textsc{z}_{\gamma}\grh_{\gamma}dt,\\
& *dt\, dx = dy +\textsc{b}_{\gamma}dt,
\end{cases}
\end{equation}
with \(\textsc{a}_{\gamma}\) being a positive function and with \(\tau\) equal to \(\frac{1}{\tau^2}\). The function
\(\textsc{a}_{\gamma}\) is constrained for the moment only to the extent that
\(\textsc{a}_{\gamma}<c_{0}^{-1}\tau^{2}\) on the support of \(\textsc{z}_{\gamma}\grh_{\gamma}\) and that
\(\textsc{a}_{\gamma}\) is independent of \(T\) on the complement in
\(\mathcal{T}_{\gamma}\) of a \(T\)-independent open set that contains the support of \(\grh_{\gamma}\)
and has compact closure in \(\mathcal{T}_{\gamma}\). This set is denoted by \(\mathcal{T}_{\gamma}'\).
This upper bound on \(\textsc{a}_{\gamma}\)
is needed so that (\ref{eq:(D.64)}) defines a positive definite metric. As for
\(\textsc{b}_{\gamma}\), it is zero on \(\mathcal{T}_{\gamma}\)
and it  is independent of \(T\) elsewhere.

\paragraph{Part 7:} This part of the subsection defines the
desired metric on \(X\) and 2-form \(w_{X}\)
on the \(s\in [-\frac{1}{2}L + 1, -\frac{1}{2} L + 5]\) part of \(X\). As done previously, these are
defined by their pull-backs via the embedding from the second bullet of
(\ref{(A.9a,11)}). The pull-back of the metric will have the form
\(ds^{2}+ \grg\) with \(\grg\) denoting an \(s\)-dependent metric on \(Y\).
Meanwhile, the pull-back of \(w_{X}\) will have the form \(ds \wedge*w_{4}+w_{4}\), with \(w_{4}\)
denoting a closed, \(s\)-dependent 2-form on \(Y\) and with \(*\) denoting
the Hodge \(*\) defined by \(\grg\). The de Rham cohomology class of \(w_{4}\)
at each \(s\) is \(c_{1}(\det(\mathbb{S}))\).

The metric \(\grg\) is independent of \(s\) for
\(s\in [-\frac{1}{2}L + 1, -\frac{1}{2}L + 2]\) and the 2-form \(w_{4}\)
is independent of \(s\) for
\(s\in [- \frac{1}{2}L + 1, -\frac{1}{2}L + 3]\). Both the metric and \(w_{4}\) are
independent of \(s\) when
\(s \in[-\frac{1}{2}L + 4, -\frac{1}{2}L + 5]\). Moreover, the restriction of both to
\(Y-(\bigcup_{\grp\in\Lambda}\mathcal{H}_\grp)\) are independent of \(s\) for all values of
\(s\). The salient difference between the \(s \leq -\frac{1}{2}
 L + 3\) version of \(w_{4}\) and the \(s \geq - \frac{1}{2}
 L + 4\) version being that the latter has nondegenerate zeros and the
 former does not. 

The construction of \(\grg\) and \(w_{4}\)
has two steps.

\paragraph{\it Step 1:} Let \(\grg_{3+}\) denote the \(-\frac{1}{2}
 L + 2\) version of the metric that is supplied in Parts 5 and 6, and let
\(w_{3+}\) denote the \(s = -\frac{1}{2}L + 2\) version of \(w_{3}\). The 2-form
\(w_{3+}\) is \(\grg_{3+}\)-harmonic but it does not vanish
transversely. By way of a reminder, the zero locus of
\(w_{3+}\) consists of the two circles in
each \(\grp \in \Lambda \) version of \(\mathcal{H}_\grp\)
where both \(u = 0\) and \(1 - 3\cos^{2}\theta = 0\).
Note in this regard that \(w_{3+}\) on
\(\mathcal{H}_\grp\)  is the 2-form

\begin{equation}\label{eq:(D.65)}
\sqrt{6} T^{-2} (- f' \cos\theta
\sin^{2}\theta \, du\, d\phi + f (1 -3\cos^{2}\theta)  \sin\theta \, d\theta \, d\phi).
\end{equation}

The construction of \(w_{4}\) starts by introducing
\(\chi_{\diamond 4}\) to denote the function on \(\mathbb{R}\) given
by \(s \mapsto \chi(s  +  L - 3)\). 
This function is equal to 1 where \(s < - 
 \frac{1}{2}L + 3\) and it is equal to 0 where \(s > - \frac{1}{2}L + 3\). The derivative of
\(\chi_{\diamond 3}\) is denoted by \(\chi_{\diamond 3}'\). Use
\(\chi_{*}\) to denote the function of \(u\) given by the
rule \(u \mapsto \chi(|u| - 1)\). This function is equal to 1 where \(|u| \leq 1\) and it is
equal to 0 where \(|u| > 2\). One last
function is needed for what follows, this denoted by
\(\chi_{\theta}\). It is a function on \([0, \pi ]\)
with values in \([0, 2]\) which has the following two properties: It is
zero near the endpoints, and has two local minima at the two values of
\(\theta\) where \(1 - 3\cos^{2}\theta = 0\). Moreover, \(\chi_{\theta}\) should appear on a neighborhood of
these minima as \(1 + (1 -3\cos^{2}\theta)^{2}\). Take \(\chi_{\theta}\) so that
\(\chi_{\theta}(\theta) =\chi_{\theta}(\pi  - \theta)\).

Fix \(z >1 \) and define the 2-form \(w_{z}\) by 
\begin{equation}\label{eq:(D.66)}\begin{split}
w_{z} & = - \big(\sqrt{6}\, f'\cos\theta\sin^{2}\theta  +
z^{-1}\cos\phi\, \chi_{\diamond 4}\, \chi_{*}\sin\theta  
\frac{\partial}{\partial\theta}\chi_{\theta}\big) \, du\, d\phi \\
&\quad +  \sqrt{6}\, f(1- 3\cos^{2}\theta)  \sin\theta \, d\theta \, d\phi  
- z^{-1}\sin\phi \, \chi_{\diamond 4}\, \chi_{*} \frac{\partial}{\partial\theta}\big(\sin\theta  \frac{\partial}{\partial\theta}\chi_{\theta}\big)\, du\, d\theta .
\end{split}
\end{equation}

This is a closed 2-form for all \(s\) that equals \(w_{3+}\) for \(s
\leq -\frac{1}{2}L + 3\) and for all \(s\) where \( |u|  > 2\).
This 2-form is independent of \(s\) when \(s \geq -\frac{1}{2}
 L + 4\). Moreover, if \(z > c_{0}\), then the \(s\)-independent
 version of \(w_{z}\) defined where \(s \geq - \frac{1}{2}
 L + 4\) has a nondegenerate zero locus, this being the four points where
\(\sin\phi = \)0, \(1 - 3\cos^{2}\theta = 0\) and \(u = 0\).

The desired 2-form \(w_{4}\) is defined to be \(w_{3+}\) on
\(Y-(\bigcup_{\grp\in\Lambda}\mathcal{H}_\grp)\) and it is defined on each \(\grp \in
\Lambda \) version of \(\mathcal{H}_\grp\) to be a \(z> c_{0}\) version of \(T^{-2}w_{z}\).

\paragraph{\it Step 2:} This step defines the metric \(\grg\). This is done by first
constructing \(\grg\) near the zero locus of \(w_{4}\) in each \(\grp
\in \Lambda \) version of \(\mathcal{H}_\grp\) and then
extending the result to the rest of \(Y\) with the help of Lemma \ref{lem:D.3}. 

Fix \(z > c_{0}\) so that \(w_{z}\) as
defined in (\ref{eq:(D.66)}) has nondegenerate zeros. The 2-form 
 \(w_{z}\) can be written as \(d\b_{z}\) where
\(\b_{z}\) is given by 
\begin{equation}\label{eq:(D.67)}\begin{split}
& \frac{\sqrt{3}}{2\sqrt{2}}z^{-1} \cos\phi 
\chi_{\diamond 4}\,( \chi_{*}'\chi_{\theta}du + \chi_{*} \chi_{\theta}'
d\theta) - \frac{\sqrt{3}}{2\sqrt{2}}z^{-1} \sin\phi
\chi_{\diamond 4}\chi_{*}\chi_{\theta} d\phi \\
&\quad + z^{-1} \sin\phi
\chi'_{\diamond 4}\, \chi_{*}\sin\theta \frac{\partial}{\partial\theta}\chi_{\theta}   du .
\end{split}
\end{equation}

Granted this formula, then \(\upsilon + \b_{z}\) has the
same zero locus as \(w_{z}\) if \(z >
c_{0}\), and it also vanishes transversely. Moreover,
\(w_{z} \wedge (\upsilon_{\diamond}+ \b_{z})\) can be written as
\(\textsc{q} \sin\theta \,  du\, 
d\theta \, d\phi \) and a calculation finds that \(\textsc{q} \geq
0\) with equality only on the joint zero locus of \(w_{z}\) and
\(\upsilon_{\diamond} + \b_{z}\). In fact, the calculation finds \(\textsc{q} \geq c_{0}^{-1}
(|u|^{2} + (1 -3\cos^{2}\theta)^{2} +z^{-2} \sin^{2}\phi
\sin^{2}\theta)\) if \(z > c_{0}\).

With \(z\) large and \(w_{4}\) defined by (\ref{eq:(D.66)}) on
\(\mathcal{H}_\grp\), the metric \(\grg\) is defined near the zeros
of (\ref{eq:(D.67)}) so that its Hodge star sends \(w_{z}\) to
\(\upsilon_{\diamond} + \b_{z}\). The
definition requires the introduction of yet another function of \(s\), this
denoted by \(\chi_{\diamond\diamond 4}\) and defined by the rule whereby
\(\chi_{\diamond\diamond 4}(s) =\chi(s + \frac{1}{2} L - 2)\). This function equals 1 where \(s < -\frac{1}{2}
 L + 2\) and it equals 0 where \(s > -\frac{1}{2}  L + 3\). The desired metric
 \(\grg\) is defined by taking its volume form to be \(\sin\theta \, du\, d\theta \, d\phi \) and its Hodge star to act as
follows:
\begin{equation}\label{eq:(D.68)}
\begin{cases}
&* \sin\theta \, d\theta \, d\phi  = \frac{1}{\sqrt{6}}
 \big(4 e^{-2R} \cosh\, (2u)  + 12 \, z^{-1} 
\sin\phi  \, \chi_{\diamond 4}'\cos\theta \sin^{2}\theta\big) du .\\
& * \sin \theta \, d\phi \, du = \frac{\sqrt{3}}{2\sqrt{2}}d\theta .\\
& * du\, d\theta   =\frac{\sqrt{3}}{2\sqrt{2}}
 \big(\chi_{\diamond\diamond 4}\sin\theta  +
 (1-\chi_{\diamond\diamond 4})\, \chi_{\theta} \,( \frac{\partial}{\partial\theta}(\sin\theta  \frac{\partial}{\partial\theta}
 \chi_{\theta}))^{-1} \big) d\phi.
\end{cases}
\end{equation}

By way of a parenthetical remark, the metric \(\grg_{3+}\) near
the zeros of \(w_{z}\) is defined by the same volume form but
with Hodge star rule given by (\ref{eq:(D.50)}). The appearance of
\(\chi_{\diamond\diamond 4}\) in the third line of (\ref{eq:(D.68)})
guarantees that \(\grg = \grg_{3+}\) where \(s\leq -\frac{1}{2}L\). 

As noted previously, \(w_{z} \wedge
(\upsilon_{\diamond} + \b_{z})> 0\) on the complement of the common zero locus of
\(w_{z}\) and \((\upsilon_{\diamond} + \b_z)\).
Having constructed \(\grg\) on a neighborhood of this locus with the desired
properties, Lemma \ref{lem:D.3} provides an extension to the whole of \(Y\) which is
independent of \(s\) where \(s < - \frac{1}{2}L + 2\), where \(s > -\frac{1}{2} L + 4\) and on
\(Y-(\bigcup_{\grp\in\Lambda}\mathcal{H}_\grp)\). This extension is such that the
2-form \(ds \wedge *w_{4} + w_{4}\) is self-dual on \([- \frac{1}{2}L + 1, - \frac{1}{2}L
+ 5] \times Y\) when self duality is defined by the metric \(ds^{2} + \grg\).

\paragraph{Part 8:} This part of the subsection supplies the input
for the definition in Part 9 of the desired metric and the 2-form
\(w_{X}\) on the \(s \in [-\frac{1}{2}L + 4, L]\) part of \(X\). The discussion in this section
refers to an auxiliary copy of the space \(X\), this denoted by
\(X_{*}\). The manifold \(X_{*}\) is the same
as \(X\), but its metric is not a metric of the sort that
is described in Parts 1-7. The eight steps that follow construct a
metric on \(X_{*}\) and a corresponding self dual 2-form
with certain desirable properties. 

\paragraph{\it Step 1:} Fix a metric in the \(Y_{\G}\) version of
\(\op{Met}^{\mathcal{N}}\). The latter with a sufficiently
large choice for \(T\) determines metrics in the set
\(\op{Met}(Y_{\G})\). This understood, choose \(T\) large
enough so that this is the case and so that two additional requirements
are met, the first being that Part 7's metric \(\grg\) and
2-form \(w_{4}\) can be constructed for any choice of \(L
> c_{T}\) with \(c_{T}\) denoting a
constant that is greater than 1 and depends only on \(T\). The second
requirement is given in Step 2. 

Let \(\grg_{-}\) and \(w_{-}\) denote the
respective \(s \in [ -\frac{1}{2}L + 4,  -\frac{1}{2} L + 5]\)
versions of \(\grg\) and \(w_{4}\), these being independent of \(s\). The metric \(\grg_{-}\) is in \(Y\)'s version
of the space \(\op{Met}_T\), so it can be used
for the metric \(\grg_{1}\) in Part 1 of Section \ref{sec:Db)}, and since
\(w_{-}\) has non-degenerate zeros, it can also be used for
the metric \(\grg_{2}\) in Part 1 of Section \ref{sec:Db)}. This part of
Section \ref{sec:Db)} uses \(w_{2}\) to denote the \(\grg_{2}\)-harmonic 2-form with de Rham cohomology class that of
\(c_{1}(\det(\mathbb{S}))\). This 2-form \(w_{2}\)
is \(w_{-}\). The 2-form \(w_{-}\) is equal to
\(w\) on \(\mathcal{Y}_{0}\) and on the \(\ir \geq\rho_{*} +\frac{5}{8}  \varepsilon \) part of \(\mathcal{N}_{\varepsilon}\)
and so it follows that \(w_{-}\) is also the 2-form that is
denoted by \(w_{3}\) in Part 2 of Section \ref{sec:Db)}. This fact
implies that the metric \(\grg_{-}\) is also a version of what
Part 2 of Section \ref{sec:Db)} denotes as \(\grg_{3T}\).
Parts 1-10 of Section \ref{sec:Dd)} will be invoked in the
upcoming steps using \(X_{*}\) and the \(\grg_{-}\)
version of \(\grg_{3T}\). These parts of Section \ref{sec:Dd)} denote the
latter version of \(\grg_{3T}\) by \(\grg_{-T}\).
What Parts 1-10 of Section \ref{sec:Dd)} denote as
\(w_{-T}\) in this case is the 2-form \(w_{-}\).

\paragraph{\it Step 2:} Let \(\grg_{\diamond}\) denote the given metric from
the \(\op{Met}(Y_{\G})\). By way of a reminder, the metric \(\grg_{\diamond}\) is
determined in part by Step 1's chosen metric from the
\(Y_{\G}\) version of \(\op{Met}^{\mathcal{N}}\) and \(T\).

As explained in Part 1 of Section \ref{sec:Db)}, a metric denoted
by \(\grg_{2}\) determines various versions of the metric \(\grg_{3T}\); and \(\grg_{\diamond}\)
can be any one of these \(\grg_{3T}\) metrics. Set
\(\grg_{+}\) to be the version of \(\grg_{2}\) that is used to construct \(\grg_{\diamond}\) and set \(\grg_{+T}\)
to denote \(\grg_{\diamond}\). What follows is the second
requirement for \(T\): It should be large enough so that the
\(Y_{-} = Y\) and \(Y_{+} =Y_{\G}\)
version of the constructions in Parts 1-10 from Section \ref{sec:Dd)} can be
invoked using \(X_{*}\) and the metrics \(\grg_{-}\)
on \(Y_{-}\) and \(\grg_{+}\) on \(Y_{+}\).

The constructions in Parts 1-8 of Section \ref{sec:Dd)} require a
closed 2-form on \(X_{*}\), this denoted by
\(\p_{X}\), whose de Rham cohomology class is
\(c_{1}(\det(\mathbb{S}))\) which has the following
additional properties: It equals \(w_{-}\) where \(s
< -102\), it equals \(w_{+}\) where \(s >
102\) and it obeys the bound in (\ref{eq:(D.20)}). Given such a
2-form, Parts 1-8 of Section \ref{sec:Dd)} supply \(L_{1} \gg 1\), a
metric on \(X_{*}\), and 
a 2-form on \(X_{*}\) with the properties listed below. The metric and 2-form are
denoted in the list and subsequently by \(\grm_{T*}\) and
\(\omega_{T*}\).

\BTitem\label{eq:(D.69)}
\item The metric \(\grm_{T*}\) obeys (\ref{eq:(A.12,15a)}) and (\ref{eq:(A.15b)}) when the version of \(L\)
in the latter is greater than \(L_{1}+ 20\). 
\item The pull-back of \(\grm_{T*}\) from the
\(s <-L_{1}-1\) part of \(X\) via the embedding from the
second bullet of (\ref{(A.9a,11)}) is
\(ds^{2}+\grg_{-}\) and the pull back of \(\grm_{T*}\) from the \(s> L_{1}+ 1\)
part of \(X_*\) by the embedding from the
third bullet of (\ref{(A.9a,11)}) is \(ds^{2}+\grg_{+}\).
\item The 2-form \(\omega_{T*}\) is self dual when self duality is
  defined by \(\grm_{T*}\). In addition, the pull-back of \(\omega_{T*}\) to any constant, \(s >1\)
  slice of \(X_{*}\) is closed.
\item The pull back of \(\omega_{T*}\) from the \(s< -L_{1}- 1\)
part of \(X_*\) by the embedding from the second bullet of
(\ref{(A.9a,11)}) is \(ds\wedge*w_{-}+w_{-}\) with \(*\) denoting here
the \(\grg_{-}\)-Hodge star.
\item The pull back of \(\omega_{T*}\) from the \(s>L_1+1\) part of \(X_*\)
via the embedding from the third bullet of (\ref{(A.9a,11)}) is
\(ds\wedge*w_{+}+w_{+}\) with \(*\) now denoting the
\(\grg_{+}\)-Hodge star and with \(w_{+}\) denoting the
\(\grg_{+}\)-harmonic 2-form with de Rham cohomology class
\(c_{1}(\det(\mathbb{S}))\).
\item The 2-form \(\omega_{T*}\) obeys the constraint in (\ref{eq:(A.13c)}).
\item The norm of \(\omega_{T*}\) and those of its \(\grm_{T*}\)-covariant
derivatives to order 10 are less than \(c_0\).
\ETitem

When comparing the notation in (\ref{eq:(D.69)}) with the notation in
Parts 1-10 of Section \ref{sec:Dd)}, keep in mind that this case
has \(\grg_{-T} = \grg_{-}\) and \(w_{-T} = w_{-}\), and \(\grg_{+T} =\grg_{\diamond}\) and \(w_{+T} =w_{\diamond}\).

The remaining steps construct a version of \(\p_{X}\) with
the required properties.

\paragraph{\it Step 3:} The construction of \(\p_{X}\) requires the three
constraints on \(\grm_{T*}\) that are described
here and a fourth constraint that is described in Step 4. The first
constraint is that imposed in Part 10 of Section \ref{sec:Dd)}.

The remaining constraints and that in Step 4 refer to the subset
\(\bigcup_{(\gamma, \textsc{z}_\gamma)\in\Theta}\mathcal{T}_{\gamma}\subset M_{\delta}\), this viewed as a subset of \(Y\)
and also as a subset of \(Y_{\G}\). The second constraint
uses the embeddings from the first and second bullets
of (\ref{(A.9a,11)}) to view the \(s < 0\) and \(s > 0\) parts of
\(X_{*}\) as \((-\infty, 0] \times Y\) and as \((0,\infty) \times  Y_{\G}\). This constraint is the
analog of that given in (\ref{eq:(D.42)}).
\begin{equation}\label{eq:(D.70)}\begin{split}
&\text{The metric \(\grm_{T*}\) on
\([-100, -96] \times\mathcal{Y}_{M}\) is the product metric
\(ds^{2}+\grg_{Y}\).}\\
&\text{ The metric \(\grm_{T*}\) on \([96,
100]\times\mathcal{Y}_{M}\) is the product metric
\(ds^{2}+\grg_{+}\).}
\end{split}
\end{equation}

By way of background for the third constraint, note that (\ref{eq:(D.43)}) holds
for \(X_{*}\), this being a consequence of what is said in
Part 1 about the ascending and descending manifolds from the critical
points of \(s\). The third constraint refers to this embedding. It also
uses \(\grm_{Y}\) and \(\grm_{+}\) to denote the metrics
\(ds^{2} + \grg_{Y}\) and \(ds^{2} +\grg_{+}\) on \(\mathbb{R} \times\bigcup_{(\gamma, \textsc{z}_\gamma)\in\Theta}\mathcal{T}_{\gamma}\).
\begin{equation}\label{eq:(D.71)}\begin{split}
&\text{There exists a \(T\)-independent constant, \(c_{*}> 1\), with}\\
&\text{the following significance: The pull-back of \(\grm_{T*}\)  from}\\
&\text{the \(s >-94\) part of \(X_{*}\) via the embedding in (\ref{eq:(D.43)}) obeys}\\
&\text{\(c_{*}^{-1}\grm_{Y}\leq \grm \leq c_{*}\grm_{Y}\) and
\(c_{*}^{-1}\grm_{+}\leq \grm \leq c_{*}\grm_{+}\).}
\end{split}
\end{equation}

This third constraint is the analog of the constraint in (\ref{eq:(D.44)}).

\paragraph{\it Step 4:} This step describes the fourth constraint on
\(\grm_{T*}\). This constraint on \(\grm_{T*}\) specifies its pull back to \([-96,
-94] \times \bigcup_{(\gamma, \textsc{z}_\gamma)\in\Theta}\mathcal{T}_{\gamma}\)
via the embedding from the second bullet of (\ref{(A.9a,11)}). The constraint
asks that this pull-back have the form \(ds^{2} + \grg\) with \(\grg\)
denoting a certain \(s\)-dependent metric on \(\bigcup_{(\gamma, \textsc{z}_\gamma)\in\Theta}
\mathcal{T}_{\gamma}\). The upcoming description of \(\grg\) refers to the depiction in (\ref{eq:(D.64)}) of
\(\grg_{-}\) on  \(\bigcup_{(\gamma, \textsc{z}_\gamma)\in\Theta} \mathcal{T}_{\gamma}\); and it refers to an analogous depiction of the metric
\(\grg_{Y}\) on  \(\bigcup_{(\gamma, \textsc{z}_\gamma)\in\Theta}\mathcal{T}_{\gamma}\). The metric \(\grg_{Y}\) on each
\(\mathcal{T}_{\gamma}\) has the same form as (\ref{eq:(D.64)}) but
with \(\grh_{\gamma}  = 0\) and with different versions of
\(\textsc{a}_{\gamma}\) and \(\textsc{b}_{\gamma}\). The \(\grg_{Y}\) versions
of these functions are denoted by
\(\textsc{a}_{Y_{\gamma}}\) and \(\textsc{b}_{Y_{\gamma}}\). Note that
\(\textsc{ a}_{Y_{\gamma}} \geq c_{0}^{-1}\). 

{The} specification of \(\grg\) uses two functions on
\(\mathbb{R}\), the first being the function \(\chi^{\mathcal{T}}_{\diamond 1}\)
given by \(\chi(s + 96)\). This function equals 1 where \(s <
-96\) and it equals 0 where \(s \geq -95\). The second function is
denoted by \(\chi^{\mathcal{T}}_{\diamond 2}\),
it is given by \(\chi(s + 95)\). The latter is equal to 1 where \(s
< -95\) and it is equal to 0 where \(s > -94\). 

The metric \(\grg\) on \(\mathcal{T}_{\gamma}\) is defined
by its volume form, this being \(\textsc{h}_{\gamma} dx dy dt\), and 
by the following Hodge star rules:
\begin{equation}\label{eq:(D.72)}\begin{split}
& *dx\, dy =\\& \quad \big(\chi^{\mathcal{T}}_{\diamond 2}\textsc{a}_{\gamma}+(1 -\chi^{\mathcal{T}}_{\diamond 2})
\textsc{a}_{Y_{\gamma}}\big)dt -\chi^{\mathcal{T}}_{\diamond
1}\textsc{a}_{\gamma}\tau^{-1}(1-\tau)\, 
\textsc{h}_{\gamma}^{-1}\textsc{z}_{\gamma}\grh_{\gamma}dx
+\textsc{b}_{\gamma}\, dy,\\
& *dy\,  dt=\\& \quad  \textsc{h}_{\gamma}^{-1}(1  +\chi^{\mathcal{T}}_{\diamond 1}\tau^{-2}\textsc{h}_{\gamma}^{-1}\textsc{a}_{\gamma}
(1-\tau)^{2}\textsc{z}_{\gamma}\grh_{\gamma}) \, dx
-\chi^{\mathcal{T}}_{\diamond 1}\textsc{a}_{\gamma}\tau^{-1}(1-
\tau)\, 
\textsc{h}_{\gamma}^{-1}\textsc{z}_{\gamma}\grh_{\gamma}dt,  \\
& *dt \, dx =\\& \quad  dy +\textsc{b}_{\gamma}dt.
\end{split}
\end{equation}

Important points to note are that \(\grg\) is independent of \(T\) and \(s\) on a
neighborhood of \(s  = -94\), that \(\grg =\grg_{-}\) on a neighborhood of \(s = -96\) and
that \(\grg =\grg_{-}\) for all \(s\) on the complement of \(\mathcal{T}_{\gamma}'\).

\paragraph{\it Step 5:} This step describes \(\p_{X}\) and says more about
the metric \(\grm_{T*}\). The 2-form
\(\p_{X}\) and the metric \(\grm_{T*}\)
on the \(s \in [-102, -98]\) part of \(X_{*}\) are described
by the analog of Step 1 in Part 11 of Section \ref{sec:Dd)} that
has \(Y\) replacing \(Y_{\G}\). By way of a summary,
\(\p_{X}\) is defined on the \(s \in [-102, -101]\) part of \(X\)
to be the 2-form \(\p_{\mathcal{N}1}\) that
is described in \(Y\)'s version of Step 3 from
Part 9 of Section \ref{sec:Dd)}. The 2-form \(\p_{X}\)
is defined on the \(s \in [-101, -100]\) part of \(X\) to be
\(Y\)'s version of the 2-form \(\p_{\mathcal{N}2 }\) that is described in
Step 4 from Part 9 of Section \ref{sec:Dd)}. The definition of
\(\p_{X}\) on the \(s \in [-100, -98]\) part of \(X\) is made by
specifying its pull back via the embedding from the second bullet of
(\ref{(A.9a,11)}). This pull back is the \(s\)-independent 2-form that equals
\(\p_{0}\) on \(\mathcal{Y}_{0}\) and
\(w_{-} - d({\sigma}_{1}\q_{1-})\) on the rest of \(Y\). The metric
\(\grm_{T*}\) on this part of \(X\) pulls back via
the embedding from the second bullet of (\ref{(A.9a,11)}) as \(ds^{2}
+ \grg\) with \(\grg\) denoting the metric given by \(\grg_{-}\) on
\(\mathcal{Y}_{M}\), the metric in (\ref{eq:(D.39)}) on \([-100, -98]
\times \mathcal{N}_{\varepsilon}\) and the metric
\(\grg_{*}\) on \([-100, -98] \times\mathcal{Y}_{0}\). Note in this regard that
\(\grm_{T*}\) is in any event described by (\ref{eq:(D.10)}).

\paragraph{\it Step 6:} This step describes \(\p_{X}\) and the metric on
the \(s \in [-98, -96]\) part of \(X\). But for one significant
difference, the description of \(\p_{X}\) here is similar to
the description of its namesake given in Step 2 from
Part 11 in Section \ref{sec:Dd)}. Both \(\p_{X}\) and
the metric on this part of \(X\) are described by their pull-backs via the
embedding from the second bullet of (\ref{(A.9a,11)}). The metric pulls back as
\(ds^{2} + \grg\) with \(\grg\) given by \(\grg_{*}\) on \(\mathcal{Y}_{0\varepsilon}\)
and by the metric in (\ref{eq:(D.39)}) on \(\mathcal{N}_{\varepsilon}\). The metric \(\grg\) on
\(\mathcal{Y}_{M}\) is the metric \(\grg_{-}\). 

As in the Step 2 from Part 11 of Section \ref{sec:Dd)}, a
1-form to be denoted by \(\q_{3-}\) is constructed with the
following properties: It obeys \(d\q_{3-} = \p -w_{-} + d(\sigma_{1}
\q_{1-})\), it vanishes on the \(r \geq\rho_{*} -\frac{1}{2}\varepsilon \) part of \(\mathcal{N}_{\varepsilon}\)
and its \(L^{2}\)-norm is bounded by \(c_{0}\). Reintroduce \(\chi_{\diamond 3}\) to
denote the function on \(\mathbb{R}\) given by \(\chi(|s| - 97)\) and use
\(\chi_{\diamond 3}'\) to denote its derivative. The 2-form \(\p_{X}\) on \([-98, -96] \times Y\)
is \(\p_{0}\) on \(\mathcal{Y}_{0}\) and given on the rest of
\(Y\) by the formula in (\ref{eq:(D.45)}). Note that
\(\p_{X}\) is \(\p_{0} + \p\) near \(\{-96\} \times Y\),
and that its \(L^{2}\)-norm on this part of \(X\) is bounded by \(c_{0}\).

To start the description of \(\q_{3-}\), let \(\gamma\)
denote a loop from a pair in the set \(\Theta\). The 2-form
\(w_{-}\) on \(\mathcal{T}_{\gamma}\) is
given by \(\tau w + (1 - \tau)\textsc{z}_{\gamma}\p_{\gamma}\) and so it can be written as
\begin{equation}\label{eq:(D.73)}
\p_{\gamma} +  \tau(\textsc{q}_{\gamma} dt -\textsc{z}_{\gamma}\grq_{\gamma} dx),
\end{equation}
where \(\textsc{q}_{\gamma}\) is a function of \(y\) and \(t\)
whose \(y\)-derivative is \(\textsc{h}_{\gamma}\), and where
\(\grq_{\gamma}\) is a function of \(x\) and \(y\) whose \(y\)-derivative
is \(\grh_{\gamma}\). Meanwhile, \(\tau=\frac{1}{T^2}\). Let \(\q_{\gamma}\) denote \(\tau(\textsc{q}_{\gamma} dt -
\textsc{z}_{\gamma}\grq_{\gamma} dx)\). Use (\ref{eq:(D.64)}) to see that the
\(\q_{\gamma} \wedge *\q_{\gamma}\) can be written as
\(\sigma\textsc{h}_{\gamma} dx\, dy \, dt\) with
\(|\sigma| \leq c_{0}\tau^{2}\textsc{a}_{\gamma}^{-1}\). Now,
\(\textsc{a}_{\gamma}\) is constrained to be positive and
less than \(c_{0}^{-1}\tau^{2}\), and these constraints are met if
\(\textsc{a}_{\gamma}\) is chosen greater than \(c_{0}^{-2}\tau^{2}\). Take
\(\textsc{a}_{\gamma}\) so that this is the case, and then
the \(L^{2}\)-norm (and pointwise norm) of \(\q_{\gamma}\) is bounded by \(c_{0}\). 

The 2-form \(w_{-} - d(\sigma \q_{1-})\) is exact on
\(\mathcal{Y}_{M}-(\bigcup_{(\gamma, \textsc{z}_\gamma)\in\Theta}\mathcal{T}_{\gamma}')\)
and on the \(\ir \leq \rho_{*}\) part of \(\mathcal{N}_{\varepsilon}\). This being the case, it
can be written as \(d\q_{*}\) on this part of \(Y\).
More to the point, Lemma \ref{lem:D.5} can be used as in the last paragraph of
Step 2 from Part 11 in Section \ref{sec:Dd)} to obtain a version
of \(\q_{*}\) that is zero where \(r \geq\rho_{*} -\frac{1}{512}  \varepsilon \) and has \(L^{2}\)-norm bounded by
\(c_{0}\) on \(\mathcal{Y}_{M}-(\bigcup_{(\gamma,
  \textsc{z}_\gamma)\in\Theta}\mathcal{T}_{\gamma}')\) and on the \(r
\leq\rho_{*}-\frac{1}{512}\varepsilon\) part of \(\mathcal{N}_{\varepsilon}\).

Let \(\gamma\) again denote a loop from a pair in \(\Theta\). The
difference \(\q_{*} - \q_{\gamma}\) on \(\mathcal{T}_{\gamma}-\mathcal{T}_{\gamma}'\)
is exact. This being the case, it follows from the Mayer-Vietoris
exact sequence and from the fact that the various loops from \(\Theta\)
freely generate \(H_{1}(M_{\delta};\mathbb{R})\) that there is a closed 1-form, \(\k\), on
\(\mathcal{Y}_{M}\) with the following three properties:
First, \(\q_{*} - \q_{\gamma} = \k\) on each \((\gamma,\textsc{z}_{\gamma})\in\Theta\) version of
\(\mathcal{T}_{\gamma}-\mathcal{T}_{\gamma}'\). Second, \(\k = 0\) near
\(\mathcal{N}_{\varepsilon}\) and on \(\bigcup_{\grp\in\Lambda}\mathcal{H}_\grp\). Finally, the
\(L^{2}\)-norm of \(\k\) is bounded by \(c_{0}\). This understood, the
sought after 1-form \(\q_{3-}\) is defined to be \(\q_{\gamma}\) on
each \((\gamma,\textsc{z}_{\gamma}) \in \Theta \) version of
\(\mathcal{T}_{\gamma}\) and to be \(\q_{*}- \k\) on \(\mathcal{Y}_{M}-\bigcup_{(\gamma,
  \textsc{z}_\gamma)\in\Theta}\mathcal{T}_{\gamma}\).

\paragraph{\it Step 7:} This step describes \(\p_{X}\)
and the metric on the \(s\in [-96, -94]\) part of \(X\). The story with
\(\p_{X}\) is simple: It is the 2-form \(\p_{0}+\p\). The metric on \(X\) is described
by its pull back to \([-96, -94]\times Y\) via the embedding
from the second bullet of (\ref{(A.9a,11)}). In particular, it pulls back as
\(ds^{2}+ \grg\) with \(\grg\) being an \(s\)-dependent metric on \(Y\). The
\(s\)-dependence involves only \(\grg\)'s restriction to \(\bigcup_{(\gamma, \textsc{z}_\gamma)\in\Theta}\mathcal{T}_{\gamma}\)
where it is given in Step 4. The metric \(\grg\) is independent of \(s\) on the rest
of \(Y\). As explained in the next paragraph, this metric on \(X\) is such
that the \(L^{2}\)-norm of \(\p_{X}\) on the \([-96, -94]\) part of \(X\) is bounded by
\(c_{0}\), a \(T\)-independent constant. 

The afore-mentioned \(L^{2}\)-norm 
bound holds for \(\p_{0}\).
To see about \(\p\), write it as \(\sum_{(\gamma,\textsc{z}_\gamma)\in\Theta}\textsc{z}_{\gamma}\p_{\gamma}\). A given version of
\(\p_{\gamma}\) has support in \(\mathcal{T}_{\gamma}\)
where the metric is given by (\ref{eq:(D.72)}). Fix \(s\in [-96, -94]\) and since
\(\p_{\gamma}=h_{\gamma}dx\,  dy\), the first bullet of (\ref{eq:(D.72)}) can be used to write
\(\p_{\gamma}\wedge*\p_{\gamma}\) as \(\textsc{p}|\grh_{\gamma}|^{2}\textsc{h}_{\gamma}dx\,  dy\, dt\) with
\(\textsc{p} =(\chi^{\mathcal{T}}_{\diamond  2}\textsc{a}_{\gamma}+(1-\chi^{\mathcal{T}}_{\diamond 2})\textsc{a}_{Y_{\gamma}})^{2}\).
Since \(\textsc{p} <c_{0}\), so the \(L^{2}\)-norm of \(\p_{\gamma}\)
at any \(s \in [-96, -94]\) slice of \([-96,-94] \times Y\) is bounded
by \(c_{0}\).

\paragraph{\it Step 8:} This last step describes \(\p_{X}\) and the metric
on the \(s \in [-94, 102]\) part of \(X\). The description of
\(\p_{X}\) starts where \(s \in [96, 102]\). The 2-form
\(\p_{X}\) here is described by the \(Y_{+} =
Y_{\G}\) version of the 2-form that is defined in Steps 1 and
2 from Part 11 of Section \ref{sec:Dd)}. The \(s \in [96,
100]\) part of the constraint in (\ref{eq:(D.70)}) and the constraint
from Part 10 of Section \ref{sec:Dd)} are needed to repeat Steps
1 and 2 from Part 11 in the case at hand. These
steps define a version of \(\p_{X}\) whose \(L^{2}\)-norm on the \(s
\in [96, 102]\) part of \(X\) is bounded by \(c_{0}\) times the \(L^{2}\) norm of
\(w_{+}\) on \(Y_{\G}\). This version of \(\p_{X}\) is \(w_{+}\) near the \(s = 102\) slice of
\(X_{*}\) and it is the 2-form \(\p_{0} + \p\) near
the \(s = 96\) slice. The 2-form \(\p_{X}\) is set equal to
\(\p_{0} + \p\) on the \(s \in [-94, 96]\) part of \(X\). Its
\(L^{2}\)-norm on the \(s \in [-94, 96]\) part of \(X\) is
bounded by \(c_{0}\), this being a consequence of (\ref{eq:(D.71)}).

\paragraph{Part 9:} Taking up where Part 8 left off, this last part
of the subsection defines the desired metric on \(X\) and
\(w_{X}\) on the \(s \in [-\frac{1}{2}L + 4, L]\) part of \(X\). To this end, fix \(T\) large and then
\(L_{1} \geq c_{T}\) so as to use the constructions in Part 8 of the metric
\(\grm_{T*}\) and \(\omega_{T*}\). With
\(L_{1}\) chosen, assume that \(L > 4L_{1}\). The metric \(\grm_{T*}\)
where \(s \in [-\frac{1}{2} L + 4, -\frac{1}{2}L + 5]\) is the same as the \(s \in [-  
\frac{1}{2} L + 4, -\frac{1}{2} L + 5]\) version of the metric from Part 8; and
\(\omega_{T*}\) on this same part of \(X\)
is the \(s \in [-\frac{1}{2} L + 4, - \frac{1}{2} L + 5]\) version of Part 8's 2-form \(w_{X}\).
This understood, the desired metric for \(X\) is taken to be
\(\grm_{T*}\) where \(s \geq - \frac{1}{2} L + 4\), and the 2-form \(w_{X}\) is taken to be
\(\omega_{T*}\) on this same part of \(X\). Here, \(\upsilon_X\) is
set to be the \(s\)-independent 1-form \(\upsilon_\diamond\); and the
bounds in items 4b) and 5c) of (\ref{eq:(A.16)}) are verified by choosing
the parameter \(m\) to be sufficiently small, as directed in Part 2 above.
\epf







\begin{thebibliography}{KM}

 \bibitem[APS]{APS} M.F. Atiyah, V. Patodi and I.M. Singer, {\em Spectral asymmetry in Riemannian geometry I}, Math. Proc. Cambridge Philos. Soc. 77 (1975), 43--69.
\bibitem[BT]{BT}R. Bott, L. Tu, Differential forms in algebraic topology. Springer, 1982.
\bibitem[BMO]{BMO}J. Bloom, T. Mrowka, P. Ozsvath, {\em A K\"unneth
    formula for monopole Floer homology}, in preparation. 
\bibitem[C]{C} H. Cartan, {\em La transgression dans un groupe de Lie
    et dans un fibr\'e principal}, in Colloque de topologie, Belgium, 1951, 57-71.

 
\bibitem[CGHH]{CGHH}V. Colin, P. Ghiggini, K. Honda, M. Hutchings,
  {\em Sutures and contact homology I }, {\tt arXiv:1004.2942}

\bibitem[D]{D}S. Donaldson, Floer homology groups in Yang-Mills
  theory, Cambridge University Press, 2004.
\bibitem[D1]{Ds}S.Donaldson, {\em The Seiberg-Witten equations and
    4-manifold topology}. Bull. Amer. Math. Soc. {\bf 33} (1996), no. 1, 45-70.
\bibitem[DK]{DK}S. Donaldson and P. Kronheimer, The Geometry of
Four-Manifolds, Oxford University Press, 1990.
\bibitem[F]{F}K. Fukaya, {\em Floer homology for connected sums of
    homology 3-spheres}, Topology {\bf 35} (1996), 89-136. 
\bibitem[GKM]{GKM}M. Goresky, R. Kottwitz, and R. MacPherson, {\em Equivariant cohomology, Koszul
duality, and the localization theorem}, Invent. Math. {\bf 131} (1998), 25-83.
\bibitem[J]{J}J. Jones, {\em Cyclic homology and equivariant homology.}
 Inventiones mathematicae {\bf 87}, 403-424.

\bibitem[Ju]{Ju}A. Juhasz, {\em Holomorphic discs and sutured
    manifolds.} Algebr. Geom. Topol. {\bf 6} (2006), 1429-1457.

\bibitem[Ho]{Ho}K. Honda, {\em Transversality theorems for harmonic
    forms.} Rocky Mountain J. Math. {\bf 34} (2004), no. 2, 629-664.





\bibitem[HT2]{HT2} M. Hutchings and C. H. Taubes, \textit{Gluing
pseudoholomorphic curves along branched covered cylinders II};
Jour. Sympl. Geom., \textbf{7} (2009) 29-133.

\bibitem[Ko]{Ko}J.-L. Koszul, {\em Sur un type d'alg\`ebres diff\'erentielles en rapport
avec la transgression} in Colloque de topologie, Belgium, 1951, 73-81.


\bibitem[KLT1]{KLT1} C. Kutluhan, Y.-J. Lee  and C. H. Taubes, {\it
    \(HF=HM \) \(I:\) Heegaard Floer 
homology and Seiberg-Witten Floer homology}; arXiv:1007.1979.

\bibitem[KLT2]{KLT2} C. Kutluhan, Y.-J. Lee  and C. H. Taubes, {\it
    \(HF=HM \) \(II:\) Reeb orbits and holomorphic curves for the ech/Heegard-Floer correspondence}; arXiv:1008.1595.

\bibitem[KLT3]{KLT3} C. Kutluhan, Y.-J. Lee  and C. H. Taubes, {\it \(HF=HM\) \(III:\) Holomorphic curves and the differential for the ech/Heegard-Floer correspondence}; arXiv:1010.3456.

\bibitem[KLT4]{KLT4} C. Kutluhan, Y.-J. Lee  and C. H. Taubes, {\it
    \(HF=HM\) \(IV:\) The Seiberg-Witten Floer homology and ech
    correspondence}; arXiv:111107.2297.

\bibitem[KMan]{KMan} P. Kronheimer, C. Manolescu, {\em Periodic Floer
    pro-spectra from the Seiberg-Witten equations}, arXiv:0203243. 

\bibitem[KM]{KM} P. Kronheimer and T. Mrowka, Monopoles and Three-Manifolds, Cambridge University Press, Cambridge UK, 2007.

\bibitem[KM1]{KMs} P. Kronheimer and T. Mrowka, {\em Knots, sutures,
    and excision.} J. Differential Geom. {\bf 84} (2010), no. 2, 301-364.

\bibitem[Lef]{Le} K. Lef\`evre-Hasegawa, {\em Sur les \(A^\infty\)-cat\`egories},
  Ph.D thesis, l'Universit\`e Paris 7, 2003.

\bibitem[L]{L} Y.-J. Lee, {\em Heegaard splittings and Seiberg-Witten
    monopoles. Geometry and topology of manifolds}, 173-202, Fields
  Inst. Commun., {\bf 47}, 2005.

\bibitem[L0]{L0} Y.-J. Lee, \textit{An extension of Taubes'
    convergence theorem}, unpublished manuscript, 2006. Superseded by
  \textit{From Seiberg-Witten to Gromov: MCE and singular symplectic
    forms}, preprint 2015.

\bibitem[L1]{L1} Y.-J. Lee, \textit{Reidemeister torsion in
    Floer-Novikov theory and counting pseudo-holomorphic tori, I.}
  Jour. Sympl. Geom. \textbf{3} (2005) 221-311.

\bibitem[LT]{LT} Y.-J. Lee and C.H. Taubes, {\it Periodic Floer
    homology and Seiberg-Witten Floer cohomology}, {\em
    J. Symp. Geom.} \textbf{10}, (2012) 81-164.

\bibitem[M]{M} C. Manolescu, {\it Seiberg-Witten-Floer stable homotopy
    type of three-manifolds with $b^1=0$}, Geom. Topol. \textbf{7}
  (2003) 889-932.


\bibitem[OS1]{OS1} P. Ozsv\'ath; Z. Szab\'o, {\em Holomorphic disks and
    topological invariants for closed three-manifolds}. Ann. of
  Math. (2) {\bf 159} (2004), no. 3, 1027-1158.I.

\bibitem[OS2]{OS2} P. Ozsv\'ath; Z. Szab\'o, {\em Holomorphic disks and
    three-manifold invariants: properties and applications}. Ann. of
  Math. (2) {\bf 159} (2004), no. 3, 1159-1245. 

\bibitem[T]{T} C. H. Taubes, {\em Seiberg-Witten
and Gromov Invariants for Symplectic 4-manifolds},
International Press, Somerville MA, 2005.

 \bibitem[T1]{T1} C. H. Taubes, \textit{Embedded contact homology and
Seiberg-Witten Floer cohomology I}, Geometry and Topology 14 (2010) 2497-2581. 

\bibitem[T2]{T2} C. H. Taubes, \textit{Embedded contact homology and
Seiberg-Witten Floer cohomology II}, 

\bibitem[T3]{T3} C. H. Taubes, \textit{Embedded contact homology and
Seiberg-Witten Floer cohomology III}, Geometry and Topology 14 (2010) 2721-2817.

\bibitem[T4]{T4} C. H. Taubes, \textit{Embedded contact homology and
Seiberg-Witten Floer cohomology IV}, Geometry and Topology 14 (2010) 2819-2960.

\bibitem[T5]{T5} C. H. Taubes, \textit{Embedded contact homology and
Seiberg-Witten Floer cohomology V}, Geometry and Topology 14 (2010)
2961-3000.

\bibitem[Ta1]{Ta1} C. H. Taubes, \(\op{SW}\Rightarrow  \op{Gr}\):
  {\em from  the Seiberg-Witten equations to pseudo-holomorphic
    curves.} J. Amer. Math. Soc. {\bf 9} (1996), no. 3, 845-918.

\bibitem[T:W1]{TW1} C. H. Taubes, \textit{The Seiberg-Witten equations and
the Weinstein conjecture}, Geometry and Topology, \textbf{11} (2007)
2117-2202. 

\bibitem[T:W2]{TW2} C. H. Taubes, \textit{The Seiberg-Witten equations and the
Weinstein conjecture II: More closed integral curves for the Reeb
vector field, }Geometry and Topology, \textbf{13} (2009) 1337-1417.

\bibitem[T:W3]{Tasf} C. H. Taubes, \textit{Asymptotic spectral flow for Dirac
operators}, Commun. Anal. and Geom. \textbf{15} (2007) 569-587.

\bibitem[TT]{TT} D. Toledo, Y.-L. Tong, {\em Duality and intersection
    theory in complex manifolds. I}, Mathematische Annalen 
 {\bf 237}, (1978) 41-77.
\end{thebibliography}
\end{document}